\documentclass{acta_2022}

\usepackage{amssymb,amsfonts,psfrag}
\usepackage{amsmath}
\usepackage{titlesec}
\usepackage{enumitem}
\usepackage{xspace}
\usepackage{multicol,multirow}
\usepackage{comment}
\usepackage{bm}
\usepackage{mathrsfs,mathtools}
\usepackage{xparse}
\usepackage{tikz}
\usepackage{gnuplot-lua-tikz}
\usepackage{hhline}
\usepackage{scalerel}
\usepackage{imakeidx}
\makeindex[columns=1, title=Index,intoc]

\usepackage[percent]{overpic}

\setlist[itemize]{leftmargin=*}
\usepackage[longnamesfirst]{natbib} \usepackage{har2nat}
\setcitestyle{square,comma,aysep={}}
\shortcites{samplepreprint} %\shortcites{} %for 5+ authors
\usepackage{amssymb,amsmath}
\usepackage{newtxtext}
\usepackage{newtxmath}
\DeclareSymbolFont{CMlargesymbols}{OMX}{cmex}{m}{n}
\DeclareMathDelimiter{(}{\mathopen} {operators}{"28}{CMlargesymbols}{"00}
\DeclareMathDelimiter{)}{\mathclose}{operators}{"29}{CMlargesymbols}{"01}
\pagerange{\pageref{firstpage}--\pageref{lastpage}}
\allowdisplaybreaks
\doi{XXXXXXXX}  %GES% 8 figures, starting with 1
\usepackage[UKenglish]{babel} %GS%
\newcommand{\ignore}[1]{}
\usepackage[twoside,
	paperheight=247mm,
	paperwidth=174mm,
	marginratio={3:5,2:3},
	left=18mm,
	top=20mm]{geometry}
\usepackage[hang]{subfigure}
\numberwithin{figure}{section}
\numberwithin{table}{section}
\usepackage{subfigure}
\usepackage{graphics,graphicx}
\usepackage{xcolor}
\usepackage{enumitem}
\setenumerate[0]{label=(\roman*)}

\usepackage{macros}

\usepackage{soul}
\usepackage[normalem]{ulem}
\newcommand{\modif}[1]{ {\textcolor{black}{#1} } }

\usepackage{todonotes}

\title[AFEM]{Adaptive Finite Element Methods}
\author[]{Andrea Bonito\\
Department of Mathematics, Texas A\&M University,\\ College Station, Texas 77843, USA\\
E-mail: {bonito@tamu.edu}\\
\and 
Claudio Canuto\\
Dipartimento di Scienze Matematiche, Politecnico di Torino, \\
Corso Duca degli Abruzzi 24, 10129 Torino, Italy\\
E-mail: {claudio.canuto@polito.it}\\
\and
Ricardo H. Nochetto\\
Department of Mathematics and Institute for Physical Science and Technology, \\
University of Maryland, College Park, MD 20742, USA\\
E-mail: {rhn@umd.edu}\\
\and
Andreas Veeser\\
Dipartimento di Matematica, Universit\`a degli Studi di Milano,\\ Via Saldini 50, 20133 Milano, Italy\\
E-mail: {andreas.veeser@unimi.it}
}
\begin{document}

\maketitle
\vspace{2.5in}
\begin{abstract}
This is a survey on the theory of adaptive finite element methods
(AFEMs), which are fundamental in modern computational science and
engineering but whose mathematical assessment is a formidable challenge.
We present a self-contained and up-to-date discussion
of AFEMs for linear second order elliptic PDEs and dimension $d>1$, 
with emphasis on foundational issues. After a brief review of functional
analysis and basic finite element theory, including piecewise polynomial
approximation in graded meshes, we present the core material for coercive problems.
We start with a novel a posteriori error analysis applicable to rough data,
which delivers estimators fully equivalent to the solution error.
They are used in the design and study of three AFEMs depending on the structure
of data. We prove linear convergence
of these algorithms and rate-optimality provided the solution and data
belong to suitable approximation classes. We also address the relation between
approximation and regularity classes. We finally extend this theory to
discontinuous Galerkin methods as prototypes of non-conforming AFEMs
and beyond coercive problems to inf-sup stable AFEMs.
\end{abstract}

\tableofcontents 

%--------------------------------------------------------------------------------

\section{Introduction: Overview of AFEMs}\label{S:intro}

This is a survey on the theory of adaptive finite element methods
(AFEMs), which are fundamental in modern computational science and
engineering. We present a self-contained and up-to-date discussion
of AFEMs for linear second order elliptic PDE in dimension $d>1$, 
with emphasis on foundational issues rather than applications of
AFEMs. This paper builds on and expands the older surveys
\cite{NoSiVe:09,NochettoVeeser:2012}. In fact, we decided to incorporate 
several new aspects to the theory described below. 

The paper develops the theory of AFEMs gradually and is meant to be
accessible to advanced students and researchers interested in learning
the fundamental aspects of adaptivity and why AFEMs outperform
classical FEMs. We quantify the superior performance of AFEMs with
precise mathematical statements rather than simulations. We present
very few numerical experiments to illustrate some key (and new) algoritmic
ideas and methods, but the paper is otherwise a tour in the numerical analysis
of  adaptive approximation of linear elliptic PDEs.

By design this paper goes deep into some foundational aspects of AFEMs
theory, provides full discussions and proofs, as well as pointers to the
main literature. We consider the following model problem on a polyhedral domain
$\Omega\subset\R^d$ with $d\ge2$ \looseness=-1
\begin{equation}\label{E:model-PDE}
L[u] := -\div (\bA\nabla u) + cu = f
\end{equation}
with general variable coefficients $(\bA,c)$, forcing $f\in H^{-1}(\Omega)$
and homogeneous Dirichlet boundary conditions $u=0$ on $\partial\Omega$ mostly, but not
exclusively. If $\V:=H^1_0(\Omega)$ and $\B:\V\times\V\to\R$ is the bilinear form associated
with \eqref{E:model-PDE}, the weak form reads:
\begin{equation}\label{E:weak-formulation}
  u\in\V: \quad
  \B[u,v] = \langle f, v \rangle \quad\forall \, v\in \V.
\end{equation}
Given a conforming and shape regular partition $\grid$ of $\Omega$, created by successive
refinement of a coarse mesh $\grid_0$, let $\V_\grid$ denote the space of continuous piecewise
polynomial functions of degree $n\ge1$ over $\grid$ vanishing on $\partial\Omega$. The Galerkin approximation $u_\grid$ of $u$ solves
\begin{equation}\label{E:galerkin-problem}
 u_\grid\in\V_\grid: \quad
 \B[u_\grid,v] = \langle f, v \rangle \quad\forall \, v\in \V_\grid.
\end{equation}
This is a conforming approximation because $\V_\grid\subset\V$.
The aim of this paper is to:

\begin{itemize}
\item
Design and analyze practical ways to
estimate the error $|u-u_\grid|_{H^1_0(\Omega)}:=\|\nabla(u-u_\grid)\|_{L^2(\Omega)}$
in terms of so-called a posteriori error
estimators, which are computable quantities depending on the discrete solution $u_\grid$
and data $\data = (\bA,c,f)$. \looseness=-1

\item
Design adaptive algorithms that equidistribute the local errors $\|\nabla(u-u_\grid)\|_{L^2(T)}$
for all elements $T\in\grid$, thereby optimizing the computational effort; this is a key step
that makes complex 3d situations accessible computationally.

\item
  Show that this strategy delivers a performance comparable with the best possible in terms of
  degrees of freedom, which is a measure of computational complexity. This is a delicate matter
  because it entails
  dealing with approximation classes and their relation to regularity classes in terms of Besov 
  and Lipschitz spaces.

\item
  Present and analyze the bisection method for mesh refinement, one of the most versatile techniques
  for local mesh refinement that guarantees shape regularity and optimal complexity; the latter is
  instrumental for the previous point. Our study includes conforming meshes as well as
  certain non-conforming meshes.

\item
  Extend the theory to a range of important problems that fail to be conforming or coercive.
  The first class is discontinuous Galerkin methods and the second one is inf-sup stable FEMs.
  The former is a notorious example of non-conforming approximation, whereas the second is
  non-coercive.
\end{itemize}
In achieving these goals we provide several new ideas and methods. We also refer to the pertinent
literature but we do not give a full list of references nor get into comparisons of various
approaches. It is not our intention to be comprehensive but rather to cover basic aspects of
adaptivity in depth at the expense of important topics we do not touch upon. Some of them are:
\begin{itemize}
\item
Adaptive eigenvalue approximation

\item
Goal oriented error analysis

\item
Non-conforming discretizations (except for discontinuous Galerkin)

\item 
Coarsening or aggregation

\item 
Anisotropic refinements

\item 
$hp$-adaptivity

\item 
Tree approximation

\item
Other PDEs: convection-diffusion equations, nonlinear and evolution equations.
\end{itemize}

We devote the rest of this introduction to provide a roadmap to the rest of the paper.
In doing so, we introduce notation that will be used later and present some topics in
their most primitive form to render an early idea about how they fit and interrelate.

%-----------------------------------------------------------------------------------
\paragraph{A posteriori error analysis.}
%-----------------------------------------------------------------------------------
%
We refer to the books \cite{AinsworthOden:00, Verfuerth:13} for the classical theory.
However, in contrast to  most of the existing literature, the current
theory deals with forcing $f \in H^{-1}(\Omega)$. This allows for rough data useful in
applications, such as line Dirac masses, but also encompasses a new approach to error estimation
that leads to error-dominated estimator and oscillation and prevents error overestimation;
this extends \cite{KreuzerVeeser:2021} to \eqref{E:model-PDE} and polynomial degree $n\ge1$.
The new twist is the construction of a projection
operator $P_\grid:H^{-1}(\Omega) \to \F_\grid$ into a space of piecewise polynomials in $\grid$
and on its skeleton $\faces$, namely the set of all internal faces. Such an operator happens to
be locally stable on stars (or patches) $\omega_z$ of $\grid$ for all vertices $z\in\vertices$ of $\grid$:
\begin{equation}\label{E:Pgrid-local-stability}
  \| P_\grid \ell \|_{H^{-1}(\omega_z)} \le \Cstab \|\ell\|_{H^{-1}(\omega_z)}
  \quad\forall \, \ell \in H^{-1}(\omega_z).
\end{equation}
An important property of $P_\grid$ and its range $\F_\grid$ is that for piecewise polynomial
coefficients $(\bA,c)$, or in short discrete coefficients, $P_\grid$ is invariant in the subspace
$L[\V_\grid]$ of $H^{-1}(\Omega)$ or equivalently
\begin{equation}\label{E:range-projection}
P_\grid \big( L[v] \big) = L[v] \in\F_\grid\quad\forall\, v\in \V_\grid.
\end{equation}
It is worth realizing that $L[v]$ is made of two distinct parts. The first one is
absolutely continuous relative to the Lebesgue
measure, namely $-\div(\bA\nabla v) + cv$ in every element $T\in\grid$. The second part is singular 
and supported in the skeleton $\faces$, namely $\jump{\bA\nabla v}\cdot\bn|_F\delta_F$ for every face
$F\in\faces$, where $\jump{\cdot}$ is the jump across $F$, $\bn$ is a unit normal to $F$, and $\delta_F$
is the Dirac mass on $F$.

These two properties of $P_\grid$ have the following crucial consequences. Let
\begin{equation}\label{E:residual}
  R_\grid := L[u-u_\grid]=f-L[u_\grid] \in H^{-1}(\Omega)
\end{equation}
be the residual of the Galerkin approximation of \eqref{E:weak-formulation}. Using \eqref{E:range-projection}
yields
\begin{equation*}
R_\grid - P_\grid R_\grid = f - P_\grid f.
\end{equation*}
This shows that $R_\grid$ decomposes into a discrete, thus finite dimensional and computable, PDE part
$P_\grid R_\grid = P_\grid f - L[u_\grid]$ and an infinite dimensional component $f-P_\grid f$, the
so-called data oscillation that depends on $f$ and can only be evaluated with additional knowledge of $f$.

The nonlocal $H^{-1}$-norm of $R_\grid$ splits into local contributions on stars, whence \footnote{Throughout this work,  $A \lesssim B$ signifies $A \leq C B$ with a constant $C$ independent of the discretization parameters, and $A \approx B$ stands for $A \lesssim B$ and $B \lesssim A$.}
\begin{equation*}
  |u-u_\grid|_{H^1_0(\Omega)}^2 \approx \|R_\grid\|_{H^{-1}(\Omega)}^2
  \approx \sum_{z\in\vertices} \|R_\grid\|_{H^{-1}(\omega_z)}^2.
\end{equation*}
The discrete nature of $P_\grid R_\grid$ allows us to derive a computable $L^2$-weighted PDE estimator
$\eta_\grid(u_\grid,z)$ equivalent to $\|R_\grid\|_{H^{-1}(\omega_z)}$, which together with the remaining
data oscillation $\osc_\grid(f,z)_{-1} := \|f-P_\grid f\|_{H^{-1}(\omega_z)}$ gives the upper bound
\begin{equation*}
|u-u_\grid|_{H^1_0(\Omega)}^2 \lesssim \sum_{z\in\vertices} \Big(\eta_\grid(u_\grid,z)^2 + \osc_\grid(f,z)_{-1}^2 \Big).
\end{equation*}
It turns out that this estimate is sharp or, in other words, that there is no overestimation of
the error. To see this important and unique property of these new estimators, we invoke
\eqref{E:Pgrid-local-stability} to write the local lower bounds
\begin{align*}
  \eta_\grid(u_\grid,z) &\approx \|P_\grid R_\grid\|_{H^{-1}(\omega_z)} \le \Cstab \|R_\grid\|_{H^{-1}(\omega_z)},
    \\
  \osc_\grid(f,z)_{-1} & = \|R_\grid - P_\grid R_\grid\|_{H^{-1}(\omega_z)} \le (1+\Cstab) \|R_\grid\|_{H^{-1}(\omega_z)}.
\end{align*}

Section \ref{S:aposteriori} constructs the operator $P_\grid$, and derives several important properties
such as its local quasi-best approximation and the above error-dominated a posteriori
bounds. The former guarantees the inequality for the local $L^2$-projection $\Pi_\grid$ \looseness=-1
\begin{equation*}
  \| f - P_\grid f\|_{H^{-1}(\omega_z)} \lesssim \| f - \Pi_\grid f\|_{H^{-1}(\omega_z)}
  \lesssim \| h(f - \Pi_\grid f) \|_{L^2(\omega_z)},
\end{equation*}
which is the typical form of data oscillation provided $f\in L^2(\Omega)$.
However, this $L^2$-weighted oscillation is not bounded above by the error and is thus responsible
for potential overestimation. Section \ref{S:aposteriori} proves further properties of $\eta_\grid(u_\grid)$
such as its reduction upon refinement and its localized discrete upper bound, as well as
quasi-monotonicity of $\osc_\grid(f)_{-1}$ upon refinement. These properties, known for the standard
$L^2$-weighted estimator and oscillation, are thus retained by the new construction. \looseness=-1

Section \ref{S:aposteriori} also deals with the alternative error estimators that result from
solving local problems,
using hierarchy, or imposing flux equilibration. We show that all of them lead, essentially, to estimators equivalent to $\|R_\grid\|_{H^{-1}(\omega_z)}$. Moreover, we present an optimal framework to deal
with non-homogeneous Dirichlet boundary conditions as well as with Robin and Neumann boundary conditions.

%-----------------------------------------------------------------------------------
\paragraph{Linear convergence of AFEMs.}
%-----------------------------------------------------------------------------------  
%
%
Local a posteriori error indicators are usually employed to mark elements (or set of
elements) with largest indicators for refinement. We are concerned with the most popular
{\it D\"orfler marking} (or bulk chasing): given a parameter $\theta\in(0,1]$,
select a set $\marked \subset\grid$ such that
\begin{equation}\label{E:doerfler}
\eta_\grid(u_\grid,\marked) \ge \theta \, \eta_\grid(u_\grid);
\end{equation}
hereafter we define
  $\eta_\grid(u_\grid,\marked)^2 := \sum_{T\in\marked} \eta_\grid(u_\grid,T)^2$ where
  $\eta_\grid(u_\grid,T)$ is the PDE indicator associated with a generic element $T\in\grid$.
Note that $\theta=1$ corresponds to uniform refinement. In Section \ref{S:convergence-coercive},
we present three AFEMs in increasing
order of complexity regarding data $\data = (\bA,c,f)$ and prove their linear convergence.

The simplest algorithm, so-called
$\GALERKIN$, works for discrete $\data$ and is the usual adaptive loop
\[
\SOLVE \to \ESTIMATE \to \MARK \to \REFINE.
\]
We assume $\SOLVE$ computes the exact Galerkin solution $u_\grid$, so we refrain from addressing
linear algebra issues. The module $\ESTIMATE$ computes the a posteriori error indicator and
the module $\MARK$ implements \eqref{E:doerfler}; in most of the paper we deal with weighted
$L^2$-error indicators $\eta_\grid(u_\grid)$ but we also address linear convergence
for alternative estimators. The module $\REFINE$ bisects marked elements
and perhaps a few more to keep meshes conforming (or $\Lambda$-admissible if they
are non-conforming). We denote by $\enorm{u-u_\grid}$ the energy error associated with the bilinear
form $\B$. This error is monotone with refinement but may stagnate. 
We thus exploit the estimator reduction property with refinement, typical of $\eta_\grid(u_\grid)$,
to show that the combined quantity
\begin{equation}\label{E:total-error}
\zeta_\grid(u_\grid)^2 := \enorm{u-u_\grid}^2 + \gamma \, \eta_\grid(u_\grid)^2
\end{equation}
contracts in every iteration of $\GALERKIN$ for a suitable scaling parameter $\gamma>0$.
This readily leads to linear convergence of
both $|u-u_\grid|_{H^1_0(\Omega)}$ and $\eta_\grid(u_\grid)$.

We next keep the coefficients $(\bA,c)$ discrete but allow for a general $f\in H^{-1}(\Omega)$.
This is to prevent the {\it multiplicative interaction} between $(\bA,c)$ and $u$ that occurs in
\eqref{E:model-PDE} if we were to approximate $(\bA,c)$. In contrast, the effect of $f$ is
{\it linear} in \eqref{E:model-PDE}. 
We show examples where $\enorm{u-u_\grid}$ may stagnate because the adaptive process is dominated by
oscillation $\osc_\grid(f)_{-1}$ (preasymptotic regime). To compensate for this fact, we design
a one-step \AFEM with switch as in \cite{Kreuzer.Veeser.Zanotti}, the so-called $\AFEMSW$, that proceeds as $\GALERKIN$ provided
$\eta_\grid(u_\grid)$ dominates and otherwise reduces $\osc_\grid(f)_{-1}$ separately. We show
that, for suitable a parameter $\gamma>0$, the combined quantity \looseness=-1
\begin{equation}\label{E:total-error-ws}
  \zeta_\grid(u_\grid)^2 := \enorm{u-u_\grid}^2 + \gamma \eta_\grid(u_\grid)^2
  +  \osc_\grid(f)_{-1}^2
\end{equation}
contracts in every loop of $\AFEMSW$. This yields linear convergence of the error $|u-u_\grid|_{H1_0(\Omega)}$
and the estimator $\est_\grid(u_\grid,f) = \eta_\grid(u_\grid) + \osc_\grid(f)_{-1}$.

The third algorithm is the two-step $\AFEM$, the so-called $\AFEMTS$,
which allows for general data $\data=(\bA,c,f)$. To handle the
aforementioned nonlinear effect of $(\bA,c)$ and also deal with general $f$, all data $\data$
are first approximated by a routine $\DATA$ to a desired level of accuracy, which is adjusted 
at every step of $\AFEMTS$,
and then fed to $\GALERKIN$ which handles discrete data. Suitably combining the accuracies of each
intermediate module leads to linear convergence and optimal complexity of $\GALERKIN$ within
each loop. The structure of $\AFEMTS$ is flexible enough to easily handle discontinuous coefficients
$(\bA,c)$ with discontinuities that may not be aligned with the mesh. This is because the
approximation of $(\bA,c)$ by discontinuous polynomials takes place in $L^p(\Omega)$ for $p < \infty$.

It is worth stressing two important points. First the approximation of $\data$ is carried out by
a $\GREEDY$ algorithm, which is shown to perform optimally starting from any refinement of
$\grid_0$. Second, the discontinuous piecewise polynomial approximations $(\wh{\bA},\wh{c})$ of
$(\bA,c)$ may not respect, for polynomial degree $\ge 1$, the positivity bounds associated with the coefficients. This requires a nonlinear correction of the output of $\GREEDY$ that restores
positivity and does not deteriorate accuracy beyond a modest multiplicative constant. We
postpone the discussion of these two delicate and technical processes to Section \ref{S:data-approx}, 
which can be omitted in a first reading. \looseness=-1

%-----------------------------------------------------------------------------------
\paragraph{Rate-optimality of AFEMs.}
%-----------------------------------------------------------------------------------
%
Showing that AFEMs outperform classical FEMs is a difficult but important matter. 
This reduces to proving
a superior relation between the required degrees of freedom (or number of elements)
for a desired accuracy; the former is in fact an acceptable
measure of complexity. Showing that AFEMs deliver a performance comparable with the best
entails the following basic ingredients:

\begin{itemize}
\item
{\it Nonlinear approximation classes}: they classify functions in terms of
the best possible algebraic decay rate of approximation $e_N(v)_X$ of a given function $v$ in a given
norm $X$ with $N$ number of elements; roughly speaking, we say $v\in\As$ if
$e_N(v)_X \lesssim N^{-s}$. These classes are related to regularity classes (Sobolev, Besov,
and Lipschitz) along Sobolev embedding lines.

\item
{\it D\"orfler marking}: If the oscillation $\osc_\grid(f)_{-1}$ is dominated by the
PDE estimator $\eta_\grid(u_\grid)$ for a given mesh $\grid$,
then any conforming refinement $\grid_* \ge \grid$ of $\grid$ that reduces $\eta_\grid(u_\grid)$
by a substantial amount induces a refined set $\mathcal{R}:=\grid\setminus\grid_*$ to modify $\grid$ into $\grid_*$
satisfying \eqref{E:doerfler}, namely $\eta_\grid(u_\grid,\mathcal{R})\ge\theta\,\eta_\grid(u_\grid)$.

\item
{\it Minimality of $\marked$}: If the subset $\marked \subset \grid$ in \eqref{E:doerfler} is minimal,
then the cardinality of $\marked$ compares favorably to the cardinality of the best mesh with a comparable error accuracy, thereby leading to rate-optimality of AFEM.
\end{itemize}

This is altogether the topic of Section \ref{S:conv-rates-coercive}.
It is important to notice that membership in $\As$ is never used explicitly by $\AFEM$ to learn
about problem \eqref{E:model-PDE} and improve its resolution. The fundamental reason behind
the superior performance of $\AFEM$ relative to \textsf{FEM} lies in nonlinear approximation
theory. We illustrate now this point with the following insightful approximation example for $d=1$ and
$X=L^\infty(0,1)$ \cite{DeVore:98,Kahane:61}: let $\Omega=(0,1)$, $\gridN=\{[x_{j-1},x_{j}]\}_{j=1}^{N}$ be a
partition of $\Omega$, with
\[
0=x_0 < x_1 < \cdots < x_j < \cdots < x_N = 1,
\]
and let $v:\Omega\to\R$ be an absolutely continuous function to be approximated
by a piecewise constant function $v_N$ over $\gridN$.
To quantify the difference between $v$ and $v_N$ we resort to the
{\it maximum norm} and study two cases depending on the regularity of $v$. We define
$v_N(x):=v(x_{j-1})$ for all $x_{j-1}\le x<x_j$ and note that
\[
| v(x)-v_N(x) | = |v(x) - v(x_{j-1})| \le \int_{x_{j-1}}^x |v'(t)| dt.
\]

\begin{itemize}
\item {\it Case 1: $W^1_\infty$-Regularity}. 
If $u\in W^1_\infty(0,1)$ and $x_{j-1}\le x<x_j$, then
\[
|v(x)-v_N(x) | \le h_j \|v'\|_{L^\infty(x_{j-1},x_j)}
\quad\Rightarrow\quad
\|v - v_N\|_{L^\infty(\Omega)} \leq \frac{1}{N} \|v'\|_{L^\infty(\Omega)}
\]
for a {\it uniform} mesh. We thus deduce a rate $N^{-1}$ using the
same integrability $L^\infty$ on both sides of the error estimate.

\item
{\it Case 2: $W^1_1$-Regularity}. Let $\|v'\|_{L^1(\Omega)}=1$ and $\grid_N$ be a {\it graded} partition so that
$
\int_{x_{j-1}}^{x_j} |v'(t)| dt = \frac{1}{N}.
$
Then, for $x\in[x_{j-1},x_j]$,
\[
|v(x) - v(x_{j-1})| \le \int_{x_{j-1}}^{x_j} |v'(t)| dt = \frac{1}{N}
\quad\Rightarrow\quad
\|v-v_N\|_{L^\infty(\Omega)} \leq \frac{1}{N} \|v'\|_{L^1(\Omega)}.
\]
\end{itemize}
We thus conclude that we could achieve
the same rate of convergence $N^{-1}$ for rougher functions with just
$\|u'\|_{L^1(\Omega)}<\infty$. Three comments are now in order. First, the contrast between Cases 1 and 2
is more dramatic for $v(x) = x^{\alpha}$ with $\alpha\in (0,1)$ because Case 1 only yields the suboptimal
rate $\|v - v_N\|_{L^\infty(\Omega)}\le N^{-\alpha}$.
Second, $\grid_N$ in Case 2 {\it equidistributes} the max-error, a concept that will permeate our
discussions later. Third, the optimal rate of Case 2 is due to the exchange of differentiability with
integrability along the critical Sobolev embedding line between the left and right-hand sides of the error estimate (nonlinear Sobolev scale),
while Case 1 relies on the linear Sobolev scale with constant integrability.

We exploit and further elaborate these concepts in Section \ref{S:conv-rates-coercive} to show
rate-optimality of the three algorithms $\GALERKIN, \AFEMSW, \AFEMTS$, discussed in
Section \ref{S:convergence-coercive}, provided $u$ and $\data$ belong to suitable
approximation classes. We also investigate the relation between these approximation
classes with regularity classes, allowing for discontinuous coefficients,
and present a rather complete discussion.

%-----------------------------------------------------------------------------------
\paragraph{Mesh refinement.}
%-----------------------------------------------------------------------------------
%
A key component of any adaptive algorithm, such as the three $\AFEM$s already described,
is the routine $\REFINE$ which refines a current mesh $\grid$ into $\grid_*$ to improve
resolution. In Section \ref{S:mesh-refinement} we study the {\it bisection method},
which is the most popular method to refine simplicial meshes in $\R^d$ for $d\ge1$. For
simplicity we focus our attention on this method, but
most results apply to other refinement strategies such as quad-trees (for quadrilaterals),
octrees (for hexagons) and red-green (for simplicial meshes). We do not insist on these
extensions but refer to \cite{bonito2010quasi} for details.

Given an initial grid $\grid_0$ with a suitable labeling, the bisection method splits a
given simplex into two children. The rules for successive cutting of simplices, for instance newest vertex
bisection for $d=2$, are such that the ensuing meshes are shape regular (with a uniform
constant only depending on $\grid_0$ and $d$). However, bisection may
not be completely local to keep conformity.
The analysis of propagation of refinement is a delicate combinatorial
problem. It is easy to see by example, that bisecting one element of large generation
(i.e.\ number of bisections needed to produce it) may require a chain of elements with length
similar to the generation. Therefore, the number of refined elements in one step cannot be
bounded by the number of marked elements. The following amazing estimate by \cite{BiDaDeV:04}
for $d=2$ and \cite{Stevenson:08} for $d>2$ shows that the cumulative effect of bisection counting
from $\grid_0$ all the marked elements $\marked_j$ is quasi-optimal:
there exists a constant $D>0$, depending on $\grid_0$
and $d$, such that
\begin{equation}\label{E:complexity-bisection}
\#\grid_k - \#\grid_0 \le D \sum_{j=0}^{k-1} \#\marked_j.
\end{equation}
This estimate is crucial for the study of rate-optimality of $\AFEM$ and is proved in
Section \ref{S:mesh-refinement} for $d=2$ and for both conforming refinement and $\Lambda$-admissible
refinement. The latter is a systematic way to handle non-conforming meshes that goes back to
\cite{BeCaNoVaVe23}. It associates a computable global index to hanging nodes and imposes a restriction
to them not to exceed
a preassigned value $\Lambda \ge 0$; if $\Lambda=0$ then the mesh is conforming.
Section \ref{S:mesh-refinement} also discusses several interesting geometric properties
of $\Lambda$-admissible meshes which turn out to be crucial for discontinuous Galerkin methods.
Since Section  \ref{S:mesh-refinement} is quite technical, it can be skipped in a first reading.

%-----------------------------------------------------------------------------------
\paragraph{Discontinuous Galerkin methods.}
%-----------------------------------------------------------------------------------
%
These methods, so-called dG, are the natural first step to investigate the role of non-conformity
in adaptivity, namely that the discrete space of discontinuous piecewise polynomials $\V_\grid$
is no longer a subspace of $H^1_0(\Omega)$. To this end, we study the symmetric interior penalty dG
method in Section \ref{S:dg} on $\Lambda$-admissible partitions $\grid$ of $\grid_0$.
Such dG methods exhibit some characteristic
and novel features with respect to conforming FEMs: the most notable one is the presence of weighted
jumps that stabilize the method and compensate for the lack of $H^1$-conformity.
We consider the formulation with lifting, which allows for minimal regularity
$u\in H^1_0(\Omega)$, and forcing $f\in H^{-1}(\Omega)$ despite that $\V_\grid$ is not a
subspace of $H^1_0(\Omega)$. The latter is possible because, within the framework of $\AFEMTS$,
$f$ is approximated by a piecewise polynomial $P_\grid f$ for which the pairing with functions in
$\V_\grid$ is meaningful.

The fact that jumps are not monotone upon refinement constitutes one of the main obstructions for studying
adaptivity for dG. To circumvent this issue we follow \cite{bonito2010quasi}, who in turn modified the
original approach of \cite{karakashian2007convergence}, and introduce the largest conforming subspace
$\V_\grid^0$ of $\V_\grid$. It turns out that, despite being coarser, $\V_\grid^0$ exhibits a local
resolution comparable with $\V_\grid$ because of key geometric properties of $\Lambda$-admissible meshes
that control the degree of non-conformity of $\grid$.
In addition, $\V_\grid^0$ is responsible for the scaled jumps to be bounded by the PDE estimator
$\eta_\grid(u_\grid)$. Exploiting properties of $\eta_\grid(u_\grid)$, similar to the conforming case,
leads to a
quasi-orthogonality estimate for the dG norm, a dG variant of the Pythagoras equality. This is
instrumental to prove a contraction property for the error plus scaled estimator and deduce
convergence for both $\GALERKIN$ and $\AFEMTS$. Moreover, we derive rate-optimality for
both algorithms provided $u$ and $\data$ belong to suitable approximation classes. Such
classes are the same as for conforming AFEMs: in fact we prove that the approximation classes for
$u$ using continuous and discontinuous piecewise polynomials on $\Lambda$-admissible meshes
coincide.

%-----------------------------------------------------------------------------------
\paragraph{Inf-sup stable AFEMs.}
%-----------------------------------------------------------------------------------
%
The convergence and optimality theories developed in Sections \ref{S:convergence-coercive} and
\ref{S:conv-rates-coercive}
rely on the bilinear form in \eqref{E:weak-formulation} being {\it coercive}.
We remove this strong restriction in Section \ref{S:conv-rates-infsup} and consider
uniformly {\it inf-sup stable} FEMs on conforming refinements $\grid_j$ of $\grid_0$.
The lack of an energy norm and its monotone behavior upon refinement has been an
obstacle for the study of this class of problems. We follow the recent work
by M. Feischl \cite{Feischl:2022}, who introduced the following form of {\it quasi-orthogonality}
between consecutive Galerkin solutions $u_j\in\V_j$, originally proposed in
\cite{Axioms:2014} as part of an abstract set of axioms of adaptivity:
\begin{equation}\label{E:quasi-optimality-feischl}
        \sum^{j+N}_{k=j} \|u_{k+1}- u_k\|^2_\V\le C(N)\|u - u_j\|^2_\V \qquad j \ge 0,
\end{equation}
where $C(N)/N\to0$ as $N\rightarrow \infty$. This is our departing point to develop a
{\it variational approach} to prove linear convergence of $u_j$ provided data $\data$ is discrete;
the latter is reflected in an equivalence property between error and estimator (without
oscillation). This is the context of a $\GALERKIN$ routine, which is next used as a
building block together with $\DATA$ for a $\AFEMTS$ that handles general data $\data$.
Moreover, we prove rate-optimality for both algorithms, thereby extending 
Sections \ref{S:convergence-coercive} and \ref{S:conv-rates-coercive}. \looseness=-1

This discussion is rather abstract. We specialize it to the Stokes equations for viscous
incompressible fluids and mixed formulations of \eqref{E:model-PDE} using Raviart-Thomas-Nedelec
and Brezzi-Douglas-Marini elements. We thereby obtain convergence and rate-optimality
for $\AFEMSW$ for the Stokes equations and $\AFEMTS$ for mixed methods with variable and
possibly discontinuous coefficients $(\bA,c)$.

We conclude with a complete proof of \eqref{E:quasi-optimality-feischl} following
\cite{Feischl:2022}. This is a tour de force in applied linear algebra and is rather
technical. It can be omitted in a first reading.

%-----------------------------------------------------------------------------------
\paragraph{Acknowledgments.}
%-----------------------------------------------------------------------------------
%
We are grateful to R. DeVore, P. Morin, and A. Salgado for discussions about regularity
classes and to Ch. Kreuzer for pointing out a gap in an earlier version of the proof of Theorem~\ref{E:zeta-contracts}. We are also thankful to G. Vacca for providing simulations and D. Fassino for
assistance with the manuscript.

This research was supported by the NSF grants DMS-2110811 (AB) and DMS-1908267 (RHN), the MUR-PRIN grants 20227K44ME (CC) and 201752HKH (AV), and  the INdAM research group GNCS (CC and AV). RHN was also partially supported by NSF grant DMS-1929284 while he was in residence at the Institute for Computational and Experimental Research in Mathematics in Providence, RI, during 
the Numerical PDEs: Analysis, Algorithms, and Data Challenges semester program.

%--------------------------------------------------------------------------------

\section{Linear Elliptic Boundary Value Problems: Basic Theory}\label{S:bvp}
%\rhn{(RHN $\longrightarrow$ CC)} %\cite{NochettoVeeser:2012}

In this section we examine the variational formulation of elliptic
partial differential equations (PDE). We start with a brief review of
Sobolev spaces and their properties and continue with two model boundary
value problems that are instrumental in our subsequent analysis. We next
present the so-called inf-sup theory that characterizes existence, uniqueness
and stability of variational problems, and apply it to coercive and saddle
point problems. These two classes will play essential roles later.

%-------------------------------------------------------------------------------------
\subsection{Sobolev spaces: scaling and embedding}\label{nv-S:sobolev-number}
%-------------------------------------------------------------------------------------
Let $\Omega\subset\R^d$ with $d>1$ be a Lipschitz and bounded domain,
and let $k\in\N, 1\le p \le \infty$. The Sobolev space $W^k_p(\Omega)$
is defined by
\[
\index{Functional Spaces!$W^k_p(\Omega)$: Sobolev spaces}
W^k_p(\Omega) := \{ v:\Omega\to\R | ~D^\alpha v\in L^p(\Omega)
~~\forall |\alpha|\le k \},
\]
and is a Banach space with the norm
\[
\|v\|_{W^k_p(\Omega)} = \left(\sum_{j=1}^k |v|_{W^j_p(\Omega)}^p \right)^{\frac{1}{p}},
\quad
|v|_{W^j_p(\Omega)} = \left( \int_\Omega |D^j v|^p \right)^{\frac{1}{p}}.
\]
The space $W^k_p(\Omega;\mathbb R^m)$ is the space $W^k_p(\Omega)$ of vector or matrix-valued functions.
If $p=2$ we write $H^k(\Omega) = W^k_2(\Omega)$ and note that this is a
Hilbert space. We denote by $H^1_0(\Omega)\subset H^1(\Omega)$ the completion of $C_0^\infty(\Omega)$
within $H^1(\Omega)$. 

Sobolev spaces $W^k_p(\Omega)$ of fractional order $k>0$ can be defined as well, by applying the real interpolation method between $W^{[k]}_p(\Omega)$ and $W^{[k]+1}_p(\Omega)$ (see \cite{BerghLofstrom:76} and \cite{AdamsFournier:03} for the details). The subsequent definitions and properties hold for Sobolev spaces of integer or fractional order. 

The {\em Sobolev number} of $W^k_p(\Omega)$ is given by
\index{Definitions!Sobolev number $\sob (W^k_p)$}
\begin{equation}\label{nv-sobolev-number}
\sob (W^k_p) := k - \frac{d}{p}.
\end{equation}
This number governs the scaling properties of the semi-norm $|v|_{W^k_p(\Omega)}$,
because rescaling variables $\wh x = \frac{1}{h} x$ for all
$x\in\Omega$, transforms $\Omega$ into $\wh\Omega$ and $v$ into 
$\wh v$, while the corresponding norms scale as
\[
|\wh v|_{W^k_p(\wh\Omega)} = h^{\sob (W^k_p)} |v|_{W^k_p(\Omega)}.
\]
%

%-------------------------------------------------------------------------------------
\subsection{Properties of Sobolev spaces}\label{S:poincare-trace}
%-------------------------------------------------------------------------------------
%
We summarize now, but not prove, several important properties of
Sobolev spaces which play a key role later. We refer to
\cite{Evans:10,GilbargTrudinger:01,Grisvard:11} for details.
We use the notation $\hookrightarrow$ to denote a continuous embedding. 
%\todo[inline]{RHN: should we introduce fractional Sobolev spaces here and the embedding theorem along with the nonlinear Sobolev line? This is a concept that repeats many times later and we need a central and early place to introduce it.}
  
%-------------------------------------------------------------------------------------
\begin{lemma}[Sobolev embedding]\label{L:embedding}
Let $m>k\ge 0$ and assume
  $\sob(W^m_q)>\sob(W^k_p)$. Then the embedding
$
    W^m_q(\Omega)\hookrightarrow W^k_p(\Omega)
$
is compact and
\[
\|v\|_{W^k_p(\Omega)} \le C \|v\|_{W^m_q(\Omega)}
\quad\forall \, v\in W^m_q(\Omega),
\]
where $C=C(m,k,q,p,\Omega,d)$.
\end{lemma}

\begin{figure}[h!]
\begin{center}
\begin{tikzpicture}[scale =0.9]
\fill[gray,opacity=0.25] plot coordinates{ (0,1.5)
(1.5+2,1.5)(3.5+2,5.5)(0,5.5)};

\draw [thick,->]plot coordinates{(0,0) (8,0)};
\draw [thick,->]plot coordinates{(0,0) (0,6.3)};
\draw [thick]plot coordinates{(3.5+2,5.5)(1.5+2,1.5)};
\draw [thick]plot[mark=o, only marks, mark size=2pt]
coordinates{(3+2,4.5)(1.5+2,1.5)};

\draw [thick, dashed]plot coordinates{(0,4.5) (3+2,4.5)};
\draw [thick, dashed]plot coordinates{(0,1.5) (1.5+2,1.5)};
\draw [thick, dashed]plot coordinates{(3+2,0) (3+2,4.5)};
\draw [thick, dashed]plot coordinates{(1.5+2,0) (1.5+2,1.5)};
%labels
\coordinate [label= {differentiability}] (z1) at (-1.5,5.6);
\coordinate [label= $1/{\text{summability}}$] (z1) at (8,-0.8);
\coordinate [label= $m$] (z1) at (-0.5,4.25);
\coordinate [label= $k$] (z1) at (-0.5,1.25);
\coordinate [label= $1/q$] (z1) at (3+2,-0.8);
\coordinate [label= $1/p$] (z1) at (1.5+2,-0.8);
\coordinate [label= $W^{k}_p(\Omega)$] (z1) at (1.3+1.5,1.5);
\coordinate [label= $W^{m}_q(\Omega)$] (z1) at (2.8+1.5,4.5);
\coordinate [label= $0$] (z1) at (-0.25,-0.5);
\coordinate [label= {$\sob(W^{m}_q(\Omega))= \text{const}$}] (z1) at
(5.5+2,5.5);
\coordinate [label= {($\text{slope} = d$) }] (z1) at (5.5+2,5);
\end{tikzpicture}
\caption{DeVore diagram \cite{DeVore:98}. The space $W^{k}_p(\Omega)$ is represented by the point $(1/p,k)$ in the first quadrant. The line $\sob(W^{m}_q)=\text{const}= \sob(W^{k}_p)$, with slope $d$, may be called the (critical) Sobolev embedding line fir $W^k_p(\Omega)$. It represents all Sobolev spaces having the same Sobolev number as $W^{k}_p(\Omega)$. Sobolev spaces corresponding to points inside the gray region and on its boundary on the vertical axis are compactly embedded in $W^{k}_p(\Omega)$. Spaces on the oblique and horizontal lines emanating from 
$W^{k}_p(\Omega)$ are generally continuously,
but not compactly embedded, in $W^{k}_p(\Omega)$ with exceptions such as $p=\infty$.
Note that indices $k$ and $m$ may take non-integer values, corresponding to fractional Sobolev spaces.}\label{F:devore-diagram}
\end{center}
\end{figure}
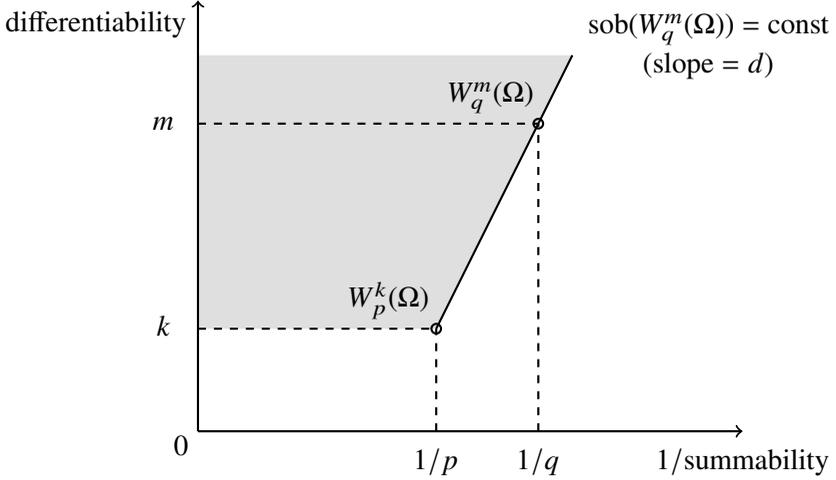

We say that two Sobolev spaces are in the same nonlinear Sobolev scale if they
have the same Sobolev number; see Figure~\ref{F:devore-diagram}. We thus note that for compactness the space
$W^m_q(\Omega)$ must be above the Sobolev scale of $W^k_p(\Omega)$, i.e. $\sob(W^m_q)>\sob(W^k_p)$.

  The assumption on the Sobolev number cannot be relaxed. To see this,
  consider $\Omega$ to be the unit ball of $\R^d$ for $d\geq 2$ and
  set $v(x)=\log\log\frac{|x|}{2}$ for
  $x\in\Omega\setminus\{0\}$. Then there holds $v\in W^1_d(\Omega) \setminus L^\infty(\Omega)$, but
  \begin{displaymath}
    \sob(W^1_d)=1-d/d=0=0-d/\infty=\sob(L^\infty).
  \end{displaymath}
  Therefore, equality cannot be expected in the embedding theorem.
  On the other hand, consider $d=1$ and the spaces $W^1_1(\Omega)$ and $L^\infty(\Omega)$.
  We see that $\sob(W^1_1) = \sob(L^\infty)=0$ but
$W^1_1(\Omega)$ is compactly embedded in $L^\infty(\Omega)$ in this
case. This shows that these two spaces are in the same nonlinear Sobolev scale
and that the above inequality between Sobolev numbers for a compact embedding
is only sufficient. 

Moreover, if $0< \alpha = \sob(W^k_p) \le 1$, then functions of $W^k_p(\Omega)$ are H\"older-$\alpha$ and
\[
|v|_{C^{0,\alpha}(\overline{\Omega})} \le C |v|_{W^k_p(\Omega)}
\quad\forall \, v \in W^k_p(\Omega).
\]
This allows for the use of the standard Lagrange
interpolation operator. We will exploit this fact later in Section \ref{S:tools}.

%-------------------------------------------------------------------------------------  
\begin{lemma}[first Poincar\'e inequality]\label{L:Poincare}
  % \todo[inline]{RHN: Is this really due to Poincar\'e? CC: Seems yes.
  % RHN (01/14/24): These are the correct quotes according to Brenner-Scott. But I would write them in 
  % $W^1_p(\Omega)$ rather than $H^1(\Omega)$. AV (14/1/24): Remark 3.32 in Ern/Guermond I gives some historical information. If I remember correctly, we would be safe if we collect them in one lemma called Poincar\'e inequalities. Doing so, I would however point out the difference in the constants.
  % RHN (01/14/24): We use these inequalities with different Lebesgue exponents on left and right-hand sides
  % (see Section 8.3.1). Should we state these lemmas more generally? What are suitable references?}
  Let $\Omega\subset\R^d$ be a bounded and Lipschitz domain. Then there is a constant
  $C_P=C_d |\Omega|^{1/d}$ such that
  \begin{equation}\label{E:Poincare}
  \index{Constants!$C_P$: Poincar\'e constant}
  \|v\|_{L^2(\Omega)} \le C_P \|\nabla v\|_{L^2(\Omega)} \quad\forall \, v\in H^1_0(\Omega).
  \end{equation}
  The same inequality is valid in $W^1_p(\Omega)$ for any $1\le p \le \infty$ provided $v$ has vanishing trace
  \cite[p.158]{GilbargTrudinger:01}.
\end{lemma}

%-------------------------------------------------------------------------------------
\begin{lemma}[second Poincar\'e inequality]\label{L:poincare-friedrichs}
  % \todo[inline]{RHN: is this the correct quote? CC: according to Wikipedia, this is the Poincare'-Wirtinger inequality - the Friedrichs inequality is like the Poincare', but with higher-order norm on the rhs AV: I would call this the Poincar\'e inequality}
  \index{Constants!$C_P$: Poincar\'e constant}
  There exists $C_P$ depending on $\Omega$ such that
  \begin{equation}\label{E:Friedrichs}
    \|v- \overline{v} \|_{L^2(\Omega)}
    \le C_P \|\nabla v\|_{L^2(\Omega)} \quad\forall \, v\in H^1(\Omega),
  \end{equation}
  where $\overline{v} := |\Omega|^{-1}\int_\Omega v$. The best constant within the class of convex domains is
$
 C_P = \frac1\pi \diam(\Omega)
$
 \cite{PayneWeinberger:60}. The same inequality is valid in $W^1_p(\Omega)$ for $1\le p \le \infty$.
 \cite{GilbargTrudinger:01}.
\end{lemma}

%-------------------------------------------------------------------------------------
\begin{lemma}[traces]\label{L:trace-id}
Let $\Omega$ be Lispchitz.
There exists a unique linear operator
 $T\colon H^1(\Omega)\to L^2(\partial\Omega)$
  such that
  \begin{equation*}
    \begin{aligned}
      \|Tv\|_{L^2(\partial\Omega)} &\le c(\Omega) \|v\|_{H^1(\Omega)}
      &\quad& \forall \, v\in H^1(\Omega),\\
      Tv &= v
      &&\forall \, v|_{\partial \Omega}\in C^0(\overline\Omega)\cap H^1(\Omega).
    \end{aligned}
  \end{equation*}
  The operator $T$ is also well defined on $W^1_p(\Omega)$ for $1\le p \le \infty$. \cite{Evans:10, Grisvard:85}
  \end{lemma}
  %
  % \todo[inline]{RHN (01/14/24): We need a reference here. Evan assumes that $\partial\Omega$ is $C^1$.
  % We also need to introduce $H^{1/2}(\partial\Omega)$. We use it in Section 8.3.1. AB (1/14/23): it is definitely in Girault-Raviart and Grisvard but I am not sure they are the appropriate refs. Maybe Adams?}
  Since $Tv=v_{|\partial\Omega}$ for continuous functions, we write $v$ for $Tv$.
  The image of $T$ is a strict
  subspace of $L^2(\partial\Omega)$, the so-called
  $H^{1/2}(\partial\Omega)$. This is a Hilbert space if equipped with the norm $\| g \|_{H^{1/2}(\partial\Omega)}:= \inf\{ \| v \|_{H^1(\Omega)} \ | \ Tv=g \}$, and $T$ is continuous with respect to this norm, 
  i.e.,
  \begin{equation}\label{E:trace-1/2}
 \|Tv\|_{H^{1/2}(\partial\Omega)} \le \|v\|_{H^1(\Omega)}
      \quad \forall \, v\in H^1(\Omega).
  \end{equation}
  The definition of $H^1_0(\Omega)$
  can be reconciled with that of traces because
\begin{equation*}
  H^1_0(\Omega) = \
  \big\{v\in H^1(\Omega) \mid v = 0 \text{ on } \partial\Omega\}.
\end{equation*}
We point out that, in view of Lemma \ref{L:Poincare}, the semi-norm $|v|_{H^1_0(\Omega)} := \|\nabla v\|_{L^2(\Omega)}$ is a norm in $H^1_0(\Omega)$. We let
$H^{-1}(\Omega)$ be the dual of $H^1_0(\Omega)$, with corresponding norm
$\|f\|_{H^{-1}(\Omega)}:=\sup_{v\in H^1_0(\Omega)} \frac{\langle f,v \rangle}{|v|_{H^1_0(\Omega)}}$.
These definitions extend to any $p\in(1,\infty)$.

\begin{lemma}[Gauss divergence theorem]\label{L:divergence-thm}
If $\bn$ is the outer unit normal to $\Omega$, then \looseness=-1
\[
    \int_\Omega \divo\vec{w} =
    \int_{\partial\Omega} \vec{w}\cdot\bn
    \qquad\forall\, \vec{w}\in W^1_1(\Omega;\R^d).
\]
\end{lemma}

\begin{lemma}[Green's formula]\label{L:Green}
The integration by parts formula holds
\begin{equation*}\label{nsv-Green2}
    \int_\Omega \divo\vec{w}\,v = -\int_\Omega \vec{w}\cdot \nabla v + 
    \int_{\partial\Omega} v\,\vec{w}\cdot\bn
    \qquad\forall \, v\in H^1(\Omega),\vec{w}\in H^1(\Omega;\R^d).
\end{equation*}
\end{lemma}

%-------------------------------------------------------------------------------------
\subsection{Examples of boundary value problems}\label{S:ex-bvp}
%-------------------------------------------------------------------------------------
%
We consider two model elliptic problems in this paper.
We start with the {\it scalar diffusion-reaction equation} with variable coefficients
\begin{equation}
\label{strong-form}
\begin{aligned}
L[u] := -\div( \bA\nabla u) + c u &= f &&\text{ in }\Omega,
\\
 u &= 0 &&\text{on }\partial\Omega,
\end{aligned}
\end{equation}
where $\Omega\subset\Rd$ is a bounded domain with Lipschitz boundary,
 $\bA\in L^\infty(\Omega;\R^{d\times d})$ is a diffusion tensor uniformly
symmetric positive definite (SPD) over $\Omega$, \ie there exist
constants $0<\alpha_1\le\alpha_2$ such that
\index{Constants!$(\alpha_1,\alpha_2)$: lower and upper bounds of the diffusion coefficient spectrum}
\begin{equation}\label{uniform_spd}
  \alpha_1|\vec{\xi}|^2\le\vec{\xi}^T \!\bA(x)\,\vec{\xi}\le\alpha_2|\vec{\xi}|^2
  \quad\fa x\in\Omega,\ \vec{\xi}\in\R^d,
\end{equation}
$c\in L^\infty(\Omega), c\ge0$ is a reaction term, and $f\in L^2(\Omega)$ is a scalar load term.

To derive the variational formulation of \eqref{strong-form} we let
$\V=H^1_0(\Omega)$ and $\V^*=H^{-1}(\Omega)$.
Since $H^1_0(\Omega)$ is the subspace of $H^1(\Omega)$ of functions
with vanishing trace,  asking for $u\in\V$ accounts for the 
homogeneous Dirichlet boundary values in \eqref{strong-form}.
We next multiply \eqref{strong-form} with a test function $v\in
H^1_0(\Omega)$, integrate over $\Omega$ and use Lemma \ref{L:Green} (Green's formula),
provided $\vec{w}=-\bA\nabla u\in H^1(\Omega;\R^d)$, to obtain
\begin{equation}\label{weak-form}
 u\in\V:
\quad
 \bilin{u}{v}
 =
 \dual{f}{v}
\quad
 \fa v\in\V.
\end{equation}
Here, $\mathcal{B}:\V\times\V \rightarrow \mathbb R$ stands for the bilinear form
\begin{equation}\label{bilinear-coercive}
    \bilin{w}{v} \definedas \int_\Omega\nabla v\cdot \bA \nabla w + cvw
    \qquad \fa v,w \in \mathbb V,
\end{equation}
and $\dual{\cdot}{\cdot}$ stands for the $L^2(\Omega)$-scalar
product and also for a duality paring $H^{-1}(\Omega)-H^1_0(\Omega)$.
Since $\bA$ is assumed to be symmetric, the bilinear form $\B$ is also symmetric; however
$\bA$ does not have to be symmetric in general.
Note that the weak form
\eqref{weak-form} allows for fluxes $\vec{w}\in L^2(\Omega;\R^d)$ and forcing $f\in H^{-1}(\Omega)$.
We examine existence, uniqueness and stability of \eqref{weak-form} in Section \ref{S:inf-sup} below.

We assume {\it homogeneous Dirichlet} boundary condition in \eqref{strong-form} for simplicity and because
this will be our basic setting later. However, we could allow a {\it non-homogeneous Dirichlet} condition
in the sense of traces
\begin{equation}\label{E:non-hom-Diirichlet}
Tu = g\quad\textrm{on }\partial\Omega,
\end{equation}
for any given function $g\in H^{\frac12}(\partial\Omega)$. To write the companion variational formulation
to \eqref{weak-form}, we first introduce the subspace $\V(g) \subset H^1(\Omega)$ of functions $v$ with trace
$Tv=g$ on $\partial\Omega$, and then rewrite \eqref{weak-form} as follows:
\begin{equation}\label{E:weak-nonhom-Dirichlet}
u \in \V(g): \quad \B[u,v] = \langle f, v \rangle \quad\forall v\in \V(0).
\end{equation}
Moreover, we could consider a {\it Robin} boundary condition for given functions $g$ and $p$ on $\partial\Omega$
\begin{equation}\label{E:Robin}
\bA \nabla u \cdot \bn + p u = g \quad\textrm{on }\modif{\partial\Omega,}
\end{equation}
\modif{where $\bn$ is the outer unit normal to $\Omega$.}
To figure out the variational formulation, we now multiply the PDE in \eqref{strong-form}
by a test function $v\in H^1(\Omega)$ and integrate by parts to find the following variant of
\eqref{weak-form} \looseness=-1
\begin{equation}\label{E:Robin-weak}
\index{Functions!$u$: solution to weak formulation \eqref{weak-form}}
u\in H^1(\Omega): \quad \B[u,v] = \ell(v) \quad\forall v\in H^1(\Omega),
\end{equation}
where for all $v,w \in H^1(\Omega)$
\begin{equation}\label{E:Robin-bilinear}
  \B[w,v] := \int_\Omega \nabla v \cdot\bA \nabla w + c v w + \int_{\partial\Omega} p vw,
  \quad \ell(v) := \langle f, v \rangle + \int_{\partial\Omega} g v.
\end{equation}
We realize that \eqref{E:Robin-bilinear} makes sense for $p\in L^\infty(\partial\Omega),  p\ge 0$ and
$g \in H^{-\frac12}(\partial\Omega)$, the dual space of $H^{\frac12}(\partial\Omega)$, whence the last integral in
\eqref{E:Robin-bilinear} means a duality pairing.
If $p=0$, then \eqref{E:Robin-bilinear} reduces to the weak form of the {\it Neumann}
boundary condition.

The second model problem is the {\it Stokes system}, which is the simplest model of a stationary viscous
incompressible fluid.
Given an external force $\bF\in L^2(\Omega;\R^d)$, let 
the velocity-pressure pair $(\bu,p)$ satisfy the momemtum and
incompressibility equations with no-slip boundary condition:
\begin{equation}\label{strong-stokes}
\begin{alignedat}{2}
-\Delta \bu + \nabla p &= \bF &\qquad&\text{in }\Omega,
\\
\div\bu & = 0 &&\text{in }\Omega,
\\
\bu &= 0 &&\textrm{on }\partial\Omega.
\end{alignedat}
\end{equation}
For the variational formulation we consider two Hilbert spaces
$\V=H^1_0(\Omega;\R^d)$ and $\Q=L^2_0(\Omega)$, where $L^2_0(\Omega)$
is the space of $L^2$ functions with zero mean value.  The space
$H^1_0(\Omega;\R^d)$ takes care of the no-slip boundary values of the
velocity. We first multiply the momentum equation in \eqref{strong-stokes} by $\vec{v}\in\V$,
assume $\vec{u}\in H^2(\Omega;\R^d)$ and use Lemma \ref{L:Green} (Green's formula)
component-wise. We next multiply the incompressibility equation in \eqref{strong-stokes} by
$q\in\Q$ and integrate over $\Omega$. We end up with the following variational formulation:
find $(\bu,p)\in\V\times\Q$ such that 
\begin{equation}\label{weak-stokes}
\begin{alignedat}{3}
  &a[\vec{u},\vec{v}] + b[p,\vec{v}] &&= \scp{\vec{f}}{\vec{v}} 
&\qquad &\forall \, \vec{v}\in\V,
  \\
  &b[q,\vec{u}] &&= 0 &&\forall \, q\in\Q.
\end{alignedat}
\end{equation}
Here the bilinear forms $a\colon\V\times\V\to\R$ and $b\colon\Q\times\V\to\R$ read
\begin{equation*}
  a[\vec{w},\vec{v}]\definedas \int_\Omega \nabla\bv : \nabla \vec{w}
  = \sum_{i=1}^d \int_\Omega \nabla v_i \cdot \nabla w_i
  \qquad\forall \, \bv,\vec{w}\in\V.
\end{equation*}
and
\begin{equation*}
  b[q,\vec{v}] \definedas -\int_\Omega q\,\divo\vec{v}
  \qquad\forall \, q\in\Q,\,\vec{v}\in\V.
\end{equation*}
We observe that $ a[\vec{w},\vec{v}]$ does not require $\vec{w} \in H^2(\Omega;\R^d)$ and that
\eqref{weak-stokes} makes sense for $\vec{f} \in H^{-1}(\Omega;\R^d)$; note that the second equation in \eqref{weak-stokes} is always satisfied for constant $q$ due to Gauss' theorem, which explains the choice $q\in \Q$.
Furthermore, \eqref{weak-stokes} can be reformulated as \eqref{weak-form}, namely
\begin{equation}\label{weak-stokes-B}
  (\vec{u},p)\in \V\times\Q: \quad
  \bilin{(\vec{u},p)}{(\vec{v},q)} =  \scp{\vec{f}}{\vec{v}} \quad\forall \,
  (\vec{v},q)\in \V\times\Q,
\end{equation}
with
\[
\bilin{(\vec{u},p)}{(\vec{v},q)} \definedas
a[\vec{u},\vec{v}] + b[p,\vec{v}]
+ b[q,\vec{u}].
\]
We discuss existence, uniqueness and stability of \eqref{weak-stokes-B}  in Section~\ref{S:inf-sup}.

We could formulate the Stokes system with other boundary conditions. First, we may allow a non-homogeneous
Dirichlet condition $\vec{u} = \vec{g}$ for a given function $\vec{g} \in H^{\frac12}(\partial\Omega;\R^d)$ satisfying the compatibility condition $\int_\Omega \vec{g}\cdot \vec{n} = 0$ imposed by Gauss' theorem,
and proceed as in the scalar case \eqref{E:weak-nonhom-Dirichlet}. Second, to deal with a Neumann boundary
condition, we introduce the stress tensor
$\vec{\sigma} (\vec{u},p) := \frac12 \big(\nabla \vec{u} + \nabla\vec{u}^T \big) - p \vec{I}$
and the symmetric part of the velocity gradient
$\vec{\eps}(\vec{u}):=\frac12 \big(\nabla \vec{u} + \nabla\vec{u}^T \big)$. Then instead of
\eqref{strong-stokes} we could write the strong form of the Neumann problem as
\begin{equation}\label{E:strong-Neumann-Stokes}
\begin{gathered}
  -\div \vec{\sigma} ( \vec{u},p) = \vec{f},
  \quad
  \div \vec{u} = 0,
  \quad \textrm{in } \Omega
  \\
  \vec{\sigma} ( \vec{u},p) \vec{n} = \vec{g} \quad \textrm{on } \partial\Omega,
\end{gathered}
\end{equation}
and its weak form as \eqref{weak-stokes} but with 
\begin{equation}\label{E: spaces_Neumann_stokes}
\V := \{ \vec{v} \in H^1(\Omega;\R^d) : \int_\Omega \vec{v} = \vec{0}\}, \qquad \Q:=L^2(\Omega),
\end{equation}
as well as bilinear form $a$ and right-hand side
\begin{equation}\label{E:momentum-Stokes}
a[\vec{u},\vec{v}] := \int_\Omega \vec{\eps} (\vec{u}) : \vec{\eps} (\vec{v}),
  \quad
 \ell(\vec{v}) := \scp{\vec{f}}{\vec{v}} + \int_{\partial\Omega} \vec{g}\cdot\vec{v}
\end{equation}
for all $\vec{v}\in \V$. Again, the last integral in \eqref{E:momentum-Stokes} is to be
interpreted as a duality pairing for $\vec{g} \in H^{-\frac12}(\partial\Omega;\R^d)$.

%-------------------------------------------------------------------------------------
\subsection{Inf-sup theory}\label{S:inf-sup}
%-------------------------------------------------------------------------------------
%
We present a functional framework for existence, uniqueness,
and stability of variational problems of the form \eqref{weak-form} or \eqref{weak-stokes-B}.
Throughout this section we let  
\begin{math}
  (\V, \scp[\V]{\cdot}{\!\!\cdot})
\end{math}
and 
\begin{math}
  (\W, \scp[\W]{\cdot}{\!\!\cdot})
\end{math}
be a pair of Hilbert spaces with induced norms
$\|\cdot\|_\V$ and $\|\cdot\|_\W$. We denote by
$\V^*$ and $\W^*$ their respective dual spaces equipped with norms
\begin{displaymath}
  \|f\|_{\V^*} = \sup_{v\in\V} \frac{\dual{f}{v}}{\|v\|_\V}
  \qquad\text{and}\qquad
  \|g\|_{\W^*} = \sup_{v\in\W} \frac{\dual{g}{v}}{\|v\|_\W}.
\end{displaymath}
We write $\Lin(\V;\W)$ for the space of all linear and continuous operators
from $\V$ into $\W$ with operator norm
\begin{displaymath}
  \|B\|_{L(\V;\W)} = \sup_{v\in\V} \frac{\|Bv\|_\W}{\|v\|_\V}.
\end{displaymath}
The following result relates a continuous bilinear form
$\B\colon\V\times\W\to\R$ with an operator $B\in\Lin(\V;\W)$
\cite{Necas:62,babuvska1971error}.

\begin{theorem}[Banach-Ne\v cas]\label{nsv-T:BanachNecasBabuska}
{\it  Let $\B:\V\times\W\to\R$ be a continuous bilinear form with norm
  \begin{equation}\label{nsv-bilinear-cont}
    \|\B\|\definedas \sup_{v\in\V}\sup_{w\in\W} 
    \frac{\bilin{v}{w}}{\|v\|_\V \|w\|_\W}.
  \end{equation}
  Then there exists a unique linear operator $B\in L(\V,\W)$ such that 
  \begin{displaymath}
    {\dual{Bv}{w}}_\W=\bilin{v}{w} \qquad\forall \, v\in\V,\,w\in\W
  \end{displaymath}	
  with operator norm 
  \begin{displaymath}
    \|B\|_{L(\V;\W)}=\|\B\|.
  \end{displaymath}
  Moreover, the bilinear form $\B$ satisfies the two conditions
  \begin{subequations}
    \begin{align}
    \index{Constants!$\alpha$: $\inf\sup$ constant}
      & i) \ \text{there exists~} \alpha>0 \text{~ such that~~}
      \alpha\|v\|_\V\le\sup_{w\in\W}\frac{\bilin{v}{w}}{\|w\|_\W} \text{ for all }
      v\in\V;
      \label{nsv-I2} \\
      & ii) \ \text{for every~ }0\neq w\in\W\text{ ~there exists~ }v\in\V
      \text{ ~such that~ }\bilin{v}{w}\neq 0; \label{nsv-I3}
    \end{align}
  \end{subequations}
  if and only if $B:\V\to\W$ is an isomorphism with 
  \begin{equation}\label{nsv-I3a}
    \|B^{-1}\|_{L(\W,\V)}\le\alpha^{-1}.
  \end{equation}
}  
\end{theorem}

We consider now the abstract variational problem
 \begin{equation}\label{var-prob-general}
    u\in\V:\quad \bilin{u}{v}=\dual{f}{v}
    \quad\forall \, v\in\W,
 \end{equation}
The following result, due to Ne\v cas \cite[Theorem~3.3]{Necas:62},
characterizes properties of the bilinear
form $\B$ that imply that \eqref{var-prob-general} is well-posed.
\begin{theorem}[Ne\v cas]\label{nsv-T:ex1-var-prob}
{\it Let $\B\colon\V\times\W\to\R$ be a continuous bilinear form.
Then the variational problem \eqref{var-prob-general}
  admits a unique solution $u\in\V$ for all $f\in\W^*$, 
  which depends continuously on $f$, if and only
  if the bilinear form $\B$ satisfies one of the equivalent
  \emph{inf-sup conditions}:
 \begin{enumerate}%[\upshape(1)]
 \setlength{\itemindent}{-0.6em}   
  \item[(1)] There exists $\alpha>0$ such that
    \begin{subequations}\label{nsv-inf-sup-1}
      \begin{align}
        & i) \qquad \sup_{w\in\W}\frac{\bilin{v}{w}}{\|w\|_\W} \ge \alpha \|v\|_\V
        \quad\text{\it for some }~\alpha>0; \qquad
        \label{nsv-inf-sup-1.a}
        % \intertext{and}
        \\
        & ii) \ \text{For every $0\ne w\in\W$ there exists $v\in\V$ such
          that $\bilin{v}{w} \ne 0$}.
        \label{nsv-inf-sup-1.b}
      \end{align}
    \end{subequations}
    \item[(2)] There holds
      \begin{equation}\label{nsv-inf-sup-2}
        \infimum_{v\in\V}\sup_{w\in\W}\frac{\bilin{v}{w}}{\|v\|_\V \|w\|_\W}>0,
        \qquad
        \infimum_{w\in\W}\sup_{v\in\V}\frac{\bilin{v}{w}}{\|v\|_\V \|w\|_\W}>0.
      \end{equation}
    \item[(3)] There exists $\alpha>0$ such that
      \begin{equation}\label{nsv-inf-sup}
        \infimum_{v\in\V}\sup_{w\in\W}\frac{\bilin{v}{w}}{\|v\|_\V \|w\|_\W}
        = \infimum_{w\in\W}\sup_{v\in\V}\frac{\bilin{v}{w}}{\|v\|_\V \|w\|_\W}
        = \alpha.
      \end{equation}
  \end{enumerate}
  In addition, the solution $u$ of \eqref{var-prob-general} satisfies the
  stability estimate
  \begin{equation}\label{cont-stab}
    \|u\|_\V \leq \alpha^{-1} \|f\|_{\W^*}.
  \end{equation}
}  
\end{theorem}

The equality in \eqref{nsv-inf-sup} might seem at first surprising but
is just a consequence of 
\begin{math}
  \|B^{-*}\|_{\Lin(\V;\W)}=\|B^{-1}\|_{\Lin(\W;\V)}.
\end{math}
In general, \eqref{nsv-inf-sup-1} is simpler to verify than
\eqref{nsv-inf-sup} and $\alpha$ of \eqref{nsv-inf-sup} is the largest
possible $\alpha$ in \eqref{nsv-inf-sup-1.a}. Since \eqref{nsv-I3a} shows that
$\|B^{-1}\|_{L(\W,\V)}$ is the best inf-sup constant $\alpha$ in \eqref{nsv-I2},
we readily obtain the following result.

\begin{corollary}[well posedness implies inf-sup]\label{C:well-posedness}
{\it Assume that
the variational problem~\eqref{var-prob-general}
  admits a unique solution $u\in\V$ for all $f\in\W^*$ so that
  \begin{displaymath}
    \|u\|_\V \leq C\|f\|_{\W^*}.
  \end{displaymath}
  Then $\B$ satisfies the inf-sup condition \eqref{nsv-inf-sup} with
  $\alpha\geq C^{-1}$.
}
\end{corollary}

We next apply these abstract results to two special but important cases. The
first class are problems with coercive bilinear form and the second
one comprises problems of saddle point type.

%----------------------------------------------------------------------
\bigskip\noindent
{\it Coercive Problems.}
An existence and uniqueness result for coercive bilinear forms was
established by Lax and Milgram
eight years prior to the result by Ne\v cas \cite{LaxMilgram:54}. Coercivity
of $\B$ is a sufficient condition for existence and uniqueness but it
is not necessary. In this case, we assume $\V=\W$.

\begin{corollary}[Lax-Milgram]\label{nsv-C:lax-milgram}
{\it  Let $\B\colon\V\times\V\to\R$ be a continuous bilinear form that
  is \emph{coercive}, namely there exists $\alpha>0$ such that
  \begin{equation}\label{nsv-coercivity}
    \bilin{v}{v} \ge \alpha \|v\|_\V^2
    \quad\forall \, v\in\V.
  \end{equation}
  Then \eqref{var-prob-general} has a unique solution that satisfies the stability estimate
  \eqref{cont-stab}.
}
\end{corollary}

If the bilinear form $\B$ is also symmetric, \ie
\begin{displaymath}
  \bilin{v}{w} = \bilin{w}{v}
  \quad\forall \, v,w\in\V,
\end{displaymath}
then $\B$ is a scalar product on $\V$. The norm induced by $\B$ is
the so-called \emph{energy norm}\looseness=-1
\begin{displaymath}
  \enorm{v}\definedas\bilin{v}{v}^{1/2}.
\end{displaymath}

For the reaction-diffusion equation \eqref{strong-form} the bilinear form given in
\eqref{bilinear-coercive} satisfies
\begin{equation}
\label{coercivity-and-continuity-constants-for_model-setting}
 0 < \alpha_1
 \leq
 \alpha
 \leq
 \| \B \| \leq \alpha_2 + \| c \|_{L^\infty(\Omega)} C_P^2,
\end{equation}
where $C_P$ is the constant in Lemma \ref{L:Poincare} (first Poincar\'e inequality).
Coercivity and continuity of $\B$, with constants $c_\B=\alpha_1$ and $C_\B=\|\B\|$,
in turn imply that the {\it energy norm} $\enorm{\,\cdot \,}$
is equivalent to the natural norm $\|\cdot\|_\V$ in $\V=H^1_0(\Omega)$:
\begin{equation}\label{nsv-enorm-Vnorm}
\index{Constants!$(c_\B,C_\B)$: norm equivalence constants}
  c_\B \|v\|_\V^2 \leq \enorm{v}^2 \leq C_\B \|v\|_\V^2
  \quad\forall \, v\in\V.
\end{equation}
Moreover, it is rather easy to show that for symmetric and coercive 
$\B$ the solution $u$ of \eqref{var-prob-general} is the unique minimizer
of the quadratic energy
\begin{displaymath}
  J[v] \definedas \frac12 \bilin{v}{v} - \dual{f}{v}
  \quad\forall \, v\in\V,
\end{displaymath}
\ie $u = \argmin_{v\in\V} J[v]$. In particular,
the energy norm and the quadratic energy play a relevant role in 
Sections~\ref{S:convergence-coercive}, \ref{S:conv-rates-coercive} and \ref{S:dg}.

This framework applies to the scalar diffusion-reaction equation \eqref{weak-form} with
{\it homogeneous Dirichlet} condition. Since the full $H^1$-norm $\|\cdot\|_{H^1(\Omega)}$ and the seminorm
$|\cdot|_{H^1_0(\Omega)}$ are equivalent in the space $\V=H^1_0(\Omega)$, according to
Lemma \ref{L:Poincare} (first Poincar\'e inequality), the bilinear form $\mathcal{B}$ in
\eqref{bilinear-coercive} is coercive and continuous in view of \eqref{uniform_spd}
and $c\ge0$. This framework applies as well to the {\it non-homogeneous Dirichlet} problem
\eqref{E:weak-nonhom-Dirichlet}, upon extending $g$ to a function $g\in H^1(\Omega)$ and rewriting
the problem in terms of $w=u-g\in H^1_0(\Omega)$ and forcing $\ell=f-L[g] \in H^{-1}(\Omega)$.

On the other hand, the bilinear form $\B$ in \eqref{E:Robin-bilinear} associated with a
{\it Robin} boundary condition
is coercive provided $p \ge p_0$ on $\partial\Omega$ (or at least on an open subset of $\partial\Omega$)
with some constant $p_0>0$. This is a consequence of the norm equivalence
\begin{equation}\label{E:norm-equiv-Robin}
\| v \|_{H^1(\Omega)}^2 \approx | v |_{H^1_0(\Omega)}^2 + \| v \|_{L^2(\partial\Omega)}^2 \quad
\forall \,v \in H^1(\Omega).
\end{equation}
For the {\it Neumann} problem, instead, we have $p=0$ in $\partial\Omega$, whence $\B$ is coercive whenever
$c > 0$ in $\Omega$ (or at least in an open subset of $\Omega$). If $c=0$ in $\Omega$, then $\B$ is only
coercive in the subspace of $H^1(\Omega)$ of functions with vanishing meanvalue, according to
Lemma \ref{L:poincare-friedrichs} (second Poincar\'e inequality).
Existence, uniqueness and stability is thus
guaranteed by Corollary \ref{nsv-C:lax-milgram} (Lax-Milgram) provided the forcing in \eqref{E:Robin-bilinear} satisfies the compatibility condition $\ell(1)=0$.

%----------------------------------------------------------------------
\bigskip\noindent
{\it Saddle Point Problems.}
We consider now an abstract problem a bit more general than \eqref{weak-stokes}, so the following
results apply to the Stokes system \eqref{strong-stokes}.

Given a pair of Hilbert spaces $(\V,\Q)$, we consider two
continuous bilinear forms
\begin{math}
  a\colon\V\times\V\to\R
\end{math}
and 
\begin{math}
  b\colon\Q\times\V\to\R.
\end{math}
If $f\in\V^*$ and $g\in \Q^*$, then we seek a pair
$(u,p)\in\V\times\Q$ solving the \emph{saddle point problem}
\begin{subequations}\label{saddle}
\begin{alignat}{3}
  &a[u,v] + b[p,v]\; &&= \scp{f}{v} &\quad &\forall \, v\in\V,
  \label{saddle.a}
  \\
  &b[q,u] &&= \scp{g}{q} && \forall \, q\in\Q.
\label{saddle.b}
\end{alignat}
\end{subequations}
Problem \eqref{saddle} is variational and can be rewritten in the form \eqref{var-prob-general}
\begin{equation}\label{saddle-B}
  (u,p)\in \V\times\Q: \quad
  \bilin{(u,p)}{(v,q)} =  \scp{f}{v}  + \scp{g}{q}\quad\forall \,
  ({v},q)\in \V\times\Q,
\end{equation}
where $\mathcal{B}$ is the bilinear form
\begin{equation}\label{bilinear-form}
\bilin{({u},p)}{({v},q)} \definedas
a[{u},{v}] + b[p,{v}]
+ b[q,{u}].
\end{equation}
Therefore, the saddle point problem \eqref{saddle} is well-posed
if and only if $\B$ satisfies the inf-sup condition
\eqref{nsv-inf-sup}. Since $\B$ is defined via the bilinear forms $a$
and $b$, and \eqref{saddle} has a degenerate structure, it
is not that simple to show \eqref{nsv-inf-sup} directly. However, the result is a
consequence of the inf-sup theorem for saddle point problems given
by Brezzi in 1974 \cite{Brezzi:74}.

\begin{theorem}[Brezzi]\label{T:Brezzi}
  The saddle point problem \eqref{saddle} has a unique solution 
  $(u,p)\in\V\times\Q$ for all data 
  $(f,g)\in\V^*\times \Q^*$, that depends continuously on data, if and
  only if there exist constants $\alpha,\beta>0$ such that
  \begin{subequations}\label{BB}
    \begin{align}
      \infimum_{v\in\V_0}\sup_{w\in\V_0} \frac{a[v,w]}{\|v\|_\V\|w\|_\V} =
      \infimum_{w\in\V_0}\sup_{v\in\V_0} \frac{a[v,w]}{\|v\|_\V\|w\|_\V}
      &= \alpha > 0,
      \label{BB.a}
      \\
      \infimum_{q\in\Q}\sup_{v\in\V} \frac{b[q,v]}{\|q\|_\Q\|v\|_\V}
      &= \beta > 0,
      \label{BB.b}
    \end{align}
  \end{subequations}
  where
  \begin{displaymath}
    \V_0 \definedas \{v\in\V \mid b[q,v]=0 \,\,\forall \, q\in\Q\}.
  \end{displaymath}
  In addition, there exists $\gamma=\gamma(\alpha,\beta,\|a\|)$ such that
  the solution $(u,p)$ is stable
  \begin{equation}\label{nsv-saddle-stab}
    \big(\|u\|_\V^2 + \|p\|_\Q^2\big)^{1/2}
    \le \gamma \big(\|f\|_{\V^*}^2 + \|g\|_{\Q^*}^2\big)^{1/2}.
  \end{equation}
\end{theorem}

Combining Theorem~\ref{T:Brezzi} (Brezzi) with Corollary~\ref{C:well-posedness}
(well posedness implies inf-up) we infer the inf-sup condition for the bilinear form $\B$ in
\eqref{bilinear-form}.

\begin{corollary}[inf-sup of $\B$]\label{C:brezzi-babuska}
{\it  Let the bilinear form $\B\colon\W\to\W$ be defined by \eqref{bilinear-form}.
  Then there holds
  \begin{displaymath}
    \infimum_{(v,q)\in\W}\sup_{(w,r)\in\W}
    \frac{\bilin{(v,q)}{(w,r)}}{\|(v,q)\|_\W \|(w,r)\|_\W}
    =
    \infimum_{(w,r)\in\W}\sup_{(v,q)\in\W}
    \frac{\bilin{(v,q)}{(w,r)}}{\|(v,q)\|_\W \|(w,r)\|_\W}
    \geq \gamma^{-1},
  \end{displaymath}
  where $\gamma$ is the stability constant from Theorem~\ref{T:Brezzi}.
 }
\end{corollary}

Assume that $a\colon\V\times\V\to\R$ is symmetric and let $(u,p)$ be
the solution to \eqref{saddle}. Then $u$ is the unique minimizer
of the energy
\begin{math}
  J[v] \definedas\frac12 a[v,v] - \dual{f}{v}
\end{math}
under the constraint $b[\cdot,u]=g$ in $\Q^*$. In view of this,
$p$ is the corresponding Lagrange multiplier and the pair $(u,p)$
is the unique saddle point of the Lagrangian
\begin{displaymath}
  L[v,q] \definedas J[v] + b[q,v] - \dual{g}{q}
  \quad\forall \, v\in\V, q\in\Q.
\end{displaymath}

\medskip\noindent
{\it Stokes system.}
Theorem \ref{T:Brezzi} (Brezzi) applies to the Stokes system \eqref{strong-stokes} and \eqref{weak-stokes}
once we verify the inf-sup property \eqref{BB.b} for the bilinear form
$b[q,\vec{v}]=-\int_\Omega q \divo \vec{v}$. This turns out to be equivalent to the following problem:
for any $q\in L^2_0(\Omega)$ there exists a $\vec{w}\in H^1_0(\Omega;\R^d)$
such that
\begin{equation}\label{E:Inverse-div}
  -\divo\vec{w} = q\quad\text{in } \Omega
  \qquad\text{and}\qquad
  \abs{\vec{w}}_{H^1(\Omega;\R^d)} \leq C(\Omega) \|q\|_{L^2(\Omega)}.
\end{equation}
This non-trivial result goes back to Ne\v cas \cite{CaDuFrGoGrNeSe:66}
% \todo[inline]{CC: incomplete ref. It says: cite{CaDuFrGoGrNeSe:66}. This is the quote, but I was not
% able to find in in mathscinet:
% Carroll, R., Duff, G., Friberg, J., Gobert, J., Grisvard, P., Neˇcas, J., Seeley, R.: Equations
% aux d´eriv´ees partielles. No. 19 in Seminaire de mathematiques superieures. Les Presses de
% l’Universit´e de Montr´eal (1966)}
and a proof can for instance be found in
\cite[Theorem~III.3.1]{Galdi:94}. This implies 
\begin{displaymath}
  \sup_{\bv\in H^1_0(\Omega;\R^d)}
  \frac{b[q,\bv]}{\abs{\bv}_{H^1(\Omega;\R^d)}} \ge
  \frac{b[q,\vec{w}]}{\abs{\vec{w}}_{H^1(\Omega;\R^d)}} 
  = \frac{\|q\|_{L^2(\Omega)}^2}{\abs{\vec{w}}_{H^1(\Omega;\R^d)}} 
  \ge
  C(\Omega)^{-1} \|q\|_{L^2(\Omega)}.
\end{displaymath}
Therefore, \eqref{BB.b} holds with $\beta\geq C(\Omega)^{-1}$.

% \todo[inline]{RHN (1/2/24): Should we mention Neumann boundary conditions? \\
% CC: do you mean mentioning the inf-sup for the appropriate spaces for these conditions?
% RHN (01/14/24): See (3.16) amd (3.17).}

The inf-sup condition is satisfied also for the spaces $\V$ and $\Q$ defined in \eqref{E: spaces_Neumann_stokes}, which are appropriate for the weak formulation \eqref{E:strong-Neumann-Stokes} of the Neumann boundary value problem. Indeed, given any $q \in L^2(\Omega)$, we can split it as $q=(q-\hat{q})+\hat{q}$ with $\hat{q} := |\Omega|^{-1}\int_\Omega q$. Let $\vec{w}_0\in H^1_0(\Omega;\R^d)$ be the function defined as in \eqref{E:Inverse-div} with $q$ replaced by $q-\hat{q}$, and let $\hat{\vec{w}} = \frac12(\hat{q}x,\hat{q}y)$.
Then, it is easily checked that the function $\vec{w}=\bar{\vec{w}}-|\Omega|^{-1}\int_\Omega \bar{\vec{w}}$ with 
$\bar{\vec{w}} = \vec{w}_0 + \hat{\vec{w}}$ belongs to $\V$ and satisfies 
\[
  -\divo\vec{w} = q\quad\text{in } \Omega
  \qquad\text{and}\qquad
  \|\vec{w} \|_{H^1(\Omega;\R^d)} \leq C \|q\|_{L^2(\Omega)}.
\]

%%%%%%%%%%%%%%%%%%%%%%%%%%%%%%%%%%%%%%%%%%%%%%%%%%%%%%%%%%%%%%%%%%%%%%%%%%%%%%%%%%
\subsection{$W^1_p$-regularity of reaction-diffusion equation}\label{S:reg-w1p}
%%%%%%%%%%%%%%%%%%%%%%%%%%%%%%%%%%%%%%%%%%%%%%%%%%%%%%%%%%%%%%%%%%%%%%%%%%%%%%%%%%

It will be useful later in Lemma \ref{L:perturbation} 
to know whether $L^\infty$-coefficients $(\bA,c)$ allow for enhanced
regularity beyond $H^1_0(\Omega)$ for solutions $u$ of \eqref{weak-form}. We can reformulate
this question as an extension of the Lax-Milgram theory that states that the solution operator
$L^{-1}$ of \eqref{strong-form}
is an isomorphism between $H^{-1}(\Omega)$ and $H^1_0(\Omega)$; see Corollary \ref{nsv-C:lax-milgram}
(Lax-Milgram).

This issue is well understood for the Laplace operator, i.e. $\vec{A}=\vec{I}$ and $c=0$. It is known
that for Lipschitz domains $\Omega \subset\R^d$, there exists
$p_0 = p_0(\Omega) > 2$ such that
\begin{equation}\label{E:Wp-reg-Laplacian}
  \| \nabla u \|_{L^p(\Omega)} \leq K \| f \|_{W^{-1}_p(\Omega)}
  \quad\forall p \in [2,p_0],
\end{equation}
where $K$ depends on $p$ (see for example \cite{JerisonKenig:95});
in particular, $p_0 > 4$ for $d=2$ and $p_0>3$ for $d=3$. Hereafter,
$W^{-1}_p(\Omega)$ denotes the dual space of $\mathring{W}^1_q(\Omega)$ -- functions in ${W}^1_q(\Omega)$
with zero trace and $q = \frac{p}{p-1}$.
For $\vec{A}\in L^\infty(\Omega,\R^{d \times d})$ and $c=0$, estimate \eqref{E:Wp-reg-Laplacian} was first derived by
\cite{Meyers:63} as a perturbation result for the Laplacian; see also \cite{BrennerScott:08}.
We now present a simple proof  for $p>2$ following \cite{BonitoDeVoreNochetto:2013}.
Let
\begin{equation}\label{E:theta-interpolation}
  \theta(p) := \frac{1/2 - 1/p}{1/2 - 1/{p_0}} \quad\forall  p \in [2,p_0],
\end{equation}
and note that $\theta(p)$ increases strictly from $0$ at $p=2$ to $1$ at $p=p_0$. Let $K_0$ be
the constant $K$ in \eqref{E:Wp-reg-Laplacian} for $p=p_0$ and, for any $t\in(0,1)$, define
\begin{equation}\label{E:p-star}
p_* (t) := \max\{ p \in  [2,p_0] \ : \ K_0^{\theta(p)}(1- t) \leq 1 \}.
\end{equation}

\begin{lemma}[$W^1_p$-regularity]\label{L:Wp-regularity}
  Let $\vec{A} \in L^\infty(\Omega,\R^{d \times d})$ satisfy \eqref{uniform_spd} with $\alpha_1\le\alpha_2$,
  $c\in L^\infty(\Omega)$ be non-negative, and let $\Omega$ be Lipschitz. 
  If $f\in W^{-1}_p(\Omega)$ for some $p \in \big[2, \min\big\{p_*\big(\frac{\alpha_1}{\alpha_2}\big),\frac{2d}{d-2}\big\}\big)$,  then the solution $u\in H^1_0(\Omega)$ of \eqref{weak-form}
  satisfies
  \begin{equation}\label{E:Wp-reg}
    \| \nabla u \|_{L^p(\Omega)} \leq C(p) \Big(1 + c_\B^{-1}C(\Omega)
    \|c\|_{L^\infty(\Omega)}  \Big)
    \| f \|_{W^{-1}_p(\Omega)}
  \end{equation}
   with constant
    $$
    C(p) = \frac{1}{\alpha_2}\frac{K_0^{\theta(p)}}{1-K_0^{\theta(p)}\big(1-\frac{\alpha_1}{\alpha_2}\big)}
    $$
    and $C(\Omega)=C_s C_P$ where $C_S$ is the constant in Lemma~\ref{L:embedding} (Sobolev embedding) and $C_P$ is the constant in Lemma~\ref{L:Poincare} (first Poincar\'e inequality).
\end{lemma}
\begin{proof}
  We first consider the principal part of the operator $L$ in \eqref{strong-form}, namely we take
  $c=0$. In fact, let the operator
  $S: \mathring{W}^1_p(\Omega) \to W^{-1}_p(\Omega)$ be defined by
  $Sv :=-\div\big(\frac{1}{\alpha_2} \vec{A}\nabla v \big)$. In order to prove \eqref{E:Wp-reg} for $S$,
  we resort to a perturbation argument for the Laplace operator $Tv :=-\Delta v$.

  The first task is to bound $K$ in \eqref{E:Wp-reg-Laplacian} in terms of $K_0$ and $p$. To this end,
  we recall that 
  the operator $T: H^1_0(\Omega)\to H^{-1}(\Omega)$ is an isomorphism with norm $\|T^{-1}\|_2=1$. Moreover,
  $T: \mathring{W}^1_p(\Omega)\to W^{-1}_p(\Omega)$ is also an isomorphism with norm $\|T^{-1}\|_{p_0}=K_0$
  according to \eqref{E:Wp-reg-Laplacian} for $p=p_0$, provided we adopt the norm $\|\nabla \cdot\|_{L^p(\Omega)}$
  in $\mathring{W}^1_p(\Omega)$. By the real method of interpolation, we know that
  $\mathring{W}^1_p(\Omega) = [H^1_0(\Omega),\mathring{W}^1_{p_0}(\Omega)]_{\theta(p),p}$
  is the interpolation space between $H^1_0(\Omega)$ and $\mathring{W}^1_{p_0}(\Omega)$ with parameter
  $\theta(p)$ given by \eqref{E:theta-interpolation}. Hence, operator interpolation theory implies that
  $T:\mathring{W}^1_p(\Omega) \to W^{-1}_p(\Omega)$ is an isomorphism with
  \[
  K = \|T^{-1}\|_p = K_0^{\theta(p)}.
  \]

  We regard $S$ as a perturbation of $T$, define the operator
  $Q:=T-S : \mathring{W}^1_p(\Omega) \to W^{-1}_p(\Omega)$, and observe that $\|Q\|_p\le 1-\frac{\alpha_1}{\alpha_2}$
  because
  \[
  \langle Qv,w \rangle = \int_\Omega \Big(I-\frac{1}{\alpha_2} \vec{A}\Big)\nabla v \nabla w
  \le \Big(1-\frac{\alpha_1}{\alpha_2}\Big) \|\nabla v\|_{L^p(\Omega)} \|\nabla w\|_{L^q(\Omega)}
  \quad w\in \mathring{W}^1_q(\Omega).
  \]
  Therefore, the operator $T^{-1}Q : \mathring{W}^1_p(\Omega) \to \mathring{W}^1_p(\Omega)$ satisfies
  \[
  \|T^{-1}Q\|_p \le \|T^{-1}\|_p \|Q\|_p \le K_0^{\theta(p)} \Big(1-\frac{\alpha_1}{\alpha_2}\Big)
  \]
  as well as $\|T^{-1}Q\|_p <1$ for any $p\in \big[2,p_*(\alpha_1/\alpha_2)\big)$ in
  view of definition \eqref{E:p-star} of $p_*(t)$. We conclude by the Neumann Theorem that the operator
  $S = T \big(I - T^{-1}Q \big) : \mathring{W}^1_p(\Omega) \to W^{-1}_p(\Omega)$ is invertible and its norm
  is bounded by
  \[
  \|S^{-1}\|_p \le \|T^{-1}\|_p \|I - T^{-1}Q\|_p \le \frac{\|T^{-1}\|_p}{1-\|T^{-1}Q\|_p}
  \le \frac{K_0^{\theta(p)}}{1- K_0^{\theta(p)} \big( 1 - \alpha_1 / \alpha_2 \big)}.
  \]
  This yields the asserted estimate for $S=-\div\big(\frac{1}{\alpha_2} \vec{A}\nabla v \big)$.

  Finally, we consider the operator  $L$ in \eqref{strong-form} with $c\ne0$. If $u\in H^1_0(\Omega)$ is the solution
  of \eqref{weak-form} given by Corollary \ref{nsv-C:lax-milgram} (Lax-Milgram), rewrite \eqref{weak-form}
  as follows
  \[
  S u = -\div \Big(\frac{1}{\alpha_2} \vec{A} \nabla u \Big) = \frac{1}{\alpha_2} (f-cu) = \frac{1}{\alpha_2} g,
  \]
  and apply the preceding estimate for $S$ to infer that
  \[
  \| S^{-1}g \|_{\mathring{W}^1_p(\Omega)} = \|\nabla u\|_{L^p(\Omega)} \le C(p) \|g\|_{W_p^{-1}(\Omega)}
  \le C(p) \big( \|f\|_{W_p^{-1}(\Omega)} + \|cu\|_{W_p^{-1}(\Omega)} \big).
  \]
  It remains to estimate the last term on the right-hand side. Using Cauchy-Schwarz in conjunction with Lemma~\ref{L:Poincare} (first Poincar\'e inequality), i.e., $\| u \|_{L^2(\Omega)} \leq C_P \| \nabla u \|_{L^2(\Omega)}$, and Lemma~\ref{L:embedding} (Sobolev embeddings), i.e., $\| w \|_{L^2(\Omega)} \leq C \| w \|_{\mathring{W}^1_q(\Omega)} \leq C_S \| \nabla w \|_{L^q(\Omega)}$ for $q=p/(p-1) \geq 2d/(d+2)$, we see that 
  \[
  \langle cu,w \rangle \le C_P C_S \|c\|_{L^\infty(\Omega)} \|\nabla u\|_{L^2(\Omega)} \|\nabla w\|_{L^q(\Omega)}
  \le \frac{C_PC_S}{c_\B}  \|c\|_{L^\infty(\Omega)} \|f\|_{H^{-1}(\Omega)} \|\nabla w\|_{L^q(\Omega)}
  \]
  because of the stability estimate \eqref{cont-stab} with constant $\alpha=c_\B$. Since
  $\|f\|_{H^{-1}(\Omega)} \le |\Omega|^{\frac{p-2}{2p}}\|f\|_{W_p^{-1}(\Omega)}$, the asserted estimate \eqref{E:Wp-reg}
  for $c\ne0$ follows immediately.
\end{proof}

%--------------------------------------------------------------------------------

\section{A Priori Approximation Theory}\label{S:tools}
 %\rhn{(RHN $\longrightarrow$ CC)} %\cite{NochettoVeeser:2012}

% {\color{brown} 
% \begin{itemize}
% \item
%   The Sobolev number: Scaling and Embedding.

% \item
%   Conforming meshes: the bisection method.
  
% \item
%   Polynomial interpolation, localized $H^1$-interpolation estimates \cite{Veeser:2016}.

% \item
%   Poincar\'e-Friedrichs estimates, trace inequalities.

% \item
%   The finite element method: inf-sup stability, best approximation, examples.

% \item
%   Adaptive approximation: greedy algorithm and approximation class.

% \item
%   Nonconforming meshes.

% \end{itemize}
% }

We devote this section to discussing basic concepts about
piecewise polynomial approximation in
Sobolev spaces over graded meshes in any dimension $d$.
We start by introducing the Petrov-Galerkin method in an abstract setting (Sect. \ref{S:Galerkin}); this motivates our interest in approximation results in Sobolev spaces. Hence, 
we briefly discuss the construction of finite element spaces in
Sect. \ref{S:fe-spaces}, along with polynomial interpolation of functions
in Sobolev spaces in Sect. \ref{S:pol-interp}. This provides local estimates
adequate for the comparison of quasi-uniform and graded meshes for $d>1$. We
exploit them in developing the so-called error equidistribution
principle in Sect. \ref{S:error-equidistribution} 
and the construction of suitably graded meshes via
a greedy algorithm in Sect. \ref{S:constructive-approximation}. We conclude that graded
meshes can deliver optimal interpolation rates for certain classes of
singular functions, and thus supersede quasi-uniform refinement.

In the second part of the section, we explore the geometric aspects of mesh refinement for conforming
meshes in Sect. \ref{S:bisection} and nonconforming meshes in Sect.
\ref{S:nonconforming-meshes}, but postpone a full and rather technical discussion until
Sect. \ref{S:mesh-refinement}. We include a statement about
complexity of the refinement procedure, which turns out to be
instrumental later and will be proved in Sect. \ref{S:mesh-refinement}.

%-------------------------------------------------------------------------------------
\subsection{The Galerkin method: best approximation}\label{S:Galerkin}
%-------------------------------------------------------------------------------------
%
In order to render the variational problem \eqref{var-prob-general} computable, we let
$\V_N\subset\V$ and $\W_N\subset\W$ be subspaces with the same dimension $N<\infty$ and
consider the Petrov-Galerkin approximation:
\begin{equation}\label{E:Galerkin}
u_N \in \V_N: \quad
\bilin{u_N}{w} = \dual{f}{w} \quad\forall \, w \in \W_N.
\end{equation}
If $\V_N=\W_N$ this is called Galerkin approximation. Since \eqref{E:Galerkin} is a square algebraic
system, existence and uniqueness of $u_N\in\V_N$ is equivalent to the kernel of the corresponding
linear discrete operator to be trivial. This leads to the following equivalent conditions for unique
solvability \cite{Necas:62,babuvska1971error}, see also  \cite[Proposition 1]{NoSiVe:09}.

\begin{lemma}[discrete inf-sup condition]\label{L:inf-sup_N}
{\it \!\! The following statements are equivalent: \looseness=-1
  \begin{enumerate} %[\upshape(1)]
  \setlength{\itemindent}{-0.6em}  
  \item[1.]  %\label{nsv-non-zero_N}
  for every $0\ne v\in\VN$ there exists $w\in\WN$ such that
    $\bilin{v}{w} \ne 0$;
\item[2.] %\label{nsv-inf-sup_N5}
    for every $0\ne w\in\WN$ there exists $v\in\VN$ such that
        $\bilin{v}{w} \ne 0$;
 \item[3.] %\label{nsv-inf-sup_N1}
 the following discrete inf-sup condition holds with a constant $\beta_N>0$
 \begin{equation}\label{E:discrete-inf-sup}
  \infimum_{v\in\VN}\sup_{w\in\WN}\frac{\bilin{v}{w}}{\|v\|_\V \|w\|_\W}
  =
  \infimum_{w\in\WN}\sup_{v\in\VN}\frac{\bilin{v}{w}}{\|v\|_\V \|w\|_\W}
  = \beta_N;
\end{equation}
\item[4.] %\label{nsv-inf-sup_N2}
    \begin{math}
      \displaystyle
      \infimum_{v\in\VN}\sup_{w\in\WN}\frac{\bilin{v}{w}}{\|v\|_\V \|w\|_\W}>0;
    \end{math}
  \item[5.] %\label{nsv-inf-sup_N3}
    \begin{math}
      \displaystyle
      \infimum_{w\in\WN}\sup_{v\in\VN}\frac{\bilin{v}{w}}{\|v\|_\V \|w\|_\W}>0.
    \end{math}
  \end{enumerate}
 } 
\end{lemma}

This is a discrete version of Theorem \ref{nsv-T:ex1-var-prob} (Ne\v{c}as) and leads to the stability bound
\begin{equation}\label{E:discrete-stability}
\|u_N\|_{\V} \leq \frac{1}{\beta_N} \|f\|_{\W^*}.
\end{equation}
Therefore, $\beta_N^{-1}$ acts as a stability constant for \eqref{E:Galerkin} and is thus 
desirable that it is uniformly bounded below away from zero. This is always the case for
{\it coercive} problems because \eqref{nsv-coercivity} is inherited within the subspaces
$\VN=\WN\subset\V$, whence $\beta_N\ge\alpha>0$. In contrast, a uniform lower bound for {\it saddle point} problems
\begin{equation}\label{E:uniform-inf-sup}
\beta_N \ge \beta >0,
\end{equation}
requires compatibility between the subspaces $\VN$ and $\WN$ \cite{BoffiBrezziFortin}. 
%We will come back to this issue in Section \ref{S:conv-rates-infsup}.

If we now subtract \eqref{E:Galerkin} from \eqref{var-prob-general}, we obtain {\it Galerkin orthogonality}
\begin{equation}\label{E:Galerkin-orthogonality}
\bilin{u-u_N}{w} = 0 \quad\forall \, w\in\WN.
\end{equation}
This relation is instrumental to derive the following best approximation property as
well to develop a posteriori error estimates in Section \ref{S:aposteriori}.

\begin{proposition}[quasi-best-approximation property]\label{P:QBAP}
  Let $\B\colon\V\times\W\to\R$ be continuous and satisfy
  \eqref{E:discrete-inf-sup}.  Then the error $u-\uN$ satisfies the bound
  \begin{equation}\label{E:best-approximation}
    \|u-\uN\|_\V \leq \frac{\|\B\|}{\beta_N} \min_{v\in\VN} \|u-v\|_\V.
  \end{equation}
\end{proposition}

\begin{proof}
We give a simplified proof, which follows Babu{\v s}ka
\cite{Babuska:71,BabuskaAziz:72} and yields the constant 
$1+\frac{\|\B\|}{\beta_N}$. The asserted constant is due to Xu and Zikatanov \cite{XuZikatanov:03}.

Combining \eqref{E:discrete-inf-sup}, \eqref{E:Galerkin-orthogonality}, and the
continuity of $\B$, we derive for all $v\in\V_N$
\[
 \beta_N\|u_N-v\|_\V
 \le
 \sup_{w\in\W_N}\frac{\bilin{u_N-v}{w}}{\|w\|_\W}
 =
 \sup_{W\in\W_N}\frac{\bilin{u-v}{w}}{\|w\|_\W}
 \le
 \|\B\| \|u-v\|_\V,
\]	
whence
\[
 \|u_N-v\|_\V\le\frac{\|\B\|}{\beta_N}\|u-v\|_\V.
\]
Using the triangle inequality yields
\[
 \|u-u_N\|_\V\le\|u-v\|_\V+\|v-u_N\|_\V
 \le
 \bigg(1+\frac{\|\B\|}{\beta_N}\bigg)\|u-v\|_\V
\]
for all $v\in\V_N$. It just remains to minimize in $\V_N$.
\end{proof}

\begin{corollary}[quasi-monotonicity]\label{C:quasi-monotonicty}
  Let $\B\colon\V\times\W\to\R$ be continuous and satisfy
  \eqref{E:discrete-inf-sup}. If $\V_M$ is a subspace of $\V_N$, then for all $v\in\V_M$
\begin{equation}\label{E:quasi-monotonicity-infsup}
\| u - u_N \|_\V \le 
\frac{\|\B\|}{\beta_N} \, 
%\Big(1+\frac{\|\B\|}{\beta_N} \Big) 
\|u-v\|_\V.
\end{equation}
Moreover, if $\V=\W$ and $\B$ is symmetric and coercive with constants $c_\B\le C_\B$,
then for all $v\in\V_M$
\begin{equation}\label{E:quasi-monotonicity-coercive}
\enorm{u-u_N} \le \enorm{u-v},
\quad \|u-u_N\|_\V \le C_\mathrm{cea} \|u-v\|_\V,
\end{equation}
\index{Constants!$C_\textrm{cea}$: best approximation constant}
where $C_\textrm{cea}:=\sqrt{\frac{C_\B}{c_\B}}$. 
\end{corollary}
\begin{proof}
Inequality \eqref{E:quasi-monotonicity-infsup} is a consequence of the previous
bound \eqref{E:best-approximation} upon taking $v\in\V_M$ instead of $\V_N$. To show the left inequality in
\eqref{E:quasi-monotonicity-coercive}, we combine \eqref{nsv-coercivity} and
\eqref{E:Galerkin-orthogonality}
\[
\enorm{u-u_N}^2 = \B[u-u_N, u-v] \le \enorm{u-u_N} \enorm{u-v}, \qquad \forall v \in \mathbb V_M.
\]
This together with the norm equivalence \eqref{nsv-enorm-Vnorm} gives the remaining
inequality.
\end{proof}

The significance of \eqref{E:best-approximation} is that we need to construct discrete spaces
$\VN$ with good approximation properties. We introduce next piecewise polynomial
approximation, which gives rise to the finite element method.

%-------------------------------------------------------------------------------------
\subsection{Finite element spaces}\label{S:fe-spaces}
%-------------------------------------------------------------------------------------
%
In this section we focus on the construction of the discrete spaces $\VN$ and $\WN$.
We consider the bilinear forms $\mathcal{B}$ introduced in Section \ref{S:ex-bvp} with emphasis on
the diffusion-reaction case \eqref{bilinear-coercive}. We assume that $\Omega$ is a
bounded polyhedral domain $\Omega\subset\R^d$ and is partitioned into a conforming or non-conforming mesh $\grid$
made of simplices $T$, which are assumed to be closed with non-overlapping interiors; thus,
\[
\overline{\Omega} = \bigcup \{T: T\in\grid\}.
\]
The reference element is denoted 
$$
\index{Meshes!$T_d$: reference element}
T_d := \big\{ x=(x_1,\dots,x_d) \in \mathbb R^d \mid 0\leq x_i \leq 1, \ i=1,..,d, \ \sum_{i=1}^d x_i \leq 1 \big\}.
$$
We will discuss the construction of conforming meshes in Section \ref{S:bisection} by the bisection
method and that of non-conforming meshes (constrained to the fulfillment of an {\em admissibility condition}) in Section \ref{S:nonconforming-meshes}, both for $d=2$. 
We will embark on a thorough discussion in Section \ref{S:mesh-refinement}. We assume for the moment 
that $\grid$
is an element of a (possibly infinite) class $\grids$ of conforming {\it shape regular} meshes. To define this
geometric concept, we denote by $\overline{h}_T$ the diameter of $T\in\grid$, by $\underline{h}_T$
the diameter of the largest ball contained in $T$, and impose the restriction
\begin{equation}\label{E:shape-regularity}
\index{Constants!$\sigma$: shape regularity constant of $\grids$}
    \sigma := \sup_{\grid\in\grids} \sup_{T\in\grid} \frac{\overline{h}_T}{\underline{h}_T} < \infty.
\end{equation}
The constant $\sigma$ is refereed as the shape regularity constant of $\grids$.

Given a shape regular mesh $\grid\in\grids$, we define the finite element space of
{\it discontinuous} piecewise polynomials of total degree up to $n\ge1$
\[
\index{Functional Spaces!${\mathbb S}^{n,-1}_\grid$: piecewise polynomials of degree $\leq n$}
{\mathbb S}^{n,-1}_\grid := \{v\in L^2(\Omega) \ | \ 
v|_T \in \P_n(T) ~~\forall T\in\grid\};
\]
and its globally {\it continuous} counterpart
\[
\index{Functional Spaces!${\mathbb S}^{n,0}_\grid$: globally continuous piecewise polynomials of degree $\leq n$}
{\mathbb S}^{n,0}_\grid := \mathbb S^{n,-1}_\grid \cap  C^0(\overline{\Omega}).
\]
%
% \todo[inline]{AB: the finite element spaces are defined on the physical element, which is equivalent to define it on the master element when affine transformations are used. However, this is not the case for quadrilaterals (not introduced yet). Shall we consider quad already here or make a comment later? \\
% CC: For better clarity, I prefer to postpone quads}
Note that ${\mathbb S}^{n,0}_\grid\subset H^1(\Omega)$ which makes it adequate for 
\eqref{weak-form}-\eqref{bilinear-coercive}.  We refer to
Braess \cite{Braess:07}, Brenner-Scott \cite{BrennerScott:08}, Ciarlet
\cite{Ciarlet:02} and Siebert \cite{Siebert:12} for a discussion on
the local construction of this space (i.e. Lagrange elements of degree $n\ge1$)
along with its properties. We denote by
\begin{equation}\label{E:def_Vgrid}
\index{Functional Spaces!$\V_\grid$: conforming finite element space}
\V_\grid := {\mathbb S}^{n,0}_\grid \cap H^1_0(\Omega)
\end{equation}
the subspace of finite element functions which vanish on
$\partial\Omega$. Note that we do not explicitly refer to the
polynomial degree, which will be clear in each context.

We focus on the conforming piecewise linear case $n=1$ (Courant
elements), but most results extend to non-conforming meshes or $n>1$. In this vein, global continuity can be simply enforced by imposing continuity at the set $\vertices$\index{Meshes!$\vertices$: set of vertices} of vertices $z$ of $\grid$, the so-called {\it nodal values}.
However, the following local construction leads to global continuity.
If $T$ is a generic simplex
of $\grid$, namely the convex hull of $\{z_i\}_{i=0}^d$,
then we associate to each vertex $z_i$ a {\it barycentric
coordinate} $\lambda_i^T$, which is the linear function in $T$ with nodal
value $1$ at $z_i$ and $0$ at the other vertices of $T$. Upon pasting
together the barycentric coordinates $\lambda_z^T$ of all simplices $T$
containing the vertex $z\in\mathcal V$, we obtain a continuous piecewise
linear function $\phi_z\in{\mathbb S}^{1,0}_\grid$ \index{Functions!$\phi_z$: Lagrange basis of $\mathbb S_\grid^{1,0}$}
as depicted in Figure \ref{nv-F:basis} for $d=2$:
\begin{figure}[h]
\begin{center}
\includegraphics[scale=0.5,trim={5cm 6cm 11cm 4cm},clip]{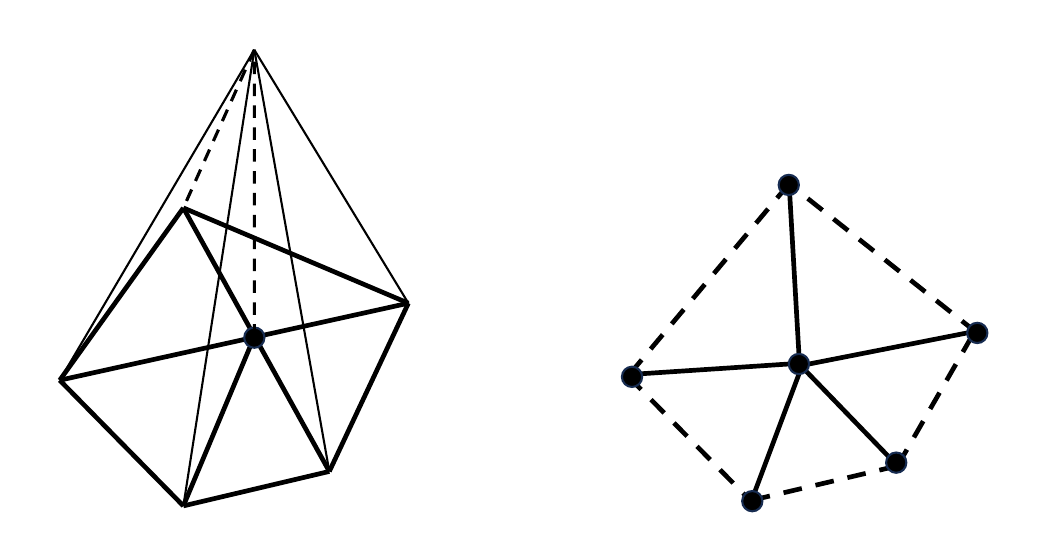}
\end{center}
\begin{picture}(0,0)
\put(117,70){$z$}
\put(132,120){$\phi_z$}
\put(210,90){$\omega_z$}
\put(250,65){$z$}
\end{picture}
\caption{(Left) Piecewise linear basis function $\phi_z$ corresponding to
  interior node $z$; (Right) Support $\omega_z$ of $\phi_z$ and skeleton $\gamma_z$ (in solid line)}
\label{nv-F:basis}
\end{figure}
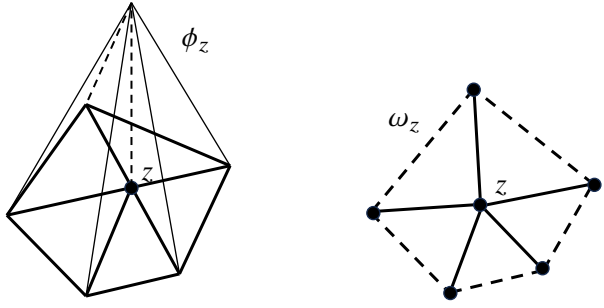

The set $\{\phi_z\}_{z\in\mathcal V}$ of all such functions
is the nodal basis of ${\mathbb S}^{1,0}_\grid$, or Courant basis.
We denote by $\omega_z := \supp (\phi_z)$ \index{Meshes!$\omega_z$: region made of elements sharing the vertex $z$} the support of
$\phi_z$, from now on called {\it star} associated to $z$, and by
$\gamma_z$\index{Meshes!$\gamma_z$: skeleton of $\omega_z$} the interior skeleton of $\omega_z$, namely the union of all the faces
containing $z$.

In view of the definition of $\phi_z$,
we have the following unique representation of any function
$v\in{\mathbb S}^{1,0}(\grid)$
\[
v(x) = \sum_{z\in\vertices} v(z) \phi_z(x).
\]
The functions $\phi_z$ are non-negative and satisfy the {\it partition of unity} property
\begin{equation}
    \sum_{z\in\vertices} \phi_z(x) = 1 \quad\forall \, x\in \Omega.
\end{equation}
If we further impose $v(z)=0$ for all $z\in\partial\Omega\cap\vertices$,
then $v\in H^1_0(\Omega)$.  

For each simplex $T\in\grid$, generated by vertices $\{z_i\}_{i=0}^d$,
the {\it dual functions} $\{\lambda_i^*\}_{i=0}^d\subset\P_1(T)$ to the
barycentric coordinates $\{\lambda_i\}_{i=0}^d$ satisfy the
bi-orthogonality relation $\int_T \lambda_i^*\lambda_j=\delta_{ij}$,
and are given by
\[
\lambda_i^* = \frac{(1+d)^2}{|T|} \lambda_i - \frac{1+d}{|T|}
\sum_{j\ne i} \lambda_j
\qquad\forall ~0\le i \le d.
\]
The {\it Courant dual basis} $\phi_z^*\in{\mathbb S}^{1,-1}(\grid)$ are the 
functions over $\grid$ given by
\[
\phi_z^* = \frac{1}{V_z} \sum_{T\ni z} (\lambda_z^T)^* \chi_T
\qquad\forall ~z\in\vertices,
\]
where $V_z\in\N$ is the valence of $z$ (number of elements of
$\grid$ containing $z$) and $\chi_T$ is the characteristic function of
$T$. These functions have the same support $\omega_z$ as the nodal
basis $\phi_z$ and satisfy the global bi-orthogonality relation
\[
\int_\Omega \phi_z^* \phi_y = \delta_{zy}
\qquad\forall ~z,y\in\mathcal V.
\]

Finally, we denote by $\nodes$\index{Meshes!$\nodes$: Lagrange nodes of order $n$} the Lagrange nodes of order $n$ of a mesh $\grid$, and by 
$\psi_z \in\mathbb{S}^{n,0}_\grid$ \index{Functions!$\psi_z$: Lagrange basis of $\mathbb S_\grid^{n,0}$} the corresponding Lagrange basis of $\mathbb{S}^{n,0}_\grid$;
hence $\mathbb{S}^{n,0}_\grid = \Span{\psi_z}_{z\in\nodes}$.

%-------------------------------------------------------------------------------------
\subsection{Polynomial interpolation in Sobolev spaces}\label{S:pol-interp}
%-------------------------------------------------------------------------------------
%
We wish to use the space $\V_\grid$ defined in \eqref{E:def_Vgrid} as the discrete space $\V_N$ in the Galerkin method \eqref{E:Galerkin}. If the bilinear form $\B$ satisfies an inf-sup condition with constant $\beta_\grid >0$, we find a discrete solution $u_\grid \in \V_\grid$ which satisfies the error bound \eqref{E:best-approximation}, i.e.
\[
 \|u-u_\grid\|_\V \leq \frac{\|\B\|}{\beta_\grid} \min_{v\in\V_\grid} \|u-v\|_\V.
\]
In turn, the minimum on the right-hand side can be bounded from above by the quantity $\|u-v\|_\V$ for any chosen $v \in \V_\grid$. This motivates the search for quasi-best approximations of $u$ in the norm of $\V$.
One classical tool to generate approximations to any given function is {\em interpolation}. Interpolation in $\V_\grid$ is discussed next.

\smallskip
If $v\in C^0(\overline{\Omega})$ we define the {\it Lagrange
  interpolant} $I_\grid v$ of $v$ as follows:
\[
I_\grid v(x) = \sum_{z\in\nodes} v(z) \psi_z(x),
\]
and note that $I_\grid v = v$ for all $v\in {\mathbb S}^{n,0}_\grid$ (i.e. $I_\grid$ is invariant
in ${\mathbb S}^{n,0}_\grid$).
For functions without point values, such as functions in $H^1(\Omega)$
for $d>1$, we need to determine nodal values by averaging. 
For any conforming mesh $\grid\in\grids$,
the averaging process extends beyond nodes and so gives rise to the
discrete neigborhood
\begin{displaymath}
\index{Meshes!$\omega_T$, $\omega_\grid(T)$: region of elements intersecting $T$}
  \omega_\grid(\elm)\definedas \bigcup_{\substack{T' \in \grid \\ T' \cap T \not = \emptyset}} T'
\end{displaymath}
for each element $\elm\in\grid$ along
with the {\it local quasi-uniformity} properties
\begin{equation*}
  \max_{\elm\in\grid}
  \#\omega_\grid(\elm) 
  \leq C(\sigma),
  \qquad
  \max_{\elm'\subset \omega_\grid(\elm)}
  \frac{\vol{\elm}}{\vol{\elm'}} \leq C(\sigma),
\end{equation*}
where $\sigma$ is the shape regularity coefficient defined in \eqref{E:shape-regularity}.
We will often write $\omega_T$ is there is no confusion about the underlying mesh $\grid$.
We shall also need the a smaller subset, namely the set of elements sharing a face with a given element $T$:
\begin{equation}
\label{tildeomegaT}
\index{Meshes!$\omega_T$, $\omega_\grid(T)2$@$\wt{\omega}_T$, $\wt{\omega}_\grid(T)$: elements sharing a face with $T$}
 \wt{\omega}_T := \wt{\omega}_\grid(T):= \bigcup_{\substack{T' \in \grid \\ T' \cap T \in \faces}} T',
\end{equation}
where $\mathcal F$ is the set of all $(d-1)$-dimensional faces of the mesh $\mathcal T$.

We introduce now one such operator $I_\grid$ due to Scott-Zhang
\cite{BrennerScott:08,ScottZhang:90}, from now on called a {\it
  quasi-interpolation operator}. We focus on polynomial degree $n=1$,
but the construction is valid for any $n$; see
\cite{BrennerScott:08,ScottZhang:90} for details.  We recall that
$\{\phi_z\}_{z\in\mathcal V}$ is the global Lagrange basis of ${\mathbb
  S}^{1,0}_\grid$, $\{\phi_z^*\}_{z\in\mathcal V}$ is the global dual
basis, and $\supp\phi_z^*=\supp\phi_z$ for all $z\in\mathcal V$.  We thus
define $I_\grid: L^1(\Omega) \to {\mathbb S}^{1,0}_\grid$ to be
\begin{equation}\label{E:quasi-interpolant}
\index{Operators!$I_\grid$: quasi-interpolation operator}
I_\grid v := \sum_{z\in\vertices} \dual{v}{\phi_z^*} \phi_z,
\end{equation}
If $0 \le s \le 2$ (integer) is a regularity index and
$1\le p\le\infty$ is an integrability index, then
we would like to prove the {\it quasi-local error estimate}
\begin{equation}
\label{E:quasi-interp-error}
 |v-I_\grid v|_{W^t_q(\elm)}
 \lesssim
 \helm^{\sob(W^s_p)-\sob(W^t_q)} |v|_{W^s_p(\omega_\elm)}
\end{equation}
for all $T\in\grid$,
provided $0\le t\le s$, $1\leq q \leq\infty$ are such that
$\sob(W^s_p) > \sob(W^t_q)$.

We first recall that $I_\grid$ is invariant in
${\mathbb S}^{1,0}_\grid$, namely,
\[
 I_\grid w=w \quad \fa w\in {\mathbb S}^{1,0}_\grid.
\]
Since the averaging process giving rise to the values of
$I_\grid v$ for each element $\elm\in\grid$ takes place in the
neighborhood $\omega_\elm$, we also deduce the local
invariance
\[
 I_\grid w|_\elm = w \quad \fa w\in \P_1(\omega_T)
\]
as well as the local stability estimate for any $1\le q \le \infty$
\[
\|I_\grid v\|_{L^q(\elm)} \Cleq \|v\|_{L^q(\omega_\elm)}.
\]
We may thus write
\[
 v - I_\grid v|_\elm
 =
 (v-w) - I_\grid (v-w)|_\elm
\quad\fa \elm\in\grid,
\]
where $w\in\P_{s-1}$ is arbitrary ($w=0$ if $s=0$).  
%
%\todo[inline]{RHN (01/09/24): Check if $s$ is used as a polynomial degree and correct.}
%
It suffices now to prove \eqref{E:quasi-interp-error} in the
reference element $\refelm$ and scale back and forth to $T$; 
the definition \eqref{nv-sobolev-number} of Sobolev
number accounts precisely for this scaling. 
We keep the notation $\elm$ for $\refelm$,
apply the inverse estimate for linear polynomials
$
 |I_\grid v|_{W^t_q(\elm)}
 \lesssim
 \|I_\grid v\|_{L^q(\elm)}
$
to $v-w$ instead of $v$, and use the above local stability
estimate, to infer that
\[
 |v-I_\grid v|_{W^t_q(\elm)}
 \lesssim
 \|v-w\|_{W^t_q(\omega_\elm)}
 \lesssim
 \|v-w\|_{W^s_p(\omega_\elm)}.
\]
The last inequality is a consequence of the inclusion
$W^s_p(\omega_T) \subset W^t_q(\omega_T)$ because
$\sob(W^s_p) > \sob(W^t_q)$ and $t\leq s$.
Estimate \eqref{E:quasi-interp-error} now follows from the
Bramble-Hilbert lemma \cite[Lemma 4.3.8]{BrennerScott:08}, \cite[Theorem 3.1.1]{Ciarlet:02}, 
\cite{DupontScott:80} or Proposition \ref{P:Bramble-Hilbert-Sobolev} below:
\begin{equation}
\label{E:local-BrambleHilbert}
 \inf_{w\in\P_{s-1}(\omega_\elm)} \| v-w \|_{W^s_p(\omega_\elm)}
 \lesssim
 |v|_{W^s_p(\omega_\elm)}.
\end{equation}
This proves \eqref{E:quasi-interp-error} for $n=1$. The construction
of $I_\grid$ and ensuing estimate \eqref{E:quasi-interp-error} are still valid
for any $n>1$ \cite{BrennerScott:08,ScottZhang:90}.

\begin{proposition}[quasi-interpolant without boundary values]
\label{P:quasi-interp}
{\it Let $s,t$ be regularity indices with $0\leq t \le s\leq n+1$, and
$1\le p,q\le\infty$  be integrability indices so that
$\sob(W^s_p) > \sob(W^t_q)$.
Then there exists a quasi-interpolation operator 
$I_\grid: L^1(\Omega)\to {\mathbb S}^{n,0}_\grid$, which is
invariant in ${\mathbb S}^{n,0}_\grid$ and satisfies
\begin{equation}
\label{E:quasi-interp-error-bis}
 |v-I_\grid v|_{W^t_q(\elm)}
 \lesssim
 \helm^{\sob(W^s_p)-\sob(W^t_q)} |v|_{W^s_p(\omega_\elm)}
 \quad\forall T\in\grid.
\end{equation}
The hidden constant in \eqref{E:quasi-interp-error}
depends on the shape coefficient of $\gridk[0]$ and $d$.
}
\end{proposition}

To impose a vanishing trace on $I_\grid v$ we may suitably modify
the averaging process for boundary nodes.
We thus define a set of dual functions with respect 
to an $L^2$-scalar product over ${(d-1)}$-subsimplices contained on 
$\partial\Omega$; see again \cite{BrennerScott:08,ScottZhang:90} 
for details. This retains the invariance
property of $I_\grid$ on ${\mathbb S}^{n,0}(\grid)$ and guarantees that
$I_\grid v$ has a zero trace if $v\in W^1_1(\Omega)$ does.
Hence, the above argument applies and
\eqref{E:quasi-interp-error-bis} follows provided $s\geq1$.

\begin{proposition}[quasi-interpolant with boundary values]
\label{P:quasi-interp-bc}
{\it Let $s,t,p,q$ be as in Proposition \ref{P:quasi-interp}.
There exists a quasi-interpolation operator 
$I_\grid: W^1_1(\Omega)\to {\mathbb S}^{n,0}_\grid$ 
invariant in ${\mathbb S}^{n,0}_\grid$ which
satisfies \eqref{E:quasi-interp-error-bis} for
$s\ge1$ and preserves the boundary values of $v$ provided they are piecewise
polynomial of degree $\le n$. In particular, if $v\in W^1_1(\Omega)$ has a
vanishing trace on $\partial\Omega$, then so does $I_\grid v$.
}
\end{proposition}

\begin{remark}[fractional regularity]\label{nv-R:fractional}
  We observe that \eqref{E:quasi-interp-error-bis} does not require the
  regularity indices $t$ and $s$ to be integer.  The proof follows along the
  same lines but replaces the polynomial degree $n$ by the greatest
  integer smaller than $s$; the generalization of
  \eqref{E:local-BrambleHilbert} can be taken from
  \cite{DupontScott:80}.
\end{remark}
%
%  \todo[inline]{RHN (01/09/24): Revise use of $s$ as a polynomial degree.}
%  

\begin{remark}[local error estimate for Lagrange interpolant]\label{nv-R:lagrange}
Let the regularity index $s$ and integrability index $1\le p\le \infty$ satisfy
$s-d/p>0$. This implies that 
$\sob(W^s_p)>\sob(L^\infty)$, whence
$W^s_p(\Omega)\subset C(\overline{\Omega})$ and
the Lagrange interpolation operator 
$I_\grid: W^s_p(\Omega)\to {\mathbb S}^{n,0}_\grid$ is well defined and
satisfies the {\it local error estimate}
\begin{equation}\label{E:error-lagrange}
 |v-I_\grid v|_{W^t_q(\elm)}
 \lesssim
 \helm^{\sob(W^s_p)-\sob(W^t_q)} | v |_{W^s_p(\elm)},
\end{equation}
provided $0\le t\le s$, $1\leq q \leq\infty$ are such that
$\sob(W^s_p)>\sob(W^t_q)$.
We point out that $\omega_\elm$ in \eqref{E:quasi-interp-error} is
now replaced by $\elm$ in \eqref{E:error-lagrange}. We also remark
that if $v$ vanishes on $\partial\Omega$ so does $I_\grid v$. The
proof of \eqref{E:error-lagrange} proceeds along the same lines as
that of Proposition \ref{P:quasi-interp} except that the nodal
evaluation does not extend beyond the element $\elm\in\grid$
and the inverse and stability estimates over the reference element
are replaced by
\[
|I_\grid v|_{W^t_q(\refelm)} \Cleq \|I_\grid v\|_{L^q(\refelm)}
\Cleq \|v\|_{L^\infty(\refelm)}
\Cleq \|v\|_{W^s_p(\refelm)}.
\]
\end{remark}

The following global interpolation error estimate builds on Proposition \ref{P:quasi-interp}
and relates to Fig.~\ref{F:devore-diagram} (DeVore diagram).

\begin{theorem}[global interpolation error estimate]\label{T:global-interp}
{\it Let $1 \le s \le n+1$, $t=0,1$, $t<s$, and $1 \le p \le q$ satisfy
$r:=\sob(W^s_p) - \sob(W^t_q) > 0$. If $v\in W^s_p(\Omega)$, then
\begin{equation}\label{E:global-interp}
|v-I_\grid v|_{W^t_q(\Omega)} \Cleq \Big(\sum_{T\in\grid} h_T^{rp} |v|_{W^s_p(\omega_T)}^p\Big)^{\frac1p}.
\end{equation}
}
\end{theorem}
\begin{proof}
Use Proposition \ref{P:quasi-interp} along with the elementary
property of series
$
\sum_n a_n\le (\sum_n a_n^{p/q})^{q/p}
$
for $0<p/q \le 1$.
\end{proof}
%

%--------------------------------------------------------------------------------
\paragraph{\it Continuous vs discontinuous approximation of gradients.}
%--------------------------------------------------------------------------------
The preceding discussion might induce us to believe that when dealing with Sobolev functions
$v\in W^1_p(\Omega)$ without pointvalues, namely $1\le p\le d$, global continuity of the quasi-interpolant
$I_\grid v$ might degrade the approximation quality relative to discontinuous approximations. The following
instrumental result shows that this is not the case \cite[Theorem~2]{Veeser:2016}.
It hinges on a new geometric concept: we say that a star $\omega_z$ is \emph{$(d-1)$-faced connected}\index{Definitions! faced-connected} if for any element $T \subset \omega_z$ and $(d-1)$-face $F \subset \omega_z$ containing $z$ there is a sequence $\{T_i\}_{i=1}^m$ such that
\begin{itemize}
    \item any $T_i$ is an element of $\omega_z$ for $0\leq i \leq m$;
    \item any intersection $T_i \cap T_{i+1}$ is a $(d-1)-$face of $\omega_z$ for $i=0\leq i \leq m-1$;
    \item $T_0$ contains $F$ and $T_m=T$.
\end{itemize}
Notice that a star $\omega_z$ is $(d-1)-$faced connected if the set $\omega_z \cap \Omega$ is connected. 

\begin{proposition}[approximation of gradients]\label{P:cont-vs-discont}
 Let $v\in W^1_p(\Omega)$ for $1\le p\le d$. Let $\grid$ be a conforming
  mesh such that its stars are $(d-1)$-face-connected.
  Then there exists a constant $C(\sigma)$ depending on the
  shape regularity coefficient $\sigma$ of \eqref{E:shape-regularity}, the dimension $d$ and the polynomial
  degree $n\ge1$, such that
  \begin{equation}\label{eq:gradient-approx}
    1 \le \frac{\min_{w\in\mathbb {S}^{n,0}_\grid} \|\nabla (v-w)\|_{L^p(\Omega)}}
    {\min_{w\in\mathbb {S}^{n,-1}_\grid} \|\nabla (v-w)\|_{L^p(\Omega; \grid)}}
    \le C(\sigma),
  \end{equation}
  where $\|\nabla w\|_{L^p(\Omega; \grid)}$ stands for the broken norm over $\grid$.
\end{proposition}

The left inequality in \eqref{eq:gradient-approx} is obvious because $\mathbb {S}^{n,0}_\grid \subset \mathbb {S}^{n,-1}_\grid$. In contrast, the right inequality is delicate and relies on examining the quasi-interpolant \eqref{E:quasi-interpolant} \cite[Theorems 1 and 22]{Veeser:2016}. An important consequence of \eqref{eq:gradient-approx} is the following localized version of \eqref{E:global-interp}: 

\begin{proposition}[localized quasi-interpolation estimate]\label{P:local-quasi-interpolation}
Let $1< s \le n+1$, $1\le p\le d$ and
$r = \sob(W^s_p) -\sob(W^1_q) > 0$. If $v\in W^1_q(\Omega)$, then 
\begin{equation}\label{eq:localized-interp}
\|\nabla(v-I_\grid v)\|_{L^q(\Omega)} 
 \Cleq  \Big( \sum_{T\in\grid} h^{pr}_T | v|_{W^{s}_p(T)}^p \Big)^{\frac{1}{p}},
\end{equation}
\end{proposition}
\begin{proof}
Since $v-I_\grid v = (v-w) - I_\grid (v-w)$ for any $w\in\mathbb {S}^{n,0}_\grid$,
combine Proposition \ref{P:cont-vs-discont} with Proposition~\ref{P:Bramble-Hilbert-Sobolev}
(Bramble-Hilbert for Sobolev spaces) to write
\begin{equation*}
\|\nabla(v-I_\grid v)\|_{L^q(\Omega)} 
\Cleq \min_{w\in\mathbb {S}^{n,-1}_\grid} \|\nabla (v-w)\|_{L^q(\Omega)}
 \Cleq  \Big( \sum_{T\in\grid} h^{pr}_T | v|_{W^{s}_p(T)}^p \Big)^{\frac{1}{p}},
\end{equation*}
This concludes the proof.
\end{proof}
The crucial difference between \eqref{eq:localized-interp} and \eqref{E:global-interp}
is that the function $v\in W^1_q(\Omega)$ does not have to belong globally to $W^s_p(\Omega)$ but rather locally, namely
$v\in W^s_p(T)$ for every $T\in\grid$, to get optimal a priori error estimates. This property
will find several applications later. A special case of \eqref{eq:localized-interp} for
$p=q=2$ and $s=2$ reads
\[
|v-I_\grid v|_{H^1(\Omega)}^2 \Cleq \sum_{T \in \mesh} h_T^2 | v |_{H^2(T)}^2,
\]
for $v \in H^2(\Omega;\grid):=\{w \in H^1(\Omega): w_{|T} \in H^2(T) \ \forall T \in \grid\}$.

%--------------------------------------------------------------------------------
%\medskip\noindent
\paragraph{Quasi-Uniform Meshes.}
We now apply Theorem \ref{T:global-interp} to {\it quasi-uniform} meshes, 
namely meshes $\grid\in\grids$ for which all its elements are of comparable
size $h$, whence
\[
h \approx (\#\grid)^{-1/d} |\Omega|^{1/d} \approx (\#\grid)^{-1/d}.
\]
\vskip -.5cm
\begin{corollary}[quasi-uniform meshes]\label{nv-C:quasi-uniform}
 {\it Let $1 < s \le n+1$ and $v\in H^s(\Omega)$. If $\grid\in\grids$ is
quasi-uniform, then
\begin{equation}\label{E:quasi-uniform-est}
\|\nabla(v-I_\grid v)\|_{L^2(\Omega)} \Cleq |v|_{H^s(\Omega)} (\#\grid)^{-(s-1)/d}.
\end{equation}
}
\end{corollary}
\begin{remark}[optimal rate]\label{nv-R:optimal-rate}
If $s=n+1$, and so $v$ has the maximal regularity $v\in
H^{n+1}(\Omega)$,
then we obtain the optimal convergence rate in a linear Sobolev scale \looseness=-1
\begin{equation}\label{nv-optimal-rate}
\|\nabla(v-I_\grid v)\|_{L^2(\Omega)} \Cleq |v|_{H^{n+1}(\Omega)} (\#\grid)^{-n/d}.
\end{equation}
The order $-n/d$ is just dictated by the polynomial degree $n$ and cannot
be improved upon assuming either higher regularity that $H^{n+1}(\Omega)$ 
or a graded mesh $\grid$. The presence of $d$ in the exponent is referred to as 
{\it curse of dimensionality}.
\end{remark}

%------------------------------------------------------------------------------------
\medskip\noindent
{\it Example (corner singularity in $2d$)}.
To explore the effect of a geometric singularity on \eqref{E:quasi-uniform-est},
we let $\Omega= (-1,1)^2 \setminus [0,1]^2$ be a L-shaped domain and $v\in H^1(\Omega)$ be
\[
v(r,\theta) = r^{\frac23} \sin (2\theta/3) - r^{2}/4.
\]
This function $v\in H^1(\Omega)$ exhibits the typical corner
singularity of the solution of $-\Delta v = f$ with suitable Dirichlet
boundary condition: $v\in H^s(\Omega)$ for $s<5/3$.
Table \ref{nv-Ta:rate} displays the best
approximation error for polynomial degree $n=1,2,3$
and a sequence of {\it uniform} refinements in the seminorm 
$|\cdot|_{H^1(\Omega)}=\|\nabla\cdot\|_{L^2(\Omega)}$. This 
gives a {\it lower} bound for the interpolation error in \eqref{E:quasi-uniform-est}.

\begin{table}[htb!]
\begin{center}%
%\begin{minipage}{1.\textwidth}
\begin{tabular}{|c|c|c|c|}
\hhline{----}
\multirow{2}{*}{$h$} & linear  & quadratic & cubic \\
         &  $(n=1)$       & $(n=2)$ & $(n=3)$\\
\hhline{----}
1/4   & 1.14 &  9.64  &  9.89     \\
1/8   & 0.74 &  0.67  &  0.67     \\
1/16  & 0.68 &  0.67  &  0.67     \\
1/32  & 0.66 &  0.67  &  0.67     \\
1/64  & 0.66 &  0.67  &  0.67     \\
1/128 & 0.66 &  0.67  &  0.67     \\
\hhline{----}
\end{tabular}
%\end{minipage}
%
\caption{Rate of convergence $s$ in term of uniform mesh-size $h$. We observe an asymptotic error decay of about $h^{2/3}$ (i.e. $s=2/3$), or equivalently $(\#\grid)^{-1/3}$,  irrespective
of the polynomial degree $n$. This provides a lower bound for 
$\|v-I_\grid v\|_{L^2(\Omega)}$ and thus shows that
\eqref{E:quasi-uniform-est} is sharp.}\label{nv-Ta:rate}
\end{center}
\end{table}
%

\begin{comment}
    \begin{figure}[ht]
\label{nv-F:uniform}
\begin{center}
\hfil\includegraphics[width=.24\textwidth,trim={6cm 6cm 6cm 13cm},clip]%
{figs/rates/mesh1}
\hfil\includegraphics[width=.24\textwidth,trim={6cm 6cm 6cm 13cm},clip]%
{figs/rates/mesh2}
\hfil\includegraphics[width=.24\textwidth,trim={6cm 6cm 6cm 13cm},clip]%
{figs/rates//mesh3}
\hfil\includegraphics[width=.24\textwidth,trim={6cm 6cm 6cm 13cm},clip]%
{figs/rates/mesh4}
\hfil
\end{center}
\caption{Sequence of consecutive uniform meshes for L-shaped domain
  $\Omega$ created by 2 bisections.}
\end{figure}
\end{comment}
%
%\todo[inline]{RHN: There is a problem with the bounding box. This picture will most likely be replaced anyway.}
%

%\todo[inline]{RHN: Should we consider Kellogg's example also to illustrate that uniform refinement goes nowhere?}

Even though $s$ is fractional,
the error estimate \eqref{E:quasi-uniform-est} is still valid
as stated in Remark \ref{nv-R:fractional}. In fact, for uniform refinement,
\eqref{E:quasi-uniform-est} can
be derived by space interpolation between $H^1(\Omega)$ and
$H^{n+1}(\Omega)$. The asymptotic rate $(\#\grid)^{-1/3}$ reported
in Table \ref{nv-Ta:rate} is consistent with
\eqref{E:quasi-uniform-est} and independent of the polynomial degree
$n$; this shows that \eqref{E:quasi-uniform-est} is sharp.
It is also suboptimal as compared with the optimal rate $(\#\grid)^{-n/2}$
of Remark \ref{nv-R:optimal-rate}.

The question arises whether the rate $(\#\grid)^{-1/3}\approx h^{2/3}$ in Table \ref{nv-Ta:rate}
is just a consequence of uniform refinement or
unavoidable.  It is important to realize that
$v\not\in H^s(\Omega)$ for $s\geq5/3$ and thus
\eqref{E:quasi-uniform-est} is not applicable.  However, the problem is not
that second-order derivatives of $v$ do not exist but rather that they
are not square-integrable.  In particular, it is true that $v\in
W^2_p(\Omega)$ if $1\le p < 3/2$.  We therefore may apply Theorem
\ref{T:global-interp} with, e.g., $n=1$, $s=2$, and $p\in[1,3/2)$ and
then ask whether the structure of \eqref{E:global-interp} can be
exploited, e.g., by compensating the local singular behavior of $v$ with
the local mesh-size $h$.  This
enterprise naturally leads to \emph{graded} meshes adapted to $v$.

%\todo[inline]{This new discussion does not use the meshsize. The minimization process is now discrete.}

%-------------------------------------------------------------------------------------
\subsection{Principle of error equidistribution.}\label{S:error-equidistribution}
%-------------------------------------------------------------------------------------
%
We investigate the relation between local interpolation error and regularity for
the design of optimal graded meshes adapted to a given function $v\in H^1(\Omega)$ for $d=2$.
We recall that $W^2_1(\Omega)$ is in the same nonlinear Sobolev scale of
$H^1(\Omega)$, namely $\sob(W^2_1)=\sob(H^1)$, but $W^2_1(\Omega) \subset C^0(\overline{\Omega})$ \cite[Lemma 4.3.4]{BrennerScott:08},
and the {\it Lagrange interpolant} $I_\grid v$ is well defined and satisfies
\begin{equation}\label{E:local-error}
\|\nabla (v- I_\grid v)\|_{L^2(T)} \le C |v|_{W^2_1(T)} =: e_\grid(v,T) \quad\forall \, T\in\grid.
\end{equation}
We formulate a discrete minimization problem on the surrogate quantity
$\be := (e_\grid(v,T))_{T\in\grid}\in\R^N$
with $N=\#\grid$: minimize the square of the total $H^1$-error $E_\grid(v)$ \looseness=-1
\[
E_\grid(v)^2 := \sum_{T\in\grid} e_\grid(v,T)^2
\]
subject to the constraint
\[
\sum_{T\in\grid} e_\grid(v,T) = C |v|_{W^2_1(\Omega)}.
\]
We idealize the problem upon allowing $e_\grid(v,T)$ to attain any nonnegative real value despite the
fact that shape regularity of $\grid$ entails geometric restrictions between adjacent elements.
We next form the Lagrangian
\[
{\mathcal L}[\be,\lambda] := \sum_{T\in\grid} e_\grid(v,T)^2 - \lambda \left(\sum_{T\in\grid} e_\grid(v,T)
- C |v|_{W^2_1(\Omega)} \right),
\]
with Lagrange multiplier $\lambda\in\R$. We thus realize that the optimality condition reads
\[
e_\grid(v,T) = \frac{\lambda}{2} \quad\forall \, T\in\grid
\]
or that $e_\grid(v,T)$ is constant over $\grid$. We rewrite this insightful conclusion as follows: \looseness=-1
\begin{equation}\label{E:equidistribution}
\begin{minipage}{0.84\linewidth}
{\em A graded mesh is quasi-optimal if the local error is equidistributed.}
\end{minipage}
\end{equation}
This calculation yields
\[
E_\grid(v)^2 = \frac{\lambda^2}{4} N,\quad
C |v|_{W^2_1(\Omega)} = \frac{\lambda}{2} N, 
\]
whence
\begin{equation}\label{E:optimal-rate}
E_\grid(v) = C |v|_{W^2_1(\Omega)} N^{-1/2}
\end{equation}
is the optimal decay rate but with regularity $v\in W^2_1(\Omega)$ rather than $H^2(\Omega)$;
this is the second instance of nonlinear approximation, namely mesh design tailored to the specific
function $v$ at hand; the first one was in Section \ref{S:intro}.
The principle of error equidistribution \eqref{E:equidistribution} was originally derived by
Babu\v{s}ka and Rheinboldt \cite{BabuskaRheinboldt:78} for $d=1$, and extended to $d=2$ by
Nochetto and Veeser \cite[Section 1.6]{NochettoVeeser:2012}, using an idealized continuous minimization
problem involving a meshsize function. The current formulation is closer to applications and does not
require a positive power of $h_T$ in \eqref{E:local-error}.

\begin{remark}[point singularities]\label{R:point-singularities}
Corner singularities \cite{Grisvard:85} as well as singularities due to
intersecting interfaces \cite{Kellogg:75} are of the form
\begin{equation}\label{E:corner-sing}
v(x) \approx r(x)^\gamma,
\quad
0<\gamma<1,
\end{equation}
for $d=2$. This implies $v\in W^2_1(\Omega)$ for all $\gamma$ and the decay rate \eqref{E:optimal-rate}
provided $\grid$ equidistributes the $H^1$-error. Babu\v ska et al \cite{BaKePi:79}
and Grisvard \cite{Grisvard:85} designed such meshes for corner singularities using weighted $H^2$-regularity.
The preceding approach is more powerful in that it does not require any characterization
of the singularities, rather than $v\in W^2_1(\Omega)$, and applies as well to line discontinuities
for $d=2$. We will come back to this point in Section \ref{S:conv-rates-coercive}.
\end{remark}

We consider now an important abstract variant of the discrete minimization process leading to
\eqref{E:equidistribution}
that will be instrumental to understand the success of greedy algorithms later.
Suppose that $0<q,p\le\infty$, $v \in L^q(\Omega)$, $X^t_p(\Omega)$ is an abstract regularity space
with $t=\frac{1}{p}-\frac{1}{q}>0$, and $E_\grid(v)_q$ and $e_\grid(v,T)_q$ are global and local
$L^q$-interpolation error indicators of $v$ that satisfy the following two abstract properties:
\begin{itemize}
\item
{\it Summability in $\ell^q$}: There exists a constant $C_1>0$ such that
\begin{equation}\label{E:lq-sum}
E_\grid(v)_q^q \le C_1^q \sum_{T\in\grid} e_\grid(v,T)_q^q;
\end{equation}

\item
{\it Summability in $\ell^p$}: There exists a constant $C_2>0$ such that
\begin{equation}\label{E:lp-sum}
\sum_{T\in\grid} e_\grid(v,T)_q^p = C_2^p |v|_{X^t_p(\Omega)}^p.
\end{equation}
\end{itemize}
We intend to find conditions on a mesh $\grid$ that
minimize the global $L^q$-error $E_\grid(v)_q$ of $v$ subject to the constraint \eqref{E:lp-sum}. We
again propose a Lagrangian
\[
{\mathcal L}[\be,\lambda] := \sum_{T\in\grid} e_\grid(v,T)_q^q -
\lambda \left(\sum_{T\in\grid} e_\grid(v,T)_q^p - C_2^p |v|_{X^t_p(\Omega)}^p \right).
\]
The optimality condition for $\be$ reads
\[
e_\grid(v,T)_q = \Big( \lambda \frac{p}{q} \Big)^{\frac{1}{q-p}} \quad\forall T\in\grid,
\]
which is a third instance of error equidistribution and is thus consistent with \eqref{E:equidistribution}. 
We now resort \eqref{E:lq-sum} and \eqref{E:lp-sum} to arrive at
\[
\sum_{T\in\grid} e_\grid(v,T)_q^q = N \Big( \lambda \frac{p}{q} \Big)^{\frac{q}{q-p}},
\quad
C_2^p |v|_{X^t_p(\Omega)}^p = N \Big( \lambda \frac{p}{q} \Big)^{\frac{p}{q-p}},
\]
whence
\begin{equation}\label{E:error-decay}
E_\grid(v)_q \le C_1 C_2 |v|_{X^t_p(\Omega)} N^{\frac{1}{q}-\frac{1}{p}}.
\end{equation}
We see that the decay rate in \eqref{E:error-decay} is $-t = \frac{1}{q}-\frac{1}{p} < 0$ and is
just dictated by the different summabilities of \eqref{E:lq-sum} and \eqref{E:lp-sum}. In the
applications of \eqref{E:error-decay} below, $t=\frac{s}{d}$ will be proportional to a
differentiability index $s$ and the condition
\[
0 = t-\frac{1}{p} + \frac{1}{q} = \frac{s}{d} -\frac{1}{p} + \frac{1}{q}
\]
will correspond to the spaces $L^q(\Omega)$ and $X^s_p(\Omega)$ being on the same nonlinear Sobolev scale.
This minimization process is an idealization that does not account for mesh regularity of $\grid$,
which in turn entail some geometric constraints in the construction of $\grid$.
A key question is whether estimates of the form \eqref{E:error-decay} can be achieved
under practical but weaker conditions than \eqref{E:equidistribution}.
In Section \ref{S:bisection} we will study the bisection method, a flexible technique for
conforming mesh refinement with optimal complexity. 
In Section \ref{S:constructive-approximation} we will present and analyze
\textsf{GREEDY}, a practical algorithm that implements these ideas and constructs quasi-optimal
conforming bisection meshes under the slightly stronger assumption
\begin{equation}\label{E:nonlinear-Sobolev-scale}
s - \frac{d}{p} + \frac{d}{q} > 0.
\end{equation}
Moreover, in Section \ref{S:nonconforming-meshes} we will extend this analysis to non-conforming meshes.

We realize from \eqref{E:error-decay} that, in order to maximize the error decay rate,
we would like to have $p$ as small as possible, even $0<p<1$. The range of $q,p$ does not matter
in the argument above and, despite the fact that $q\ge 1$ in all applications below, the range of $p$
is only limited by that of $s$, which in turn depends on the polynomial degree $n\ge1$ in that
$0 < s \le n+1$.

We now return to the special case \eqref{E:local-error}, namely $q=2, p=1$ and $\nabla v \in L^2(\Omega)$.
As already shown in \eqref{E:local-error}, in the nonlinear Sobolev scale
\[
\sob(W^2_1) - \sob(H^1) = \Big( 2 - \frac{2}{1} \Big) - \Big(1 - \frac{2}{2} \Big)
= 0
\]
we expect the best error decay
\[
\|\nabla(v-I_\grid v)\|_{L^2(\Omega)} \Cleq |v|_{W^2_1(\Omega)} (\#\grid)^{-\frac{1}{2}},
\]
whereas the linear Sobolev scale yields the reduced order
\[
\|\nabla(v-I_\grid v)\|_{L^2(\Omega)} \Cleq |v|_{H^s(\Omega)} (\#\grid)^{-(s-1)/2}
\]
for $s < 1+\gamma < 2$ and $v$ satisfying \eqref{E:corner-sing}, where $I_\grid$ is the Lagrange
interpolation operator. The nonlinear Sobolev scale entails
a trade of differentiability with integrability:
we gain up to differentiability $s=2$ at the expense of lower integrability $p=1$ for
polynomial degree $n=1$.  This
trade-off is at the heart of the optimal estimate \eqref{E:optimal-rate} and is represented
in the so-called {\it DeVore diagram} in Fig. \ref{F:devore-diagram}.

If the polynomial degree is $n\ge2$, then the largest differentiability index is
$s=n+1$, which for $d=2$ leads to integrability index $p<1$:
\begin{equation}\label{E:p<1}
\Big(s-\frac{2}{p}\Big) - \Big(1 - \frac{2}{2}\Big) = 0 \quad\Rightarrow\quad  p = \frac{2}{n+1} < 1.
\end{equation}
To measure regularity of $v$, the corresponding Sobolev space
must be replaced by the Besov space $B^{n+1}_{p,p}(\Omega)$
or the Lipschitz space ${\it Lip}^{n+1}_p(\Omega)$. We will introduce and study these spaces in 
Section \ref{S:Besov}.

%-------------------------------------------------------------------------------------
\subsection{Conforming meshes: the bisection method}\label{S:bisection}
%-------------------------------------------------------------------------------------
%\todo[inline]{AB: we need to add a statement about conforming meshes and generations (cite N-S-V)}

In order to approximate functions in $W^k_p(\Omega)$ by piecewise
polynomials, we decompose $\Omega$ into simplices. We briefly
discuss the {\em bisection} method, an elegant and versatile technique
for subdividing $\Omega$ in any dimension into a conforming mesh. 
We also discuss briefly 
nonconforming meshes in \S \ref{S:nonconforming-meshes}. 
We present complete proofs, especially of the complexity of bisection, later in 
\S \ref{S:mesh-refinement}. 

We focus on $d=2$ and follow \cite{BiDaDeV:04}, but the results carry over to any dimension $d>2$ 
\cite{Stevenson:08}. We refer to \cite{NoSiVe:09} for a rather complete discussion for
$d\ge2$. 

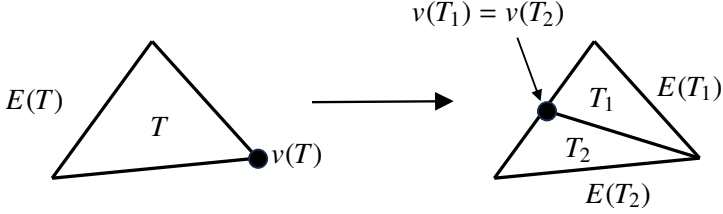
\begin{figure}[ht]
	\begin{center}
	    \includegraphics[scale=0.7,trim={5cm 12cm 15cm 2.5cm},clip]{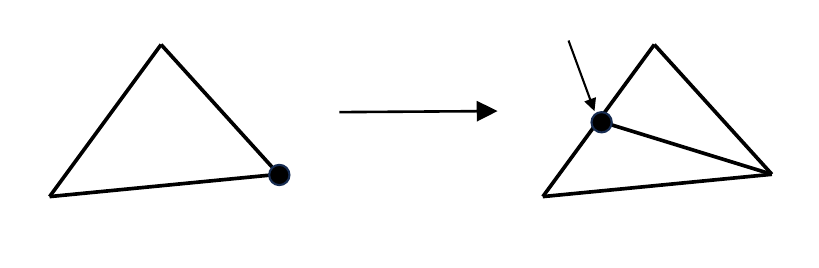}
    \begin{picture}(0,0)
        \put(-280,50){$E(T)$}
        \put(-225,40){$T$}
        \put(-180,30){$v(T)$}
        \put(-69,32){$T_2$}
        \put(-60,51){$T_1$}
        \put(-62,15){$E(T_2)$}
        \put(-35,55){$E(T_1)$}
        \put(-127,83){$v(T_1)=v(T_2)$}
    \end{picture}
 
        \caption{Triangle $T\in\grid$ with vertex $v(T)$ and opposite
          refinement edge $E(T)$. The bisection rule for $d=2$ consists of
          connecting $v(T)$ with the midpoint of $E(T)$, thereby giving rise to
          children $T_1,T_2$ with common vertex $v(T_1)=v(T_2)$, the
          newly created vertex, and 
          opposite refinement edges $E(T_1), E(T_2)$.\label{nv-F:edge-vertex}}
    \end{center}
\end{figure}

Let $\grid$ denote a {\it mesh} (triangulation or grid) made of
simplices $T$, and let $\grid$ be {\it conforming} (edge-to-edge). Each
element is labeled, namely it has an edge $E(T)$ assigned for
refinement (and an opposite vertex $v(T)$ for $d=2$); 
see Figure \ref{nv-F:edge-vertex}.

\begin{figure}[t!]
  \centering
  \includegraphics[width=0.8\hsize]{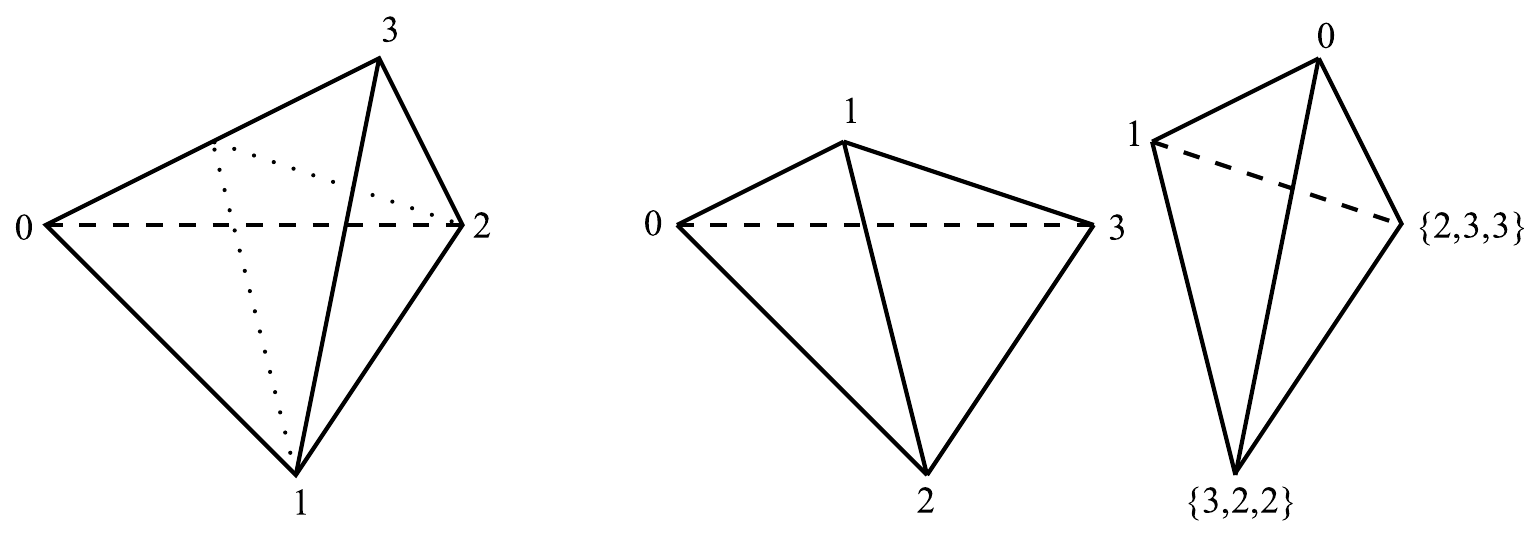}
  \caption{Refinement of a single tetrahedron $\elm$ of type $t$. The child
    $\elm_1$ in the middle has the same node ordering regardless of
    type. In contrast, for the child $\elm_2$ on the right
    a triple is appended to two nodes. The local vertex index is given
    for these nodes by the $t$-th component of the triple.
    \label{nsv-F:bisect-3d}}
\end{figure}

The bisection method consists of a suitable {\it labeling} of the initial mesh
$\grid_0$ and a rule to assign the refinement edge to the two
children. For $d=2$ we consider the {\it newest vertex bisection}
as depicted in Figure \ref{nv-F:edge-vertex}.
For $d>2$ the situation is more complicated and one needs the
concepts of type and vertex order \cite{NoSiVe:09, Stevenson:08}.
More precisely, we identify a simplex $\elm$ with the set of
its \emph{ordered vertices} and its \emph{type} $t$ by
\begin{displaymath}
  \elm = \{z_0,z_1,\dots,z_d\}_{t} \,,
\end{displaymath}
with $ t \in \{0,\dots,d-1\}$. Given such a $d$-simplex $\elm$ we use the following bisection rule to
split it in a unique fashion and to impose both vertex order and type
to its children. The edge $\overline{z_0z_d}$ connecting the first and
last vertex of $\elm$ is the \emph{refinement edge} of $\elm$ and its
midpoint $\bar{z}=\frac{z_0+z_d}2$ becomes the
new vertex.  Connecting the new vertex $\bar{z}$ with the vertices of
$\elm$ other than $z_0,z_d$ determines the common face $S = \{\bar z,
z_1,\dots,z_{d-1}\}$ shared by the two children $\elm_1,\elm_2$ of
$\elm$. The \emph{bisection rule} dictates the following vertex order
and type for $\elm_1,\elm_2$
\begin{equation}\label{nsv-bisect}
  \begin{aligned}
    \elm_1&\definedas
    \{z_0,\bar{z},\underbrace{z_1,\dots,z_t}_{\rightarrow},
    \underbrace{z_{t+1},\dots,z_{d-1}}_{\rightarrow}\}_{(t+1)\!\!\!\!\mod \!d},\\[-2mm]
    \elm_2&\definedas
    \{z_d,\bar{z},\underbrace{z_1,\dots,z_t}_{\rightarrow}
    ,\underbrace{z_{d-1},\dots,z_{t+1}}_{\leftarrow}\}_{(t+1)\!\!\!\!\mod \!d},
  \end{aligned}
\end{equation}
with the convention that arrows point in the direction of increasing
indices and $\{z_1,\dots,z_0\}=\emptyset$,
$\{z_{d},\dots,z_{d-1}\}=\emptyset$. 
%
%In 2d the bisection rule does not depend on the element type
%and we get for
%\begin{math}
%  \elm = \{z_0,z_1,z_2\}
%\end{math}
%the two children
%\begin{displaymath}
%  \elm_1 = \{z_0,\bar{z},z_1\}\qquad\text{and}\qquad
%  \elm_2 = \{z_2,\bar{z},z_1\}.
%\end{displaymath}
%\begin{figure}[htbp]
%  \centering
%  \includegraphics[width=0.47\hsize]{bisect_2d-1}
%  \hfill
%  \includegraphics[width=0.47\hsize]{bisect_2d-2}
%
%  \caption{Refinement of a single triangle $T=\{z_0,z_1,z_2\}$
%    and its reflected triangle $T_R =\{z_2,z_1,z_0\}$.
%    \label{nsv-F:bisect-2d}}
%\end{figure}
%%
%As depicted in Fig.~\ref{nsv-F:bisect-2d}, the refinement edge of the
%two children is opposite to the new vertex $\bar{z}$, whence this
%procedure coincides with \emph{the newest vertex bisection} for
%$d=2$. For $d\geq 3$ the bisection of an element does depend on its
%type, and, as we shall see below, this is important for preserving
%shape regularity. 
For instance, in 3d the children of
$T=\{z_0,z_1,z_2,z_3\}_t$ are (see Fig.~\ref{nsv-F:bisect-3d})
\begin{displaymath}
  \begin{aligned}
    t&=0: &\qquad T_1&=\{z_0,\bar{z},z_1,z_2\}_1 &&\text{and}&
    T_2&=\{z_3,\bar{z},z_2,z_1\}_1,\\
    t&=1: & T_1&=\{z_0,\bar{z},z_1,z_2\}_2 &
    &\text{and}& T_2&=\{z_3,\bar{z},z_1,z_2\}_2,\\
    t&=2: & T_1&=\{z_0,\bar{z},z_1,z_2\}_0 &
    &\text{and}& T_2&=\{z_3,\bar{z},z_1,z_2\}_0.
  \end{aligned}
\end{displaymath}
Note that the vertex labeling of $\elm_1$ is type-independent, whereas
that of $\elm_2$ is the same for type $t=1$ and $t=2$. To account for this
fact the vertices $z_1$ and $z_2$ of $\elm$ are tagged $(3,2,2)$ and
$(2,3,3)$ in Fig.~\ref{nsv-F:bisect-3d}. The type of $\elm$ then
dictates which component of the triple is used to label the vertex.
%

%Any different labeling of an element's vertices does not change its
%geometric shape but applying the above bisection rule it does change the
%shape and vertex order of its two children. This holds true for any
%relabeling except one. An element with this special relabeling of vertices
%is called \emph{reflected element}. We state next its precise definition.
%
%%\begin{definition}[Reflected Element]\label{nsv-D:ref-elm}
%Given an element $T=\{z_0,\cdots,z_d\}_t$, the \emph{reflected element} is
%given by
%%
%\begin{displaymath}
%  \elm_R\definedas \{z_d,\underbrace{z_1,\dots,z_t}_{\rightarrow}
%  ,\underbrace{z_{d-1},\dots,z_{t+1}}_{\leftarrow},z_0\}_{t}.
%\end{displaymath}
%%
%%\end{definition}
%
%Fig.~\ref{nsv-F:bisect-2d} depicts for 2d $\elm = \{z_0,z_1,z_2\}$ and
%$\elm_R = \{z_2,z_1,z_0\}$. It shows that the children of $\elm$ and
%$\elm_R$ are the same. This property extends to $d\ge3$; compare with
%Problem~\ref{nsv-E:bisect-reflected}. Any other
%relabeling of vertices leads to different shapes of the children, in
%fact as many as $\frac12(d+1)!$

Bisection creates a \textit{unique} master
forest $\mathbb{F}$ of binary trees with infinite depth, where each
node is a simplex (triangle in 2d), its two successors are the two
children created by bisection, and the roots of the binary trees are
the elements of the initial conforming partition $\grid_0$. It is
important to realize that, no matter how an element arises in the
subdivision process, its associated newest vertex is unique and only
depends on the labeling of $\grid_0$: so the edge $E(T)$ assigned for 
refinement (and the opposite vertex $v(T)$ for $d=2$)  are
independent of the order of the subdivision process for all
$T\in\mathbb{F}$; see Lemma \ref{L:unique-labeling} in Sect. \ref{S:mesh-refinement}.
Therefore, $\mathbb{F}$ is unique.

A finite subset $\mathcal{F}\subset\mathbb{F}$ is called a
\textit{forest} if $\grid_0\subset\mathcal{F}$ and the nodes of
$\mathcal{F}$ satisfy
\begin{itemize}
	\item[$\bullet$] all nodes of $\mathcal{F}\setminus\grid_0$ have a predecessor;
	\item[$\bullet$] all nodes in $\mathcal{F}$ have either two successors or none.
\end{itemize}
Any node $T\in\mathcal{F}$ is thus uniquely connected with a node
$T_0$ of the initial triangulation $\grid_0$, i.e. $T$ belongs to the
infinite tree ${\mathbb{F}}(T_0)$ emanating from $T_0$.
Furthermore, any forest may have
\textit{interior nodes}, i.e. nodes with successors, as well as
\textit{leaf nodes}, i.e. nodes without successors. The set of leaves
corresponds to a mesh (or triangulation, grid, partition) 
$\grid=\grid(\mathcal{F})$ of $\grid_0$ which may not
be conforming or edge-to-edge. 

We thus introduce the set  $\grids$ of all {\it conforming} refinements of 
$\gridk[0]$:
\begin{displaymath}
\index{Meshes!$\grids$: set of all conforming refinements of $\grid_0$}
  \grids := \{\grid = \grid(\mathcal{F}) \mid \mathcal{F}\subset\mathbb{F} 
  \text{ is finite and }\grid(\forest) \text{ is conforming}\}.
\end{displaymath}
If $\gridk[*]=\grid(\forest_*)\in\grids$ is a conforming refinement of
$\grid=\grid(\forest)\in\grids$, we write $\gridk[*]\ge\grid$ and understand this
inequality in the sense of trees, namely $\forest\subset\forest_*$.

\begin{figure}[ht]
	\centering
	\includegraphics[scale=0.75]{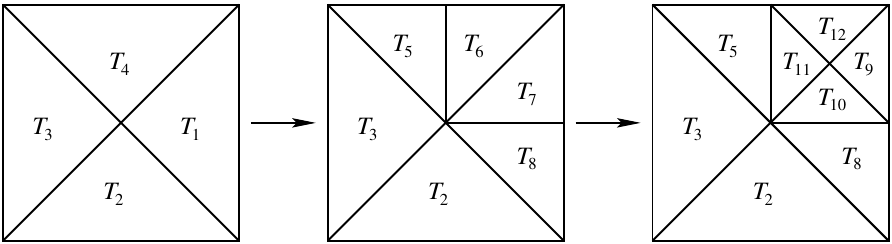}
        \caption{Sequence of bisection meshes $\{\gridk[k]\}_{k=0}^2$ starting
        from the initial mesh $\gridk[0]=\{T_i\}_{i=1}^4$ with longest
        edges labeled for bisection. Mesh
        $\gridk[1]$ is created from $\gridk[0]$ upon bisecting $T_1$
        and $T_4$, whereas mesh $\gridk[2]$ arises from $\gridk[1]$
        upon refining $T_6$ and $T_7$. The bisection rule is described
        in Figure \ref{nv-F:edge-vertex}.\label{nv-F:sequence}}
\end{figure}
\begin{figure}[ht]
	\centering
	\includegraphics[scale=0.65]{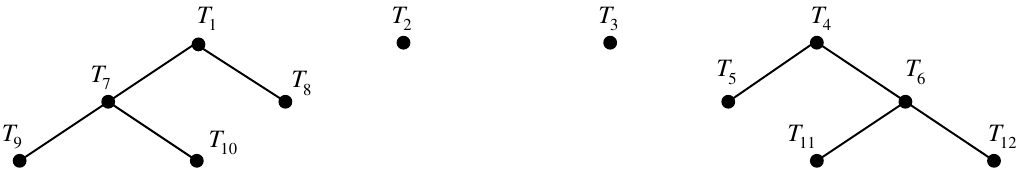}
        \caption{Forest $\forest_2$ corresponding to the grid sequence
          $\{\gridk[k]\}_{k=0}^2$ of Figure \ref{nv-F:sequence}. The
          roots of $\forest_2$ form the initial mesh $\gridk[0]$
          and the leaves of $\forest_2$ constitute the conforming
          bisection mesh $\gridk[2]$. Moreover, each level of
          $\forest_2$ corresponds to all elements with generation 
          equal to the level.\label{nv-F:tree}}
\end{figure}
\medskip\noindent
{\it Example.}
Consider $\grid_0=\{T_i\}_{i=1}^4$ and the longest edge
to be the refinement edge. Figure \ref{nv-F:sequence}
displays a sequence of conforming meshes $\gridk[k]\in\grids$
created by bisection.

\noindent
Each element $T_i$ of $\grid_0$ is a root of a finite tree
emanating from $T_i$, which together form the forest $\forest_2$
corresponding to mesh $\gridk[2]=\grid(\forest_2)$.
Figure \ref{nv-F:tree} displays $\forest_2$, whose leaf nodes are
the elements of $\gridk[2]$. 

\medskip\noindent
%----------------------------------------------------------------------------
{\it Properties of Bisection.}
We now discuss several crucial geometric properties of bisection.
We start by recalling the concept of shape regularity. For any $T\in\grid$,
we define

\begin{minipage}[b]{0.4\linewidth}
\includegraphics[scale=0.45,trim={4cm 9cm 19cm 4cm},clip]{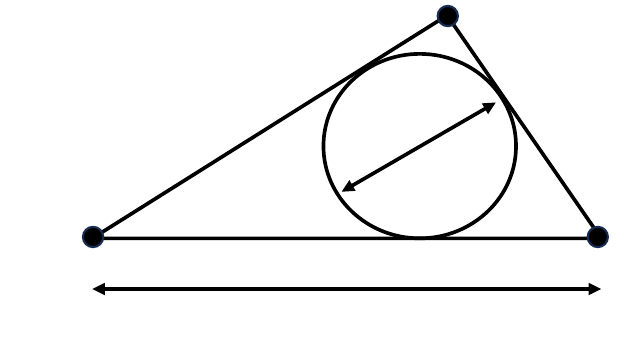}
\begin{picture}(0,0)
\put(-50,35){$\underline{h}_T$}
\put(-70,0){$\overline{h}_T$}
\end{picture}
\end{minipage}
\begin{minipage}[b]{0.6\linewidth}
\begin{align*}
	&\overline{h}_T:=\diam(T)\\
	&h_T:=|T|^{1/d}\\
	&\underline{h}_T:=2\sup\{r>0\,|\, B(x,r)\subset T \text{ for }x\in T\}.
\end{align*}
\vspace{0.35cm}
\end{minipage}

Then
\[
\underline{h}_T \le h_T \le \overline{h}_T\le \sigma \underline{h}_T
\qquad\forall T\in\grid,
\]
where $\sigma>1$ is the shape regularity constant of \eqref{E:shape-regularity}. 
The next lemma guarantees that bisection keeps $\sigma$ bounded.

\begin{figure}[htbp]
  \centering
  \includegraphics[width=0.9\hsize]{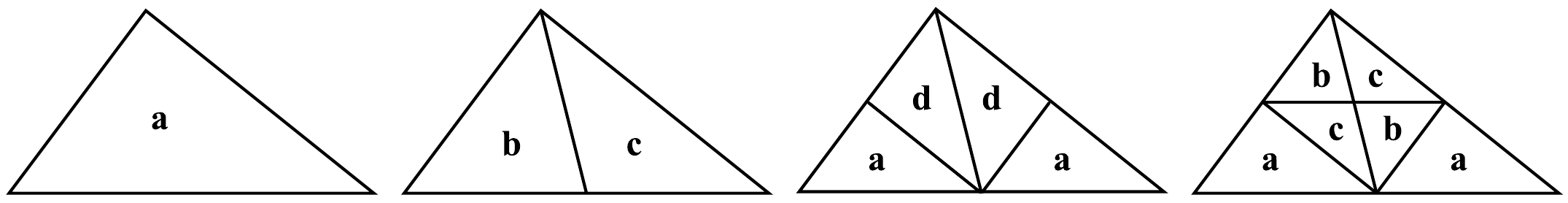}
  \caption{Bisection produces at most 4 similarity classes for any
    triangle.\label{nv-F:2dsim-gen}}
\end{figure}
\begin{lemma}[shape regularity]\label{nv-L:shape-reg}
\it The partitions $\grid$ generated by newest vertex bisection satisfy a
uniform minimal angle condition, or equivalently $\sigma$ is uniformly
bounded, only depending on the initial partition $\gridk[0]$.
\end{lemma}
\begin{proof}
Each $T\in\gridk[0]$ gives rise to a fixed number of similarity
classes, namely $4$ for $d=2$ according to Figure \ref{nv-F:2dsim-gen}.
This, combined with the fact that $\#\gridk[0]$ is finite, yields
the assertion.
\end{proof}

We define the {\it generation (or level)} $g(T)$ \index{Meshes!$g(T)$: generation of $T$} of an 
element $T\in\grid$ as the number of bisections needed to create $T$
from its ancestor $T_0\in\gridk[0]$. Since bisection splits an element
into two children with equal measure, we realize that
\begin{equation}\label{nv-mesh-gen}
h_T = 2^{-g(T)/2} h_{T_0} \qquad\forall T\in\grid.
\end{equation}
Referring to Figure \ref{nv-F:tree} we observe that
the leaf nodes $T_9,T_{10},T_{11},T_{12}$ have generation $2$, whereas
$T_5,T_8$ have generation $1$ and $T_2,T_3$ have generation $0$.

The following geometric property is a simple consequence of
\eqref{nv-mesh-gen}.
\begin{lemma}[element size vs generation]\label{nv-L:diam-gen}
\it There exist constants $0<D_1<D_2$, only depending on $\gridk[0]$, such
that
\begin{equation}\label{nv-diam-gen}
D_1  2^{-g(T)/2} \le h_T < \overline{h}_T \le D_2  2^{-g(T)/2}
\qquad\forall T\in\grid.
\end{equation}
\end{lemma}

\begin{figure}[h]
	\centering
	\includegraphics[scale=0.75]{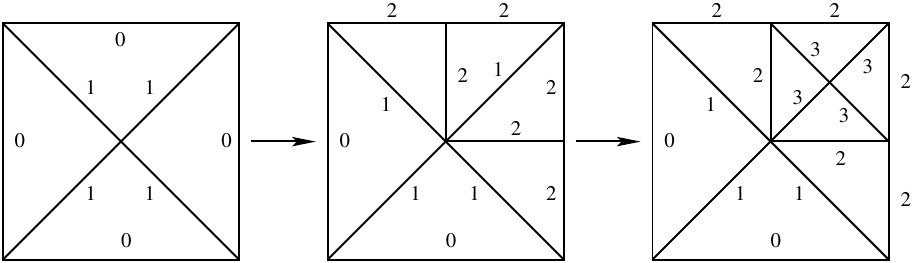}
        \caption{Initial labeling and its evolution for the sequence
          of conforming refinements
          $\gridk[0]\le\gridk[1]\le\gridk[2]$
          of Figure \ref{nv-F:sequence}. \label{nv-F:labeling}}
\end{figure}

%-----------------------------------------------------------------------------------
\medskip\noindent
{\it Labeling and Bisection Rule.}
Whether the recursive application of bisection does not lead to
inconsistencies depends on a suitable initial labeling of edges and a 
bisection rule. For $d=2$ they are simple to state  \cite{BiDaDeV:04}.
%but for $d>2$ we refer to Condition (b) of Section 4
%of \cite{Stevenson:08}.
Given $T\in\grid$ with generation $g(T)=i$, we assign the label $(i+1,i+1,i)$ to $T$ 
with $i$ corresponding to the refinement edge $E(T)$. 
The following rule dictates how the labeling changes with refinement: 
the side $i$ is bisected and both new sides as well as the bisector 
are labeled $i+2$ whereas the remaining labels do not change.
To guarantee that the label of an edge is independent of the elements
sharing this edge, we need a special labeling for $\grid_0$, due
\cite[Theorem 2.9]{Mitchell:89} and \cite[Lemma 2.1]{BiDaDeV:04} for $d=2$:
\begin{equation}\label{nv-initial-label}
\begin{minipage}{0.75\linewidth}
\emph{edges of $\grid_0$ have labels $0$ or $1$ and all elements
  $T\in\grid$ have exactly two edges with label $1$ and one with label
  $0$}.
\end{minipage}
\end{equation}
There is a variant for $d>2$ due to \cite[Section 4]{Stevenson:08}.
It is not obvious that labeling \eqref{nv-initial-label} exists, but if it does then all
elements of $\gridk[0]$ can be split into pairs of compatibly
divisible elements.
We refer to Figure \ref{nv-F:labeling} for an example of initial
labeling of $\gridk[0]$ satisfying \eqref{nv-initial-label}
and the way it evolves for two successive
refinements $\gridk[2]\ge\gridk[1]\ge\gridk[0]$ corresponding to
Figure \ref{nv-F:sequence}. 
%
%------------------------------------------------------------------------------

\begin{figure}[ht]
	\centering
	\includegraphics[scale=0.8]{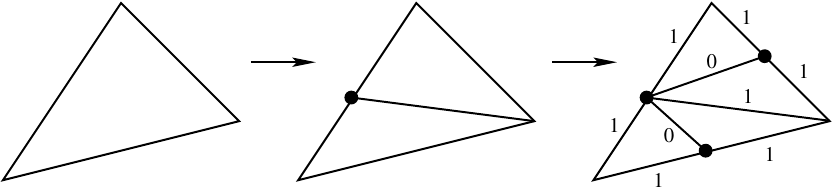}
        \caption{Bisecting each triangle of $\gridk[0]$ twice and
          labeling edges in such a way that all boundary edges have
          label 1 yields an initial mesh satisfying
          \eqref{nv-initial-label}. \label{nv-F:initial-label-2d}}
\end{figure}

To guarantee \eqref{nv-initial-label} we can proceed as follows:
given a coarse mesh of elements $T$ we can bisect twice each $T$ and
label the four grandchildren, as indicated in Figure
\ref{nv-F:initial-label-2d} 
for the resulting mesh $\grid_0$ to satisfy the initial labeling
\cite{BiDaDeV:04}.

For $d \geq 3$  a general strategy of initial labeling is due to 
\cite[Section 4 - Condition (b)]{Stevenson:08}, who in turn improves
upon \cite{Maubach:95} and \cite{Traxler:97} and shows
how to impose it upon further refining each element of $\grid_0$. 
We refer to the survey \cite{NoSiVe:09} for a discussion of this condition:
a key consequence is \looseness=-1
\begin{equation}\label{E:uniform-refinement}
\begin{minipage}{0.8\linewidth}
\emph{every uniform refinement of $\grid_0$ gives
a conforming bisection mesh}.
\end{minipage}
\end{equation}
Condition \eqref{nv-initial-label} is still valid, and a construction by successive bisections similar, 
but much trickier, than the one described for $d=2$ can be performed to fulfil it; yet 
for $d=3$ the number of elements increases by an order of magnitude, which indicates
that \eqref{nv-initial-label} is a severe restriction in practice.
Finding alternative, more practical, conditions is an important problem.
%
% \todo[inline]{Is there now an alternative? For quadrilateral meshes with handing nodes
%   this is not an issue, and these meshes are easier to handle for $d>2$ Should
%   we comment on this here?
% }

\medskip\noindent
{\it Initialization of arbitrary triangulations.} A novel initialization procedure that can be applied to {\em any} conforming triangulation $\gridk[0]$ has been recently proposed by \cite{Diening2023}; hereafter we present a short account of it.

The key concept is that of {\em coloring} the vertices of $\gridk[0]$. A colored initial triangulation in $\mathbb{R}^d$ is a pair $(\gridk[0],c)$, where $c: {\cal V}_{\gridk[0]} \to \{0, \dots, d \}$ is such that the colors of all vertices of each $T \in \gridk[0]$ are distinct. The color map $c$ allows one to sort the vertices of each initial element $T=\{ z_0, \dots, z_d\}_t \in \gridk[0]$ so that 
$$
c(z_j) = j \qquad j \in \{0, \dots, d \}\,.
$$
To refine a marked $T \in \gridk[0]$, one applies the Maubach bisection rule leading to \eqref{nsv-bisect}, and possibly adds a recursive closure, which is proven to terminate, to guarantee the conformity of the final triangulation. The coloring property is conserved in this process, and the conclusion of Theorem \ref{nv-T:complexity-refine} below holds true, starting from any initially colored triangulation $\gridk[0]$.

Unfortunately, not every initial triangulation can be colored. For instance, consider in dimension $d=2$ a patch of triangles sharing a common vertex. If color $0$ is assigned to such inner vertex, then the outer vertices must take successively colors 1 and 2, but if the number of triangles in the patch is odd there will be a vertex which is not colorable. 

To overcome this obstruction, \cite{Diening2023} propose to use more colors, and introduce the concept of {\em generalized coloring}: a pair $(\gridk[0],c)$ is a $(N+1)$-colored triangulation if there exists an integer $N\geq d$ and a mapping $c: {\cal V}_{\gridk[0]} \to \{0, \dots, N \}$ such that the colors of all vertices of each $T \in \gridk[0]$ are distinct. Any initial $\gridk[0]$ can be colored in this generalized sense: indeed, after the initialization $c(z)=+\infty$ for all $z \in {\cal V}_{\gridk[0]}$, one defines
$$
c(z):=\min \, ( \mathbb{N}_0 \setminus \{c(w):[z,w] \text{ is an edge of }\gridk[0] \} ), \qquad z \in {\cal V}_{\gridk[0]},
$$
as the smallest color not already attained by a neighboring vertex. Then, one sets $N:=\max \, \{c(z) : z \in {\cal V}_{\gridk[0]} \}$, and notes that $N$ is bounded by the maximal number of edges connected to a vertex of $\gridk[0]$.

A generalized $(N+1)$-colored triangulation $(\gridk[0],c)$ in $\mathbb{R}^d$ can be seen as a collection of $d$-simplices contained in a virtual, colored triangulation $\gridk[0]^+$ in $\mathbb{R}^N$. It suffices to add $N-d$ virtual nodes to each simplex in $\gridk[0]$, so that it becomes a $N$-simplex, and attribute to these nodes the remaining $N-d$ colors. Note that these virtual simplices are only connected via their $d$-subsimplices belonging to $\gridk[0]$. 
In the example mentioned above of a patch of triangles sharing a vertex, a $(3+1)$-colored triangulation is defined as follows: a tetrahedron is built on top of each triangle; the previously uncolorable vertex takes the new color 3, whereas color 2 is attributed to the new vertices of the two tetrahedra sharing that vertex; the new vertex of any other tetrahedron takes the color 3.

With the new triangulation $\grid_0^+$ at hand, one could apply the Mauback bisection rule to it, which as a by-product would refine the initial triangulation $\grid_0$. However, \cite{Diening2023} suggests a short-cut that refines directly $\grid_0$, by invoking an algorithm that bisects a $k$-simplex in dimension $m>k$. A further round of recursive refinements may be needed to guarantee conformity. Diening, Gehring, and Storn  prove that the recursion terminates. In addition, for any $(N+1)$-colored initial triangulation, the conclusion of Theorem \ref{nv-T:complexity-refine} below holds true also in this case, with a constant $\Ccompl$ satisfying $\Ccompl \lesssim N^d$.

%----------------------------------------------------------------------------------
\medskip\noindent
{\it The procedure} \REFINE.
Given $\grid\in\grids$ and a selected subset $\marked\subset\grid$ (the set of {\it marked}
elements), the procedure
\[
[\gridk[*]] = \REFINE (\grid,\marked)
\]
creates a new conforming refinement $\gridk[*]$ of $\grid$ by
bisecting all elements of $\marked$ at least once and perhaps
additional elements to keep conformity.

Conformity is a constraint in the refinement procedure that prevents
it from being completely local. The propagation of refinement beyond
the set of marked elements $\marked$ is a rather delicate matter,
which we discuss later in Sect. \ref{S:mesh-refinement}. For
instance, we show that a naive estimate of the form
\[
\#\gridk[*] - \#\grid \le \Ccompl ~\#\marked
\]
is {\it not} valid with an absolute constant $\Ccompl$ independent of the
refinement level. This can be repaired upon considering the cumulative
effect for a sequence of conforming bisection meshes
$\{\gridk\}_{k=0}^\infty$.
This is expressed in the following crucial complexity result due to
\cite{BiDaDeV:04} for $d=2$ and \cite{Stevenson:08} for $d>2$. We present a complete
proof later in Sect. \ref{S:mesh-refinement}.

\begin{theorem}[complexity of $\REFINE$]\label{nv-T:complexity-refine}
{\it If $\gridk[0]$ satisfies the initial labeling
\eqref{nv-initial-label} for $d=2$, or that in \cite[Section 4]{Stevenson:08}
for $d>2$,
then there exists a constant $\Ccompl>0$ only depending on
$\gridk[0]$ and $d$ such that for all $k\ge 1$
}
\index{Constants!$\Ccompl$: complexity of \REFINE constant}
\[
\#\gridk - \#\gridk[0] \le \Ccompl \sum_{j=0}^{k-1} \# \markedk[j].
\]
\end{theorem}
If elements $T\in\marked$ are to be bisected $b\ge1$ times, then the
procedure \REFINE can be applied recursively, and 
Theorem \ref{nv-T:complexity-refine} remains valid with $\Ccompl$
also depending on $b$. \looseness=-1

%------------------------------------------------------------------------------------
\medskip\noindent
{\it Mesh Overlay.}
%------------------------------------------------------------------------------------
For the subsequent discussion it will be convenient to merge (or superpose) two conforming
meshes $\grid_1,\grid_2\in\grids$, thereby giving rise to the 
so-called \emph{overlay} $\grid_1\oplus\grid_2$\index{Meshes!$\grid_1\oplus\grid_2$: mesh overlay}. This operation
corresponds to the union in
the sense of trees \cite{CaKrNoSi:08,Stevenson:07}. 
We next bound the
cardinality of $\grid_1\oplus\grid_2$ in terms of that of $\grid_1$
and $\grid_2$.

\begin{lemma}[mesh overlay]\label{nv-L:overlay}
Let $\grid_1,\grid_2\in\grids$.
	The overlay $\grid=\grid_1\oplus\grid_2\in\grids$ is conforming and 
	\begin{equation}\label{nv-overlay}
		\#\grid\le\#\grid_1+\#\grid_2-\#\grid_0.
	\end{equation}
\end{lemma}
For a proof we refer to \cite[Lemma~3.7]{CaKrNoSi:08} and to Proposition~\ref{P:mesh-overlay-Lambda} 
below for a more general situation.

%---------------------------------------------------------------------------------
\subsection{Constructive approximation}\label{S:constructive-approximation}
%---------------------------------------------------------------------------------
We now construct graded bisection meshes $\grid$ for
$n=1,d=2$ that achieve the optimal
decay rate $(\#\grid)^{-1/2}$ of \eqref{E:optimal-rate}
under the global regularity assumption 
\begin{equation}\label{E:weaker-reg}
    v\in W^2_p(\Omega),\quad p>1.
\end{equation}
Therefore, $W^2_p(\Omega)$ is strictly above the Sobolev line for the space $H^1(\Omega)$:
$\sob(W^2_p)=2-\frac{2}{p}>0=\sob(H^1)$. Note that $s=1, p>1$ and $q=2$
obey the restriction \eqref{E:nonlinear-Sobolev-scale} for $\nabla v \in L^2(\Omega)$. In particular, 
$W^2_p(\Omega)$ is compactly embedded into $H^1(\Omega)$ according to
Lemma \ref{L:embedding} (Sobolev embedding).

Following \cite{BiDaDeVPe:02} and \cite{GaspozMorin:14}, we use a
greedy algorithm that is based on the knowledge of the element
errors and on bisection.
The algorithm hinges on \eqref{E:equidistribution}:
if $\delta>0$ is a given tolerance, the element error is
equidistributed and within tolerance  $e_\grid(v,T)\approx\delta$, and the global error
decays with maximum rate $(\#\grid)^{-1/2}$, then
\[
\delta^2 \#\grid
 \approx
 \sum_{T\in\grid} e_\grid(v,T)^2
 =
 |v-I_\grid v|^2_{H^1(\Omega)}
 \lesssim
 (\#\grid)^{-1}
\]
whence $\#\grid\lesssim\delta^{-1}$; here $I_\grid$ stands for the Lagrange interpolation operator.
With this in mind, we impose $e_T(v)\le\delta$ as a threshold to stop refining and expect
$\#\grid\Cleq\delta^{-1}$.
The following algorithm implements this idea.

\begin{algo}[greedy algorithm] \label{algo:greedy}
Given a
tolerance $\delta>0$ and a conforming mesh $\grid_0$, \textsf{GREEDY} finds
a conforming refinement $\grid\ge\grid_0$ of $\grid_0$ by bisection
such that $e_\grid(v,T)\le\delta$ for all $T\in\grid$: let $\grid=\grid_0$
and
\begin{algotab}
  \> [\grid] = $\textsf{GREEDY}(\grid,\delta,v)$\\
  \>  \quad while $\marked :=\{T\in\grid\,|\,e_\grid(v,T)>\delta\}\ne\emptyset$\\
  \> \> \quad $\grid := \REFINE(\grid,\marked)$\\
  \> \quad return $\grid$
\end{algotab}
\end{algo} 

Since $W^2_p(\Omega)\subset
C^0(\overline{\Omega})$, because $p>1$, we can use the Lagrange interpolant
and local estimate \eqref{E:error-lagrange} with 
$r=\sob(W^2_p)-\sob(H^1)=2-2/p>0$.
We deduce
\begin{equation}\label{E:termination}
	e_\grid(v,T)\lesssim h_T^r\,\|D^2 v\|_{L^p(T)}.
\end{equation}

We assess the quality of the resulting mesh in a slightly more general setting, following
\cite[Proposition~1 and Corollary 1]{BonitoCascon:2016},
needed later in Sections~\ref{S:conv-rates-coercive} and \ref{S:data-approx} for solution and
data approximation for any polynomial degree.

%----------------------------------------------------------------------------------------
\medskip\noindent
{\it An abstract greedy algorithm.}
%----------------------------------------------------------------------------------------
%
We consider a generic (possibly vector-valued) function $v\in L^q(\Omega, \mathbb R^M)$, with $M\geq 1$ and
$1\le q \le \infty$ and denote by $e_\grid(v,T) = \|v-\Pi_\grid v\|_{L^q(T)}$
the abstract $L^q$-local error for $T\in\grid$ used in the
\textsf{GREEDY} procedure and by $E_\grid(v) = \|v-\Pi_\grid v\|_{L^q(\Omega)}$ the global
$L^q$-interpolation error by either continuous of discontinuous piecewise polynomials (the definition
of $\Pi_\grid v$ is irrelevant now). We formulate the following assumptions:
\begin{itemize}
\item
{\it Summability in $\ell^q$}: The errors $\{ e_\grid(v,T) \}_{T \in \grid}$ satisfy
\begin{equation}\label{E:lq-summability}
\|v\|_{L^q(\Omega)}^q \lesssim \sum_{T\in\grid} e_\grid(v,T)^q.
\end{equation}
\end{itemize}

Rather than \eqref{E:weaker-reg} we assume that $v$ belongs to an abstract
space $X^s_p(\Omega;\grid_0)$ of functions with differentiability index $s\in (0,n]$ and 
integrability index $p \in (0,\infty]$ piecewise over $\grid_0$ with two crucial properties:
\begin{itemize}
\item
{\it Local error estimate}: For $r = s - \frac{d}{p} + \frac{d}{q} > 0$ and all $T\in\mesh$
\begin{equation}\label{A:termination_abstract}
e_\grid(v,T) \lesssim h_T^r | v |_{X^s_p(T)}
\end{equation}
\item
{\it Norm subadditivity}: For $p<\infty$, and obvious modification for $p=\infty$,
\begin{equation}\label{E:subadditive-norm}
\sum_{T \in \grid} | v |_{X^s_p(T)}^p \lesssim | v |_{X^s_p(\Omega;\grid_0)}^p.
\end{equation}
\end{itemize}
The space $X^s_p(\Omega;\grid_0)$ will later be either a Sobolev space $W^s_p(\Omega;\grid_0)$, 
a Besov space $B^s_{p,p}(\Omega;\grid_0)$ or a Lipschitz space ${\it Lip}^s_p(\Omega;\grid_0)$, with 
piecewise regularity over $\grid_0$; the latter two will
allow $0<p<1$. For the moment we do not need to be specific and just rely on the two properties above.

\begin{proposition}[abstract greedy error] \label{P:abstract_greedy}
{\it Let $\gridk[0]$ be an initial subdivision of $\Omega\subset\R^d$
satisfying the initial labeling property
\eqref{nv-initial-label} for $d=2$, or its variant for $d>2$.
Let $M\ge 1, 0< q,p \le \infty$ and $s-\frac{d}{p}+\frac{d}{q}>0$. Let 
$v\in L^q(\Omega,\R^M)$ satisfy \eqref{E:lq-summability},
\eqref{A:termination_abstract} and \eqref{E:subadditive-norm}.
Then {\rm\GREEDY}\!$(\grid_0,\delta,v)$ terminates in a finite number of iterations with
local errors verifying $e_\grid(v,T) \le \delta$ for all $T\in \grid$,
and there is a constant
$C = C(p,q,s,d,\Omega,\grid_0)$ such that the output $\grid\in\grids$ satisfies 
\begin{equation}\label{E:greedy_complexity_tools}
\|v-\Pi_\grid v\|_{L^q(\Omega)} \le C |v|_{X^s_p(\Omega;\grid_0)} \big(\#\grid-\#\grid_0\big)^{-\frac{s}{d}}.
\end{equation}
}
\end{proposition}
\begin{proof}
We proceed in several steps.

\step{1} {\it Termination}.
Since $h_T$ decreases monotonically to $0$ with bisection, so does $e_\grid(v,T)$ in view of
\eqref{A:termination_abstract}. Consequently, \textsf{GREEDY} terminates in finite number $k\ge 1$ of
iterations. Upon termination, the local errors satisfy $e_\grid(v,T) \le \delta$ for all $T\in\grid$
by construction, whence \eqref{E:lq-summability} implies
\[
\|v-\Pi_\grid v\|_{L^q(\Omega)} \lesssim \delta (\#\grid)^{\frac{1}{q}}.
\]

\step{2} {\it Counting}.
Let $\marked=\marked_0\cup\dots\cup\marked_{k-1}$ 
be the set of marked elements. We organize the elements in $\marked$ 
by size in  a way that allows for a counting argument. Let 
${\mathcal P}_j$ be the set of elements $T$ of $\marked$ with size
\[
 2^{-(j+1)}
 \le
 |T|<2^{-j}
\quad\Rightarrow\quad
 2^{-\frac{j+1}{d}}\le h_T<2^{-\frac{j}{d}},
\]
because $h_T=|T|^{1/d}$ for shape regular meshes $\grid\in\grids$.

We first observe that all $T$'s in ${\mathcal P}_j$ are
\emph{disjoint}. This is because if $T_1,\,T_2\in {\mathcal P}_j$ and $\mathring{T}_1\cap\mathring{T}_2\neq\emptyset$, then one of them is contained in the other, say $T_1\subset T_2$, due to the bisection procedure which works in any dimension $d\ge1$; see Section \ref{S:mesh-refinement}. Hence,
\[
	|T_1|\le\frac{1}{2}\,|T_2|
\]
contradicting the definition of ${\mathcal P}_j$. This implies the first bound
\begin{equation}\label{nv-small}
 2^{-(j+1)}\,\#{\mathcal P}_j
 \le
 |\Omega|
\quad\Rightarrow\quad
 \#{\mathcal P}_j\le|\Omega|\, 2^{j+1}.
\end{equation}
In light of \eqref{A:termination_abstract}, we have for $T\in {\mathcal P}_j$
\[
	\delta\le e_\grid(v,T) \lesssim 2^{-\frac{jr}{d}}|v|_{X^s_p(T)}.
\]
Therefore, accumulating these quantities in $\ell^p$ and invoking \eqref{E:subadditive-norm} yields
\[
\delta^{p}\,\#{\mathcal P}_j\lesssim 2^{-\frac{jrp}{d}}
\sum_{T\in{\mathcal P}_j}|v|_{X^s_p(T)}^p\lesssim 2^{-\frac{jrp}{d}}\, |v|_{X^s_p(\Omega;\grid_0)}^p
\]
and gives rise to the second bound 
\begin{equation}\label{nv-large}
	\#{\mathcal P}_j\lesssim \delta^{-p}\,2^{-\frac{jrp}{d}}\,|v|^p_{X^s_p(\Omega;\grid_0)}.
\end{equation}
\step{3} {\it Cardinality}. The two bounds for $\#{\mathcal P}$ in \eqref{nv-small} and
\eqref{nv-large} are complementary. The first one is good for $j$ small
whereas the second is suitable for $j$ large (think of $\delta\ll
1$). The crossover takes place for $j_0$ such that 
\[
2^{j_0+1}|\Omega| \approx \delta^{-p}\,2^{-\frac{j_0rp}{d}}|v|^p_{X^s_p(\Omega;\grid_0)}
\quad\Rightarrow\quad 
2^{j_0}\approx \left(|\Omega|^{-1} \delta^{-p} |v|_{X^s_p(\Omega;\grid_0)}^p\right)^{\frac{d}{d+rp}}.
\]
We now compute 
\[
\#\marked=\sum_j\#{\mathcal P}_j\lesssim
\sum_{j\le j_0}2^j|\Omega|+\delta^{-p}\,|v|^p_{X^s_p(\Omega)}
\sum_{j>j_0}(2^{-\frac{rp}{2}})^j.
\]
Since 
\[
\sum_{j\le j_0}2^j\approx 2^{j_0},
\qquad \sum_{j>j_0}(2^{-\frac{rp}{d}})^j\lesssim 2^{-\frac{rpj_0}{d}},
\]
we can write 
\[
\#\marked\lesssim |\Omega|^{1-\frac{d}{d+rp}}\big(\delta^{-1} | v |_{X^s_p(\Omega;\grid_0)}\big)^{\frac{dp}{d+rp}}.
\]
We finally apply Theorem \ref{nv-T:complexity-refine} (complexity of \REFINE) to arrive at 
\[
\#\grid-\#\grid_0\lesssim \#\marked\lesssim |\Omega|^{\frac{rp}{d+rp}} \big(\delta^{-1} | v |_{X^s_p(\Omega;\grid_0)}\big)^{\frac{dp}{d+rp}},
\]
or equivalently
\[
\delta \lesssim |\Omega|^{\frac{r}{d}} |v|_{X^s_p(\Omega;\grid_0)} \big(\#\grid - \#\grid_0 \big)^{-\frac{d+rp}{dp}}.
\]

\step{4} {\it Total error}. Since $\frac{d+rp}{dp} = \frac{s}{d} + \frac{1}{q}$ we deduce from Step 1
\[
\|v-\Pi_\grid v\|_{L^q(\Omega)} \lesssim \delta (\#\grid)^{\frac{1}{q}} \lesssim |\Omega|^{\frac{r}{d}}
|v|_{X^s_p(\Omega;\grid_0)} (\#\grid-\#\grid_0)^{-\frac{s}{d}},
\]
which is the desired estimate.
\end{proof}

The output mesh $\grid$ of $\textsf{GREEDY}(\grid_0,\delta,v)$ starting from $\grid_0$ satisfies
$\#\grid\ge c_0\#\grid_0$ for some $c_0>1$, whence
$\#\grid - \#\grid_0 \ge \big(1-\frac{1}{c_0} \big) \#\grid$ and \eqref{E:greedy_complexity_tools} yields
\begin{equation}\label{E:greedy_grid}
\|v-\Pi_\grid v\|_{L^q(\Omega)} \lesssim C |v|_{X^s_p(\Omega;\grid_0)} \big(\#\grid\big)^{-\frac{s}{d}},
\end{equation}
where $C$ depends on $c_0$. It will be convenient in many applications of $\textsf{GREEDY}$,
to be discussed later in Sections \ref{S:conv-rates-coercive} and \ref{S:data-approx}, that
the starting mesh be a conforming refinement of $\grid_0$ to enhance its efficiency. We
will prove in Section \ref{S:greedy} that \eqref{E:greedy_grid} remains valid.

It is instructive to realize that \textsf{GREEDY} is a practical algorithm that hinges on the
different summabilities of \eqref{E:lq-summability} and \eqref{E:subadditive-norm}, and delivers
a global $L^q$-error consistent with \eqref{E:error-decay} of Section \ref{S:error-equidistribution}.
Moreover, the outcome graded grid $\grid$ is quasi-optimal but may not equidistribute the error,
not even approximately.

We are now in a position to show that \textsf{GREEDY} constructs optimal graded meshes for the interpolation error in $H^1(\Omega)$ alluded to at the beginning of this section. To this end, we let $I_\grid$
be the Lagrange interpolation operator for $d=2$. \looseness=-1

\begin{corollary}[optimal $H^1$-convergence rate]\label{C:optimal-mesh}
{\it If $v\in H^1(\Omega)\cap W^2_p(\Omega)$ for $1<p\leq 2$ and $d=2$, then {\rm \GREEDY} 
yields graded bisection meshes $\grid$ so that}
$$
| v - I_\grid v|_{H^1(\Omega)} \lesssim | \Omega |^{1-1/p} \|D^2 v\|_{L^p(\Omega)} (\#\grid)^{-1/2}.
$$
\end{corollary}
\begin{proof}
We invoke Proposition~\ref{P:abstract_greedy} (abstract greedy error) and equation \eqref{E:greedy_grid}
for $\nabla v \in L^2(\Omega,\R^2)$
with $\Pi_\grid \nabla w = \nabla I_\grid v$ and $s=1$, $q=2$, $p>1$ whence $s-\frac{d}{p}+\frac{d}{q}>0$.
\end{proof}

\begin{remark}[piecewise $W^2_p$-smoothness]\label{R:pw-smooth}
Since \eqref{E:termination} is completely local for $d=2$, we see from \eqref{E:subadditive-norm}
that it suffices for $v\in H^1(\Omega)$ to be piecewise in $W^2_p$ over the initial partition $\grid_0$,
namely $W^2_p(\Omega;\grid_0)$. It turns out that
this statement is valid for any dimension $d\ge2$ in view of
Proposition \ref{P:cont-vs-discont} (approximation of gradients). We will revisit this
issue in Section \ref{S:Besov}.
\end{remark}

\begin{remark}[case $p<1$]\label{R:p<1}
We consider now polynomial degree $n\ge1$. The integrability $p$ 
corresponding to differentiability $n+1$ results from equating
Sobolev numbers:
\[
n+1-\frac{d}{p}=\sob(H^1)=1-\frac{d}{2}
\quad\Rightarrow\quad p=\frac{2d}{2n+d}.
\]
Depending on $d\ge2$ and $n\ge1$, this may lead to $0<p<1$, in which 
case $W^{n+1}_p(\Omega)$ is to be replaced by the Besov space $B^s_{p,p}(\Omega)$ for $s<n+1$
or the Lipschitz space
${\it Lip}^{n+1}_p(\Omega)$ \cite{DeVore:98}. We will discuss this matter in Section \ref{S:Besov}
and make the abstract greedy setting precise.
\end{remark}

\begin{remark}[isotropic vs anisotropic elements]\label{R:anisotropic}
Since geometric singularities are of the form \eqref{E:corner-sing} for $d=2$,
Corollary \ref{C:optimal-mesh} (optimal $H^1$-convergence rate) shows
that isotropic graded meshes are able to deliver optimal convergence rates for $d=2$.
Unfortunately, this is no longer the case for $d>2$ and
anisotropic meshes are necessary for optimal meshes. This topic is delicate is several respects.
Deriving reliable and efficient a posteriori error estimators is largely open for anisotropic meshes;
this is the subject of Section \ref{S:aposteriori} for isotropic meshes. Even having such estimators,
building a theory of adaptivity is open; this is the subject of Sections \ref{S:convergence-coercive},
\ref{S:conv-rates-coercive} and \ref{S:conv-rates-infsup} for isotropic meshes. Finally, constructing
anisotropic meshes based on a posteriori information alone and that easily allow for refinement and
coarsening is problematic.
For these reasons we do not dwell on anisotropic refinement in this survey.
\end{remark}

%-------------------------------------------------------------------------------------
\subsection{Nonconforming meshes}\label{S:nonconforming-meshes}
%-------------------------------------------------------------------------------------

\begin{figure}[ht]
\begin{center}
\includegraphics[width=0.27\textwidth]{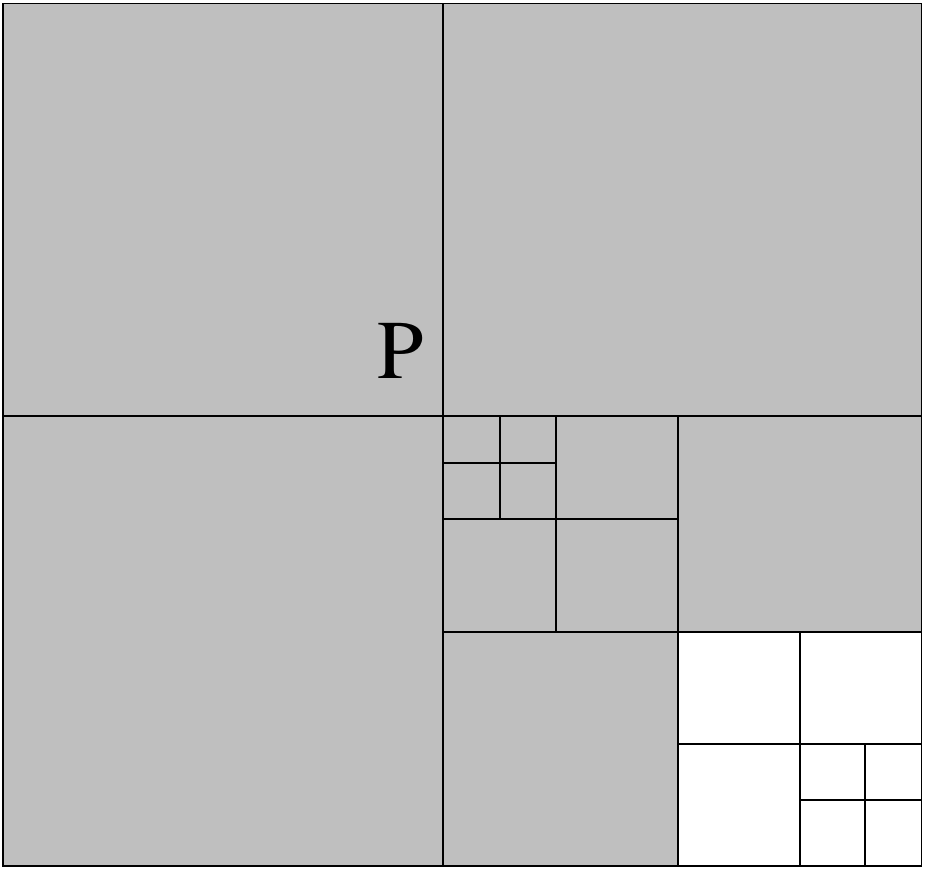}
\hskip0.3cm
\includegraphics[width=0.32\textwidth]{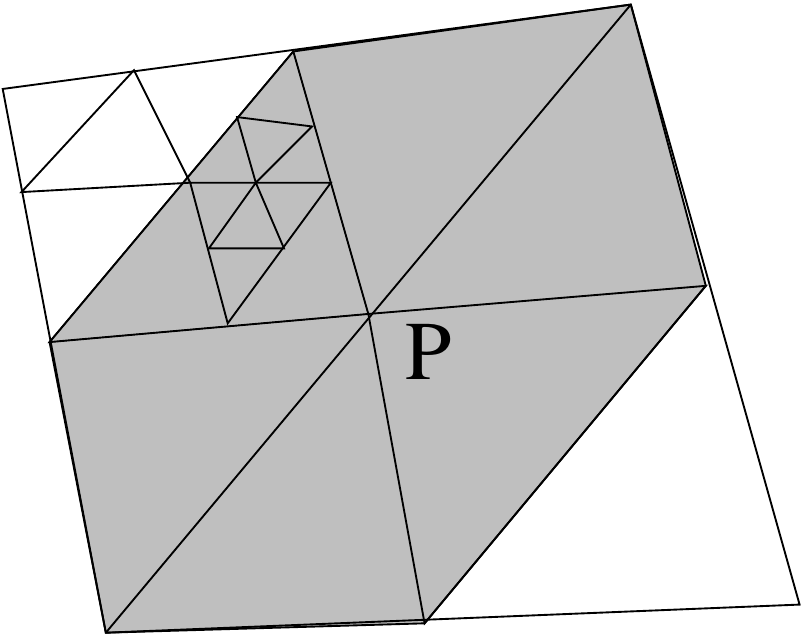}
\hskip0.2cm
\includegraphics[width=0.32\textwidth]{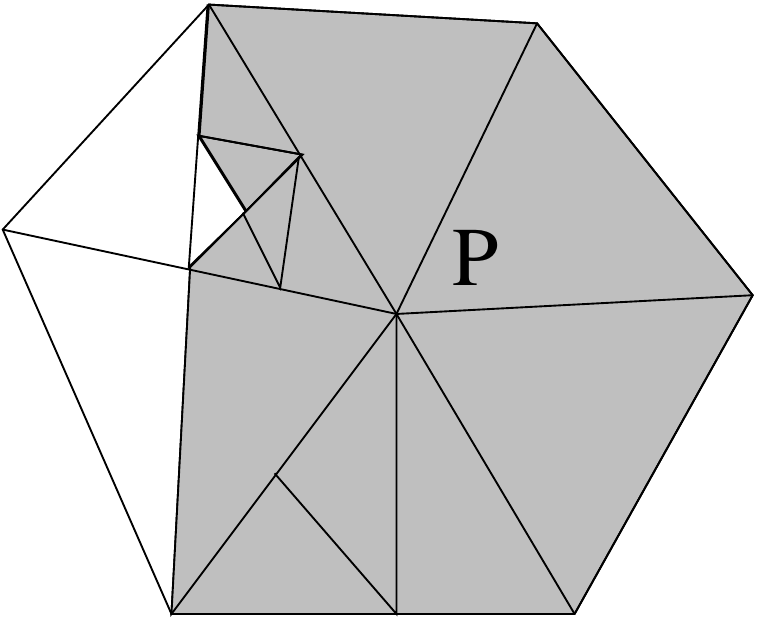}
\caption{Nonconforming meshes made of quadrilaterals (a), triangles with 
red refinement (b), and triangles with bisection (c). The shaded
regions depict the domain of influence of a proper or conforming node $P$.}
\label{nv-F:nonconforming-meshes}
\end{center}
\end{figure}

More general subdivisions of $\Omega$ than those in 
\S \ref{S:bisection} are used in practice.
If the elements of $\gridk[0]$ are quadrilaterals for $d=2$, or their
multidimensional variant for $d>2$, then it is natural to allow for improper
or {\it hanging nodes} for the resulting refinements $\grid$ to be
graded; see Figure \ref{nv-F:nonconforming-meshes} (a).
On the other hand, if $\gridk[0]$ is made of triangles for $d=2$,
or simplices for $d>2$, then red refinement without green
completion also gives rise to graded meshes with hanging nodes;
see Figure \ref{nv-F:nonconforming-meshes} (b). In both cases, the
presence of hanging nodes is inevitable to enforce mesh grading.
Finally, bisection may produce meshes with hanging nodes, as depicted
in Figure \ref{nv-F:nonconforming-meshes} (c), if the completion
process is incomplete. All three refinements maintain shape
regularity, but for both practice and theory, they
cannot be arbitrary: we need to restrict the level of nonconformity.
We discuss this next, starting with the case of polynomial degree $n=1$
\cite{bonito2010quasi,BeCaNoVaVe23}. \looseness=-1

We say that a node $P$ of $\grid$ is a {\it proper} (or conforming) node if
it is a vertex of all elements containing $P$; otherwise, we say that
$P$ is an {\it improper} (nonconforming or hanging) node. The set ${\cal N}$ of all nodes of $\grid$ is thus partitioned into the set ${\cal P}$ \index{Meshes!$\cal P$: proper nodes} of proper nodes, and the set ${\cal H}={\cal N}\setminus {\cal P}$ of hanging nodes. 

\begin{figure}[!h]\label{f-globalindex}
\begin{center}
\begin{overpic}[scale=0.25]{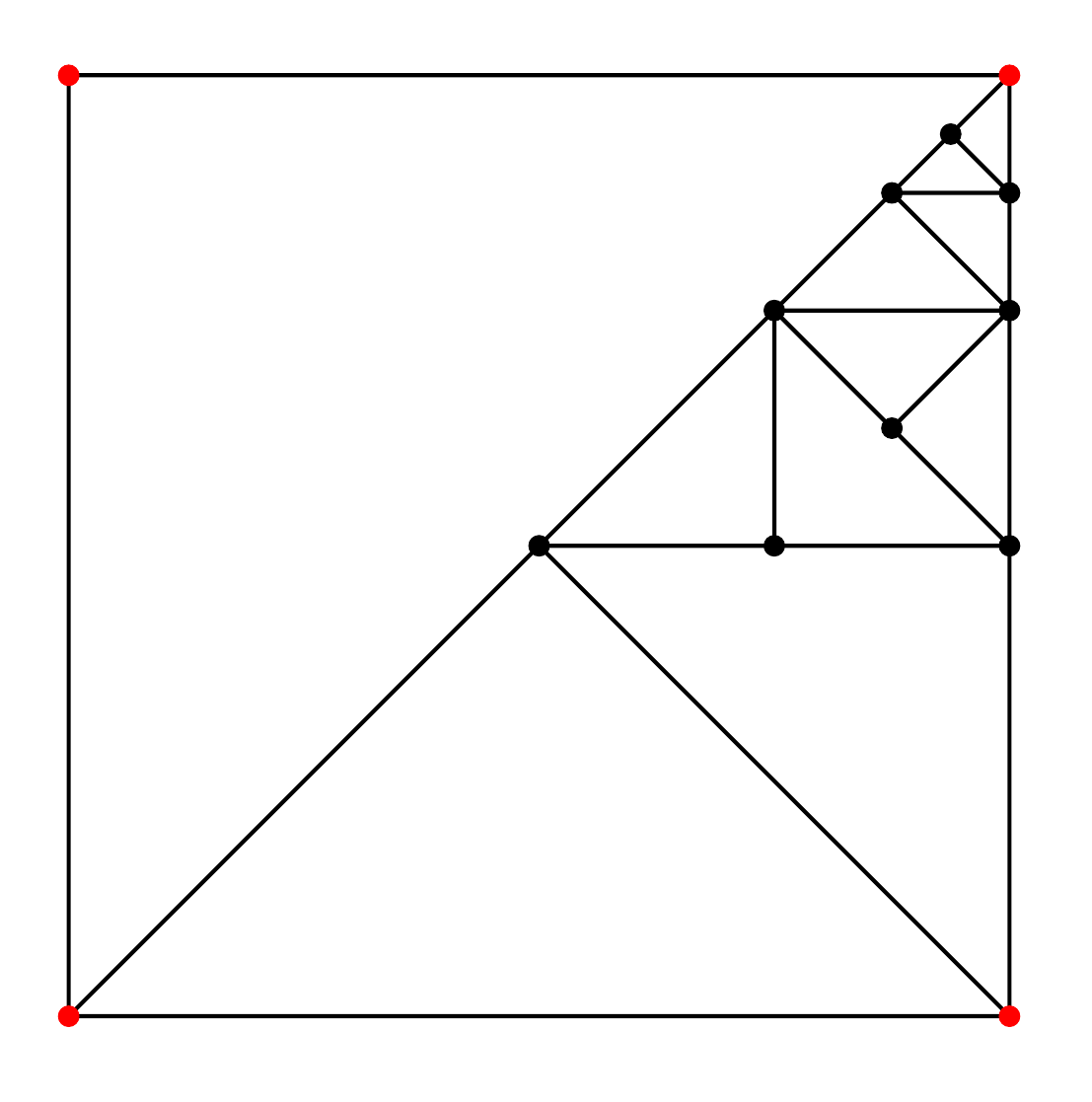}
\put( 2, 2){$0$}
\put( 2,95){$0$}
\put(92, 2){$0$}
\put(92,95){$0$}
\put(45,52){$1$}
\put(93,52){$1$}
\put(65,70){$2$}
\put(69,44){$2$}
\put(93,70){$2$}
\put(75,80){$3$}
\put(79,54){$3$}
\put(93,80){$3$}
\put(82,87){$4$}
\end{overpic}
\end{center}
\caption{Example of distributions of proper nodes (red) and hanging nodes (black), with associated global indices $\lambda$. }
\end{figure}
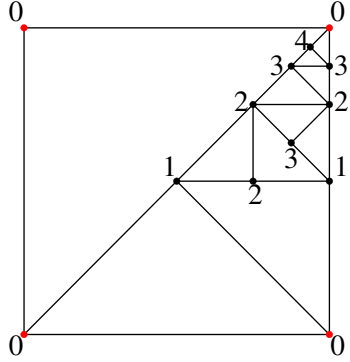

A useful notion in dealing with hanging nodes is the {\em global index} of a node, introduced in \cite{BeCaNoVaVe23}: it measures the number of nonconforming refinements needed to generate a hanging node from proper nodes. To define it, for any $x \in {\cal H}$ which has been generated by the bisection of an edge $[x',x'']$, let us set ${\cal B}(x)=\{x',x''\}$.
\begin{definition}[global index of a node]\label{d-globalindex}
\index{Definitions!Global index of a node}
The global index $\lambda(x)$ \index{Meshes!$\lambda(x)$: global index of a node $x\in \mathcal N$} of a node $x \in {\cal N}$ is defined recursively as follows:
\begin{itemize}
\item if $x \in {\cal P}$, set $\lambda(x)=0$;
\item if $x \in {\cal H}$ and ${\cal B}(x)=\{x',x''\}$, set $\lambda(x)=\max(\lambda(x'), \lambda(x'')) +1 $.
\end{itemize}
\end{definition}
The set of all nodes of $\grid$ is thus partitioned according to the value of the global index: for any integer $l \geq 0$, we set ${\cal H}_l = \{ x \in {\cal N} : \lambda(x)=l\}$. Note that ${\cal H}_0 = {\cal P}$. An example of distribution of global indices for $d=2$ is shown in Figure \ref{f-globalindex}.

We define the global index of the triangulation $\grid$ by %setting
$\lambda(\grid) := \max_{x \in {\cal N}} \lambda(x)$. The level of nonconformity of the triangulations we are dealing with is controlled by the following condition of admissibility.
\begin{definition}[$\Lambda$-admissibility]\label{d-LambdaAdmissibility} Let $\Lambda \geq 0$ be an integer. A refinement $\grid$ of $\grid_0$ is $\Lambda$-admissible if
\begin{equation}\label{nv-level-nonconformity}
\index{Constants!$\Lambda$: $\Lambda$-admissibility constant}
\lambda(\grid) \leq \Lambda \,.
\end{equation}
If $\lambda(\grid)\ge1$, then $\grid$ is nonconforming, but otherwise $\grid\in\grids$ is conforming
if $\lambda(\grid)=0$.
The collection of all $\Lambda$-admissible partitions is denoted by $\grids^\Lambda$.
\index{Meshes!$\grids^\Lambda$: set of all $\Lambda$-admissibility refinements of $\grid_0$}
\end{definition}

$\Lambda$-admissibility has the following basic implications. 
\begin{proposition}[properties of $\Lambda$-admissible partitions]\label{P:property-Lambda-adm}
Let $T$ be any element of a $\Lambda$-admissible partition $\grid$.
\begin{enumerate}
\item If $e \subset \partial T$ is an edge of  $T$, then $e$ may contain at most $2^\Lambda-1$ hanging nodes. 
\item If $e \subset \partial T$ is an edge of some other element $T'$, then $h_{T'} \simeq h_T$, where the hidden constants only depend on the shape of the initial triangulation $\grid_0$ and possibly on $\Lambda$.
\end{enumerate}
\end{proposition}
\begin{proof} {\it (i)} stems from the fact that the edge may contain at most $2^{k-1}$ hanging nodes of level $k$ for $1 \leq k \leq \Lambda$. To prove {\it (ii)} we observe that the length ratio $\frac{|\bar{e}|}{|e|}$, where $\bar{e}$ is the edge of $T$ containing $e$, is at most $2^\Lambda$, and we conclude invoking the shape regularity of the partition.  
\end{proof}

\begin{figure}\label{f-basis-function}
\begin{center}
\includegraphics[width=0.45\textwidth]{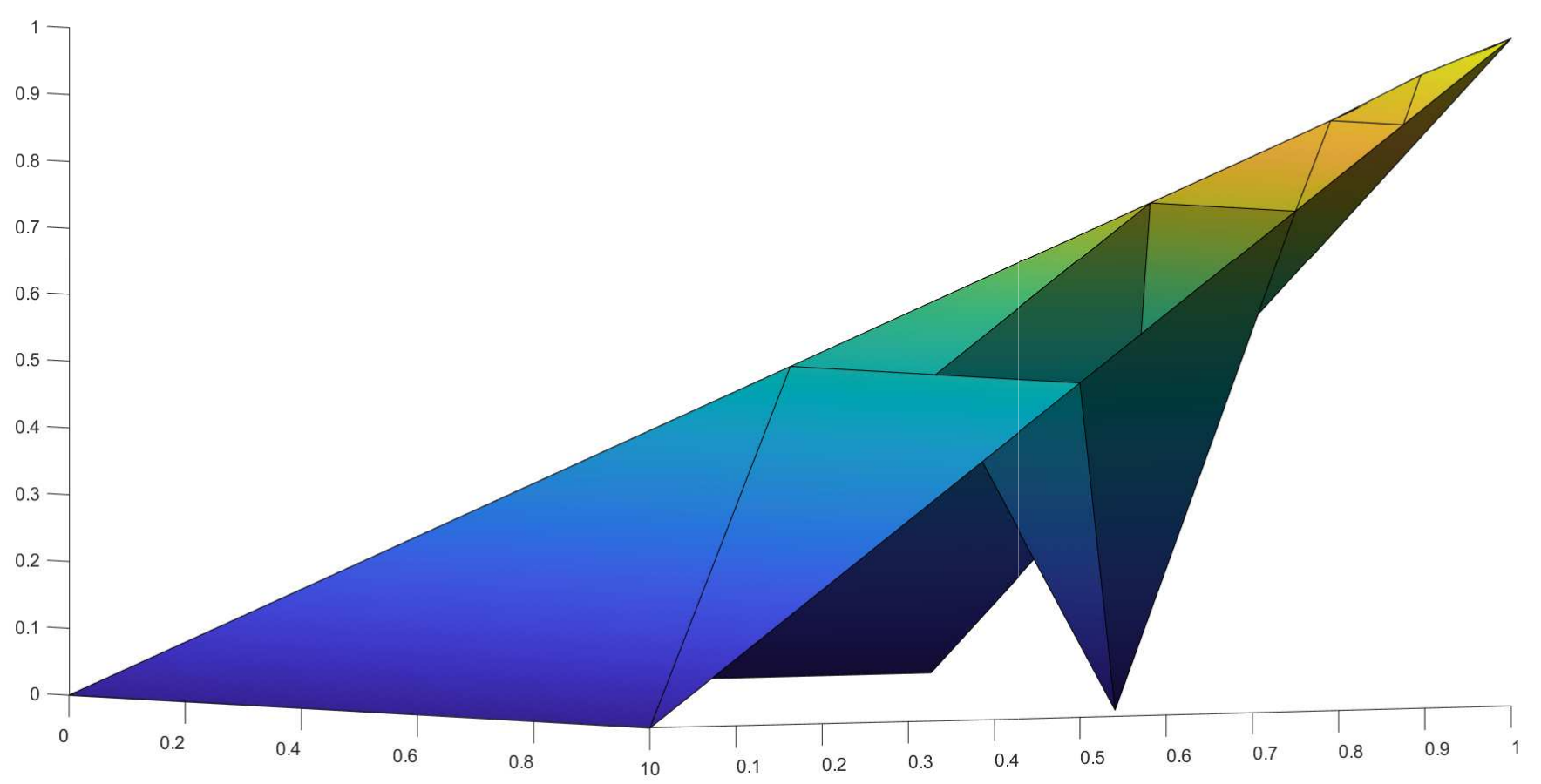}
\end{center}
\caption{Plot of the basis function $\phi_P$ on the nonconforming triangulation shown in Figure \ref{f-globalindex}, for $P$ equal to the upper right corner of the domain, after using bisection to convert the lowest hanging nodes with global index $3$ into a proper node.} 
\end{figure}

In the space $\V_\grid$ of continuous piecewise-linear maps over $\grid$, functions  are uniquely defined by their values at the proper nodes of $\grid$. So  it is natural to introduce
the canonical continuous piecewise-linear basis functions $\phi_P$
associated with any proper node $P$. They satisfy
\begin{equation}\label{e:representation_proper}
v = \sum_{P\in {\cal P}} v(P) \phi_P \qquad\forall ~v\in\V_\grid\,,
\end{equation}
and are defined by the conditions: $\phi_P \in \V_\grid$ and
\begin{itemize}
\item $\phi_P(z)=1$ if $z=P$,  \ $\phi_P(z)=0$ if $z \in {\cal P}\setminus \{P\}$.
\end{itemize}
The values of $\phi_P$ at the hanging nodes, hence everywhere in the domain, can be reconstructed by linear interpolation as follows: assuming that $\phi_P$ has been defined at all nodes of global index $< l$, if $z \in {\cal H}_l$ and ${\cal B}(z)=\{z',z''\}$, then 
$$
\phi_P(z)= \frac12\big(\phi_P(z')+\phi_P(z'')\big)\,.
$$
An example of basis function $\phi_P$ on a nonconforming triangulation is provided in Figure \ref{f-basis-function}.

The {\em domain of influence} of a proper node $P$ is the set
\index{Meshes!$\omega_\grid(P)$: domain of influence of a proper node $P$}
\[
\omega_\grid(P) = \supp (\phi_P)\,,
\]
highlighted in grey in Figure \ref{nv-F:nonconforming-meshes}; this notion was introduced in
\cite{BabuskaMiller:87} in the context of $K$-meshes; see also \cite{bonito2010quasi}.  
To identify elements $T \in \grid$ contained in $\omega_\grid(P)$,  we introduce for any node $x \in {\cal N}$ the set ${\cal P}(x)$ of the {\em proper nodes influencing $x$}, which is defined recursively as follows:
\begin{itemize}
\item initialize ${\cal P}(x)= \{x\}$; 
\item while ${\cal P}(x) \cap {\cal H} \not = \emptyset$, if $y \in {\cal P}(x) \cap {\cal H}$ replace ${\cal P}(x)$ by $({\cal P}(x) \setminus \{y\})\cup {\cal B}(y)$.
\end{itemize} 
Then, $ T \subseteq \omega_\grid(P)$ if and only if $P$ influences some vertex of $T$, i.e., there exists a vertex $v$ of $T$ such that $P \in {\cal P}(v)$.

One of the consequences of the $\Lambda$-admissibility assumption of $\grid$ is the following result, which says that all elements $T$ contained in $\omega_\grid(P)$ have comparable size.
\begin{proposition}[size of the domain of influence]\label{p-size-domain-influence}
There exists a positive constant $C=C(\grid_0,\Lambda)$, only depending on the shape of the initial triangulation $\grid_0$ and possibly on $\Lambda$, such that for any $P \in {\cal P}$
\[
{\rm diam \,} \omega_\grid(P) \leq C \, h_T  \qquad \forall T \in \grid, \ T \subseteq  \omega_\grid(P)\,.
\]
\end{proposition}
\begin{proof}
Elements in $\omega_\grid(P)$ having $P$ as a vertex share in pairs an edge or a portion of an edge, hence -- as noted above -- $\Lambda$-admissibility implies the existence of a characteristic size, say $h_P$, which is comparable to the diameter of each of them.
On the other hand, any $T \subset \omega_\grid(P)$ not containing $P$ has at least one vertex $v_T \in {\cal H}$ such that  $P \in {\cal P}(v_T)$. Thus, there exists a sequence $\{y_k : 0 \leq k \leq K\}$ of vertices satisfying $y_0=v_T$, $y_K=P$, and $y_{k+1} \in {\cal B}(y_k)$ for $0 \leq k < K$; since $\lambda(y_k) \geq \lambda(y_{k+1})+1$, necessarily $K \leq \Lambda$. Correspondingly, we can find a chain of at most $K$ elements, starting at $T$ and ending at an element containing $P$, which share in pairs an edge or a portion of an edge. We deduce that $h_T \simeq h_P$, and ${\rm dist}(T,P) \simeq h_T$, where the hidden constants may depend only on the shape of the initial triangulation and on $\Lambda$. The conclusion easily follows from these results.
\end{proof}

We now turn to the case of polynomial degree $n>1$; we refer to \cite{CaFa23} for more details. The concept of hanging node is no longer solely related to the geometry of the mesh, but also to the distribution of degrees of freedom along the edges of the elements. For instance, consider a full edge $e$ shared by two triangles $T$ and $T'$, and bisect $T'$ to create two new elements $T_1$ and $T_2$ having $e$ as a vertex. If we use quadratic Lagrangian elements, the midpoint $x$ of $e$ carries a degree of freedom for the three elements that share it, so we do not consider it as a hanging node; on the other hand, the nodes at distance $\frac14 |e|$ and $\frac34 |e|$ from an endpoint of $e$ are hanging nodes (despite they are vertices of no triangle) since they do not carry a degree of freedom for the element $T$. If we move to cubic Lagrangian elements, then $x$ becomes a hanging node, together with the nodes at distance $\frac16 |e|$ and $\frac56 |e|$ from an endpoint of $e$, whereas the nodes at distance $\frac13 |e|$ and $\frac23 |e|$ are not hanging nodes, since they carry a degree of freedom for each triangle they belong to.

In general, for a partition $\grid$ made of classical affine Lagrangian or Hermitian elements, the hanging nodes are defined as follows. 
\begin{itemize}
    \item 
    Given an element $T\in\grid$, the set ${\cal P}_T$ of the {\em proper
nodes of $T$} is made of all images of the reference $n$-lattice via the
affine transformation.
%    Given an element $T\in\grid$, the set ${\cal P}_T$ of the {\em proper nodes of $T$} is made of those images of the reference $n$-lattice that sit on the boundary $\partial T$. 
The set ${\cal H}_T$ of the {\em hanging nodes of $T$} collects the points of $\partial T$ that are not proper nodes of $T$, but are proper nodes of some other contiguous element $T'$. The set of all nodes of $T$ is ${\cal N}_T := {\cal P}_T \cup {\cal H}_T$.
    \item At the global level, if ${\cal N}=\bigcup_{T \in \grid} {\cal N}_T$ is the set of all nodes of $\grid$, the set ${\cal P} \subseteq {\cal N}$ of the {\em proper nodes of $\grid$} \index{Meshes!$\cal P$: proper nodes} contains those nodes that are proper nodes for all elements they belong to. The complementary set ${\cal H}:={\cal N}\setminus {\cal P}$ is the set of the {\em hanging nodes of $\grid$}.    
\end{itemize}
In other words, a hanging node of $\grid$ is a point that carries a degree of freedom for some, but not all elements it belong to. With this definition of proper nodes, representation \eqref{e:representation_proper} of continuous piecewise-linear maps extends to $n>1$.

The {\em global index} $\lambda(x)$ of a node
$x \in {\cal N}$ is precisely defined as in Definition \ref{d-globalindex}. The set ${\cal B}(x) \subset {\cal N}$ collects the endpoints of an interval $[x', x'']$, contained in the skeleton of $\grid$, that has been bisected when $x$ has been created, and contains no other node inside.

\begin{figure}\label{F:global_index_n}
\begin{center}
\includegraphics[width=0.8\textwidth]{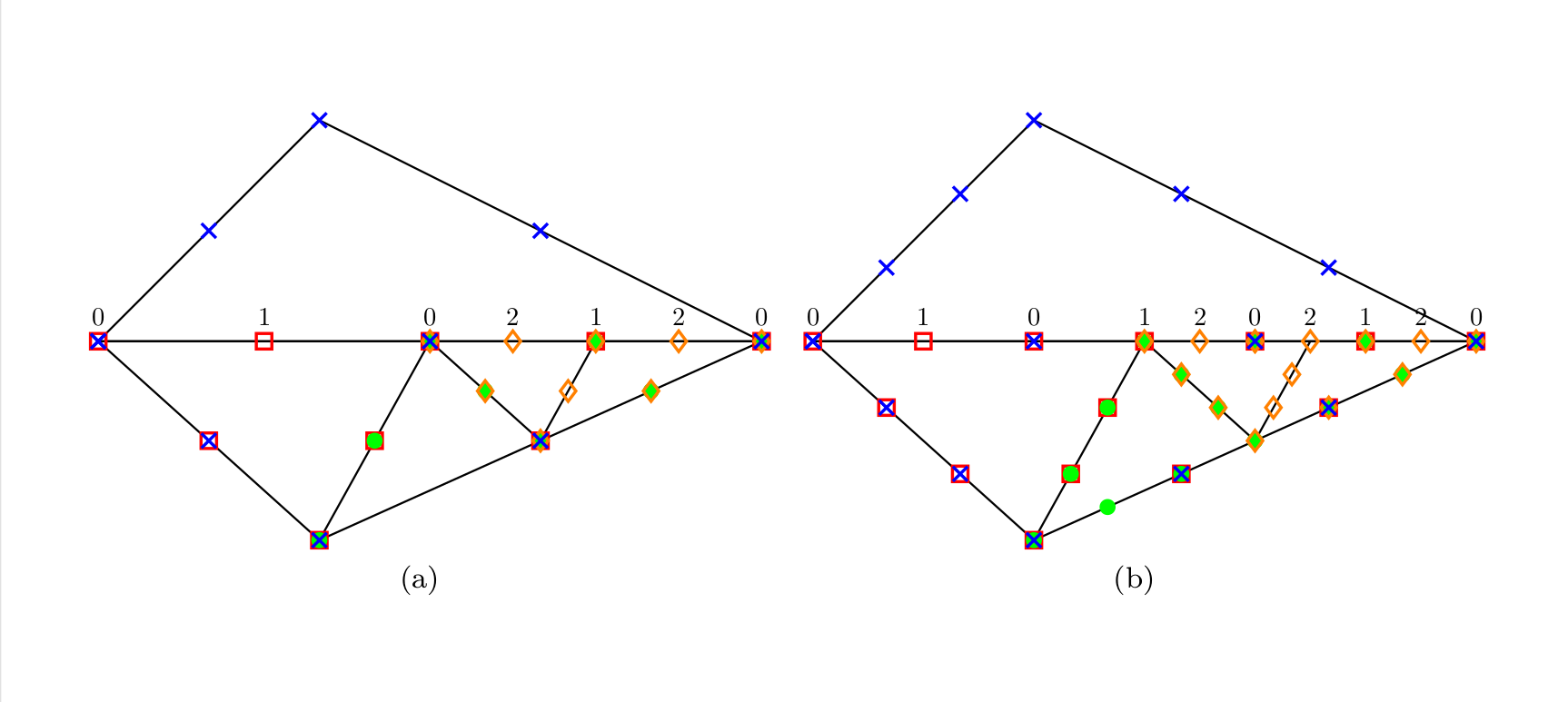}
\end{center}
\caption{Triangulation after the three refinements in the case $n=2$ (a) and in the case $n=3$ (b). Blue crosses represent the original degrees of freedom on the initial conforming mesh. Red squares, green circles and orange triangles are used for the degrees of freedom of the first, second and third refinement, respectively. All nodes are proper, except those on the horizontal line, whose global index is reported.}
\end{figure}

Figure \ref{F:global_index_n} provides two examples,  for $n=2$ and $n=3$, of distributions of hanging nodes and corresponding global indices, created by successive bisections starting from an initial conforming partition.

The concept of $\Lambda$-admissibility, given in Definition \ref{d-LambdaAdmissibility}, remains unchanged for $n>1$. The statements in Proposition~\ref{P:property-Lambda-adm}, too, extend to the higher order case; the maximum number of hanging nodes on an edge being now $O(n2^\Lambda)$. Consequently, the conclusion of Proposition~\ref{p-size-domain-influence} remains valid when $n>1$ as well: there is a constant $C=C(\grid_0,\Lambda)$ such that
\begin{equation}\label{e:domain_influence_control}
{\rm diam \,} \omega_\grid(P) \leq C \, h_T  \qquad \forall T \in \grid, \ T \subseteq  \omega_\grid(P)\,.
\end{equation}

\begin{remark}[quadrilateral and hexahedral partitions]\label{R:Lambda-quad}
It is readily seen that the definitions of global index of a node and $\Lambda$-admissible partition extend seamlessly to shape-regular meshes made of quadrilaterals refined by a quadtree strategy ($d=2$) (see Fig. \ref{nv-F:nonconforming-meshes} (a) for an example), or by hexahedra refined by an octree strategy ($d=3$). The same holds for heterogeneous partitions made of a combination of simplices and hexahedra. All results reported above are valid for such partitions.
We refer to \cite{bonito2010quasi} for details.
\end{remark}

\paragraph{$\Lambda$-admissible meshes under refinement.}
Given a $\Lambda$-admissible grid $\grid$, a subset $\marked$ 
of elements marked for refinement, and a desired number $b \ge1$ of
subdivisions to be performed in each marked element, the procedure
\[
\gridk[*] = \REFINE (\grid,\marked, \Lambda)
\]
creates a minimal $\Lambda$-admissible mesh $\gridk[*]\ge\grid$ such that all
the elements of $\marked$ are subdivided at least $b$ times. In order
for $\gridk[*]$ to be $\Lambda$-admissible, perhaps other elements not in
$\marked$ must be partitioned. Despite the fact that admissibility is a 
constraint on the refinement procedure weaker than conformity,
it cannot avoid the propagation of refinements beyond $\marked$.
The complexity of $\REFINE$ is again an issue which we discuss
in \S \ref{SS:nonconforming}: we show that 
Theorem \ref{nv-T:complexity-refine} extends to this case.

\begin{theorem}[complexity of $\REFINE$ for $\Lambda$-admissible meshes]\label{T:nonconforming-meshes}
Let $\gridk[0]$ be an arbitrary conforming partition of $\Omega$, 
except for the bisection algorithm in which case $\gridk[0]$ satisfies the labeling
\eqref{nv-initial-label} for $d=2$ or its higher dimensional counterpart
\cite{Stevenson:08}. Then the estimate
\[
\#\gridk - \#\gridk[0] \le \Ccompl \sum_{j=0}^{k-1} \# \markedk[j]
\quad\forall k\ge 1
\]
\index{Constants!$\Ccompl$: complexity of \REFINE constant}
holds with a constant $\Ccompl$ depending on $\gridk[0]$, $d$, $n$ and $\Lambda$.
\end{theorem}

The following result about uniform refinements of $\Lambda$-admissible partitions will be used in the sequel. {\it The uniform refinement $\grid_*$ of a partition $\grid \in \grids^\Lambda$ is the partition obtained by bisecting $d$ times each element of $\grid$}. This implies, in particular, that each edge of $\grid$ is bisected once.
\begin{proposition}[$\Lambda$-admissibility of uniform refinements]\label{P:Lambda-uniform}
If $\grid \in \grids^\Lambda$ is a $\Lambda$-admissible partition and $\grid_*$ is its uniform refinement, then $\grid_*$ is $\Lambda$-admissible.
\end{proposition}
\begin{proof}
A simple recursion argument on the global index of the hanging nodes of $\grid$ shows that after refinement each such node either becomes a proper node or its global index is reduced by 1. At the same time, new nodes are created by the refinement, whose global index is at most 1 plus the maximal global index of the pre-existing nodes. In both cases, the maximal global index of $\grid_*$ cannot exceed $\Lambda$. \looseness=-1
\end{proof}

A simple consequence is the following result, which is useful to control the meshsize between consecutive refinements.
\begin{corollary}[bound on the refinements]\label{C:Lambda-refined}
{\rm \REFINE} with $b=1$ never refines an element of a $\Lambda$-admissible partition $\grid$ more than $d$ times.
\end{corollary}
\begin{proof}
\REFINE gives the smallest $\Lambda$-admissible mesh $\grid_*$ such that all the marked elements of $\grid$ have been refined. Since the uniform refinement of $\grid$ remains $\Lambda$-admissible, the minimality of $\grid_*$ implies that no element of the marked set can be refined more than d times.    
\end{proof}

We conclude by emphasizing that the polynomial interpolation and adaptive
approximation theories of Sects. \ref{S:pol-interp}
and \ref{S:constructive-approximation} extend to nonconforming meshes with fixed
level of incompatibility as well.

%--------------------------------------------------------------------------------

%
%
%
\section{A Posteriori Error Analysis}\label{S:aposteriori}
%
%
% intro: a posteriori error analysis and its goals
Numerical solutions to a boundary value problem serve to approximate its unknown exact solution. In such a context, it is of interest
\begin{itemize}
\item to quantify the error of the numerical solution,
m\item to gain information for \emph{adapting} the discretization to the exact solution
\end{itemize}
in a computationally accessible manner. These are the two goals of an \emph{a~posteriori error analysis}, where the adjective `a~posteriori' hints to the fact that the numerical solution itself can be involved. To achieve the two goals, the a~posteriori analysis individuates so-called \emph{error estimators} that, ideally, are computable, split into local contributions called \emph{indicators}, and bound the error from above and below.

% this section
This section exemplifies such an analysis, considering the numerical solution of the boundary value problem \eqref{strong-form}, i.e.\
\begin{equation*}
 -\operatorname{div}(\bA\nabla u) + cu = f \text{ in }\Omega,
\quad
 u = 0 \text{ on }\partial\Omega, 
\end{equation*}
with Lagrange elements of arbitrary fixed order $n \geq 1$. Throughout this section, we adopt the notations and assumptions of previous sections for this model setting. In particular, the \emph{exact solution} $u \in H^1_0(\Omega)$ solves the variational problem \eqref{weak-form} and, given a simplicial conforming mesh $\grid \in \grids$ of $\Omega$ and finite element space
\index{Functions!$u_\grid$: Galerkin approximation}
\begin{equation*}
\V_\grid
 \definedas
 \big\{ v \in \mathbb S^{n,0}_\grid \mid v|_{\Omega} = 0 \big\}
 \subset
 H^1_0(\Omega),
\end{equation*}
the \emph{Galerkin approximation} solves
\begin{equation}
\label{setting-apost:Galerkin}
u_\grid \in \V_\grid: \quad  \B[u_\grid,w] = \langle f, w \rangle \quad \forall \, w \in \V_\grid,
\end{equation}
with the bilinear form $\B$ from \eqref{bilinear-coercive}.

% key feature of given approach 
We stress that the analysis will be conducted under the regularity assumptions
\begin{equation}
\label{setting-apost;nat-regularity}
 \bA \in L^\infty\big( \Omega; \R^{d \times d} \big),
\quad
 c \in L^\infty(\Omega),
\quad
 f \in H^{-1}(\Omega) = H^1_0(\Omega)^*,
\end{equation}
used in Sect.\ \ref{S:inf-sup} to establish existence and uniqueness of the exact solution. This fact distinguishes the approach below, which builds on \cite{KreuzerVeeser:2019}, from most other ones requiring additional regularity; cf., e.g., \cite{Verfuerth:13}. Notably, this difference not only allows for covering more examples but is also related to strengthening the relationship between error and estimator to a true equivalence on any admissible mesh $\grid \in \grids$.

% roles of forcing vs coefficients
It is useful to recall two differences between the forcing $f$ and the coefficients $(\bA,c)$. First, while the exact solution $u$ depends linearly on the forcing $f$, it depends nonlinearly on the diffusion tensor $\bA$ and on the reaction coefficient $c$. To state the second difference, let $u \in H^1_0(\Omega)$ and note that the assumptions \eqref{setting-apost;nat-regularity} on $(\bA,c)$ imply the ``missing'' one $f \in H^{-1}(\Omega)$. On the other hand, the assumptions on $(A,f)$ imply only $cu \in H^{-1}(\Omega)$ while the assumptions on $(c,f)$ imply only $-\operatorname{div}(\bA \nabla u) \in H^{-1}(\Omega)$. These conditions are weaker than the ``missing'' one $u \in H^1_0(\Omega)$,
and are due to the multiplicative role of $(\bA,c)$ in the differential equation. \looseness=-1

% plan
In order to elucidate the new twists allowing for \eqref{setting-apost;nat-regularity}, this section is organized as follows. We start with steps of the a posteriori analysis that are common to the `classical' and the new approach. We then illustrate the classical approach with the \emph{standard residual estimator}, and afterwards develop the new approach resulting in a modification of the standard residual estimator, called \emph{modified residual estimator}. Finally, we conclude by adapting the new approach to other techniques of a~posteriori error estimation and boundary conditions. 

% notation
In what follows, the notation may or may not indicate the dependencies of a given quantity. We shall balance readability and the importance of the dependence in the given context. For example, in \eqref{setting-apost:Galerkin}, the discrete solution depends not only on the mesh $\grid$ but also on the data $\Omega$, $\bA$, $c$, and $f$ in problem \eqref{strong-form}. We write $\grid$ explicitly
because of its more prominent role in the a posteriori analysis. Let $\faces:=\faces_\grid$
\index{Meshes!$\faces$, $\faces_\grid$: interior faces} denote the set of all interior $(d-1)$-dimensional faces of $\grid$. The letter $C$ will be used for a generic constant, with possibly different values at each occurrence. If not stated otherwise, it may depend on the shape regularity coefficient $\sigma$ from \eqref{E:shape-regularity}, the dimension $d$, and the polynomial degree $n$ in $\V_\grid$.
  
%
%-----------------------------------------------------------------------------------------
\subsection{Error, residual and localization of residual norm}
\label{S:err-res}
%-----------------------------------------------------------------------------------------
%
% intro
This section starts the a~posteriori analysis by establishing that a suitable norm of the so-called \emph{residual}
\begin{itemize}
\item is equivalent to the error $\| \nabla (u-u_\grid) \|_{L^2(\Omega)}$
%$\| \nabla e_\grid \|_{L^2(\Omega)}$, $e_\grid \definedas u-u_\grid$
and
\item admits a \emph{localization} in the sense that it splits into suitable local contributions depending on accessible quantities, i.e.\ on data $\data=(\bA,c,f)$ and the discrete solution $u_\grid$.
\end{itemize}
We do not consider computability yet --- this important aspect will be addressed in the following sections.

% residual
Replacing in the weak form \eqref{weak-form} the exact solution $u$ by its approximation $u_\grid$, we define the \emph{residual} $R_\grid \in H^{-1}(\Omega)$:
\begin{equation*}
\index{Error Estimators!$R_\grid$: residual in $H^{-1}(\Omega)$}
\langle R_\grid, w \rangle
 =
 \langle f, w \rangle
 -
 \B[u_\grid,w]
\quad
\forall \, w \in H^1_0(\Omega).
\end{equation*}
We thus have a quantity that depends only on data $\data$ and the approximate solution $u_\grid$ and relates to the error function $u-u_\grid$
%$e_\grid$
as follows:
\begin{equation}
\label{err-res-id}
 \langle R_\grid, w \rangle
 =
 \B[u-u_\grid %e_\grid
 ,w] \quad \forall w \in H^1_0(\Omega).
\end{equation}
Continuity and coercivity of the bilinear form $\B$ then provide a quantitative relationship
between error and residual.

\begin{lemma}[error and residual]
\label{L:err-residual}
The error of the approximation $u_\grid$ is equivalent to the residual norm. More precisely, \begin{equation*}
 \frac{1}{\| B \|} \| R_\grid \|_{H^{-1}(\Omega)}
 \leq
 \| \nabla(u-u_\grid) \|_{L^2(\Omega)}
 \leq
 \frac{1}{\alpha} \| R_\grid \|_{H^{-1}(\Omega)},
\end{equation*}
where $\|\B\|\geq\alpha>0$ are, respectively, the continuity and coercivity constants of the bilinear form $\B$.
\end{lemma}

\begin{proof}
The error-residual relationship \eqref{err-res-id} yields the lower bound,
\begin{align*}
 \| R_\grid \|_{H^{-1}(\Omega)}
 &=
 \sup_{w \in H^1_0(\Omega)} \frac{\langle R_\grid, w \rangle}{\| \nabla w \|_{L^2(\Omega)}}
 =
 \sup_{w \in H^1_0(\Omega)} \frac{\B[u-u_\grid, w]}{\| \nabla w \|_{L^2(\Omega)}}
\\
 &\leq
 \| \B \| \| \nabla (u-u_\grid) \|_{L^2(\Omega)},
\end{align*}
while the choice $w = u - u_\grid$ therein gives
\begin{align*}
 \alpha \| \nabla(u-u_\grid) \|_{L^2(\Omega)}^2
 &\leq
 \B[u-u_\grid,u-u_\grid]
 =
 \langle R_\grid, u - u_\grid \rangle
\\
 &\leq
 \| R_\grid \|_{H^{-1}(\Omega)} \| \nabla (u-u_\grid) \|_{L^2(\Omega)}
\end{align*}
and thus the upper bound.
\end{proof}

\begin{remark}[role of forcing vs role of coefficients]
\label{R:role-of-f-vs-A,c}
In addition to the two differences between right-hand side $f$ and coefficients $(\bA,c)$ mentioned in the introduction of this section, a third one implicitly arises in the proof of Lemma~\ref{L:err-residual}: the coefficients defining the bilinear form $\B$ are fixed, while the right-hand side $f$ is replaced by the residual $R_\grid$ in \eqref{err-res-id}, which varies with the mesh $\grid$.
\end{remark}

\begin{remark}[local lower estimate for the error]
\label{R:local-lower-bounds}
The proof of Lemma~\ref{L:err-residual} shows that the lower bound of the error hinges on the continuity of the bilinear form $\B$. Since the evaluation of $\B$ involves only local operators, one might expect that there are also local lower bounds. This however depends on the interplay of the underlying differential operator and the choice of the test space norm. Indeed, in the case of the Poisson problem, i.e.\ $\bA = \vec{I}$, $c=0$, and the test space norm $\| \nabla \cdot \|_{L^2(\Omega)}$, one easily sees that
\begin{equation*}
  \| R_\grid \|_{H^{-1}(\omega)}
  \leq
 \| \nabla(u-u_\grid) \|_{L^2(\omega)}
\end{equation*}
for any subdomain $\omega \subset \Omega$. This local lower bound however does not carry over to the general case with $c\neq0$ as the error function itself is bounded by its gradient only through the global inequality in Lemma~\ref{E:Poincare} (first Poincar\'e inequality). On the other side, endowing the test space $H^1_0(\Omega)$ with the full $H^1$-norm $\|\cdot\|_{H^1(\Omega)}$ yields
\begin{equation*}
 \| R_\grid \|_{(H^1(\omega))^*}
 \leq
 \max \big\{ \alpha_1, \| c \|_{\modif{L^\infty(\Omega)}} \big\} \| u - u_\grid \|_{H^1(\omega)}
\end{equation*}
for any subdomain $\omega \subset \Omega$.

In line with \cite{Axioms:2014}, we shall not invoke local lower bounds to derive convergence and rate optimality for the error of $\AFEM$, although they might appear useful or even crucial in other settings.
\end{remark}

\begin{remark}[constants in error-residual relationship]
\label{R:constants-for-error-res-rel}
Given the norm measuring the error, that is the norm of the trial space, the choice of the test space norm is important for the ensuing constants in the error-residual relationship;
%. Indeed, the bilinear form $\B$ is continuous and coercive for test space norms from a class of %equivalent norms. For our purposes, a perfect one is given by
%\begin{equation*}
% \| w \|
% \definedas
% \sup_{v \in H^1_0(\Omega)} \frac{\B[v,w]}{\| \nabla v \|_{L^2(\Omega)}}
%\end{equation*}
%as it implies that inf-sup and continuity constants satisfy $\beta = 1 = \| \B \|$. However, using this %norm, or a close variant, would yield the following analysis more involved;
see, e.g., \cite[Sects. 4.3 and 4.6]{Verfuerth:13}.
%for analyses in this direction.
To avoid related additional technicalities and difficulties, the test space is endowed with the straightforward norm $\|\nabla\cdot\|_{L^2(\Omega)}$.  
\end{remark}

\medskip Lemma~\ref{L:err-residual} establishes the first goal that was set out at the beginning of this section. We now turn to the second one, that is, we split the residual norm $\| R_\grid \|_{H^{-1}(\Omega)}$ into local contributions. Note that the nature of the dual norm $\|\cdot\|_{H^{-1}(\Omega)}$ makes this task less obvious than for integral norms as in the error $\| \nabla(u - u_\grid) \|_{L^2(\Omega)}$.

We start by recalling that the definition of the Galerkin approximation $u_\grid$ implies that its residual is orthogonal to the discrete trial space $\V_\grid = \mathbb{S}^{n,0}_\mesh \cap H^1_0(\Omega)$:
\begin{equation*}
%\label{orth-of-residual}
 \langle R_\grid, w \rangle = 0 \quad  \forall \,w \in \V_\grid.
\quad
\end{equation*}
Denote by $\vertices$ the set of vertices of $\grid$ and by $\phi_z \in \mathbb{S}^{1,0}_\grid$ the hat function with $\phi_z(y) = \delta_{yz}$ for all vertices $y \in \vertices$. In what follows, the \emph{partial} orthogonality
\begin{equation}
\label{partial-orth-of-residual}
 \langle R_\grid, \phi_z \rangle = 0 \quad \forall \, z \in \vertices % = \mathbb{S}^{n,0}_\mesh \cap H^1_0(\Omega)
\quad
\end{equation}
will be crucial to split the nonlocal norm of the residual into local contributions. The latter ones will be formulated in terms of the supports of the hat functions and, thus, we associate
to each vertex $z \in \vertices$ the following subset and submesh:
\begin{equation}
\label{stars}
\index{Meshes!$\omega_z$: region made of elements sharing the vertex $z$}
\omega_z \definedas  \operatorname{\supp} \phi_z = \bigcup_{T \in \grid_z} T
\quad\text{with}\quad
\index{Meshes!$\grid_z$: star of elements sharing the vertex $z$}
 \grid_z \definedas \{ T \in \grid \mid T \ni z \}.
\end{equation}
These subsets, called \emph{stars},  form a \emph{subdomain covering} of $\Omega$, viz.\ each interior $\mathring{\omega}_z$ is a domain and $\overline{\Omega} = \cup_{z\in\vertices} \; {\omega}_z$. The \emph{overlapping index} $\esssup_{x\in\Omega} \#\{ z \in \vertices \mid \omega_z \ni x \}$ of this covering is bounded by $(d+1)$.

\begin{lemma}[localization of $H^{-1}$-norm]
\label{L:loc-of-H^-1-norm}
Let $\ell \in H^{-1}(\Omega)$ be an arbitrary linear functional on $H^1_0(\Omega)$.
\begin{enumerate}
\item If $\langle \ell, \phi_z \rangle = 0$ for all interior vertices $z \in \vertices \cap \Omega$, then
\begin{equation*}
\index{Constants!$C_\mathrm{loc}$: localization constant}
 \| \ell \|_{H^{-1}(\Omega)}^2
 \leq
 (d+1) C_\mathrm{loc}^2 \sum_{z \in \vertices} \| \ell \|_{H^{-1}(\omega_z)}^2
\end{equation*}
where $C_\mathrm{loc}$ depends only on the shape regularity coefficient $\sigma$ from \eqref{E:shape-regularity} and $d$. \looseness=-1
\item For any subdomain covering  $(\omega_i)_{i\in I}$ of $\Omega$ with finite overlapping index $C_\mathrm{ovrl} \definedas \esssup_{x \in \Omega} \#\{ i \in I \mid \omega_i \ni x \}$, we have
\begin{equation*}
\index{Constants!$C_\mathrm{ovrl}$: overlay constant}
  \sum_{i \in I} \| \ell \|_{H^{-1}(\omega_i)}^2
  \leq
   C_\mathrm{ovrl} \| \ell \|_{H^{-1}(\Omega)}^2.
\end{equation*}
\end{enumerate}
\end{lemma}

\cite[Theorem~3.5]{Blechta.Malek.Vohralik:2020}
%\todo{AV: Problem to get the accent on the `i' in `Vohralik'. AB: cannot make it work, let's leave this to the editors}
generalizes Lemma~\ref{L:loc-of-H^-1-norm} to the $W^{-1}_p$-norm, $1<p<\infty$. Lemma~\ref{L:alt-loc-of-H^-1-norm} below provides an alternative localization with different local norms.

\begin{proof}
\fbox{\scriptsize 1} We start by showing statement (i). Thanks to the orthogonality of $\ell$, we may write
\begin{equation*}
 \langle \ell, w \rangle
 =
 \langle \ell, w - \sum_{z\in\vertices} c_z \phi_z \rangle,
\end{equation*}
where any $c_z$ is an arbitrary constant if $z \in \vertices \cap \Omega$ is an interior vertex and $0$ if  $z \in \vertices \cap \partial\Omega$ is a boundary vertex. Using the partition of unity $\sum_{z\in\vertices}\phi_z = 1$ on $\Omega$, we split the new test function
\begin{equation*}
 w - \sum_{z\in\vertices} c_z \phi_z
 =
 \sum_{z\in\vertices} (w-c_z)\phi_z,
\end{equation*}
into local contributions $(w-c_z)\phi_z \in H^1_0(\omega_z)$, $z \in \vertices$. The constant $c_z$ allows us to counter the gradient generated by the cut-off with $\phi_z$. Indeed, the product rule, $0\leq\phi_z\leq 1$, and $|\nabla \phi_z| \leq C(d) \sigma h_T^{-1}$ on an element $T \in \grid$ lead to
\begin{align*}
 \| \nabla \big( (w-c_z) &\phi_z \big) \|_{L^2(\omega_z)}
 \leq
 \| \phi_z \nabla (w-c_z) \|_{L^2(\omega_z)}
 +
 \| (w-c_z) \nabla \phi_z \|_{L^2(\omega_z)}
\\
 &\leq
 \| \phi_z \|_{L^\infty(\omega_z)} \| \nabla w\|_{L^2(\omega_z)}
 +
 \| \nabla \phi_z \|_{L^\infty(\omega_z)} \| w-c_z \|_{L^2(\omega_z)}
\\
 &\leq
 \| \nabla w\|_{L^2(\omega_z)} + C(d) \sigma\left( \max_{T\subset\omega_z} h_T^{-1} \right)
  \| w-c_z \|_{L^2(\omega_z)}.
 \end{align*}
If we choose $c_z = \fint_{\omega_z} w$ for interior vertices $z \in \vertices \cap \Omega$, then Lemma~\ref{L:poincare-friedrichs} (second Poincar\'e inequality) on reference stars implies
\begin{equation*}
 \| w-c_z \|_{L^2(\omega_z)}
 \Cleq
 \diam \omega_z \| \nabla w \|_{L^2(\omega_z)}.
\end{equation*}
The same inequality follows for boundary vertices $z \in \vertices \cap \partial\Omega$ thanks to the fact that $w$ vanishes at least on one face of $\partial \omega_z \cap \partial\Omega$. Combing this with $\diam \omega_z \Cleq h_T$ for $T \subset \omega_z$, we thus obtain, for all local contributions, the stability bound
\begin{equation}
\label{def-Cloc}
 \| \nabla \big((w-c_z)\phi_z\big) \|_{L^2(\omega_z)}
 \leq
 C_\mathrm{loc} \| \nabla w \|_{L^2(\omega_z)},
\end{equation}
where the constant $C_\mathrm{loc}$ depends only on $d$ and $\sigma$. Hence,
\begin{equation*}
 \langle \ell, w \rangle
 =
 \langle \ell, w - \sum_{z\in\vertices} c_z \phi_z \rangle
 =
 \sum_{z\in\vertices} \langle \ell, (w-c_z)\phi_z \rangle
\end{equation*}
gives
\begin{align*}
 |\langle \ell, w \rangle|
% &\leq
% \sum_{z\in\vertices} \| \ell \|_{H^{-1}(\omega_z)} \| \nabla \big((w-c_z)\phi_z\big) \|_{L^2(\omega_z)}
%\\
 &\leq
 C_\mathrm{loc} \sum_{z\in\vertices} \| \ell \|_{H^{-1}(\omega_z)} \| \nabla w \|_{L^2(\omega_z)}
\\
 &\leq
 \sqrt{d+1} C_\mathrm{loc} \left( \sum_{z\in\vertices} \| \ell \|_{H^{-1}(\omega_z)}^2 \right)^{1/2}
   \| \nabla w \|_{L^2(\Omega)}
\end{align*}
and (i) is proven.

\smallskip\fbox{\scriptsize 2}
We verify statement (ii). For each index $i \in I$, define $v_i \in H^1_0(\omega_i) \subset H^1_0(\Omega)$ by
\begin{equation*}
 \int_{\omega_i} \nabla v_i \cdot \nabla w
 =
 \langle \ell, w \rangle, \qquad  \forall w \in H^1_0(\omega_i).
\quad
\end{equation*}
We obtain $ \langle \ell, v_i \rangle = \| \nabla v_i \|_{L^2(\omega_i)}^2 = \| \ell \|_{H^{-1}(\omega_z)}^2$ by arguments similar to the ones in the proof of Lemma \ref{L:err-residual} (error and residual). The sum $v \definedas \sum_{i \in I} v_i$ is in $H^1_0(\Omega)$ with
\begin{align*}
 \| \nabla v \|_{L^2(\Omega)}^2
 &\leq
 \int_\Omega \Big| \sum_{i \in I_x} \nabla v_i(x) \Big|^2 \, dx
 \leq
 \int_\Omega \#I_x \sum_{i \in I_x} |\nabla v_i(x)|^2 \, dx
\\
 &\leq
 C_\mathrm{ovrl} \sum_{i \in I} \|\nabla v_i \|_{L^2(\omega_i)}^2
 =
 C_\mathrm{ovrl} \sum_{i \in I} \| \ell \|_{H^{-1}(\omega_z)}^2,
\end{align*}
where we denote the set of active indices in $x \in \Omega$ by $ I_x \definedas \{ i \in I \mid \omega_i \ni x \}$. Inserting this in 
\begin{equation*}
 \sum_{i \in I}  \| \ell \|_{H^{-1}(\omega_z)}^2
 =
 \sum_{i \in I}   \langle \ell, v_i \rangle
 =
 \langle \ell, v \rangle
 \leq
 \| \ell \|_{H^{-1}(\omega_z)} \| \nabla v \|_{L^2(\Omega)}
\end{equation*}
establishes the desired inequality.
\end{proof}

Thanks to the partial orthogonality \eqref{partial-orth-of-residual} and the properties of the star covering $\omega_z$, $z \in \vertices$, we readily obtain the following statement.
\begin{corollary}[star localization of residual norm]
\label{C:residual-loc}
The $H^{-1}$-norm of the residual can be split into local contributions on stars:
\begin{equation*}
 \frac{1}{d+1} \sum_{z \in \vertices} \| R_\grid \|_{H^{-1}(\omega_z)}^2
 \leq
 \| R_\grid \|_{H^{-1}(\Omega)}^2
 \leq
 (d+1) C_\mathrm{loc} \sum_{z \in \vertices} \| R_\grid \|_{H^{-1}(\omega_z)}^2,
\end{equation*}
where $C_\mathrm{loc}$ depends only on $d$ and the shape regularity coefficient $\sigma$.
\end{corollary}

The upper bound of the global residual norm in Corollary~\ref{C:residual-loc} employs the stars  $\omega_z$, $z\in\vertices$, as local domains. The next remark assesses this choice by discussing conceivable alternatives in terms of elements and domains of the type
\begin{equation}
\label{pairs}
\index{Meshes!$\omega_F$: region of elements containing the face $F$}
 \omega_F \definedas \bigcup_{T \in \grid_F} T
\quad\text{with}\quad
 \grid_F \definedas \{ T \in \grid \mid T \supset F \}, 
\end{equation}
where $F \in \faces$ is an interior face of $\grid$.

\begin{remark}[star localization is minimal for $d\geq2$]
\label{R:star-loc-is-minimal}
The use of stars in the localization of the global residual norm is a sort of minimal choice, except for the special case $d=1$ where elements can be used:
\begin{itemize}
\item If $d=1$, point values are defined for functions in $H^1(\Omega)$. This allows an upper bound with elements instead of stars as local domains. In fact, choosing $c_z = w(z)$ for all interval endpoints, the function $\sum_{z\in\vertices} c_z \phi_z$ amounts to the Lagrange interpolant $I_\grid w \in {\mathbb S}^{1,0}_\grid \cap H^1_0(\Omega)$, and we have $(w-I_\grid w)|_I \in H^1_0(I)$ with $\| \nabla (w-I_\grid w) \|_{L^2(I)} \Cleq \| \nabla w \|_{L^2(I)}$ for any interval $I$ of the mesh $\grid$. Arguing as in the proof of Lemma~\ref{L:loc-of-H^-1-norm} (i) then gives
\begin{equation*}
 \| R \|_{H^{-1}(\Omega)}^2
 \Cleq
 \sum_{I \in \grid} \| R \|_{H^{-1}(I)}^2.
\end{equation*}
\item An upper bound where the stars are replaced by elements $T\in\grid$ cannot hold in general because it does not account for face-supported residual contributions. For example, consider our setting with
\begin{equation*}
\begin{gathered}
 d\geq2, \quad \bA=\vec{I}, \quad c=0,
\quad
 \langle f, w \rangle = \int_F q w, \; w \in H^1_0(\Omega) 
\\
 \text{where $F \in \faces$ and $q\neq0$ is $L^2$-orthogonal to $\P_n(F)$.}
\end{gathered}
\end{equation*}
 Then we have $u \neq 0 = u_\grid$ and therefore $\| R_\grid \|_{H^{-1}(\Omega)} > 0$ but $\| R_\grid \|_{H^{-1}(T)} = 0$ for any $T \in \grid$.
\item An upper bound with pairs $\omega_F$, $F \in \faces$, instead of stars cannot hold in general. This can be shown by considering our setting with
\begin{equation*}
\begin{gathered}
 d=2, \quad \bA=\vec{I}, \quad c=0,
\quad
 \langle f_\eps, w \rangle = \frac{1}{\pi \eps^2} \int_{B_\eps(z)} w, \; w \in H^1_0(\Omega)
\\
 \text{where $z \in \vertices$ is a vertex of a suitable triangulation $\grid$}
\end{gathered}
\end{equation*}
for $\eps \searrow 0$. The limiting right-hand side is the Dirac measure in $z$, which, formally, is not be seen by any $\| \cdot \|_{H^{-1}(\omega_F)}$, $F\in\faces$. We thus have $\| R_\grid \|_{H^{-1}(\Omega)} \to \infty$ but $\sum_{F\in\faces} \| R_\grid \|_{H^{-1}(\omega_F)} \Cleq 1$; cf.\ \cite{Tantardini.Veeser.Verfuerth}.
\end{itemize}
\end{remark}

The bisection method for mesh refinement is element-oriented. It is therefore advantageous to dispose of an element-indexed reformulation of the localization  in Corollary~\ref{C:residual-loc}. For that purpose, we recall the notion of \emph{patches}
\begin{equation}
\index{Meshes!$\omega_T$, $\omega_\grid(T)$: region of elements intersecting $T$}
\omega_T := \bigcup_{\substack{T' \in \grid \\ T' \cap T \neq \emptyset}} T'
\end{equation}
and may use the following equivalence.

\begin{lemma}[localization re-indexing]\label{L:local-re-indexing}
For any functional $\ell \in H^{-1}(\Omega)$, we have
\begin{equation*}
  \sum_{z \in \vertices} \| \ell \|_{H^{-1}(\omega_z)}^2
  \approx
  \sum_{T \in \grid} \| \ell \|_{H^{-1}(\omega_T)}^2 \,,
\end{equation*}
where the hidden constants depend on $d$ and the shape regularity coefficient $\sigma$.
\end{lemma}

Note that, in contrast to the localization itself, its re-indexing does not require any orthogonality like $\langle \ell, \phi_z \rangle = 0$ for all $z \in \vertices \cap \Omega$.

\begin{proof}
For any vertex $z \in \vertices$, there is an element $T \in \grid$ containing $z$. Then the inclusion $\omega_z \subset \omega_T$ yields the inequality $\| \ell \|_{H^{-1}(\omega_z)} \leq \| \ell \|_{H^{-1}(\omega_T)}$. Hence,
\begin{equation*}
 \sum_{z\in\vertices} \| \ell \|_{H^{-1}(\omega_z)}^2
 \leq
 \sum_{T\in\grid} \| \ell \|_{H^{-1}(\omega_T)}^2.
\end{equation*}
To show the converse inequality, let $T\in\grid$ be any element and $w \in H^1_0(\omega_T)$. Given any vertex $z \in \vertices \cap \omega_T$, Lemma~\ref{L:Poincare} (first Poincar\'e inequality) on $\omega_T$ implies the stability bound
\begin{align*}
 \| \nabla (w \phi_z) \|_{L^2(\omega_z \cap \omega_T)}
 &\leq
 \| \phi_z \nabla w \|_{L^2(\omega_z \cap \omega_T)}
 +
 \| w \nabla \phi_z \|_{L^2(\omega_z \cap \omega_T)}
\\
 &\leq
 \| \nabla w \|_{L^2(\omega_z \cap \omega_T)}
 +
 C(d,\sigma) \max_{T \subset \omega_z \cap \omega_T} h_T ^{-1} \, \| w \|_{L^2(\omega_T)}
\\
 &\Cleq
 \| \nabla w \|_{L^2(\omega_T)}.
\end{align*}
We thus derive 
\begin{align*}
 \langle \ell, w \rangle
 &=
 \sum_{z \in \vertices \cap \omega_T} \langle \ell, w\phi_z \rangle
 \leq
 \sum_{z \in \vertices \cap \omega_T}
  \| \ell \|_{H^{-1}(\omega_z \cap \omega_T)} \| \nabla (w\phi_z) \|_{L^2(\omega_z \cap \omega_T)}
\\
 &\Cleq
  \left( \sum_{z \in \vertices \cap \omega_T} \| \ell \|_{H^{-1}(\omega_z)} \right)
  \| \nabla w \|_{L^2(\omega_T)}
\end{align*}
and, since $\#(\vertices \cap \omega_T)$ is bounded in terms of the shape regularity coefficient $\sigma$,
\begin{equation*}
 \| \ell \|_{H^{-1}(\omega_T)}^2
 \Cleq
 \sum_{z \in \vertices \cap \omega_T} \| \ell \|_{H^{-1}(\omega_z)}^2.
\end{equation*}
Summing over $T\in\mesh$ and taking into account that $\#\{ T \in \mesh \mid \omega_T \ni z \}$ is again bounded in terms of $\sigma$, concludes the proof.
\end{proof}

\subsection{Standard residual estimator and its flaws}
\label{S:std-res-est}
%
%
% intro
Exploiting the results of Sect.\ \ref{S:err-res}, we derive an \emph{a~posteriori upper bound} of the error in terms of the standard residual estimator and discuss its flawed sharpness. This discussion will serve as the starting point for an improved a~posteriori analysis in the following sections.

% definition of standard residual estimator
\smallskip The \emph{standard residual estimator} needs the additional regularity
\begin{equation}
\label{extra-reg-for-std-res-est}
 f \in L^2(\Omega)
\quad\text{and}\quad
 \bA \in W^{1}_{\infty}(\Omega; \mathbb R^{d\times d})
\end{equation}
for the data in our model problem \eqref{strong-form}. Given the Galerkin approximation $u_\grid$ from \eqref{setting-apost:Galerkin}, it may be defined as follows, cf., e.g., \cite{Verfuerth:13}:
\begin{subequations}
\label{std-res-est}
\begin{equation}
\label{std-res-est;tot}
\index{Error Estimators!$\est_\grid^{\text{std}}(u_\grid,\data)$: standard residual estimator}
 \est_\grid^\text{std}
 \definedas
 \est_\grid^\text{std}(u_\grid,\data)
 \definedas
 \left(
  \sum_{T\in\grid} \est_\grid^\text{std}(u_\grid,\data,T)^2
\right)^{\frac{1}{2}}
\end{equation}
with the \emph{local indicators}
\begin{equation}
\label{loc-ind-of-std-res-est}
\index{Error Estimators!$\est_\grid^{\text{std}}(u_\grid,\data,T)$: standard local indicators}
%\begin{aligned}
 \est_\grid^\text{std}(u_\grid, \data, T)^2
 \definedas
% &\definedas
 h_T \| j(u_\grid) \|_{L^2(\partial T \setminus \partial\Omega)}^2
 +
 h_T^2 \| r(u_\grid) \|_{L^2(T)}^2,
% h_T \| [A\nabla u_\grid]\cdot\vec{n} \|_{L^2(\partial T \setminus \partial\Omega)}^2
%\\
% &\qquad +  h_T^2  \| f - c u_\grid + \div(\bA\nabla u_\grid) \|_{L^2(T)}^2.
%\end{aligned}
\end{equation}
\end{subequations}
where
\begin{itemize}
\item the \emph{scaling factor} $h_T = |T|^{1/d}$ measures the size of the element $T \in \grid$,
\item $j(v)=j_\grid(v)$ is the \emph{jump residual} given face-wise for $v \in \mathbb V_\grid$ by 
\begin{align*}
\index{Error Estimators!$j(u_\grid)$: jump residual}
 j(v)|_{F}
 &\definedas
 \big( \jump{\bA\nabla v}\cdot \vec{n}_{T_1} \big)|_{F}
 \definedas
 \big( (\bA \nabla v)|_{T_1} - (\bA \nabla v)|_{T_2} \big) \cdot \vec{n}_{T_1}
\\
 &=
 (\bA \nabla v)|_{T_1} \cdot \vec{n}_{T_1}  + (\bA \nabla v)|_{T_2} \cdot \vec{n}_{T_2}
\end{align*}
where $F \in \faces$, $T_1, T_2 \in \grid$ are such that $F = T_1 \cap T_2$, $\vec{n}_{T_i}$ denotes the outer normal of $\partial T_i$, $i=1,2$, and
\item 
\index{Error Estimators!$r_\grid(u_\grid)$, $r(u_\grid)$: element residual}
$r(v)=r_\grid(v)$ is the \emph{element residual}, a function given for $v \in \mathbb V_\grid$ by 
\[
r(v)|_{T} \definedas \big(f - c v + \div(\bA\nabla v) \big)|_{T}
\]
on any element $T\in\grid$.
\end{itemize}
Note that the definition itself already uses the extra regularity \eqref{extra-reg-for-std-res-est}. For notational simplicity, we shall write $j$ and $r$ instead of $j(u_\grid)$ and $r(u_\grid)$ for the rest of this section.
Also, for any interior face $F \in \faces$, we have $F = T_1 \cap T_2$ with $T_1,T_2 \in \grid$. If $\vec{n}_{T_i}$ denotes the outer normal of $\partial T_i$, $i=1,2$, we set $\vec{n}_F=\vec{n}_{T_1}$\index{Meshes!$\vec{n}_F$: normal to the face $F$}. This particular choice of $\bn_F$ is irrelevant as it does not affect the following definition of normal jump of any vector-valued field $\bg$ with well defined trace on $F$
$$
\index{Meshes!$\jump{\bg} \cdot \vec{n}_F$: normal jump across $F$}
\jump{\bg }\cdot \vec{n}_{F} :=  \bg|_{T_1} \cdot \vec{n}_{T_1} + \bg|_{T_2}  \cdot \vec{n}_{T_2}.
$$

\begin{theorem}[upper bound with standard residual estimator]
\label{T:up-bd-with-std-res-est}
Suppose the additional regularity \eqref{extra-reg-for-std-res-est} holds. Then the error is bounded by the standard residual estimator:
\begin{equation*}
 \| \nabla (u-u_\grid) \|_{L^2(\Omega)}
 \Cleq
 \est_\grid^\text{std},
\end{equation*}
where the hidden constant depends on the coefficients $(\bA,c)$, the shape regularity coefficient $\sigma$, and $d$.
\end{theorem}

\begin{proof}
As Lemma~\ref{L:err-residual} (error and residual) and Corollary~\ref{C:residual-loc} (star localization of residual norm) imply
\begin{equation*}
 \| \nabla(u-u_\grid) \|_{L^2(\Omega)}^2
 \Cleq
 \| R_\grid \|_{H^{-1}(\Omega)}^2
 \Cleq
 \sum_{z \in \vertices} \| R_\grid \|_{H^{-1}(\omega_z)}^2
\end{equation*}
and $\#\{z \in \vertices \mid \omega_z \supset T \}=d+1$, it suffices to establish
\begin{equation}
\label{loc-up-bd-with-std-res-est}
 \| R_\grid \|_{H^{-1}(\omega_z)}^2
 \Cleq
 \sum_{T \subset \omega_z} \est_\grid^\text{std}(u_\grid,f,T)^2
\end{equation}
for any vertex $z \in \vertices$. To this end, let $w \in H^1_0(\omega_z)$. The extra regularity $f \in L^2(\Omega)$ and $\bA \in W^{1}_{\infty}(\Omega;\mathbb R^{d\times d})$ allows for piecewise integration by parts, which leads to the following $L^2$-representation of the residual:
\begin{equation*}
 \langle R_\mesh, w \rangle
 =
 \sum_{F \ni z} \int_F jw %[A\nabla u_\grid] \cdot \vec{n} w
 +
 \sum_{T \ni z} \int_T rw. %\big(f - cu_\grid + \operatorname{div} (\bA \nabla u_\grid) \big) w.
\end{equation*}
In order to bound the right-hand side suitably, we use
the scaled trace theorem
\begin{equation}
\label{scaled-trace-thm-for_simplex-faces}
 \| w \|_{L^2(F)}^2
 \leq
 \frac{|F|}{|T|} \| w \|_{L^2(T)}^2
 +
 \frac{2}{d} \frac{|F|\diam T}{|T|} \| w \|_{L^2(T)} \| \nabla w \|_{L^2(T)}
\end{equation}
for any face $F \subset \partial T$, see e.g.\ \cite[Corollary~4.5]{Veeser.Verfuerth:09}, the inequality
\begin{equation*}
 \| w \|_{L^2(\omega_z)}
 \leq
 \diam \omega_z \| \nabla w \|_{L^2(\omega_z)}
\end{equation*}
from Lemma~\ref{L:Poincare} (first Poincar\'e inequality), and the two geometric relationships
\begin{gather*}
 \diam \omega_z \Cleq h_T \text{ whenever } T \subset \omega_z,
\quad
 |F| \diam T \Cleq |T| \text{ for } F \subset \partial T. 
\end{gather*}
We thus obtain
\begin{equation*}
%\label{flawed-triangle-ineqs}    
 | \langle R_\grid, w \rangle |
 \Cleq
 \left(
 \sum_{T \ni z}  h_T \| r \|_{L^2(T)}
  +
  h_T^{1/2} \sum_{F \ni z, \ F \subset T} \| j \|_{L^2(F)}
 \right) \| \nabla w \|_{L^2(\omega_z)}.
\end{equation*}
As the number of faces and elements in the star $\omega_z$ is bounded in terms of the shape regularity coefficient $\sigma$, we arrive at the desired bound \eqref{loc-up-bd-with-std-res-est} and the proof is finished.
\end{proof}

\begin{remark}[alternative derivation of upper bound]
\label{R:alt-derivation-of-std-up-bd}
The upper bound for the standard residual estimator in Theorem \ref{T:up-bd-with-std-res-est} is often derived with a suitable interpolation operator, by-passing the localization of the $H^{-1}$-norm in Lemma~\ref{L:loc-of-H^-1-norm}. That approach is useful for the proof of Theorem~\ref{T:ubd-corr} below and is presented therein. Here we opted for using the localization of the $H^{-1}$-norm in order to facilitate the comparison with the following subsections. The approach at hand is also convenient to keep the ensuing constants small; cf.\ \cite{Veeser.Verfuerth:09}. 
\end{remark}

% discussion of sharpness of upper bound
An important question is the sharpness of the upper bound in Theorem \ref{T:up-bd-with-std-res-est}. The so-called \emph{a~posteriori lower bounds} provide some answer by trying to bound the estimator in terms of the error. For many estimators, there arises however additional terms of \emph{oscillatory} nature. The following remark justifies the presence of such terms for the case at hand.

\begin{remark}[nonasymptotic overestimation]
\label{R:nonasym-overest}
The lower bound
\begin{equation*}
 \est_\grid^\text{std} \Cleq \| \nabla(u-u_\grid) \|_{L^2(\Omega)},
\end{equation*}
which would imply equivalence of error and estimator, \emph{cannot hold} in general for the following reason.

Fix a mesh $\grid$ and a functional $f \in H^{-1}(\Omega) \setminus L^2(\Omega)$ and consider a sequence $(f_n)_n$ of functions in $L^2(\Omega)$ with $\lim_{n\to\infty} \| f - f_n\|_{H^{-1}(\Omega)} = 0$. Then the sequences $(u_n)_n$ and $(u_{\grid,n})_n$ of exact and Galerkin solutions on a fixed mesh $\grid$ remain bounded. The error sequence $\big( \| \nabla (u_n - u_{\grid,n}) \|_{L^2(\Omega)}\big)_n$ is therefore also bounded, while the standard residual estimator $\est_\grid^\text{std}(u_{\grid,n},f_n) \to \infty$ becomes unbounded. Note that, in the special case $f = -\operatorname{div}(\bA\nabla v) + cv$ with $v \in \V_\grid$, we even have for the error $\lim_{n\to\infty} \| \nabla (u_n - u_{\grid,n}) \|_{L^2(\Omega)} = 0$.

In other words: in certain cases, the standard residual estimator bounds almost $0$ by almost $\infty$ and a lower bound has to involve an additional term that cannot be bounded by the error in general. 
\end{remark}

% oscilllations
We shall define these additional terms with the help of the following local best approximations. Let $K$ be an element or face of $\grid$ and $m \in \N_0$ a polynomial degree. Given $v \in L^2(K)$, denote by $\Pi_Kv \definedas \Pi_K^m v$ \index{Operators!$\Pi_K$, $\Pi_K^m$: $L^2$ projection onto $\P_m(K)$} the best approximation in $\P_m(K)$ with respect to the norm $\|\cdot\|_{L^2(K)}$. It is convenient to allow also for $m=-1$ with $\P_{-1}(K) = \{0\}$ and $\Pi_K^{(-1)} v = 0$. Writing $\data=(\bA,c,f)$ for the data in problem \eqref{strong-form}, the $(m_1,m_2)$-oscillation for the standard residual estimator is then given by
\begin{subequations}
\label{std-osc}
\begin{equation}
\index{Error Estimators!$\osc_\grid^\text{std}(u_\grid,\data)$: standard oscillation}
 \osc_\grid^\text{std}(u_\grid,\data)^2
 \definedas
% \left(
  \sum_{T \in \grid} \osc_\grid^\text{std}(u_\grid,\data,T)^2,
% \right)^{1/2}
\end{equation}
with the local indicators
\begin{equation}
\label{loc-ind-of-std-osc}
\index{Error Estimators!$\osc_\grid^\text{std}(u_\grid,\data,T)$: standard local oscillation}
 \osc_\grid^\text{std}(u_\grid,\data,T)^2
 \definedas
 h_T^2 \, \| r - \Pi_T^{m_2} r \|_{L^2(T)}^2
 +
 h_T \!\!\!\! \sum_{F\subset \partial T \setminus \partial\Omega}
  \| j - \Pi_F^{m_1} j \|_{L^2(F)}^2.
\end{equation}
\end{subequations}

\begin{proposition}[partial lower bound]
\label{P:lw-bd-with-std-res-est}
If $f \in L^2(\Omega)$ and $\bA \in W^{1}_{\infty}(\Omega;\mathbb R^{d\times d})$, the standard residual estimator is bounded by error and oscillation:
\begin{equation*}
 \est_\grid^\text{std}
 \Cleq
 \| \nabla (u-u_\grid) \|_{L^2(\Omega)}
 +
 \osc_\grid^\text{std}(u_\grid,\data),
\end{equation*}
where the hidden constant depends on $d$, the coefficients $\bA$ and $c$, the shape regularity coefficient $\sigma$ as well as the oscillation degrees $(m_1,m_2)$.
\end{proposition}

\begin{proof}
In light of Lemma \ref{L:err-residual} (error and residual) and Corollary \ref{C:residual-loc} (localization of residual norm), we may establish the claimed bound by bounding each indicator with a corresponding local residual norm. To this end, we shall consider here only the case of %the Poisson problem with $A=I$ and $c=0$, linear finite elements, viz.\ $n=1$, and and
the oscillation degrees $(m_1,m_2)=(0,0)$. The general case can be verified along the same lines with additional technicalities and is treated in the proof of Lemma~\ref{L:local-stability-of-P} below in a slightly different context.

\smallskip\fbox{\scriptsize 1} We start by bounding an arbitrary element residual $h_T \| r \|_{L^2(T)}$, $T \in \grid$, in terms of some local residual norm. To this end, we may try to invert the following consequence of Lemma~\ref{L:Poincare} (first Poincar\'e inequality):
\begin{equation*}
 \| R_\grid \|_{H^{-1}(T)}
 =
 \sup_{w \in H^1_0(T)} \frac{|\langle R_\grid, w \rangle|}{\| \nabla w \|_{L^2(T)}}
 =
 \sup_{w \in H^1_0(T)} \frac{\int_T rw}{\| \nabla w \|_{L^2(T)}}
 \Cleq
 h_T \| r \|_{L^2(T)},
\end{equation*}
the residual norm of which avoids involving the jump residual. We thus actually ask for an equivalence of two different smoothness norms.  Such an equivalence can hold only for special $r$, e.g., from a finite-dimensional space. Furthermore, writing $\|r\|_{L^2(T)}^2 = \int_T r (r\chi_T)$ suggests the choice $w=r\chi_T$, which is however not admissible for the residual $R_\grid$ as both $r$ and the characteristic function $\chi_T$ do not belong to $H^1_0(\Omega)$. We shall overcome these issues by replacing $r$ with its mean value $\Pi_T^0 r$ and $\chi_T$ with the element bubble
\begin{equation}
\label{elm-bubble}
\phi_T \definedas (d+1)^{(d+1)} \prod_{z\in\vertices \cap T} \phi_z.
\end{equation}
Thanks to $\int_T \phi_T = C_d |T|$ and the inverse estimate $\| \nabla w \|_{L^2(T)} \Cleq h_T^{-1} \| w \|_{L^2(T)}$ for $w = (\Pi_T^0 r) \phi_T \in H^1_0(T) \cap \P_{d+1}(T)$, we derive 
\begin{align*}
 \| \Pi_T^0 r \|_{L^2(T)}
 &\Cleq
 \int_T (\Pi_T^0 r) w
 \leq
 \| (\Pi_T^0 r) \chi_T \|_{H^{-1}(T)} \| \nabla w \|_{L^2(T)}
\\
 &\Cleq
 h_T^{-1} \| (\Pi_T^0 r) \chi_T \|_{H^{-1}(T)} \| w \|_{L^2(T)}
 \leq
 h_T^{-1} \| (\Pi_T^0 r) \chi_T \|_{H^{-1}(T)} \| \Pi_T^0r \|_{L^2(T)},
\end{align*}
whence
\begin{equation}
\label{lw-bd-with-discretized-elm-residual}
 h_T \| \Pi_T^0 r \|_{L^2(T)}
 \Cleq
 \| (\Pi_T^0 r) \chi_T \|_{H^{-1}(T)}.
\end{equation}
This implies the desired partial lower bound for the element residual by a perturbation argument and the inequality $\| r - \Pi_T^0 r  \|_{H^{-1}(T)} \Cleq h_T \| r - \Pi_T^0 r \|_{L^2(T)}$, which follows from another application of Lemma~\ref{L:Poincare} (first Poincar\'e inequality):
\begin{equation}
\label{partial-low-bd-for-elm-res}
\begin{aligned}
 h_T \| r \|_{L^2(T)}
 &\leq
 h_T \| \Pi_T^0 r \|_{L^2(T)} + h_T \| r - \Pi_T^0 r \|_{L^2(T)}
\\
 &\Cleq
 \| (\Pi_T^0 r) \chi_T \|_{H^{-1}(T)}  + h_T \| r - \Pi_T^0 r \|_{L^2(T)}
\\
 &\lesssim
 \| r \chi_T \|_{H^{-1}(T)}  + h_T \| r - \Pi_T^0 r \|_{L^2(T)}
\\
 &=
 \| R_\grid \|_{H^{-1}(T)}  + h_T \| r - \Pi_T^0 r \|_{L^2(T)}.
\end{aligned}
\end{equation}

\smallskip\fbox{\scriptsize 2} We bound an arbitrary jump residual $\| j \|_{L^2(F)}$, $F \in \faces$, in a similar manner. Note that here an interference of the element residual is unavoidable because the support of nontrivial test functions have nonempty interior. We thus may try to invert %and we may partially invert
\begin{equation*}
 \| R_\grid \|_{H^{-1}(\omega_F)}
 =
 \sup_{w \in H^1_0(\omega_F)} \frac{ \int_F jw + \int_{\omega_F} rw}{\| \nabla w \|_{L^2(\omega_F)}}
 \Cleq
 h_F^{\frac{1}{2}} \| j \|_{L^2(F)} + \sum_{T \subset \omega_F} h_T \| r \|_{{L^2(T)}},
\end{equation*}
where also the scaled trace theorem \eqref{scaled-trace-thm-for_simplex-faces} is used. To this end, we write  $\delta_F$ for the Dirac measure of the face $F$,
\begin{equation}
\label{face-bubble}
 \phi_F 
 \definedas
 d^d \prod_{z \in \vertices \cap F} \phi_z
\end{equation}
for the face bubble of $F$ and choose the test function $w = (\Pi_F^0j)\phi_F \in H^1_0(\omega_F)$. Using in addition $\| w \|_{L^2(\omega_F)} \Cleq h_F^{\frac{1}{2}} \| w \|_{L^2(F)}$ and \eqref{lw-bd-with-discretized-elm-residual}, we deduce
\begin{align*}
 \| \Pi_F^0 j \|_{L^2(F)}^2
 &\Cleq
 \int_F (\Pi_F^0j) w
 +
 \sum_{T\subset\omega_F} \int_T (\Pi_T^0 r)w
 -
 \sum_{T\subset\omega_F} \int_T(\Pi_T^0 r)w
\\
 &\leq
 \left\|
  (\Pi_F^0j) \delta_F + \sum_{T\subset\omega_F} (\Pi_T^0 r)\chi_T
 \right\|_{H^{-1}(\omega_F)}
  \| \nabla w \|_{L^2(\omega_F)}
\\ 
 & \qquad +
 \sum_{T \subset \omega_F} \| \Pi_F^0 r \|_{L^2(T)} \| w \|_{L^2(T)}
\\
 &\Cleq
  \left\|
  (\Pi_F^0j) \delta_F + \sum_{T\subset\omega_F} (\Pi_T^0 r)\chi_T
 \right\|_{H^{-1}(\omega_F)}
 h_F^{-\frac{1}{2}} \| \Pi_F^0 j \|_{L^2(F)},
\end{align*}
whence
\begin{equation}
\label{lw-bd-with-discrete-jump-res}
 h_F^{\frac{1}{2}} \| \Pi_F^0 j \|_{L^2(F)}^2
 \Cleq
 \left\|
  (\Pi_F^0j) \delta_F + \sum_{T\subset\omega_F} (\Pi_T^0 r)\chi_T
 \right\|_{H^{-1}(\omega_F)}.
\end{equation}
Passing to the proper jump residual $j$, we arrive at the partial lower bound for the jump residual:
\begin{equation}
\label{partial-low-bd-for-jump-res}
\begin{aligned}
  h_F^{\frac{1}{2}} \| \Pi_F^0 j \|_{L^2(F)}
 &\Cleq
 \| R_\grid \|_{H^{-1}(\omega_F)}
\\
 & \quad +
 h_F^{\frac{1}{2}}\| j - \Pi_F^0 j \|_{L^2(F)}
 +
 \sum_{T \subset \omega_F} h_T \| r - \Pi_T^0 r \|_{L^2(T)}.
\end{aligned}
\end{equation}

\smallskip\fbox{\scriptsize 3} We square the bounds \eqref{partial-low-bd-for-elm-res} and \eqref{partial-low-bd-for-jump-res} from the previous steps and sum them, respectively, over all elements and faces to conclude the claimed partial lower bound with the help of Lemma~\ref{L:loc-of-H^-1-norm}(ii) (localization of $H^{-1}-$norm) and Lemma~\ref{L:err-residual} (error and residual).  
\end{proof}

The significance of Proposition~\ref{P:lw-bd-with-std-res-est} (partial lower bound) strongly depends on the choice of the polynomial degrees $(m_1,m_2)$ in the oscillation from \eqref{std-osc}. The following two remarks address this important aspect.

\begin{remark}[oscillation degrees - asymptotics]
\label{R:choosing-osc-degs;asymptotics}
It is desirable that, under refinement, the oscillation in Proposition~\ref{P:lw-bd-with-std-res-est} (partial lower bound) converges to $0$ at least as fast as the error.
The maximal convergence order of the error under uniform refinement is $\| \nabla (u-u_\grid) \|_{L^2(\Omega)} = O(h^{n})$ as $h\to0$.  In view of the scaling factors and derivative orders appearing in jump and element residual, we are thus led to require
\begin{equation*}
 m_1 \geq n-1
\quad\text{and}\quad
 m_2 \geq n-2.
\end{equation*}
One might hope that strict inequalities lead to higher order. Note however that, since $\osc^\text{std}$ involves in general both discrete solution $u_\grid$ and data $\data=(\bA,c,f)$, this will not be guaranteed without additional assumptions. Furthermore, increasing $m_1$ and $m_2$ entails bigger hidden constants in the lower bounds \eqref{lw-bd-with-discretized-elm-residual} and \eqref{lw-bd-with-discrete-jump-res}, as these bounds cannot hold for arbitrary $L^2$-functions. Consequently, a potentially higher asymptotic speed of the oscillation $\osc^\text{std}$ comes with a bigger constant in front of it and, therefore, with diminished non-asymptotic significance.
\end{remark}

\begin{remark}[oscillation degrees - data oscillation reduction]
\label{R:choosing-osc-degs;reduction-to-data-osc}
In the particular case of the Poisson equation, i.e.\ $\bA=\vec{I}$ and $c=0$, and linear elements, i.e.\ $n=1$, the oscillation  with the degrees $(m_1,m_2)=(0,0)$ reduces to the \emph{data oscillation}
\begin{equation*}
 \osc_\grid^\text{std}(u_\grid,\data)^2
 =
 \sum_{T\in\grid} h_T^2 \| f - \Pi_T^{0} f \|_{L^2(T)}^2;
\end{equation*}
it depends only on the data, here the right-hand side $f$. Note also that here the regularity of $f$ is determined by the regularity of the exact solution $u$.

For elements with degree $n\geq2$, the choices $(m_1,m_2)=(n-1,n-2)$ ensure that for $F\in \faces$
\begin{equation}
\label{pw-reproduction-of-Laplacian}
 \Pi_F^{n-1} \big( \jump{\nabla u_\grid}|_{F} \cdot \vec{n}_F \big)
 = 
 \jump{\nabla u_\grid}|_{F} \cdot \vec{n}_{F} 
\quad\text{and}\quad
 \Pi_T^{n-2} (\Delta u_\grid{}|_{T})
 =
 \Delta u_\grid{}|_{T}
\end{equation}
and so, again, oscillation reduces to data oscillation in $f$:
\begin{equation*}
 \osc_\grid^\text{std}(u_\grid,\data)^2
 =
 \sum_{T\in\grid} h_T^2 \| f - \Pi_T^{n-2} f \|_{L^2(T)}^2.
\end{equation*}
If we add a reaction term, viz.\ we consider $\bA=\vec{I}$ and $c=1$, we can again obtain the reduction to data oscillation by increasing $m_2$ to $n$.

For a more general operator with piecewise polynomial coefficients
\begin{equation*}
 \bA \in (\mathbb{S}_\grid^{n_{\bA},-1})^{d \times d}
\quad\text{and}\quad
 c \in \mathbb{S}_\grid^{n_c,-1},
\end{equation*}
the choice
\begin{equation}
\label{osc-degrees-for-pw-discrete-coefficients}
 (m_1,m_2)
 =
 \big( n_A + n - 1, \max\{n_c + n, n_A + n - 2\} \big)
\end{equation}
again reduces $\osc^\text{std}$ to data oscillation in $f$.

Finally, for a general operator without piecewise polynomial coefficients $(\bA, c)$, a reduction to data oscillation with piecewise polynomial best approximations as before is not possible. The argument in Remark~\ref{R:choosing-osc-degs;asymptotics} suggests approximating the general coefficients with piecewise polynomial coefficients  satisfying
\begin{equation*}
 n_{\bA} = n-1
\quad\text{and}\quad
 n_c = n-1.
\end{equation*}
As we shall see below in Sect.\ \ref{S:bds-for-corr}, the choice \eqref{osc-degrees-for-pw-discrete-coefficients} with these values allows us to bound $\osc^\text{std}$  in terms of $\| \nabla u_\grid \|_{L^2(\Omega)}$, which is controlled by stability, and data oscillation terms involving $f$ and the coefficients $\bA$ and $c$. Note however that the nature of these data oscillation terms differs from the preceding reductions: e.g., the regularity of the coefficients $\bA$ and $c$ is not determined by the exact solution $u$.
\end{remark}

In the light of Remark~\ref{R:choosing-osc-degs;asymptotics} (oscillation degrees - asymptotics), one might hope that the overestimation described in Remark~\ref{R:nonasym-overest} (nonasymptotic overestimation) disappears under refinement. This can be ensured under suitable regularity assumption but is not guaranteed in general as the following remark reveals.
\begin{remark}[asymptotic overestimation]
\label{R:asymp-overest}
Considering a variant of the standard residual estimator that allows for $f \in H^{-1}(\Omega)$ and adaptive refinement, \cite[Sect. 6.4]{CohenDeVoreNochetto:2012} give an example where the error converges asymptotically faster than the estimator; see also \cite[Lemma~21]{KreuzerVeeser:2021}.
\end{remark}

After having recognized the above flaws of the standard residual estimator, let us conclude with an observation that will be the departure point of an improved analysis.

\begin{corollary}[equivalence for discrete data]
\label{C:equiv-pw-poly-data}
Suppose all data $\data=(\bA,f,c)$ of problem \eqref{strong-form} are piecewise polynomial, i.e.\  there are $n_{\bA}, n_c, n_f \in \N_0$ such that
\begin{equation*}
 \bA \in (\mathbb{S}_\grid^{n_{\bA},-1})^{d \times d},
\quad
 c \in \mathbb{S}_\grid^{n_c,-1},
\quad\text{and}\quad
 f \in \mathbb{S}_\grid^{n_f,-1}.
\end{equation*}
Then error and standard residual estimator are equivalent:
\begin{equation*}
 \| \nabla(u-u_\grid) \|_{L^2(\Omega)}
 \approx
 \est_\grid^\text{std},
\end{equation*}
where the hidden constants depend only on $d$, the coefficients $\bA$ and $c$, the shape regularity coefficient $\sigma$, and the degrees $n_{\bA}$, $n_c$, and $n_f$.
\end{corollary}

\begin{proof}
The upper bound follows from Theorem~\ref{T:up-bd-with-std-res-est} (upper bound with standard residual estimator), while Proposition~\ref{P:lw-bd-with-std-res-est} (partial lower bound) with
\begin{equation*}
 (m_1,m_2)
 =
 \big( n_{\bA} + n - 1, \max\{ n_{\bA} + n - 2, n_c+n,n_f\} \big).
\end{equation*}
yields the lower bound. 
\end{proof}

Motivated by the above discussion, one may define a variant of the standard residual estimator, characterized by a splitting into two different parts. More precisely, choosing $(m_1,m_2)$ according to Remark~\ref{R:choosing-osc-degs;reduction-to-data-osc} (oscillation degrees - data oscillation reduction), one may replace the local indicators in \eqref{loc-ind-of-std-res-est} by
\begin{subequations}
\label{loc-ind-of-std-res-est;alt}
\begin{equation}
\index{Error Estimators!$\est_\grid^\text{std}(u_\grid,f,T)$: standard local indicators}
 \est_\grid^\text{std}(u_\grid,f,T)^2
 \definedas
 \eta_\grid^\text{std}(u_\grid,T)^2 + \osc^\text{std}_\grid(u_\grid,\data,T)^2,
\end{equation}
where the first part, the so-called \emph{PDE indicator}, is given by
\begin{equation}
\index{Error Estimators!$\eta_\grid^\text{std}(u_\grid,T)$: standard local PDE indicators}
 \eta_\grid^\text{std}(u_\grid,T)^2
 \definedas
 h_T^2  \| \Pi_T^{m_2} r \|_{L^2(T)}^2
 +
 h_T \sum_{F \subset \partial T \setminus \partial\Omega}
  \| \Pi_F^{m_1} j \|_{L^2(F)}^2 \,,
\end{equation}
\end{subequations}
while the second part corresponds to the local oscillation from \eqref{loc-ind-of-std-osc}; compare with \cite[Theorems~1.5 and 4.7]{Verfuerth:13}. In this way,
\begin{itemize}
\item the PDE indicators are computable (in terms of the Galerkin approximation $u_\grid$ and the local projections),
\item the oscillation indicators typically have to be approximated by numerical quadrature,
\item both types of indicators are, in general, not dominated by the error.

\end{itemize}
\subsection{Discrete functionals and a posteriori error analysis}
\label{S:discretize-functionals}
%
% intro
This section introduces the notion of \emph{discrete functionals} and individuates properties in their approximation that are useful in a~posteriori error analysis.  The realization of these properties distinguishes the subsequent approach, which is adapted from \cite{KreuzerVeeser:2021} and \cite{Kreuzer.Veeser.Zanotti}.

% discrete functionals and meshed subdomains
\smallskip The notion of discrete functionals and its local counterparts are of interest for at least two reasons. The first one is that their $H^{-1}$-norm can be rather easily quantified, as we shall see in Corollary \ref{C:quantify-discrete-dual-norms} below. This property is related to Corollary \ref{C:equiv-pw-poly-data}  (equivalence for discrete data), which can be read in the following way: the standard residual estimator is equivalent to the error whenever the residual is a discrete functional. The second reason lies in the observation that an important part of the residual, namely the application of the differential operator to a discrete function, is itself of discrete nature. This feature is partially captured by the following definition of discrete functionals with polynomial densities and is discussed in Remark \ref{R:diff-op-and-discrete-functionals}. 

\begin{definition}[discrete functionals and meshed subdomains]
\label{D:discrete-functionals}
Given $m_1 \in \N_0$, $m_2 \in \N_0 \cup \{-1\}$, let $\F_\grid \definedas \F(\grid) \definedas \F_{m_1,m_2}(\grid)$
\index{Functional Spaces!$\F_\grid$, $\F(\grid)$: discrete functionals}
denote the subspace 
\begin{multline*}
 \Big\{ \ell \in H^{-1}(\Omega) \mid 
   \forall w \in H^1_0(\Omega)
\;
    \langle \ell, w \rangle
    =
    \sum_{F\in \faces} \int_F q_F w
    +
    \sum_{T \in \grid} \int_T q_T w
\\
 \text{ with fixed }  q_F \in \P_{m_1}(F),\,  q_T \in \P_{m_2}(T)
 \Big\}
\end{multline*}
of \emph{discrete functionals}, i.e., functionals that are given by piecewise polynomial densities over elements and interior faces. We call $(m_1,m_2)$ the \emph{degrees of the discrete functionals}.

A set $\omega$ is a \emph{$\grid$-meshed subdomain}\index{Definitions!$\grid$-meshed subdomain} if it is a subdomain of $\Omega$ and it is triangulated by a submesh $\grid_\omega \subset\grid$ \index{Meshes!$\grid_\omega$: triangulated submesh}, i.e.\ we have $\overline{\omega} = \cup_{T \in \grid_\omega} T$. A functional $\ell \in H^{-1}(\Omega)$ is then \emph{discrete in the meshed subdomain $\omega$} whenever $\ell|_{H^1_0(\omega)} \in \F(\grid_\omega)$. Here the faces 
$$
\faces_\omega:= \{ F \in \faces \mid \ F \cap \omega \not = \emptyset, \ F \not \subset \partial \omega \}
\index{Meshes!$\faces_\omega$: faces interior to $\omega$}
$$
involved in $\F(\grid_\omega)$ are interior to $\omega$; e.g., the subspaces $\F(\{T\})$, $T \in \grid$, do not involve any faces. In accordance with \eqref{stars}, we use the abbreviations $\grid_z$\index{Meshes!$\grid_z$, $\grid_{\omega_z}$: elements forming $\omega_z$} and $\faces_z$\index{Meshes!$\faces_z$, $\faces_{\omega_z}$: faces interior to $\omega_z$} for $\grid_{\omega_z}$ and $\faces_{\omega_z}$.
\end{definition}

% local \F as restriction of global one
Alternatively, the local space $\F(\grid_\omega)$ can be obtained from the global space $\F(\grid)$ by restriction:
\begin{equation}
\label{local-F-by-restriction}
 \index{Functional Spaces!$\F(\grid_\omega)$: local discrete functionals}
\F(\grid_\omega)
 =
 \F(\grid)|_{H^1_0(\omega)}
 \definedas
 \left\{ \ell|_{H^1_0(\omega)} \mid \ell \in \F(\grid) \right\}.
 \end{equation}

\begin{remark}[differential operator and discrete functionals]
\label{R:diff-op-and-discrete-functionals}
The image of the finite element space $\V_\grid$ under the linear differential operator $-\operatorname{div} (\bA \nabla \cdot) + c\cdot$ is again a finite-dimensional space. For differential operators with piecewise polynomial coefficients $\bA$ and $c$, the above notion captures this by the property that the application of such operators to discrete functions $v \in \V_\grid$ yields discrete functionals. Indeed, if
\begin{equation*}
 \bA \in (\mathbb{S}_\grid^{n_{\bA},-1})^{d \times d},
\quad\text{and}\quad
 c \in \mathbb{S}_\grid^{n_c,-1}
\end{equation*}
with $n_{\bA}, n_c \in \N_0$, piecewise integration by parts gives the representation
\begin{equation}
\label{diff-op-on-fe-fct}
\begin{multlined}
    \int_\Omega \bA \nabla v \cdot \nabla w + c v w
\\
    =
    \sum_{F \in \faces} \int_F \jump{\bA \nabla v} \cdot \vec{n}_F w
    +
    \sum_{T \in \grid} \int_T \big( c v - \div (\bA \nabla v) \big) w,
\end{multlined}
\end{equation}
where, for any interior face $F\in \faces$ and any element $T \in \grid$,
\begin{equation*}
    \jump{\bA \nabla v} \cdot \vec{n}_F \in \P_{m_1}(F),
\quad
    c v - \div (\bA \nabla v) \in \P_{m_2}(T)
\end{equation*}
with $m_1 = n_{\bA} + n - 1$ and $m_2 = \max\{n_{\bA} + n - 2, n_c + n\}$. 
Note however that not every functional in $\F_{m_1,m_2}(\mesh)$ can be written in the form of \eqref{diff-op-on-fe-fct}. In fact, as the representation of a discrete functional is made up of $L^2$-scalar products on domains that are mutually disjoint or of different dimension, we have
\begin{equation}
\label{dimFmesh=}
 \dim \F_{m_1,m_2}(\mesh) = \#\faces \dim \P_{m_1} + \#\mesh \dim \P_{m_2},
\end{equation}
which is strictly greater than $\dim \V_\mesh$. This enlargement, which is implicitly used in the proof of Lemma \ref{P:lw-bd-with-std-res-est} (partial lower bounds), turns out to be convenient also in the constructive approximation of discrete functionals.
\end{remark}

% approximation of discrete functionals by standard approach
\smallskip In view of the aforementioned properties of discrete functionals, we may split the residual into a discrete and a non-discrete part. Splitting the standard residual estimator in the alternative local indicators \eqref{loc-ind-of-std-res-est;alt} is in a similar spirit.
To see this, we introduce $\Pi_\grid \ell \in H^{-1}(\Omega)$ given by
\begin{equation}
\label{E:Pigrid}
 \langle \Pi_\grid \ell, w \rangle
 \definedas
 \sum_{F \in \faces} \int_F \big( \Pi_F^{m_1}g \big) w
 +
 \sum_{T \in \grid} \int_T \big( \Pi_T^{m_2}f \big) w,
\quad
 w \in H^1_0(\Omega),
\end{equation}
for all $\ell \in H^{-1}(\Omega)$ admitting the representation
\begin{equation*}
 \langle \ell, w \rangle
 =
 \sum_{F \in \faces} \int_F g w
 +
 \sum_{T \in \grid} \int_T f w,
\quad
 w \in H^1_0(\Omega)
\end{equation*} 
with suitable density functions $g$ and $f$.
% Indeed, defining $\Pi_\grid$ as the operator given by the local $L^2$-projections appearing in the oscillation \eqref{std-osc}, it corresponds to writing
Then the splitting of the alternative indicators \eqref{loc-ind-of-std-res-est;alt} corresponds to writing
\begin{equation}
\label{splitting-residual-with-Pigrid}
 R_\grid
 =
 \Pi_\grid R_\grid + (I-\Pi_\grid) R_\grid.
\end{equation}
Moreover, Remark~\ref{R:choosing-osc-degs;reduction-to-data-osc} (oscillation degrees - data oscillation reduction) discusses in particular conditions for the identity
\begin{equation}
\label{red-to-osc(f);Pigrid}
 (I-\Pi_\grid) R_\grid = f-\Pi_\grid f,
\end{equation}
which follows from the property that $\Pi_\grid$ reproduces the functionals in Remark~\ref{R:diff-op-and-discrete-functionals} (differential operator and discrete functionals); compare with \eqref{pw-reproduction-of-Laplacian}, which in terms of $\Pi_\grid$ reads $\Pi_\grid(\Delta u_\grid) = \Delta u_\grid$, where $\Delta$ is now the distributional Laplacian.

% towards new approach
The fact that the definition of $\Pi_\grid$ requires the extra regularity $f \in L^2(\Omega)$ and $\bA \in W^{1}_{\infty}(\Omega;\mathbb R^{d\times d})$ not only excludes applications but, in light of Remark \ref{R:nonasym-overest}, (nonasymptotic overestimation) entails overestimation. To circumvent this flaw, we therefore aim at constructing a new approximation operator $P_\grid$ that is defined for all functionals $\ell \in H^{-1}(\Omega)$. Furthermore, we want this operator to be a projection onto $\F_\grid$ so that the counterpart
\begin{equation*}
 (I-P_\grid) R_\grid
 =
 f - P_\grid f 
\end{equation*}
of \eqref{red-to-osc(f);Pigrid} holds under the same conditions.

% summary of plan
To summarize, our plan is to develop a quasi-optimal a posteriori error analysis by constructing a locally computable, linear projection
\begin{equation*}
 P_\mesh:H^{-1}(\Omega) \to \F_\mesh \subset H^{-1}(\Omega)
\end{equation*}
onto the discrete functionals that induces a splitting 
\begin{equation}
\label{splitting-residual-with-Pgrid}
 R_\grid
 =
 P_\grid R_\grid + (I-P_\grid) R_\grid
\end{equation}
of the residual into a \emph{discretized residual} $P_\grid R_\grid$, which can be easily quantified, as well as an \emph{oscillatory residual} $(I-P_\grid) R_\grid$, which under the conditions of Remark \ref{R:diff-op-and-discrete-functionals} (differential operator and discrete functionals) reduces to an oscillation of the right-hand side $f$. 

% avoiding overestimation
The proof of an upper bound of the error will then involve a triangle inequality applied to the right-hand side of \eqref{splitting-residual-with-Pgrid}. The following remark provides criteria to prevent overestimation in such a context, and is followed by a comparison of the two approaches represented by \eqref{splitting-residual-with-Pigrid} and \eqref{splitting-residual-with-Pgrid}.

\begin{remark}[avoiding overestimation]
\label{R:avoiding-overestimation}
Overestimation can be often avoided by ensuring two relatively simple conditions. In order to discuss them informally, consider the model inequality
\begin{equation}
\label{model-ineq}
 |\cdot| \leq |\cdot|_1 + |\cdot|_2,
\end{equation}
where $|\cdot|$, $|\cdot|_i$, $i=1,2$, are seminorms and denote the domain and kernel of $|\cdot|$, respectively, by $\operatorname{dom} |\cdot|$ and $\operatorname{ker} |\cdot|$ etc.

% kernel criterion
The first condition, the \emph{kernel condition}, is that zero is not overestimated:
\begin{subequations}
\label{cond-for-acc-ineq}
\begin{equation}
\label{kernels-ok}
 \operatorname{ker} |\cdot|
 \subset
 \operatorname{ker} |\cdot|_1  \cap \operatorname{ker} |\cdot|_2 \,.
\end{equation}
% domain criterion 
The second condition, the \emph{domain condition}, is that a finite value is never bounded by $\infty$, or in other, still informal, words: if the evaluation of the left-hand side is (or can be uniquely defined to be) a finite value, the same holds for the right-hand side:
\begin{equation}
\label{domains-ok}
 \operatorname{dom} | \cdot |
 \subset
 \operatorname{dom} | \cdot |_1
 \cap
 \operatorname{dom} | \cdot |_2.
\end{equation}
\end{subequations}
% point to rigorous discussion
\cite{Kreuzer.Veeser.Zanotti} provide a precise version of the domain condition \eqref{domains-ok}, show that, given inequality \eqref{model-ineq}, the two conditions \eqref{cond-for-acc-ineq} are also sufficient for equivalence, and discuss further applications of this viewpoint. 
\end{remark}

In order to illustrate the application of Remark~\ref{R:avoiding-overestimation}, let us consider only the special case of the Poisson equation, i.e.\ $\bA = \vec{I}$, $c = 0$, and linear elements, i.e.\ $n=1$.  We start with the upper bound in terms of the standard residual estimator in Theorem~\ref{T:up-bd-with-std-res-est} and view it as a function of the right-hand side $f$. Then the domain condition is violated as the left-hand side is defined for any $f \in H^{-1}(\Omega)$, while the right-hand side is defined only for $f \in L^2(\Omega)$. Also the kernel condition is not verified: the left-hand side vanishes whenever $f = -\Delta v$ for some $v \in \V_\grid$, while the right-hand side vanishes only for $f=0$. The splitting in the alternative local indicators \eqref{loc-ind-of-std-res-est;alt} does not worsen this situation, i.e.\ it does not add further instances in which kernel and domain condition are missed. Note however that the oscillation indicators alone are in conflict with the domain condition and, therefore, another PDE indicator cannot cure the overestimation. Finally, for the outlined approach, the splitting \eqref{splitting-residual-with-Pgrid} and the required properties for the operator $P_\grid$ ensure both kernel and domain condition.

\subsection{Testing discrete functionals}
\label{S:approx-of-discrete-functionals}
%
% intro
The $H^{-1}$-projection $P_\grid$ onto the discrete functionals $\F_\grid$ will be defined by means of a Petrov-Galerkin-type approach. This section prepares its definition by individuating a suitable \emph{test space} $\V^+_\mesh$. The key property of $\V^+_\mesh$ is that the dual pairing $\langle\cdot,\cdot\rangle$ in $H^{-1}(\Omega)$ is nondegenerate on the product $\F_\mesh \times \V^+_\mesh$.
%, i.e., $\ell \in \F_\mesh$ with $\langle \ell, w \rangle = 0$ for all $w \in \V^+_\mesh$ implies $\ell=0$
Doing so, the degrees $(m_1,m_2)$ of the discrete functionals will be parameters
% in $\N_0 \times (\N_0 \cap \{-1\})$
that are omitted in the notation. The construction of the test space $\V^+_\grid$ proceeds in two steps. First, we locally associate to the degrees of freedom in $\F_\mesh$ certain functions on $\Omega$. For the degrees of freedom on the skeleton, this will inolve a suitable \emph{extension operator}.  Second, we turn the ensuing functions into admissible test functions with the help of a \emph{cut-off}.

% first step: extensions
\smallskip The degrees of freedom in $\F_\mesh$ are given by density polynomials over element and faces. For an element $T\in\mesh$, if we extend such a density polynomial $q_T$ by $0$ off $T$, it is already a function on $\Omega$. For a polynomial $q_F$ associated with a face $F\in\faces$, we employ the following extension operator $E_F$ mapping a function $v$ on $F$ to a function on $\omega_F$, the union of all elements $T$ containing $F$.

%Given such an element $T \subset \omega_F$, write $z_0,\dots,z_d$ for its vertices, $z_d$ being the one %opposite to $F$, and set
%\begin{equation}\label{e:E_F}
%\begin{split}
% \big(E_F v\big)(x)
% &\definedas
% v \left( z_0 + \sum_{i=1}^{d-1} \phi_{z_i}(x)(z_i-z_0) \right)
% \\
%  &=
%  v \left( \left( \phi_{z_0}(x) + \phi_{z_d}(x) \right) z_0 + \sum_{i=1}^{d-1} \phi_{z_i}(x) z_i \right),
%\qquad
% x \in T,
% \end{split}
%\end{equation}
Given such an element $T \subset \omega_F$, write $z_0,\dots,z_d$ for its vertices, $z_d$ being the one opposite to $F$, denote by $b_F \definedas \frac{1}{d} \sum_{i=0}^{d-1} z_i$ the barycenter of $F$, and set
\begin{equation*}
\label{e:E_F}
%\begin{split}
 \big(E_F v\big)(x)
 \definedas
 v \left( \phi_d(x) b_F + \sum_{i=0}^{d-1} \phi_{z_i}(x) z_i \right)
\qquad
 x \in T,
 %\end{split}
\end{equation*}
and extend by 0 off $\omega_F$. Note that the definition of $E_F$ is affine invariant and does not depend on the enumeration of the vertices of $F$. The next lemma collects two useful properties of this extension operator. 

\begin{lemma}[extending from faces]
\label{L:extending-from-faces}    
Let $F \in \faces$ be a face. For any function $v \in L^2(F)$, we have
\begin{align*}
 \| E_F v \|_{L^2(\omega_F)}
 \lesssim
 h_F^{\frac{1}{2}} \| v \|_{L^2(F)},
\end{align*}
where $h_F$ stands for the diameter of $F$ and the hidden constant depends only on $d$ and the shape regularity coefficient $\sigma$. Furthermore, if $v$ is a polynomial, then $E_F v$ is a continuous piecewise polynomial of the same degree.
\end{lemma}

\begin{proof}
\fbox{\scriptsize 1} In view of $(E_Fv)^2 = E_F(v^2)$, we may show the inequality by verifying
\begin{equation}
\label{L1-stability-of-EF}
 \int_{\omega_F} E_Fw \lesssim h_F \int_F w
\end{equation}
for any positive function $w:F\to\mathbb{R}$, which amounts to $L^1$-stability. To this end, we shall use a standard argument involving the following reference situation, which slightly differs from the common one with $T_d = \big\{ x=(x_1,\dots,x_d) \in \mathbb R^d \mid 0\leq x_i \leq 1, \ \sum_{i=1}^d x_i \leq 1 \big\}$  and $b_d \definedas (d+1)^{-1}(1,\dots,1) \in \R^d$. Let the reference face $\wh{F} \definedas T_{d-1} - b_{d-1} \subset \R^{d-1}$ be a translation of $T_{d-1}$ and let the reference simplex $\wh{T} \subset \R^d$ be the convex hull of $\wh{F} \times \{ 0 \}$ and the canonical basis vector $e_d = (0,\dots,0,1) \in \R^d$. The barycenter of $\wh{F}\times\{0\}$ is then the origin in $\R^d$ and the barycentric coordinate of the vertex $e_d$ of $\wh{T}$ is $x_d$. Fixing an element $T$ with $T \subset \omega_F$, denote
%by $z_0,\dots,z_d$ the vertices of $T$ and
by $G_T:\wh{T}\to T$ a bi-affine map sending vertices of $\wh{F} \times \{0\}$ into vertices of $F$ and $e_d$ into the vertex of $T$ opposite to $F$, and write $G_F:\wh{F}\times\{0\} \to F$ for the restriction $G_T{}|_{\wh{F}\times\{0\}}$. The pullbacks of $E_Fw$ and $w$ satisfy
\begin{equation*}
 G_T^*(E_Fw)(x',x_d) = G_F^*w(x',0)
\end{equation*}
for all $x = (x',x_d) \in \wh{T} = \big\{ y=(y',y_d) \in \wh{F} \times \mathbb{R} \mid 0 \leq y_d \leq 1-|y'+b_{d-1}|_1 \big\}$, where $|z'|_1 = \sum_{i=1}^{d-1} |z'_i|$ stands for the $\ell_1$-norm in $\mathbb{R}^{d-1}$. Consequently, the transformation rule, the fact that the Jacobians of $G_T$ and $G_F$ are constant, the Fubini theorem, $w \geq 0$, and $|\wh{F}|/|\wh{T}| = |T_{d-1}|/|T_d| = d$ yield
\begin{align*}
 \int_T E_Fw
 &=
 \frac{|T|}{|\wh{T}|} \int_{\wh{T}} G_T^*(E_Fw)
 =
 \frac{|T|}{|\wh{T}|} \int_{\wh{F}} \int_{0}^{1-|x'+b_{d-1}|_1}
  G_F^*w(x',0) \,\mathrm{d}x_d \,\mathrm{d}x'
\\
 &\leq
 \frac{|T|}{|\wh{T}|} \int_{\wh{F}} G_F^*w(\cdot,0)
 =
 \frac{|T| |\wh{F}|}{|\wh{T}| |F|}  \frac{}{} \int_F w
 =
 d \frac{|T|}{|F|} \int_F w.
\end{align*}
Since the hidden constant in $|T|/|F| \lesssim h_F$ depends only on the shape regularity coefficient $\sigma$, this implies the $L^1$-stability bound \eqref{L1-stability-of-EF} and so also the claimed $L_2$-stability is proven.

\smallskip\fbox{2} The second statement for polynomial arguments of $E_F$ is a direct consequence of its definition.
\end{proof}

% second step: cut-off
In view of the above different treatment of elements and faces, we need two types of cut-off functions: one for elements denoted by $\phi_T$ and another one for faces denoted by $\phi_F$.  Possible choices are the element and face bubbles from \eqref{elm-bubble} and \eqref{face-bubble}. Since other choices will be useful in Sect.\ \ref{S:bds-for-corr} below, we shall rely henceforth only on the following properties.

\begin{assumption}[abstract cut-off]
\label{A:abstract-cut-off}
\index{Assumptions!Abstract cut-off}
The \emph{cut-off functions} $\phi_T$, $T \in \grid$, and $\phi_F$, $F \in \faces$, satisfy
\begin{gather*}
 \mathop{\mathrm{supp}} \phi_T = T,
\quad
 0 \leq \phi_T \leq 1,
\quad
 \mathop{\mathrm{supp}} \phi_F = \omega_F,
\quad
 0 \leq \phi_F \leq 1,
\end{gather*}
and act in an affine-equivalent manner on the element level: there exists a finite-dimensional linear space $\mathbb{S}^+ \subset L^\infty(T_d)$ of functions defined on the reference element $T_d$ such that
$G_T^*\phi_T$ does not depend on $T$, $G_F^*\phi_F$ does not depend on $F$, and
\begin{align*}
 &\forall T \in \mesh , \  \forall q \in \P_{m_2}(T)
\quad
 G_T^*\left( q\phi_T \right) \in \mathbb{S}^+,
\\
  &\forall F \in \faces , \  \forall q \in \P_{m_1}(F), \  \forall T \in \mesh 
\quad
 G_T^*\Big( \big(\big(E_Fq\big)\phi_F \big)|_{T} \Big)\in \mathbb{S}^+,
\end{align*}
where $G_T$ is a bi-affine map from the reference element $T_d$ to the generic element $T$ and $G_T^*(v) = v\circ G_T$ denotes the pullback of a function $v:T\to\mathbb{R}$ via $G_T$. 
\end{assumption}

In the case of the bubble functions \eqref{elm-bubble} and \eqref{face-bubble}, Assumption \ref{A:abstract-cut-off} holds with $\mathbb{S}^+ = \P_{\max\{m_1 + d-1, m_2 + d\} }(T_d)$ as the extension operators $E_F$, $F \in \faces$, preserve the polynomial degree.

\begin{lemma}[properties of cut-off]
\label{L:cut-off-prop}
If the cut-off functions $\phi_T$, $T \in \mesh$, and $\phi_F$, $F \in \faces$, satisfy Assumption \ref{A:abstract-cut-off}, then we have
\begin{align*}
 \| q \|_{L^2(T)} &\Cleq \| q \phi_T^{1/2} \|_{L^2(T)}
\quad\text{and}\quad
 \| \nabla(q\phi_T) \|_{L^2(T)} \Cleq h_T^{-1} \| q\phi_T \|_{L^2(T)}
\end{align*}
for all $q \in \P_{m_2}(T)$ as well as
\begin{align*}
 \| q \|_{L^2(F)} &\Cleq \| q \phi_F^{1/2} \|_{L^2(F)}
\quad\text{and}\quad
 \| \nabla\big( (E_Fq)\phi_F\big) \|_{L^2(\omega_F)}
 \Cleq
 h_F^{-1} \| (E_Fq)\phi_F \|_{L^2(\omega_F)}
\end{align*}
for all $q \in \P_{m_1}(F)$. The hidden constants depend only on $d$, the shape regularity coefficient $\sigma$, the degrees $(m_1,m_2)$ of the discrete functionals, and the space $\mathbb{S}^+$. 
\end{lemma}

\begin{proof}
\fbox{\scriptsize 1} To verify the first claimed inequality, we start by noting that, thanks to $\supp \phi_T = T$, we have $\phi_{T_d} \definedas G_T^*(\phi_T)>0$ in the interior of $T_d$. Hence, $\|\cdot\|_{L^2(T_d)}$ and $\|\,\cdot\,\phi_{T_d}^{1/2}\|_{L^2(T_d)}$ are norms on $\P_{m_2}(T_d)$ and, thanks to $\dim \P_{m_2}(T_d) < \infty$, are equivalent. A standard round trip to the reference element and $G_T^*q \, G_T^*\phi_T^{1/2} = G_T^*(q\phi_T^{1/2})$ thus yield
\begin{align*}
 \| q \|_{L^2(T)}
 &\Cleq
 h_T^{d/2} \| G_T^*q \|_{L^2(T_d)}
 \Cleq
  h_T^{d/2} \| G_T^*q \, G_T^* \phi_T^{1/2}\|_{L^2(T_d)}
\\
 &\Cleq
 \| q \phi_T^{1/2}\|_{L^2(T)},
\end{align*}
and the first claimed inequality is established. The third one is proved along the same lines, but with a round trip to the reference face.

\smallskip\fbox{\scriptsize 2} For the other claimed inequalities note that $\|\nabla\,\cdot\,\|_{L^2(T_d)}$ and $\inf_{c \in \mathbb{R}}\|\,\cdot\,-c\|_{L^2(T_d)}$ are equivalent norms on the finite-dimensional quotient space $\mathbb{S}^+/\mathbb{R}$. Consequently, further round trips to the reference element give 
\begin{align*}
 \| \nabla(q\phi_T) \|_{L^2(T)}
 &\Cleq
 h_T^{-1+d/2} \| \nabla G_T^*(q\phi_T) \|_{L^2(T_d)}
 \Cleq
 h_T^{-1+d/2} \inf_{c \in \mathbb{R}} \| G_T^*(q\phi_T) - c \|_{L^2(T_d)}
\\
 &\Cleq
 h_T^{-1+d/2} \| G_T^*(q\phi_T) \|_{L^2(T_d)}
 \Cleq
  h_T^{-1} \| q \phi_T \|_{L^2(T)}
\end{align*}
and
\begin{align*}
 \| \nabla \big( E_F(q)\phi_T \big) \|_{L^2(\omega_F)}^2
 &=
 \sum_{T \subset \omega_F}
 \| \nabla \big( E_F(q)\phi_T \big) \|_{L^2(T)}^2
\\
 &\Cleq
 h_F^{-1} \sum_{T \subset \omega_F} \| E_F(q) \phi_T \|_{L^2(T)}^2
 =
 h_F^{-1} \| E_F(q) \phi_T \|_{L^2(\omega_F)}^2
\end{align*}
and the proof is completed.
\end{proof}

These preparations lead to the following test space for discrete functionals.
\begin{definition}[test space for discrete functionals]
\label{D:V^+}
Using the cut-off functions from Assumption \ref{A:abstract-cut-off}, we associate to the space $\F_\grid$  of discrete functionals the following \emph{test space}:
\begin{equation*}
%\label{Def-V^+}
\index{Functional Spaces!$\V^+(\grid)$, $\V^+_\grid$: test space for discrete functionals}
\begin{multlined}
    \V^+_\mesh
    \definedas
    \V^+(\mesh)
    \definedas
    \operatorname{span} \left(
     \big\{ q_T \phi_T \mid q_T \in \P_{m_2}(T), \; T \in \grid \big\}
     \phantom{\bigcup} \right.
\\
    \left.  \bigcup
   \big\{
     E_F(q_F) \phi_F \mid q_F \in \P_{m_1}(F),\; F \in \faces
    \big\}
    \right).
\end{multlined}
\end{equation*}
If $\omega$ is a subdomain of $\Omega$ meshed by $\grid_\omega$, then $\V^+(\grid_\omega)$ is the test space for $\F(\grid_\omega)$.
\end{definition}

Similarly as for the space $\F_\grid$ of discrete functionals, the test space over a subdomain $\omega$ meshed by $\grid_\omega$ can be obtained from the global test space, namely 
\begin{equation}
\label{localV+-by-restriction}
\index{Functional Spaces!$\V^+(\grid_\omega)$: local test space for discrete functionals}
 \V^+(\grid_\omega) = \V^+(\grid) \cap H^1_0(\omega).
\end{equation}

\subsection{A projection onto discrete functionals}
\label{S:proj-onto-discrete-functionals}
Having the test space $\V_\grid^+$ from Definition \ref{D:V^+} at our disposal, we are now ready to construct a $H^{-1}$-\emph{projection} $P_\grid$ as suggested in Sect.\ \ref{S:discretize-functionals}. Like in the previous section, the degrees $(m_1,m_2)$ of the discrete functionals in $\F_\grid$ are hidden parameters.

\begin{definition}[projection onto discrete functionals]
\label{D:Pgrid}
Given the discrete functionals $\F_\grid$ and the test space $\V^+_\grid$, we define a \emph{projection} $P_\grid:H^{-1}(\Omega) \to \F_\mesh$ by \looseness=-1
\begin{equation}\label{E:defPmesh}
    \langle P_\grid \ell, w \rangle
    =
    \langle \ell, w \rangle \quad \forall w \in \V^+_\mesh.
 \index{Operators!$P_\grid$: projection operator from $H^{-1}(\Omega)$ into $\F_\grid$}   
\end{equation}
The well-posedness of this definition and algebraic properties of $P_\grid$ are verified in the following Lemma~\ref{L:algebra-of-P}. Moreover, a representation of $P_\grid$ in form of a quasi-interpolation operator is given in Corollary~\ref{C:proj-as-interpolators} below.

The \emph{polynomial densities} of $P_\mesh\ell$ are denoted by $P_T\ell \definedas P_{\grid, T}\ell$, $T\in\mesh$, and $P_F \ell \definedas P_{\grid,F}$, $F\in\faces$, so that
\begin{equation}
\label{def-densities-of-P}
\index{Operators!$P_T$, $P_F$: polynomial densities of $P_\grid$}
 \langle P_\mesh \ell, w \rangle
 =
 \sum_{T \in \mesh} \int_T P_T\ell \, w + \sum_{F \in \faces} \int_F P_F\ell \, w.
\end{equation}
\end{definition}

In the next lemma we show in particular that $P_\grid$ is a local operator. In order to formulate this, we shall use \emph{$\grid$-meshed local subdomains}, i.e.\ $\grid$-meshed subdomains $\omega$ for which there exists a mesh element $T\in\grid$ with $\omega \subset \omega_T$.
For the next lemma addressing algebraic properties of the operator $P_\grid$, recall the notation $\omega_F = \cup_{T\in\grid_F} T$ from \eqref{pairs} for an interior face $F \in \faces$.

\begin{lemma}[algebraic properties]% of $P_\grid$]
\label{L:algebra-of-P}
The operator $P_\grid$ is a local linear projection onto the subspace $\F_\mesh$ of discrete functionals. More precisely, for any local subdomain $\omega$ meshed by $\grid$, there is a linear projection $P_\omega:H^{-1}(\omega) \to \F(\grid_\omega)$ such that \looseness=-1
\begin{equation*}
 P_\grid \ell|_{H^{-1}(\omega)}
 =
 P_\omega \left( \ell|_{H^1_0(\omega)} \right)
 \in \F(\grid_\omega)
\end{equation*} 
for all $\ell \in H^{-1}(\Omega).
$\end{lemma}

\begin{proof}
% linear independence of DOFs in V^+
\fbox{\scriptsize 1} We first show that the degrees of freedom in Definition \ref{D:V^+} of $\V^+_\mesh$ are linearly independent. To this end, we fix an element $T\in\mesh$, denote  by $F_1, \dots, F_l$, $l\leq d+1$ its faces that are in $\F_\grid$ and write $\phi_0 \definedas \phi_T$, $\phi_i \definedas \phi_{F_i}$ and $E_i \definedas E_{F_i}$ for $i=1,\dots,l$. We then claim that, for all $q_0 \in \P_{m_2}(T)\setminus\{0\}$ and all $q_i \in \P_{m_1}(F_i)\setminus\{0\}$, $i=1,\dots,l$, we have 
\begin{equation}
\label{lin-indep-of-dofs-of-V+}
 \alpha_0 q_0 \phi_0 + \sum_{i=1}^l \alpha_i E_i(q_i) \phi_i = 0 \text{ in } T
 \implies
 \alpha_0 = \dots = \alpha_l = 0.
\end{equation}
In light of $\supp\phi_T = T$ and $\supp \phi_F = \omega_F$, we observe $\phi_0{}|_{\partial T} = 0$ and $\phi_i{}|_{\partial T \setminus F_i} = 0$ for $i=1,\dots,l$. We thus evaluate the hypothesis of \eqref{lin-indep-of-dofs-of-V+} first on the faces $F_1,\dots,F_l$ and then in the element $T$. This gives $\alpha_i = 0$ for $i=0,\dots,l$ and \eqref{lin-indep-of-dofs-of-V+} is verified.

% well-posedness of P
\smallskip\fbox{\scriptsize 2} Next, we discuss the well-posedness of \eqref{E:defPmesh}. The linear independence \eqref{lin-indep-of-dofs-of-V+} and \eqref{dimFmesh=} lead to
\begin{equation}
\label{dimV+mesh=dimFmesh}
 \dim \V^+_\mesh
 =
 \#\mesh \dim \P_{m_2} + \#\faces \dim \P_{m_1}
 =
 \dim \F_\mesh.
\end{equation}
Thus, it suffices to show the implication
\begin{equation}
\label{pairing-nondeg-on-Fmesh-V+mesh;2}
  \ell \in \F_\mesh \, :  \  \langle \ell, w \rangle = 0 \,\, \forall w \in \V^+_\mesh
 \quad \implies \quad 
 \ell = 0.
\end{equation}
For that purpose, let $\ell \in \F_\mesh$ and let $q_T$, $T \in \mesh$, and $q_F$, $F \in \faces$ be the determining polynomials. For any element $T\in\mesh$, the choice $w = q_T\phi_T \in H^1_0(T)$ implies
\begin{equation*}
 0 = \langle \ell, w \rangle = \int_T (q_T)^2 \phi_T,
 \quad\text{i.e.}\quad
 q_T = 0.
\end{equation*}
Thus, $\ell$ does not have contributions from elements. Regarding faces, any choice $w = \big(E_F q_F\big) \phi_F \in H^1_0(\omega_F)$, $F\in\faces$, therefore gives
\begin{equation*}
 0 = \langle \ell, w \rangle = \int_F (q_F)^2 \phi_F,
\quad\text{i.e.}\quad
 q_F = 0,
\end{equation*}
and implication \eqref{pairing-nondeg-on-Fmesh-V+mesh;2} is established. Combining \eqref{dimV+mesh=dimFmesh} and \eqref{pairing-nondeg-on-Fmesh-V+mesh;2}, we have that the dual pairing in $H^{-1}(\Omega)$ is nondegenerate on $\F_\mesh \times \V^+_\mesh$, which in turn ensures that $P_\mesh$ is well-defined. The Petrov-Galerkin character of the definition \eqref{E:defPmesh} then ensures that $P_\mesh$ is a linear projection onto $\F_\mesh$.

\smallskip\fbox{\scriptsize 3} It remains to show that $P_\mesh$ is a local operator. Given any local subdomain $\omega$ meshed by $\grid_\omega$, we can apply the preceding proof to $\grid_\omega$ instead of $\grid$. This shows that the dual pairing in $H^{-1}(\omega)$ is nondegenerate on $\F(\grid_\omega) \times \V^+(\grid_\omega)$ and ensures a local projection operator $P_\omega:H^{-1}(\omega) \to \F(\grid_\omega)$. Taking into account \eqref{local-F-by-restriction} and \eqref{localV+-by-restriction}, we note
\begin{equation*}
 P_\grid{}|_{H^{-1}(\omega)} = P_\omega,
\end{equation*}
which completes the proof.
\end{proof}

The verification of \eqref{pairing-nondeg-on-Fmesh-V+mesh;2}, which amounts to a proof of uniqueness, suggests the following approach to compute $P_\mesh\ell$, $\ell \in H^{-1}(\Omega)$.

\begin{remark}[local computation]
\label{R:Computation-of-P}
Let $\ell \in H^{-1}(\Omega)$. Recalling \eqref{def-densities-of-P}, the polynomials $P_T\ell$, $T \in \mesh$, and $P_F\ell$, $F \in \faces$ can be computed by solving first
\begin{equation}\label{e:Computation-of-PT_local}
 \int_T P_T\ell \, q \phi_T = \langle \ell, q \phi_T \rangle
 \quad
 \forall T \in \mesh, q \in \P_{m_2}(T),
\end{equation}
and then
\begin{equation}\label{e:Computation-of-PF_local}
 \int_F P_F\ell \, q \phi_F
 =
 \langle \ell, q \phi_F \rangle - \sum_{T\subset\omega_F} \int_T P_T\ell \, q \phi_F
 \quad
 \forall F \in \faces, q \in \P_{m_1}(F).
\end{equation}
This amounts to two block-diagonal linear systems with, respectively, $\#\mesh$ blocks of size $\dim \P_{m_2}$ and $\#\faces$ blocks of size $\dim \P_{m_1}$. Each block and each corresponding right-hand side arises from local computations.
\end{remark}

\begin{remark}[star localization vs locality of $P_\grid$]
\label{R:stars-vs-pairs}
Stars $\omega_z$, $z \in \vertices$, are meshed local subdomains. Lemma \ref{L:algebra-of-P} thus shows that, for any vertex $z \in \vertices$, there is a linear projection $P_z:H^{-1}(\omega_z) \to \F(\grid_z)$ such that
\begin{equation*}
 P_\grid \ell|_{H^{-1}(\omega_z)}
 =
 P_z \left( \ell|_{H^1_0(\omega_z)} \right)
 \in \F(\grid_z)
\end{equation*} 
for all $\ell \in H^{-1}(\Omega)$. The stars also appear in the localizing upper bound of the global residual norm in Corollary~\ref{C:residual-loc}. As they are minimal subdomains therein, cf.\ Remark \ref{R:star-loc-is-minimal}, it may appear that a finer localization with smaller domains cannot be exploited in a~posteriori analysis. Although this is true in the context of upper bounds, the increased locality of $P_\grid$ is useful in the context of lower bounds; see the reduced lower bound \eqref{E:local-lower-reduced}, which follows from the interior vertex property introduced in Definition \ref{D:interior-node}, and is crucial to derive the contraction result \eqref{E:simplest-contraction}.
\end{remark}

% local stability properties of Pgrid
We have already mentioned that we shall use $P_\grid$ to split the residual. In light of the bounds for the residual norm in Corollary~\ref{C:residual-loc} (star localization of residual norms), this should be done in a locally stable manner. In order to formulate and employ the local stability properties of $P_\grid$, the following notation is useful. Given a local subdomain $\omega$ meshed by $\grid_\omega$, we define the \emph{$\V^+(\grid_\omega)$-discrete dual norm} by
\begin{equation}
\label{def-|.|_V^+(omega)'}
 \| \ell \|_{\V^+(\grid_\omega)^*}
 \definedas
 \sup_{w \in \V^+(\grid_\omega), \| \nabla w \|_{L^2(\omega)}=1} \langle \ell, w \rangle
\quad
 \ell \in H^{-1}(\omega).
\end{equation}
In view of $\V^+(\grid_\omega) \subset H^1_0(\omega)$, we have $\| \ell \|_{\V^+(\grid_\omega)^*} \leq \| \ell \|_{H^{-1}(\omega)}$ for all $\ell \in H^{-1}(\omega)$.

\begin{lemma}[local $H^{-1}$-stability]% properties of $P_\grid$]
\label{L:local-stability-of-P}
The projection $P_\grid$ is locally $H^{-1}$-stable: for any local subdomain $\omega$ meshed by $\grid_\omega$, we have
\begin{align*}
 \| P_\grid \|_{\mathcal{L}(H^{-1}(\omega))} = 
 \sup_{\ell \in \F(\grid_\omega)} \frac{\|\ell\|_{H^{-1}(\omega)}}{\|\ell\|_{\V^+(\grid_\omega)^*}}
 \leq
 C_\mathrm{lStb},
\end{align*}
where $C_\mathrm{lStb} = C_\mathrm{lStb}(d,\sigma, m_1, m_2)$ depends only on $d$, the shape regularity coefficient $\sigma$ from \eqref{E:shape-regularity}, the degrees $(m_1,m_2)$ of the discrete functionals and the space $\mathbb{S}^+$.

%Furthermore, the pairs $\big( \F(\grid_\omega, \V^+(\grid_\omega) \big)$ are uniformly inf-%sup stable with respect to the pairing $\langle \cdot, \cdot\rangle $ of  $H^{-1}(\omega) %\times H^1_0(\omega)$ for all local meshed subdomains $\omega$.
\end{lemma}
\begin{proof}
\fbox{\scriptsize 1} We start by verifying the '$\leq$'-part of the claimed identity for the operator norm. The definition of the operator norm leads to
\begin{equation*}
  \| P_\omega\|_{\mathcal{L}(H^{-1}(\omega))} =
  \sup_{\ell \in H^{-1}(\omega)} \frac{\|P_\omega\ell\|_{H^{-1}(\omega)}}{\|\ell\|_{H^{-1}(\omega)}}
  \le 
 \sup_{\ell \in H^{-1}(\omega)} \frac{\|P_\omega \ell\|_{H^{-1}(\omega)}}{\|\ell\|_{\V^+(\grid_\omega)^*}}.
\end{equation*}
We now notice that $\|\ell\|_{\V^+(\grid_\omega)^*} = \|P_\omega\ell\|_{\V^+(\grid_\omega)^*}$ in view of
\eqref{def-|.|_V^+(omega)'} and \eqref{E:defPmesh}. Hence,
\[
\| P_\omega\|_{\mathcal{L}(H^{-1}(\omega))} \le
\sup_{\ell \in H^{-1}(\omega)} \frac{\|P_\omega \ell\|_{H^{-1}(\omega)}}{\|P_\omega\ell\|_{\V^+(\grid_\omega)^*}}
= \sup_{\ell \in \F(\grid_\omega)} \frac{\|\ell\|_{H^{-1}(\omega)}}{\|\ell\|_{\V^+(\grid_\omega)^*}}
\]
because the projection $P_\omega$ is onto $\F(\grid_\omega)$. 

\smallskip\fbox{\scriptsize 2} Next, we show that $\| P_\grid \|_{\mathcal{L}(H^{-1}(\omega))}$ is uniformly bounded.
Let $\ell \in \F(\grid_\omega)$ be a discrete functional, namely
\begin{equation*}
 \langle\ell, w\rangle = \sum_{T\in\grid_\omega} \int_T q_Tw + \sum_{F \in \faces_\omega} \int_F q_Fw
 \quad\forall w\in H^1_0(\omega).
\end{equation*}
We proceed in two steps that are quite similar to the classical standard residual estimates. Arguing as in \eqref{loc-up-bd-with-std-res-est}, and using Lemma \ref{L:Poincare} (first Poincar\'e inequality)
in the domain $\omega$ of diameter about $h_T$ for $w \in H^1_0(\omega)$, we obtain
\begin{align*}
 |\langle \ell, w \rangle|
% &\leq
% \sum_{T \in \mesh_z} \left|\int_T q_T w \right|
%  + \sum_{F\in\faces_z} \left| \int_F q_F w \right|
%\\
 &\Cleq
 \left(
  \sum_{T \in \grid_\omega} h_T^2\| q_T \|_{L^2(T)}^2
   + \sum_{F \in \faces_\omega} h_F\| q_F \|_{L^2(F)}^2
\right)^{1/2} \| \nabla w \|_{L^2(\omega)}
\end{align*}
whence
\begin{equation}
\label{Pstability;upper-bd}
 \| \ell \|_{H^{-1}(\omega)}^2
 \Cleq
 \sum_{T \in \mesh_\omega} h_T^2 \| q_T \|_{L^2(T)}^2
   + \sum_{F \in \faces_\omega} h_F \| q_F  \|_{L^2(F)}^2.
\end{equation}
Here we do not exploit that $\ell$ is discrete in $\omega$; this will be crucial in the second step, when we bound each term on the right-hand side, in a manner reminding the derivation of classical lower bounds. For any $T \in \mesh_\omega$, we write $\V^+(T)^*$ as shorthand for $\V^+(\{T\})^*$ and exploit Lemma \ref{L:cut-off-prop} (properties of cut-off) to deduce
\begin{equation*}
\begin{aligned}
 \|q_T\|_{L^2(T)}^2
 &\Cleq
 \int_T \left( q_T \right)^2 \phi_T
 =
 \langle \ell, q_T \phi_T \rangle
 \leq
 \| \ell \|_{\V^+(T)^*} \| \nabla (q_T\phi_T) \|_{L^2(T)}
\\
 &\Cleq
 \| \ell \|_{\V^+(T)^*} h_T^{-1} \| q_T\phi_T \|_{L^2(T)}
 \leq
 \| \ell \|_{\V^+(T)^*} h_T^{-1} \| q_T \|_{L^2(T)},
\end{aligned}
\end{equation*}
whence
\begin{equation}
\label{elm-den<=disc-dual-norm}
 h_T \|q_T\|_{L^2(T)}
 \Cleq
 \| \ell \|_{\V^+(T)^*}.
\end{equation}
For an interior face $F \in \faces_\omega$, we proceed similarly, taking into account also Lemma~\ref{L:extending-from-faces} and that $\V^+(T) \subset \V^+(\grid_F)$ entails $\| \ell \|_{\V^+(T)^*} \leq \| \ell \|_{\V^+(\grid_F)}$ for $T \in \grid_F$. We thus obtain
\begin{equation*}
\begin{aligned}
 \|q_F &\|_{L^2(F)}^2
 \Cleq
 \int_F \left( q_F \right)^2 \phi_F
 =
 \langle \ell, \big(E_Fq_F\big) \phi_F \rangle
 -
 \sum_{T \subset \omega_F} \int_T q_T \big(E_Fq_F\big) \phi_F
\\
 &\Cleq
 \| \ell \|_{\V^+(\grid_F)^*} \| \nabla \big((E_Fq_F)\phi_T\big) \|_{L^2(\omega_F)}
 +
 \sum_{T \subset \omega_F} \| q_ T \|_{L^2(T)} \| (E_Fq_F)\phi_T \|_{L^2(T)}
\\
 &\Cleq
 \| \ell \|_{\V^+(\grid_F)^*} h_T^{-1} \| \big( E_Fq_F\big) \phi_T \|_{L^2(\omega_F)}
 \Cleq
 \| \ell \|_{\V^+(\grid_F)^*} h_T^{-\frac{1}{2}} \| q_F \|_{L^2(F)},
\end{aligned}
\end{equation*}
i.e.\ 
\begin{equation}
\label{face-den<=disc-dual-norm}
 h_F^{\frac{1}{2}} \|q_F\|_{L^2(F)}
 \Cleq
 \| \ell \|_{\V^+(\grid_F)^*}.
\end{equation}
The number of elements and interior faces in the local subdomain $\omega$ is uniformly bounded by $d$ and the shape regularity coefficient $\sigma$. Hence, inequalities \eqref{elm-den<=disc-dual-norm} and \eqref{face-den<=disc-dual-norm} together with the inclusions $\V^+(T) \subset \V^+(\grid_\omega)$ for $T \in \grid_\omega$ and  $\V^+(\grid_F) \subset \V^+(\grid_\omega)$ for $F \in \faces_\omega$ imply
\begin{equation}
\label{Pstability;lower-bd}
 \sum_{T \in \mesh_\omega} h_T^2 \| q_T \|_{L^2(T)}^2 + \sum_{F \in \faces_\omega} h_F \| q_F  \|_{L^2(F)}^2
 \Cleq
   \| \ell \|_{\V^+(\grid_\omega)^*}^2.
\end{equation}
Combing \eqref{Pstability;upper-bd} and \eqref{Pstability;lower-bd} shows that the ratio
$\|\ell\|_{H^{-1}(\omega)} / \|\ell\|_{\V^+(\grid_\omega)^*}$ for $\ell \in \F(\grid_\omega)$ is bounded by a universal constant depending on $d,\sigma,m_1,m_2$ and $\mathbb S^+$.

\smallskip\fbox{\scriptsize 3} It remains to complete the proof of the claimed identity for the operator norm. To this end, we first introduce an operator $Q_\omega: H^1_0(\omega) \to \V^+(\grid_\omega)$ by
\begin{equation*}
 \langle \ell, Q_\omega w \rangle
 =
 \langle \ell, w \rangle \quad \forall \ell \in \F(\grid_\omega).
\end{equation*}
As the one for $P_\omega$, this definition is well-posed because the pair $\big( \F(\grid_\omega), \V^+(\grid_\omega) \big)$ is nondegenerate for the dual pairing of $H^{-1}(\omega)$; cf.\ \eqref{dimV+mesh=dimFmesh} and \eqref{pairing-nondeg-on-Fmesh-V+mesh;2} of the proof of Lemma~\ref{L:algebra-of-P} (algebraic properties). By the Petrov-Galerkin character of the definition, $Q_\omega$ is a linear projection onto $\V^+(\grid_\omega)$. Given arbitrary $\ell \in H^{-1}(\omega)$ and $w \in H^1_0(\omega)$, the definitions of $Q_\omega$ and $P_\omega$ imply
\begin{equation*}
 \langle P_\omega \ell, w \rangle
 =
 \langle P_\omega \ell, Q_\omega w \rangle
 =
 \langle \ell, Q_\omega w \rangle,
\end{equation*}
that is $Q_\omega = P_\omega^*$ is the (Hilbert) adjoint to $P_\omega$.  In other words: the adjoint $P_\omega^*$ is a projection onto $\V^+(\grid_\omega)$. With this, we can prove the missing inequality. Let $\ell \in \F(\grid_\omega)$ be discrete. In fact,
\[
\langle \ell, w \rangle = \langle P_\omega \ell, w \rangle = \langle \ell, P_\omega^* w \rangle
\quad\Rightarrow\quad
\|\ell\|_{H^{-1}(\omega)} =
\sup_{w\in H^1_0(\omega)} \frac{\langle \ell, P_\omega^* w \rangle}{\|\nabla w\|_{L^2(\omega)}}
\]
leads to
\[
\|\ell\|_{H^{-1}(\omega)} \le \|\ell\|_{\V^+(\grid_\omega)^*} \|P_\omega^*\|_{\mathcal{L}(H^1_0(\omega))}
= \|\ell\|_{\V^+(\grid_\omega)^*} \| P_\omega\|_{\mathcal{L}(H^{-1}(\omega))}.
\]
This concludes the proof.
\end{proof}

\begin{remark}[failing global $H^{-1}$-stability]
\label{R:missing-global-H^-1-stabilty}
For Lebesgue norms, local stability of linear operators in terms of shape regularity entails that their respective global stability is uniform under mesh refinement. The fact that the first part of Lemma~\ref{L:loc-of-H^-1-norm} (localization of $H^{-1}$-norm) needs a condition to be true, may lead one to suspect that this implication might not be true in general for the $H^{-1}$-norm. This suspicion is confirmed by Example \ref{E:glb-instab-of-proj} below, where we show that $\| P_\grid \|_{\mathcal{L}(H^{-1}(\Omega))}$ can tend to $\infty$ under mesh refinement.  
\end{remark}

The proof of Lemma \ref{L:local-stability-of-P} provides all nontrivial ingredients to allow the approximate computation of $\| \ell \|_{H^{-1}(\omega)}$ whenever $\ell \in H^{-1}(\Omega)$ is discrete in $\omega$.

\begin{corollary}[quantifying $H^{-1}$-norms of discrete functionals]
\label{C:quantify-discrete-dual-norms}
Let $\omega \subset \Omega$ be a local subdomain meshed by $\grid_\omega$ and $\ell \in H^{-1}(\Omega)$ be discrete in $\omega$, given by the polynomials $q_T$ for $T \in \grid_\omega$, and $q_F$ for $F\in\faces_\omega$, where $F\in\faces_\omega$ are the interior faces in $\omega$. We then have
\begin{equation*}
    \| \ell \|_{H^{-1}(\omega)}^2
    \approx
    \sum_{T \in \grid_\omega} h_T^2 \| q_T \|_{L^2(T)}^2
    +
    \sum_{F \in \faces_\omega} h_F \| q_F \|_{L^2(F)}^2,
\end{equation*}
where the hidden constants depend on $d$, the shape regularity coefficient $\sigma$, the degrees $(m_1,m_2)$ of the discrete functionals, and the space $\mathbb{S}^+$ appearing in Assumption~\ref{A:abstract-cut-off}.
\end{corollary}
\begin{proof}  
This is a consequence of \eqref{Pstability;upper-bd} and \eqref{Pstability;lower-bd}, the
latter requiring $\ell$ to be discrete,
along with the fact that $\| \ell \|_{\V^+(\grid_\omega)^*} \le \| \ell \|_{H^{-1}\omega)^*}$
for any $\ell \in H^{-1}(\omega)$.
\end{proof}

\begin{corollary}[local near-best approximation]% of $P_\grid$]
\label{C:local-near-best-approx-of-P}
The projection $P_\grid$ yields local near-best approximations: for any functional $\ell \in H^{-1}(\Omega)$ and any local subdomain  $\omega$ meshed by $\grid_\omega$, we have
\begin{equation*}
    \| \ell - P_\grid \ell \|_{H^{-1}(\omega)}
    \leq
    C_\mathrm{lStb} \inf_{\chi \in \F(\grid_\omega)} \| \ell - \chi \|_{H^{-1}(\omega)},
\end{equation*}
where $C_\mathrm{lStb}$ is the constant of Lemma \ref{L:local-stability-of-P} (local $H^{-1}$ stability).
\end{corollary}

\begin{proof}
Fix a local subdomain $\omega$ meshed by $\grid_\omega$ and let $\chi \in \F(\grid_\omega)
$ be arbitrary. Thanks to Lemma \ref{L:algebra-of-P} (algebraic properties), we have $P_\omega \chi = \chi$ and 
\begin{equation*}
 \left( \ell - P_\grid\ell \right)|_{H^1_0(\omega)}
 =
 \left( I - P_\omega \right)\ell|_{H^{1}_0(\omega)}
 =
 \left( I - P_\omega \right) \left(\ell - \chi \right)|_{H^1_0(\omega)}.
\end{equation*}
As $P_\omega$ is a nontrivial projection on the Hilbert space $H^{-1}(\omega)$, \cite{Szyld:06} ensures
\begin{equation}
\label{local-boundedness-of-P-andI-P}
 \| I - P_\omega \|_{\mathcal{L}(H^{-1}(\omega))}
 =
 \| P_\omega \|_{\mathcal{L}(H^{-1}(\omega))}
 \leq
 C_\mathrm{lStb}.
\end{equation}
Hence
\begin{equation*}
  \| \ell - P_\grid \ell \|_{H^{-1}(\omega)}
  \leq
  C_\mathrm{lStb} \| \ell - \chi \|_{H^{-1}(\omega)}
\end{equation*}
concludes the proof because $\chi \in \F(\grid_\omega)$ is arbitrary.
\end{proof}

% \todo[inline]{RHN: We could use this estimate to get upper bounds for the oscillation term under various assumption on $f$. For instance, if $f\in L^2(\Omega)$ we may take $\chi$ with zero deltas and local $L^2$-projections onto polynomials of degree $n-1$ over each element. This would give the standard weighted $L^2$-oscillation. This is to be discussed in Section 8. AB: I think it would be good already here to mention what happens when $f \in L^2(\Omega)$. E.G: $P_\grid f$ is not necessarily in $L^2$ but thanks to the above result, we have the approximation property we expect. Same story when $f$ is a dirac on the skeleton, it would be good to discuss already here what we can say about $\| l - P_\grid l\|$. This would make $P_\grid$  more tangible.}

We illustrate the approximation of possible parts of the residual with the projection $P_\grid$ in a series of three remarks. For that purpose, the approximation quality is to be measured with a local $H^{-1}(\omega)$-norm, and it is instructive to compare with the operator $\Pi_\grid$ from \eqref{E:Pigrid}. %\eqref{splitting-residual-with-Pigrid}
Recall that the operator $\Pi_\grid$ is used implicitly in the standard approach (see Sect.\ \ref{S:std-res-est}) to approximate the discrete functionals $\F_\grid$. 

\begin{remark}[approximating functions]
\label{R:std-vs-new;approx-functions}
For functions, the local error with $P_\grid$ is uniformly dominated by the one with $\Pi_\grid$.
More precisely, if $m_2 \geq 0$ and $\ell \in H^{-1}(\Omega)$ satisfies $\langle \ell, w \rangle = \int_\Omega fw$ where $f \in L^p(\Omega)$ with $p > 2 d/ (2+d)$, then Corollary \ref{C:local-near-best-approx-of-P} (local near-best approximations) and $\Pi_\grid f  = \sum_{T \in \grid} (\Pi_T^{m_2} f) \chi_T \in \F(\grid_\omega)$ imply, for any local meshed subdomain $\omega$,
\begin{equation*}
 \| \ell - P_\grid \ell \|_{H^{-1}(\omega)}
 \Cleq
 \| f - \Pi_\grid f \|_{H^{-1}(\omega)}.
\end{equation*}
Observe that, although $\ell$ is a function, $P_\grid \ell$ is typically not a function. This property might look undesirable but it is crucial for an advantage of $P_\grid$ over $\Pi_\grid$ and closely related to the fact that the opposite inequality does not hold; cf.\ Remark~\ref{R:std-vs-new;stability} about stability.

% comparison with classical oscillation
Furthermore, supposing $f \in L^2(\Omega)$ and combining the preceding inequality with Lemma~\ref{L:Poincare} (first Poincar\'e inequality) gives
\begin{equation}
\label{bound-by-class-osc}
 \| \ell - P_\grid \ell \|_{H^{-1}(\omega_T)}^2
 \Cleq
 \sum_{T' \subset \omega_T} h_T^2 \| f - \Pi_{T'}^{m_2} f \|_{L^2(\omega_{T'})}^2,
\end{equation}
which establishes that the local $P_\grid$-oscillation of functions is uniformly dominated by its classical $\Pi_\grid$-counterpart but not vice versa.

In Sect.\ \ref{ss:data_f_Lp} this case is considered in the context of adaptive approximation.
\end{remark}

\begin{remark}[approximating admissible functionals]
\label{R:std-vs-new;admissible-functionals}
For functionals allowing for the application of $\Pi_\grid$, the local error with $P_\grid$ is again uniformly dominated by the one with $\Pi_\grid$.
In view of the previous remark, let us consider only $\ell \in H^{-1}(\Omega)$ such that $\langle \ell, w \rangle = \int_\Sigma gw$ where $g \in L^p(\Sigma)$ with $\Sigma \definedas \cup_{F\in\faces} F$ and $p > 2 (d-1)/d$. Note that we have again  $\Pi_\grid \ell \in \F(\grid_\omega)$ as $\langle \Pi_\grid \ell, w \rangle = \sum_{F \in \faces} \int_F  (\Pi_F^{m_1}g) \, w$ for all $w \in H^1_0(\Omega)$. Corollary~\ref{C:local-near-best-approx-of-P} (local near-best approximations) thus ensures, for any local meshed subdomain $\omega$,
\begin{equation*}
 \| \ell - P_\grid \ell \|_{H^{-1}(\omega)}
 \Cleq
 \| \ell - \Pi_\grid \ell \|_{H^{-1}(\omega)}.
\end{equation*}
Moreover, supposing $g \in L^2(\Sigma)$ and combining the scaled trace theorem \eqref{scaled-trace-thm-for_simplex-faces} with Lemma~\ref{L:Poincare} (first Poincar\'e inequality) yields
\begin{equation*}
 \| \ell - P_\grid \ell \|_{H^{-1}(\omega_T)}^2
 \Cleq
 \sum_{\substack{F \subseteq \overline{\omega}_T \\ F \not \subseteq \partial \omega_T} } h_F \| g - \Pi_F^{m_1} g \|_{L^2(F)}^2.
\end{equation*}
Also this case will be revisited in the context of adaptive approximation, namely in Sect.\ \ref{ss:data_f_div}.
\end{remark}

\begin{remark}[stability of approximation]
\label{R:std-vs-new;stability}
The error with $P_\grid$ is stable, while the one with $\Pi_\grid$ is not.
To see this by example, we restrict to $(m_1,m_2)=(0,0)$, fix some interior face $F \in \faces$ and, for $\eps>0$ sufficiently small, consider
\begin{align*}
 \langle \ell_\eps, w \rangle
 \definedas
 \int_\Omega f_\eps w
 =
 \frac{1}{2\eps} \int_{-\eps}^\eps \int_F w(y + s\vec{n}_F) \, dy \, ds,
\quad
 w \in H^1_0(\Omega),
\end{align*}
where $ f_\eps = (2\eps)^{-1} \chi_{F_\eps}$ is a multiple of the characteristic function of $ F_\eps \definedas \{ x + s\vec{n}_F \mid x \in F, -\eps < s < \eps \}$. As
\begin{align*}
 \langle \ell_\eps - \delta_F, w \rangle
 &=
 \frac{1}{2\eps} \int_{-\eps}^\eps \int_F \int_0^s \partial_{\vec{n}_F} w(y+t\vec{n}_F) \, dt \, dy \, ds
\\
 &\leq
 \frac{1}{2\eps} \int_{-\eps}^\eps \int_{F_\eps} | \nabla w(x) | \, dx \, ds
 \leq
 |F_\eps|^{1/2} \| \nabla w \|_{L^2(\Omega)},
\end{align*}
the functions $(\ell_\eps)_{\eps>0}$ tend to the proper functional $\delta_F$:
\begin{equation*}
 \ell_\eps \to \delta_F \text{ in }H^{-1}(\Omega).
\end{equation*}
Combining the convergence with the local stability of $P_\grid$, see \eqref{local-boundedness-of-P-andI-P}, yields
\begin{align*}
 \Big| \| \ell_\eps - P_\grid \ell_\eps \|_{H^{-1}(\omega_F)}
   &- \| \delta_F - P_\grid \delta_F \|_{H^{-1}(\omega_F)} \Big|
 \leq
 \| (I-P_\grid) (\ell_\eps - \delta_F) \|_{H^{-1}(\omega_F)}
\\
 &\Cleq
 \| \ell_\eps - \delta_F \|_{H^{-1}(\omega_F
 )}
\to 0,
\end{align*}
the stability of the error with $P_\grid$. Furthermore, since $P_\grid\delta_F=\delta_F$, the stability entails here
$\| \ell_\eps - P_\grid \ell_\eps \|_{H^{-1}(\omega_F)}\to0$. For $\Pi_\grid$ however, the approximation on the skeleton and in the volume are independent of each other. Hence, combining $\Pi_\grid \delta_F = \delta_F$, which follows from \eqref{E:Pigrid}, with $\lim_{\eps\to0} \Pi_T^0 f_\eps = |F| / ( 2|T| )$ for the two elements $T \in \grid$ containing $F$, leads to
\begin{equation*}
  \| \delta_F - \Pi_\grid \delta_F \|_{H^{-1}(\omega_F)}
  = 0 <
  \lim_{\eps\to0} \| \ell_\eps - \Pi_\grid \ell_\eps \|_{H^{-1}(\omega_F)}.
\end{equation*}
Measuring the error in weighted $L^2$-norms instead of the $H^{-1}$-norm results in a more dramatic instability. Indeed, denoting by $1_K$ and $0_K$ the constant functions on a simplex $K$ equal to $1$ or $0$, the two sides translate in
\begin{gather*}
 h_F^{\frac{1}{2}} \| 1_F - \Pi_F^0 1_F \|_{L^2(F)}
 +
 \sum_{T \in \grid_F} h_T \| 0_T - \Pi_T^0 0_T \|_{L^2(T)} = 0,
\\
 \lim_{\eps\to 0} \left(
  h_F^{\frac{1}{2}} \| 0_F - \Pi_F^0 0_F \|_{L^2(F)}
  +
  \sum_{T \in \grid_F} h_T \| f_\eps - \Pi_T^0 f_\eps \|_{L^2(T)}
\right) = \infty.
\end{gather*}
Note that such a transformation of volume contributions into contributions on the skeleton may occur by perturbation in the right-hand side or, in the opposite direction, by an improvement of the Galerkin approximation thanks to refinement.

In view of this instability of $\Pi_\grid$, the inequalities in the preceding Remarks~\ref{R:std-vs-new;approx-functions} and \ref{R:std-vs-new;admissible-functionals} cannot be reversed - a fact that can be inferred also from Remark~\ref{R:avoiding-overestimation}.

The above perturbations of $\delta_F$ are in the domain of $\Pi_\grid$. For the functionals
\begin{equation*}
 \langle \wh{\ell}_\eps, w \rangle
 \definedas
 \int_F w(y + \eps \vec{n}_F) \, dy,
\quad
 w \in H^1_0(\Omega)
\end{equation*}
however, it is not clear how to directly apply $\Pi_\grid$ for $\eps\neq0$. To the contrary, the approximations $P_\grid\wh{\ell}_\eps$ are defined and stable around $0$. Noteworthy, $P_\grid \wh{\ell}_\eps$ uses volume contributions to compensate for the displacement in the representation of the singular contribution. 
\end{remark}

\subsection{Discretized and oscillatory residual}
\label{S:dis-osc-res}
%
% intro
%AB I AM HERE
We now turn to the proper a~posteriori analysis, that is we shall derive \emph{upper} and \emph{lower bounds of the error}, implementing the following plan, which is motivated %and partially outlined
in Sect.\ \ref{S:discretize-functionals}. We use the projection $P_\grid$ onto discrete functionals $\F_\grid$ to split the residual into discretized and oscillatory parts. Then the quantification of the \emph{oscillatory residual} is reduced to \emph{data} oscillation through suitable choices of the degrees $(m_1,m_2)$ of the discrete functionals. The \emph{discretized residual} can be quantified in various manners; see Sections \ref{S:std-res-est} and \ref{S:other_estimators} below.

% splitting into discretized and oscillatory residual
We start by introducing indicators reflecting the announced splitting of the residual into discretized and oscillatory parts. They are vertex-indexed and, given $z \in \vertices$, defined by
\begin{equation}
\label{abstract_indicator}
\begin{aligned}
\index{Error Estimators!$\est^\text{abs}_\grid(z)$: abstract total estimator}
 \est^\text{abs}_\grid(z)^2
 &\definedas
\index{Error Estimators!$\eta^\text{abs}_\grid(z)$: abstract PDE estimator}
 \eta^\text{abs}_\grid(z)^2
 +
\index{Error Estimators!$\osc^\text{abs}_\grid(R_\grid,z)$: abstract oscillation}
 \osc_\grid(R_\grid,z)^2
 \quad\text{with}
\\
 \eta^\text{abs}_\grid(z)
 &\definedas \| P_\grid R_\grid \|_{H^{-1}(\omega_z)}
\text{ and }
 \osc_\grid(R_\grid,z)
 \definedas
 \| (I - P_\grid) R_\grid \|_{H^{-1}(\omega_z)}.
\end{aligned}
\end{equation}
Note that these quantities are not proper indicators: they still need to be quantified in a computable manner. By using `abs' (shorthand for `abstract'), we hint at the fact that $\eta^\mathrm{abs}_\grid$ can be quantified by various approaches.
\begin{lemma}[splitting of local residual norm]
\label{L:abstract_bds}
For any vertex $z\in\vertices$, the local residual norm is equivalent to the abstract indicator from \eqref{abstract_indicator}:
\begin{equation*}
  \frac{1}{\sqrt{2}C_\mathrm{lStb}} \est^\text{abs}_\grid(z)
  \leq
  \| R_\grid \|_{H^{-1}(\omega_z)}
  \leq
  \sqrt{2} \est^\text{abs}_\grid(z),
\end{equation*}
where $C_\mathrm{lStb}$ is the constant of Lemma \ref{L:local-stability-of-P} (local $H^{-1}$ stability).
\end{lemma}

\begin{proof}
As announced, we use the linear projection $P_\grid$ in order to split the residual into a discretized and an oscillatory part:
%In view of $-cu_\grid + \operatorname{div}(\bA\nabla u_\grid) \in \F_\grid = \operatorname{Im} P_\grid$, this leads to
\begin{equation}
\label{residual_splitting}
 R_\grid
 =
 P_\grid R_\grid + (I-P_\grid) R_\grid.
% =
% R^\text{apr}_\grid
% +
% \big( f - P_\grid f \big),
\end{equation}
% where $R^\text{apr}_\grid \definedas P_\grid f - cu_\grid + \operatorname{div}(\bA\nabla u_\grid)$ is a discrete approximation of the residual and $f - P_\grid f$ is an error in approximating $f$.
The upper bound of the local residual norm then readily follows from the triangle inequality:
\begin{equation*}
 \| R_\grid \|_{H^{-1}(\omega_z)}
 \leq
 \eta^\text{abs}_\grid(z) 
 +
 \osc_\grid(R_\grid,z) 
 \leq
 \sqrt{2} \est^\text{abs}_\grid(z).
\end{equation*}
To show the lower bound, we exploit the local stability of $P_\grid$, cf.\ \eqref{local-boundedness-of-P-andI-P}, to obtain
\begin{equation*}
  \eta^\text{abs}_\grid(z)
  =
  \| P_\grid R_\grid \|_{H^{-1}(\omega_z)}
  \leq
  C_\mathrm{lStb}  \| R_\grid \|_{H^{-1}(\omega_z)}
\end{equation*}
and
\begin{equation*}
 \osc_\grid(R_\grid,z)
 =
 \| (I-P_\grid)  R_\grid \|_{H^{-1}(\omega_z)}
 \leq
 C_\mathrm{lStb}  \| R_\grid \|_{H^{-1}(\omega_z)}.
\end{equation*}
Squaring both inequalities, summing them, and then taking the square root finishes the proof.
\end{proof}

In many a~posteriori analyses, this lemma is replaced by steps breaking a possible true equivalence between error and estimator. Therefore, the following remark points out the key ingredients.
\begin{remark}[ensuring proper equivalence]
\label{R:Ingredients-for-proper-equiv}
The fact that the projection $P_\grid$ and so also $I-P_\grid$ are linear and locally bounded operators precludes overestimation; see also Remark~\ref{R:avoiding-overestimation}. Comparing with Sect.\ \ref{S:std-res-est} and $\Pi_\grid$ in \eqref{splitting-residual-with-Pigrid}, we see that the local stability in $H^{-1}$ is crucial to that end and, in view of Remark~\ref{R:std-vs-new;stability}, requires discrete functionals with contributions on the skeleton.
\end{remark}

% reduction to data oscillation
Next, we want to simplify the residual oscillations $\osc_\grid(R_\grid,z)$, $z \in \vertices$, in the spirit of Remark~\ref{R:choosing-osc-degs;reduction-to-data-osc}. This will be dependent on the coefficients $\bA$ and $c$ of the differential operator and involve the following `polynomial degrees':
\begin{subequations}
\label{n_A-and-n_c}
\begin{gather}
 n_{\bA}
 \definedas
 \min \left\{ k \in \N_0 \mid 
  \bA \in \big( \mathbb{S}_\grid^{k,-1} \big)^{d \times d}
 \right\},
\\
  n_c
 \definedas
 \min \left\{ k \in \N_0 \cup \{-1\} \mid 
  c \in \mathbb{S}_\grid^{k,-1}
 \right\},
\end{gather}
\end{subequations}
where we use the convention $\min \emptyset = \infty$. We shall say that the differential operator $-\operatorname{div}(\bA \nabla \cdot) + c(\cdot)$ in \eqref{strong-form} \emph{has discrete coefficients} whenever $\max\{n_{\bA},n_c\} < \infty$, otherwise it has \emph{nondiscrete coefficients}.

\begin{lemma}[data oscillation reduction for discrete coefficients]
\label{L:data-oscillation-reduction}
If the coefficients $A$ and $c$ are discrete, the choices
\begin{align*}
 m_1 &= n_{\bA} + n - 1,
\\
 m_2 &= \max \{ n-2+n_{\bA}, \wt m_c \}
\quad\text{with}\quad
 \wt m_c = \begin{cases}
   n + n_c, &\text{if }c\neq0,
   \\ 0, &\text{otherwise}
  \end{cases} 
\end{align*}
ensure that the oscillatory residual reduces to data oscillation of the right-hand side:
\begin{equation*}
 (I-P_\grid) R_\grid
 =
 f - P_\grid f.
\end{equation*}
\end{lemma}

\begin{proof}
The choices for $m_1$ and $m_2$ yield, for any face $F \in \faces$ and any element $T \in \grid$,
\begin{equation*}
 \jump{\bA \nabla u_\grid} \cdot \vec{n}_F \in \P_{m_1}(F)
\quad\text{and}\quad
 \operatorname{div}(\bA \nabla u_\grid)|_{T} \in \P_{m_2}(T).
\end{equation*}
Furthermore, if $c\neq0$, we also have $cu_\grid{}|_{T} \in \P_{m_2}(T)$ and the claimed identity follows from $-\operatorname{div}(\bA\nabla u_\grid) + c u_\grid \in \F_\grid$.
\end{proof}

\begin{remark}[Poisson equation with linear elements]
\label{Poisson-with-Courant}
In the case of the Poisson equation with linear elements, the choices in Lemma \ref{L:data-oscillation-reduction} lead to $m_1=0$ and $m_2=0$. Alternatively, one may use $m_1=0$ and $m_2=-1$ (recall we have set $\mathbb{P}_{-1}(T)=\{0\}$), cf.\ \cite{DieningKreuzerStevenson:16} or \cite{Siebert.Veeser:2007-unilaterally-constrained-min-with-afem}. The choice here leads to an oscillation for which the standard oscillation indicators $h_T \| f - \Pi_T f \|_T$, $T \in \grid$, can be used as a surrogate; see also Remark \ref{R:data-osc-surrogates} about surrogates.
\end{remark}

% surrogate reduction for nondiscrete coefficients
If one of the coefficients, $\bA$ or $c$, is nondiscrete, the range of the finite element space $\V_\grid$ under the differential operator $-\operatorname{div}(\bA\nabla\cdot) + c(\cdot)$ consists of functionals whose densities are not piecewise polynomial. Consequently, the oscillatory residual cannot be reduced to the oscillation $f-P_\grid f$, or to any other oscillation of $f$ involving discrete functionals with piecewise polynomial densities. The next result illustrates the idea of a non-perfect remedy, namely bounding the residual oscillation defined in \eqref{abstract_indicator} in terms of data oscillation and discrete stability. For its formulation, we define global $H^{-1}$-oscillations by
\begin{equation}
\label{global-osc}
 \osc_\grid(\ell)^2
 \definedas
 \sum_{z\in\vertices} \osc_\grid(\ell,z)^2,
\quad
 \ell \in H^{-1}(\Omega),
\end{equation}
which, in contrast to $\| (I-P_\grid) \ell \|_{H^{-1}(\Omega)}$, is bounded in terms of $\ell$; cf.\ Remark~\ref{R:missing-global-H^-1-stabilty} (failing global $H^{-1}$-stability).

\begin{lemma}[surrogate data oscillation reduction]
\label{L:surrogate-osc}
Let
$%\begin{equation*}
 m_{\bA} \definedas \min\{ n_{\bA}, n-1 \}
$ and $%\quad\text{and}\quad
 m_{c} \definedas \min\{ n_c, n-1 \}
$ %\end{equation*}$
and define $m_1$ and $m_2$ as in Lemma~\ref{L:data-oscillation-reduction}, but replacing $n_{\bA}$ and $n_c$, respectively, with $m_{\bA}$ and $m_c$. Given any approximations
\begin{equation*}
 \wh{\bA} \in \big( \mathbb{S}^{m_{\bA},-1}_\grid \big)^{d \times d}
\quad\text{and}\quad
 \wh{c} \in \mathbb{S}^{m_c,-1}_\grid,
\end{equation*}
we then have, for all vertices $z \in \vertices$,
\begin{align*}
 \osc_\grid(R_\grid,z)
 &\leq
 \osc_\grid(f,z)
 +
 C_\mathrm{lStb} C(d,\sigma)\| \bA-\wh{\bA} \|_{L^\infty(\omega_z)} \| \nabla u_\grid \|_{L^2(\omega_z)}
\\
 &\quad +
 C_\mathrm{lStb} C(d,\sigma) \|h(c-\wh{c})\| _{L^\infty(\omega_z)} \| u_\grid \|_{L^2(\omega_z)}
\end{align*}
and thus
\begin{align*}
 \osc_\grid(R_\grid)^2
 &\leq
 3 \osc_\grid(f)^2
\\
 &+
 \frac{3(d+1)}{\alpha^2} \| f \|_{H^{-1}(\Omega)}^2 \left(
  \|\bA-\wh{\bA}\|_{L^\infty(\Omega)}^2
  %\max_{z\in\vertices} \|\bA-\wh{\bA}\|_{L^\infty(\omega_z)}^2
  +
  C_P^2 C(d,\sigma) \|h(c-\wh{c})\|_{L^\infty(\Omega)}^2
  \right)
\end{align*}
where $C_\mathrm{lStb}$ is the constant from Lemma~\ref{L:local-stability-of-P} (local $H^{-1}$-stability), $h$ the meshsize function defined by $h|_{T} = h_T$ for all $T\in\grid$,
%$C(d,\sigma)$ depends only on $d$ and the shape regularity coefficient $\sigma$,
$\alpha$ is the coercivity constant from \eqref{coercivity-and-continuity-constants-for_model-setting}, and $C_P$ is the constant in Lemma~\ref{L:Poincare} (first Poincar\'e inequality).
\end{lemma}

The bounds of Lemma~\ref{L:surrogate-osc} are obviously not convenient if $\bA$ or $c$ are not continuous. We therefore implement the underlying idea in Sect.\ \ref{S:general-data} differently.

\begin{proof}
\fbox{\scriptsize 1} To verify the local bound, let $z \in \vertices$ be any vertex. By linearity of $P_\grid$, we obtain
\begin{align*}
 \osc_\grid(R_\grid,z)
 & \leq
 \osc_\grid(f,z)
 \\& +
 \|(I-P_\grid) \big( -\operatorname{div}( \bA\nabla u_\grid ) \big)\|_{H^{-1}(\omega_z)}
+
 \|(I-P_\grid) (c u_\grid)\|_{H^{-1}(\omega_z)}
\end{align*}
and it remains to bound appropriately the two terms involving the coefficients $\bA$ and $c$. As the definitions of $m_1$ and $m_2$ ensure $-\operatorname{div}(\wh{\bA}\nabla u_\grid) \in \F_\grid$, Corollary~\ref{C:local-near-best-approx-of-P} (local near-best approximation), the scaled trace theorem \eqref{scaled-trace-thm-for_simplex-faces} and Lemma~\ref{L:Poincare} (first Poincar\'e inequality) give
\begin{align*}
 \|(I-P_\grid) &\big(-\operatorname{div}( \bA\nabla u_\grid ) \big)\|_{H^{-1}(\omega_z)}
 \leq
 C_\mathrm{lStb} \|-\operatorname{div}\big( (\bA-\wh{\bA}) \nabla u_\grid \|_{H^{-1}(\omega_z)}
\\
 &\lesssim
 C_\mathrm{lStb} \| \bA-\wh{\bA} \|_{L^\infty(\omega_z)} \| \nabla u_\grid \|_{L^2(\omega_z)}.
\end{align*}
As $\wh{c}u_\grid \in \F_\grid$ thanks to the definition of $m_2$, a similar argument using again Lemma~\ref{L:Poincare} and $\diam \omega_z \leq C h_z$ on $\omega_z$ provides
\begin{align*}
 \|(I-P_\grid) (cu_\grid)\|_{H^{-1}(\omega_z)}
 &\leq
 C_\mathrm{lStb} \|(c-\wh{c})u_\grid \|_{H^{-1}(\omega_z)}
\\
 &\leq
 C C_\mathrm{lStb} \| h(c-\wh{c}) \|_{L^\infty(\omega_z)}
  \| u_\grid \|_{L^2(\omega_z)},
\end{align*}
and the local bound is verified.

\smallskip\fbox{\scriptsize 2} To show the global bound, we square the local bound and sum it over all vertices $z \in \vertices$ to obtain
\begin{align*}
 \osc_\grid(R_\grid)^2
 &\leq
 3 \sum_{z\in\vertices} \osc_\grid(f,z)^2
 +
 3(d+1)C \|\bA-\wh{\bA}\|_{L^\infty(\Omega)}^2 \| \nabla u_\grid \|_{L^2(\Omega)}^2
\\
 & \quad+
 3(d+1) C \|h(c-\wh{c})\|_{L^\infty(\Omega)}^2  \| u_\grid \|_{L^2(\Omega)}^2.
\end{align*}
Hence, Lemma~\ref{L:Poincare} (first Poincar\'e inequality) on $\Omega$ and discrete stability,
\begin{equation*}
 \| u_\grid \|_{L^2(\Omega)}
 \leq
 C_P \| \nabla u_\grid \|_{L^2(\Omega)}
 \leq
 \frac{C_P}{\alpha} \| f \|_{H^{-1}(\Omega)},
\end{equation*}
finish the proof.
\end{proof}

The following remarks set Lemma~\ref{L:abstract_bds} (splitting of local residual norm) and the accompanying results Lemma~\ref{L:data-oscillation-reduction} and Lemma~\ref{L:surrogate-osc} on the reduction to data oscillation in the context of adaptive algorithms. 
\begin{remark}[structure of splitting]
\label{R:abstract-estimator-splitting}
Combing Lemma~\ref{L:abstract_bds} (splitting of local residual norms) with Lemma~\ref{L:data-oscillation-reduction} or Lemma~\ref{L:surrogate-osc} about reduction to data oscillation thus provides an  abstract estimator with the following two global parts:
\begin{equation*}
 \eta^\text{abs}_\grid(u_\grid)^2
 \definedas
 \sum_{z\in\vertices} \eta^\text{abs}_\grid(u_\grid,z)^2
\end{equation*}
and, writing $\data=(\bA,c,f)$ for the data of the partial differential equation,
\begin{align*}
\osc^\mathrm{abs}_\grid(\data)^2
 &\definedas \osc_\grid(f)^2,
 \\
 \osc^\mathrm{abs}_\grid(\data)^2
 &\definedas
 \osc_\grid(f)^2 + C_1 \max_{z\in\vertices} \|\bA-\wh{\bA}\|_{L^\infty(\omega_z)}^2
 + C_2 \max_{z\in\vertices} \|h(c-\wh{c})\|_{L^\infty(\omega_z)}^2,
\end{align*}
the latter provided $(\bA,c)$ are not discrete.
%\begin{equation*}
% \osc^\mathrm{abs}_\grid(\data)^2
% \definedas
% \begin{cases}
%  \osc_\grid(f)^2 \ \text{ or}
% \\
% \osc_\grid(f)^2 + C_1 \max_{z\in\vertices} \|\bA-\wh{\bA}\|_{L^\infty(\omega_z)}^2 
%  \\
%  \qquad {} + C_2 \max_{z\in\vertices} \|h(c-\wh{c})\|_{L^\infty(\omega_z)}^2.
% \end{cases}
%\end{equation*}
%
It is important to note the different nature of these two parts. The first part $\eta^\text{abs}_\grid(u_\grid)$, the abstract \emph{PDE indicator}, 
\begin{itemize}
\item is strictly related to the structure of the underlying PDE,
\item involves only discrete functionals from $\F_z$, and
\item the evaluation of its local indicators $\eta^\text{abs}_\grid(u_\grid,z)$ requires the global computation of the discrete solution $u_\grid$.
\end{itemize}
In contrast, the second part $\osc^\text{abs}_\grid(\data)$, the \emph{oscillation (indicator)},
\begin{itemize}
\item depends only on the data $\data$ of the differential operator,
\item involves non-discrete functionals, and
\item the evaluation of its local indicators $\osc_\grid(f,z)$, $\| \bA - \wh{\bA} \|_{L^\infty(\omega_z)}$, $\| h(c-\wh{c}) \|_{L^\infty(\omega_z)}$, $z\in\vertices$, is completely local.
\end{itemize}
The respective properties `discrete nature' and `local dependence' of the two parts are the key advantage over the whole local residual indicators $\| R_\grid \|_{H^{-1}(\omega_z)}$, $z \in \vertices$, and will be instrumental in the algorithmic design in Sections \ref{S:convergence-coercive} and \ref{S:conv-rates-coercive} below. 
\end{remark}

\begin{remark}[minimal regularity and regularizing $P_\grid$]
\label{R:no-extra-regularity}
It is worth noting that the results in this section do not involve any regularity beyond \eqref{setting-apost;nat-regularity} and that the projection $P_\grid$ has a regularizing effect. In particular, we have
\begin{equation*}
 \operatorname{Im} P_\grid
 =
 \F_\grid
 \subset
 H^{-\frac{1}{2}-\eps}(\Omega)
\quad
 \text{ for any small }\eps>0
\end{equation*}
thanks to the trace theorem in fractional Sobolev spaces. As a consequence, most techniques for a~posteriori error estimation can be directly applied to the discretized residual $P_\grid R_\grid$,  without any special twisting and under natural regularity assumptions. 
\end{remark}

\begin{remark}[reduction vs surrogate reduction]
\label{R:surrogate-reduction}
The kernel condition of Remark~\ref{R:avoiding-overestimation} (avoiding overestimation) is not verified for the bounds in Lemma~\ref{L:surrogate-osc} (surrogate data oscillation reduction). These bounds may thus exhibit overestimation and cannot be reversed. If we use the right-hand side of an overestimating bound as a part of an estimator, we shall call that part a \emph{surrogate}. This label marks a crucial difference between the cases represented by Lemma~\ref{L:surrogate-osc} (surrogate data oscillation reduction) and Lemma~\ref{L:data-oscillation-reduction} (data oscillation reduction), which is free of any overestimation.
\end{remark}

\begin{remark}[surrogate data oscillation]
\label{R:data-osc-surrogates}
Surrogates for data oscillation indicators can be useful in order to provide more direct access for computation. For example, if $f \in L^2(\Omega)$, the bound \eqref{bound-by-class-osc} by the classical $\Pi_\grid$-oscillation can be approximated by numerical integration. In such a context, it is useful to take the following points into account:
\begin{itemize}
\item Computable surrogates, i.e.\ computable upper bounds, for data oscillation indicators are in general impossible. In fact, generic data from an infinite-dimensional space will not be completely seen by the finite information available at any stage of a computation; cf.\ also \cite[Lemma~2 and Corollary~5]{KreuzerVeeser:2021} illustrating this fact for $\osc_\grid(f)$ with the help of orthogonality. Hence, computable surrogates will hinge on additional a~priori information on the given data. We postpone a discussion of examples to Section \ref{S:right-hand-side}.
\item As a general rule, surrogates should be applied last. This avoids that other parts of the estimators are also affected by overestimation; see Remark~\ref{R:std-vs-mod-res-est} (standard vs modified residual estimator) below.
\end{itemize}
\end{remark}

%
%-----------------------------------------------------------------------------------
\subsection{Modified residual estimation}
\label{S:mod-res-est}
%-----------------------------------------------------------------------------------
% intro

In view of the splitting into PDE and oscillation indicators and the discussion on the computability of the latter, it remains to quantify the abstract PDE indicators $\eta^\text{abs}_\grid(u_\grid,z)$, $z\in \vertices$.  To this end, we can employ Corollary~\ref{C:quantify-discrete-dual-norms} (quantifying $H^{-1}$-norms of discrete functionals), resulting in a modification of the standard residual estimator $\est^\mathrm{std}_\grid(u_\grid,\data)$ from Section \ref{S:std-res-est}. Alternative quantifications by other techniques of a~posteriori error estimation are discussed in Section \ref{S:other_estimators} below. Doing so, for simplicity, we consider only the case given by the following assumption.

\begin{assumption}[discrete coefficients and discrete functionals]
\label{A:discrete-coefficients}
\index{Assumptions!Discrete coefficients and discrete functionals}
\, Suppose that the coefficients $\bA$ and $c$ in \eqref{strong-form} are discrete and choose the degrees $(m_1,m_2)$ of the discrete functionals in $\F_\grid$ according to Lemma~\ref{L:data-oscillation-reduction} (data oscillation reduction for discrete coefficients).
\end{assumption}
For nondiscrete coefficients, one essentially has to invoke  Lemma~\ref{L:surrogate-osc} (surrogate data oscillation reduction) instead of Lemma~\ref{L:data-oscillation-reduction} in order to reduce to data oscillation.

% re-indexing by elements for bisection
\smallskip We shall employ the bisection method in order to refine the mesh. Since this method is based upon the subdivision of elements, it is convenient to split the estimator into contributions associated with elements and not with vertices as in Sect.\ \ref{S:dis-osc-res}.

% definition of modified standard residual estimator 
To define the modified residual estimator, we 
recall the representation \eqref{def-densities-of-P} of the $H^{-1}$-projection $P_\grid$, and we use Assumption \ref{A:discrete-coefficients} (discrete coefficients and discrete functionals) to set
\begin{equation}
\label{mod-res-est}
\begin{aligned}
\index{Error Estimators!$\est_\grid$: total estimator}
 \est_\grid^2
 \definedas
 \sum_{T\in\grid} \est_\grid(T)^2
\quad&\text{with}\quad
\\ 
\index{Error Estimators!$\est_\grid(T)$, $\est_\grid(u_\grid,f,T)$: local total estimator}
 \est_\grid(T)^2
 \definedas
 \est_\grid(u_\grid,f,T)^2
 &\definedas
 \eta_\grid(u_\grid,T)^2 + \osc_\grid(f,T)^2,
%\text{ with }
\\
\index{Error Estimators!$\eta_\grid(u_\grid,T)$: PDE local estimator}
 \eta_\grid(u_\grid,T)^2
 &\definedas
  h_T \sum_{F \subset \partial T \setminus \partial\Omega}
  \| \jump{\bA \nabla u_\grid} \cdot \vec{n}_F - P_F f \|_{L^2(F)}^2
\\
 &\qquad
  + h_T^2  \| P_T f - c u_\grid + \div(\bA\nabla u_\grid) \|_{L^2(T)}^2,
\\
\index{Error Estimators!$\osc_\grid(f,T)$, $\osc_\grid(f,T)_{-1}$: local oscillation for the load function}
 \osc_\grid(f,T)^2
 &\definedas
% \osc_\grid(f,T)_{-1}^2
% \definedas
  \| f - P_\grid f \|_{H^{-1}(\omega_T)}^2. 
\end{aligned}
\end{equation}
Clearly, this is a variant of the standard residual estimator in \eqref{std-res-est}, where the main differences are given by the corrections $P_Ff$, $F \in \faces$, of the jump residual and the replacement of $f|_{T}$ by $P_Tf$, $T \in \grid$, in the PDE indicator.  As shown by the following theorem and remarks, the modification leads to more accurate a~posteriori bounds.

\begin{theorem}[modified residual estimator]
\label{T:modified-estimator}
Under Assumption~\ref{A:discrete-coefficients}, the modified residual estimator \eqref{mod-res-est} is equivalent to the error: more precisely, we have
\begin{equation*}
\index{Constants!$(C_L,C_U)$: a-posteriori lower and upper bounds constants}
 C_L \est_\grid
 \leq
 \| \nabla (u - u_\grid) \|_{L^2(\Omega)}
 \leq
 C_U \est_\grid, %\left(
%  \sum_ {T \in \grid} \est_\grid(T)^2
% \right)^{\frac{1}{2}}
\end{equation*}
%and, for any element $T \in \grid$,
%\begin{equation*}
% C_L \est_\grid(T)
% \leq
% \| \nabla(u_\grid - u) \|_{L^2(\omega_T)}.
%\end{equation*}
where the constants $C_U \geq C_L > 0$ depend only on the coefficients $(\bA,c)$, the shape regularity coefficient $\sigma$ from \eqref{E:shape-regularity}, the polynomial degree $n$, and $d$.
\end{theorem}

\begin{proof}
% upper bound
To derive the upper bound, we use the ones of Lemma~\ref{L:err-residual} (error and residual), Corollary~\ref{C:residual-loc} (star localization of residual norm), Lemma~\ref{L:abstract_bds} (splitting of local residual norms), Corollary~\ref{C:quantify-discrete-dual-norms} (quantifying $H^{-1}$-norms of discrete functionals) with stars and obtain
\begin{align*}
 \| \nabla (u-u_\grid) \|_{L^2(\Omega)}^2
 &\Cleq
 \| R_\grid \|_{H^{-1}(\Omega)}^2
 \Cleq
 \sum_{z \in \vertices} \| R_\grid \|_{H^{-1}(\omega_z)}^2
\\
 &\Cleq
 \sum_{z \in \vertices} \est^\text{abs}_\grid(z)^2
 =
 \sum_{z \in \vertices} \eta^\text{abs}_\grid(u_\grid,z)^2
  + \sum_{z \in \vertices} \osc_\grid(f,z)^2
\\
  &\Cleq
 \sum_{z \in \vertices} \sum_{T \in \grid_z} \eta_\grid(u_\grid,T)^2
  + \sum_{z \in \vertices} \osc_\grid(f,z)^2
\end{align*}
with $\grid_z = \{ T\in\grid \mid T \ni z\}$. As a given mesh element appears in the star meshes $\grid_z$ for at most $d+1$ vertices, we have
\begin{equation*}
 \sum_{z \in \vertices} \sum_{T \in \grid_z} \eta_\grid(u_\grid,T)^2
 \leq
 (d+1) \sum_{T \in \grid} \eta_\grid(u_\grid,T)^2
\end{equation*}
for the first sum and Lemma~\ref{L:local-re-indexing} (localization re-indexing) yields
\begin{equation*}
 \sum_{z \in \vertices} \osc_\grid(f,z)^2
 \Cleq
 \sum_{T \in \mesh} \osc_\grid(f,T)^2
\end{equation*}
for the second sum. Inserting the last two inequalities in the previous one, we conclude the upper bound:
\begin{equation*}
 \| \nabla (u-u_\grid) \|_{L^2(\Omega)}^2
 \Cleq
 \sum_{T \in \grid} \eta_\grid(u_\grid,T)^2 + \sum_{T \in \mesh} \osc_\grid(f,T)^2
 =
 \sum_ {T \in \grid} \est_\grid(T)^2.
\end{equation*}

% lower bounds
% AB I AM HERE
To show the lower bound, fix a mesh element $T \in \grid$. Applying the local lower bounds in Corollary~\ref{C:quantify-discrete-dual-norms}, Lemma~\ref{L:local-stability-of-P} (local $H^{-1}$ stability) and Lemma \ref{L:err-residual} on the local meshed subdomain $\wt{\omega}_T$ defined in \eqref{tildeomegaT} yields for the PDE indicator
\begin{equation}
\label{lower-bd;PDE-ind}
 \eta_\grid(u_\grid,T)
 \Cleq
% \| P_\grid R_\grid \|_{\V^+(\wt{\grid}_T)'}
% \leq
 \| P_\grid R_\grid \|_{H^{-1}(\wt{\omega}_T)}
 \Cleq
 \| R_\grid \|_{H^{-1}(\wt{\omega}_T)}.
% \leq
% \| \nabla(u - u_\grid) \|_{L^2(\wt{\omega}_T)}.
\end{equation}
In the case of the oscillation indicator, we exploit $-cu_\grid + \operatorname{div}(\bA\nabla u_\grid) \in \F_\grid$ with the help of Lemma~\ref{L:algebra-of-P} (algebraic properties) and apply Lemma~\ref{L:local-stability-of-P} on the local meshed subdomain $\omega_T$:
\begin{equation*}
%\begin{aligned}
 \osc_\grid(f,T)
 =
 \| f - P_\grid f \|_{H^{-1}(\omega_T)}
 =
 \| (I-P_\grid) R_\grid \|_{H^{-1}(\omega_T)}
\\
 \Cleq
 \| R_\grid \|_{H^{-1}(\omega_T)}.
% \Cleq
% \| \nabla (u-u_\grid) \|_{L^2(\omega_T)}.
%\end{aligned}
\end{equation*}
Thanks to $\wt{\omega}_T \subset \omega_T$, combining the last two inequalities gives the desired local lower bound
\begin{equation}\label{e:equivalence_tmp}
\begin{aligned}
 \est_\grid(T)^2
 &=
 \eta_\grid(u_\grid,T)^2 + \osc_\grid(f,T)^2
\\
 &\Cleq
 \| R_\grid \|_{H^{-1}(\wt{\omega}_T)}
 +
 \| R_\grid \|_{H^{-1}(\omega_T)}
 \Cleq
 \| R_\grid \|_{H^{-1}(\omega_T)}.
% \| \nabla(u - u_\grid) \|_{L^2(\wt{\omega}_T)}^2
%  +  \| \nabla (u-u_\grid) \|_{L^2(\omega_T)}^2
% \Cleq
% \| \nabla (u-u_\grid) \|_{L^2(\omega_T)}^2.
\end{aligned}
\end{equation}
As the number of patches $\omega_T$, $T \in \grid$, containing a given mesh element is uniformly bounded by $d$ and the shape regularity coefficient $\sigma$, summing this bound over all mesh elements yields the global lower bound
\begin{equation*}
 \sum_{T \in \grid} \est_\grid(T)^2
 \Cleq
 \sum_{T \in \grid} \| R_\grid \|_{H^{-1}(\omega_T)}^2
 \Cleq
 \| \nabla(u-u_\grid) \|_{L^2(\Omega)}^2
\end{equation*}
with the help of Lemma~\ref{L:loc-of-H^-1-norm} (localization of $H^{-1}$ norm) and  Lemma~\ref{L:err-residual}. Thus the equivalence of error and estimator is established.
\end{proof}

% \todo[inline]{RHN: These estimates could be used to present and discuss a one-step AFEM with $H^{-1}$-oscillation and study its convergence and optimality properties. We could also replace the $H^{-1}$-oscillation by an
% $L^2$-weighted oscillation and get the standard AFEM. In both cases we would be dealing with discrete coefficients. This would justify the two-step AFEM as an algorithm that takes variable coefficients $(A,c)$ into account.}

A detailed comparison of the modified residual estimator with the standard one is in order.
\begin{remark}[modified vs standard residual estimator]
\label{R:std-vs-mod-res-est}
We compare the modified residual estimator \eqref{mod-res-est} with the standard one given by \eqref{std-res-est;tot} and the local split indicators \eqref{loc-ind-of-std-res-est;alt}. As a common characterizing feature, both residual estimators use properly scaled $L^2$-norms of jump and element residual, ready for numerical integration. However, we observe the following differences:
\begin{itemize}
\item while the modified estimator $\est_\grid$ is defined under the natural regularity assumptions \eqref{setting-apost;nat-regularity}, the standard estimator $\est^\mathrm{std}_\grid$ requires $\bA \in W^{1}_{\infty}(\Omega;\mathbb R^{d\times d})$ and $f \in L^2(\Omega)$ in addition;
\item while the modified estimator $\est_\grid$ is truly equivalent to the error, the standard estimator may may overestimate it, limited however by Proposition~\ref{P:lw-bd-with-std-res-est} (partial lower bound).
\end{itemize}
By the domain test in Remark~\ref{R:avoiding-overestimation} (avoiding overestimation), we know that these two points are interrelated. However, also the kernel test is at play in the overestimation. Indeed, revisiting the proof of Theorem~\ref{T:up-bd-with-std-res-est} (upper bound with standard residual estimator), we can replace the scaled $L^2$-norms of the element residuals %$\sum_{T\subset\omega_z} h_T^2 \| r \|_{L^2(T)}^2$
on a star $\omega_z$ by $\| r \|_{H^{-1}(\omega_z)}$ and the resulting vertex-oriented variant of the residual estimator with unsplit local indicators is defined for all $f \in H^{-1}(\Omega)$. Overestimation can, however, still occur non-asymptotically as well as asymptotically; cf.\ \cite{CohenDeVoreNochetto:2012}. Indeed, in the case of the Poisson equation and linear finite elements, the kernel test is obviously not satisfied. This shows that the splitting in jump and element residual is quite delicate and highlights the crucial role of the modifications of the standard residual estimator: not only do they allow for stability in line with Remark~\ref{R:std-vs-new;stability} (stability of approximation) but also imply the kernel test.

To conclude this comparison, let us illustrate the second point of Remark~\ref{R:data-osc-surrogates} (surrogate data oscillation), namely that surrogates should be applied last. Using \eqref{bound-by-class-osc} in Remark~\ref{R:std-vs-new;approx-functions} (approximating functions), we may replace in the modified residual estimator the $H^{-1}$-oscillation $\osc_\grid(f)$ by the standard oscillation $\osc^\mathrm{std}_\grid(f)$, which can be readily approximated with numerical integration. Doing so, we first split the residual with $P_\grid$ and then apply $\Pi_\grid$ to obtain the surrogate. Note however that, if we apply $\Pi_\grid$ earlier to split the residual, the crucial modifications will not appear and, therefore, also the PDE indicator of the standard residual estimator exhibits overestimation. 
\end{remark}

%-------------------------------------------------------------------------------------
\subsection{Bounds for corrections and reduction of PDE estimator}
\label{S:bds-for-corr}
%-------------------------------------------------------------------------------------
In the following sections we shall use the modified residual estimator $\est_\grid$ from \eqref{mod-res-est} in adaptive algorithms. In their convergence analyses, not only its relationship with the error is important, but also its relationship with the norm $\| \nabla (u_{\grid_*}-u_\grid)\|_{L^2(\Omega)}$ of (possible) \emph{corrections}, where $u_{\grid_*}$ is the Galerkin approximation %in $\V_{\grid_*}$
to $u$ over some refinement $\grid_*$ of $\grid$. This section establishes corresponding upper and lower bounds, as well as related results about the \emph{global PDE indicator}
\begin{equation}
\label{global-PDE-estimator}
\index{Error Estimators!$\eta_\grid(u_\grid,f)$: PDE estimator}
\eta_\grid(u_\grid,f)^2
 \definedas
 \sum_{T\in\mesh} \eta_\grid(u_\grid,T)^2
\end{equation}
and the \emph{global oscillation}
\begin{equation}
\label{global-osc-with-elms}
\index{Error Estimators!$\osc_\grid(f)$, $\osc_\grid(f)_{-1}$: oscillation for the load function}
 \osc_\grid(f)^2
 \definedas
 \sum_{T \in \grid} \osc_\grid(f,T)^2.
\end{equation}
When it is important to indicate that the oscillations are measured in $H^{-1}$, we use the notation
$$
\index{Error Estimators!$\osc_\grid(f)$, $\osc_\grid(f)_{-1}$: oscillation for the load function}
\index{Error Estimators!$\osc_\grid(f,T)$, $\osc_\grid(f,T)_{-1}$: local oscillation for the load function}
\osc_\grid(f)_{-1} \qquad  \textrm{and} \qquad \osc_\grid(f,T)_{-1}.
$$

\smallskip Let $\grid_*$ be a conforming mesh that is a \emph{refinement} of $\mesh$, i.e.\ for any element $T \in \mesh$, there exists a submesh $\grid_{*,T}$ of $\grid_*$ such that $T = \cup \{T_* : T_* \in \grid_{*,T}\}$.
The Galerkin approximation in $\V_{\grid_*}$ is characterized by
\begin{equation*}
 u_{\grid_*} \in \V_{\grid_*}:
\quad
 \B[u_{\grid_*},w] = \langle f, w \rangle \quad\forall w \in \V_{\grid_*}.
\end{equation*}
Hence, the discrete solution $u_\grid$ on the original mesh $\grid$ is not only a Galerkin approximation to the exact solution $u$ satisfying \eqref{weak-form} but also to $u_{\grid_*}$. The norm $\| \nabla(u_{\grid_*}-u_\grid) \|_{L^2(\Omega)}$ of the correction therefore can be viewed as the error in approximating $u_{\grid_*}$ on the mesh $\grid$. This viewpoint suggests considering the variant
\begin{equation}
\label{res-corr-id}
 \langle R_\grid, w \rangle
 =
 \B[u_{\grid_*} - u_\grid, w] \quad  \forall w \in \V_{\grid_*}
\end{equation}
of the error-residual identity \eqref{err-res-id} and introducing the discrete dual norm 
\begin{equation}
\label{residual-norm-of-corr}
 \| R_\grid \|_{(\V_{\grid_*})^*}
 \definedas
 \sup_{w \in \V_{\grid_*}} \frac{\langle R_\grid, w \rangle}{\| \nabla w \|_{L^2(\Omega)}}
\end{equation}
of the residual as a counterpart of $\| R_\grid \|_{H^{-1}(\Omega)}$. Arguing as in the proof of Lemma~\ref{L:err-residual} (error and residual), we thus readily obtain the following quantitative relationship between the correction and the residual.

\begin{lemma}[correction and residual]
\label{L:corr-res}
If $\grid_*$ is a refinement of the mesh $\grid$, the norm of the correction $u_{\grid_*}-u_\grid$ is equivalent to the discrete residual norm. More precisely,
\begin{equation*}
 \frac{1}{\|\B\|} \| R_\grid \|_{(\V_{\grid_*})^*}
 \leq
 \| \nabla(u_{\grid_*}-u_\grid) \|_{L^2(\Omega)}
 \leq
 \frac{1}{\alpha} \| R_\grid \|_{(\V_{\grid_*})^*}
\end{equation*}
where $\| \B \| \geq \alpha >0$, are, respectively, the continuity and coercivity constant of the bilinear form $\B$.
\end{lemma}

% towards sharp upper bd of correction
We first exploit the upper bound in Lemma~\ref{L:corr-res}. As the inclusion $\V_{\grid_*} \subset H^1_0(\Omega)$ implies
\begin{equation}
\label{discrete-smaller-that-org-dualnorm}
 \| R_\grid \|_{(\V_{\grid_*})^*}
 \leq
 \| R_\grid \|_{H^{-1}(\Omega)},
\end{equation}
Theorem~\ref{T:modified-estimator} (modified residual estimator) immediately yields the upper bound
\begin{equation}
\label{simple-upper-bd-for-corr}
 \| \nabla (u_{\grid_*}-u_\grid) \|_{L^2(\Omega)}
 \leq
 C_U \est_\grid(u_\grid,f).
\end{equation}
This bound however appears to be not accurate in view of the use of \eqref{discrete-smaller-that-org-dualnorm}. We shall sharpen it by following the line of its proof but exploiting the \emph{full} orthogonality
\begin{equation}
\label{full-ortogonality-of-residual;corr}
 \langle R_\grid, w \rangle
 = 0 \quad  \forall w \in \V_\grid
\end{equation}
with suitably tuned Scott-Zhang interpolation \cite{ScottZhang:90}.

% tuned Scott-Zhang interpolation
In order to prepare the use of this interpolation, denote by $\mathcal{N}$ the Lagrange nodes of order $n$ of the mesh $\mesh$ and by $\faces$ and $\faces_*$, respectively, the $(d-1)$-dimensional faces of $\grid$ and $\grid_*$, including boundary ones. Given a node $z \in \nodes$, fix a face $F_z \in \faces$ such that $F_z$ contains $z$ and the following conditions are met:
%\begin{subequations}
\begin{gather*}
 z \in \partial\Omega \implies F_z \subset \partial\Omega,
\\
 \{ F \in \faces \cap \faces_* \mid F \ni z \} \neq \emptyset
 \implies
 F_z \in \faces_*.
\end{gather*}
%\end{subequations}
Furthermore, denote by $\psi_z^*$ the polynomial in $\P_n(F_z)$ satisfying
\begin{equation*}
 \forall y \in \nodes
\quad
 \int_{F_z} \psi_z^* \psi_y = \delta_{yz},
\end{equation*}
where $\{\psi_y\}_{y\in\nodes}$ is the Lagrange basis of ${\mathbb S}^{n,0}_\grid$, and define
\begin{equation}\label{d:SZ}
 I_\grid w
 =
 \sum_{z\in\nodes} \left( \int_{F_z} \psi_z^* w \right) \psi_z.
\end{equation}
The two conditions on the fixed face $F_z$ then ensure, respectively,
\begin{subequations}
\label{approx-prop-of-SZ}
\begin{align}
 w \in H^1_0(\Omega) &\implies I_\grid w \in \V_\grid,
\\ \label{local-invariance-of-tuned-Scott-Zhang}
 w \in \V_{\grid_*} \text{ and } T \in \grid \cap \grid_*
 &\implies
 I_\grid w = w \text{ on }T.
\end{align}
\end{subequations}
In particular, if $w \in \V_{\grid_*}$, its approximation $I_\grid w \in \V_\grid$ is an admissible test function and coincides with $w$ whenever possible. Finally, $I_\grid$ has the following stability and approximation properties, where the hidden constants depend only on $d$, $n$, and the shape regularity coefficient $\sigma$: for any element $T \in \grid$ and any face $F \in \faces$, % with $F \subset T$,
\begin{subequations}
\label{SZ-estimates}
\begin{align}
\label{SZ-H1-stability}
 \| \nabla I_\grid w \|_{L^2(T)}
 &\Cleq
 \| \nabla w \|_{L^2(\omega_T)},
\\ \label{SZ-error-estimate-on-elm}
 \| w - I_\grid w \|_{L^2(T)}
 &\Cleq
 h_T \| \nabla w \|_{L^2(\omega_T)},
\\ \label{SZ-error-estsimate-on-face}
 \| w - I_\grid w \|_{L^2(F)}
 &\Cleq
 h_F^{\frac{1}{2}} \| \nabla w \|_{L^2(\omega_F)}.
\end{align}
\end{subequations}

% proper upper bound for correction
The sharpening of the simple upper bound \eqref{simple-upper-bd-for-corr} lies in the fact that only a part of the estimator in \eqref{mod-res-est} will be invoked. To formulate this, we define
\begin{equation}
\label{estimator-of-submesh}
 \est_\grid(u_\grid,f,\wt{\grid})
 \definedas
 \left(
  \sum_{T\in \wt{\grid}} \est_\grid(u_\grid,f,T)^2
 \right)^{1/2}
\end{equation}
%I AM HERE
where $\wt \grid \subset \grid$ is a subset of elements in $\grid$. In the same vein, we shall denote $\eta_\grid(u_\grid,f,\wt \grid)$ and $\osc_\grid(f,\wt \grid)$. 

\begin{theorem}[upper bound for corrections]
\label{T:ubd-corr}
Let Assumption~\ref{A:discrete-coefficients} hold and let $\grid_*$ be a refinement of the mesh $\grid$. The correction $u_{\grid_*} - u_\grid$ is bounded in terms of the indicators of the refined elements $\grid \setminus \grid_*$:
\begin{equation*}
\index{Constants!$\wt{C}_U$: localized upper bound constant}
 \| \nabla(u_\grid - u_{\grid_*}) \|_{L^2(\Omega)}
 \leq
 \wt{C}_U \est_\grid \big(u_\grid, f, \grid \setminus \grid_* \big),
\end{equation*}
where the constant $\wt{C}_U>0$ depends only on the dimension $d$, the coefficients $\bA$ and $c$, the polynomial degree $n$, and the shape regularity coefficient $\sigma$ from \eqref{E:shape-regularity}.
\end{theorem}
    
\begin{proof}
\fbox{\scriptsize 1} \, {\it Localization and splitting of the residual norm.} In light of Lemma~\ref{L:corr-res}, it suffices to bound the discrete residual norm $\| R_\grid \|_{(\V_{\grid_*})^*}$. Given $w \in \V_{\grid_*}$, we prepare the localization of the residual by full orthogonality \eqref{full-ortogonality-of-residual;corr} and split it with help of the projection $P_\grid$ on discrete functionals:
\begin{align*}
 |\langle R_\grid, w \rangle|
 &=
 |\langle R_\grid, w - I_\grid w \rangle|
\\
 &\leq
 |\langle P_\grid R_\grid, w - I_\grid w \rangle|
  +
% \langle R - P_\grid R, w - I_\grid w \rangle|
%\\
% &=
% |\langle P_\grid f  - cu_\grid + \operatorname{div}(\bA \nabla u_\grid), w - I_\grid w \rangle|
%  +
 |\langle f - P_\grid f, w - I_\grid w \rangle|,
\end{align*}
where we used the identity $R-P_\grid R = f - P_\grid f$ in the last step. In light of \begin{equation*}
 \est_\grid(u_\grid,f,\grid \setminus \grid_*)^2
 =
 \eta_\grid(u_\grid,f,\grid \setminus \grid_*)^2
  +
 \osc_\grid(f,\grid \setminus \grid_*)^2,
\end{equation*}
it remains to bound the two terms with discretized residual $P_\grid R_\grid$ and the oscillation of $f$ appropriately.

\fbox{\scriptsize 2} {\it Bounding the discretized residual.}  We adopt the notation \eqref{def-densities-of-P} for the densities of $P_\grid$, and
exploit the piecewise nature of the discretized residual and the local invariance \eqref{local-invariance-of-tuned-Scott-Zhang} of $I_\grid$ to deduce
\begin{multline*}
 \langle P_\grid R_\grid, w - I_\grid w \rangle
\\
 =
 \sum_{T\in\grid\setminus\grid_*} 
 \left(\int_T (P_T R_\grid) (w -I_\grid w)
 + \frac{1}{2}
 \sum_{F \subset \partial T \setminus \partial\Omega} \int_F (P_F R_\grid) (w -I_\grid w) \right).
\end{multline*}
Invoking the local approximation properties \eqref{SZ-error-estimate-on-elm}  and \eqref{SZ-error-estsimate-on-face} of $I_\grid$ leads to the desired bound for the discretized residual:
\begin{align*}
 |\langle P_\grid R_\grid, w - I_\grid w \rangle|
 &\Cleq
 \sum _{T \in \grid \setminus \grid_*} \eta_\grid(u_\grid,f,T) \, \| \nabla w \|_{L^2(\omega_T)}
\\
 &\Cleq
 \eta_\grid \big( u_\grid,f,\grid \setminus \grid_* \big) \, \| \nabla w \|_{L^2(\Omega)}.
\end{align*}

\fbox{\scriptsize 3} {\it Bounding the oscillation.} We need to split the oscillation into suitable local contributions and first proceed similarly to the proof of Lemma~\ref{L:loc-of-H^-1-norm} (localization $H^{-1}$ norm) (i). Writing
\begin{equation*}
 w - I_\grid w
 =
 \sum_{z \in \vertices} (w-I_\grid w) \phi_z
\quad\text{and}\quad
 \Omega_0 \definedas \bigcup_{T \in \grid \setminus \grid_*} T,
\end{equation*}
we have $(w- I_\grid w)\phi_z \in H^1_0(\omega_z \cap \Omega_0)$ thanks to \eqref{local-invariance-of-tuned-Scott-Zhang} and, for any $T \subset \omega_z \cap \Omega_0$
\begin{align*}
 \| \nabla \big( (w-I_\grid w)\phi_z \big) \|_{L^2(T)}
 &\leq
 \| \phi_z \nabla (w-I_\grid w) \|_{L^2(T)}
 +
 \| (w-I_\grid w) \nabla \phi_z \|_{L^2(T)}
\\
&\leq
 \| \nabla (w-I_\grid w) \|_{L^2(T)}
 +
 C(d) \sigma \| \nabla w \|_{L^2(\omega_T)}
\\
 &\Cleq
 \| \nabla w \|_{L^2(\omega_T)}
\end{align*}
by means of $0\leq\phi_z\leq1$, $|\nabla\phi_z| \leq C(d) \sigma h_T^{-1}$, \eqref{SZ-H1-stability}, and \eqref{SZ-error-estimate-on-elm}. Hence, we get 
\begin{align*}
 |\langle f - P_\grid f &, w - I_\grid w \rangle|
 \leq
 \sum_{z\in\vertices} |\langle f - P_\grid f, (w - I_\grid w)\phi_z \rangle|
\\
 &\leq
 \sum_{z\in\vertices} \| f - P_\grid f \|_{H^{-1}(\omega_z \cap \Omega_0)}
  \| \nabla \big( (w - I_\grid w)\phi_z \big) \|_{L^2(\omega_z \cap \Omega_0)}
\\
 &\Cleq
 \sum_{z\in\vertices} \| f - P_\grid f \|_{H^{-1}(\omega_z \cap \Omega_0)}
  \| \nabla w \|_{L^2(\cup_{T \subseteq \omega_z \cap \Omega_0} \omega_T)}
\\
 &\Cleq
 \left( \sum_{z\in\vertices} \| f - P_\grid f \|_{H^{-1}(\omega_z \cap \Omega_0)}^2 \right)^{1/2}
  \| \nabla w \|_{L^2(\Omega)} .
\end{align*}
Since 
\begin{equation*}
  \sum_{z\in\vertices} \| f - P_\grid f \|_{H^{-1}(\omega_z \cap \Omega_0)}^2
  \leq
  \sum_{T \in \grid \setminus \grid_*} \| f - P_\grid f \|_{H^{-1}(\omega_T)}^2,
\end{equation*}
the oscillation of $f$ is therefore bounded by
\begin{equation*}
 |\langle f - P_\grid f, w - I_\grid w \rangle|
 \leq
 \osc_\grid(f,\grid \setminus \grid_*) \| \nabla w \|_{L^2(\Omega)}
\end{equation*}
and the proof is complete. 
\end{proof}

% lower bounds by interior vertices
Proposition~\ref{P:lw-bd-with-std-res-est} (partial lower bound) as well as Lemma~\ref{L:abstract_bds} (splitting of local residual norm) illustrate that the test space $\V_\grid^+$ is closely related to lower bounds for the error. This observation suggests establishing lower bounds for the correction $\| \nabla ( u_\grid - u_{\grid_*} ) \|_{L^2(\Omega)}$ by ensuring conditions like
\begin{equation*}
 \V^+(\grid_\omega) \subset \V(\grid_*), 
\end{equation*}
where $\omega$ is a $\grid$-mesh subdomain. Inspecting the construction of $\V^+_\grid$, we realize that such conditions can be achieved if $ \max \{ m_1, m_2 \} \leq n-1 $ and the cut-off is implemented with hat functions of a virtual refinement of $\grid$.

\begin{lemma}[cut-off by refined hat functions]
\label{L:cut-off-with-hat}    
Let $\grid_+$ be the minimal bisection refinement of $\grid$ such that the relative interior of each element $T \in \grid$ and each face $F \in \faces$ of the original mesh $\grid$ contains at least one vertex from $\grid_+$. Then there exist hat functions  $\phi_T$, $T \in \grid$, and $\phi_F$, $F \in \faces$, in $\mathbb{S}^{1,0}(\grid_+)$ satisfying Assumption~\ref{A:abstract-cut-off} if $\max\{m_1,m_2\} \leq n-1$.
\end{lemma}

\begin{proof}
The details of the proof depend on bisection and we therefore restrict to the case $d=2$; for $d > 2$, the following reference situation used to define the hat functions is replaced by several ones with ``tagged'' reference simplices. Let $\wh{T}=T_2$ be the reference element in $\R^2$ with the standard enumeration of its vertices $\wh{z}_0 = 0$, $\wh{z}_1 = e_1$,  and $\wh{z}_2 = e_2$. Furthermore, let $\wh{\grid}_+$ be the mesh obtained by applying 5 bisections so that vertices in the interiors of $\wh{T}$ and of its faces are generated. Denote by $\wh{\phi}_{\wh{T}}$, $\wh{\phi}_{F'}$, $F'\subset\wh{T}$, the four hat functions in $\mathbb{S}^{1,0}(\wh{\grid}_+)$ associated with these generated vertices.  Given an arbitrary element $T \in \grid$, let $H_T$ denote the bi-affine map $T \to \wh{T}$ preserving the numbering of the vertices for bisection and define the pull-backs
\begin{equation*}
 \phi_T
 \definedas
 H_T^* \big( \wh{\phi}_{\wh{T}} \big),
\qquad
 \phi_F{}|_{T}
 \definedas
 H_T^* \big( \wh{\phi}_{H_T(F)} \big),
\quad
 F \subset T,
\end{equation*}
and extend by 0 off $T$ or $\omega_F$. As the extension operators $E_F$ preserve the polynomial degree, see Lemma~\ref{L:extending-from-faces} (extending from faces),  $\max\{ m_1, m_2 \} \leq n-1$, and
\begin{equation*}
\big\{ H_T^{-1}(\wh{T}_+) \mid \wh{T}_+ \in \wh{\grid}_+ \big\}
 =
 \big\{ T_+ \in \grid_+ \mid T_+ \subset T \big\},
\end{equation*}
the hat functions $\phi_T$, $T \in \grid$, and $\phi_F$, $F \in \faces$, then satisfy Assumption \ref{A:abstract-cut-off} with $G_T = H_T^{-1}$ and $\mathbb{S}^+ = \mathbb{S}^{n,0}(\wh{\grid}_+)$.
\end{proof}

% lower bounds by interior vertices

\begin{definition}[interior vertex property]\label{D:interior-node}
\index{Definitions!Interior vertex property}
A mesh element $T \in \grid$ satisfies the \emph{interior vertex property}
 with respect to $\grid_* \geq \grid$ whenever each interior face $F \subset \partial T \setminus \partial\Omega$ of $T$ and each element in $\wt{\omega}_T$ (defined in \eqref{tildeomegaT}) have in their relative interiors at least one vertex from $\grid_*$. 

A set $\marked \subset \grid$ satisfies the interior vertex property with respect to a refinement
$\grid_*\ge\grid$ if each element $T\in\marked$ satisfies the interior vertex property. 
\end{definition}

The interior vertex property is valid upon enforcing a fixed number $b$ of bisections ($b=3,6$ for $d=2,3$).
An immediate consequence is the following lower bound for corrections.

\begin{theorem}[lower bound for corrections]
\label{T:lubd-corr}
Suppose $\bA$ is piecewise constant over $\grid$ and $c=0$, define $P_\grid$ with the help of the cut-off functions in Lemma~\ref{L:cut-off-with-hat}, and denote by $\mathcal{M}$ the subset of elements in $\grid$ satisfying the interior vertex property with respect to $\grid_*$. Then
\begin{equation*}
 \sum_{T \in \mathcal{M}} \est_\grid(u_\grid, f, T)^2
 \leq
 C \| \nabla(u_\grid - u_{\grid_*}) \|_{L^2(\Omega)}^2
 +
 \sum_{T \in \mathcal{M}} \| f - P_\grid f \|_{H^{-1}(\omega_T)}^2,
\end{equation*}
where $C$ depends on $d$, the shape regularity coefficient $\sigma$, the coefficients $(\bA,c)$, and on the polynomial degree $n$.
\end{theorem}

\begin{proof}
\fbox{\scriptsize 1} We first show a local bound with the PDE indicator $\eta_\grid(u_\grid,T)$. In view of $\bA \in {\mathbb S}^{0,-1}_\grid$ and $c=0$, we choose $m_1=n-1$ and $m_2=n-2$ as the degrees for the discrete functionals. We can thus apply Lemma~\ref{L:cut-off-with-hat} and construct $P_\grid$ with the refined hat functions. Let $T \in \mathcal{M}$ and so, using the interior vertex property and the notation associated with $\wt{\omega}_T$ in \eqref{tildeomegaT}, we deduce $\V^+(\wt{\grid}_T) \subset \V(\grid_*,T) \definedas \V(\grid_*) \cap H^1_0(\wt{\omega}_T)$, where $\wt \grid_T:= \{ T \in \grid \ | \ T \subset \wt \omega_T\}$. Combining this with inequality \eqref{Pstability;lower-bd} and Definition~\ref{D:Pgrid} (projection onto discrete functionals), we conclude
\begin{equation*}
 \eta_\grid(u_\grid,T)
 \Cleq
 \| P_\grid R_\grid \|_{\V^+(\wt{\grid}_T)^*}
 =
 \| R_\grid \|_{\V^+(\wt{\grid}_T)^*}
 \leq
 \| R_\grid \|_{\V(\grid_*,T)^*}.
\end{equation*}

\smallskip\fbox{\scriptsize 2} To collect the local bounds of the first step, we first show that, for any $\ell \in \V(\grid_*)^*$,
\begin{equation*}
 \sum_{T\in\grid} \| \ell \|_{_{\V(\grid_*,T)^*}}^2
 \leq
 (d+2)  \| \ell \|_{\V(\grid_*)^*}^2.
\end{equation*}
To this end, we just repeat the proof of Lemma~\ref{L:loc-of-H^-1-norm} (localization of $H^{-1}$-norm), replacing the spaces $H^1_0(\omega_i)$ and $H^1_0(\Omega)$, respectively, with $\V(\grid_*,T)$ and $\V(\grid_*)$. Hence, squaring and summing the bound of the first step as well as using Lemma~\ref{L:err-residual} (error and residual) yield
\begin{equation}
\label{marked-PDE-ind-bdd-by-corr}
 \sum_{T \in \mathcal{M}} \eta_\grid(u_\grid,T)^2
 \Cleq
 \sum_{T \in \mathcal{M}}  \| R_\grid \|_{\V(\grid_*,T)^*}^2
 \Cleq
 \| \nabla (u_{\grid_*} - u_\grid) \|_{L^2(\Omega)}^2.
\end{equation}

\smallskip\fbox{\scriptsize 3} We finally prove the claimed bound by simply inserting \eqref{marked-PDE-ind-bdd-by-corr}:
\begin{align*}
 \sum_{T \in \mathcal{M}} \est_\grid(u_\grid,f,T)^2
 &=
 \sum_{T \in \mathcal{M}} \big( \eta_\grid(u_\grid,T)^2 + \osc_\grid(f,T)^2 \big)
\\
 &\leq
 C  \| \nabla (u_{\grid_*} - u_\grid) \|_{L^2(\Omega)}^2
 +
 \sum_{T \in \mathcal{M}} \osc_\grid(f,T)^2
\end{align*}
and the proof is finished.
\end{proof}

\begin{remark}[oscillation and correction]
\label{R:osc-corr}
In general, fixing first the finer mesh $\grid_*$, it is impossible to bound oscillation indicators $\osc_\grid(f,T)$ by some suitable correction. Indeed, these indicators can contain contributions to $f$ and so to $R_\grid$ of ``arbitrarily high frequency'', while the correction can control only contributions of the residual $R_\grid$ with frequencies representable over $\grid_*$; cf.\ \eqref{residual-norm-of-corr}.
\end{remark}

Monotonicity properties of the error $|u-u_\grid|_{H^1_0(\Omega)}= \| \nabla (u-u_\grid) \|_{L^2(\Omega)}$ and the PDE
error estimator $\eta_\grid(u_\grid) = \eta_\grid(u_\grid,f)$ with respect to $\grid$ would be useful but
fail to hold. To investigate this issue, we consider two admissible meshes $\grid,\grid_*\in\grids$,
the latter being a refinement of the former $\grid_* \ge \grid$, and a third admissible mesh
$\wh{\grid} \le \grid$. We further assume that data $\data=(\bA,c,f)$ is discrete over $\wh{\grid}$
in the sense that $\data\in\D_{\wh{\grid}}$ where
\begin{equation*} %\label{E:space-discrete-data}
\D_{\wh\grid} := \big[ {\mathbb S}^{n-1,-1}_{\wh{\grid}} \big]^{d\times d}
\times {\mathbb S}^{n-1,-1}_{\wh{\grid}} \times \F_{\wh{\grid}}
\end{equation*}
and $\data$ does not change in the transition from $\grid$ to $\grid_*$ irrespective of
the degree of local refinement; in particular $f=P_{\wh\mesh} f \in \F_{\wh\mesh}$. We will later denote discrete data as
$\wh{\data}=(\wh{\bA},\wh{c},\wh{f})$ to distinguish it from exact data $\data$, and to study their
discrepancy, but we prefer to keep the simple notation $\data=\wh{\data}$ now because there is no 
reason for confusion. In particular, this implies that the bilinear form in \eqref{bilinear-coercive}
and forcing function
are the same for both Galerkin solutions $u_\grid\in\V_\grid$ and $u_{\grid_*}\in\V_{\grid_*}$, whence
the energy errors are monotone according to \eqref{E:quasi-monotonicity-coercive}
\[
\enorm{u-u_{\grid_*}} \le \enorm{u-u_{\grid}},
\]
but not $|u-u_\grid|_{H^1_0(\Omega)}$.
Moreover, $\eta_\grid(u_\grid)$ is not monotone because the discrete solution
$u_\grid\in\V_\grid$ changes with the mesh. It is thus useful to quantify
the behavior of $\eta_\grid(u_\grid)$ in terms of $\grid$ and $u_\grid$
following \cite{CaKrNoSi:08}; see also \cite{MoSiVe:08}. We do this next.

The first lemma exploits the structure of the PDE residual 
estimator, namely the presence of a positive power of the local meshsize,
and expresses the {\it reduction} of $\eta_{\grid_*}(v,f)$ relative to $\eta_{\grid}(v,f)$ for
fixed functions $v\in\V_\grid$ and $f\in\F_\mesh$. This quantitative property is instrumental in studying
convergence of {\AFEM}s for coercive problems in Section \ref{S:conv-rates-coercive}
as well as discontinuous Galerkin methods in Section \ref{S:dg} and
inf-sup stable problems in Section \ref{S:conv-rates-infsup}.

\begin{lemma}[reduction property of the estimator]\label{L:strict-est}
If the elements of $\marked\subset\grid$ are bisected at least 
$b\ge1$ times to refine $\grid$ into $\grid_*$, and $\lambda=1-2^{-b/d}$, then
\begin{equation}\label{E:strict-est}
\rhn{\eta_{\gridk[*]} (v,f,\gridk[*])^2 \le \eta_{\grid} (v,f,\grid)^2
-\lambda \, \eta_{\grid} (v,f,\marked)^2
\qquad\forall v\in\V_\grid, f\in\F_\grid.}
\end{equation}
\end{lemma}
\begin{proof}
Given $T\in\grid$, we rewrite \eqref{mod-res-est} as follows
\[
\eta_\grid(v,T)^2 = h_T j_\grid(v,T)^2 + h_T^2 r_\grid(v,T)^2
\]
\rhn{with $\eta_\grid(v,T) = \eta_\grid(v,f,T)$ and}
\begin{align*}
j_\grid(v,T)^2 &= \rhn{j_\grid(v,f,T)^2 =} \sum_{\substack{F \in \faces\\ F\subset \partial T}} \| \jump{\bA\nabla v}\cdot \vec{n}_F - P_\grid f\|_{L^2(F)}^2,
\\
r_\grid(v,T) &= \rhn{r_\grid(v,f,T) =} \| P_\grid f + \div(\bA\nabla v) -cv \|_{L^2(T)},
\end{align*}
where $f = P_\grid f = P_{\grid_*}f \in \F_\grid$ does not change from $\grid$ to $\grid_*$.
We readily have
\[
\sum_{T_* \in \grid_*, \ T_* \subseteq T}\eta_{\grid_*}(v,T_*)^2 \le \eta_\grid(v,T)^2.
\]
because $h_{T_*} \le h_T$ for all $T_*\subset T$ and $T_*\in\grid_*$. If, in addition,
$T$ is bisected at least $b$ times, then any such $T_*$ satisfies $h_{T_*} \le 2^{-\frac{b}{d}} h_T$,
whence
\[
\sum_{T_* \in \grid_*, \ T_* \subseteq T}\eta_{\grid_*}(v,T_*)^2  \le 2^{-\frac{b}{d}} \eta_\grid(v,T)^2.
\]
Therefore, adding over $T\in\grid$ we obtain
\[
\eta_{\grid_*}(v)^2 = \sum_{T\in\grid}\sum_{T_* \in \grid_*, \ T_* \subseteq T} \eta_{\grid_*}(v,T_*)^2
\le 2^{-\frac{b}{d}} \sum_{T\in\marked} \eta_{\grid}(v,T)^2 + \sum_{T\in\grid\backslash\marked}
\eta_{\grid}(v,T)^2
\]
which implies the assertion \eqref{E:strict-est}.
\end{proof}

The next result complements Lemma \ref{L:strict-est} in that it
expresses the Lipschitz continuity of \rhn{$\eta_\grid(v,f)$ with respect to the
argument $v\in\V_\grid$ for fixed $\grid$ and $f\in\F_\grid$.} \looseness=-1

\begin{lemma}[Lipschitz property of the estimator]\label{L:lipschitz}
\rhn{Let $\grid$ and $f\in\F_\grid$ be fixed.}
\index{Constants!$\Clip$: estimator Lipschitz property constant}
There exists a constant $\Clip$ proportional to $\|A\|_{L^\infty(\Omega)}+\|c\|_{L^\infty(\Omega)}$
such that 
  \begin{equation}\label{E:lipschitz-estimator}
    \rhn{|\eta_\grid(v,f)-\eta_\grid(w,f)|} \le \Clip |v-w|_{H^1_0(\Omega)}
    \quad\forall v,w\in\V_\grid.
  \end{equation}
\end{lemma}
\begin{proof}
Since \rhn{$\eta_\grid(v) = \eta_\grid(v,f)$} is the $\ell^2$-norm of the vector $(\eta_\grid(v,T))_{T\in\grid}\in\R^{\#\grid}$,
applying the triangle inequality gives
\begin{align*}
\big| \eta_\grid(v) - \eta_\grid(w) \big|^2
&\le \sum_{T\in\grid}\big| \eta_\grid(v,T) - \eta_\grid(w,T) \big|^2
\\
& \le \sum_{T\in\grid} h_T \big|j_\grid(v,T) - j_\grid(w,T) \big|^2
+ h_T^2 \big|r_\grid(v,T) - r_\grid(w,T) \big|^2.
\end{align*}
We first consider the jump terms and apply again the triangle inequality
followed by an inverse estimate to find that
\begin{equation*}
    \begin{split}
\big|j_\grid(v,T) - j_\grid(w,T) \big|^2 &\le \sum_{\substack{F \in \faces \\ F \subset \partial T}} \| \jump{\bA \nabla(v-w) } \cdot \vec{n}_F \|_{L^2(F)}^2\\
&\lesssim h_T^{-1}\|\bA\|_{L^\infty(\Omega)}^2 \| \nabla(v-w) \|_{L^2(\omega_T)}^2,        
    \end{split}
\end{equation*}
where $\omega_T$ is the patch of $T$.
A similar reasoning for the element residuals yields
\begin{align*}
\big|r_\grid(v,T) - r_\grid(w,T) \big|^2 & \lesssim \|c(v-w)\|_{L^2(T)}^2 + \|\div(\bA \nabla(v-w))\|_{L^2(T)}^2
\\ & \lesssim \|c\|_{L^\infty(\Omega)}^2 \|(v-w)\|_{L^2(T)}^2 + h_T^{-2} \|\bA\|_{L^\infty(\Omega)}^2 \|\nabla(v-w)\|_{L^2(T)}^2,
\end{align*}
because $A$ is piecewise polynomial.
Finally, adding over $T\in\grid$ and applying Lemma \ref{L:Poincare} (first Poincar\'e inequality)
concludes the proof.
\end{proof}

\rhn{
Since the estimator $\eta_\grid(v,f)$ depends explicity on $P_\grid f$, and $P_\grid f$ may change with $\grid$, it is crucial to account for the variations of $\eta_\grid(v,f)$ while keeping $\grid$ and $v\in\V_\grid$ fixed. This is the purpose of our next result.

\begin{lemma}[estimator dependence on discrete forcing]\label{L:est-forcing}
  Let $\grid$ and $v\in\V_\grid$ be fixed. Then there exists a constant $\Clip$ 
such that 
  \begin{equation}\label{E:est-forcing}
    |\eta_\grid(v,f)-\eta_\grid(v,g)| \le \Clip \left( \sum_{T\in\grid} \|f-g\|_{H^{-1}(\omega_T)}^2 \right)^{1/2}
    \quad\forall f,g\in\F_\grid.
  \end{equation}
\end{lemma}
\begin{proof}
We proceed elementwise, as in Lemma~\ref{L:lipschitz}, except that after applying the triangle inequality we end up with the weighted $L^2$-norms
$$
h_T^2 \| f - g \|_{L^2(T)}^2 + h_T \| f - g\|_{L^2(\partial T)}^2, \qquad \forall T\in \grid. 
$$
Extending these norms to patches $\omega_T$ and its interior faces $\sigma_T$, and appealing to Corollary~\ref{C:quantify-discrete-dual-norms} (quantifying $H^{-1}$-norms of discrete functionals), we deduce
$$
h_T^2 \| f - g \|_{L^2(\omega_T)}^2 + h_T \| f - g\|_{L^2(\sigma_T)}^2 \approx \| f - g \|_{H^{-1}(\omega_T)}^2,
$$
Adding over $T\in\grid$ finishes the proof. 
\end{proof}

In the subsequent applications of Lemma \ref{L:lipschitz} the discrete coefficients $(\bA,c)$
may change with
the change of the supporting mesh $\wh{\grid}$ but they will always be uniformly bounded in
$L^\infty(\Omega)$; hence the constant $\Clip$ is uniformly bounded as well.
Upon combining Lemmas \ref{L:strict-est}, \ref{L:lipschitz} and \ref{L:est-forcing}
we obtain the following crucial property.

\begin{proposition}[estimator reduction]\label{P:est-reduction}
Given $\grid\in\grids$ and a subset $\marked\subset\grid$
of elements marked for refinement, let {\rm \REFINE} be the procedure discussed in
Section \ref{S:bisection} that bisects the elements of $\marked$ at least $b$ times
and $\grid_*=\REFINE\big(\grid,\marked\big)$ be the resulting conforming mesh. Let the
coefficients $(\bA,c)$ be discrete and fixed. Then for $\lambda=1-2^{-b/d}$, for all $v\in\V_\grid$, $v_*\in\V_{\grid_*}$, $f\in\F_\grid, f_*\in\F_{\grid_*},$ and any $\delta>0$
\index{Constants!$\Clip$: estimator Lipschitz property constant}
\begin{align*}
      \eta_{\grid_*}(v_*,f_*,\grid_*)^2 &\le
      (1+\delta)\big(\eta_\grid(v,f,\grid)^2 - \lambda\,\eta_\grid(v,f,\marked)^2\big)
      \\ &+
      2(1+\delta^{-1})\,\Clip^2 \, \left(|v_*-v|_{H^1_0(\Omega)}^2 +
      \sum_{T_*\in\grid_*} \|f_*-f\|_{H^{-1}(\omega_{T_*})}^2 \right),
\end{align*}
where $\Clip$ is the constant in Lemmas~\ref{L:lipschitz} and \ref{L:est-forcing}.
\end{proposition}
\begin{proof}
For any $\delta>0$, write
\[
\eta_{\grid_*}(v_*,f_*,\grid_*)^2 \leq (1+\delta) \eta_{\grid_*} (v,f,\grid_*)^2 +
\big(1 + \delta^{-1} \big)\big(  \eta_{\grid_*} (v_*,f_*,\grid_*) -  \eta_{\grid_*} (v,f,\grid_*) \big)^2 
\]
and apply Lemma \ref{L:strict-est} to the first term and Lemmas \ref{L:lipschitz} and \ref{L:est-forcing} to the second one combined with a triangle inequality.
\end{proof}
}

We finish this section by investigating the behavior of the global oscillation under refinement.
\begin{lemma}[quasi-monotonicity of oscillation]
\label{L:quasi-mono-osc}
If $f \in H^{-1}(\Omega)$ and  $\grid,\grid_* \in \grids$ with $\grid_* \geq \grid$, then
\begin{equation*}
\index{Constants!$C_{\mathrm{osc}}$: oscillation quasi-monotonicity constant}
%\label{e:quasi-mono-osc}
\osc_{\grid_*}(f) \leq C_{\mathrm{osc}} \osc_{\grid}(f),
\end{equation*}
where $C_\mathrm{osc}$ depends only on the shape regularity coefficient $\sigma$ and $d$.
\end{lemma}

\begin{proof}
Given $T\in \grid$, let $T_* \in \grid_*$ such that $T_* \subset T$. Since $\grid_*$ is a refinement of $\grid$, this implies that the patch $\omega_{\grid_*}(T_*)$ in $\grid_*$ around $T_*$ is contained in the patch $\omega_\grid(T)$ in $\grid$ around $T$. Thanks to Lemma~\ref{C:local-near-best-approx-of-P} (local near-best approximation), we derive 
\begin{equation*}
 \osc_{\grid_*}(f,T_*)^2
 =
 \| (I-P_{\grid_*}) f \|_{H^{-1}(\omega_{\grid_*}(T_*))}^2
 \leq
 C_\mathrm{lStb}^2  \| (I-P_{\grid}) f \|_{H^{-1}(\omega_{\grid_*}(T_*))}^2
\end{equation*}
and therefore, with the  help of (ii) of Lemma~\ref{L:loc-of-H^-1-norm} (localization of $H^{-1}$-norm),
\begin{equation*}
 \sum_{T_* \subset T} \osc_{\grid_*}(f,T_*)^2
 \leq
 C_\mathrm{lStb}^2 C_\mathrm{ovrl} \| (I-P_{\grid}) f \|_{H^{-1}(\omega_{\grid}(T))}^2
 =
 C_\mathrm{lStb}^2 C_\mathrm{ovrl} \osc_\grid(f,T)^2,
\end{equation*}
where $C_\mathrm{ovrl}$ is bounded in terms of the shape regularity coefficient $\sigma$ and $d$. Hence, summing over $T \in \grid$ yields
\begin{equation*}
 \osc_{\grid_*}(f)^2
 =
 \sum_{T \in \grid} \sum_{T_* \subset T} \osc_{\grid_*}(f,T_*)^2
 \leq
 C_\mathrm{lStb}^2 C_\mathrm{ovrl} \osc_\grid(f)^2
\end{equation*}
and the proof is finished.
\end{proof}
%-------------------------------------------------------------------------------------
\subsection{Alternative estimators}
\label{S:other_estimators}
%-------------------------------------------------------------------------------------

% intro
In Sect.\ \ref{S:mod-res-est} we used the $H^{-1}$-projection $P_\grid$ to derive a posteriori bounds for the error in the spirit of the standard residual estimator. The goal of this section is to illustrate that the approach with $P_\grid$ can be combined also with other techniques of a~posteriori error estimation, generalizing and expanding the discussion in \cite[Sect. 4]{KreuzerVeeser:2021} with the $H^{-1}$-projection $P_\grid$. 

Alternative techniques have been developed with the desire to reduce or even circumvent that constants spoil the relationship between error and estimator. In the framework of the aforementioned approach, we shall see that the various techniques based upon
\begin{itemize}
\item \emph{local (discrete) problems},
\item \emph{hierarchy},
\item \emph{flux equilibration}
\end{itemize}
amount to different ways of quantifying a local norm of the discretized residual $P_\grid R_\grid$. This observation is useful for comparing the techniques  and for a common treatment in the following sections about adaptive algorithms. 

% general setting
\smallskip As in Sect.\ \ref{S:mod-res-est} on modified residual estimation, we shall consider only the case given by Assumption~\ref{A:discrete-coefficients} (discrete coefficients and discrete functionals). For the hidden constants in the results of this section, it is useful to keep in mind Remark~\ref{R:constants-for-error-res-rel} (constants in error-residual relationship).

% vertex-indexed variant of modified residual estimator for comparison
Theorem~\ref{T:modified-estimator} (modified residual estimator) analyzed an element-indexed version of the residual estimator. For the sake of simplicity, we shall refrain here from such an element-indexed setting and remain in the vertex-indexed setting of the abstract analysis of Sect.\ \ref{S:dis-osc-res}. In order to facilitate the comparison with the other estimators below, we offer the following vertex-indexed variant of Theorem~\ref{T:modified-estimator}.
For $z \in \vertices$, we set $\faces_z:= \{ F \in \faces \ | \ z \in F \}$ and
$\grid_z:= \{ T \in \grid \ | \ z \in T \}$. 
Given the Galerkin approximation $u_\grid$ from \eqref{setting-apost:Galerkin}, define the PDE indicator by
\begin{subequations}
\label{vertex-mod-res-est}
\begin{equation}
\label{vertex-mod-rest-est;pde-ind}
\begin{aligned}
 \eta^\mathrm{res}_\grid(u_\grid)
 &\definedas
 \sum_{z \in \vertices}  \eta^\mathrm{res}_\grid(u_\grid,z)^2
\quad\text{with}
\\
 \eta^\mathrm{res}_\grid &(u_\grid,z)^2
 \definedas
 \sum_{F \in \faces_z} h_F \| P_F R_\grid \|_{L^2(F)}^2
 +
 \sum_{T \in \grid_z} h_T^2 \| P_T R_\grid \|_{L^2(T)}^2,
\end{aligned}
\end{equation}
where $R_\grid = f + \operatorname{div}(\bA \nabla u_\grid) - c u_\grid \in H^{-1}(\Omega)$ is the residual and $P_F$, $F\in\faces$ and $P_T$, $T \in \grid$ yield the polynomial densities of $P_\grid$; see Definition \ref{D:Pgrid}. The vertex-indexed modified residual estimator is then
\begin{equation}
\label{vertex-mod-rest-est;total}
 \est^\mathrm{res}_\grid
 \definedas
 \est^\mathrm{res}_\grid(u_\grid,f)^2
 \definedas
 \eta^\mathrm{res}_\grid(u_\grid)^2
 +
 \osc_\grid(f)^2.
\end{equation}
\end{subequations}
\begin{theorem}[vertex-indexed modified residual estimator]
\label{T:vertex-mod-res-est}
Suppose Assumption~\ref{A:discrete-coefficients}. The modified residual estimator \eqref{vertex-mod-rest-est;total} is equivalent to the error: 
\begin{equation*}
 \frac{ \min\{ 1, C_\mathrm{L,res} \} }{ C_\mathrm{lStb}^\star C_d } \est^\mathrm{res}_\grid
 \Cleq
 \| \nabla(u-u_\grid) \|_{L^2(\Omega)}
 \Cleq
 \max\{ 1, C_\mathrm{U,res} \} C_d C_\mathrm{loc} \est^\mathrm{res}_\grid, 
\end{equation*}
while its PDE indicator \eqref{vertex-mod-rest-est;pde-ind} is locally equivalent to the discretized residual: for all vertices $z \in \vertices$,
\begin{equation*}
 C_\mathrm{L,res} \eta^\mathrm{res}_\grid(u_\grid,z)
 \leq
 \| P_\grid R_\grid \|_{H^{-1}(\omega_z)}
 \leq
 C_\mathrm{U,res} \eta^\mathrm{res}_\grid(u_\grid,z).
\end{equation*}
Here, $C_\mathrm{L,res}$ and $ C_\mathrm{U,res}$ are the hidden constants of Corollary \ref{C:quantify-discrete-dual-norms} on stars, $C_\mathrm{lStb}^\star$ is the stability constant of $P_\grid$ on stars from Lemma~\ref{L:local-stability-of-P}, $C_\mathrm{loc}$ is the constant from Corollary~\ref{C:residual-loc}, $C_d = \sqrt{2(d+1)}$, and the hidden constants depend only on the error-residual relationship in Lemma~\ref{L:err-residual}.
\end{theorem}

\begin{proof}
Local equivalence is a reformulation of Corollary \ref{C:quantify-discrete-dual-norms} (quantifying $H^{-1}$-norms of discrete functionals) on stars. The global bounds follow by combining local equivalence with Lemma~\ref{L:err-residual} (error and residual), Corollary~\ref{C:residual-loc} (star localization of residual norm), and Lemma~\ref{L:abstract_bds} (splitting of local residual norms).
\looseness=-1
\end{proof}

%
%-----------------------------------------------------------------------------------
\subsubsection*{Adjoint projection}
%-----------------------------------------------------------------------------------
%
The projection $P_\grid$ relates residual $R_\grid$ and discretized residual $P_\grid R_\grid$. In order to exploit this relationship on the test space $H^1_0(\Omega)$, we shall need the adjoint $P_\grid^*$ to the projection $P_\grid$. Curiously, operators employed in this vein appeared first; see, e.g., \cite{Morin.Nochetto.Siebert:2003-Local-probs} and \cite{Veeser:2002-nonlinear-Laplacian}.

\smallskip Given $w \in H^1_0(\Omega)$, the function $P_\grid^* w$ can be directly defined by requiring 
\begin{equation}
\label{adjoint-of_Pgrid}
P_\grid^* w \in \V_\grid^+: \quad
 \langle \ell, P_\grid^* w \rangle
 =
 \langle \ell, w \rangle \quad  \forall \ell \in \F_\grid
\end{equation}
This definition is well-posed thanks to Lemma~\ref{L:inf-sup_N} (discrete inf-sup condition) and Lemma~\ref{L:algebra-of-P} (algebraic properties), especially \eqref{dimV+mesh=dimFmesh} and \eqref{pairing-nondeg-on-Fmesh-V+mesh;2}. Clearly, $P_\grid^*$ is a linear projection onto the finite-dimensional subspace $\V_\grid^+ \subset H^1_0(\Omega)$. A representation as interpolation operator will be derived in Corollary~\ref{C:proj-as-interpolators} below. Using both definitions of $P_\grid$ and $P_\grid^*$, we see that they are actually adjoint:
\begin{equation}
\label{Pgrid*-is-adjoint}
 \langle P_\grid \ell, w \rangle
 =
 \langle P_\grid \ell, P_\grid^* w \rangle
 =
 \langle \ell, P_\grid^* w \rangle \quad
 \forall \ell \in H^{-1}(\Omega),\,  w \in H^1_0(\Omega).
\end{equation}
Consequently, Lemmas~\ref{L:algebra-of-P} and \ref{L:local-stability-of-P} (local $H^{-1}$ stability) show that $P_\grid^*$ is a local operator with
\begin{equation}
\label{local-stability-of-Pgrid*}
 \| P_\grid^* \|_{\mathcal{L}(H^1_0(\omega))}
 =
 \| P_\grid \|_{\mathcal{L}(H^{-1}(\omega))}
 \leq
 C_\mathrm{lStb}.
\end{equation}
The choice $\ell = R_\grid$ in \eqref{Pgrid*-is-adjoint} leads to
\begin{equation*}
 \langle P_\grid R_\grid, w \rangle
 =
 \langle P_\grid R_\grid, P_\grid^* w \rangle
 =
 \langle R_\grid, P_\grid^* w \rangle \quad
 \forall w \in H^1_0(\Omega),
\end{equation*}
where the two identities show that the discretized residual $P_\grid R_\grid$ can be analyzed with discrete test functions in $\V_\grid^+$ only, cf.\ the norm equivalence in Lemma \ref{L:local-stability-of-P}. Restricting to discrete test functions in $\V_\grid^+ = \operatorname{Im} P_\grid^*$, we obtain:
\begin{equation}
\label{testing-with-Vgrid+}
 \langle P_\grid R_\grid, w \rangle
 =
 \langle R_\grid, w \rangle
 \quad \forall w \in \V_\grid^+.
\end{equation}

\subsubsection*{An estimator based upon local problems}
Local dual norms can be quantified by solving local problems. Requiring computability of these solutions leads to \emph{finite-dimensional} or \emph{discrete} local problems. In other words, we lift the residual to local and finite-dimensional extensions of the finite element space. Starting with \cite{BabuskaRheinboldt:78}, this idea was used to soften the impact of constants in the relationship between error and estimator; cf.\ \cite[Remark~1.22]{Verfuerth:13} for more references.

\smallskip Within the approach of Sections \ref{S:err-res} and \ref{S:dis-osc-res}, we can use local discrete problems to quantify the local $H^{-1}$-norms of the discretized residual $P_\grid R_\grid$. In this manner, constants arise only due to the localization of the residual norm and to the splitting into discretized and oscillatory residual by the $H^{-1}$-projection $P_\grid$.

% definition of indicator
We start by introducing the vertex-oriented PDE indicator. Given the Galerkin approximation $u_\grid$ from \eqref{setting-apost:Galerkin}, set
\begin{subequations}
\label{lpb-est}
\begin{equation}
\label{lpb-est;pde-ind}
 \eta^\mathrm{lpb}_\grid(u_\grid)
 \definedas
 \sum_{z \in \vertices}  \eta^\mathrm{lpb}_\grid(u_\grid,z)^2
\quad\text{with}\quad
 \eta^\mathrm{lpb}_\grid(u_\grid,z)
 \definedas
 \| \nabla v_z \|_{L^2(\omega_z)},
\end{equation}
where $v_z \in \V^+(\grid_z)$ is the solution of the local problem
\begin{equation}
\label{lpb-est;lpb}
 \int_{\omega_z} \nabla v_z \cdot \nabla w 
 =
 \langle R_\grid, w \rangle \quad \forall w \in \V^+(\grid_z).
\end{equation}
Note that this problem is discrete for $\dim \V^+(\grid_z) < \infty$ and therefore can be solved up to machine precision. The resulting estimator is then
\begin{equation}
 \est^\mathrm{lpb}_\grid
 \definedas
 \est^\mathrm{lpb}_\grid(u_\grid,f)^2
 \definedas
 \eta^\mathrm{lpb}_\grid(u_\grid)^2
 +
 \osc_\grid(f)^2.
\end{equation}
\end{subequations}

\begin{theorem}[estimator based on local problems]
\label{T:lpb-est}
Under Assumption~\ref{A:discrete-coefficients} the estimator \eqref{lpb-est} based on local problems is equivalent to the error, while its PDE indicator is locally equivalent to the discretized residual with constant $1$ in the lower bound so that
\begin{equation*}
 \frac{1}{C_\mathrm{lStb}^\star C_d} \est^\mathrm{lpb}_\grid
 \Cleq
 \| \nabla(u-u_\grid) \|_{L^2(\Omega)}
 \Cleq
 C_\mathrm{lStb}^\star C_d C_\mathrm{loc} \est^\mathrm{lpb}_\grid
\end{equation*}
and, for all vertices $z \in \vertices$,
\begin{equation*}
 \eta^\mathrm{lpb}_\grid(u_\grid,z)
 \leq
 \| P_\grid R_\grid \|_{H^{-1}(\omega_z)}
 \leq
 C_\mathrm{lStb}^\star \eta^\mathrm{lpb}_\grid(u_\grid,z)
\end{equation*}
where $C_\mathrm{lStb}^\star$ is the stability constant of $P_\grid$ on stars from Lemma~\ref{L:local-stability-of-P}, $C_\mathrm{loc}$ from Corollary~\ref{C:residual-loc}, $C_d = \sqrt{2(d+1)}$ and the hidden constants depend only on the error-residual relationship in Lemma~\ref{L:err-residual}.
\end{theorem}

\begin{proof}
It suffices to show the local equivalence for the PDE indicator; %$\eta^\textrm{lpb}_\grid(u_\grid)$;
cf.\ Theorem~\ref{T:vertex-mod-res-est} (vertex-indexed modified residual estimator) and note $C_\mathrm{lStb}^\star \geq 1$.  Let $z \in \vertices$ be any vertex. In view of \eqref{testing-with-Vgrid+}, the definition of $v_z \in \V^+(\grid)$ readily implies
\begin{equation*}
 \eta^\mathrm{lpb}_\grid(u_\grid,z)
 =
 \| \nabla v_z \|_{L^2(\omega_z)}
 =
 \| P_\grid R_\grid \|_{\V^+(\grid)^*}\,.
\end{equation*}
Hence Lemma~\ref{L:local-stability-of-P} on the local $H^{-1}$-stability of $P_\grid$ yields the asserted local equivalence of PDE indicator and discretized residual $P_\grid R_\grid$.
\end{proof}

\subsubsection*{A stable biorthogonal system for $\F_\grid \times \V^+_\grid$}
%
% intro
Stable biorthogonal systems induce linear bounded projections, which enjoy near-best approximation thanks to the Lebesgue lemma. Supposing Assumption~\ref{A:abstract-cut-off} (abstract cut-off), we now outline the construction of such a system for the finite-dimensional product $\F_\grid \times \V^+_\grid$. The constructed system will induce both projections $P_\grid$ and its adjoint $P_\grid^*$. This generalizes the bi-orthogonal system in \cite[Sect. 3.4]{KreuzerVeeser:2021} to arbitrary degrees of the discrete functionals and provides an alternative approach to $P_\grid$, its local stability as well as its computation. Furthermore, we use it for devising an hierarchical estimator.

% construction of reference biorthogonal system 
\smallskip The construction is implemented in an affine equivalent manner and our first step consists in setting-up a suitable \emph{reference biorthogonal system}. Let $\wh{T} \definedas T_d$ be the reference element, $\wh{F} \definedas T_{d-1} \times \{0\} \subset \wh{T}$ be the reference face, and denote the polynomial degrees in $\F_\grid$ by $m_1 \in \N_0$ and $m_2 \in \N_0$. Writing
\begin{equation*}
 K_1 \definedas \dim \P_{m_1}(\wh{F}),
\qquad
 K_2 \definedas \dim \P_{m_2}(\wh{T}),    
\end{equation*}
assume that we are given orthonormal bases $q_{(\wh{F},1)},\dots,q_{(\wh{F},K_1)} \in \P_{m_1}(\wh{F})$ and $q_{(\wh{T},1)},\dots,q_{(\wh{T},K_2)} \in \P_{m_2}(\wh{T})$ in the sense that
\begin{equation}
\label{ref-orthonormal-bases}
\begin{aligned}
% \forall k,l = 1,\dots,K_1
%\quad
 \int_{\wh{F}} q_{(\wh{F},k)} q_{(\wh{F},l)} \phi_{\wh{F}}
 =
 \delta_{kl},
\qquad %\\
% \forall k,l = 1,\dots,K_2
%\quad
 \int_{\wh{T}} q_{(\wh{T},k)} q_{(\wh{T},l)} \phi_{\wh{T}}
 =
 \delta_{kl}.
\end{aligned}
\end{equation}
for all admissible $k,l$, i.e.\ $k,l \in \{1,\dots, K_1\} \text{ or } \{1,\dots,K_2\}$ depending on the underlying domain. These bases induce the \emph{reference functionals}
\begin{equation*}
%\label{ref-functionals}
%\begin{aligned}
 \wh{\ell}_{(\wh{F},k)}(w)
 \definedas
 \int_{\wh{F}} q_{(\wh{F},k)} w,
\qquad
% k = 1,\dots,K_1,
%\\
 \wh{\ell}_{(\wh{T},k)}(w)
 \definedas
 \int_{\wh{T}} q_{(\wh{T},k)} w,
%\quad
% k = 1, \dots, K_2
%\end{aligned}
\end{equation*}
on $H^1(\wh{T})$, which in turn span the reference space $\wh{\F}$. In order to define complementing test functions, %associated with the reference face $\wh{F}$
let $\wh{E}$ be the extension operator \eqref{e:E_F} associated with the reference face $\wh{F}$ adapted to the current situation with the only element $\wh{T}$,
%write 
%\begin{equation*}
% \P_{m_2}(\phi_{\wh{T}})
% \definedas
% \{ q \phi_{\wh{T}} \mid q \in \P_{m_2} \}
%\end{equation*}
and, given some $v \in L^2(\wh{T})$, define $\wh{Q}v \in \P_{m_2}$ by
\begin{equation}
\label{ref-Q}
 \int_{\wh{T}} q \big( \wh{Q}v \big) \phi_{\wh{T}}
 =
 \int_{\wh{T}} q v 
 \quad \forall q \in \P_{m_2}.
\end{equation}
We thus define the \emph{reference test functions}
\begin{equation}
\label{ref-test-functions}
\begin{aligned}
 \wh{w}_{(\wh{T},k)}
 &\definedas
 q_{(\wh{T},k)} \phi_{\wh{T}},
\quad
 k = 1,\dots,K_2,
\\
 \wh{w}_{(\wh{F},k)}
 &\definedas
 \widehat{v}_{(\wh{F},k)} - \big( \wh{Q}\widehat{v}_{(\wh{F},k)} \big) \phi_{\wh{T}}, \quad k=1,...,K_1,
\end{aligned}
\end{equation}
with $\widehat{v}_{(\wh{F},k)}
 \definedas
 \big( \wh{E} q_{(\wh{F},k)} \big) \phi_{\wh{F}}$.
Note $\wh{w}_{(\wh{F},k)} \neq 0$. Writing $\wh{I} \definedas \{ (\wh{F},k) \mid k=1,\dots, K_1 \} \cup
 \{ (\wh{T},k) \mid k = 1, \dots, K_2 \}$,
we then have
\begin{equation*}
 \wh{w}_i \in \wh{\V}^+
\definedas
\left\{
 \big( \wh{E}q_1 \big) \phi_{\wh{F}} + q_2 \phi_{\wh{T}}
\mid
 q_1 \in \P_{m_1}, q_2 \in \P_{m_2}
 \right\}
\end{equation*}
for all $i \in \wh{I}$ and the biorthogonality
\begin{equation}
\label{ref-biorthogonality}
 \forall i,j \in \wh{I}
\quad
 \langle \wh{\ell}_i, \wh{w}_j \rangle = \delta_{ij}
\end{equation}
thanks to \eqref{ref-orthonormal-bases} and \eqref{ref-Q}. We thus dispose of a biorthogonal system in the reference product $\wh{\F} \times \wh{\V}^+$.

% construction of global biorthogonal system
Using pull-backs with some minor tweaks, this reference biorthogonal system induces a global biorthogonal one. To this end, we employ bi-affine maps $G_F$, $G_T$, and $G_{(T,F)}$. Here, e.g., given a  pair $(T,F) \in \grid \times \faces$ with $F \subset T$, the map $G_{(T,F)}$ is bi-affine and sends vertices into vertices such that $G_{(T,F)}(\wh{T}) = T$ and $G_{(T,F)}(\wh{F}) = F$. The fact that these maps are only unique up to some renumbering of the vertices is irrelevant as all objects in the reference situation on $(\wh{T},\wh{F})$ are invariant under such renumberings.
We denote the respective inverse maps of $G_F$, $G_T$, and $G_{(T,F)}$ by $H_F$, $H_T$, and $H_{(T,F)}$. Motivated by the transformation rule, we introduce the scaled pull-backs, for $F \in \faces, T \in \grid, k \text{ admissible}$,
\begin{equation}
\label{orthonormal-bases-on-faces-and-elms}
\begin{aligned}
 q_{(F,k)}
 \definedas
 \left( \frac{|\wh{F}|}{|F|} \right)^{\tfrac{1}{2}} H_F^*q_{(\wh{F},k)}, % \in \P_{m_1}(F),
\quad
 q_{(T,k)}
 \definedas
 \left( \frac{|\wh{T}|}{|T|} \right)^{\tfrac{1}{2}} H_T^*q_{(\wh{T},k)}, % \in \P_{m_w}(T)
\end{aligned}
\end{equation}
of the reference orthonormal bases in \eqref{ref-orthonormal-bases}. These lead to the basis
\begin{equation}
\label{basis-for-Fgrid}
%\begin{aligned}
 \ell_{(F,k)}(w)
 \definedas
 \int_{F} q_{(F,k)} w,
\qquad
% k = 1,\dots,K_1,
%\\
 \ell_{(T,k)}(w)
 \definedas
 \int_{T} q_{(T,k)} w,
%\quad
% k = 1, \dots, K_2
%\end{aligned}
\end{equation}
of $\F_\grid$, while the associated test functions are given again through pull-backs:
\begin{equation}
\label{basis-for-Vgrid+}
  w_{(T.k)}
 \definedas
 \left( \frac{|\wh{T}|}{|T|} \right)^{\tfrac{1}{2}} H_T^* \wh{w}_{(\wh{T},k)},
\quad
 w_{(F.k)}{}|_{T}
 \definedas
 \left( \frac{|\wh{F}|}{|F|} \right)^{\tfrac{1}{2}} H_{(T,F)}^* \wh{w}_{(\wh{F},k)}
\end{equation}
for all $T \in \grid$, $F \in \faces$ with $F \subset T$ and all admissible $k$. Note that $w_{(F.k)} \in H^1_0(\Omega)$. Finally, we introduce the index set
\begin{equation*}
 I
 \definedas
 \big( \faces \times \{1, \dots, K_1 \} \big)
 \cup
 \big( \grid \times \{1, \dots, K_2 \} \big).
\end{equation*}
and observe that $w_i \in \V_\grid^+$ for all $i \in I$.
\begin{lemma}[biorthogonal system]
\label{L:biorthogonal-system}
The pairs $(\ell_i,w_i)$, $i \in I$, provide a stable biorthogonal system of the product $\F_\grid \times \V_\grid^+$: indeed, 
%\begin{equation*}
% \forall i,j \in I
%\qquad
$\langle \ell_i, w_j \rangle = \delta_{ij}$ for all $i,j \in I$
%\end{equation*}
and, writing $I_z \definedas \{ (S,k) \in I \mid S \ni z \}$ for all $z \in \vertices$,
\begin{equation*}
 \sum_{i \in I_z} \| \ell_i \|_{H^{-1}(\omega_z)} \| \nabla w_i \|_{L^2(\omega_z)}
 \leq
 C_\mathrm{bOS}^\star,
\end{equation*}
where the constant $C_\mathrm{bOS}^\star$ depends only on $d$, $m_1$, $m_2$ and the shape regularity coefficient $\sigma$ from \eqref{E:shape-regularity}.
\end{lemma}

\begin{proof}
\fbox{\scriptsize 1} We first establish the biorthogonality. Thanks to the transformation rule, the scaled pull-backs \eqref{orthonormal-bases-on-faces-and-elms} indeed form local 
orthonormal bases of $P_{m_1}(F)$, $F \in \faces$, and $P_{m_2}(T)$, $T \in \grid$: 
\begin{equation}
\label{orthonormality-on-faces-and-elms}
 \int_{F} q_{(F,k)} q_{(F,l)} \phi_{F}
 =
 \int_{\wh{F}} q_{(\wh{F},k)} q_{(\wh{F},l)} \phi_{\wh{F}}
 =
 \delta_{kl},
\quad
 \int_{T} q_{(T,k)} q_{(T,l)} \phi_{T}
 =
 \delta_{kl}
\end{equation}
for all admissible $k$ and $l$. This orthonormality, combined with the local supports of the pairs $(\ell_i,w_i)$, $i \in I$, shows the biorthogonality, except for the cases when $(T,F) \in \grid \times \faces$ with $F \subset T$ and $k,l$ are admissible. Here transformation rule and the definition of $\wh{Q}$ imply
\begin{align*}
 \langle \ell_{(T,k)}, w_{(F,l)} \rangle
 &=
 \int_T q_{(T,k)} w_{(F,l)}
 =
 \left( \frac{|T| \, |\wh{F}|}{ |\wh{T}| \, |F|}\right)^{\frac{1}{2}}
  \int_{\wh{T}} q_{(\wh{T},k)} \wh{w}_{(\wh{F},l)}
\\
 &=
 \left( \frac{|T| \, |\wh{F}|}{ |\wh{T}| \, |F|}\right)^{\frac{1}{2}}
 \int_{\wh{T}} q_{(\wh{T},k)}
  \left( \wh{v}_{(\wh{F},l)} - \big( \wh{Q}\wh{v}_{(\wh{F},l)} \big) \phi_{\wh{T}} \right)
 =
 0
\end{align*}
and biorthogonality is verified.

\smallskip\fbox{\scriptsize 2} It remains to show the stability bound for any vertex $z \in \vertices$. Given $i \in I_z$, we have either $i=(T,k)$ with $T \in\grid$ or $i=(F,k)$ with $F \in \faces$. On the one hand, the functional $\ell_i$ satisfies
\begin{equation*}
 \| \ell_{(T,k)} \|_{H^{-1}(\omega_z)}
 \Cleq
 h_T
\quad\text{or}\quad
 \| \ell_{(F,k)} \|_{H^{-1}(\omega_z)}
 \Cleq
 h_F^{\frac{1}{2}}
\end{equation*}
since, by passing to the reference element and using orthonormality \eqref{orthonormality-on-faces-and-elms}, (a variant of the) Poincar\'e inequality \eqref{E:Poincare} and the trace inequality on $\wh{T}$, we have  
\begin{align*}
 \langle \ell_{(T,k)}, w \rangle
 &=
 \int_T q_{(T,k)} w
 =
 \left( \frac{|T|}{|\wh{T}|} \right)^{\frac{1}{2}} \int_{\wh{T}} q_{(\wh{T},k)} G_T^* w
\\
 &\leq
 \left( \frac{|T|}{|\wh{T}|} \right)^{\frac{1}{2}} \| q_{(\wh{T},k)} \|_{L^2(\wh{T})}
  \| G_T^* w \|_{L^2(\wh{T})}
\\
 &\Cleq
  \left( \frac{|T|}{|\wh{T}|} \right)^{\frac{1}{2}}
  \left( \int_{\wh{T}} |q_{(\wh{T},k)}|^2 \phi_{\wh{T}} \right)
  \| \nabla (G_T^* w) \|_{L^2(\wh{T})}
 \Cleq
 h_T \| \nabla w \|_{L^2(T)}
\end{align*}
or, with $T\in\grid$ such that $T \supset F$,
\begin{align*}
 \langle \ell_{(F,k)}, w \rangle
 &=
 \int_F q_{(F,k)} w
 =
 \left( \frac{|F|}{|\wh{F}|} \right)^{\frac{1}{2}} \int_{\wh{F}} q_{(\wh{F},k)} G_{(T,F)}^* w
\\
 &\leq
   \left( \frac{|F|}{|\wh{F}|} \right)^{\frac{1}{2}} \| q_{(\wh{F},k)} \|_{L^2(\wh{F})}
  \| G_{(T,F)}^* w \|_{L^2(\wh{F})}
\\
 &\Cleq
 \left( \frac{|F|}{|\wh{F}|} \right)^{\frac{1}{2}}
 \left( \int_{\wh{F}} |q_{(\wh{F},k)}|^2 \phi_{\wh{F}} \right)
 \| \nabla \big( G_{(T,F)}^* w \big) \|_{L^2(\wh{T})}
\\
 &\Cleq
 \left( \frac{|F| \, |\wh{T}|}{|\wh{F}| \, |T|} \right)^{\frac{1}{2}}
 h_T \| \nabla w \|_{L^2(T)}
 \Cleq
 h_F^{\frac{1}{2}} \| \nabla w \|_{L^2(T)}.
\end{align*}
On the other hand, we obtain that the function $w_i$ verifies
\begin{equation*}
 \| \nabla w_{(T,k)} \|_{L^2(T)}
 =
 \left( \frac{|\wh{T}|}{|T|} \right)^{\frac{1}{2}}
  \| \nabla \big( H_T^*\wh{w}_{(\wh{T},k)} \big) \|_{L^2(T)}
 \Cleq
 h_T^{-1} \| \wh{w}_{(\wh{T},k)} \|_{L^2(\wh{T})}
 \Cleq
 h_T^{-1}
\end{equation*}
or
\begin{align*}
 \| \nabla w_{(F,k)} \|_{L^2(T)}
 &=
 \left( \frac{|\wh{F}|}{|F|} \right)^{\frac{1}{2}}
  \| \nabla \big( H_{(T,F)}^*\wh{w}_{(\wh{T},k)} \big) \|_{L^2(T)}
\\
 &\Cleq
 h_T^{-1} \left( \frac{|T|}{|F|} \right)^{\frac{1}{2}} \| \wh{w}_{(\wh{T},k)} \|_{L^2(\wh{T})}
 \Cleq
 h_F^{-\frac{1}{2}}.
\end{align*}
Using these inequalities to bound $\| \ell_i \|_{H^{-1}(\omega_z)} \| \nabla w_i \|_{L^2(\omega_z)}$, summing over all $i \in I_z$ then establishes the stability bound as the cardinality $\#I_z$ is uniformly bounded in terms of the shape regularity coefficient $\sigma$. 
\end{proof}

\begin{corollary}[projections as interpolation operators]
\label{C:proj-as-interpolators}
The biorthogonal system $(\ell_i,w_i)$, $i \in I$, induces the $H^{-1}$-projection $P_\grid$ from Definition~\ref{D:Pgrid} and its adjoint $P_\grid^*$. Indeed, we have
\begin{equation*}
 P_\grid \ell
 =
 \sum_{i \in I} \langle \ell, w_i \rangle \ell_i
\quad\text{and}\quad
 P_\grid^* w
 =
 \sum_{i \in I} \langle \ell_i, w \rangle w_i
\end{equation*}
for all $\ell \in H^{-1}(\Omega)$ and $w \in H^1_0(\Omega)$. The stability of the biorthogonal system then provides an alternative proof of the $H^{-1}$-stability on stars of both projection $P_\grid$ and $P_\grid^*$, entailing $C_\mathrm{lStb}^\star \leq C_\mathrm{bOS}^\star$.
\end{corollary}

\begin{proof}
\fbox{\scriptsize 1} We only show the identity for $P_\grid$; the one for $P_\grid^*$ can be verified along the same lines. The biorthogonality in Lemma~\ref{L:biorthogonal-system} readily implies 
\begin{equation*}
 \left\langle \sum_{i \in I} \langle \ell, w_i \rangle \ell_i, w_j \right\rangle
 =
 \sum_{i \in I} \langle \ell, w_i \rangle \langle \ell_i, w_j \rangle
 =
 \langle \ell, w_j \rangle
 \quad \forall j \in I.
\end{equation*}
As $w_j, \j \in I$ is a basis of $\mathbb V_\grid^+$, we conclude the claimed identity for $P_\grid$.

\smallskip\fbox{\scriptsize 2} To verify the stability statement, we again restrict ourselves to the case of the projection $P_\grid$. Observe first that the proof of the stability of the biorthogonal system does invoke the local stability of $P_\grid$. Thanks to the representation of $P_\grid$ and the stability of the biorthogonal system in Lemma~\ref{L:biorthogonal-system}, we have
\begin{equation*}
 \langle P_\grid \ell, w \rangle
 =
 \sum_{i \in I_z} \langle \ell, w_i \rangle \langle \ell_i, w \rangle
%\\
% &\leq
% \sum_{i \in I} \| \ell \|_{H^{-1}(\omega_z)} \| \nabla w_i \|_{L^2(\omega_z)}
%  \|\ell_i\|_{H^{-1}(\omega_z)} \| \nabla w \|_{L^2(\omega_z)}
 \leq
 C_\mathrm{bOS}^\star  \| \ell \|_{H^{-1}(\omega_z)} \| \nabla w \|_{L^2(\omega_z)},
\end{equation*}
where $I_z:= \{ (S,k) \in I \ | S \ni z \}$. 
The proof is finished.
\end{proof}

The following two remarks illustrate the practical and theoretical usefulness of the representation formulae.
\begin{remark}[alternative computation of $P_\grid$]
\label{R:alt-computation-of-Pgrid}
A by-product of Corollary~\ref{C:proj-as-interpolators} is a way of computing $P_\grid\ell$ for a given functional $\ell \in H^{-1}(\Omega)$ that ``diagonalizes'' the approach in Remark~\ref{R:Computation-of-P}. In fact, given reference orthonormal bases as in \eqref{ref-orthonormal-bases}, we can compute the functionals $\ell_i$, $i\in I$, and test functions $w_i$, $i \in I$, by means of the formulae \eqref{ref-test-functions}, \eqref{orthonormal-bases-on-faces-and-elms}, \eqref{basis-for-Fgrid}, \eqref{basis-for-Vgrid+} whence, evaluating $\langle \ell, w_i \rangle$, $i \in I$, everything in the representation of $P_\grid\ell$ in Corollary~\ref{C:proj-as-interpolators} is at our disposal. 
\end{remark}

\begin{example}[global instability of $P_\grid$ and $P_\grid^*$]
\label{E:glb-instab-of-proj}
While the projections $P_\grid$ and $P_\grid^*$ are locally stable, both may become globally unbounded under mesh refinement. To see this, recall \eqref{local-stability-of-Pgrid*}, note $\| P_\grid \|_{\mathcal{L}(H^{-1}(\Omega))} = \| P_\grid^* \|_{\mathcal{L}(H^1_0(\Omega))}$ and, following the spirit of an example in \cite{Tantardini.Veeser.Verfuerth}, consider
\begin{equation*}
 w \definedas \sum_{z\in\vertices \cap \Omega} \phi_z \in H^1_0(\Omega).
\end{equation*}
Then, for all {\it quasi-uniform} meshes $\grid$ with shape regularity coefficient $\sigma$, there is a constant $C$ depending on $\sigma$ and quasi-uniformity such that
\begin{equation}
\label{lw-bd-norm-of-Pgrid*}
 \| P_\grid^* \|_{\mathcal{L}(H^1_0(\Omega))}^2
 \geq
 \frac{ \| \nabla (P_\grid^* w) \|_{L^2(\Omega)}^2 }{ \| \nabla w \|_{L^2(\Omega)}^2 }
 \geq
 C \frac{ \#\{ T \in \grid \mid T \cap \partial\Omega = \emptyset \} }
  { \#\{ T \in \grid \mid T \cap \partial\Omega \neq \emptyset \} }.
\end{equation}
Obviously, the last term tends to $\infty$ under uniform refinement.

\smallskip\noindent To prove \eqref{lw-bd-norm-of-Pgrid*}, we proceed in several steps, mostly hiding constants depending on quasi-uniformity of $\grid$ and, as usual, the shape regularity coefficient $\sigma$ and $d$.

\smallskip\fbox{\scriptsize 1}  We first bound $\| \nabla w \|_{L^2(\Omega)}$ from above. Noting that
\begin{equation*}
 w = 1
\quad\text{on }
 \bigcup_{T \cap \partial\Omega = \emptyset} T,
\end{equation*}
the bound $\|\nabla\phi_z\|_{L^\infty(T)} \Cleq h_T^{-1}$ readily implies
\begin{equation}
\label{nablawleq}
 \| \nabla w \|_{L^2(\Omega)}^2
 =
 \sum_{T \cap \partial\Omega \neq \emptyset} \| \nabla w \|_{L^2(T)}^2
% \leq
% C  \sum_{T \cap \partial\Omega \neq \emptyset} h_T^{-2} |T|
 \Cleq
 \# \{ T \in \grid \mid T \cap \partial\Omega \neq \emptyset \} \, h_\grid^{d-2},
\end{equation}
where $h_\grid$ stands for the meshsize of $\grid$.

\smallskip\fbox{\scriptsize 2} The lower bound for $\| \nabla (P_\grid^* w) \|_{L^2(\Omega)}$ is more involved. We start by showing the following representation for any $T \in \grid$ with $T \cap \partial\Omega = \emptyset$:
\begin{equation}
\label{Pgrid*|T=}
 P_\grid^* w|_{T}
 =
 H_T^*\wh{v}
\end{equation}
with the fixed function
\begin{equation*}
 \wh{v}
 \definedas
 \sum_{(\wh{T},k) \in \wh{I}} \left( \int_{\wh{T}} q_{(\wh{T},k)} \right)  \wh{w}_{(\wh{T},k)}
 %q_{(\wh{T},k)} \phi_{\wh{T}}
 +
 \sum_{(F',k)}
 \left( \frac{|F'|}{|\wh{F}|} \right)^{\frac{1}{2}}
 \left( \int_{\wh{F}} q_{(\wh{F},k)} \right) \wh{w}_{(F',k)}
 \not\in
 \P_0,
\end{equation*}
where the indices of the second sum vary according to $F' \subset \wh{T}$, $k=1,\dots,K_1$ and $\wh{w}_{(F',k)}$ is given by \eqref{basis-for-Vgrid+} with the transformation $H_{(\wh{T},F')}$.
Note first that, thanks to $w = 1$ on $T$ and \eqref{orthonormal-bases-on-faces-and-elms}, the coefficients in the expansion of $P_\grid^* w|_{T}$ satisfy
\begin{equation*}
  \langle \ell_{(T,k)}, w \rangle
  =
  \int_T q_{(T,k)}
  =
  \left( \frac{|T|}{|\wh{T}|} \right)^{1/2} \int_{\wh{T}} q_{(\wh{T},k)}
\end{equation*}
and, for any $F \subset T$,
\begin{equation*}
  \langle \ell_{(F,k)}, w \rangle
  =
  \int_F q_{(F,k)}
  =
  \left( \frac{|F|}{|\wh{F}|} \right)^{1/2} \int_{\wh{F}} q_{(\wh{F},k)}.
\end{equation*}
Combining these identities with \eqref{basis-for-Vgrid+} yields the claimed identity \eqref{Pgrid*|T=} and it remains to verify $\wh{v} \not\in \P_0$. Suppose $\wh{v} = c \in \mathbb R$. As a consequence, for any face $F' \subset \wh{T}$ and $k \in \{1, \dots, K_1\}$, we have
\begin{equation*}
 c
 =
 \wh{w}_{(F',k)}
 =
 |\wh{F}|^{1/2} |F'|^{-1/2} \wh{w}_{(\wh{F},k)}
 =
 |\wh{F}|^{1/2} |F'|^{-1/2} c.
\end{equation*}
As not all faces of the reference simplex have the same volume, this yields $c=0$. From \eqref{ref-test-functions} and \eqref{ref-orthonormal-bases}, we infer that the coefficients in the definition of $\wh{v}$ vanish. In particular, $\int_{\wh{T}} q_{(\wh{T},k)} = 0$ for all $k=1,\dots,K_2$ means $\wh{Q}1=0$, where $\wh{Q}$ is the operator given in \eqref{ref-Q}. This however is a contradiction because the restriction of $\wh{Q}$ to $\P_{m_2}$ is injective. Hence, $\wh{v} \not\in \P_0$ is proven.

\smallskip\fbox{\scriptsize 3} We are ready to show the bound for $\| \nabla(P_\grid^*w) \|_{L^2(\Omega)}$. Given any element $T \in \grid$ with $T \cap \partial\Omega = \emptyset$,  we pass to the reference element to exploit the previous step and obtain
\begin{align*}
 \| \nabla(P_\grid^*w) \|_{L^2(T)}
 =
 \| \nabla (H_T^* \wh{v}) \|_{L^2(T)}
 \gtrsim
 h_T^{d/2-1} \| \nabla \wh{v} \|_{L^2(\wh{T})}
\end{align*}
with $\| \nabla \wh{v} \|_{L^2(\wh{T})} > 0$ independent of $T$. Consequently, 
\begin{equation*}
 \| \nabla(P_\grid^*w) \|_{L^2(\Omega)}^2
 \geq
 \sum_{T \cap \partial\Omega = \emptyset} \| \nabla(P_\grid^*w) \|_{L^2(T)}^2
 \gtrsim
 \#\{ T \in \grid \mid T \cap \partial\Omega = \emptyset \} \, h_\grid^{d-2}
\end{equation*}
because $\grid$ is quasi-uniform. Combining this lower bound with the upper bound \eqref{nablawleq} of the first step, we conclude \eqref{lw-bd-norm-of-Pgrid*}.
\end{example}

\subsubsection*{A hierarchical estimator}
%
% intro
Like estimators based upon local problems, hierarchical estimators aim at softening the impact of constants in the lower bound, with the difference that they are explicit. While global higher order extensions were used originally, \cite{Bornemann.Erdmann.Kornhuber:1996_aposteriori} use an extension tailored to the residual structure and derive an upper bound with indicators testing the residual with a basis of the extension. One may expect that such explicit indicators come at the price of increased constants in the upper bound. For the following example, this expectation is confirmed by the inequality $C_\mathrm{lStb}^\star \leq C_\mathrm{bOS}^\star$.

Given the Galerkin approximation $u_\grid$ from \eqref{setting-apost:Galerkin}, the hierarchical PDE indicator is defined by
\begin{subequations}
\label{hier-est}
\begin{equation}
\label{hier-est;pde-ind}
 \eta^\mathrm{hier}_\grid(u_\grid)
 \definedas
 \sum_{z \in \vertices}  \eta^\mathrm{hier}_\grid(u_\grid,z)^2
\quad\text{with}\quad
 \eta^\mathrm{hier}_\grid(u_\grid,z)
 \definedas
 \max_{i \in I_z}
  \frac{|\langle R_\grid, w_i \rangle|}{\| \nabla w_i \|_{L^2(\omega_z)}}
\end{equation}
with $I$ and $I_z$ as in Lemma~\ref{L:biorthogonal-system} (biorthogonal system). Note that the test functions $w_i$, $i\in I$, are available, cf.\ Remark \ref{R:alt-computation-of-Pgrid}, and therefore $\eta^{\mathrm{hier}}_\grid(u_\grid)$ is explicit. The resulting estimator is then
\begin{equation}
 \est^\mathrm{hier}_\grid
 \definedas
 \est^\mathrm{hier}_\grid(u_\grid,f)^2
 \definedas
 \eta^\mathrm{hier}_\grid(u_\grid)^2
 +
 \osc_\grid(f)^2.
\end{equation}
\end{subequations}

\begin{theorem}[hierarchical estimator]
\label{T:hier-est}
Suppose the coefficients $\bA$ and $c$ are discrete. The hierarchical estimator \eqref{hier-est} is equivalent to the error, while its PDE indicator is locally equivalent to the discretized residual with constant 1 in the lower bound:
\begin{equation*}
 \frac{1}{C_\mathrm{lStb}^\star C_d} \est^\mathrm{hier}_\grid
 \Cleq
 \| \nabla(u-u_\grid) \|_{L^2(\Omega)}
 \Cleq
 C_\mathrm{bOS}^\star C_d C_\mathrm{loc} \est^\text{hier}_\grid
\end{equation*}
and, for all vertices $z\in\vertices$,
\begin{equation*}
 \eta^\text{hier}_\grid(u_\grid,z)
 \leq
 \| P_\grid R_\grid \|_{H^{-1}(\omega_z)}
 \leq
 C_\mathrm{bOS}^\star \eta^\text{hier}_\grid(u_\grid,z),
\end{equation*}
where $C_\mathrm{lStb}^\star$ is the stability constant of $P_\grid$ on stars from Lemma~\ref{L:loc-of-H^-1-norm}, $C_\mathrm{loc}$ from Corollary~\ref{C:residual-loc}, $C_d = \sqrt{2(d+1)}$, and the hidden constants depend only on the error-residual relationship in Lemma~\ref{L:err-residual}.
\end{theorem}

\begin{proof}
It suffices to verify the local equivalence for the PDE indicator; cf.\ Theorem~\ref{T:vertex-mod-res-est} (vertex-indexed modified residual estimation). Its lower bound simply follows from \eqref{testing-with-Vgrid+}: for all $i \in I_z$, we have
\begin{equation*}
 |\langle R_\grid, w_i \rangle|
 =
 |\langle P_\grid R_\grid, w_i \rangle|
 \leq
 \| P_\grid R_\grid \|_{H^{-1}(\omega_z)} \| \nabla w_i \|_{L^2(\omega_z)}.
\end{equation*}
To show its upper bound, let $w \in H^1_0(\omega_z)$ and, with the help of Corollary \ref{C:proj-as-interpolators} (projections as interpolation operators) and Lemma~\ref{L:biorthogonal-system} (biorthogonal system), we derive
\begin{align*}
 \langle P_\grid R_\grid, w \rangle
 &=
 \sum_{i \in I_z} \langle R_\grid, w_i \rangle \langle \ell_i, w \rangle
\\
 &\leq
 \sum_{i \in I_z}
  \frac{\langle R_\grid, w_i \rangle}{\| \nabla w_i \|_{L^2(\omega_z)}}
  \| \nabla w_i \|_{L^2(\omega_z)} \| \ell_i \|_{H^{-1}(\omega_z)} \| \nabla w \|_{L^2(\omega_z)}
\\
 &\leq
 C_\mathrm{bOS}^\star \, \eta^\mathrm{hier}_\grid(u_\grid,z) \, \| \nabla w \|_{L^2(\omega_z)}
\end{align*}
and the local equivalence is established.
\end{proof}

\begin{remark}[different test functions]
\label{R:different-test-functions-for-hier-est}
The hierarchical estimator \eqref{hier-est} does not generalize the one in \cite{Bornemann.Erdmann.Kornhuber:1996_aposteriori} as it uses slightly different test functions for edges. The given framework however applies to their variant, too; cf.\ \cite[Sect. 4.1]{KreuzerVeeser:2021}.  
\end{remark}

\subsubsection*{Alternative localization and residual splitting}
%
% intro
Lemma~\ref{L:loc-of-H^-1-norm} (localization of $H^{-1}$-norm) is not well suited for reducing or avoiding constants in the upper bounds. The following modification however allows this.

\smallskip We replace the local spaces $H^1_0(\omega_z)$, $z \in \vertices$, with
\begin{equation*}
 \W_z
 \definedas
 \begin{cases}
  \big\{ w \in H^1(\omega_z) \mid \int_{\omega_z} w = 0 \big\},
  &\text{if }z \in \vertices \cap \Omega,
 \\
  \big\{ w \in H^1(\omega_z) \mid w = 0 \text{ on } \partial\omega_z \cap \partial\Omega \big\},
  &\text{if }z \in \vertices \cap \partial\Omega,
 \end{cases}
\end{equation*}
endow them with the norm $\| \nabla \cdot \|_{L^2(\omega_z)}$, and denote by $\W_z^*$ the respective dual spaces endowed in turn with
\begin{equation}
 \| \ell \|_{\W_z^*}
 \definedas
 \sup \big\{
  \langle \ell, w \rangle \mid w \in \W_z,\,  \| \nabla w \|_{L^2(\omega_z)} \leq 1
 \big\}.
\end{equation}

\begin{lemma}[alternative localization of $H^{-1}$-norm]
\label{L:alt-loc-of-H^-1-norm}
Let $\ell \in H^{-1}(\Omega)$ be any linear functional.
\begin{enumerate}
\item If $\langle \ell, \phi_z \rangle = 0$ for all interior vertices $z \in \vertices \cap \Omega$, then
\begin{equation*}
 \| \ell \|_{H^{-1}(\Omega)}^2
 \leq
 (d+1) \sum_{z \in \vertices} \| \phi_z \ell \|_{\W_z^*}^2.
\end{equation*}
\item We have
\begin{equation*}
  \sum_{z \in \vertices} \| \phi_z\ell \|_{\W_z^*}^2
  \leq
 (d+1) C_\mathrm{loc}^2 \| \ell \|_{H^{-1}(\Omega)}^2,
\end{equation*}
where $C_\mathrm{loc}$ is the constant in Lemma~\ref{L:loc-of-H^-1-norm}~(i).
\end{enumerate}
\end{lemma}

\begin{proof}
The proof is essentially a regrouping of the arguments in Lemma~\ref{L:loc-of-H^-1-norm}, where the constant $C_\mathrm{loc}$ in the stability bound \eqref{def-Cloc} now arises in the proof of the lower bound from the following argument: we have
\begin{equation}
\label{H-1-is-subset-of-W_z^*}
 \| \phi_z \ell \|_{\W_z^*}
 \leq
 C_{\mathrm{loc}} \| \ell \|_{H^{-1}(\omega_z)}.
\end{equation}
thanks to
\begin{equation*}
\begin{aligned}
 \langle \phi_z \ell, w \rangle
 &=
 \langle \ell, \phi_z w \rangle
 \leq
 \| \ell \|_{H^{-1}(\omega_z)} \| \nabla (w\phi_z) \|_{L^2(\omega_z)}
\\
 &\leq
 C_{\mathrm{loc}} \| \ell \|_{H^{-1}(\omega_z)} \| \nabla w \|_{L^2(\omega_z)},
\end{aligned}
\end{equation*}
for all $w \in \W_z$.
\end{proof}

The question arises whether the inequality \eqref{H-1-is-subset-of-W_z^*} between the two local dual norms can be reversed. The following lemma reveals that this is only partially possible, covering discrete functionals as arguments.
%Its proof, as well as the definition of the estimator based on flux equilibration, employ the operator

%
\begin{lemma}[partial equivalence for local dual norms]
\label{L:partial-equiv-local-dual-norms}
If $z \in \vertices \cap \Omega$ is an interior vertex, the functional $\ell = \phi_z^{-1}$ satisfies
\begin{equation*}
 \| \phi_z \ell \|_{\W_z^*} = 0
\quad\text{and}\quad
 \| \ell \|_{H^{-1}(\omega_z)} > 0.
\end{equation*}
Furthermore, for any vertex $z\in\vertices$,
\begin{equation*}
 \| \ell \|_{H^{-1}(\omega_z)}
 \leq
 C_{\F} \| \phi_z \ell \|_{\W_z^*}
 \quad \forall \ell \in \F(\grid_z),
\end{equation*}
where the constant $C_{\F}$ depends only on $d$, the shape regularity coefficient $\sigma$, and the degrees $m_1$ and $m_2$ of the discrete functionals.
\end{lemma}

\begin{proof}
\fbox{\scriptsize 1} We show the claims on the functional $\ell = \phi_z^{-1}$ for an interior vertex $z \in \vertices \cap \Omega$. By the definition of $\W_z$, we have, for all $w \in \W_z$,
\begin{equation*}
 \langle \phi_z \ell, w \rangle
 =
 \int_{\omega_z} w = 0,
\end{equation*}
whence $\phi_z \ell \in \W_z^*$ with $ \| \phi_z \ell \|_{\W_z^*} = 0$.

To verify that $\ell = \phi_z^{-1} \in H^{-1}(\omega_z)$, we write $d_z \definedas \operatorname{dist}(\cdot,\partial\omega_z)$ for the distance function of the star boundary and shall use the weighted Poincar\'e inequality
\begin{equation*}
 \forall w \in H^1_0(\omega_z)
\quad
 \| w d_z^{-1} \|_{L^2(\omega_z)}
 \Cleq
 \| \nabla w \|_{L^2(\omega_z)},
\end{equation*}
which follows from the Hardy inequality; cf.\ \cite[Lemma~3.6]{SacchiVeeser:06}. Consequently, exploiting also $d_z \leq \phi_z$ on $\omega_z$, we obtain, for all $w \in H^1_0(\Omega)$,
\begin{align*}
 \langle \ell, w \rangle
 =
 \int_{\omega_z} (\phi_z^{-1} d_z) (w d_z^{-1})
 \leq
 |\omega_z|^{1/2} \| w d_z^{-1} \|_{L^2(\omega_z)}
 \Cleq
|\omega_z|^{1/2} \| \nabla w \|_{L^2(\omega_z)}.
\end{align*}
This and $\langle \ell,\phi_z \rangle = |\omega_z|$ ensure $\ell \in H^{-1}(\omega_z)$ with $\| \ell \|_{H^{-1}(\omega_z)} > 0$. % \geq |\omega_z|/\| \nabla \phi_z \|_{L^2(\omega_z)}$.

\smallskip\fbox{\scriptsize 2} We start the proof of the asserted inequality by checking that $\| \cdot \|_{H^{-1}(\omega_z)}$ is a norm on $\F(\grid_z)$. To this end, consider $q_F \in \P_{m_1}(F)$, $F \in \faces_z$ and $q_T \in \P_{m_2}(T)$, $T \in \grid_z$ such that, for all $w \in H^1_0(\omega_z)$,
\begin{equation*}
 0
 =
 \langle \ell, w \rangle
 \definedas
 \sum_{F \in \mathcal{F}_z} \int_F q_F w + \sum_{T \in \grid_z} \int_T q_T w.
\end{equation*}
We need to show $\ell = 0$. Testing with $w \in H^1_0(T)$, $T \in \grid_z$, the fundamental lemma of the calculus of variations yields $q_T = 0$ for all $T \in \grid_z$. Similarly, testing now with $w \in H^1_0(\omega_F)$, $F \in \faces_z$, gives $q_F = 0$ for all $F \in \grid_z$. Thus, $\ell = 0$ holds.

\smallskip\fbox{\scriptsize 3} Next, we check that also $\| \phi_z \cdot \|_{\W_z^*}$ is a norm on $\F(\grid_z)$. This time,  consider $q_F \in \P_{m_1}(F)$, $F \in \faces_z$ and $q_T \in \P_{m_2}(T)$, $T \in \grid_z$ such that, for all $w \in \W_z$,
\begin{equation*}
 0
 =
 \langle \phi_z\ell, w \rangle
 \definedas
 \sum_{F \in \mathcal{F}_z} \int_F \phi_z q_F w + \sum_{T \in \grid_z} \int_T \phi_z q_T w,
\end{equation*}
and again, we need to conclude $\ell = 0$. If $z \in \vertices \cap \partial\Omega$ is a boundary node, we obtain $\ell = 0$ by the arguments of the previous step. We are thus left with the case $z \in \vertices \cap \Omega$ of interior nodes. Given $w \in H^1(\omega_z)$, we set $c_w \definedas \fint_{\omega_z} w$ and $c_\ell = |\omega_z|^{-1} \langle \phi_z \ell, 1 \rangle$, and observe 
\begin{equation*}
 0
 =
 \langle \phi_z \ell, w - c_w \rangle
 =
 \langle \phi_z \ell - c_\ell, w - c_w \rangle
 =
 \langle \phi_z \ell - c_\ell, w \rangle.
\end{equation*}
Hence, testing with $w \in H^1_0(T)$, $T \in \grid_z$, we deduce $\phi_z q_T = c_\ell$ on each $T \in \grid_z$. This is however only possible if $c_\ell = 0$ and $q_T = 0$ for all $T \in \grid_z$. Therefore, testing with $w \in H^1_0(\omega_F)$, $F \in \faces_z$, yields $q_F = 0$ for all $F \in \faces_z$ and $\ell = 0$ is established in general.

\smallskip\fbox{\scriptsize 4} To conclude the asserted inequality, note that $\F(\grid_z)$ has finite dimension and, for a fixed polynomial degrees $m_1$ and $m_2$, is invariant under continuous piecewise affine transformations. Furthermore, both norms scale in the same manner. We therefore can pass to reference stars and use there the equivalence of norms in finite-dimensional spaces. Transforming the inequality back from the reference star then finishes the proof.
\end{proof}

The alternative localization entails that we need to adapt Lemma~\ref{L:abstract_bds} (splitting of local residual norms). Relying on the local $H^{-1}$-stability of $P_\grid$, Lemma~\ref{L:partial-equiv-local-dual-norms} (partial equivalence for local dual norms) reveals that the adaptation has to be global.
\looseness=-1

\begin{lemma}[alternative splitting]
\label{L:alt-abstract_bds}
Using the local norms $\|\cdot\|_{\W_z^*}$, $z \in \vertices$, the residual can be split into discretized and oscillatory residual:
\begin{multline*}
  \frac{1}{(C_\mathrm{lStb}^\star){}^2 C_d^2 C_\mathrm{loc}^2} \sum_{z\in\vertices}
   \left( \| \phi_z P_\grid R_\grid \|_{\W_z^*}^2
   +
   \| \phi_z (I - P_\grid) R_\grid \|_{\W_z^*}^2 \right)
\\ %&
 \leq
  \| R_\grid \|_{H^{-1}(\Omega)}^2
  \leq
  C_d^2 \sum_{z\in\vertices} \left(
   \| \phi_z P_\grid R_\grid \|_{\W_z^*}^2
   +
   \| \phi_z (I - P_\grid) R_\grid \|_{\W_z^*}^2
   \right),
\end{multline*}
where $C_\mathrm{lStb}^\star$ is the stability constant of $P_\grid$ on stars from Lemma \ref{L:local-stability-of-P}, and $C_d = \sqrt{2(d+1)}$.
\end{lemma}

\begin{proof}
Combine the localization in Lemma~\ref{L:alt-loc-of-H^-1-norm} with the proof of Lemma~\ref{L:abstract_bds}, replacing the local norm $\|\cdot\|_{H^{-1}(\omega_z)}$ in most places but apply \eqref{H-1-is-subset-of-W_z^*} before using the local $H^{-1}$-stability of $P_\grid$. 
\end{proof}

\subsubsection*{An estimator based on flux equilibration}
%
% intro
Estimators based on flux equilibration have been designed with the goal to obtain constant $1$ in the upper bound. The principal obstruction that computation can access only a finite-dimensional part of infinite-dimensional objects like the residual norm is overcome by means of the \emph{Prager-Synge theorem}. Realizations of this approach can be found, e.g., in \cite{Ainsworth:2010}, \cite{Braess.Pillwein.Schoeberl:2009_equilibrated}, \cite{Ern.Smears.Vohralik:2017:discrete-Hidvliftings}, and \cite{Luce.Wohlmuth:2004_equlilbrated-fluxes-estimator}. 

% PDE indicator
\smallskip The definition of the PDE indicator needs some preparation. Let $d \in \{2,3\}$, as in the aforementioned works, and let $z \in \vertices$ be vertex. Given the operator $\pi_z:\{  \phi_z \ell \mid \ell \in H^{-1}(\Omega) \} \to \W_z^*$ defined by
\begin{equation*}
 \pi_z(\phi_z \ell)
 \definedas
 \begin{cases}
  \phi_z\ell - \displaystyle \frac{\langle \ell, \phi_z \rangle}{|\omega_z|},
   &\text{if } z \in \vertices \cap \Omega,
 \\
  \phi_z \ell, &\text{if } z \in \vertices \cap \partial\Omega.
 \end{cases}
\end{equation*}
and
\begin{equation*}
 \gamma_z
 \definedas
 \begin{cases}
  \partial\omega_z, &\text{if } z \in \vertices \cap \Omega,
 \\
  \partial\omega_z \setminus \partial\Omega &\text{if } z \in \vertices \cap \partial\Omega,
 \end{cases}
\end{equation*}
we introduce the local space $\mathbb{D}_z \neq \emptyset$
\begin{equation*}
 \mathbb{D}_z
 \definedas
 \big\{ \vec{\xi} \in L^2(\omega_z;\R^d) \mid
  \operatorname{div} \vec{\xi} = \pi_z(\phi_z P_\grid R_\grid) \text{ and }
   \vec{\xi} \cdot \vec{n}_F = 0 \text{ on } F, \ \forall F \subseteq \gamma_z
 \big\},
\end{equation*}
and its discretization
\begin{equation*}
 \mathbb{D}_z(\grid)
 \definedas
 \big\{ \vec{\xi} \in \mathbb{D}_z \mid
 \; \vec{\xi} \in \mathrm{RTN}_m(T) \;\; \forall T \in \grid_z
 \big\}
\end{equation*}
with the Raviart-Thomas-N\'ed\'elec elements 
\begin{equation*}
 \mathrm{RTN}_m(T)
 =
 \left\{ \vec{\xi}:T\to\R^d \mid
  \vec{\xi}(x) = \vec{q}(x) + q(x)x \text{ with } \vec{q} \in (\P_m )^d, q \in \P_m
 \right\}
\end{equation*}
of order $m \definedas \max\{m_1,m_2\}+1$. Given the Galerkin approximation $u_\grid$ from \eqref{setting-apost:Galerkin}, the PDE indicator is then given by
\begin{subequations}
\label{feq-est}
\begin{equation}
 \eta^\mathrm{feq}_\grid(u_\grid)^2
 \definedas
 \sum_{z \in \vertices} \eta^\mathrm{feq}_\grid(u_\grid)^2
\quad\text{with}\quad
 \eta^\mathrm{feq}_\grid(u_\grid)
 \definedas
 \min_{\vec{\xi} \in \mathbb{D}_z(\grid)} \| \vec{\xi} \|_{L^2(\omega_z)}
\end{equation}
and the total estimator by
\begin{equation}
 \est^\mathrm{feq}_\grid
 \definedas
 \est^\mathrm{feq}_\grid(u_\grid,f)^2
 \definedas
 \eta^\mathrm{feq}_\grid(u_\grid)^2
 +
 \| \phi_z(f - P_\grid f) \|_{\W_z^*}^2.
\end{equation}
Note that the local PDE indicators $\eta^\mathrm{feq}_\grid(u_\grid,z)$ are computable up to machine precision.
\end{subequations}

\begin{theorem}[estimator based on flux equilibration]
\label{T:feq-est}
Suppose that the coefficients $\bA$ and $c$ are discrete and that $d \in \{2,3\}$. The estimator \eqref{feq-est} based on flux equilibration is equivalent to the error, while its PDE indicator is locally equivalent to the discretized residual with constant $1$ in the upper bound for the $\| \cdot \|_{\W_z^*}$-norm, so that
\begin{equation*}
 \frac{ C_\mathbb{D}}{C_\mathrm{lStb}^\star C_d C_\mathrm{loc}}
  \est^\mathrm{feq}_\grid
 \Cleq
 \| \nabla(u-u_\grid) \|_{L^2(\Omega)}
 \Cleq
 C_d \est^\mathrm{feq}_\grid
\end{equation*}
and, for all vertices $z\in\vertices$,
\begin{equation*}
 C_\mathbb{D} \eta^\mathrm{feq}_\grid(u_\grid,z)
 \leq
 \| \phi_z P_\grid R_\grid \|_{\W_z^*}
 \leq
 \eta^\mathrm{feq}_\grid(u_\grid,z)
\end{equation*}
as well as
\begin{equation*}
 \frac{C_\mathbb{D}}{C_\mathrm{loc}} \eta^\mathrm{feq}_\grid(u_\grid,z)
 \leq
 \| P_\grid R_\grid \|_{H^{-1}(\omega_z)}
 \leq
 C_{\F} \, \eta^\mathrm{feq}_\grid(u_\grid,z),
\end{equation*}
where $C_\mathbb{D}$ depends on $d$ and the shape regularity coefficient $\sigma$, $C_\mathrm{lStb}^\star$ is the stability constant of $P_\grid$ on stars from Lemma~\ref{L:local-stability-of-P}, $C_\mathrm{loc}$ comes from Lemma~\ref{L:loc-of-H^-1-norm}, $C_d = \sqrt{2(d+1)}$, $C_{\F}$ from Lemma~\ref{L:partial-equiv-local-dual-norms}. and the hidden constants depend only on the error-residual relationship in Lemma~\ref{L:err-residual}.
\end{theorem}

\begin{proof}
\fbox{\scriptsize 1} We start by verifying the local equivalence for the $\| \cdot \|_{\W_z^*}$-norm. Let $z \in \vertices$ be any vertex. The Prager-Synge theorem on the star $\omega_z$ implies
\begin{equation*}
 \| \phi_z P_\grid R_\grid \|_{\W_z^*}
 =
 \| \pi_z(\phi_z P_\grid R_\grid) \|_{\W_z^*}
 =
 \min_{\vec{\xi} \in \D_z} \| \vec{\xi} \|_{L^2(\omega_z)};
\end{equation*}
cf., e,g, \cite[Proposition~1.40]{Verfuerth:13}. Hence the upper bound with constant $1$ readily follows the inclusion $\mathbb{D}_z(\grid) \subset \mathbb{D}_z$, while the lower bound is a consequence of the nontrivial inequality
\begin{equation*}
 C_{\mathbb{D}}  \min_{\vec{\xi} \in \D_z(\grid)} \| \vec{\xi} \|_{L^2(\omega_z)}
 \leq
  \min_{\vec{\xi} \in \D_z} \| \vec{\xi} \|_{L^2(\omega_z)},
\end{equation*}
where $C_{\mathbb{D}}$ depends only on $d$ and the shape regularity coefficient $\sigma$; cf., e.g., \cite[Theorem~7]{Braess.Pillwein.Schoeberl:2009_equilibrated} and 
\cite[Theorem~1.1]{Ern.Smears.Vohralik:2017:discrete-Hidvliftings}.

\smallskip\fbox{\scriptsize 2} We verify the local equivalence for the $\|\cdot\|_{H^{-1}(\omega_z)}$-norm. On the one hand, combining the first equivalence with \eqref{H-1-is-subset-of-W_z^*}, we obtain
\begin{equation*}
 C_{\mathbb{D}} \eta^\mathrm{feq}_\grid(u_\grid,z)
 \leq
 \| \phi_z P_\grid R_\grid \|_{\W_z^*}
 \leq
 C_\mathrm{loc}  \| P_\grid R_\grid \|_{H^{-1}(\omega_z)}.
\end{equation*}
On the other hand, using Lemma~\ref{L:partial-equiv-local-dual-norms} instead of \eqref{H-1-is-subset-of-W_z^*} yields
\begin{equation*}
 \| P_\grid R_\grid \|_{H^{-1}(\omega_z)}
 \leq
 C_{\F} \| \phi_z P_\grid R_\grid \|_{\W_z^*}
 \leq
 C_{\F} \eta^\mathrm{feq}_\grid(u_\grid,z)
\end{equation*}
and the equivalence for the $\|\cdot\|_{H^{-1}(\omega_z)}$ is verified, too.

\smallskip\fbox{\scriptsize 3} The global bounds follow by combining Lemmas~\ref{L:err-residual} (error and residual), \ref{L:alt-abstract_bds} (alternative splitting in discretized and oscillatory residual), and \ref{L:data-oscillation-reduction} (data oscillation reduction), as well as the first local equivalence.
\end{proof}

\begin{remark}[improved upper bound]
\label{R:alt-up-bd}
Applying the Prager-Synge theorem on $\Omega$, we can improve the upper bound in Theorem~\ref{T:feq-est} to
\begin{equation}
\label{feq;improved-up-bd}
 \| \nabla (u-u_\grid) \|_{L^2(\Omega)}
 \Cleq
  \| \vec{\xi}_\Omega \|_{L^2(\Omega)}
 +
 \sqrt{d+1} \left(
  \sum_{z \in \vertices} \| \pi_z\big( \phi_z (P_\grid f - f)\big) \|_{\W_z^*}^2
 \right)^{1/2}
\end{equation}
with $\vec{\xi}_\Omega \definedas \sum_{z \in \vertices} \vec{\xi}_z$, where $\vec{\xi}_z \definedas \argmin_{\xi \in \mathbb{D}_z(\grid)} \| \vec{\xi} \|_{L^2(\omega_z)}$ are the minimizing vector fields associated with the PDE indicators, extended  by $0$ off $\omega_z$.

\smallskip To see this, we derive, thanks to the partial orthogonality \eqref{partial-orth-of-residual} of the residual and Lemma~\ref{L:data-oscillation-reduction} (data oscillation reduction),
\begin{align*}
 \operatorname{div}\vec{\xi}_\Omega
 &=
 \sum_{z \in \vertices} \pi_z\big( \phi_z P_\grid R_\grid \big)
 =
 \sum_{z \in \vertices} \pi_z\big( \phi_z R_\grid \big)
 +
 \sum_{z \in \vertices} \pi_z\big( \phi_z (P_\grid R_\grid - R_\grid)\big)
\\
 &=
 \sum_{z \in \vertices} \phi_z R_\grid
 +
 \sum_{z \in \vertices} \pi_z\big( \phi_z (P_\grid f - f)\big)
 =
 R_\grid + \delta_\grid
\end{align*}
with $\delta_\grid \definedas \sum_{z \in \vertices} \pi_z\big( \phi_z (P_\grid f - f)\big)$. Hence, 
\begin{equation*}
 \| \nabla (u-u_\grid) \|_{L^2(\Omega)}
 \Cleq
 \| R_\grid \|_{H^{-1}(\Omega)}
 \leq
 \| R_\grid + \delta_\grid \|_{H^{-1}(\Omega)} + \| \delta_\grid \|_{H^{-1}(\Omega)},
\end{equation*}
inserting $\vec{\xi}_\Omega$ in the Prager-Synge theorem on $\Omega$ and Lemma~\ref{L:alt-loc-of-H^-1-norm} (alternative localization of $H^{-1}$-norm) establish the claimed bound.

\smallskip In view of the bound \eqref{feq;improved-up-bd}, the alternative local PDE indicators $\| \phi_z^{1/2} \vec{\xi}_\Omega \|_{L^2(\omega_z)}$, $z \in \vertices$, may be used in an adaptive context. Note however that this alternative does not necessarily strengthen the link with the local residual as the definition of $\vec{\xi}_\Omega$ suggests an increased overlapping in the lower bound.
\end{remark}

%-------------------------------------------------------------------------------------
%\subsection{A posteriori error estimates for other boundary conditions}
\subsection{Other boundary conditions}
\label{L:other-bc-apost}
%-------------------------------------------------------------------------------------
%
This section illustrates that the preceding analysis of homogeneous Dirichlet conditions can be adapted to other boundary conditions. In particular, we discuss
\begin{itemize}
\item \emph{Robin and Neumann boundary conditions}, as an example for variationally formulated boundary conditions,
\item the \emph{pure Neumann problem}, with its global solvability constraint,
\item \emph{non-homogeneous Dirichlet boundary conditions}, formulated in an essential manner.
\end{itemize}
\emph{Mixed boundary conditions}, suitably discretized, give rise to a posteriori error estimators combining in a straightforward manner the indicators of, for instance, the first and third of the above groups. We therefore omit further details for such a setting. 

%\medskip\noindent
%-------------------------------------------------------------------------------------
\subsubsection*{Robin and Neumann boundary conditions}
%-------------------------------------------------------------------------------------
%
% intro and error-residual relationship
The {\it Robin} bilinear form in \eqref{E:Robin-bilinear} is coercive and continuous in $\V := H^1(\Omega)$ provided its coefficient $p\ge p_0$ on an open subset of $\partial\Omega$ for some constant $p_0>0$, according to the norm equivalence \eqref{E:norm-equiv-Robin}. Consequently, \eqref{E:Robin-weak} admits a unique solution $u\in\V$. If $\V_\grid = \mathbb{S}^{n,0}_\grid$ is the subspace of $\V$ of continuous piecewise polynomial functions of degree $\le n$, then the Galerkin counterpart of \eqref{setting-apost:Galerkin} reads
\[
u_\grid\in\V_\grid: \quad \B[u_\grid, v] = \ell(v) \quad\forall v\in \V_\grid,
\]
with $\ell=f + g\delta_{\partial\Omega}\in\V^*$; %, the dual of $\V$;
cf.\ \eqref{E:Robin-bilinear}. Its residual $R_\grid \in \V^*$ is defined as
\[
\langle R_\grid,w \rangle := \ell(w) - \B[u_\grid,w] \quad w \in \V,
\]
and $\|R_\grid\|_{\V^*}$ is equivalent to the error $\|u-u_\grid\|_{H^1(\Omega)}$ due to
Lemma \ref{L:err-residual} (error and residual), whose proof easily extends to $\V$.

% localization
The global norm $\|R_\grid\|_{\V^*}$ also localizes to all stars $\omega_z$ because Galerkin orthogonality $\langle R_\grid,\phi_z\rangle = 0$ is now valid also for boundary vertices $z\in\vertices \cap \partial\Omega$. Indeed, the proof of Lemma~\ref{L:loc-of-H^-1-norm} (localization of $H^{-1}$-norm) extends with minor modifications, where the local spaces for boundary vertices $z \in \vertices \cap \partial\Omega$ are now $\{v \in H^1(\omega_z) \mid v = 0 \text{ on } \partial\omega_z \setminus \partial\Omega \}$. Also the proof of Lemma~\ref{L:alt-loc-of-H^-1-norm} (alternative localization of $H^{-1}$-norm) is easily modified, using the local space $\{ v \in H^1(\omega_z) \mid \int_{\omega_z} v  = 0 \}$ at the boundary, too.

The next key step is the construction of a projection %operator
$P_\grid:\V^*\to\F_\grid$ that mimics the projection operator $P_\grid$ of Section \ref{S:approx-of-discrete-functionals}. For that purpose, the space of discrete functionals $\F_\grid$ has to include boundary face Dirac masses $q_F \delta_F$ with densities $q_F\in \P_{m_1}(F)$ for $F\subset\partial\Omega$. Consequently, $g$ can be approximated on $\partial\Omega$ similarly to the forcing $f$ in $\Omega$, while the coefficient $p$ is at play like the coefficient $c$. Indeed, considering for simplicity only the case of discrete coefficients $(\bA,c,p)$, the condition $m_1 \geq n_p + n$ arises in addition to those in Remark~\ref{R:choosing-osc-degs;reduction-to-data-osc}. With these caveats, the tools developed in Sections \ref{S:discretize-functionals}, \ref{S:approx-of-discrete-functionals} and \ref{S:proj-onto-discrete-functionals} give rise to a suitably adapted projection $P_\grid$ to split the residual $R_\grid$ into a discretized residual $P_\grid R_\grid$ and an oscillatory residual $(f - P_\grid f) + (g\delta_{\partial\Omega} - P_\grid(g\delta_{\partial\Omega}))$. Here  $g\delta_{\partial\Omega} - P_\grid(g\delta_{\partial\Omega})$ is supported only on $\partial\Omega$ and therefore contributes only to the oscillation indicators based upon the aforementioned new local spaces for boundary stars. This modified oscillation $\osc^\mathrm{Rob}_\grid(\data)$ with $\data=(\bA,c,p,f,g)$ and be combined with any of the presented PDE indicators, but we focus on residual estimation. In fact,  the new discrete residual $P_\grid R_\grid$ leads to a definition of the PDE estimator $\eta^\mathrm{Rob}_\grid(u_\grid)$ as in \eqref{mod-res-est}, but with additional contributions related to the boundary faces. Given any boundary face $F$ of $\grid$, such a contribution reads
\[
h_F  \| \jump{\bA\nabla_\grid u_\grid} \cdot\bn_F + p \, u_\grid - P_F g \|_{L^2(F)}^2
\]
and measures the discretized Robin residual. Combining as usual the PDE estimator $\eta^\mathrm{Rob}_\grid(u_\grid)$ and oscillation $\osc^\mathrm{Rob}_\grid(\ell)$ yields the total estimator $\est^\mathrm{Rob}_\grid(u_\grid,\ell)$, whence the following variant of Theorem~\ref{T:modified-estimator} (modified residual estimator) follows:
{\it for discrete coefficients $(\bA,c,p)$, the $H^1$-error and $\est^\mathrm{Rob}_\grid(u_\grid,\ell)$ are equivalent}
\[
C_L \est^\mathrm{Rob}_\grid(u_\grid,\ell)
\le
\|u-u_\grid\|_{H^1(\Omega)}
\le
C_U \est^\mathrm{Rob}_\grid(u_\grid,\ell).
\]
The estimates in Section \ref{S:bds-for-corr} for corrections and estimator reduction extend as well.

\subsubsection*{Pure Neumann problem}
%
% intro: the new guy ...
Neumann conditions are already covered by the previous section, except for the case of the pure Neumann problem with $p = 0$ on $\partial\Omega$ in \eqref{E:Robin-weak} requiring, as key novelty, the solvability constraint $\ell(1_\Omega) = 0$, i.e.\ right-hand side applied to the constant function equal to $1$ gives $0$. For such problems, unique exact and discrete solutions exist provided we choose $\V$ to be the subspace of $H^1(\Omega)$ of functions with zero mean value and $\V_\grid$ its natural 
finite element counterpart of degree $n$.

% ... and the choice of the residual norm
The residual $R_\grid$ is defined on all $H^1(\Omega)$ and satisfies $\langle R_\grid, 1_\Omega \rangle = 0$. Combining this fact with Lemma~\ref{L:poincare-friedrichs} (second Poincar\'e inequality) and $\inf_{c \in \R} \| v - c \|_{L^2(\Omega)} = \| w \|_{L^2(\Omega)}$ with $w = v - \fint_\Omega v \in \mathbb V$, we derive
\begin{align*}
 \| R_\grid \|_{\V^*}
 \!=\!
 \sup_{w \in \V} \frac{\langle R_\grid, w \rangle}{\| \nabla w \|_{L^2(\Omega)}}
 \!\approx\!
 \sup_{w \in \V} \frac{\langle R_\grid, w \rangle}{\| w \|_{H^1(\Omega)}}
 \!=\!
 \sup_{v \in H^1(\Omega)} \frac{\langle R_\grid, v \rangle}{\| v \|_{H^1(\Omega)}}
 \!=\! \| R_\grid \|_{H^1(\Omega)^*}.
\end{align*}
Consequently, localizing $\| R_\grid \|_{H^1(\Omega)^*}$ as in the previous section, we can derive a~posteriori error estimators with suitable contributions from the boundary $\partial\Omega$. 

% ... and AFEM-TS
However, the projection $P_\grid$ from the previous section cannot be used to generate discrete data in some auxiliary problem because $\langle \ell, 1_\Omega \rangle = 0$ does not imply $ \langle P_\grid\ell, 1_\Omega \rangle = 0$ in general. Further, a simple modification like $P_\grid \ell - \langle P_\grid\ell, 1_\Omega \rangle \langle 1_\Omega,1_\Omega\rangle^{-1} 1_\Omega$ with a global correction destroys the crucial local approximation properties.

% a suitable modification of the projection
To address this issue, we modify the projection $P_\grid$ such that the new projection $\wt{P}_\grid$ enforces locally 
$\langle \wt{P}_\grid \ell, 1_\Omega \rangle = \langle \ell, 1_\Omega \rangle $ in the spirit of the construction of the Lagrange multiplier in \cite{Fierro.Veeser:2003_reg-total-var}. To this end, recall that $P_\grid$ is now defined on $H^1(\Omega)$ and its range, the discrete functionals $\F(\grid)$, includes also boundary face Dirac masses, and that the first localization involves the local spaces $\V_z \definedas \{ v \in H^1(\omega_z) \mid v = 0 \text{ on } \partial\omega_z \setminus \partial\Omega \}$, $z \in \vertices$. Given $\ell \in H^1(\Omega)^*$, set
\begin{equation}
\label{def-modification-of-projection}
 \wt{P}_\grid \ell
 \definedas
 \sum_{z \in \vertices} \phi_z \wt{P}_z\ell
\quad\text{with}\quad 
 \wt{P}_z \ell
 \definedas
% \begin{cases}
 P_\grid \ell
 - \displaystyle
 \frac{\langle P_\grid \ell - \ell, \phi_z \rangle}{\int_{\omega_z} \phi_z} 1_{\omega_z}.
% &\text{if }z \in \Omega,
% \\
% P_\grid \ell &\text{if }z \in \partial\Omega,
% \end{cases}
\end{equation}
%Note the increase of the polynomial degrees for the global operator $\wt{P}_\grid$. 
%However, $\wt{P}_\grid$ possesses enhanced properties.

%
\begin{lemma}[new projection]
\label{L:modification-of-projection}
The operator \eqref{def-modification-of-projection} is linear, local, and satisfies 
\begin{equation*}
 \langle \wt{P}_z \ell, \phi_z \rangle = \langle \ell, \phi_z \rangle
\quad
 \forall z \in \vertices
\quad\text{and}\quad
 \langle \wt{P}_\grid \ell, 1 \rangle = \langle \ell, 1 \rangle.
\end{equation*}
%as well as
%\begin{equation*}
% \wt{P}_\grid \ell = \ell
%\quad
% \forall \ell \in F_{m_1,m_2}(\grid).
%\end{equation*}
Furthermore, $\wt{P}_z$ provides near-best approximation in $\F(\grid)|_{\V_z}$ and
\begin{equation*}
 \| \ell - \wt{P}_\grid \ell \|_{H^1(\Omega)^*}^2
 \leq
 C_\mathrm{loc} \sum_{z \in \vertices} \| \ell - \wt{P}_z \ell \|_{\V_z^*}^2.
\end{equation*} 
\end{lemma}

\begin{proof}
\fbox{\scriptsize 1} We start with the algebraic properties. By the definition of $\wt{P}_z$, we have the local relationships
$ %\begin{equation*}
 \langle \ell - \wt{P}_z \ell, \phi_z \rangle = 0
$, %\end{equation*}
viz.\  $ \langle \phi_z (\ell - \wt{P}_z\ell), 1 \rangle = 0$ for all vertices $z\in\vertices$.
Summing over all vertices immediately yields the global $\langle \ell - \wt{P}_\grid \ell, 1 \rangle = 0$. 
%To verify the invariance for $\ell \in \F_{m_1,m_2}(\grid)$, note that, for such functionals, %we have $\ell = P_\grid \ell$, whence $\wt{P}_z \ell = P_\grid \ell$ on $H^1(\omega_z)$ and %$\wt{P}_\grid \ell = P_\grid \ell = \ell$.

\smallskip\fbox{\scriptsize 2} To show that $\wt{P}_z$ is near-best approximating in $\F(\grid)|_{\V_z}$, we bound its error in terms of the one of $P_\grid$. The triangle inequality readily gives
\begin{equation*}
 \| \ell - \wt{P}_z \ell \|_{\V_z^*}
 \leq
 \| \ell - P_\grid \ell \|_{\V_z^*}
 +
 \left\|
  \frac{\langle P_\grid \ell - \ell, \phi_z \rangle}{\int_{\omega} \phi_z}  1_{\omega_z}
 \right\|_{\V_z^*},
\end{equation*}
while a variant of Lemma~\ref{L:Poincare} (first Poincar\'e inequality) and the properties of $\phi_z$ deliver
\begin{multline*}
 \left\|
  \frac{\langle P_\grid \ell - \ell, \phi_z \rangle}{\int_{\omega} \phi_z}  1_{\omega_z}
 \right\|_{\V_z^*}
 \Cleq
 |\omega_z|^{-1} h_z \| \langle P_\grid \ell - \ell, \phi_z \rangle 1_{\omega_z}\|_{L^2(\omega_z)}
\\
 \Cleq
 |\omega_z|^{-1/2} h_z \| \ell - P_\grid \ell \|_{\V_z^*}
  \| \nabla \phi_z \|_{L^2(\omega_z)}
%\\ &
 \Cleq
 \| \ell - P_\grid \ell \|_{\V_z^*}.
\end{multline*}
Hence, the error of $\wt{P}_z$ is dominated by the one of $P_\grid$,
\begin{equation}
\label{tildePgrid-locally-near-best}
 \| \ell - \wt{P}_z \ell \|_{\V_z^*}
 \Cleq
 \| \ell - P_\grid \ell \|_{\V_z^*}.
\end{equation}
and the near-best approximation of $\wt{P}_z$ follows from Corollary~\ref{C:local-near-best-approx-of-P} (local near-best approximation), adapted to the setting at hand.

\smallskip\fbox{\scriptsize 3} It remains to prove the claimed inequality. Given $w \in H^1(\Omega)$, the definition of $\wt{P}_\grid$ and the first stp yield the following identity 
\begin{align*}
 \langle \ell - \wt{P}_\grid \ell, w \rangle
 &=
 \sum_{z \in \vertices} \langle \ell, w \phi_z \rangle - \langle \phi_z \wt{P}_z \ell, w \rangle
\\
 &= \sum_{z \in \vertices} \langle \ell - \wt{P}_z \ell, w \phi_z \rangle
 =
 \sum_{z \in \vertices} \langle \ell - \wt{P}_z \ell, (w - c_z)\phi_z \rangle,
\end{align*}
with $c_z = \fint_{\omega_z} w$. Proceeding as in Lemma~\ref{L:loc-of-H^-1-norm} (localization of $H^{-1}$ norm) establishes the desired inequality and concludes the proof.
\end{proof}
The operator $\wt{P}_\grid$ possesses additional enhanced global properties, which are not needed here. Lemma~\ref{L:modification-of-projection} and \eqref{tildePgrid-locally-near-best} allow us to solve auxiliary pure Neumann problems with discrete data $\wt{P}_\grid \ell$, with the option of replacing in the local indicators the restrictions of $P_\grid$ with $\wt{P}_z$.

%-------------------------------------------------------------------------------------
\subsubsection*{Non-homogeneous Dirichlet boundary conditions}
%-------------------------------------------------------------------------------------    
%
%
Let $\V:=H^1(\Omega)$ and $\V_\grid \definedas \mathbb{S}^{n,0}_\grid$ be the subspace of $\V$ of continuous piecewise
polynomials of degree $\le n$. Given Dirichlet boundary data $g\in H^{1/2}(\Gamma)$,
where $\Gamma:=\partial\Omega$ for
simplicity, recall that $u\in\V(g)=\{ v \in \V \mid v = g \text{ on }\Gamma \}$ satisfies \eqref{E:weak-nonhom-Dirichlet}. Let
$g_\grid \in \mathbb{S}_\grid^{n,0}$ be a continuous finite element 
approximation of $g$ on $\Gamma$ and 
$\V_\grid(g_\grid)$ be the subspace of $\V_\grid$ of discrete functions %$\mathbb{S}_\grid^{n,0}$
with trace $g_\grid$. The Galerkin approximation of $u$ satisfies
\begin{equation*}
u_\grid\in\V_\grid(g_\grid): \quad \B[u_\grid,v] = \langle f, v \rangle\quad\forall v\in\V_\grid(0).
\end{equation*}
The error $e_\grid=u-u_\grid$ obviously satisfies Galerkin orthogonality
\[
\B[e_\grid,v]=0 \quad \forall v\in\V_\grid(0),
\]
but in general $e_\grid = g - g_\grid \ne 0$ on $\Gamma$. We follow \cite{SacchiVeeser:06}
to derive %computable
a posteriori bounds of $\|e_\grid\|_{H^1(\Omega)}$ using minimal regularity
$g\in H^{1/2}(\Gamma)$. %\looseness=-1

% orthogonal decomposition
We start with an orthogonal decomposition of the error $e_\grid$ arising from the two equations of the problem. Let $R_G=R_G(u_\grid,f)\in H^{-1}(\Omega)$ be
the Galerkin residual already introduced in Section \ref{S:err-res}, namely
\[
\langle R_G, v \rangle = \langle f, v \rangle - \B [u_\grid,v] \quad\forall v \in \V(0)=H^1_0(\Omega),
\]
and define the Galerkin error $e_G$ as its representation in $H^1_0(\Omega)$:
\[
e_G\in H^1_0(\Omega): \quad \B[e_G,v] = \langle R_G, v \rangle \quad\forall v \in \V(0).
\]
Furthermore, let $R_D = R_D (g) = g - g_\grid \in H^{1/2}(\Gamma)$ be the \emph{Dirichlet residual}, represented by the \emph{Dirichlet error} $e_D$ defined by
\[
e_D \in \V(g-I_\grid g): \quad \B[e_D,v] = 0 \quad\forall v \in \V(0).
\]
Then $e_\grid=e_G+e_D$ and the orthogonality $\B[e_D,e_G]=0$ yields %in the energy norm of $\B$:
\[
\enorm{e_\grid}^2 = \enorm{e_G}^2 + \enorm{e_D}^2,
\]
while the derivation for homogeneous Dirichlet conditions readily provides %the equivalence
\[
\enorm{e_G} \approx \|\nabla e_G\|_{L^2(\Omega)} \approx \est_\grid(u_\grid,f),
\]
where the Galerkin estimator $\est_\grid(u_\grid,f)$ is defined by \eqref{mod-res-est}, or any other estimator from Sect.\ \ref{S:other_estimators}. It thus remains to clarify whether $\enorm{e_\grid}$ is definite in the sense $\enorm{e_\grid} = 0 \implies e_\grid =0$ and to derive suitable lower and upper bounds for $\enorm{e_D}$.

% discretization of Dirichlet values
To this end, we need to be more specific about the choice of $g_\grid$. Let
$g_\grid = I_\grid g$ be the Scott-Zhang quasi-interpolant of $g$, which is defined locally using boundary values of $g$ exclusively \cite{ScottZhang:90,BrennerScott:08} and satisfies \looseness=-1
\begin{subequations}\label{E:scott-zhang}
\begin{align}
\label{E:invariance}
\qquad\qquad
 v\in\V_\grid|_\Gamma &\Rightarrow I_\grid v = v \quad \text{on }\Gamma,  && (\textrm{invariance})
\qquad\qquad
\\
\label{E:stability}
 \|I_\grid v\|_{L^2(\Gamma)} &\lesssim \| v\|_{L^2(\Gamma)} \quad \forall v \in \V. && (\textrm{stability})
\end{align}
\end{subequations}
These two properties ensure a variant of the equivalence $\| \cdot \|_{H^1(\Omega)} \approx \| \nabla \cdot \|_{L^2(\Omega)}$ for functions with zero trace on $\Gamma$.
\begin{lemma}[equivalence for vanishing discretized trace]
\label{L:nonzero-trace}
There exists a constant $C$ depending only on the shape regularity  of $\grids$ and $\Omega$ such that
\[
  \|v\|_{H^1(\Omega)} \le C \|\nabla v\|_{L^2(\Omega)} 
\quad
 \forall v \in \V \text{ with } I_\grid v = 0 \text{ on }\Gamma.
\]
\end{lemma}
\begin{proof}
Note that the core of the claimed inequality amounts to a variant of the first Poincar\'e inequality. In view of the norm equivalence \eqref{E:norm-equiv-Robin}, it suffices to prove that $\|v\|_{L^2(\Gamma)} \lesssim \|\nabla v\|_{L^2(\Omega)}$. Letting $\bar{v}_\Omega \definedas \fint_\Omega v$, and using \eqref{E:invariance} yields $v=v-I_\grid v = (v-\bar{v}_\Omega) - I_\grid(v-\bar{v}_\Omega)$ on $\Gamma$. Consequently, \eqref{E:stability} implies
\begin{equation}
\label{bounding-L2Gamma}
 \|v\|_{L^2(\Gamma)}
 \lesssim
 \| v - \bar{v}_\Omega \|_{L^2(\Gamma)}
 \lesssim
 \| v - \bar{v}_\Omega \|_{H^1(\Omega)}
 \lesssim
 \|\nabla v\|_{L^2(\Omega)},
\end{equation}
because of Lemma~\ref{L:trace-id} (traces), and Lemma~\ref{L:poincare-friedrichs} (second Poincar\'e inequality). 
\end{proof}
Observing $I_\grid e_\grid = I_\grid g - I_\grid^2 g = 0$ on $\Gamma$, we can apply Lemma~\ref{L:nonzero-trace} to get
\begin{equation*}
%\label{measuring-egrid}
 \| e_\grid \|_{H^1(\Omega)}
 \leq
 C \| \nabla e_\grid \|_{L^2(\Omega)}
 \leq
 \frac{C}{\alpha_1} \enorm{e_\grid}
 \leq
 \frac{ C \max\{\alpha_2,\| c \|_{L^\infty(\Omega)}\} } {\alpha_1}  \| e_\grid \|_{H^1(\Omega)},
\end{equation*}
establishing in particular that $\enorm{e_\grid}$ is definite. In the same vein, we derive
%the equivalence
\begin{equation*}
 \enorm{e_D}
 \approx
 \| e_D \|_{H^1(\Omega)}
 \approx
 \| \nabla e_D \|_{L^2(\Omega)}
\end{equation*}
for the Dirichlet error.

% passing to equivalent, intrinisic description 
With the intent to achieve directly computable bounds for the Dirichlet error, we next establish the equivalence $\|e_D\|_{H^1(\Omega)} \approx \| g-I_\grid g \|_{H^{1/2}(\Gamma)}$, where the \emph{intrinsic} $H^{1/2}$-norm combines the $L^2(\Gamma)$-norm  with the seminorm
\begin{equation*}
 |v|_{H^{1/2}(\Gamma)}^2
 =
 \int_\Gamma\int_\Gamma \frac{|v(x)-v(y)|^2}{|x-y|^d} \, dx \, dy.
\end{equation*}
This equivalence follows with the help of the trace and extension theorems for $H^{1/2}(\Gamma)$, see, e.g., \cite[Theorem~6.2.40]{Hackbusch:92}. In fact, on the one hand, that trace theorem immediately gives
$ %\begin{equation*}
 \| g-I_\grid g \|_{H^{1/2}(\Gamma)}
 \Cleq
 \| e_D \|_{H^1(\Omega)}.
$ %\end{equation*}
On the other hand, let $\chi \in H^1(\Omega)$ denote the extension of $g - I_\grid g$ from \cite[Theorem~6.2.40]{Hackbusch:92}. Then
\begin{equation*}
 \| e_D \|_{H^1(\Omega)}
 \Cleq
 \enorm{e_D}
 \le
 \enorm{\chi}
 \Cleq
 \| \chi \|_{H^1(\Omega)}
 \Cleq
 \| g - I_\grid g\|_{H^{1/2}(\Gamma)},
\end{equation*}
where the second inequality is thanks to $\B[e_D,e_D-\chi]=0$.

% localization and Dirichlet indicators
We are left with the issue that the $H^{1/2}(\Gamma)$-seminorm is {\it nonlocal}. To handle this delicate matter, we invoke its localization due to \cite{Faermann:00,Faermann:02}
\begin{align*}
 |v|_{H^{1/2}(\Gamma)}^2
 \le
 \sum_{F \in \faces_\Gamma} \left(
  \int_F\int_{\omega_F} \frac{|v(x)-v(y)|^2}{|x-y|^d} \, dx \, dy
  +
  \frac{C}{h_F} \|v\|_{L^2(F)}^2
 \right),
\end{align*}
where $F\in\faces_\Gamma$ is a generic face of $\grid$ lying on $\Gamma$ and $\omega_F$ is the patch on $\Gamma$ associated with $F$. The last term seems problematic. However, applied to $v=g-I_\grid g$, we can mimic the steps of \eqref{bounding-L2Gamma} with local variants of \eqref{E:scott-zhang}, but using in the last step the second Poincar\'e inequality in $H^{1/2}$, see, e.g., \cite[Lemma~3.2]{SacchiVeeser:06}:
\[
 \|v\|_{L^2(F)}^2
 =
 \|v-I_\grid v\|_{L^2(F)}^2
 \lesssim
 \|v-\bar{v}_F\|_{L^2(\omega_F)}^2
 \lesssim
 h_F |v|_{H^{1/2}(\omega_F)}^2
\]
where $\bar{v}_F$ is the mean value of $v$ on $\omega_F$. Note that this bounds also means that the $L^2$-part in $\| g - I_\grid g \|_{H^{1/2}(\Gamma)}$ is (locally) controlled by its seminorm. Altogether, this leads to defining the \emph{Dirichlet oscillation} with the following local indicators:
\begin{equation}
\label{E:Dirichlet-osc}
\begin{aligned}
 \osc_\grid(g)^2_{1/2} &\definedas \sum_{F\in\faces_\Gamma}\osc_\grid (g,F)^2_{1/2},
\\
 \osc_\grid (g,F)^2_{1/2}
 &\definedas
 \int_{\omega_F}\int_{\omega_F}
  \frac{|(g-I_\grid g)(x)-(g-I_\grid g)(y)|^2}{|x-y|^d} \, dx \,dy.
\end{aligned}
\end{equation}
We observe that $\osc_\grid (g,F)$ is a double singular integral but computationally accessible, for instance, by using suitable quadrature provided $g$ is continuous \cite[Section 4.1]{SacchiVeeser:06}.

\begin{proposition}[Dirichlet oscillation]
\label{P:Dirichlet}
There exist constants $D_1 \ge D_2>0$ depending on the shape regularity of $\grids$ and geometry of $\Gamma$, such that
\[
 D_2 \osc_\grid(g)_{1/2}
 \le
 \|\nabla e_D\|_{L^2(\Omega)}
 \le
 D_1 \osc_\grid(g)_{1/2}.
\]
\end{proposition}
\begin{proof}
The preceding derivation verifies the upper bound. For the lower one, note that for any $v \in H^1(\Omega)$ such that $v =g - I_\grid g$ on $\Gamma$
%
%  \[
%  \|\nabla e_D\|_{L^2(\Gamma)}^2 \approx |g-I_\grid g|_{H^{1/2}(\Gamma)}^2
%    \lesssim \sum_{F\in\faces_\Gamma}\osc_\grid (g,F)^2_{1/2} = \osc_\grid (g)^2_{1/2}.
%  \]
%  % 
%  On the other hand, if $v=g-I_\grid g$ we see that
%
\begin{align*}
 \osc_\grid (g)^2_{1/2}
 =
 \sum_{F\in\faces_\Gamma}\osc_\grid (g,F)^2_{1/2}
 &\le
 \sum_{F\in\faces_\Gamma}
  \int_{\omega_F}\int_{\Gamma} \frac{|v(x)-v(y)|^2}{|x-y|^d} \, dx \, dy
\\
 &\lesssim
 \int_{\Gamma}\int_{\Gamma} \frac{|v(x)-v(y)|^2}{|x-y|^d} \, dx \, dy
 =
 |v|_{H^{1/2}(\Gamma)}^2,
\end{align*}
because the patches $\omega_F$, $ F \in \faces_\Gamma$, possess a uniform overlapping property due to shape regularity of $\grids$. Applying this to $v=e_D$ finishes the proof.
\end{proof}
For suitable settings, local lower a posteriori estimates for the Dirichlet error $e_D$ can be derived; see \cite[Theorem 3.2]{SacchiVeeser:06}.

Combining the Dirichlet oscillation with some Galerkin estimator $\est_\grid(u_\grid,f)$ by
\begin{equation*}
 \est^\mathrm{Dir}_\grid(u_\grid,f,g)^2
 \definedas
 \est_\grid(u_\grid,f)^2 + \osc_\grid(g)_{1/2}^2,
\end{equation*}
the preceding discussion is summarized by the following result.

\begin{theorem}[estimators for general Dirichlet condition]
\label{T:Dirichlet}
If Assumption~\ref{A:discrete-coefficients} (discrete coefficients) is valid, then there exist constants $C_L \le C_U$ depending on $(\bA,c)$, $\Omega$, $\Gamma$, and the shape regularity of $\grids$ such that
\[
 C_L \est^\mathrm{Dir}_\grid(u_\grid,f,g)  
 \le
 \|\nabla(u-u_\grid)\|_{L^2(\Omega)} \le
 C_U \est^\mathrm{Dir}_\grid(u_\grid,f,g).
\]
\end{theorem}

%--------------------------------------------------------------------------------

\section{Convergence of AFEM for Coercive Problems}\label{S:convergence-coercive}
% \rhn{(RHN $\longrightarrow$  CC)}

% {\color{brown}
%\begin{itemize}
%\item
%  Definition of AFEM, D\"orfler marking \cite[Section 4.1]{NochettoVeeser:2012}.
  
%\item
%  Contraction property for discrete data. Pythagoras equality \cite[Section 5]{NochettoVeeser:2012}.

%\item\index{Algorithms!\GALERKIN: Procedure that iterates \SOLVE, \ESTIMATE, \MARK, \REFINE}
%  Variable data: perturbation theory, two-step algorithm \cite{BonitoDeVoreNochetto:2013}.

%\item
%  Modules to handle data approximation: definition and properties.
  
%\item
%  Conditional contraction and complexity of second step.

%\item
%  Numerical experiments.
%\end{itemize}
%}

In this section we consider the coercive problem \eqref{strong-form} with the intent to design and analyze three $\AFEM$s in increasing order of
complexity and applicability, depending on properties of data $\data$. Our basic regularity
assumption on data reads $\data = (\bA,c,f) \in\D$, where 
\begin{equation}\label{E:space-data}
\D := L^\infty(\Omega; \R^{d\times d})  \times L^\infty(\Omega) \times H^{-1}(\Omega).
\index{Functional Spaces!$\D$: data}
\end{equation}
We approximate $\data$ with discrete data $\wh{\data} = (\wh{\bA}, \wh{c}, \wh{f})\in \D_{\wh\grid}$, where
\begin{equation}\label{E:space-discrete-data}
\D_{\wh\grid} := \big[ {\mathbb S}^{n-1,-1}_{\wh{\grid}} \big]^{d\times d}
\times {\mathbb S}^{n-1,-1}_{\wh{\grid}} \times \F_{\wh{\grid}}
\index{Functional Spaces!$\D_{\wh\grid}$: discrete data subordinate to $\wh \grid$}
\end{equation}
is subordinate to a partition $\wh\grid \in \grids$. We will often assume that {\it data is discrete},
meaning precisely that $\data = \wh{\data}$.

We start with the {\it one-step $\AFEM$}, hereafter called $\GALERKIN$,
which is the standard SEMR loop
\[
\SOLVE \,\,
\longrightarrow \,\,
\ESTIMATE \,\,
\longrightarrow \,\,
\MARK \,\,
\longrightarrow \,\,
\REFINE
\]
introduced in \cite{Dorfler:96} and further developed in \cite{MoNoSi:00,MoNoSi:02,CaKrNoSi:08}.
This is the simplest algorithm in that it requires data $\data = (\bA,c,f)$ to be discrete, but it is
a building block for the other two methods. After reviewing a few crucial properties of
error and estimator in Section \ref{S:properties-AFEM},
we fully discuss $\GALERKIN$ in Section \ref{S:discrete-data}.

The second algorithm is the {\it one-step $\AFEM$ with switch}, which still assumes the coefficients
$(\bA,c)$ to be discrete but allows for general forcing $f \in H^{-1}(\Omega)$. This is a new contribution
of this survey that, similarly to \cite{Kreuzer.Veeser.Zanotti}, exploits the structure of the error estimator $\est_\grid(u_\grid,f)$ of
Section \ref{S:aposteriori}
\[
\est_\grid(u_\grid,f)^2 = \eta_\grid (u_\grid)^2 + \osc_\grid(f)_{-1}^2,
\]
and its equivalence to the energy error.
The PDE estimator $\eta_\grid (u_\grid)$ relies on the discrete forcing $P_\grid f \in \F_\grid$
and is fully computable, whereas the data oscillation $\osc_\grid(f)_{-1}$ encodes the infinite
dimensional nature of $f$ and could be estimated in important cases of practical interest
further discussed in Section \ref{S:right-hand-side}. The quantity $\osc_\grid(f)_{-1}$ measures the deviation of $f$ from being discrete and  may dictate the pre-asymptotic
regime of $\AFEM$. Therefore, $\osc_\grid(f)_{-1}$ must be handled separately
from $\eta_\grid(u_\grid)$; hence the name of the new method, hereafter called $\AFEMSW$. Assuming
that $\osc_\grid(f)_{-1}$ is computable, the module
\[
[\wh{\grid}] = \DATA \, (\grid,f,\tau)
\]
deals with $\osc_\grid(f)_{-1}$ whenever it is large relative to $\est_\grid(u_\grid,f)$. In fact, it
creates an admissible refinement $\wh{\grid}$ of the input mesh $\grid$ such that
$\osc_{\wh{\grid}} (f)_{-1}$ is below the desired tolerance $\tau$, i.e. $\osc_{\wh{\grid}} (f)_{-1} \le \tau$.
%Given the generality
%of $f$, this is achieved by a tree algorithm\todo{AB (1/22/24): make sure this is the case} of \cite{BinevDeVore:2004,Binev:2018,BinevFierroVeeser:2023}.
We explain the role of data oscillation for error analysis, %review tree approximation,
design $\AFEMSW$ and prove its linear
convergence in Section \ref{S:one-step-switch}.

The third algorithm deals with variable data $\data$ and various degrees of regularity of $\data$,
and is able to handle discontinuous coefficients $(\bA,c)$ not aligned with admissible meshes
$\grid\in\grids$ emanating from $\grid_0$. To handle the multiplicative structure of $(\bA,c)$
in the model problem \eqref{strong-form}, we consider the following {\it two-step {\rm \AFEM}}:
\medskip
\begin{algo}[\AFEMTS]\label{algo:AFEM-Kernel}
\index{Algorithms!\AFEMTS: two step AFEM successively approximating the data and the Galerkin solution with approximate data}
Given an initial mesh $\grid_0$, an initial tolerance $\eps_0$, and a parameter $\omega$ sufficiently small to be determined later, iterate
\begin{algotab}
 \>  $\AFEMTS \, (\mesh_0,\varepsilon_0,\omega)$ \\
  \> \>  $k=0$\\
  \> \>  $[\widehat\mesh_{k},\widehat\data_{k}]= \DATA \, (\mesh_k, \data, \omega \, \varepsilon_k)$ \\
  \> \>  $[\mesh_{k+1},u_{k+1}]= \GALERKIN \, (\widehat{\mesh}_{k},\widehat\data_{k},\varepsilon_k)$ \\
  \> \>  $\varepsilon_{k+1}=\tfrac12 {\varepsilon_k}; ~ k \leftarrow k+1$
\end{algotab}
\end{algo}

This structure was first proposed in \cite{Stevenson:08} and further explored in
\cite{BonitoDeVoreNochetto:2013,CohenDeVoreNochetto:2012,BonitoCascon:2016,BonitoCascon:2013,BonitoDevaud:2015}. 
The three components of
data $\data = (\bA,c,f)\in\D$ are first approximated by discrete data
$\wh{\data}=(\wh{\bA},\wh{c},\wh{f})\in\D_{\wh\grid}$, as defined in \eqref{E:space-data} and
\eqref{E:space-discrete-data}, within the module
\[
[\wh\mesh,\wh\data]= \DATA \, (\mesh, \data, \tau)
\index{Algorithms!$\DATA$: procedure to approximates the data $\data=(\bA,c,f)$}
\]
to accuracy $\tau = \omega\eps$ significantly smaller than $\eps$.
This is achieved by an algorithm similar to
Algorithm \ref{algo:greedy} (greedy algorithm), which is fully discussed along with applications
to $\data$ in Section \ref{S:data-approx}.
The resulting admissible refinement $\wh{\grid}$ of $\grid$
and discrete data $\wh{\data}$ over $\wh{\grid}$ are next taken by $\GALERKIN$ to reduce the
PDE error to the desired tolerance $\eps$, namely the module
\[
[\grid,u_\grid]= \GALERKIN \, (\wh{\grid},\wh\data,\eps)
\index{Algorithms!\GALERKIN: procedure that iterates \SOLVE, \ESTIMATE, \MARK, \REFINE}
\]
constructs a refinement $\grid$ of $\wh{\grid}$ with discrete data $\wh\data$ over $\wh{\grid}$
such that $\eta_\grid(u_\grid) \le \eps$. We point out that if data is discrete, i.e. $\data=\wh{\data}$,
then $\DATA$ is skipped and $\AFEMTS$ reduces to $\GALERKIN$.
We tackle $\AFEMTS$ in Section \ref{S:general-data}, where we prove a perturbation estimate
with respect to $\data$ and next discuss convergence properties of $\AFEMTS$.
We will extend this approach to
discontinuous FEMs in Section~\ref{S:dg} and to mixed FEMs for \eqref{strong-form} as well as the Stokes system \eqref{strong-stokes}
in Section \ref{S:conv-rates-infsup}.
%
 
%%%%%%%%%%%%%%%%%%%%%%%%%%%%%%%%%%%%%%%%%%%%%%%%%%%%%%%%%%%%%%%%%%%%%%%%%%%%%%%%
\subsection{Properties of error and estimator}\label{S:properties-AFEM}
%%%%%%%%%%%%%%%%%%%%%%%%%%%%%%%%%%%%%%%%%%%%%%%%%%%%%%%%%%%%%%%%%%%%%%%%%%%%%%%%
%
We follow Casc\'on, Kreuzer, Nochetto, and Siebert
\cite{CaKrNoSi:08} and summarize some basic properties of \GALERKIN that
emanate from the symmetry of the differential operator (i.e. of $\bA$)
and features of the modules. In doing this, any explicit constant or hidden constant
in $\Cleq$ will only
depend on the uniform shape-regularity of $\grids$, the dimension $d$,
the polynomial degree $n$, and the (global) eigenvalues of $\bA$, but
not on a specific grid $\grid\in\grids$, except if explicitly stated. 

We recall that the bilinear form $\B$ in \eqref{bilinear-coercive} with continuous coefficients
$(\bA,c)$ is symmetric, coercive and continuous in the space $H^1_0(\Omega)$ (see \eqref{nsv-enorm-Vnorm}),
namely $\enorm{v} = \B[v,v]^{\frac12}$\index{Norms!$\enorm{.}$: energy norm with exact coefficients}
 is a norm equivalent to $|v|_{H^1_0(\Omega)}$
with equivalence constants
$0  <c_\B \le C_\B$
\index{Constants!$(c_\B,C_\B)$: norm equivalence constants}
\begin{equation}\label{E:exactB}
c_\B |v|_{H^1_0(\Omega)}^2 \le \enorm{v}^2 \le C_\B |v|_{H^1_0(\Omega)}^2
\quad\forall v \in H^1_0(\Omega).
\end{equation}
The module $\DATA$ approximates $(\bA,c)$ over
a mesh $\grid$ by piecewise polynomial coefficients $(\wh{\bA},\wh{c})$ obeying
side constraints so that the corresponding perturbed bilinear form $\wh{B}$ still defines a
uniform scalar product in $H^1_0(\Omega)$ 
\begin{equation}\label{E:scalar-prod-wh}
    \enorm{v} = \wh{\B}[v,v]^{\frac12} \quad\forall v \in H^1_0(\Omega),
\index{Norms!$\enorm{.}$: energy norm with perturbed coefficients}
\end{equation}
which satisfies \eqref{E:exactB} with constants 
$0 < c_{\wh{\B}} \le C_{\wh{\B}}$ independent of $\grid$. We hope this slight abuse of notation will
not create confusion because we will always refer to the energy norm in \eqref{E:scalar-prod-wh} when
dealing with $\wh\B$.
We denote by $\wh{u}=u(\wh \data) \in H^1_0(\Omega)$\index{Functions!$\wh u$: solution to the perturbed problem \eqref{E:perturbed-weak-form}} the solution of \eqref{weak-form}
with coefficients $(\wh{\bA},\wh{c})$ and forcing function either $\wh{f}=f\in H^{-1}(\Omega)$ or its projection
$\wh{f}=P_\grid f \in \F_\grid$ defined in \eqref{def-densities-of-P}, namely
\begin{equation}\label{E:perturbed-weak-form}
\wh{\B} [\wh{u},v] = \langle \wh{f},v \rangle
\quad\forall v\in H^1_0(\Omega),
\end{equation}

In the sequel, we will often compare discrete functions on different meshes. Given
$\grid\in\grids$, we denote by $\grid_*\in\grids$ an admissible refinement of $\grid$ and write
\begin{equation}\label{E:grid-gridstar}
\grid \le \grid_* \quad\Leftrightarrow\quad
\grids(\grid) \subset \grids(\grid_*)
\index{Meshes!$\grid \leq \grid_*$: refinement relation}
\end{equation}
in the sense that the supporting tree of $\grid$ is contained in that of $\grid_*$.
For any $\grid_*\ge\grid$, we have the following crucial property.

\begin{lemma}[Pythagoras]\label{L:Pythagoras}
{\it  Let $\grid_*\ge\grid \ge \wh \grid$ and let
  $\wh{u} \in H^1_0(\Omega)$ be the solution of \eqref{E:perturbed-weak-form} with discrete
  \rhn{coefficients $(\wh{\bA},\wh{c})$ over $\wh{\grid}$. The corresponding
  Galerkin solutions $u_\grid \in \V_\grid$ and $u_{\grid_*}\in\V_{\grid_*}$ with coefficients
  $(\wh{\bA},\wh{c})$ and forcing $f\in H^{-1}(\Omega)$}
  satisfy the orthogonality property
\begin{equation}\label{E:Pythagoras}
	\enorm{\wh{u}-v_\grid}^2=\enorm{\wh{u}-u_{\grid_*}}^2+\enorm{u_{\grid_*}-v_\grid}^2
        \quad\forall v_\grid \in \V_\grid.
\end{equation}  
}
\end{lemma}
\begin{proof}
Exploit the nestedness property $\V_\grid \subset \V_{\grid_*}$, and the Galerkin orthogonality property
$\wh{\B}[\wh{u}-u_{\grid_*},v_\grid - u_{\grid_*}]=0$ in $\V_{\grid_*}$ for the scalar product induced by $\wh{\B}$.
\end{proof}

Property \eqref{E:Pythagoras} is very restrictive: it relies on space nestedness and is valid
exclusively for the energy norm. However, it is instrumental for the subsequent analysis in
the energy norm or the equivalent norm $|\cdot|_{H^1_0(\Omega)}$,
but it does not extend to other, perhaps more practical, norms such as the maximum norm. This is an
important open problem and a serious limitation of this theory.

We recall that the residual a posteriori error analysis of Section \ref{S:aposteriori} relies on the
projection operator $P_\grid:H^{-1}(\Omega)\to\F_\grid$, with element and face components
$P_\grid f|_T = P_T f$ for $T\in\grid$ and $P_\grid f|_F = P_F f$. The {\it full local error indicator}
\[
\est_\grid(u_\grid,f,T)^2 = \eta_\grid (u_\grid,T)^2 + \osc_\grid (f,T)_{-1}^2
\]
splits into a computable {\it PDE error indicator} with discrete coefficients $(\wh{\bA},\wh{c})$
\begin{equation}\label{E:pde-estimator}
	\eta_\grid(v,T)^2=h_T^2\,\|r(v)\|^2_T+h_T\,\|j(v)\|^2_{\partial T}\quad\forall T\in\grid,
\end{equation}
where the \emph{interior} and \emph{jump residuals} are given by
\begin{equation}\label{E:residuals}
\begin{alignedat}{2}
	& r(v)|_T = P_T f + \divo(\wh{\bA}\nabla v)-\wh{c}v \quad &&\fa T\in\grid\\
	& j(v)|_{F} = [\wh{\bA}\nabla v] \cdot\bn \,|_F - P_F f  \quad &&\fa F\in\faces,
\end{alignedat}
\end{equation}
and $j(v)|_F=0$ for boundary faces $F$, and {\it data oscillation}
\begin{equation}\label{E:data-osc}
\osc_\grid(f,T)_{-1}^2 = \|f - P_\grid f\|_{H^{-1}(\omega_T)}^2
\quad\forall T\in\grid,
\end{equation}
where $\omega_T$ is the patch associated with $T$.
The corresponding global quantities are
\begin{equation}\label{E:global-est-osc}
\begin{gathered}
\est_\grid(u_\grid,f)^2 = \sum_{T\in\grid} \est(u_\grid,f,T)^2,
\\
\eta_\grid(u_\grid)^2 = \sum_{T\in\grid} \eta_\grid(u_\grid,T)^2,
\quad
\osc_\grid(f)_{-1}^2 = \sum_{T\in\grid} \osc_\grid(f,T)_{-1}^2,
\end{gathered}
\end{equation}
and have the following a posteriori error estimates proved in Theorem \ref{T:modified-estimator}
(modified residual estimator) for the $H^1_0$-norm.

\begin{proposition}[a posteriori error estimates]\label{L:ApostEnergy}
{\it Let $\wh{u} \in H^1_0(\Omega)$ be the solution of \eqref{E:perturbed-weak-form} with
discrete coefficients $(\wh{\bA},\wh{c})$ over $\grid\in\grids$ but general forcing $f\in H^{-1}(\Omega)$.
Then, there exist constants $0<C_L\le C_U$, depending on the shape regularity of $\grids$
such that the Galerkin solution $\uG\in\V_\grid$ satisfies
\begin{equation}\label{E:aposteriori-bounds-H1}
\index{Constants!$(C_L,C_U)$: a-posteriori lower and upper bounds constants}
    C_L\,\estG(\uG,f) \le |\wh{u}-\uG|_{H^1_0(\Omega)} \le C_U\,\estG(\uG,f),
  \end{equation}
}
\end{proposition}

Moreover, if $\enorm{\,\cdot \,}$ stands for the energy norm in \eqref{E:scalar-prod-wh} with equivalence
constants $c_{\wh{\B}}\le C_{\wh{\B}}$ satisfying \eqref{E:exactB}, then \eqref{E:aposteriori-bounds-H1}
yields
\begin{equation}\label{E:aposteriori-bounds}
    C_2\,\estG(\uG,f) \le \enorm{\wh{u}-\uG} \le C_1\,\estG(\uG,f),
  \end{equation}
with $C_1 = \sqrt{C_{\wh{\B}}} C_U$ and $C_2 = \sqrt{c_{\wh{\B}}} C_L$.

There is a fundamental difference between \eqref{E:aposteriori-bounds-H1} and
earlier versions of a posteriori error estimates, 
going back to the seminal paper \cite{BabuskaMiller:87}; see also 
\cite{AinsworthOden:00,Braess:07,NoSiVe:09,Verfuerth:13}. It is about the role of data
oscillation $\osc_\grid(f)_{-1}$, which is now dominated by the error $|\wh{u}-\uG|_{H^1_0(\Omega)}$
and does not spoil the lower bound. This is due to the fact that $\osc_\grid(f)_{-1}$ is evaluated in
the natural space $H^{-1}(\Omega)$ and quantifies the discrepancy between $f$ and a suitable
projection $P_\grid f$ which gives rise to a quasi-best local approximation of $f$.
We refer to \cite{NoSiVe:09} and \cite{KreuzerVeeser:2021} for a
discussion of data oscillation.

Suppose now that we have two conforming meshes $\grid, \grid_* \in \grids$ with $\grid_* \ge \grid$. Let
\begin{equation}\label{E:refined-set}
\Res:=\refined:=\grid\backslash\grid_*
\end{equation}
be the subset of refined elements of $\grid$, namely, those elements in $\grid$ 
that are no longer in $\grid_*$.
We stress that the upper bound in \eqref{E:aposteriori-bounds-H1} cannot be local due to the nonlocal nature
of the error $|\wh{u} - u_\grid|_{H^1_0(\Omega)}$. However, in view of Theorem \ref{T:ubd-corr} (upper
bound for corrections), the following remarkable local upper bound for Galerkin solutions
$u_\grid \in V_\grid, u_{\grid_*} \in V_{\grid_*}$ holds
\begin{equation}\label{E:local-upper-bound}
\enorm{u_\grid - u_{\grid_*}} \le C_1 \est_\grid(u_\grid,f,\Res),
\end{equation}
where for $\cal S \subset \grid$, $\est_\grid(u_\grid,f,{\cal S}) := \big( \sum_{T\in {\cal S}} \est_\grid(u_\grid,f,T)^2 \big)^{1/2}$ is the
error estimator restricted to $\cal S$. Consequently, solely the elements where $\grid$ and $\grid_*$
differ, namely the set $\Res$, account for the discrepancy between $u_\grid$ and $u_{\grid_*}$. This turns
out to be consistent with \eqref{E:aposteriori-bounds} because $\grid$ has to be refined everywhere to
get to $\wh{u}$, whence $\Res=\grid$.

In contrast to the upper bound in \eqref{E:aposteriori-bounds-H1}, the corresponding lower bound is
local according to Theorem \ref{T:modified-estimator} (modified residual estimator). This is due
to the local nature of the PDE \eqref{strong-form}. However, when comparing $u_\grid$ and $u_{\grid_*}$,
this bound is not valid unless the {\em interior vertex property} (given in Definition~\ref{D:interior-node}) is satisfied \cite{MoNoSi:00}; in fact, we
present a counterexample later in Example \ref{EX:interior-node} taken from \cite{MoNoSi:00}.

The interior vertex property is valid upon enforcing a fixed number $b$ of bisections ($b=3,6$ for $d=2,3$).
An immediate consequence, proved in Theorem \ref{T:lubd-corr} (local lower
bound for corrections), is the {\it discrete} lower a posteriori bound for piecewise constant
diffusion coefficient $\bA$ and reaction coefficient $c=0$ on $\grid_0$,
\begin{equation}\label{E:local-lower-discrete}
{C}_{L,1} \, \est_\grid(u_\grid,f,\marked) \le \enorm{u_\grid - u_{\grid_*}}
+ {C}_{L,2} \osc_\grid (f, \omega(\marked))_{-1},
\end{equation}
where $\omega(\marked) := \cup\{\omega_T: T\in\marked\}$ is the union of all patches of elements in
$\marked$ and $\osc_\grid (f,\omega(\marked))_{-1}^2=\sum_{T \in \omega(\marked)}\osc_\grid(f,T)_{-1}^2$; we refer to \cite{MoNoSi:00,MoNoSi:02}. We stress that if $f=P_\grid f$ is discrete,
then $\osc_\grid(f)_{-1}=0$ and \eqref{E:local-lower-discrete} reduces to
\begin{equation}\label{E:local-lower-reduced}
C_2 \, \eta_\grid(u_\grid,\marked) \le \enorm{u_\grid - u_{\grid_*}}.
\end{equation}

One serious difficulty in dealing with AFEM is that one has access to
the energy error $\enorm{\wh{u}-\uG}$, or equivalently to $|\wh{u}-\uG|_{H^1_0(\Omega)}$,
only through the full error estimator $\estG(\uG,f)$. Lemma \ref{L:Pythagoras} (Pythagoras) implies
monotonicity of the energy error with respect to $\grid$, namely for $\grid_* \ge \grid$
\[
\enorm{\wh{u}-u_{\grid_*}} \le \enorm{\wh{u}-u_\grid}.
\]
However, the PDE estimator $\eta_\grid(u_\grid)$ fails to be monotone for fixed
discrete coefficients $(\wh{\bA},\wh{c})$ 
because it depends on the discrete solution $\uG\in\VG$ that changes with the mesh.
The following estimate, proved in Proposition \ref{P:est-reduction} (estimator reduction), 
quantifies the deviation
of $\eta_\grid(u_\grid)$ from monotonicity: there exists $\lambda>0$ such that for any $\delta>0$, $v \in \V_\grid$ and $v_* \in \V_{\grid_*}$,
\rhn{
\begin{align*}
\index{Constants!$\Clip$: estimator Lipschitz property constant}
      \eta_{\grid_*}(v_*,\grid_*)^2 &\le
      (1+\delta)\big(\eta_\grid(v,\grid)^2-\lambda\,\eta_\grid(v,\marked)^2\big)
      \\ & \quad +
      2\big(1+\delta^{-1}\big)\,\Clip \left( \enorm{v_*-v}^2 
      + \sum_{T\in\grid_*} \|P_\grid f - P_{\grid_*} f\|_{H^{-1}(\omega_T)}^2 \right),
\end{align*}
where $\Clip$ depends on $(\bA,c)$ and the shape regularity constant of $\grids$.
We refer to \cite{CaKrNoSi:08,MoSiVe:08} for the case 
$P_\grid f = P_{\grid_*} f = f \in L^2(\Omega).$
}
 
%%%%%%%%%%%%%%%%%%%%%%%%%%%%%%%%%%%%%%%%%%%%%%%%%%%%%%%%%%%%%%%%%%%%%%%%%%%%%%%%
\subsection{Convergence for discrete data: one-Step {\rm \AFEM}}\label{S:discrete-data}
%%%%%%%%%%%%%%%%%%%%%%%%%%%%%%%%%%%%%%%%%%%%%%%%%%%%%%%%%%%%%%%%%%%%%%%%%%%%%%%%

We now present the four basic modules of $\GALERKIN$, the one-step $\AFEM$ within
Algorithm \ref{algo:AFEM-Kernel} (\AFEMTS), namely
\begin{equation}\label{E:adaptive-loop}
\SOLVE \,\,
\longrightarrow \,\,
\ESTIMATE \,\,
\longrightarrow \,\,
\MARK \,\,
\longrightarrow \,\,
\REFINE,
\end{equation}
discuss their main properties, and prove a contraction property between
consecutive iterates of $\GALERKIN$. According to Algorithm \ref{algo:AFEM-Kernel},
given discrete data $\wh\data$ over a conforming mesh $\wh\grid$, created by $\DATA$,
and a desired tolerance $\eps>0$, the module
\begin{equation}\label{E:galerkin-eps}
[\mesh,u_\mesh]=\GALERKIN \, (\widehat{\mesh},\widehat{\data},\varepsilon)
\end{equation}
stops the loop \eqref{E:adaptive-loop} as soon as the error
tolerance $\eps$ is reached, i.e. as soon as 
\begin{equation}\label{E:eta-eps}
\eta_\mesh(u_\mesh) \leq \varepsilon.
\end{equation}

Since the data never change within $\GALERKIN$ and is always discrete, we assume in this section that $\data \in \D_{\grid}$ and do not use the hat symbol to indicate quantities defined using the (discrete) data.  

%-------------------------------------------------------------------------------
\subsubsection{Modules of {\rm \GALERKIN}}\label{S:GALERKIN}
%-------------------------------------------------------------------------------

\paragraph{Module {\rm\SOLVE}.}\index{Algorithms!\SOLVE: construct FEM approximation}
If $\grid\in\grids$ is a conforming refinement of $\grid_0$, and $\V_\grid$
is the finite element space of $C^0$ piecewise polynomials of degree $\le n$, then 
\[
 [\uG] = \SOLVE \, (\grid)
 \]
determines the Galerkin FEM solution \emph{exactly}, namely without algebraic error,
\begin{equation}\label{E:galerkin}
    \uG \in \VG:\quad \B[\uG,v] = \int_\Omega \nabla v \cdot \bA \nabla\uG + c vu
    = \langle f,v \rangle,
\end{equation}   
where $f\in H^{-1}(\Omega)$. However,
if $f\in\F_\grid$ is discrete as defined in \eqref{def-densities-of-P}, then 
\begin{equation*}
 \langle f,v \rangle =    
      \sum_{T\in \grid} \int_T q_T v + \sum_{F\in\faces} \int_F q_F v
      \quad\forall v\in\VG.
\end{equation*}
The assumption of exact solvability is made
for simplicity. The algebraic error committed in solving \eqref{E:galerkin} by iterative solvers
can be accommodated within the forthcoming theory. We refer to \cite{Stevenson:07,DanielVohralik:2023}
for details about how to relate the algebraic and PDE errors.
%

%-------------------------------------------------------------------------------------
\paragraph{Module {\rm \ESTIMATE}.}
\index{Algorithms!\ESTIMATE: compute the element error indicators and element data oscillations}
Given a conforming mesh $\grid\in\grids$ and the Galerkin solution
$\uG\in\VG$, the output of
\[
	[\{\eta_\grid(\uG,T),\osc_\grid(f,T)_{-1}\}_{T\in\grid}]=\ESTIMATE \, (\uG,\grid,\data)
\]
are the element error indicators $\eta_\grid(\uG,T)$ defined in \eqref{E:pde-estimator} with the discrete data $\data$, namely
\begin{equation*}
\eta_\grid (u_\grid,T)^2 = h_T^2 \| r(u_\grid) \|_T^2 + h_T \| j(u_\grid) \|_{\partial T}^2
\quad T\in\grid,
\end{equation*}
and element data oscillation $\osc_\grid(f,T)_{-1}$ defined in \eqref{E:data-osc}, namely
\begin{equation*}
\osc_\grid(f,T)_{-1} = \| f - P_\grid f \|_{H^{-1}(\omega_T)}.
\end{equation*}
We observe that for discrete forcing $f=P_\grid f$, global data oscillation vanishes
\begin{equation}\label{E:no-osc}
\osc_\grid (f)_{-1} = \| f - P_\grid f \|_{H^{-1}(\Omega)} = 0;
\end{equation}
this property is always valid within $\GALERKIN$. In this case, the output of $\ESTIMATE$ reduces
to just the PDE error indicators. Given $\cal S \subset\grid$, we denote
\[
\eta_\grid(v,{\cal S})^2:= \sum_{T\in {\cal S}}\eta_\grid(v,T)^2,
\quad
\eta_\grid(v)=\eta_\grid(v,\grid) \quad v\in \V_\grid.
\]
%

%---------------------------------------------------------------------------
\paragraph{Module {\rm \MARK}.}
\index{Algorithms!\MARK: D\"orfler marking}
%---------------------------------------------------------------------------
Given $\grid\in\grids$, the Galerkin solution $\uG\in\V_\grid$,
and element error indicators $\{\eta_\grid(\uG,\elm)\}_{\elm\in\grid}$,
the module \MARK selects elements for refinement using
\emph{D\"orfler Marking} (or bulk chasing) \cite{Dorfler:96,MoNoSi:00,NoSiVe:09,NochettoVeeser:2012},
\ie given a parameter $\theta\in(0,1]$ the output $\marked$ of
\begin{displaymath}
\index{Constants!$\theta$: D\"orfler marking parameter}
  [\marked]=\MARK\big(\{\eta_\grid(\uG,\elm)\}_{\elm\in\grid},\grid, \theta\big)
\end{displaymath}
satisfies
\begin{equation}\label{E:theta}
  \eta_\grid(\uG,\marked)\ge\theta\,\eta_\grid(\uG,\grid).
\end{equation}
This marking guarantees that $\marked$ contains a substantial part of
the total (or bulk) error, thus its name. The choice of $\marked$
does not have to be minimal at this stage, that is, the marked elements
$T\in\marked$ do not necessarily must be those with largest
indicators.

%---------------------------------------------------------------------------
\paragraph{Module {\rm \REFINE}.}
\index{Algorithms!\REFINE: refine all marked elements $b$ times and others necessary to produce a
conforming mesh}
%---------------------------------------------------------------------------
Let $b\in\N$ be the number of desired bisections per marked element.
Given $\grid\in\grids$ and a subset $\marked$ of marked elements, the
output $\gridk[*]\in\grids$ of
\begin{displaymath}
  [\grid_*]=\REFINE\big(\grid,\,\marked\big)
\end{displaymath}
is the smallest admissible refinement $\gridk[*]$ of $\grid$ so that
all elements of $\marked$ are at least bisected $b$ times. Therefore,
we have $h_{\grid_*}\leq h_\grid$ and the strict reduction property 
\begin{equation}\label{nv-h-reduction}
h_{\grid_*}|_T \leq 2^{-b/d} h_\grid|_T \qquad\fa T\in\marked,
\end{equation}
where $h_\grid:\Omega\to\R^+$ is a piecewise constant meshsize function that coincides with
$h_T=|T|^{1/d}$ on every $T\in\grid$.
We finally let
\begin{equation*}
\Res:=\refined:=\grid\backslash\grid_*
\end{equation*}
be the subset of refined elements of $\grid$ and note that
$
\marked\subseteq\Res.
$

Concatenating these four modules we get the standard SEMR one-step $\AFEM$.

%------------------------------------------------------------------------------------
\begin{algo}[$\GALERKIN$]\label{A:GALERKIN}
\index{Algorithms!\GALERKIN: procedure that iterates \SOLVE, \ESTIMATE, \MARK, \REFINE}
Let $\grid \ge \grid_0$ be a conforming refinement of a suitable initial mesh $\grid_0$.
Let data $\data=(\bA,c,f) \in D_{\grid}$ be discrete on $\grid$ and $\eps>0$ be
a stopping tolerance. The following one-step {\rm AFEM} creates a conforming refinement
$\grid_*\ge \grid$ and Galerkin solution $u_{\grid_*}\in\V_{\grid_*}$ for data $\data$
such that $\eta_{\grid_*}(u_{\grid_*})\le\eps$:
\begin{algotab}
  \> $[\grid_*,u_{\grid_*}] = \GALERKIN \, (\grid, \data, \eps)$ \\
  \> \> $\text{set } j=0, \grid_0 = \grid$\\
  \> \> $\text{do}$ \\
  \> \> \> $[u_{j}] = \SOLVE \, (\grid_j)$\\
  \> \> \> $[\{\eta_j(u_{j},T)\}_{T\in\grid_j}] = \ESTIMATE(u_{j},\grid_j,\data)$\\
  \> \> \>  $\text{if } \eta_j(u_{j}) \le \eps$ \\
  \> \> \> \> \text{return } $\grid_j,u_{j}$ \\
  \> \> \>  $[\marked_j] = \MARK \, \big(\{\eta_j(u_{j},T)\}_{T\in\grid_j},\grid_j,\theta\big)$\\
  \> \> \>  $[\grid_{j+1}] = \REFINE(\grid_j,\marked_j)$ \\
  \> \> \> $j\leftarrow j+1$ \\
  \> \> $\text{while true}$
\end{algotab}
\end{algo}

%-------------------------------------------------------------------------------
\subsubsection{Contraction property of {\rm \GALERKIN}}\label{S:contraction}
%-------------------------------------------------------------------------------

A key question to ask is what is (are) the quantity(ies) that $\GALERKIN$ may
contract. In light of \eqref{E:Pythagoras},
an obvious candidate is the energy error $\enorm{u-u_j}$\,, where $u_j\in \V_j = \V_{\grid_j}$ solves the problem
\begin{equation}\label{E:galerkin-j}
\B[u_j,w]= \langle f, w \rangle \qquad \forall w \in \V_j \,.    
\end{equation}
We now show that this is in fact the
case for discrete data $\data\in\D_{\grid}$ provided the {\it discrete local estimate} \eqref{E:local-lower-reduced} holds. The latter is a consequence of the {\it interior vertex property} of Definition \ref{D:interior-node} whenever $\bA$ is piecewise constant and $c=0$ in $\grid$ and data oscillation vanishes
$\osc_{\grid}(f)_{-1}=0$ \cite{MoNoSi:00,MoNoSi:02}.

\begin{lemma}[contraction property with discrete lower bound]\label{L:contraction-discrete-data}
Let data $\data\in\D_{\grid}$ be discrete and $u=u(\data)\in H^1_0(\Omega)$ be the corresponding exact solution.
If the subset $\marked_j \subset \grid_j$ of elements marked by {\rm \MARK} satisfies the discrete local estimate \eqref{E:local-lower-reduced} with respect to $\grid_{j+1} \ge \grid_j$, then
for $\alpha := \big(1 - \big(\theta\frac{C_2}{C_1}\big)^2\big)^{1/2} < 1$ the Galerkin solutions
$u_j\in\V_j, u_{j+1}\in\V_{j+1}$ of \eqref{E:galerkin-j} satisfy
\begin{equation}\label{E:simplest-contraction}
\index{Constants!$\alpha$: contraction constant}
\enorm{u-u_{j+1}} \le \alpha \enorm{u-u_j},
\end{equation}
where $0<\theta<1$ is the parameter in \eqref{E:theta}
and $C_1 \ge C_2$ are the constants in \eqref{E:aposteriori-bounds} and \eqref{E:local-lower-reduced}
respectively.
\end{lemma}
\begin{proof}
For convenience, we use the notation 
\begin{displaymath}
  e_j=\enorm{u-u_j},\ E_j=\enorm{u_{j+1}-u_j},\ \eta_j=\eta_j(u_j,\grid_j),\
  \eta_j(\marked_j)=\eta_j(u_j,\marked_j)
\end{displaymath}
and recall that $\estG[\grid_j](u_j,f) = \eta_j$ because $\osc_{\grid_j}(f)_{-1}=0$.
The key idea is to use the Pythagoras equality
\eqref{E:Pythagoras},
namely $e_{j+1}^2 = e_j^2 - E_j^2$,
and show that $E_j$ is a significant portion of $e_j$. 
Since \eqref{E:local-lower-reduced} implies
\[
C_2 \eta_j(\marked_j) \le E_j,
\]
applying D\"orfler marking \eqref{E:theta} and the upper
bound in \eqref{E:aposteriori-bounds},
we deduce
\begin{equation*}
E_j^2 \ge C_2^2 \theta^2 \eta_j^2 \ge \Big(\theta\frac{C_2}{C_1}\Big)^2 e_j^2.
\end{equation*}
This is the desired property of $E_j$ and leads to \eqref{E:simplest-contraction}.
\end{proof}

The contraction property \eqref{E:simplest-contraction} is very special and only valid for the
energy norm. For the $H^1_0$-norm we have the following simple but useful consequence.

\begin{corollary}[linear convergence]\label{C:linear-convergence}
If $c_\B \le C_\B$ are the constants in \eqref{E:exactB}, then
\begin{equation*}
|u-u_k|_{H^1_0(\Omega)} \le \sqrt{\frac{C_\B}{c_\B}} \, \alpha^{k-j} \, |u-u_j|_{H^1_0(\Omega)}
\quad k\ge j \ge 0.
\end{equation*}
\end{corollary}

We wonder whether or not the interior vertex property is necessary for \eqref{E:local-lower-reduced}, and thus
for \eqref{E:simplest-contraction}. We present an example, introduced
in \cite{MoNoSi:00,MoNoSi:02} to justify such a property for
constant data and $n=1$.

\begin{example}[lack of strict error monotonicity]\label{EX:interior-node}
\rm
Let $\Omega=(0,1)^2$, $\bA=\bI, c=0, f=1$ (constant data), 
and consider the sequences of
meshes depicted in Figure \ref{F:interior-node}.
If $\phi_0$ denotes the basis function associated with the only
interior vertex of the initial mesh $\grid_0$, then 
$u_0=u_1=\frac{1}{12}\,\phi_0$ and $u_2\neq u_1$.
\begin{figure}[h]
	\centering	\includegraphics[width=0.8\textwidth]{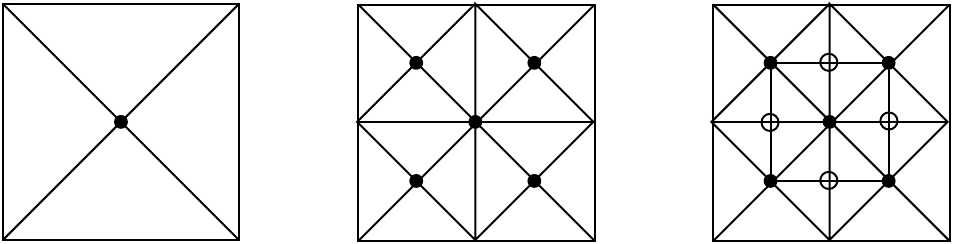}
  \caption{Grids $\gridk[0]$, $\gridk[1]$, and $\gridk[2]$ of
    Example \ref{EX:interior-node}. The mesh $\gridk[1]$ has nodes in the
    middle of edges of $\gridk[0]$, but only $\gridk[2]$ has nodes in
    the interior of elements of $\gridk[0]$. Hence, $\gridk[2]$
    satisfies the interior vertex property of Definition
    \ref{D:interior-node} with respect to $\gridk[0]$ whereas
    $\grid_1$ does not.
    \label{F:interior-node}}
\end{figure}

The mesh $\grid_1\ge\grid_0$ is produced by a standard $2$-step
bisection $(b=2)$ in $2d$. Since $u_0=u_1$ we conclude that the energy
error does not change 
$\enorm{u-u_0}=\enorm{u-u_1}$, whence \eqref{E:local-lower-reduced} fails,
between two consecutive steps of $\GALERKIN$ for $b=d=2$. This is no longer true provided
an interior vertex in each marked element is created, because then
Lemma \ref{L:contraction-discrete-data} (contraction property with discrete lower bound) holds.
\end{example}

\paragraph{Circumventing the discrete lower bound.}
Enforcing \eqref{E:local-lower-reduced} requires a minimal number $b_*$ of bisections, say $b_*=3,6$
for $d=2,3$, to guarantee the interior vertex property. This can be quite taxing, especially for $d=3$,
and relies on the strong assumption of $\bA$ being piecewise constant and $c=0$ on $\grid$.
It is clear from the preceding discussion that the energy error alone cannot
be expected to contract between consecutive iterates.
We explore next what quantity to monitor instead of the energy error in the
analysis, with the aim of avoiding \eqref{E:local-lower-reduced} and building a theory applicable
to general discrete coefficients $(\bA,c)$. This exploits the special structure of residual estimators
and does not directly extend to non-residual estimators.

%------------------------------------------------------------------------------------
\paragraph{Heuristics.} According to \eqref{E:Pythagoras}, the energy
error is monotone 
$\enorm{u-u_{j+1}}\le\enorm{u-u_j}$,	
but the previous Example shows that strict inequality may
fail. However, if $u_{j+1}=u_j$, estimate \eqref{E:strict-est} reveals
a strict estimator reduction 
$\eta_{j+1}(u_{j+1})<\eta_j(u_j)$. 
We thus expect that, for a suitable
scaling factor $\gamma>0$, the so-called \emph{quasi error}
\begin{equation}\label{E:quasi-error}
	\zeta_j^2(u_j) := \enorm{u-u_j}^2+\gamma\,\eta_j^2(u_j)
\end{equation}
may contract. This heuristics illustrates a distinct aspect of
{\rm AFEM} theory, the interplay between continuous quantities such the
energy error $\enorm{u-u_j}$ and discrete ones such as the estimator
$\eta_j(u_j)$: no one alone has the requisite properties to
yield a contraction between consecutive adaptive steps.
This result was originally proven in \cite{CaKrNoSi:08}.

\begin{theorem}[general  contraction property]\label{T:contraction}
  Let $\data \in\D_{\grid}$ be discrete data. Let $\theta\in (0,1]$ be the D\"orfler marking parameter,
  and $\{\grid_j,\V_j,u_j\}_{j=0}^\infty$
  be a sequence of conforming meshes, finite element spaces and
  discrete solutions $u_j\in\V_j$ created by {\rm \GALERKIN} for the model problem
  \eqref{E:galerkin-j}. 
  If $u=u(\data) \in H^1_0(\Omega)$ is the exact solution of \eqref{E:perturbed-weak-form},
  then there exist constants $\gamma>0$ and $0<\alpha<1$, additionally
  depending on the number $b\ge1$ of bisections and $\theta$, such that
  for all $j\ge 0$
  \begin{equation}\label{E:contraction-prop}
    \enorm{u-u_{j+1}}^2+\gamma\,\eta_{j+1}^2(u_{j+1})
    \le \alpha^2\,\Big(\enorm{u-u_j}^2+\gamma\,\eta_j^2(u_j)\Big).
  \end{equation}
\end{theorem}
\begin{proof}
We split the proof into four steps and use the notation in
Lemma \ref{L:contraction-discrete-data} (contraction property with discrete lower bound)

\medskip
\step{1} The error orthogonality \eqref{E:Pythagoras} 
reads
\begin{equation}\label{E:contract-1}
  e_{j+1}^2=e_j^2-E_j^2.
\end{equation}
Employing Proposition~\ref{P:est-reduction} (estimator reduction)
with $\grid=\grid_j$, $\grid_*=\grid_{j+1}$, 
\rhn{$v=u_j, v_*=u_{j+1}$ and $f=f_*\in\F_j$} gives
\begin{equation}\label{E:contract-2}
  \eta_{j+1}^2\le(1+\delta)\big(\eta_j^2-\lambda
  \,\eta_j^2(\marked_j)\big)+(1+\delta^{-1})\, \Clip^2 \,E_j^2.
\end{equation}
After multiplying
\eqref{E:contract-2} by $\gamma>0$, to be determined later, we
add \eqref{E:contract-1} and \eqref{E:contract-2} to obtain
\begin{align*}
  e_{j+1}^2+\gamma\,\eta_{j+1}^2&\le e_j^2
  +\big(\gamma\,(1+\delta^{-1})\,\Clip^2 -1\big)\,E_j^2
  +\gamma\,(1+\delta)\,\big(\eta_j^2-\lambda\,\eta_j^2(\marked_j)\big).
\end{align*}
\step{2}
We now choose the parameters $\delta, \gamma$: let $\delta$ satisfy
\[
(1+\delta)\big(1-\lambda\theta^2\big) 
= 1 -\frac{\lambda\theta^2}{2},
\]
and $\gamma$ verify
\[
\gamma\,(1+\delta^{-1})\,\Clip^2 =1.
\]
Note that this choice of $\gamma$ yields
\begin{equation}\label{e:monoton_galerkin_total_error-A}
  e_{j+1}^2+\gamma\,\eta_{j+1}^2\le e_j^2
  +\gamma\,(1+\delta)\,\big(\eta_j^2-\lambda\,\eta_j^2(\marked_j)\big).
\end{equation}
\step{3}
We next employ D\"orfler Marking \eqref{E:theta}, namely 
$\eta_j(\marked_j)\ge\theta\eta_j$, to deduce
\begin{equation*}
  e_{j+1}^2+\gamma\,\eta_{j+1}^2
  \le
  e_j^2
  + \gamma (1+\delta)(1-\lambda\theta^2)\eta_j^2.
\end{equation*}
This, in conjunction with the choice of $\delta$, gives
\begin{equation}\label{e:monoton_galerkin_total_error}
  e_{j+1}^2+\gamma\,\eta_{j+1}^2 
  \le 
  e_j^2
  + \gamma \bigg(1-\frac{\lambda\theta^2}{2}\bigg)\eta_j^2
\end{equation}
which we write 
\begin{equation*}
  e_{j+1}^2+\gamma\,\eta_{j+1}^2 
  \le e_j^2 - \frac{\gamma\lambda\theta^2}{4}\eta_j^2
  + \gamma \bigg(1-\frac{\lambda\theta^2}{4}\bigg)\eta_j^2.
\end{equation*}

\step{4} Finally, the upper bound in \eqref{E:aposteriori-bounds}, namely
\begin{math}
  e_j\le C_1\,\eta_j,
\end{math}
implies that
\begin{displaymath}
  e_{j+1}^2+\gamma\,\eta_{j+1}^2 \le 
  \bigg(1-\frac{\gamma\lambda\theta^2}{4C_1^2}\bigg) e_j^2
  + \gamma \bigg(1-\frac{\lambda\theta^2}{4}\bigg) \eta_j^2.
\end{displaymath}
This in turn leads to
\begin{displaymath}
  e_{j+1}^2+\gamma\,\eta_{j+1}^2 
  \le  \alpha^2
  \big(e_j^2 +\gamma\,\eta_j^2\big),
\end{displaymath}
with
$\alpha^2 \definedas \max \big\{1-\frac{\gamma\lambda\theta^2}{4C_1^2},
  1-\frac{\lambda\theta^2}{4} \big\}<1$,
and thus concludes the proof of theorem. 
\end{proof}

\begin{remark}[basic ingredients]
  This proof solely uses D\"orfler
  marking \eqref{E:theta}; 
  Pythagoras identity \eqref{E:Pythagoras}; the a posteriori
  upper bound in \eqref{E:aposteriori-bounds}; and Proposition~\ref{P:est-reduction} (estimator
  reduction). The proof circumvents the use altogether of the lower bound in \eqref{E:aposteriori-bounds}
  and the discrete lower bound \eqref{E:local-lower-reduced}.
\end{remark}

The contraction property \eqref{E:contraction-prop} is valid for a suitable combination of
the energy norm $\enorm{u-u_j}$ and the PDE estimator $\eta_j(u_j)$. We cannot expect this type of
result for the underlying space norm $|u-u_j|_{H^1_0(\Omega)}$. We have instead the following
statement, whose structure reflects the possible stagnation of $|u-u_j|_{H^1_0(\Omega)}$ during the
refinement process, as documented in Example \ref{EX:interior-node}.

\begin{corollary}[linear convergence of error]\label{C:linear-convergence-error}
If the assumptions of Theorem \ref{T:contraction} are valid, and $0<\alpha<1, \gamma>0$ are the
constants in \eqref{E:contraction-prop}, then there holds
\begin{equation}\label{E:linear-convergence-error}
|u-u_k|_{H^1_0(\Omega)} \le C_* \alpha^{k-j} |u-u_j|_{H^1_0(\Omega)} \quad\forall \, k \ge j \ge 0,
\end{equation}
with $C_* = \big(\frac{C_\B}{c_\B} \big( 1 + \frac{\gamma}{C_2^2} \big)\big)^{1/2}>1$ and
constants $C_\B\ge c_\B>0$ and $C_2>0$ given in \eqref{E:exactB} and \eqref{E:aposteriori-bounds}
respectively. 
\end{corollary}
\begin{proof}
Simply concatenate \eqref{E:exactB}, \eqref{E:contraction-prop} and \eqref{E:aposteriori-bounds}
to obtain
\begin{align*}
c_\B |u-u_k|_{H^1_0(\Omega)}^2 &\le \enorm{u-u_k}^2 + \gamma \, \eta_k(u_k)^2
\\ & \le \alpha^{2(k-j)} \Big( \enorm{u-u_j}^2 + \gamma \, \eta_j(u_j)^2 \Big)
\\ &
\le \alpha^{2(k-j)} \left(C_\B \Big(1+\frac{\gamma}{C_2^2} \Big) \right) |u-u_j|_{H^1_0(\Omega)}^2.
\end{align*}
This implies \eqref{E:linear-convergence-error} and concludes the proof.
\end{proof}

We stress that, in contrast to \eqref{E:contraction-prop}, \eqref{E:linear-convergence-error}
does rely on the lower bound in \eqref{E:aposteriori-bounds}. This is  not the case if we
express linear convergence in terms of the PDE estimator. The proof is similar to the
preceding one and is omitted.

\begin{corollary}[linear convergence of estimator]\label{C:linear-convergence-est}
If the assumptions of Theorem \ref{T:contraction} are valid, and $0<\alpha<1, \gamma>0$ are the
constants in \eqref{E:contraction-prop}, then there holds \looseness=-1
\begin{equation}\label{E:linear-convergence-est}
\eta_k(u_k) \le C_\# \alpha^{k-j} \eta_j(u_j) \quad\forall \, k \ge j \ge 0,
\end{equation}
with $C_\# = \big(1+\frac{C_1^2}{\gamma}\big)^{1/2}>1$ and $C_1$ given in \eqref{E:aposteriori-bounds}.
\end{corollary}

\begin{remark}[stopping]
In view of \eqref{E:linear-convergence-est}, \eqref{E:aposteriori-bounds-H1}, we realize that \GALERKIN requires $j\leq J$ iterations until the stopping criteria $\eta_j \leq \eps$ is satisfied and delivers the error $|u-u_j|_{H^1_0(\Omega)} \le C_U\eps$,
where
\[
J \leq 1 + \frac{\log{\frac{\eps}{C_\#\eta_0}}}{\log\alpha}.
\]
\end{remark}

%-----------------------------------------------------------------------------------
\subsubsection{Discontinuous coefficients: Kellogg's example}
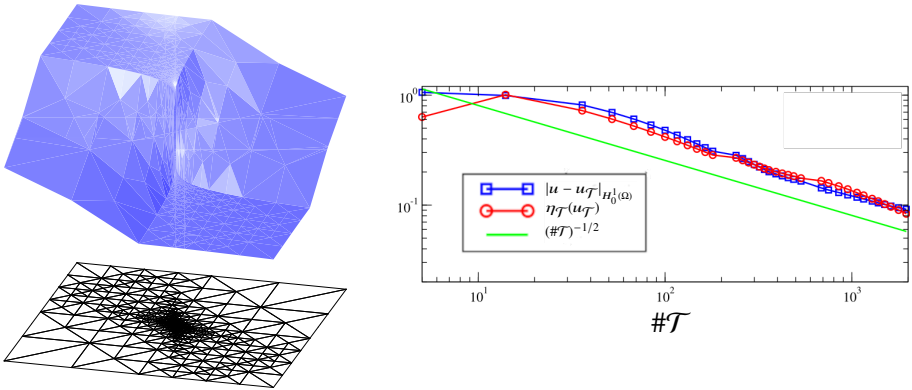
\begin{figure}[h!]
\begin{minipage}{0.4\linewidth}
\centering
%\begin{figure}%[ht]
  \includegraphics[width=0.9\textwidth]{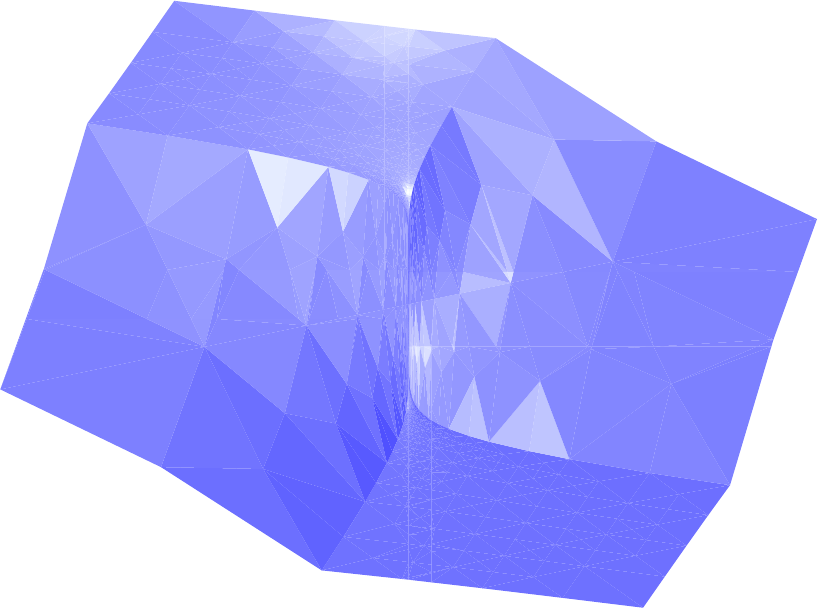}

  \includegraphics[width=0.9\textwidth]{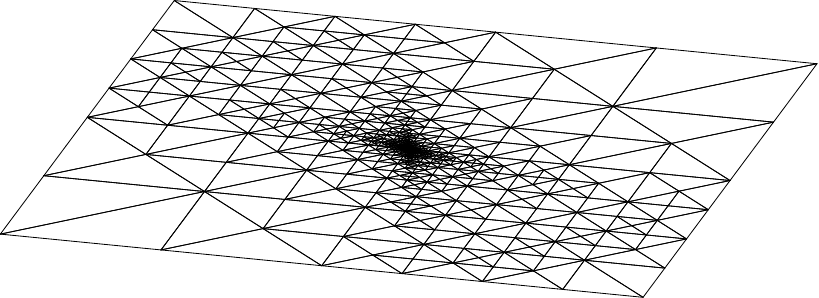}
%\end{figure}
\end{minipage}
\begin{minipage}{0.6\linewidth}
\centering
%\begin{figure}
\includegraphics[width=0.9\textwidth]{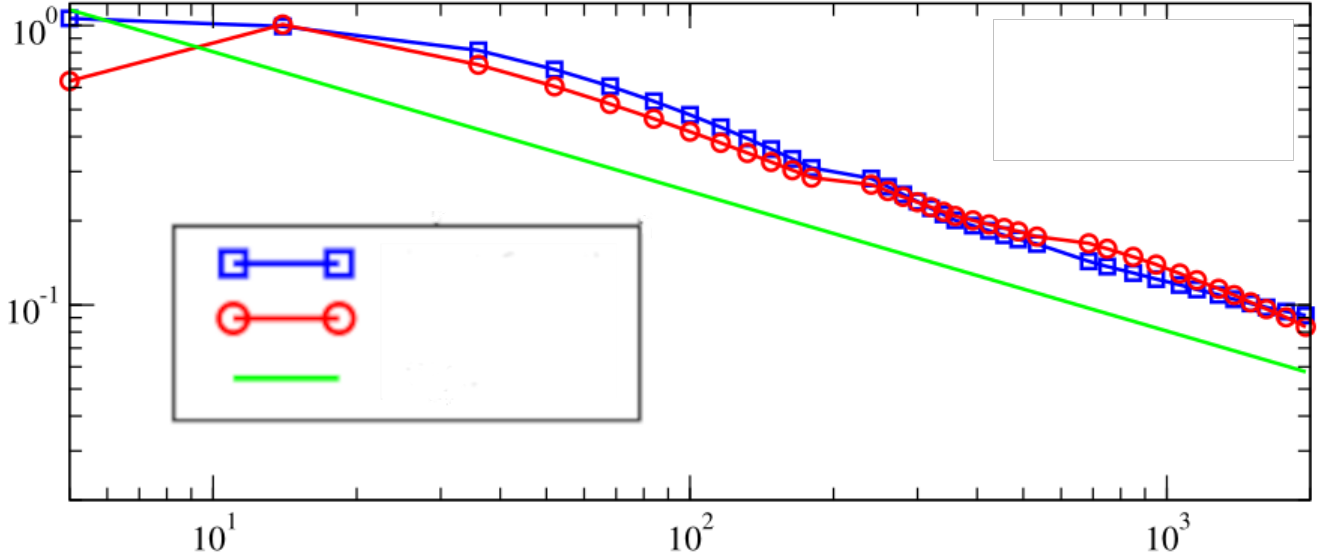}
\begin{picture}(0,0)
\put(-140,42){\tiny $|u-u_\grid|_{\scaleto{H^1_0(\Omega)}{6pt}}$}
\put(-140,35){\tiny $\eta_\grid(u_\grid)$}
\put(-140,26){\tiny $(\#\grid)^{\scaleto{-1/2}{4pt}}$}
\put(-100,-10){$\#\grid$}
\end{picture}
%\end{figure}
\end{minipage}
\caption{Discontinuous coefficients in checkerboard pattern: (left) graph
   of the discrete solution $u$, which is $u\approx r^{0.1}$, and
   underlying strongly graded grid $\grid$ towards the origin 
   (notice the steep gradient of $u$ at the origin);
   (right) estimate and true error in terms of $\#\grid$
   (the optimal decay for piecewise linear
   elements in $2d$ is indicated by the green line with slope $-1/2$).
\label{F:kellogg-1}}
\end{figure}

%-----------------------------------------------------------------------------------
%
We examine a simple yet quite demanding example with piecewise constant coefficients in
checkerboard pattern for $d=2$ due to \cite{Kellogg:75}, and used in
\cite{MoNoSi:00,MoNoSi:02} as a 
benchmark for $\GALERKIN$. We consider 
$\Omega=(-1,1)^2$, $\bA = a_1 {\vec{I}}$ in the first and third
quadrants, and $\bA = a_2 {\vec{I}}$ in the second and fourth
quadrants. This checkerboard pattern is the worst for the
regularity of the solution $u$ at the origin. For $f=c=0$, a function
of the form $u(r,\theta) =  r^\gamma\mu(\theta)$ in polar coordinates
solves \eqref{strong-form} with nonvanishing Dirichlet condition
for suitable $0<\gamma<2$ and $\mu$
\cite{MoNoSi:00,MoNoSi:02,NoSiVe:09}. We choose $\gamma=0.1$, which 
leads to $u\in H^s(\Omega)$ for $1\le s<1.1$ and piecewise in $W^2_p$ for
some $p>1$. This corresponds to diffusion coefficients
$a_1 \approx 161.44$ and $a_2 = 1$,
which can be computed via Newton's method; the closer $\gamma$ is to
$0$, the larger is the ratio $a_1/a_2$. The solution $u$ and a
sample mesh are depicted in Figure \ref{F:kellogg-1}(left).

\begin{figure}[ht]
  \centering
  \begin{tabular}{cccc}
  \includegraphics[width=0.2\textwidth]{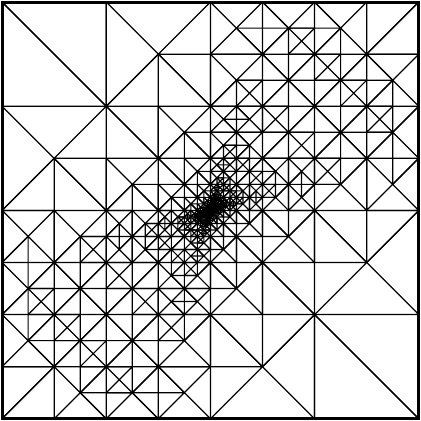}%
  &\includegraphics[width=0.2\textwidth]{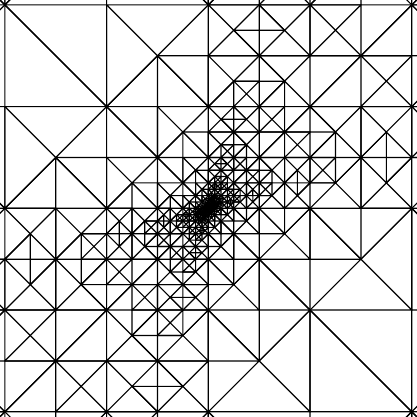}%
  &\includegraphics[width=0.2\textwidth]{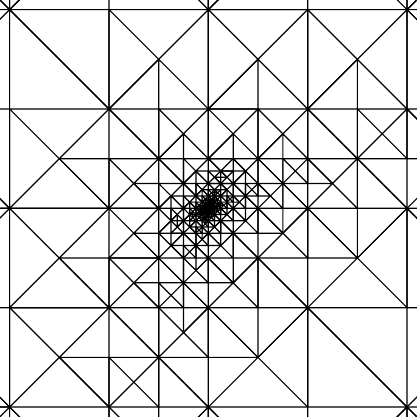}%
  &\includegraphics[width=0.2\textwidth]{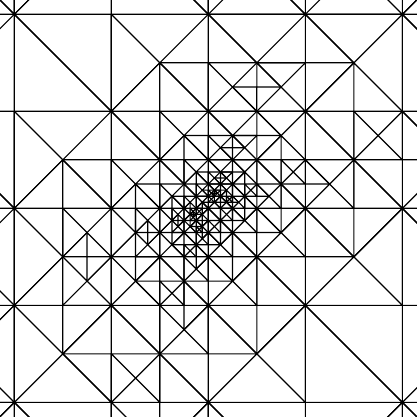}\\
  (a) & (b) & (c) & (d)
  \end{tabular}
  \caption{Discontinuous coefficients in checkerboard pattern: (a) final
    grid $\grid$ highly graded towards the origin
    with cardinality $\#\grid \approx 2000$; (b) zoom to
    $(-10^{-3},10^{-3})^2$; (c) zoom to $(-10^{-6},10^{-6})^2$;
    (d) zoom to $(-10^{-9},10^{-9})^2$. For a similar resolution, a
    uniform grid $\grid$ would require cardinality $\#\grid\approx 10^{20}$.
    \label{F:kellogg-2}}
\end{figure}

Figure \ref{F:kellogg-1} (right) documents the optimal performance of
$\GALERKIN$: both the energy error and estimator exhibit optimal decay
$(\#\grid)^{-1/2}$ in terms of the cardinality $\#\grid$ of the
underlying mesh $\grid$ for piecewise linear finite elements.
On the other hand, Figure \ref{F:kellogg-2} displays a strongly
graded mesh $\grid$ towards the origin generated by $\GALERKIN$ using bisection, 
and three zooms which reveal a selfsimilar structure. It is worth
stressing that the meshsize is of order $10^{-10}$ at the origin
and that $\#\grid \approx 2\times 10^3$, whereas to reach a similar resolution with
a uniform mesh $\grid$ we would need $\#\grid\approx 10^{20}$.
This example clearly reveals that adaptivity can restore optimal
performance even with modest computational resources. 

Classical FEMs with quasi-uniform meshes $\grid$ require regularity $u\in
H^2(\Omega)$ to deliver an optimal convergence rate
$(\#\grid)^{-1/2}$ with polynomial degree $n=1$. Since $u\notin H^s(\Omega)$ for any $s>1.1$, this
is not possible for the example above. However, the problem is not
quite the lack of second derivatives, but rather the fact that they
are not square integrable. In fact, the function $u$ is in $W^2_p$ for
$p>1$ in each quadrant, and so over the initial mesh $\grid_0$, namely $u\in
W^2_p(\Omega;\grid_0)$. The computational rate of convergence $(\#\grid)^{-1/2}$ is
consistent with Corollary \ref{C:optimal-mesh}. We will prove that $\GALERKIN$ delivers this
rate in Section \ref{S:conv-rates-coercive}.
  
%------------------------------------------------------------------------------------
\subsection{Data oscillation: one-Step {\rm AFEM} with switch}\label{S:one-step-switch}
%------------------------------------------------------------------------------------

In Section \ref{S:discrete-data} we assumed that the full data $\data=(\bA,c,f)\in\D_{\grid}$ is discrete,
and in particular $f = P_{\grid} f \in \F_{\grid}$.
The finite dimensional nature of $\wh \data$ allowed us to develop a rather simple theory of
convergence for $\GALERKIN$, the one-step AFEM, that hinges exclusively on the PDE local error indicator
$\eta_\grid (u_\grid,T)$ defined in \eqref{E:pde-estimator}.
We now keep $(\bA,c)$ discrete, whence the elliptic operator
in \eqref{strong-form} includes the Laplacian, but explore the role of
a general forcing $f \ne P_\grid f$. Therefore, in contrast with \eqref{E:no-osc},
we now investigate the effect of data oscillation \eqref{E:global-est-osc}
\begin{equation*} %\label{E:data-osc}
\osc_\grid (f)_{-1}^2 = \sum_{T\in\grid} \| f - P_\grid f \|_{H^{-1}(\omega_T)}^2
\end{equation*}
for any $\grid\in\grids$, and present a linear convergence theory. We recall from Theorem
\ref{T:modified-estimator} (modified residual estimator) that the total error estimator
$\est_\grid(u_\grid,f)^2 = \eta_\grid (u_\grid)^2 + \osc_\grid (f)_{-1}^2$ is equivalent to the
$H^1$-error, namely
\begin{equation}\label{E:error-equiv-estimator}
C_L \est_\grid(u_\grid,f) \le \| \nabla(u-u_\grid) \|_{L^2(\Omega)} \le C_U \est_\grid(u_\grid,f).
\end{equation}

As in the previous section, to simplify notation
we do not use the hat symbol to indicate quantities defined with the discrete data $(\bA,c)$.

%------------------------------------------------------------------------------
\subsubsection{Role of data oscillation}\label{S:role-osc}
%------------------------------------------------------------------------------
%
At a first sight, it might seem that Example~\ref{EX:interior-node} (lack of strict error monotonicity)
is very special and can only occur at the beginning of the refinement process. We now show that this
situation can happen at any stage and that even an interior vertex property may not guarantee error or data
oscillation decrease.

\begin{example}[interior vertex] \label{EX:chekerboard-f}
\rm
Let the polynomial degree be $n=1$,
fix $m \in \N$ and consider \eqref{E:galerkin} with $\bA={\vec{I}}$ the identity matrix, $c=0$, $\Omega=(0,1)^2$
and checkerboard $f$ given by the following expression and depicted in Figure~\ref{F:example-checkerboard} 
 (left)
\[
f(x) = \begin{cases}
      \phantom{-}1 &\text{if }x \in (i\,2^{-m},(i+1)\,2^{-m}) \times
                         (j\,2^{-m},(j+1)\,2^{-m}) \text{ and }
                   i+j\text{ odd} \\
      -1 &\text{otherwise}.
	 \end{cases}
\]
We start with the same mesh $\grid_0$ with four elements as in Example~\ref{EX:interior-node},
and construct recursively grids $\grid_{k+1}\in\grids$, $k\ge0$, as a conforming refinement of
$\grid_k\in\grids$ via two newest-vertex bisections of every triangle of $\grid_k$; see Figure~\ref{F:example-checkerboard} (right). Since $f$
is $L^2$-orthogonal to every piecewise linear basis function of the space
$\V_{\grid_k}=\mathbb{S}_{\grid_k}^{1,0}$ for $0 \le k \le m-1$, we deduce that $u_{\grid_k}=0$
and the energy error does not change
\begin{equation}\label{E:error=error0}
\enorm{u-u_{\grid_k}} = \enorm{u-u_{\grid_0}} \quad 0 \le k \le m-1.
\end{equation}
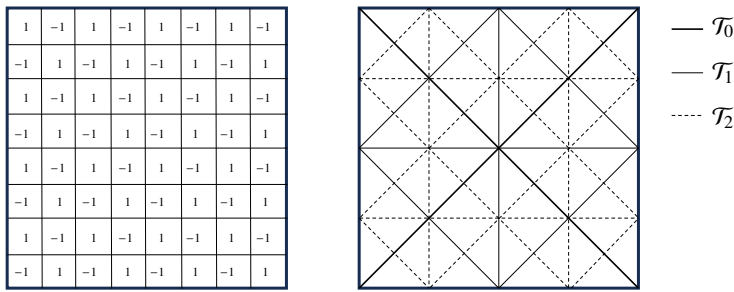
\begin{figure}[htbp]
\begin{center}
\scalebox{0.85}{
\includegraphics[trim={3.5cm 3cm 2.5cm 4cm},clip,width=0.9\textwidth]{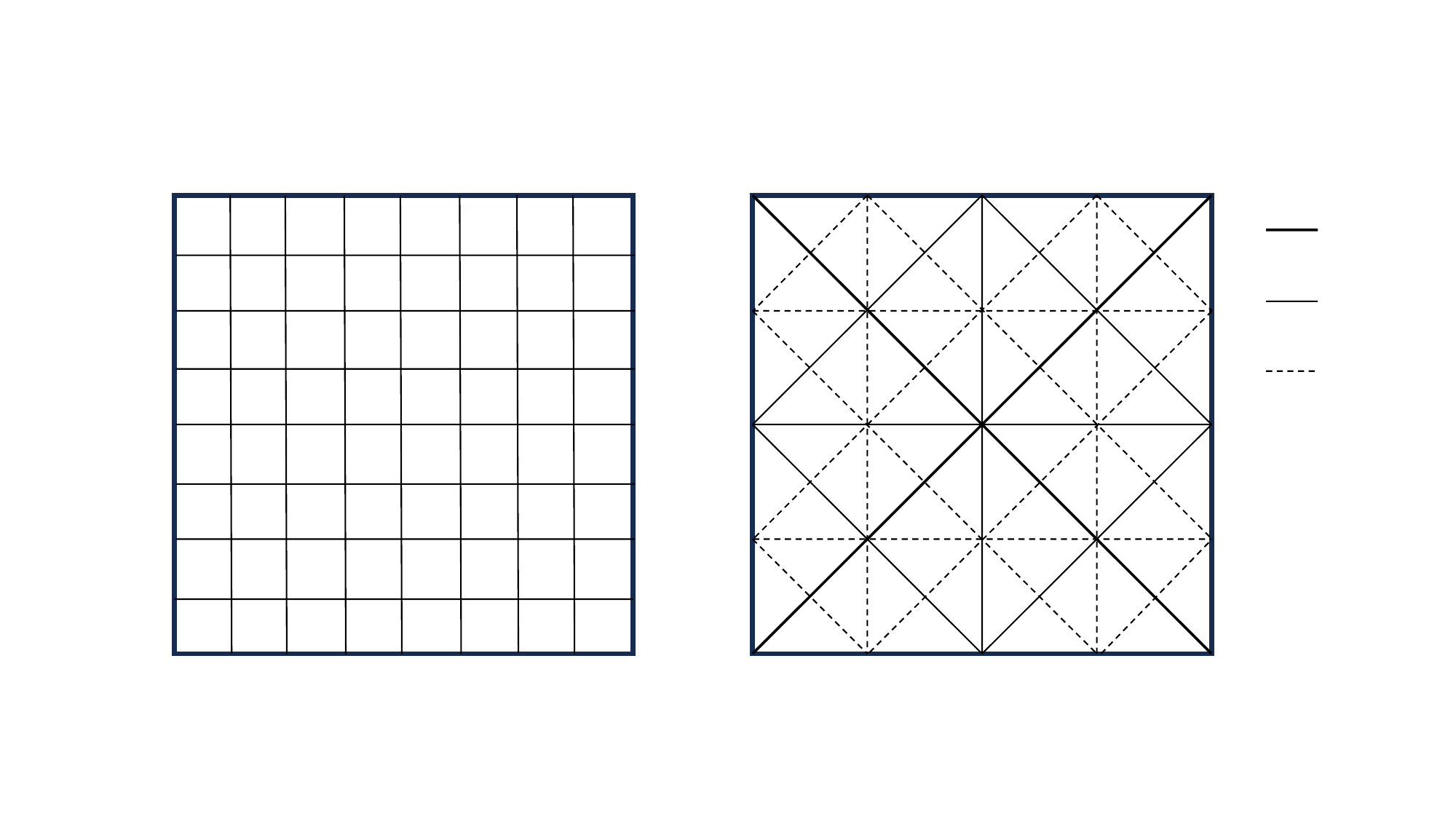}
\begin{picture}(0,0)
\put(-7,122){$\grid_0$}
\put(-7,102){$\grid_1$}
\put(-7,82){$\grid_2$}
% first column
\put(-312,123){\tiny $1$}
\put(-316,107){\tiny $-1$}
\put(-312,92){\tiny $1$}
\put(-316,76){\tiny $-1$}
\put(-312,61){\tiny $1$}
\put(-316,46){\tiny $-1$}
\put(-312,30){\tiny $1$}
\put(-316,15){\tiny $-1$}
% second column
\put(-300,123){\tiny $-1$}
\put(-297,107){\tiny $1$}
\put(-300,92){\tiny $-1$}
\put(-297,76){\tiny $1$}
\put(-300,61){\tiny $-1$}
\put(-297,46){\tiny $1$}
\put(-300,30){\tiny $-1$}
\put(-297,15){\tiny $1$}
% third column
\put(-282,123){\tiny $1$}
\put(-286,107){\tiny $-1$}
\put(-282,92){\tiny $1$}
\put(-286,76){\tiny $-1$}
\put(-282,61){\tiny $1$}
\put(-286,46){\tiny $-1$}
\put(-282,30){\tiny $1$}
\put(-286,15){\tiny $-1$}
% 4th column
\put(-270,123){\tiny $-1$}
\put(-267,107){\tiny $1$}
\put(-270,92){\tiny $-1$}
\put(-267,76){\tiny $1$}
\put(-270,61){\tiny $-1$}
\put(-267,46){\tiny $1$}
\put(-270,30){\tiny $-1$}
\put(-267,15){\tiny $1$}
% 5th column
\put(-252,123){\tiny $1$}
\put(-256,107){\tiny $-1$}
\put(-252,92){\tiny $1$}
\put(-256,76){\tiny $-1$}
\put(-252,61){\tiny $1$}
\put(-256,46){\tiny $-1$}
\put(-252,30){\tiny $1$}
\put(-256,15){\tiny $-1$}
% 6th column
\put(-240,123){\tiny $-1$}
\put(-237,107){\tiny $1$}
\put(-240,92){\tiny $-1$}
\put(-237,76){\tiny $1$}
\put(-240,61){\tiny $-1$}
\put(-237,46){\tiny $1$}
\put(-240,30){\tiny $-1$}
\put(-237,15){\tiny $1$}
% 7th column
\put(-221,123){\tiny $1$}
\put(-225,107){\tiny $-1$}
\put(-221,92){\tiny $1$}
\put(-225,76){\tiny $-1$}
\put(-221,61){\tiny $1$}
\put(-225,46){\tiny $-1$}
\put(-221,30){\tiny $1$}
\put(-225,15){\tiny $-1$}
% 8th column
\put(-209,123){\tiny $-1$}
\put(-206,107){\tiny $1$}
\put(-209,92){\tiny $-1$}
\put(-206,76){\tiny $1$}
\put(-209,61){\tiny $-1$}
\put(-206,46){\tiny $1$}
\put(-209,30){\tiny $-1$}
\put(-206,15){\tiny $1$}
\end{picture}
}
\end{center}
\caption{Representation of the checkerboard function $f$ of Example~\ref{EX:chekerboard-f} for
         $m=3$~(left), and grids $\grid_k$ for $k=0,1,2$~(right).
         \label{F:example-checkerboard}}
\end{figure}

We see that this procedure creates three interior vertices in every triangle of $\grid_k$ after two
refinement steps, namely in $\grid_{k+2}$ as long as $k+2 \le m$. Since the error does not change,
we conclude that the interior vertex
property is necessary for error reduction but is not sufficient in the presence of data oscillation
$\osc_{\grid_k}(f)_{-1} \ne 0$. We conclude
\begin{equation}\label{E:oscillation-high-order}
\begin{minipage}{0.85\linewidth}
\emph{Data oscillation $\osc_\grid(f)_{-1}$ is not generally of higher order than the error,
especially in the early stages of the adaptive process}. 
\end{minipage}
\end{equation}
On the other hand, for $k=m$ the discrete solution $u_{\grid_m}$ does no longer vanish globally, but is still
zero along the lines where $f$ changes sign due to the symmetry of the
problem, and the same happens with $u_{\grid_{m+1}}$. Therefore, the behavior of $u_{\grid m}$ and
$u_{\grid_{m+1}}$ in a fixed square, where $f$ is constant, is exactly the same as in
Example~\ref{EX:interior-node}. This implies that $u_{\grid_m} = u_{\grid_{m+1}}$, and illustrates that the
rather special situation of Example~\ref{EX:interior-node} can occur at any stage of the
refinement process. 
\end{example}

\begin{example}[vanishing of $P_\grid f$ for $n=1$] \label{EX:projection-Pgrid}
\rm
Since $P_\grid f$ is constructed locally upon testing $f$ against cubic and quadratic bubbles
(see Remark~\ref{R:Computation-of-P} (local computation)), and $f$ of Example~\ref{EX:chekerboard-f}
is highly oscillatory, we realize
that $P_{\grid_k} f$ is rather small relative to $f$ in $H^{-1}(\Omega)$, but it is not zero. This is
due to the lack of complete symmetry of the checkerboard pattern and the triangular grid. Suppose that
each square of Fig.~\ref{F:example-checkerboard}, where $f=\pm 1$, is further split across the diagonals
into four triangles, and that $f$ is assigned the alternating values $\pm 1$ and $\mp 1$ in each triangle
depending on whether $f$ was originally $1$ or $-1$ in that square;
this configuration is displayed in Fig.~\ref{F:symmetric-chekerboard}. Suppose further that the
coefficients $(\bA,c)$ of the operator \eqref{strong-form} are piecewise constant, as it happens for
the Laplacian, the polynomial degree is $n=1$, and the definition \eqref{e:Computation-of-PT_local}
of $P_T$ over a triangle $T\in\grid$ uses $q \in \P_0$ rather than $\P_1$. In light of
\eqref{e:Computation-of-PT_local} and \eqref{e:Computation-of-PF_local}, symmetry yields
for all $T\in\grid$ and $F\in\faces$
\begin{equation}\label{E:Pgrid=0}
\int_T f \phi_T = 0 \quad\Rightarrow\quad P_T f = 0,
\qquad
\int_F f \phi_F = 0 \quad\Rightarrow\quad P_F f = 0,
\end{equation}
whence $P_\grid f = 0$. Since also $u_\grid = 0$ because $f$ is orthogonal to all basis functions of
$\V_\grid$, we deduce $\est_\grid(u_\grid,f) = 0$ and all the information about the error $\enorm{u-u_\grid}\ne0$
resides in the data oscillation $\osc_\grid(f)_{-1}\ne0$. Moreover, the fact that $P_\grid f = 0$ for several
iterations reveals the important property that $\osc_\grid(f)_{-1}$ may not change upon refinement because
\begin{equation}\label{E:osc=f}
\osc_\grid(f)_{-1}^2 = \sum_{T\in\grid} \| f \|_{H^{-1}(\omega_T)}^2.
\end{equation}
Since $\enorm{u-u_\grid} \approx \osc_\grid(f)_{-1}$, according to
\eqref{T:modified-estimator}, a special care
must be exercised to reduce data oscillation when it dominates. This justifies the structure of
the algorithm one-step AFEM with switch below.
\begin{figure}[hpbt]
\begin{center}
\includegraphics[trim={3cm 3.5cm 7cm 4cm},clip,width=.5\textwidth]{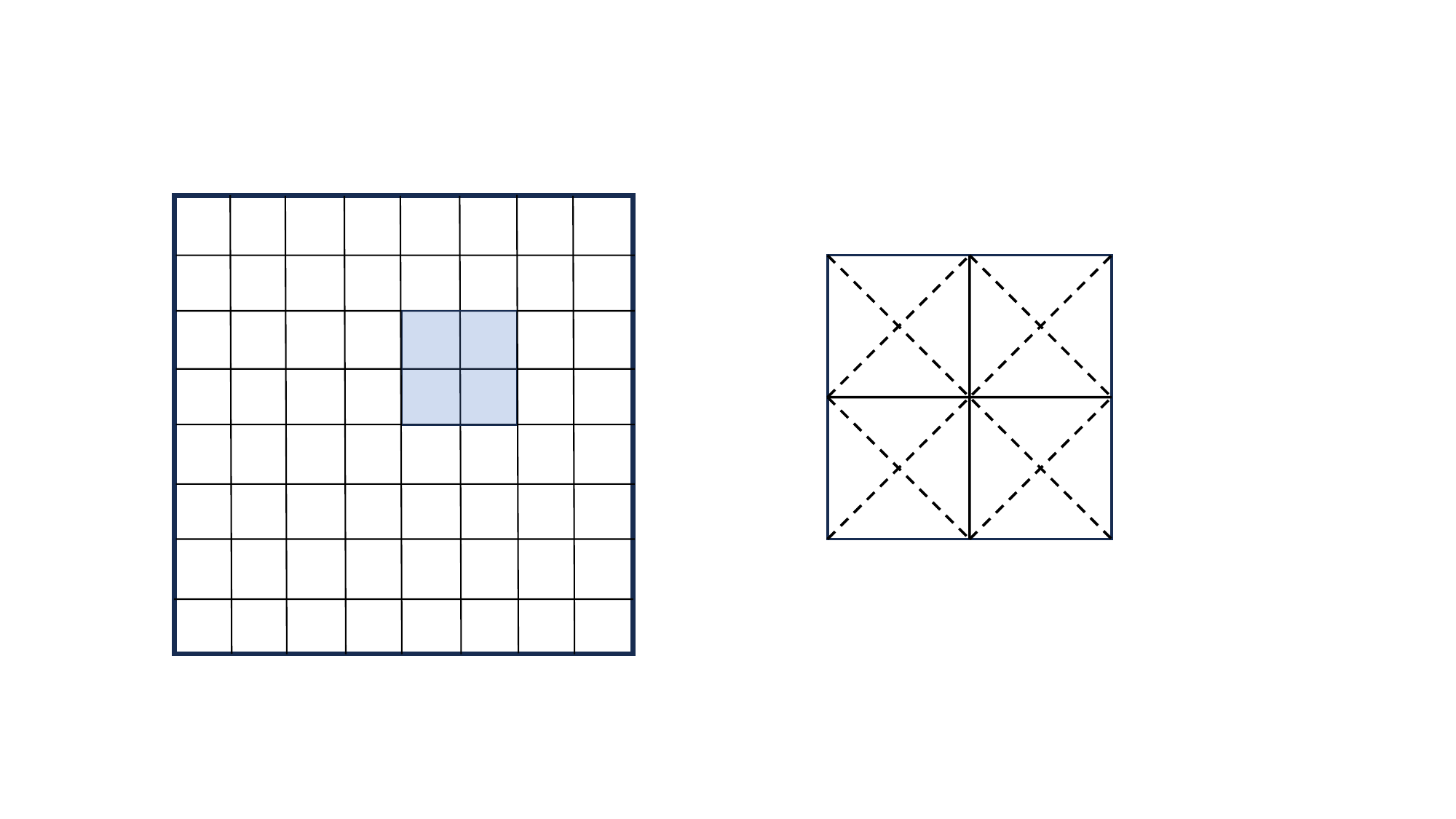}
\begin{picture}(0,0)
\put(-131,56){\tiny $1$}
\put(-122,56){\tiny $-1$}
\put(-133,46){\tiny $-1$}
\put(-120,46){\tiny $1$}
\put(-56,58){\tiny $1$}
\put(-42,58){\tiny $1$}
\put(-33,58){\tiny $-1$}
\put(-19,58){\tiny $-1$}
\put(-58,33){\tiny $-1$}
\put(-43,33){\tiny $-1$}
\put(-31,33){\tiny $1$}
\put(-17,33){\tiny $1$}
\put(-51,66){\tiny $-1$}
\put(-51,51){\tiny $-1$}
\put(-49,41){\tiny $1$}
\put(-49,26){\tiny $1$}
\put(-25,66){\tiny $1$}
\put(-25,51){\tiny $1$}
\put(-28,41){\tiny $-1$}
\put(-28,26){\tiny $-1$}
\put(-85,45){$\Longrightarrow$}
\end{picture}
\caption{
    Refinement of the shaded according to the process described in Example~\ref{EX:projection-Pgrid}.
    \label{F:symmetric-chekerboard}}
\end{center}
\end{figure}
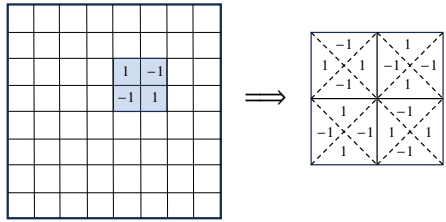
\end{example}

\begin{example}[vanishing of $P_\grid f$ for $n>1$] \label{EX:Pgrid-n>1}
Given $n\geq 1$ a polynomial degree and $\grid_k$, $k=1,...,m$, uniform refinements of $\grid_0$, there are finitely many conditions to verify for $f \in H^{-1}(\Omega)$
to be orthogonal to $\V_{\grid_k}$ and to $\F_{\grid_k}$. Since $\dim H^{-1}(\Omega) = \dim L^2(\Omega) = \infty$ there are infinitely many loads $f\in H^{-1}(\Omega)$ as well as
in $L^2(\Omega)$ that yield $u_{\grid_k} = P_{\grid_k} f = 0$, which implies \eqref{E:error=error0}. Moreover, $\eta_{\grid_k}(u_{k})=0$ and $\est_{\grid_k}(u_{\grid_k})=\osc_{\grid_k}(f)_{-1}$ satisfies  \eqref{E:osc=f}. 
One explicit example is as follows.

Given an initial mesh $\grid_0$, suppose that $f$ consists of line Dirac masses supported on the
skeleton of $\grid_0$ with densities $g_F$ on $F\in\faces_{\grid_0}$ made of piecewise polynomials of
degree $2n+1$. We further assume that the $g_F$ are orthogonal to $\P_{2n}$ over $F$ as well as over all sub-faces obtained from $m \ge 1$ uniform refinements of $\grid_0$;
see Fig.~\ref{F:example-checkerboard}. In such a situation,
\eqref{E:Pgrid=0} applies and $u_{\grid_k}=P_{\grid_k} f = 0$, whence \eqref{E:error=error0} and
\eqref{E:osc=f} are valid for $0 \le k \le m$.
\end{example}

These three examples reveal the following crucial and novel feature about the interplay of
the energy error $\enorm{u-u_\grid}$ and data oscillation $\osc_\grid(f)_{-1}$: 
\begin{equation}\label{E:role-of-oscillation}
\begin{minipage}{0.85\linewidth}
\emph{Data oscillation $\osc_\grid(f)_{-1}$ may be responsible for the energy error $\enorm{u-u_\grid}$
to stagnate, even with the interior vertex property, and may entirely dominate it relative to the
error estimator $\est_\grid(u_\grid,f)$ over many mesh refinements unless it is reduced. }
\end{minipage}
\end{equation}

%------------------------------------------------------------------------------
\subsubsection{Reducing data oscillation}\label{S:reducing-osc}
%------------------------------------------------------------------------------
%
The PDE error estimator $\eta_\grid(u_\grid)$ in \eqref{E:pde-estimator} is fully discrete and thus
computable. In contrast, the computation, or rather estimation, of $\osc_\grid(f)_{-1}$ hinges on a priori
knowledge of $f$ and cannot be assessed in general. Assuming that the local indicators 
introduced in Lemma \ref{L:local-re-indexing} (localization re-indexing)
\begin{equation}\label{E:local-omegaT}
\osc_\grid(f,T)_{-1} = \| f - P_\grid f\|_{H^{-1}(\omega_T)} \quad T \in \grid.
\end{equation}
are computable without further regularity than $f\in H^{-1}(\Omega)$, it is natural to think of
{\it tree approximation} as the algorithm of choice to reduce $\osc_\grid(f)_{-1}$
\cite{BinevDeVore:2004,BinevFierroVeeser:2023,Binev:2018}. However, this optimal algorithm
is not readily applicable because of the lack of a suitable sub-additivity property.

On the other hand, greedy algorithms, such as that in Section \ref{S:constructive-approximation}
(constructive approximation), do not work under minimal regularity. 
In Section \ref{S:right-hand-side}, we present practical examples of rough $f$ for which 
$\osc_\grid(f)_{-1}$ can be replaced by a larger computable surrogate estimator
$\wt\osc_\grid(f)_{-1}$. The latter splits into element contributions and is amenable to a greedy
strategy. Since this is specialized and technical, we prefer to postpone the full discussion
to Section \ref{S:right-hand-side} and assume now the existence of a module $\DATA$ with the 
following property: given a tolerance $\tau>0$ and a conforming mesh $\grid\in\grids$, $\DATA$
constructs a conforming refinement $\grid_*\in\grids$
\[
[\grid_*] = \DATA \, (\grid, f, \tau),
\]
such that $\osc_{\grid_*} (f)_{-1} \le \tau$. The complexity of $\DATA$ depends on the decay 
rate of the best approximation error $ \min_{\grid\in\grids_N} \wt\osc_\grid(f)_{-1}$ of $f$ 
with $N$ degrees of freedom. We address this important issue in Section \ref{S:right-hand-side}
for each example separately.

%------------------------------------------------------------------------------
\subsubsection{Linear convergence}\label{S:linear-convergence}
%------------------------------------------------------------------------------
%
The following algorithm, $\AFEMSW$, a one-step \AFEM with switch, is a minor, but essential, modification of $\GALERKIN$ in that the call to
the modules $\MARK$ and $\REFINE$ is conditional to the size of $\osc_\grid(f)_{-1}$ relative to
$\est_\grid(u_\grid,f)$. This structure is consistent with the heuristic discussion in
\cite[Section 6]{CaKrNoSi:08} to avoid separate marking. A similar algorithm is
being developed in \cite{Kreuzer.Veeser.Zanotti}.

\begin{algo}[\AFEMSW: one-step AFEM with switch]\label{A:one-step-switch}
\index{Algorithms!\AFEMSW: \AFEM with switch}
Let $\grid_0$ be a suitable initial mesh, the coefficients $(\bA,c)$ be discrete over $\grid_0$,
and $\eps>0$ be a stopping tolerance.
Given parameters $0<\theta,\omega,\xi<1$, $\AFEMSW$ iterates the following loop until
$\est_\grid(u_\grid,f) \le \eps$:

\medskip
\begin{algotab}
  \> $[\grid,u_\grid] = \AFEMSW \, (\grid_0, \data, \eps)$
  \\
  \> \> set $j=0$\\
  \> \> do\\
  \> \> \> $[ u_{\grid_j} ] = \SOLVE \, (\grid_j)$
  \\
  \> \> \> $[\eta_{\grid_j}(u_{\grid_j}), \osc_{\grid_j}(f)_{-1} ] = \ESTIMATE \, (u_{\grid_j}, \grid_j, \data)$ 
  \\
  \> \> \> if $\est_{\grid_j}(u_{\grid_j},f) \le \eps$\\
  \> \> \> \> return $\grid_j, u_{\grid_j}$\\
  \> \> \> else if $\osc_{\grid_j}(f)_{-1} \le \sigma_j := \omega \est_{\grid_j}(u_{\grid_j},f)$
  \\
  \> \> \> \> $[\marked_j] = \MARK \, (\{\eta_{\grid_j}(u_{\grid_j},T)\}_{T\in\grid_j}, \grid_j, \theta)$
  \\
  \> \> \> \> $[ \grid_{j+1} ] = \REFINE \, (\grid_j, \marked_j)$
  \\
  \> \> \> else
  \\
  \> \> \> \> $[ \grid_{j+1} ] = \DATA \, (\grid_j,f,\xi\sigma_j)$
  \\
  \> \> \> $j \leftarrow j+1$
  \\
  \> \> \textrm{while true}
\end{algotab}
\medskip
\AB{
Note that $\SOLVE$ computes the Galerkin approximation using the exact right-hand side $f \in H^{-1}(\Omega)$ (not necessarily in $\mathbb F_{\grid_j}$), thereby preserving the Galerkin orthogonality property. Moreover,  $\ESTIMATE$ is now responsible for computing the PDE estimator 
$$
\eta_{\grid_j}(u_{\grid_j})=\eta_{\grid_j}(u_{\grid_j},f,\grid_{j})
$$ 
using $P_{\grid_j}f\in \mathbb F_{\grid_j}$,
as well as data oscillation $\osc_{\grid_j}(f)_{-1}$, which together give
\[
\est_{\grid_j}(u_{\grid_j},f) = \Big( \eta_{\grid_j}(u_{\grid_j})^2 + \osc_{\grid_j}(f)_{-1}^2 \Big)^{1/2},
\]
and $\MARK$ consists of D\"orfler marking \eqref{E:theta} with parameter $\theta$.
}
\end{algo}

\AB{
We proceed as in Section~\ref{S:contraction}
 to prove linear convergence of $\AFEMSW$. We
first show a contraction property for the {\it quasi-error}, which instead of \eqref{E:quasi-error} reads
\begin{equation}\label{E:quasi-error-sw}
\zeta_{\grid_j}(u_{\grid_j},f)^2 : = \enorm{u-u_{\grid_j}}^2 + \gamma \eta_{\grid_j}(u_{\grid_j})^2
+  \osc_{\grid_j} (f)_{-1}^2,
\end{equation}
where $u=u(\bA,c,f)$ is the Galerkin solution with $(\bA,c)$ discrete but $f$ exact and the scaling parameter satisfies $0<\gamma\le1$.

\begin{theorem}[contraction property of $\AFEMSW$]\label{T:contraction-sw}
Let $(\bA,c)$ be discrete coefficients over $\grid_0$ and let $f \in H^{-1}(\Omega)$. Let $\theta \in (0,1]$ be the D\"orfler parameter
and $(\grid_j,\V_j,u_j)$ be the sequence of conforming meshes $\grid_j$, finite element spaces $\V_j$,
and Galerkin solutions $u_j\in\V_j$ produced by {\rm \AFEMSW}. There exist parameters $0< \omega_0 <1$ 
sufficiently small and $0<\gamma \leq 1$ and $0<\alpha<1$ such that for any $\omega \leq \omega_0$ 
and $\xi \le 1/2$ the quasi-error $\zeta_{\grid_j}$ in
\eqref{E:quasi-error-sw} contracts \looseness=-1
\begin{equation}\label{E:zeta-contracts}
\index{Constants!$\alpha$: contraction constant}
\zeta_{\grid_{j+1}}(u_{\grid_{j+1}},f) \le \alpha \zeta_{\grid_j}(u_{\grid_j},f)  \quad j\ge 0.
\end{equation}
\end{theorem}
\begin{proof}
We argue as in Theorem~\ref{T:contraction} (general contraction property) upon distinguishing the two possible cases
within Algorithm~\ref{A:one-step-switch}.
But first, we must account for a crucial difference: the discrete forcing function $P_{\grid_j}f$ used in the definition of the estimator $\est_{\grid_j}(u_{\grid_j},f)$ changes in each iteration.
We use the same notation as in Theorem~\ref{T:contraction} along with $\osc_j:=\osc_{\grid_j}(f)_{-1}$, $\est_j^2:= \eta_j^2 + \osc_j^2$, and $P_j:=P_{\grid_j}$.

\medskip
\step{1} {\it Estimator reduction property}. 
In view of Proposition \ref{P:est-reduction} (estimator reduction), we need to estimate the discrepancy between discrete forcing functions
\[
\! \sum_{T\in\grid_{j+1}} \|P_{j+1} f - P_j f\|_{H^{-1}(\omega_T)}^2 
\le 2 \!\! \sum_{T\in\grid_{j+1}} \!\! \left(\| f - P_{j+1} f \|_{H^{-1}(\omega_T)}^2 
\!+\! \| f - P_j f \|_{H^{-1}(\omega_T)}^2 \right).
\]
For the first term we recall Lemma~\ref{L:quasi-mono-osc} (quasi-monotonicity of oscillation) to 
write
\[
\sum_{T\in\grid_{j+1}} \| f - P_{j+1} f \|_{H^{-1}(\omega_T)}^2 
= \osc_{j+1}^2 
\le C_{\textrm{osc}}^2 \osc_j^2.
\]
For the second term, instead, we combine the projection property $P_{j+1}(P_j f) = P_j f$ with 
Lemma \ref{L:loc-of-H^-1-norm} (localization of $H^{-1}$-norm) and 
Corollary~\ref{C:local-near-best-approx-of-P} (local near-best approximation), and the fact that
$\grid_{j+1}$ is a refinement of $\grid_j$, to see that
\begin{align*}
\sum_{T' \subset \omega_T} \| f - P_{j+1}(P_j f) \|_{H^{-1}(\omega_{T'})}^2 
&\le C_{\textrm{lStb}}^2 \sum_{T' \subset \omega_T} \| f - P_j f \|_{H^{-1}(\omega_{T'})_{-1}}^2
\\ &
\le C_{\textrm{lStb}}^2 C_{\textrm{ovrl}}^2 \osc_{\grid_j} (f,\omega_T)^2
\quad\forall T\in\grid_j.
\end{align*}
Adding over $T$ and recalling Proposition \ref{P:est-reduction} we end up with
\begin{equation}\label{E:est-reduction-forcing}
\begin{aligned}
      \eta_{\grid_{j+1}}(u_{j+1},f,\grid_{j+1})^2 &\le
      (1+\delta)\big(\eta_{\grid_j}(u_j,f,\grid_j)^2 - \lambda\,\eta_{\grid_j}(u_j,f,\marked_j)^2\big)
      \\ &+
      (1+\delta^{-1})\,\Clip^2 \, \big( |u_j - u_{j+1}|_{H^1_0(\Omega)}^2 +
      \osc_j^2 \big),
\end{aligned}
\end{equation}
for a constant $\Clip$ large enough to absorb all preceding constants, and any $\delta>0$.

\medskip
\step{2} {\it Case $\osc_j \le \omega \est_j$.} We first observe that $\eta_j^2\ge (1-\omega^2) \est_j^2$ and
$\osc_j^2 \le \frac{\omega^2}{1-\omega^2} \eta_j^2$. We then proceed as in Theorem \ref{T:contraction} with the
quantity $e_j^2 + \gamma \eta_j^2$, and observe that the choices of $\delta$ and $\gamma$
\begin{equation}
    \label{e:constraint_d_gamma}
\delta \leq -1 + \frac{1 - \frac{\lambda\theta^2}{2}}{1-\lambda\theta^2} =
\frac{\lambda\theta^2}{2(1-\lambda\theta^2)},
\qquad
\gamma \leq \frac{\delta}{4\Clip^2} \leq \frac{1}{2\big(1+\delta^{-1}\big) \Clip^2},
\end{equation}
imply $\gamma(1+\delta^{-1})\Clip^2 \leq 1/2$. This, together with \eqref{E:Pythagoras}
and \eqref{E:est-reduction-forcing}, leads to
\[
e_{j+1}^2 + \gamma \eta_{j+1}^2 \le e_j^2 + \gamma \Big(1-\frac{\lambda\theta^2}{2}\Big) \eta_j^2+\frac 1 2 \osc_j^2;
\]
compare with \eqref{e:monoton_galerkin_total_error}.
We invoke the upper bound in \eqref{E:aposteriori-bounds} to write 
$$
\eta_j^2 \ge (1-\omega^2)\est_j^2 \ge (1-\omega^2)\frac{e_j^2}{C_1^2} \ge \frac{e_j^2}{2C_1^2}
$$ 
provided $\omega^2 \leq 1/2$, whence
\begin{equation*}\label{E:error+eta}
e_{j+1}^2 + \gamma \eta_{j+1}^2 \le \Big(1 - \frac{\gamma\lambda\theta^2}{8C_1^2}\Big) e_j^2 +
\gamma \Big(1-\frac{\lambda\theta^2}{8} \Big) \eta_j^2 - \frac{\gamma\lambda\theta^2}{8} \eta_j^2+\frac 1 2 \osc_j^2. 
\end{equation*}
We next consider the data oscillation, for which we invoke Lemma~\ref{L:quasi-mono-osc} (quasi-monotonicity of oscillation) %have the {\it quasi-monotonicity} property \eqref{e:quasi-mono-osc}
\begin{equation*}\label{E:quasi-monotonicty}
\osc_{j+1} \le C_\textrm{osc} \osc_j, \qquad   \osc_j^2 \le C_\textrm{osc}^2 \frac{\omega^2}{1-\omega^2} \eta_j^2 \leq 2 C_\textrm{osc}^2 \omega^2 \eta_j^2.
\end{equation*}
Adding the two preceding inequalities yields
\begin{align*}
\zeta_{j+1}^2 = e_{j+1}^2 + \gamma \eta_{j+1}^2 +  \osc_{j+1}^2
\le &  \Big(1 - \frac{\gamma\lambda\theta^2}{8C_1^2}\Big) e_j^2
\\
& + \Big(1-\frac{\lambda\theta^2}{8} \Big)  \big(\gamma \eta_j^2 + \osc_j^2 \big)
\\
& + \Big[ -\gamma \frac{\lambda\theta^2}{8} + 2\big(C_\textrm{osc}^2 - 1 + \frac{\lambda\theta^2}{8}
+ \frac 1 2\big) \omega^2 \Big] \eta_j^2.
\end{align*}
We drop the term $- \frac 12 + \frac{\lambda\theta^2}{8} \leq 0$ and  let $\gamma= \frac{\delta}{4\Clip^2}$, which is consistent with \eqref{e:constraint_d_gamma}.
We seek conditions on $\omega$ that make the factor of $\eta_j^2$ non-positive. Imposing
\begin{equation}\label{E:condition-1}
\omega^2 \le \frac{\gamma\lambda\theta^2}{16C_\textrm{osc}^2} = \frac{\lambda\theta^2}{64 C_\textrm{osc}^2 \Clip^2}\delta
\end{equation}
yields
$$
\zeta_{j+1}^2 \leq \alpha_1^2 \zeta_j^2
$$
with
\[
\alpha_1^2 := \max \Big\{ 1 - \frac{\delta \lambda\theta^2}{32C_1^2\Clip^2}, 1 - \frac{\lambda\theta^2}{8} \Big\} < 1.
\]

\step{3} {\it Case $\osc_j > \omega\est_j$}. The module $\DATA$ with input parameter $\xi \le 1/2$ gives
\[
\osc_{j+1} \le \xi \omega \est_j < \xi \osc_j.
\]
We now exploit the contraction of $\osc_j$ to compensate the moderate increase of
$\eta_j^2$ and presence of $\osc_j^2$, both governed by \eqref{E:est-reduction-forcing}. 
In fact, $\gamma(1+\delta^{-1})\Clip^2 \leq \frac 1 2$ yields
\[
e_{j+1}^2 + \gamma \eta_{j+1}^2 \le e_j^2 + \gamma (1+\delta) \eta_j^2+\frac 1 2 \osc_j^2.
\]
We add $\osc_{j+1}^2$ to both sides and rewrite the right-hand side to arrive at
\begin{align*}
\zeta_{j+1}^2 = e_{j+1}^2 + \gamma \eta_{j+1}^2 +  \osc_{j+1}^2  \le &
e_j^2 - \frac{1-2\xi^2}{8} \osc_j^2 \\
&+(1-\delta) \gamma \eta_j^2  +
\left( \frac{1+2\xi^2}{4}+ \frac 1 2 \right) \osc_j^2 \\
&+ 2 \delta \gamma \eta_j^2 - \frac{1-2\xi^2}{8} \osc_j^2.
\end{align*}
Our next task is to find conditions on $\omega$ for the last line to be non-positive.
To this end, we resort to the upper bound $\eta_j^2< \omega^{-2} \osc_j^2$ and $\xi \le 1/2$ to obtain
$$
2 \delta \gamma \eta_j^2 -  \frac{1-2\xi^2}{8}\osc_j^2 < \left( \frac{2\delta \gamma}{\omega^2}- \frac{1}{16}\right) \osc_j^2 \leq 0
$$
provided we impose the relation
\begin{equation}\label{E:condition-2}
\omega^2 \geq 32 \delta \gamma = \frac{8}{\Clip^2}\delta^2.
\end{equation}
We next use the upper bound $e_j \leq C_1 \mathcal E_j \leq C_1 \omega^{-1} \osc_j$ to write
$$
e_j^2 - \frac{1-2\xi^2}{8}\osc_j^2 \leq \left( 1- \frac{\omega^2}{16 C_1^2}\right) e_j^2,
$$
whence we end up with
\[
\zeta_{j+1}^2 \leq \alpha_2^2 \zeta_j^2
\]
provided we define
\begin{equation*}
\alpha_2^2 := \max\Big\{1-\frac{\omega^2}{16 C_1^2}, 1-\delta, \frac{3+2\xi^2}{4}  \Big\}<1.
\end{equation*}

\step{4} {\it Choosing the parameters}. We see that the asserted estimate
\eqref{E:zeta-contracts} is valid with $\alpha = \max\{\alpha_1,\alpha_2\}<1$ provided the constraints
\eqref{E:condition-1} and \eqref{E:condition-2} are compatible, i.e., 
$$
\frac{8}{\Clip} \delta^2 \le \omega^2 \le \frac{\lambda \theta^2}{64 C_\textrm{osc}^2 \Clip^2}\delta.
$$
We choose $\delta_0 = \frac{\lambda\theta^2}{512 C_{\textrm{osc}}^2 \Clip}$ and 
$\omega_0 = \frac{\lambda \theta^2}{128 C_\textrm{osc}^2 \Clip \sqrt{2\Clip}}$. Then,
for all $\omega\le\omega_0$, there exists $\delta \le \delta_0$ that satisfies the previous
inequalities as well as $\gamma = \frac{\delta}{4\Clip^2} \leq 1$, perhaps upon reducing $\delta_0$.
This completes the proof of Theorem~\ref{E:zeta-contracts}.
\end{proof}
}

Note that
we could replace the conditional $\osc_{\grid_j}(f)_{-1} \le \omega \, \est_{\grid_j}(u_{\grid_j},f)$ by
$\osc_{\grid_j}(f)_{-1} \le \omega \, \eta_{\grid_j}(u_{\grid_j})$, but the tolerance $\tau$ of $\DATA$ cannot
be
\[
\tau = \xi \, \omega \, \eta_{\grid_j}(u_{\grid_j})
\]
because the algorithm may not terminate when $\eta_{\grid_j}(u_{\grid_j})=0$, see e.g. Examples \ref{EX:chekerboard-f} - \ref{EX:Pgrid-n>1}.
In fact, the tolerance
$\tau = \xi \omega\est_{\grid_j}(u_{\grid_j},f)$ is {\it dynamic} and relative to $\est_{\grid_j}(u_{\grid_j},f)$.
This avoids separate marking, which was shown in \cite[Section 6]{CaKrNoSi:08} to give
non-optimal convergence rates. In contrast, we will prove in Section \ref{S:conv-rates-coercive}
that Algorithm \ref{A:one-step-switch} is rate optimal.

\rhn{It turns out that} Theorem~\ref{T:contraction-sw} yields linear convergence of error and estimator.

\begin{corollary}[linear convergence of error]\label{C:linear-convergence-error-sw}
\AB{For $0<\alpha<1$ and $0<\omega \leq \omega_0$, $\xi\le\frac{1}{2}$ as in Theorem \ref{T:contraction-sw},
and $C_* = \big( 1 +C_2^{-1} \big)^{1/2}$ with $C_2$ as in \eqref{E:aposteriori-bounds},}
there holds
\[
|u-u_{\grid_k}|_{H^1_0(\Omega)} \le C_* \, \alpha^{k-j} |u-u_{\grid_j}|_{H^1_0(\Omega)} \quad\forall k\ge j \ge 0.
\]
\end{corollary}
\begin{proof}
We again use the same notation as in Lemma~\ref{L:contraction-discrete-data} and Theorem~\ref{T:contraction-sw}. In view of the definition \eqref{E:quasi-error-sw} of quasi-error $\zeta_j:=\zeta_{\grid_j}(u_{\grid_j},f)$, we thus have $e_j\le\zeta_j$
and
\AB{
\begin{align*}
\zeta_j^2 \le e_j^2 +  \eta_j^2 + \osc_j^2  
\le \big( 1 +C_2^{-1} \big) e_j^2
\end{align*}
}
because $C_2 \est_j\le e_j$ from \eqref{E:aposteriori-bounds}. This implies
\[
e_j \le \zeta_j \le C_* e_j \quad\forall j\ge 0,
\]
and invoking Theorem \ref{T:contraction-sw} (contraction property for $\AFEMSW$)
\[
e_k^2 \le \zeta_k^2 \le \alpha^{2(k-j)} \zeta_j^2 \le \alpha^{2(k-j)} C_*^2 e_j
\]
gives the desired estimate.
\end{proof}

We stress that Corollary \ref{C:linear-convergence-error-sw} relies on the lower bound in
\eqref{E:aposteriori-bounds} whereas Corollary \ref{C:linear-convergence-est-sw} uses only the upper bound.
Its proof is similar and thus omitted.

\begin{corollary}[linear convergence of estimator]\label{C:linear-convergence-est-sw}
\AB{For $0<\alpha<1$ and $0<\omega \leq \omega_0$, $\xi \le \frac 1 2$ as in Theorem \ref{T:contraction-sw},
and $C_\# = \big( (1+C_1^2)\gamma^{-1}\big)^{1/2}$ with $C_1$ as in \eqref{E:aposteriori-bounds},}
there holds
\[
\est_{\grid_k}(u_{\grid_k},f) \le C_\# \, \alpha^{k-j} \, \est_{\grid_j}(u_{\grid_j},f) \quad\forall k\ge j \ge 0.
\]
\end{corollary}

%%%%%%%%%%%%%%%%%%%%%%%%%%%%%%%%%%%%%%%%%%%%%%%%%%%%%%%%%%%%%%%%%%%%%%%%%%%%%%%%
\subsection{Convergence for general data: two-step {\rm AFEM}}\label{S:general-data}
%%%%%%%%%%%%%%%%%%%%%%%%%%%%%%%%%%%%%%%%%%%%%%%%%%%%%%%%%%%%%%%%%%%%%%%%%%%%%%%%

We now remove the restriction of Sections \ref{S:discrete-data} and \ref{S:one-step-switch}
of discrete data and allow for general data $\data = (\bA,c,f) \in\D$ as defined in
\eqref{E:space-discrete-data}.
The current goal is to study Algorithm \ref{algo:AFEM-Kernel} ($\AFEMTS$), which 
concatenates the modules $\DATA$ and $\GALERKIN$.
We start with the study of continuous dependence with respect to data $\data$. We next discuss
the approximation of $\data$ within the module $\DATA$, the computational cost of $\GALERKIN$,
and eventually the convergence of Algorithm \ref{algo:AFEM-Kernel}.

%-------------------------------------------------------------------------------
\subsubsection{Perturbation theory}\label{S:perturbation}
%-------------------------------------------------------------------------------
%
We start with a brief discussion of data perturbation. Given constants $0<\alpha_1\leq \alpha_2$ and $0 \le c_1 \leq c_2$, we define the constrained spaces for the diffusion and reaction coefficients by
\index{Constants!$(\alpha_1,\alpha_2)$: lower and upper bounds of the diffusion coefficient spectrum}
\index{Constants!$(c_1,c_2)$: lower and upper bounds of the reaction coefficient}%
\begin{equation}\label{d:DiffusionSpace}
\index{Functional Spaces!$M(\alpha_1,\alpha_2)$: Admissible set for $\bA$}
\begin{split}
    M(\alpha_1,\alpha_2) := \big\{\bA\in L^\infty(\Omega; \R^{d\times d}_{\textrm{sym}}):
\ & 0<\alpha_1\le\lambda_j(\bA(x))\le \alpha_2 \\
&\textrm{for a.e. }x\in\Omega, 1\le j \le d \big\},
\end{split}
\end{equation}
where $\lambda_j(\bA(x))$ denotes the $j$-th eigenvalue of $\bA$ at
$x\in\Omega$
and
\begin{equation}\label{d:ReactionSpace}
\index{Functional Spaces!$R(c_1,c_2)$: Admissible set for $c$}
R(c_1,c_2) := \big\{ c \in L^\infty(\Omega): c_1 \le c(x) \le c_2 \quad \textrm{for a.e. }x\in\Omega \big\}.
\end{equation}

The coefficients $(\bA,c)$ are assumed to satisfy the structural assumption
\begin{equation}\label{E:structural-assumption}
\index{Assumptions!Structural assumption for exact data}
\bA \in M(\alpha_1,\alpha_2), \quad c \in R(c_1,c_2);
\end{equation}
see \eqref{uniform_spd}. 
This guarantees coercivity and continuity of the bilinear form
$\B$ in \eqref{bilinear-coercive}, and thus unique solvability of \eqref{weak-form}.

Regarding the discrete coefficients, $(\wh{\bA},\wh{c})$ will ultimately be piecewise
polynomials in a grid $\wh\grid\in\grids$. The side constraints in \eqref{d:DiffusionSpace} and \eqref{d:ReactionSpace} are generally violated by any linear projection onto piecewise polynomials of degree $n-1 \ge 1$,
e.g. the $L^2$-projection, and require a nonlinear correction maintaining high order accuracy. This is a crucial but delicate matter addressed later in Section
\ref{ss:positivity}.
For the moment, we simply assume that the 
discrete coefficients $(\wh{\bA},\wh{c})$ satisfy 
\begin{equation}\label{E:structural-assumption-wh}
\index{Assumptions!Structural assumption for discrete data}
\index{Constants!$(\halpha_1,\halpha_2)$: lower and upper bounds of the approximate diffusion coefficient spectrum}
\index{Constants!$(c_1,c_3)$@$(\hc_1,\hc_2)$: lower and upper bounds of the approximate reaction coefficient}
\hbA \in M(\halpha_1,\halpha_2), \quad \hc \in R(\hc_1,\hc_2),
\end{equation}
with \looseness=-1
\begin{equation}\label{E:wh-alpha}
\index{Constants!$\Cctr$: constrain upper bound amplification constant}
\frac{\alpha_1}2 \leq \halpha_1 \leq \halpha_2 \leq \Cctr \alpha_2, \qquad 
-\frac{\alpha_1}{4C_P^2} \leq \hc_1 \leq \hc_2 \leq \Cctr (\alpha_1+c_2)
\end{equation}
where $C_P>0$ is the Poincar\'e constant in \eqref{E:Poincare} and $\Cctr \geq 1$ is a constant; see \eqref{e:halphas} and \eqref{e:def_hc}.
This implies coercivity and continuity of the perturbed bilinear form 
\begin{equation}\label{E:perturbed-B}
\wh\B[v,w]:= \int_\Omega \nabla v \cdot \hbA \nabla w + \hc v w, \qquad \forall v,w \in H^1_0(\Omega),
\end{equation}
because for all $v,w \in H^1_0(\Omega)$
\[
\wh\B[v,v] \ge \halpha_1 \int_\Omega  |\nabla v|^2 -\frac{\alpha_1}{4C_P^2} \int_\Omega |v|^2
\ge  \frac{\alpha_1}{4} |v|_{H^1_0(\Omega)}^2
\]
and
\[
\big|\wh\B[v,w] \big| \le \int_\Omega \halpha_2 |\nabla v| \, |\nabla w| + \wh{c}_2 |v| \, |w|
\le \big(\halpha_2 + \wh{c}_2 C_P^2 \big) |v|_{H^1_0(\Omega)} |w|_{H^1_0(\Omega)}.
\]
Therefore, the energy norm $\enorm{v}^2 = \wh\B[v,v]$ is equivalent to the $H^1_0$-seminorm
\begin{equation}\label{E:enorm-Vnorm-wh}
\index{Constants!$(c_{\wh\B},C_{\wh\B})$: norm equivalence constants for the perturbed problem}
c_{\wh\B} |v|_{H^1_0(\Omega)}^2 \le \enorm{v}^2 \le C_{\wh\B} |v|_{H^1_0(\Omega)}^2,
\end{equation}
where $c_{\wh\B} = \frac{\alpha_1}{4}$ and $C_{\wh\B} = \halpha_2 + \wh{c}_2 C_P^2$.
Hence, the Lax-Milgram Theorem guarantees the
existence of a unique solution $\hu = u(\wh{\data})\in H^1_0(\Omega)$ of the perturbed
problem \eqref{E:perturbed-weak-form} defined using the discrete data $\wh \data=(\wh{\bA},\wh c,\wh f)$. 

We now quantify the effect of perturbing data from $\data$ to $\wh\data$ in the space
\begin{equation}\label{E:space-DO}
\index{Functional Spaces!$D(\Omega)2$@$\wh{D}(\Omega)$: temporary metric space for the data perturbation}
\wh{D}(\Omega) := L^r(\Omega; \R^{d\times d}) \times W^{-s}_q(\Omega) \times H^{-1}(\Omega),
\end{equation}
where $2\le r \le \infty$ and $0 \le s \le 1, \frac{d}{2-s} < q \le \infty$; $W^{-s}_q(\Omega)$\index{Functional Spaces!$W^{-s}_q(\Omega)$: dual of $W^{s}_{q*}(\Omega)$ with $q^*=\frac{q}{q-1}$} is the dual of $W^{s}_{q*}(\Omega)$ with $q^*=\frac{q}{q-1}$.  
The use of $r=\infty$ for $\bA$ entails the further assumption
\begin{equation}\label{E:assAunicont}
\text{\em $\bA$ is piecewise uniformly continuous
over a generic mesh $\grid\in\grids$,}
\end{equation}
which turns out to be rather restrictive but customary in the theory of $\AFEM$. Our present
approach allows for $r<\infty$ and thus for discontinuous coefficients $(\bA,c)$ non aligned with
$\grid$, which is important in practice. However, it requires the following slightly stronger regularity
property of the solution $u \in H^1_0(\Omega)$ of \eqref{weak-form}:
\begin{equation}\label{E:reg-w1p}
\|\nabla u\|_{L^p(\Omega)} \le C_p \|f\|_{W^{-1}_p(\Omega)}
\quad
2< p \le p_0.
\end{equation}
We refer to Lemma \ref{L:Wp-regularity} ($W^1_p$-regularity) 
that shows the existence of $C_p>0$ and
$p_0>2$ that only depend on $\Omega, \alpha_1, \alpha_2$ and $c_2$.

\begin{lemma}[continuous dependence on data] \label{L:perturbation}
Let $\data=(\bA,c,f) \in \D$ be such that $\bA\in M(\alpha_1,\alpha_2)$ and $c\in R(c_1,c_2)$. Let
$\wh\data=(\wh\bA,\wh{c},\wh{f})\in\D$ be an approximation of $\data$ such that
$\wh{\bA} \in M(\wh{\alpha}_1,\wh{\alpha}_2)$ and $\wh{c} \in R(\wh{c}_1,\wh{c}_2)$.
Let $2\le r \le \infty, 2 \le r_* = \frac{2r}{r-2} \le p_0$ be so that $f \in W^{-1}_{r_*}(\Omega)$.
If $u = u(\data), \wh{u} = u(\wh{\data}) \in H^1_0(\Omega)$ are the solutions of \eqref{weak-form}
and \eqref{E:perturbed-weak-form} with data $\data, \wh{\data}$, respectively, and $u$ satisfies
\eqref{E:reg-w1p} with $p = r_*$ in case $r<\infty$, then
for any $0\le s \le 1$ and $\frac{d}{2-s} < q \le \infty$ there holds
\begin{equation}\label{E:perturbation}
 \| \nabla (u- \hu) \|_{L^2(\Omega)} \le C(\data,\Omega) \, \|\data-\hdata\|_{\wh{D}(\Omega)},
\end{equation}
where the constant $C(\data,\Omega)$ depends on $\data$, $\Omega$, $p_0$, $q$ and $s$, and blows up
as $q\to\frac{d}{2-s}$ for $d=2$ while it remains bounded for $d>2$.
\end{lemma}
\begin{proof}
Subtracting the weak formulations \eqref{weak-form} for $u$ and \eqref{E:perturbed-weak-form} for $\wh{u}$,
and reordering, we easily obtain for any $v \in H^1_0(\Omega)$
$$
\int_\Omega \nabla v \cdot \hbA \nabla (u-\hu) + \hc v (u-\hu) =
\int_\Omega \nabla v \cdot (\hbA-\bA) \nabla u + (\hc-c) v u + \langle f-\hf,v \rangle.
$$
We choose $v = u-\hu \in H^1_0(\Omega)$ and invoke \eqref{E:enorm-Vnorm-wh} to deduce
\begin{equation*}\label{e:perturb_coerc}
c_{\wh\B} \| \nabla v \|_{L^2(\Omega)}^2 \leq \int_\Omega \nabla v \cdot (\hbA-\bA) \nabla u + (\hc-c) v u + \langle f-\hf,v \rangle.
\end{equation*}
We estimate each term separately, starting with the first and last terms
\begin{equation}\label{e:perturb_a}
\int_\Omega \nabla v \cdot (\hbA-\bA) \nabla u \leq \| \hbA-\bA  \|_{L^r(\Omega)} \| \nabla u \|_{L^{r^*}(\Omega)}
\| \nabla v \|_{L^2(\Omega)}
\end{equation}
with $2 \le r_* = \frac{2r}{r-2} \le p_0$, as well as
\begin{equation*}\label{e:perturb_f}
\langle f-\hf,v \rangle \leq \| f - \hf \|_{H^{-1}(\Omega)} \| \nabla v \|_{H^1(\Omega)}.
\end{equation*}
For the reaction term, which is more delicate, we invoke the duality pairing $W^s_{q'}-W^{-s}_{q}$ for any $0\le s \le 1$ and  $q' = \frac{q}{q-1} \geq 1$ to obtain
\begin{equation*}\label{e:perturb_reac}
\int_\Omega (\hc-c) v u \le  \| \hc -c \|_{W^{-s}_q(\Omega)} | vu |_{W^{s}_{q'}(\Omega)}.
\end{equation*}
We now estimate $| vu |_{W^{s}_{q'}(\Omega)} \lesssim | vu |_{W^{1}_{p'}(\Omega)}$ where
$1/p' = \min\{ 1, (1-s)/d + 1/q'\}$ guarantees that $W^{1}_{p'}(\Omega) \subset W^{s}_{q'}(\Omega)$ \cite[Theorem 14.32]{Leoni:2009}.
Recalling that $q>d/(2-s)$, we deduce
\[
\frac{1-s}{d} + \frac{1}{q'}=\frac{1-s}{d} + 1 - \frac{1}{q} > \frac{d-1}{d} \ge \frac{1}{2},
\]
whence $1/p' > 1/2$ and there exists $t < \infty$ satisfying $1/t + 1/2 = 1/p'$ and 
\begin{equation*}\label{e:prod_ws}
| vu |_{W^{s}_{q'}(\Omega)} \lesssim  \| \nabla v\|_{L^2(\Omega)} \| u \|_{L^{t}(\Omega)} + \| v \|_{L^t(\Omega)} \| \nabla u \|_{L^{2}(\Omega)}.
\end{equation*}
Using the definition of $p'$, we obtain the explicit expression $t = \max\{2, t_0\}$, where
\[
t_0 = \frac{2dq}{q \big(2(1-s) + d\big) - 2d}.
\]
Moreover, for the Sobolev embedding $H^1(\Omega) \hookrightarrow  L^t(\Omega)$ we require
\[
1 - d \Big( \frac{1}{2} - \frac{1}{t} \Big) > 0
\quad\Rightarrow\quad
q > \frac{d}{2-s},
\]
which is our assumption on $q$. Therefore, as $q\to\frac{d}{2-s}$, we see that
$t_0\to\frac{2d}{d-2}$ and the limit is
infinite for $d=2$ but finite  and larger than $2$ for $d>2$.
Sobolev embedding together with the first Poincar\'e inequality \eqref{E:Poincare} gives the estimate
\[
\| vu \|_{W^{s}_{q'}(\Omega)} \le C(\Omega,t) \| \nabla v\|_{L^2(\Omega)} \| \nabla u \|_{L^{2}(\Omega)}.
\]
where $C(\Omega,t)$ is proportional to $t$ for $d=2$.
We finally observe that the factors $\| \nabla u \|_{L^{r_*}(\Omega)}$ and $\| \nabla u \|_{L^2(\Omega)}$
appear in the estimates of the coefficients $\bA$ and $c$, thereby reflecting the multiplicative
nature of these terms. Since $2\le r_* \le p_0$, they can be further bounded in terms
of $\|f\|_{W^{-1}_{r_*}(\Omega)}$ according to \eqref{E:reg-w1p}.
This in conjunction with the preceding estimates yields the assertion \eqref{E:perturbation}.
\end{proof}

A natural and rather popular choice of parameters $(r,q,s)$ in Lemma \ref{L:perturbation}
(continuous dependence on data) is $r=q=\infty$ and $s=0$, but this would prevent the
coefficients $(\bA,c)$ from being discontinuous within elements; see \eqref{E:assAunicont}.
We will explore this matter further in Section \ref{S:data-approx} (data approximation).

\begin{remark}[$L^2$-approximation of $\bA$]\label{R:L2-A}
\rm
It is appealing to estimate the distortion $\bA - \hbA$ in $L^2(\Omega)$ rather than in $L^r(\Omega)$
because it is a simpler norm to deal with.
Since $\| \bA\|_{L^\infty(\Omega)}\le\alpha_2$, $\|\hbA \|_{L^\infty(\Omega)} \leq \halpha_2$,
and $2 \le r \le \infty$, we deduce
\[
\| \bA - \hbA \|_{L^r(\Omega)} \leq 
\| \bA - \hbA \|_{L^\infty(\Omega)}^{1-\frac{2}{r}} \| \bA - \hbA \|_{L^2(\Omega)}^{\frac{2}{r}}
\lesssim \| \bA - \hbA \|_{L^2(\Omega)}^{\frac{2}{r}}.
\]
However, this may be sub-optimal in general.
One important situation where this is sharp corresponds to $\bA$ being
piecewise constant with jump discontinuities across a Lipschitz hypersurface $\gamma$ and $\hbA=\bA$ on every element $T \in \grid$ not intersecting $\gamma$. In that case, the equivalence
$$
\| \bA - \hbA \|_{L^{p}(\Omega)} \approx \big| \{ x \in \Omega \ : \ \bA(x) \not = \hbA(x) \} \big|^{\frac{1}{p}}, 
$$
is valid for $1 \le p \leq \infty$, whence
$$
\| \bA - \hbA \|_{L^r(\Omega)} \approx \| \bA - \hbA \|_{L^2(\Omega)}^{\frac{2}{r}}.
$$
\end{remark}

%-------------------------------------------------------------------------------
\subsubsection{Approximation of $\data$: module {\rm \DATA}}\label{S:DATA}
%-------------------------------------------------------------------------------
%

In this section we briefly discuss the structure of $\DATA$, which is the module of Algorithm
\ref{algo:AFEM-Kernel} (two-step $\AFEM$) responsible for data approximation.

In the sequel we will no longer rely on the Banach space $\wh{D}(\Omega)$ defined in \eqref{E:space-DO}
and used in Lemma \ref{L:perturbation} (continuous dependence on data). We rather restrict
the error notion to the following stronger Banach space
\begin{equation}\label{E:space-DO-redefine}
\index{Functional Spaces!$D(\Omega)$: metric space for the data perturbation}
D(\Omega) := L^r(\Omega; \R^{d\times d}) \times L^q(\Omega) \times H^{-1}(\Omega),
\end{equation}
where $q=2$ for $d<4$ or $q > \frac{d}{2}$ for $d\ge 4$; we justify the choice of $q$ below.
Let $\D$ and $\D_\grid$ be the spaces defined in \eqref{E:space-data} and \eqref{E:space-discrete-data}
for a conforming mesh $\grid\in\grids$. Given $\data=(\bA,c,f)\in\D$, let 
$\delta_\grid (\data)$ be the best approximation error of $\data$ within $\D_\grid$
measured in the space $D(\Omega)$, namely
\begin{equation}\label{E:data-oscillation}
\delta_\grid (\data) := \inf_{\data_\grid \in \D_\grid} \| \data - \data_\grid \|_{D(\Omega)}.
\end{equation}
This quantity characterizes the approximation quality of $\D_\grid$, thereby having theoretical value.
Since $\delta_\grid (\data)$ is hard to access in
view of the norms involved in the definition of $D(\Omega)$, the module $\DATA$ computes the
surrogate quantity 
$$
\index{Error Estimators!$\osc_\grid(\data)$: surrogate for the data error}
\osc_\grid(\data):= \|\data - \wh\data\|_{D(\Omega)} 
$$
for some approximation $\wh \data \in \D_\grid$ to be specified below. 

\begin{assumption}[properties of $\DATA$]\label{A:prop-data}
\index{Assumptions!Properties of \DATA}
Given a conforming mesh $\grid\in\grids$ and a tolerance $\tau>0$, the call
\begin{equation}\label{E:data}
[\wh\grid,\wh\data] = \DATA \,(\grid,\data,\tau)
\end{equation}
creates an admissible refinement $\wh\grid$ of $\grid$ and discrete data
$\wh\data=\data_{\wh\grid}\in\D_{\wh\grid}$ such that for a constant $C_\textrm{data}$\looseness=-1
\begin{equation}\label{E:bound-data}
\index{Constants!$\Cdata$: \DATA approximation constant}
\osc_{\wh\grid}(\data) := \|\data - \wh\data\|_{D(\Omega)} \le C_\textrm{data}\tau
\end{equation}
as well as the structural conditions \eqref{E:structural-assumption-wh} are achieved
in a finite number of iterations that depends on the regularity of $\data$, and such that
\begin{equation}\label{E:almost-best}
\index{Constants!$\Ldata$: \DATA quasi-optimality constant}
\osc_{\wh\grid} (\data) \le \Ldata \, \delta_{\wh\grid} (\data)
\end{equation}
with $\Ldata \ge 1$ depending only on the
shape regularity of $\grids$, the polynomial degree $n$ and the Lebesgue exponents in the space $D(\Omega)$.
\end{assumption}

In view of Lemma \ref{L:perturbation} (continuous dependence on data), there exists a constant
$C_D>0$ depending on $\data, \Omega$, and the shape regularity of $\grids$, such that the exact solutions
$u = u(\data)$ and $\hu = u(\wh{\data})$ of \eqref{strong-form} and \eqref{E:perturbed-weak-form},
corresponding to data $\data$ and $\wh\data$ respectively, satisfy the error estimate
\begin{equation}\label{E:error-DATA}
\index{Constants!$C_D$: \DATA constant}
|u-\hu|_{H^1_0(\Omega)} \le C_D \tau.
\end{equation}
A brief discussion follows about computing $\osc_{\wh\grid}(\data)$, where $\wh\grid$ remains fixed and
is replaced by $\grid$ to simplify the notation. Specific details
are given later in Assumptions \ref{A:approx-data} and \ref{A:optim-data} of Section \ref{A:approx-data}
and especially in Section \ref{S:data-approx}.

\paragraph{Approximating the coefficients.}

We now construct approximations $(\hbA,\hc)$ using  local $L^2$-projections and emphasize that this does not enforce the side constraints in the structural assumption \eqref{E:structural-assumption-wh}. We propose in Section~\ref{S:data-approx} a nonlinear correction satisfying the side constraints without sacrificing the accuracy. 

Given $T\in\grid$, and $v \in L^p(T)$ with $1\le p \le \infty$, we denote by $\Pi_T v:=\Pi_T^{n-1}v$
\index{Operators!$\Pi_K$, $\Pi_K^m$: $L^2$ projection onto $\P_m(K)$}
the $L^2$-projection
of $v$ onto the space $\P_{n-1}$ of polynomials of degree $\le n-1$, namely
\begin{equation}\label{E:L2-projection-def}
\Pi_T v\in \P_{n-1}: \quad \int_T \Pi_T v \, w = \int_T v \, w \quad\forall w\in \P_{n-1}.
\end{equation}

\begin{lemma}[$L^p$-stability of $\Pi_T$]\label{L:Lp-stab}
For every $1\le p \le \infty$ and $v\in L^p(T)$, 
there exists a constant $C$ depending on $p,n$ and the shape regularity of $\grids$ such that
\begin{equation}\label{E:Lp-stab}
\| \Pi_T v \|_{L^p(T)} \le C \| v \|_{L^p(T)} \quad\forall T\in\grid.
\end{equation}
\end{lemma}
\begin{proof}
It is trivial to see that $\| \Pi_T v \|_{L^2(T)} \le \| v \|_{L^2(T)}$. Let $2<p\le\infty$ and
combine an inverse estimate with a H\"older inequality to write
\[
\| \Pi_T v \|_{L^p(T)} \le C h_T^{\frac{d}{p}-\frac{d}{2}}\| \Pi_T v \|_{L^2(T)} 
\le C h_T^{\frac{d}{p}-\frac{d}{2}} \| v \|_{L^2(T)} \le C \| v \|_{L^p(T)}.
\]
For $1\le p < 2$ we proceed by duality. Let $\varphi\in L^q(T)$ with $q=\frac{p}{p-1}$. Then
\[
\int_T \Pi_T v \varphi = \int_T v \Pi_T \varphi \le \|v\|_{L^p(T)} \|\Pi_T \varphi\|_{L^q(T)}
\le C \|v\|_{L^p(T)} \|\varphi\|_{L^q(T)},
\]
which implies \eqref{E:Lp-stab} and concludes the proof.
\end{proof}

We immediately have the following simple consequence of Lemma \ref{L:Lp-stab}.

\begin{corollary}[best approximation of $\Pi_T$]\label{C:best-approx}
For every $1\le p \le \infty$ and $v\in L^p(T)$, 
there exists a constant $\Cba\ge1$ depending on $p,n$ and the shape regularity of $\grids$ such that
\begin{equation}\label{E:best-approx}
\index{Constants!$\Cba$: best approximation constant of $\Pi_T$}
\| v - \Pi_T v \|_{L^p(T)} \le \Cba \inf_{w\in \P_{n-1}} \| v - w \|_{L^p(T)}.
\end{equation}
\end{corollary}
\begin{proof}
We combine the invariance of $\Pi_T$ on $\P_{n-1}$, i.e. $\Pi_T w = w$ for $w \in\P_{n-1}$,
with \eqref{E:Lp-stab} to see that
\[
\|v-\Pi_T v\|_{L^p(T)} = \|(v-w)-\Pi_T (v-w)\|_{L^p(T)} \le C \|v-w\|_{L^p(T)}.
\]
This implies \eqref{E:best-approx} as asserted.
\end{proof}

The $L^2$-projection is easily computable because it entails solving the linear system
\eqref{E:L2-projection-def}. However, this flexibility comes at the expense of a best approximation
constant $\Cba>1$ in \eqref{E:best-approx} for $p\ne2$. The best $L^p$-approximation of $v$ in $T$ 
is also computable, because it boilds down to a convex minimization problem, and would result in
$\Cba=1$. This excellent property is superseded by the simplicity of \eqref{E:L2-projection-def},
which makes $\Pi_T v$ the approximation of choice.

\begin{corollary}[quasi-monotonicity of $\Pi_T$]\label{C:quasi-monotonicity}
Let $\grid,\grid_*\in\grids$ be so that $\grid\le\grid_*$, and let $T\in\grid,T_*\in\grid_*$
satisfy $T_*\subset T$. If $\Cba$ is the constant in \eqref{E:best-approx}, then
\begin{equation}\label{E:quasi-monotonicity}
\index{Constants!$\Cba$: quasi-monotonicity constant of $\Pi_T$}
\|v-\Pi_{T_*} v\|_{L^p(T_*)} \le \Cba \|v-\Pi_T v\|_{L^p(T)}
\end{equation}
for all $1\le p \le \infty$, and $\Cba=1$ for $p=2$.
\end{corollary}
\begin{proof}
Simply use \eqref{E:best-approx} to write
\[
\|v-\Pi_{T_*} v\|_{L^p(T_*)} \le \Cba  \|v-\Pi_T v\|_{L^p(T_*)} \le \Cba  \|v-\Pi_T v\|_{L^p(T)}.
\]
This is the desired bound.
\end{proof}

We are now ready to define the discontinuous $\P_{n-1}$-approximation $\wh{v}$ of $v\in L^p(\Omega)$.
Inequality \eqref{E:quasi-monotonicity} with $\Cba>1$ is fine for most
instances except Lemma \ref{l:greedy_start_grid} below. Therefore, we introduce a nonlinear
modification of the obvious choice $\wh{v}$ for $T\in\grid$,
namely $\wh{v}=\Pi_T v$. We give a recursive (and computable) definition as follows:
If $T\in\grid_0$, then $\wh{v}|_T := \Pi_T v$; if $T\in\grid$ let $\wh{v}|_{P(T)}\in \P_{n-1}$
be the approximation of $v$ in the parent element $P(T)$ of $T$, and set
\begin{equation}\label{E:wh-v}
\wh{v}|_T :=
\begin{cases}
\Pi_T v \quad &\textrm{if } \|v-\Pi_Tv\|_{L^p(T)} \le \|v-\wh{v}|_{P(T)}\|_{L^p(T)},
\\
\wh{v}|_{P(T)}  \quad &\textrm{if } \|v-\Pi_Tv\|_{L^p(T)} > \|v-\wh{v}|_{P(T)}\|_{L^p(T)}.
\end{cases}
\end{equation}
We then define
\begin{equation}\label{E:osc-v}
\index{Error Estimators!$\osc_\grid(v,T)_p$: generic surrogate for data error}
\osc_\grid(v,T)_p:= \|v-\wh{v}\|_{L^p(T)} \quad\forall T\in\grid.
\end{equation}
Since the chain of elements emanating from $\grid_0$ and culminating with $T$ is unique,
the notion $\osc_\grid(v,T)_p$ is well defined and independent of $\grid$. The following result is an immediate consequence of \eqref{E:wh-v}.

\begin{lemma}[monotonicity of oscillation]\label{L:monotonicity-osc}
For all $1\le p \le \infty$, $\grid,\grid_*\in\grids$ with $\grid\le\grid_*$, and
$T_*\in\grid_*, T\in\grid$ so that $T_*\subset T$, there holds
\begin{equation}\label{E:monotoncicity-osc}
\osc_{\grid_*}(v,T_*)_p \le \osc_\grid (v,T)_p.
\end{equation}
\end{lemma}

Consequently, for any $n\ge1$ and $T\in\grid$, let $\wh{\bA} \in [\mathbb{S}^{n-1,-1}_\grid]^{d\times d}$, $\wh{c}\in\mathbb{S}^{n-1,-1}_\grid$
be defined locally via \eqref{E:wh-v}, and let 
the surrogate element error indicators of $(\bA,c)$ be given by
\begin{equation}\label{E:zeta-Ac-loc}
\index{Error Estimators!$\osc_\grid(\bA,T)_r$: local surrogate for the diffusion coefficient approximation error}
\index{Error Estimators!$\osc_\grid(c,T)_q$: local surrogate for the reaction coefficient approximation error}
\osc_{\mesh}(\bA, T)_r := 
\Vert \bA - \hbA\Vert_{L^r(T)} \,,
\quad
\osc_{\mesh}(c, T)_q  := 
\|  c-\wh{c} \|_{L^q(T)}\,,
\end{equation}
for some $2 \le r \le \infty$ and $\frac{d}{2} < q \le \infty$ according to \eqref{E:perturbation}
for $s=0$. The simplest choice $q=2$ yields $\wh{c}_T=\Pi_T c$ in \eqref{E:wh-v}, but
requires the restriction $d < 4$, which is fine in practice.

For $n=1$ the situation is a bit special on two counts. First, $\Pi_T v$
reduces to meanvalues of $v$, namely
\begin{equation}\label{E:meanvalues}
\Pi_T \bA := \frac{1}{|T|} \int_T \bA,
\quad
\Pi_T c := \frac{1}{|T|} \int_T c,
\quad\forall T\in \grid .
\end{equation}
for $\bA \in M(\alpha_1,\alpha_2), c\in R(c_1,c_2)$ defined in \eqref{E:structural-assumption}.
Hence, $\hbA \in M(\halpha_1,\halpha_2)$ with $\wh\alpha_1=\alpha_1$, $\wh\alpha_2 = \alpha_2$ and $\wh c \in R(\hc_1,\hc_2)$ with $\wh{c}_1 = c_1, \wh{c}_2 = c_2$, i.e., the $L^2$-projections \eqref{E:meanvalues} on piecewise constants over $\grid$ as well as
$\wh{\bA}$ and $\wh{c}$ satisfy the side conditions
in \eqref{E:structural-assumption-wh} without changing the original range of parameters.
In addition, instead of \eqref{E:zeta-Ac-loc}, we can exploit {\it superconvergence}
in $W^{-1}_q(\Omega)$ with $q>\frac{d}{2-s}=d$ in \eqref{E:perturbation}. In fact, we utilize
the orthogonality of $\Pi_T$ in conjunction with \eqref{E:best-approx}
and \eqref{E:quasi-interp-error-bis} to obtain for an arbitrary function $w\in W^1_{q^*}(\Omega)$ 
and $q^*=\frac{q}{q-1}$
\begin{align*}
\int_T (c-\Pi_T c ) w = \int_T (c-\Pi_T c ) (w - \Pi_T w)
\lesssim h_T^{t} \| c-\Pi_T c \|_{L^r(T)} |w|_{W^1_{q^*}(\omega_T)}
\end{align*}
where $t=1-\frac{d}{q^*}+\frac{d}{r^*}=1+\frac{d}{q}-\frac{d}{r}>0$ and $r^*=\frac{r}{r-1}$. 
We consider two cases: $r=2,\infty$. If $r=2$ and $s=1$, then $q>d$ results in $0<t<2-\frac{d}{2}$ 
and entails the restriction $d<4$. This 
implies $\|c-\wh{c}\|_{W^{-1}_q(\Omega)} \lesssim \osc_\grid (c)_2$ where
\begin{equation}\label{E:c_super_conv}
\index{Error Estimators!$\osc_\grid(c,T)_2$: super-convergent local surrogate for the reaction coefficient approximation error}
\osc_{\grid} (c,T)_2 := h_T^t\|c-\Pi_T c\|_{L^2(T)}.
\end{equation}
If $r=\infty$ and $s=1$, then $q=\infty$ yields $t=1$ and 
$\|c-\wh{c}\|_{W^{-1}_q(\Omega)} \lesssim \osc_\grid (c)_\infty$ where
\begin{equation}\label{E:c_super_conv-infty}
\index{Error Estimators!$\osc_\grid(c,T)_\infty$: super-convergent local surrogate for the reaction coefficient approximation error}
\osc_{\grid} (c,T)_\infty := h_T\|c-\Pi_T c\|_{L^\infty(T)}.
\end{equation}

\paragraph{Approximating the load.} Dealing with $f\in H^{-1}(\Omega)$ is trickier for several reasons.
First the norm in $H^{-1}(\Omega)$ is nonlocal, so its localization is non-obvious. We recall
the definition \eqref{mod-res-est} of local oscillation $\osc_{\grid}(f,T)_{-1}$ for $T\in\grid$ and Corollary
\ref{C:local-near-best-approx-of-P} (local near-best approximation) to deduce
\begin{equation}\label{E:zeta-star}
\index{Error Estimators!$\osc_\grid(f,T)$, $\osc_\grid(f,T)_{-1}$: local oscillation for the load function}
\osc_{\mesh}(f,T)_{-1} := 
\Vert f - P_{\grid} f\Vert_{H^{-1}(\omega_T)} \le C_{\textrm{lStb}} \inf_{\chi\in\F_{\grid_{\omega_T}}} \| f - \chi\|_{H^{-1}(\omega_T)}
\end{equation}
where $C_{\textrm{lStb}}$ is the constant in Lemma \ref{L:local-stability-of-P} (local $H^{-1}$-stability);
equivalently, $\osc_\grid(f,T)_{-1}$ delivers a near best approximation of $f$ in $H^{-1}(\omega_T)$.
The second issue at stake is that without further assumptions on $f$, it is not possible to
evaluate or bound the left hand side of \eqref{E:zeta-star}. In Section \ref{S:data-approx} we
will consider several classes of loads amenable to computation and yet relevant in practice.

A popular variant of this approach for $f\in L^2(\Omega)$ replaces $\chi$ in \eqref{E:zeta-star}
by the $L^2$-projection $\Pi_\grid$ onto discontinuous piecewise polynomials of degree $n-1$,
and sets $\wh{f} = \Pi_\grid f$. This leads to the standard local weighted $L^2$-element error indicator
\looseness=-1
\begin{equation}\label{E:weighted-L^2}
\wt{\osc}_{\mesh}(f,T)_{-1} := 
h_T \Vert f - \wh f\Vert_{L^2(T)}
\quad\forall T\in\grid.
\end{equation}

\paragraph{Data error estimators.} They are the following quantities for the coefficients $(\bA,c)$
\looseness=-1
\begin{equation}\label{E:zeta-Ac-global}
\index{Error Estimators!$\osc_{\mesh}(\bA)_r$: oscillation for the diffusion coefficient}
\index{Error Estimators!$\osc_{\mesh}(c)_q$: oscillation for the reaction coefficient}
\begin{aligned}
\osc_{\mesh}(\bA)_r &:= \bigg(\sum_{T\in\grid} \osc_{\mesh}(\bA,T)_r^r\bigg)^{\frac{1}{r}},
\\
\osc_{\mesh}(c)_q &:= \bigg(\sum_{T\in\grid} \osc_{\mesh}(c,T)_q^q\bigg)^{\frac{1}{q}},
\end{aligned}
\end{equation}
which accumulate in $\ell^r$ and $\ell^q$ for $2\le r \le \infty$ and $\frac{d}{2}<q$; recall
that $q=2$ is an admissible choice provided $d<4$. In contrast, the global error estimator for $f$
\begin{equation}\label{E:zeta-f-global}
\index{Error Estimators!$\osc_{\mesh}(f)$, $\osc_{\mesh}(f)_{-1}$: oscillation for the load function}
\osc_{\mesh}(f)_{-1} := \bigg(\sum_{T\in\grid} \osc_{\mesh}(f,T)_{-1}^2 \bigg)^{\frac{1}{2}}
\end{equation}
accumulates in $\ell^2$. The total data error estimator satisfies \eqref{E:almost-best} and reads
\begin{equation}\label{E:total-data-est}
\index{Error Estimators!$\osc_{\grid}(\data)$: total data error estimator}
\osc_{\mesh}(\data) := 
\osc_{\mesh}(\bA)_r +
\osc_{\mesh}(c)_q +
\osc_{\mesh}(f)_{-1} \,.
\end{equation}

\paragraph{The module {\rm \DATA}.}
%-------------------------------------------------------------------------------------
%
This module reduces the oscillation of data $\data = (\bA,c,f)$ sequentially.
It consists of a linear approximation followed by a nonlinear correction.

Given a coefficient $v = \bA,c$, a mesh $\grid\in\grids$, a tolerance $\tau$, an
accumulation index $1\le p \le \infty$,
and a number of bisections $b\ge1$ per marked element, the call
\[
[\wt\grid,\wt{v}] = \GREEDY \big(v,\grid,\tau,p,b\big)
\]
returns a conforming refinement $\wt\grid$ of $\grid$ and a piecewise polynomial approximation
$\wt{v}$ of $v$ over $\wt\grid$ such that the oscillation computed with $v-\wt{v}$ satisfies
\[
\osc_{\wt\grid} (v)_p \le  \tau.
\]
For the load function $f$, since the computation of $\osc_\grid(f)_{-1}$ is impossible without
further assumptions on $f$, we will consider three surrogate estimators 
$\wt \osc_{\grid}(f)_{-1}$ in Section \ref{S:right-hand-side} that also accumulate in $\ell^p$ such that, for all $\grid \in \grids$,
$$
\index{Constants!$\Cdata$: \DATA approximation constant}
\osc_\grid(f)_{-1} \leq  \Cdata \wt \osc_{\grid}(f)_{-1},
$$
where $\Cdata \geq 1$.
$\GREEDY$ applied to the surrogate estimator constructs $\wt \grid \geq \grid$ satisfying
\begin{equation}\label{e:data_osc_f_red}
\wt \osc_{\wt \grid}(f)_{-1} \leq  \tau \quad  \Rightarrow \quad  \osc_{\wt \grid}(f)_{-1} \leq  \Cdata \tau.
\end{equation}
In all cases, the routine $\GREEDY$ is similar to that in Algorithm \ref{algo:greedy} (greedy algorithm)
with several important
distinctions: it accumulates the local error indicators in the $\ell^p$-norm and starts from any
mesh $\grid \ge \grid_0$ to save computational work.

Finally, the structure of the module $\DATA$ is as follows:
it concatenates $\GREEDY$  with
$\APPLYCONSTRA$ and $\APPLYCONSTRC$ in order to satisfy Assumption
\ref{A:prop-data} (properties of $\DATA$).
The routine $\GREEDY$ deals with pure approximation without constraints: called with tolerance $\tau/3$, it sequentially reduces 
the oscillation for $\bA, c, f$ with the most recent updated mesh to reduce their errors so that
\[
\osc_{\wh\grid} (\bA)_r \le   \tau/3, \quad \osc_{\wh\grid} (c)_q \le \tau/3, \quad \wt \osc_{\wh\grid} (f)_{-1} \leq  \tau/3
\]
on a conforming refinement $\wh\grid \geq \grid$. This
is discussed in detail in Section \ref{S:greedy}.

From \eqref{e:data_osc_f_red}, we get $\osc_{\wh \grid}(f)_{-1} \leq C_{\textrm{data}} \tau/3$. On the other hand, the resulting coefficients $(\wt{\bA},\wt{c})$ most likely do not satisfy the constraints
\eqref{E:structural-assumption-wh} for $n>1$. This requires a further nonlinear correction
\[
[\wh{\bA}] = \APPLYCONSTRA (\wh{\grid},\wt{\bA}),
\qquad
[\wh{c}] = \APPLYCONSTRC (\wh{\grid},\wt{c}),
\]
that enforces \eqref{E:structural-assumption-wh} on the same grid $\wh{\grid}$ without compromising the
accuracy gain produced by $\GREEDY$: there exists a constant $\ge 1$, still denoted $\Cdata$ for simplicity, such that
\[
\index{Constants!$\Cdata$: \DATA approximation constant}
\osc_{\wh\grid} ( \hbA )_r \leq \Cdata  \tau/ 3, \quad \osc_{\wh\grid} (\wh{c} )_q \le \Cdata \tau /3
\quad\Rightarrow\quad
\osc_{\wh\grid} (\data) \le \Cdata \tau.
\]
For instance, for a fixed parameter $L\ge2$, we get $\halpha_1 = \frac 1 2 \alpha_1$ and
$\halpha_2 = (1+4L)\frac {\alpha_2} 2$ for the parameters in \eqref{E:structural-assumption-wh}.
We give details in Sections \ref{ss:positivity}, \ref{S:right-hand-side}, and \ref{S:data}.

The optimality properties of $\DATA$ hinge on the performance of $\GREEDY$ and
the regularity of $\data$. Since this is not necessary for the present convergence assessment,
we discuss it later in Section \ref{S:data-approx}.

%-------------------------------------------------------------------------------
\subsubsection{Computational cost of {\rm \GALERKIN}}\label{S:cost-galerkin}
%-------------------------------------------------------------------------------
%
The output pair $(\widehat{\mesh},\widehat\data)$ of $\DATA$
is next taken by $\GALERKIN$, the one-step AFEM of Algorithm \ref{A:GALERKIN} in Section \ref{S:GALERKIN},
to run an inner loop of the form \eqref{E:adaptive-loop}
with fixed discrete data $\widehat\data$ and initial mesh
$\widehat{\mesh}$. The call \eqref{E:galerkin-eps} of $\GALERKIN$ stops as soon as the error
tolerance $\varepsilon$
is reached, which takes a finite number of iterations because $\GALERKIN$
is a contraction between
consecutive iterates, and creates the next mesh-solution pair $(\mesh,\uG)$.
It is worth noticing that, in the absence of
this stopping test, the Galerkin solution $\uG$ would converge to the solution
$\hu=u(\wh\data)$ of \eqref{E:perturbed-weak-form}, which is not the desired solution $u=u(\data)$ of
\eqref{strong-form}.

We stress that, in view of \eqref{E:bound-data} and \eqref{E:error-DATA},
the relative resolution of the modules $\DATA$ and $\GALERKIN$ is critical for the discrepancy between
the exact and perturbed solutions $u$ and $\hu$. This is ultimately
responsible for the performance of $\AFEMTS$ and is studied in Section
\ref{S:conv-rates-coercive}.

We now investigate the number of iterations within $\GALERKIN$, which dictate its computational
cost. We point out that at iteration $k-1\ge0$ of $\AFEMTS$, the output $(\grid_k,u_k)$ of $\GALERKIN$, and thus of
$\AFEMTS$, satisfies
\begin{equation}\label{E:output-galerkin}
\eta_k(u_k) = \eta_{\grid_k}(u_k) \le \eps_{k-1}
\quad\Rightarrow\quad
|u_k - \hu_{k-1}|_{H^1_0(\Omega)} \le C_U \eps_{k-1}
\end{equation}
according to \eqref{E:aposteriori-bounds-H1}. We recall that $\hu_{k-1}=\hu_{k-1}(\wh{\data}_{k-1})\in H^1_0(\Omega)$ is the exact solution with discrete data $\wh\data_{k-1}$, and that $\est_{\grid_k}(u_k,f)$ is defined with discrete data $\wh{\data}_{k-1}$ and satisfies  $\est_{\grid_k}(u_k,f) = \eta_{\grid_k}(u_k)$ because
data oscillation $\osc_{\grid_k}(f)_{-1} = 0$. The next iteration $k$ of $\AFEMTS$ calls $\DATA$,
which in turn refines the mesh $\grid_k$ to $\wh\grid_k$ and updates the data approximation from
$\wh\data_{k-1}$ to $\wh\data_k$ over $\wh\grid_k$. The pair $(\wh\grid_k,\wh\data_k)$ determines
the first Galerkin solution $u_{k,0}\in\V_{k,0}=\V_{\wh\grid_k}$ of $\GALERKIN$
and corresponding estimator $\eta_{k,0}(u_{k,0})$ with $\grid_{k,0}=\wh\grid_k$, which must satisfy
\begin{equation}\label{E:galerkin-active}
\eta_{k,0} (u_{k,0}) > \eps_k
\end{equation}
for $\GALERKIN$ to be executed.
The reduction of $\eta_{\grid_{k,j}}(u_{k,j})$ for $j\ge 1$ dictates the
number of iterations of $\GALERKIN$. We examine this next.

\begin{proposition}[computational cost of $\GALERKIN$]\label{P:cost-galerkin}
\AB{If the assumptions of Theorem \ref{T:contraction} are valid, then} for any $k\in\N$, the number of subiterations $J_k$ inside a call to {\rm \GALERKIN} at
iteration $k$ of {\rm \AFEMTS} is bounded independently of $k$.
\end{proposition}
\begin{proof}
The $j$-th error $e_{k,j} := |\wh{u}_k-u_{k,j}|_{H^1_0(\Omega)}$ within $\GALERKIN$ converges linearly
in view of Corollary \ref{C:linear-convergence-error} (linear convergence of error) because the
discrete data
$\wh{\data}_k$ is fixed in these inner iterations. Exploiting the lower bound
$C_L \eta_{k,j}(u_{k,j}) \le e_{k,j}$ stated in \eqref{E:aposteriori-bounds-H1}, we thus
deduce
\[
\eta_{k,j}(u_{k,j}) \le C_L^{-1} e_{k,j}\le C_L^{-1} C_* \alpha^{j-i}
e_{k,i}, \quad j\ge i \ge 0,
\]
whence $\eta_{k,j}(u_{k,j}) \le C_\# \alpha^j e_{k,0}$ with  $C_\#:=C_L^{-1} C_*$. The number of iterations of $\GALERKIN$
depends on the size of $\eta_{k,0}(u_{k,0})$ relative to $\eps_k$. We assume that
$\eta_{k,0}(u_{k,0})>\eps_k$ according to \eqref{E:galerkin-active}.
We first prove that $\eta_{k,0}(u_{k,0})\lesssim\eps_k$ and next argue that $J_k$
is bounded uniformly in $k$. We proceed in two steps.

\smallskip
\step{1} {\it Bound on $|\wh{u}_k-u_{k,0}|_{H^1_0(\Omega)}$}: Since $u_k\in\V_k\subset\V_{k,0}=\V_{\wh\grid_k}$,
and the Galerkin solution $u_{k,0}\in\V_{k,0}$ minimizes the error $\enorm{u_{k,0}-\hu_k}$ in $\V_{k,0}$,
relative to the energy norm induced by the bilinear
form $\wh{\B}$ with discrete data $\wh{\data}_k$, we deduce
\[
\enorm{u_{k,0}-\hu_k} \le \enorm{u_k-\hu_k}
\le \sqrt{C_{\wh\B}} \Big( |u_k-\hu_{k-1}|_{H^1_0(\Omega)} + |\hu_{k-1}-\hu_k|_{H^1_0(\Omega)} \Big),
\]
where the last inequality uses \eqref{E:exactB} for $\wh{\B}$.
Invoking the a posteriori upper bound \eqref{E:aposteriori-bounds} and the termination condition
of $\GALERKIN$ at step $k-1$, we obtain
\[
|u_k-\hu_{k-1}|_{H^1_0(\Omega)} \le C_U \est_{\grid_k}(u_k,f) = C_U \eta_k(u_k) \le C_U \eps_{k-1} = 2C_U \eps_k.
\]
On the other hand, using \eqref{E:error-DATA} with $\tau=\omega\eps_k$
and $0<\omega\le 1$, we arrive at
\[
|u-\hu_k|_{H^1_0(\Omega)} \le C_1 \eps_k,
\]
with $C_1=\omega C_D$. The triangle inequality thus yields
\[
|\hu_{k-1}-\hu_k|_{H^1_0(\Omega)} \le |u-\hu_{k-1}|_{H^1_0(\Omega)} + |u-\hu_k|_{H^1_0(\Omega)}
\le C_1(\eps_{k-1}+\eps_k)
= 3C_1 \eps_k,
\]
whence
\[
e_{k,0}= |u_{k,0}-\hu_k|_{H^1_0(\Omega)} \le \sqrt{\frac{C_{\wh{\B}}}{c_{\wh{\B}}}}\big(2C_U + 3C_1 \big) \eps_k
=: C_2 \eps_k.
\]

\step{2} {\it Bound on $J_k$}: We observe that $\GALERKIN$ stops once $\eta_{k,j}(u_{k,j}) \le\eps_k$.
Since the smallest such $j$ is $J_k$, we see that 
\[
\eps_k < \eta_{k,J_k-1}(u_{k,J_k-1}) \le C_\# \alpha^{J_k-1} e_{k,0} \le  C_\#  C_2 \eps_k \alpha^{J_k-1}.
\]
This implies the asserted bound $J_k \le 1 + \frac{\log (C_\# C_2)}{\log \alpha^{-1}}$ uniform in $k$.
\end{proof}

%-------------------------------------------------------------------------------
\subsubsection{Realization of $\AFEMTS$}\label{S:realization-AFEM}
%-------------------------------------------------------------------------------
%
We now make the two step \AFEM algorithm precise.

\begin{algo}[\AFEMTS]\label{algo:AFEMTS}
\index{Algorithms!\AFEMTS: two step AFEM successively approximating the data and the Galerkin solution with approximate data}
Given an initial tolerance $\eps_0>0$, a target tolerance $\tol$ and initial mesh $\mesh_0$,
as well as a safety parameter $\omega \in (0,1]$, {\rm AFEM} consists of the two-step algorithm:

\medskip
\begin{algotab}
  \>  $[\mesh, u_\mesh]=\AFEMTS \, (\mesh_0, \varepsilon_0, \omega, \tol)$
  \\
  \>  \> $\text{set } k= 0 \text{ and do }$ 
  \\
  \>  \> \> $[\widehat\mesh_{k},\widehat{\data}_{k}]=\DATA \, (\mesh_k, \data, \omega \, \varepsilon_k) $ 
  \\
  \>  \> \> $[\mesh_{k+1},u_{k+1}]=\GALERKIN \, (\widehat{\mesh}_{k},\widehat{\data}_{k},\varepsilon_k)$ 
  \\
  \>  \> \> $\varepsilon_{k+1}=\tfrac12 {\varepsilon_k}$
  \\
  \>  \> \> $k \leftarrow k+1$
  \\
  \>  \> $\text{while } \eps_{k-1}>\tol$
  \\
  \>  \> $\text{return} \, \mesh_k, u_k$
\end{algotab}
\end{algo}

\begin{proposition}[convergence of $\AFEMTS$]\label{P:convergence-AFEM}
For each $k\ge0$ the modules {\rm \DATA} and {\rm \GALERKIN} converge in a finite number of iterations, the
latter independent of $k$.
Moreover, there exists a constant $C_*$ depending on $\mesh_0, \Omega, d, n$, the Lebesgue exponents
$r,q$ in $D(\Omega)$, the parameters $\alpha_1, \alpha_2, c_1, c_2$ in \eqref{d:DiffusionSpace}
and \eqref{d:ReactionSpace}, and the shape regularity constant of $\grids$,
such that the output of the $(k+1)$-th iteration
$[\mesh_{k+1},u_{k+1}]=\GALERKIN \, (\widehat{\mesh}_{k},\widehat{\data}_{k},\varepsilon_k)$ satisfies
$|u -u_{k+1}|_{H^1_0(\Omega)} \leq C_* \varepsilon_k$ for all $k \geq 0$.    
Therefore, {\rm \AFEMTS} stops after $K<2+\frac{\log\frac{\eps_0}{\tol}}{\log2}$ iterations and delivers
\begin{equation*}
  | u -u_{K} |_{H^1_0(\Omega)} \leq C_* \tol.
\end{equation*}   
\end{proposition}
\begin{proof}
In view of Assumption \ref{A:prop-data} (properties of $\DATA$) the module $\DATA$ iterates a finite number
of steps to reach tolerance $\tau=\omega\varepsilon_k$ for every $k\ge0$.
Moreover, the number of iterations of $\GALERKIN$ is independent of $k$
due to Proposition \ref{P:cost-galerkin} (computational cost of $\GALERKIN$), whence we deduce
that each loop of $\AFEMTS$ requires finite
iterations. Thus, the output $u_{k+1}$ of the $(k+1)$-th loop satisfies
\[
|u -u_{k+1}|_{H^1_0(\Omega)} \leq | u - \widehat{u}_k |_{H^1_0(\Omega)} + | \widehat{u}_k -u_{k+1} |_{H^1_0(\Omega)}
\le \big( \omega C_D + C_U \big) \varepsilon_k = C_* \eps_k,
\]
according to \eqref{E:error-DATA} with $\tau\le\omega\varepsilon_k$ and \eqref{E:output-galerkin} for all
$k \geq 0$. Finally, $\AFEMTS$ terminates after $K$
loops, where $K$ satisfies $\frac{1}{2}\tol < \varepsilon_{K-1}  \leq \tol$, and the
asserted estimate holds.
\end{proof}

This elementary proof gives no insight whether the $H^1_0$-error decays optimally in
terms of degrees of freedom.
We assess this fundamental question in Sections \ref{S:conv-rates-coercive} and \ref{S:data-approx},
but investigate it computationally in Section \ref{S:AFEM-experiments}.

A two-step algorithm similar to $\AFEMTS$ was first proposed in \cite{Stevenson:08},
and further explored in \cite{BonitoDeVoreNochetto:2013,CohenDeVoreNochetto:2012}. Note that other quantities, such as the number of degrees of freedom, could be employed to stop $\AFEMTS$ instead.   
It is also worth realizing that the structure of the algorithm is independent of the size of tolerance {\tt tol}. In this vein, a user could take $\eps_0 = {\tt tol}$, provided {\tt tol} is affordable by the computational resources at hand. With such a choice, the modules $\DATA$ and $\GALERKIN$ run only once, in sequence: data are approximated to the desired accuracy in one shot, then fed to the PDE solver which produces the approximate solution.
%   
% It is worth realizing that such a structure is independent of the size of tolerance {\tt tol} and that perhaps other quantities, such as the number of degrees of freedom, could be employed to stop $\AFEMTS$ instead. 
% In this vein, a user could take $\eps_0 = {\tt tol}$ provided {\tt tol} is affordable by the computational resources at hand. In this way, the modules $\DATA$ and $\GALERKIN$ would run only once without interacting with each other. 
Since the quasi-optimality theory in Section \ref{S:conv-rates-coercive} would also hold for this choice of $\eps_0$, one might wonder why not using this simpler strategy. We stress that iterating over $\eps_k$ has the following advantages:
\begin{itemize}%[$\bullet$]
\item
{\it Restarts:} Dynamical shrinking of {\tt tol}, for instance to account for the user decision to improve the accuracy, does not entail a restart of $\AFEMTS$ but rather a continuation from the previous computed solution. In this sense, the resulting iteration would be similar to the proposed structure of $\AFEMTS$.

\medskip
\item
{\it Computational resources:} $\AFEMTS$ allows for ``balanced investment" of computational resources between the modules $\DATA$ and $\GALERKIN$. If the stopping criterion, either accuracy or number of degrees of freedom, is unrealistic for the problem at hand, $\AFEMTS$ would still produce a discrete solution with equilibrated data and solution errors.

\medskip
\item
{\it Nonlinear problems:} The interleaving approach of $\AFEMTS$ appears to be better suited for treating nonlinear problems for which data $\data$ may depend on the solution. Therefore a call to $\GALERKIN$, and corresponding solution update, must precede a call to $\DATA$.

\medskip
\item
{\it Iterative solvers:} If an efficient iterative solver is adopted within $\SOLVE$, then the previous discrete solution of $\GALERKIN$ could be taken as initial iterate, thereby making $\SOLVE$ fast because $\frac{\epsilon_{k+1}}{\epsilon_k}= \frac12$. If instead one computes with $\DATA$ alone until the fixed tolerance {\tt tol} is reached, then $\GALERKIN$ would  work directly on fine meshes, which are not adapted to the geometric domain singularities, and without good initial guess. This would lead to fewer but heavier iterations of $\GALERKIN$, which is detrimental from a linear algebra perspective.
\end{itemize}

%%%%%%%%%%%%%%%%%%%%%%%%%%%%%%%%%%%%%%%%%%%%%%%%%%%%%%%%%%%%%%%%%%%%%%%%%%%%%%%%%%%%%%%%%%%%
\subsubsection{Computational assessment of {\rm \AFEMTS}}\label{S:AFEM-experiments}
%%%%%%%%%%%%%%%%%%%%%%%%%%%%%%%%%%%%%%%%%%%%%%%%%%%%%%%%%%%%%%%%%%%%%%%%%%%%%%%%%%%%%%%%%%%%
%
In this section we explore computationally the relative performance of $\GALERKIN$ and $\DATA$,
for the two-step $\AFEM$, and elucidate the behavior of data and coefficient oscillations within
$\DATA$. Our observations motivate the rigorous study of Section \ref{S:conv-rates-coercive}, which
provides theoretical support to our experiments. The numerical computations are made with the help of
\cite{FunkenPraetorius:11}.

We consider problem \eqref{strong-form} in the L-shaped domain
$\Omega = (-1, 1)^2\setminus([0, 1]\times [-1, 0])$,  with diffusion term $\bA = a \vec{I}$, where
\begin{align*}
%a(x,y) = 1 + \exp\bigl( -50 |(x + 0.5, y + 0.5)|^2 \bigr) + \exp\bigl( -50|(x+0.5,y-0.5)|^2\bigr)
%\\
a(x, y) = 1 &+ \exp\bigl(-50 ((x + 0.5)^2 + (y+0.5)^2) \bigr)
+ \exp\bigl(-50 ((x + 0.5)^2 + (y-0.5)^2) \bigr) \,,
\end{align*}
and reaction term
\[
\begin{aligned}
c(x, y) &= 1 + \exp\bigl(-50 ((x + 0.5)^2 + y^2) \bigr)
+ \exp\bigl(-50 (x^2 + (y-0.5)^2) \bigr) \,;
\end{aligned}
\]
note that the Gaussians in the definition of $a$ and $c$ have the same intensity but are located
in different places within $\Omega$. The load term $f$ and the Dirichlet
boundary conditions are chosen in accordance with the analytical solution
\[
u(x, y) = r^{\frac{2}{3}} \sin \bigl(  2\alpha/3\bigr) + 
\exp\bigl(-1000 ((x - 0.5)^2 + (y-0.5)^2) \bigr) \,,
\]
where $(r, \alpha)$ are the polar coordinates around the origin.
Notice that the exact solution $u$ is singular at the reentrant corner: it belongs to the Sobolev spaces
$H(\Omega)^{\frac{5}{3}-\eps}$ with $\eps > 0$ and $W^2_p(\Omega)$ with $p>1$. It also exhibits a rapid
transition of order $10^{-3/2}$ around the point $(0.5,0.5)$ due to the presence of a very narrow Gaussian.
The Gaussians are meant to test the performance of the module $\DATA$,
while in addition the corner singularity of the solution tests the execution of the module $\GALERKIN$.
 
\begin{figure}[ht]
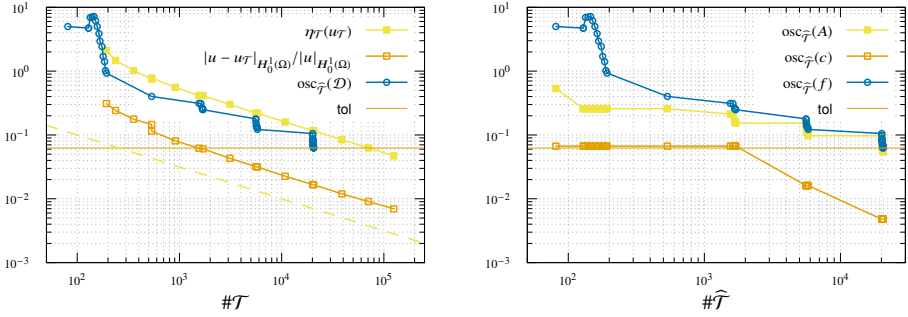

  \centering
\scalebox{0.7}{\input{GNUPLOT/AFEM-Convergence}
\input{GNUPLOT/DATA-Convergence}}
  \caption{{Left: 
estimator $\eta_\grid(u_\grid)$,
data error $\osc_{\wh\mesh}(\data)$, and relative $H^1$-error 
obtained with the algorithm $\AFEM$ performing $b=3$ bisections per marked element. 
The optimal decay is indicated by the dashed line with slope $-0.5$.
Right: Diffusion error $\osc_{\wh\grid}(\bA)$,
reaction error $\osc_{\wh\grid}(c)$,
load error $\osc_{\wh\grid}(f)_{-1}$,
obtained with the algorithm $\AFEM$.}}
\label{F:avem}
\end{figure}

\begin{figure}
    \centering    \includegraphics[trim={2.25cm 8cm 2cm 6.5cm},clip,scale=0.3]{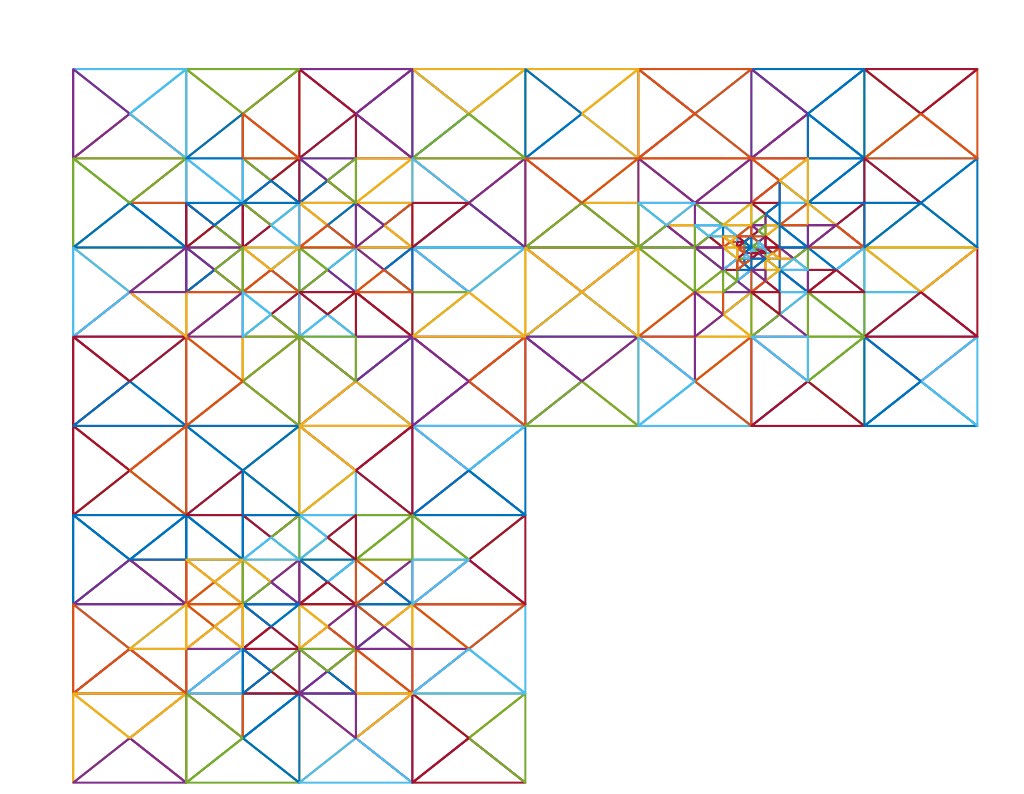}\includegraphics[trim={2.25cm 8cm 2cm 6.5cm},clip,scale=0.3]{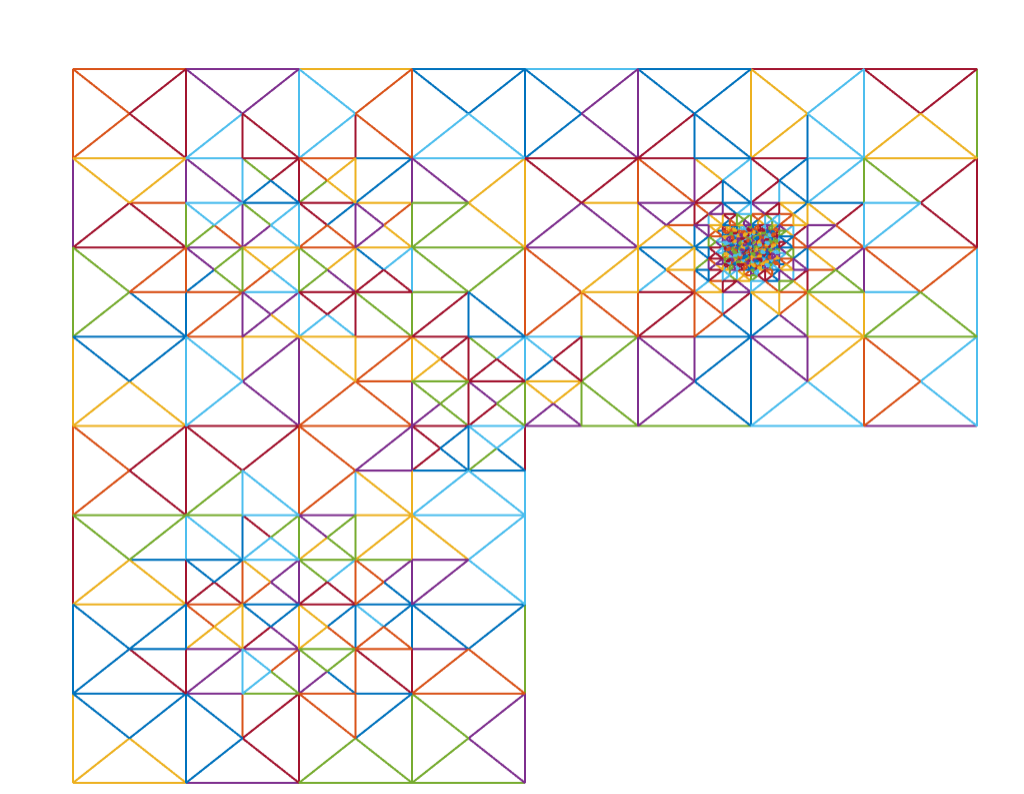}
    \caption{Resulting meshes after the first iteration of \DATA (left) and after the first iteration of \GALERKIN (right). \DATA marked 29 elements for refinement while \GALERKIN marked 308 elements. Refer to Table~\ref{t:n_marked} for more details.}
    \label{f:mesh_conv}
\end{figure}

We utilize the following parameters in the numerical test
\[
\theta = 0.5, \quad \omega = 1,
\quad \texttt{tol}=2^{-4},
\quad h_0=0.125,
\quad
\eps_0=1.
\]
Notice that the number of iterations of the algorithm $\AFEM$ is $K=\log_2(\epsilon_0/\texttt{tol}) = 4$.
We compute the relative $H^1$-error between the exact solution $u$ and
the FEM solution $u_\grid$ and notice that its decay rate is $(\#\grid)^{-1/2}$ in Fig. \ref{F:avem} 
 (left). This rate is
consistent with that of the PDE estimator $\eta_\grid(u_\grid)$ and data estimator
$\osc_\grid(\data)$. In Fig. \ref{F:avem} (right) we display the component of the data error $\osc_{\wh\grid}(\bA)$, $\osc_{\wh\grid}(c)$,
$\osc_{\wh\grid}(f)_{-1}$ defined in \eqref{E:zeta-Ac-global} and \eqref{E:zeta-f-global} with
local contributions defined in \eqref{E:zeta-Ac-loc} for $\bA$ with $r=\infty$ , in \eqref{E:c_super_conv-infty} for $c$ with $t=1$ and \eqref{E:weighted-L^2} for $f$. Recall that at each iteration $k$, \DATA circles through $\osc_{\wh\grid}(\bA)$, $\osc_{\wh\grid}(c)$, and $\osc_{\wh\grid}(f)_{-1}$ reducing each of these oscillations to 1/3 of the iteration tolerance $\varepsilon_k=2^{-k}$. The presence of the weight $h^t$ in $\osc_{\wh\grid}(c)$ considerably reduces the influence of the approximation of $c$, which is below threshold from the start and thus never generates any refinement (see Table~\ref{t:n_marked}). The local oscillation for $f$ also includes a weight vanishing as $h \to 0$ but $\osc_{\wh\grid}(f)_{-1}$ is above the desired tolerance, which would in principle generate refinements. However, since at each iteration \DATA considers $\osc_{\wh\grid}(\bA)$ first and the regions refined to reduce $\osc_{\wh\grid}(f)_{-1}$ are included in the regions needed to be refined to reduce $\eta_\grid(u_\grid)$ and $\osc_{\wh\grid}(\bA)$, the \GREEDY routine applied to $f$ does not refine any element except during the first iteration when the Galerkin error has not yet been reduced by the algorithm. Overall, the reduction of the Galerkin error is driving most of the refinements. The number of marked element to reduce the approximation errors of $\bA$, $c$, $f$, and the residual estimator are reported in Table~\ref{t:n_marked} along with those when $b=1$ refinement is used per marked element. In Figure~\ref{f:mesh_conv}, we provide the resulting meshes after the first iteration of \DATA and \GALERKIN. 

\begin{table}[htb!]
    \centering
    \begin{tabular}{|c||c|c||c|c||c|c||c|c|}
    \hhline{---------}
       \multirow{2}{*}{$k$} 
       & \multicolumn{2}{c||}{$\osc_{\wh\grid}(\bA)$} & \multicolumn{2}{c||}{$\osc_{\wh\grid}(c)$} & \multicolumn{2}{c||}{$\osc_{\wh\grid}(f)_{-1}$} & \multicolumn{2}{c|}{$\eta_\grid(u_\grid)$} \\
    \hhline{|~||--||--||--||--}
       & b=1 & b=3 & b=1 & b=3 & b=1 & b=3 & b=1 & b=3\\
    \hhline{-||--||--||--||--}
       1  &  32 & 16 &  0 & 0  & 26 & 13 & 363 & 308 \\
       2  &  16 & 16 &  0 & 0  &  0 & 0   & 1636 & 1138  \\
       3  &  120 &43 &  0 & 0  &  0 & 0   & 7447 & 4227  \\
       4  &  123 &62 &  0 & 0  &  0 & 0  &  42792 & 15268 \\
       5  &  82 & 138&  0 & 0  &  0 & 0  & 144345 & 102350  \\
    \hhline{---------}
    \end{tabular}
    \caption{Number of marked elements to reduce the data and Galerkin errors at each iterations $k=1,2,3,4$ of \AFEMTS when using $b=1$ and $b=3$ refinements per marked element. Regardless of the value used for $b$, the reduction of the Galerkin error is driving most of the refinements followed by the error in the approximation of the diffusion coefficient $\bA$. The approximation of $f$ is subordinate to the approximation of $u$ and $\bA$ arising earlier in the adaptive loop and thus does not generate any refinement except during the first iteration when the Galerkin error has not yet been tackled by the algorithm. The approximation of $c$ is below the final tolerance from the start and does not generate any refinement.}
    \label{t:n_marked}
\end{table}

%%%%%%%%%%%%%%%%%%%%%%%%%%%%%%%%%%%%%%%%%%%%%%%%%%%%%%%%%%%%%%%%%%%%%%%%%%%%%%%%%%%%%%%%%%%%
\subsection{Convergence for other boundary conditions}
\label{S:conv-other-bdry-cond}
%%%%%%%%%%%%%%%%%%%%%%%%%%%%%%%%%%%%%%%%%%%%%%%%%%%%%%%%%%%%%%%%%%%%%%%%%%%%%%%%%%%%%%%%%%%%

We consider first the variational problem \eqref{E:Robin-bilinear} with {\it Robin} boundary condition.
We approximate data $\data = (\bA,c,p,f,g)$ by piecewise polynomials 
$\wh{\data} = (\wh{\bA},\wh{c},\wh{p},\wh{f},\wh{g})$,
The only difference with respect to \eqref{E:space-discrete-data} is that the new functions
$(p,g)$ are approximated on $\partial\Omega$ by discontinuous polynomials $(\wh{p},\wh{g})$
of degree $n-1$ and $2n-1$.
The projection operator $P_\grid$ approximates $g\delta_{\partial\Omega}$ by 
$\wh{g} \delta_{\partial\Omega}=P_\grid (g\delta_{\partial\Omega})$ without component 
in the bulk because $g\delta_{\partial\Omega}$ is a line Dirac mass aligned with 
the mesh. Discrete functions $(\wh{p},\wh{g})$ must be produced by $\DATA$, subject to a sign constraint on $p$.
The approximate bilinear form $\wh{\B}$ and linear functional $\wh{\ell}$ read
\begin{equation}\label{E:B-robin}
\wh\B[w,v] := \int_\Omega \nabla v \cdot \wh{\bA} \nabla w + \wh{c} v w +\int_{\partial\Omega} \wh{p} vw,
\quad
\wh{\ell} (v) := \langle \wh{f},v \rangle + \int_{\partial\Omega} \wh{g} \, v.
\end{equation}
The a posteriori error estimates of Section \ref{S:aposteriori} extend to this pair $(\wh{\B},\wh{\ell})$.
The algorithms $\GALERKIN$, $\AFEMSW$, and $\AFEMTS$ are similar to those above and possess a
similar supporting convergence theory. The {\it Neumann} boundary condition is a 
particular case with $p=0$. We do not pursue this any further.

However, the {\it pure Neumann} boundary condition is special because of the global
compatibility condition $\wh{\ell}(1)= \langle \wh{\ell},1 \rangle =0$. In Section \ref{L:other-bc-apost} we introduce a new
projection operator $\wt{P}_\grid$, a modification of $P_\grid$, with the requisite properties
of local approximation and global compatibility $\langle \wt{P}_\grid \ell, 1\rangle =0$ provided 
$\ell\in H^1(\Omega)^*$ satisfies $\langle \ell,1 \rangle=0$. We thus set $\wt{\ell} = \wt{P}_\grid \ell$
to solve the Galerkin problems and use $\wt{P}_z$ in the local indicators. We do not explore
this matter further.

For a {\it non-homogeneous Dirichlet} boundary data $g \in H^{1/2}(\partial\Omega)$,
$\DATA$ must produce a continuous piecewise polynomial approximation $\wh{g}$ of degree $n$, thereby consistent
with the Galerkin solution $u_\grid$. The Dirichlet oscillation $\osc_\grid(g)_{1/2}$ is defined in
\eqref{E:Dirichlet-osc} and is locally computable. Data oscillation now becomes
\[
\osc_\grid(\ell) = \osc_\grid(f)_{-1} + \osc_\grid(g)_{1/2}
\]
and added to the PDE estimator $\eta_\grid(u_\grid)$ for $g=0$ gives a full estimator equivalent to the error,
according to Theorem \ref{T:Dirichlet} (estimators for general Dirichlet conditions). With these minor
modifications, the convergence theory for $\GALERKIN, \AFEMSW$, and $\AFEMTS$ extends to this 
case. We do not provide any further details.

%%%%%%%%%%%%%%%%%%%%%%%%%%%%%%%%%%%%%%%%%%%%%%%%%%%%%%%%%%%%%%%%%%%%%%%%%%%%%%%%%%%%%%%%%%%%
\subsection{Convergence for alternative estimators}
\label{S:conv-non-residual-est}
%%%%%%%%%%%%%%%%%%%%%%%%%%%%%%%%%%%%%%%%%%%%%%%%%%%%%%%%%%%%%%%%%%%%%%%%%%%%%%%%%%%%%%%%%%%%

We have so far developed a convergence theory for the residual estimator $\est_\grid(u_\grid,f)$.
The purpose of this section is to extend this theory to the three alternative estimators
discussed in Section \ref{S:other_estimators}, namely
\begin{itemize}
\item
$\est^\mathrm{lpb}_\grid(u_\grid,f)^2 = \eta^\mathrm{lpb}_\grid(u_\grid)^2 + \osc_\grid(f)_{-1}^2$:
estimator based on {\it local  problems};

\item
$\est^\mathrm{hier}_\grid(u_\grid,f)^2 = \eta^\mathrm{hier}_\grid(u_\grid)^2 + \osc_\grid(f)_{-1}^2$:
{\it hierarchical} estimator

\item
$\est^\mathrm{feq}_\grid(u_\grid,f)^2 = \eta^\mathrm{feq}_\grid(u_\grid)^2 + \osc_\grid(f)_{-1}^2: $
estimator based on {\it flux equilibration}.
\end{itemize}
They are all computed on stars $\omega_z$ with $z\in\vertices$ and possess a similar structure.
The first term is the PDE estimator, from now on called $\zeta_\grid(u_\grid)$ to refer to any of
them, and is locally equivalent to the discrete residual $P_\grid R_\grid$
\begin{equation}\label{E:star-equiv-res}
\zeta_\grid (u_\grid,z) \approx \| P_\grid R_\grid \|_{H^{-1}(\omega_z)}
\quad \forall \, z\in\vertices;
\end{equation}
see Theorems \ref{T:lpb-est}, \ref{T:hier-est}, and \ref{T:feq-est}. In fact, they are all
different mechanisms to extract information from $P_\grid R_\grid$. Since the vertex-indexed
residual PDE indicator $\eta_\grid(u_\grid,z):=\eta^\mathrm{res}_\grid(u_\grid,z)$, defined in \eqref{vertex-mod-rest-est;pde-ind}, is also proven to be equivalent to $\| P_\grid R_\grid \|_{H^{-1}(\omega_z)}$
in Theorem \ref{T:vertex-mod-res-est} (vertex-indexed modified residual estimator), we deduce the existence of two equivalence constants $\CLequiv \le \CUequiv$ such that
\begin{equation}\label{E:star-equiv-est}
\index{Constants!$(\CLequiv,\CUequiv)$: lower and upper estimators equivalence constant}
\CLequiv \eta_\grid (u_\grid,z) \le \zeta_\grid (u_\grid,z) \le \CUequiv \eta_\grid (u_\grid,z)
\quad \forall \, z\in\vertices.
\end{equation}

Following \cite{KreuzerSiebert:11}, we will exploit this property to prove convergence of
$\AFEM$ driven by $\zeta_\grid(u_\grid)$. An obstruction for a direct convergence theory is that
our preceding results rely heavily on Lemma \ref{L:strict-est} (reduction property of the estimator),
which is not necessarily valid for any of the alternative estimators. We refer to
\cite{CasconNochetto:12} who present a direct approach based on the local lower
bound for discrete solutions of Theorem \ref{T:lubd-corr} (local lower bound for corrections).
The latter is guaranteed by Definition \ref{D:interior-node} (interior vertex property) for operators
with coefficients $\bA$ piecewise constant and $c=0$, and any polynomial degree $n\ge1$, but we do
not know its validity for more general coefficients $(\bA,c)$.

The key for convergence is imposing a D\"orfler marking. 
We say that a set of vertices $\marked_{\cal V}$ satisfies a D\"orfler property with parameter $\theta\le1$
if
\begin{equation}\label{E:doerfler-property}
\zeta_\grid(u_\grid,\marked_{\cal V})^2 := \sum_{z\in\marked_{\cal V}} \zeta_\grid(u_\grid, z)^2
\ge \theta^2 \sum_{z\in\vertices} \zeta_\grid(u_\grid, z)^2 =: \zeta_\grid(u_\grid)^2.
\end{equation}
Let $\marked$ be the collection of elements contained in the stars $\omega_z$ with $z\in\marked_{\cal V}$. Then $\MARK$ marks  all elements in $\marked$, and $\REFINE$ bisects them $b\ge1$ times. This gives rise to a star-driven $\GALERKIN$ procedure.

\begin{lemma}[D\"orfler property]\label{L:doerfler-property}
If the set of vertices $\marked_{\cal V}$ satisfies a D\"orfler property with parameter $\theta$
for $\zeta_\grid(u_\grid)$, then $\marked$ satisfies a D\"orfler property with parameter $\overline{\theta} = \frac{\CLequiv}{\CUequiv}\theta$ for $\eta_\grid(u_\grid)$.
\end{lemma}
\begin{proof}
Simply use \eqref{E:star-equiv-est} to derive \eqref{E:doerfler-property} for $\eta_\grid(u_\grid)$
with parameter $\overline{\theta}$.
\end{proof}

Hence, star-driven procedures for $\zeta_\grid(u_\grid)$ lead to the corresponding counterparts for $\eta_\grid(u_\grid)$. It turns out that algorithms $\GALERKIN$, $\AFEMSW$, and $\AFEMTS$ can be reformulated for vertex-indexed indicators $\{\eta_\grid(u_\grid,z)\}_{z\in\vertices}$
as  defined in \eqref{vertex-mod-rest-est;pde-ind}, without changing their essential properties. We may thus wonder about them driven by
$\{\zeta_\grid(u_\grid,z)\}_{z\in\vertices}$ instead. Since 
these algorithms hinge on the D\"orfler property \eqref{E:doerfler-property},
Lemma \ref{L:doerfler-property} gives rise to similar convergence properties for
$\zeta_\grid(u_\grid)$-driven algorithms provided \eqref{E:doerfler-property} is enforced. 
We state this next without proof.

\begin{corollary}[convergence of $\GALERKIN$]\label{C:galerkin-zeta}
If the coefficients $(\bA,c,f) \in \AB{\D_{\grid}}$, then there exist
$0<\overline{\alpha}<1$ and $C_*,C_\#>0$ such that the solution-estimator
pairs $(u_j,\zeta_j(u_j))$ of {\rm \GALERKIN} converge linearly, namely for all $k\ge j \ge0$
\[
|u-u_k|_{H^1_0(\Omega)} \le C_* \, \overline{\alpha}^{k-j} |u-u_j|_{H^1_0(\Omega)},
\quad
\zeta_k(u_k) \le C_\# \, \overline{\alpha}^{k-j} \zeta_j(u_j).
\]
\end{corollary}

\begin{corollary}[convergence of $\AFEMSW$]\label{C:afem-ws-zeta}
If the coefficients $(\bA,c)$ are discrete and $f\in H^{-1}(\Omega)$,
then for \AB{$0<\omega \leq \omega_0$, $\xi\le\frac{1}{2}$ as in Theorem \ref{T:contraction-sw}},
there exist $0<\overline{\alpha}<1$ and $C_*,C_\#>0$ such that the solution - estimator
pairs $\big(u_j,\est_j(u_j,f)\big)$ of {\rm \AFEMSW}, where $\est_j(u_j,f)^2 = \zeta_j(u_j)^2 + \osc_j(f)_{-1}^2$, converge linearly:  for all $k\ge j \ge0$
\[
|u-u_k|_{H^1_0(\Omega)} \le C_* \, \overline{\alpha}^{k-j} |u-u_j|_{H^1_0(\Omega)},
\quad
\est_k(u_k,f) \le C_\# \, \overline{\alpha}^{k-j} \est_j(u_j,f).
\]
\end{corollary}

Both $\GALERKIN$ and $\AFEMSW$ converge under restrictions on data $\data = (\bA,c,f)$.
For arbitrary data $\data$, $\AFEMTS$ concatenates $\GALERKIN$ and $\DATA$, the later being unrelated
to $\zeta_\grid(u_\grid)$. Therefore, Corollary \ref{C:galerkin-zeta} and
Proposition \ref{P:cost-galerkin} (computational cost of $\GALERKIN$) yield the following
extension of Proposition \ref{P:convergence-AFEM} (convergence of $\AFEMTS$).

\begin{corollary}[convergence of $\AFEMTS$]\label{C:afem-ts-zeta}
The algorithm {\rm \AFEMTS} driven by $\zeta_\grid(u_\grid)$ stops after 
$K<2+\frac{\log\frac{\eps_0}{\tol}}{\log2}$ iterations and delivers the error
\begin{equation*}
  | u -u_{K} |_{H^1_0(\Omega)} \leq C_* \tol.
\end{equation*}
The number of iterations of {\rm \GALERKIN} is bounded uniformly for all outer loops.
\end{corollary}

%--------------------------------------------------------------------------------

\section{Convergence Rates of AFEM for Coercive Problems}\label{S:conv-rates-coercive}
% \rhn{(RHN $\longrightarrow$ CC)}

% {\color{brown}
%\begin{itemize}
%\item
%  Approximation classes for solution and data, relation to Besov spaces.

%\item
%  $\epsilon$-approximation of order $s$ \cite{BonitoDeVoreNochetto:2013}.

%\item
%  Optimal cardinality of second step (GALERKIN): optimal marking, cardinality of
%  marked set.

%\item
%  Optimal cardinality of AFEM: concatenation of first and second steps.

%\item
%  Numerical experiments: investigate relative cost of each module.  
%\end{itemize}
%}

The ultimate goal of $\AFEM$ is to produce a \emph{quasi-best} approximation $u_\grid\in\V_\grid$ to
the solution $u\in\V$ of \eqref{weak-form}
with error measured in $\V = H^1_0(\Omega)$. The performance of $\AFEM$ is measured by the
size of the error $|u-u_\grid|_{H^1_0(\Omega)}$ relative to the cardinality $\#\grid$ of $\grid$.
The latter usually reflects the total computational cost of implementing $\AFEM$. As a benchmark,
it is useful to compare the
performance of $\AFEM$ with the best approximation of $u\in\V$ and $\data = (\bA,c,f) \in \D$, provided we
have full knowledge of them. This is the main purpose of this section.

Under suitable assumptions on the solution $u$ and data $\data$, we prove the existence of constants $C(u,\data) >0$ and $s \in (0, \frac{n}{d}]$ such that
\begin{equation}\label{E:optimality-bound}
|u-u_{\grid_k}|_{H^1_0(\Omega)} \leq C(u,\data) \, \big( \#\mesh_k \big)^{-s} \, ,
\end{equation}
provided $s$ is the best decay rate with meshes in $\grids$ with a comparable number of degrees
of freedom. The upper bound $\frac{n}{d}$ of $s$ is dictated by the best decay rate with polynomials
of degree $n\ge1$ in dimension $d$ unless $u$ is degenerate (for instance, $u$ belongs to a finite
element space $\V_\grid$ with $\grid\in\grids$).
The dependence on $\data$ of the constant $C(u,\data)$ accounts for the multiplicative
structure of the interaction between the coefficients $(\bA,c)$ and $u$, and cannot be avoided in general.

A crucial insight for the simplest scenario, the Laplacian and piecewise constant
forcing $f$, is due to \cite{Stevenson:07}. It has been extended to operators with variable coefficients
by \cite{CaKrNoSi:08}
and later expressed in terms of the estimator by \cite{Axioms:2014}. It reads as follows:
\begin{equation}\label{E:stevenson}
\begin{minipage}{0.85\linewidth}
\emph{if a marking strategy reduces the PDE estimator $\eta_\grid(u_\grid)$ to
a fraction of its current value, then the refined set of elements ${\cal R}$
inherits an error indicator $\eta_\grid(u_\grid,{\cal R})$ comparable to $\eta_\grid(u_\grid)$,
hence a D\"orfler marking}.
\end{minipage}
\end{equation}
This allows to compare meshes produced by \AFEM 
with optimal ones and to conclude a quasi-optimal error decay.
To this end, we introduce in Section \ref{S:approx-classes} approximation classes for functions in $\V$ and $\D$, tailored to the decomposition of $\Omega$ into conforming refinements of an initial conforming partition $\grid_0$, the root of $\grids$. We will assume that $u=u(\data)\in\V$ and $\data=(\bA,c,f)\in\D$ belong to these classes which, however, are not characterized in terms of regularity of $u$ and $\data$. In Section \ref{S:eps-approx}, we investigate the approximability properties of perturbations $\wh{u}=u(\wh{\data})$ of the exact solution $u$, namely exact solutions of \eqref{E:perturbed-weak-form} with perturbed data $\wh\data$. Next, in Section \ref{S:optim-mark}, we consider a conforming refinement $\grid_*\in\grids$ of a partition $\grid\in\grids$, and give conditions under which an optimal D\"orfler marking property holds. We first apply this in Section \ref{S:quasi-opt-one-step} to study and derive rate-optimality of $\GALERKIN$ and $\AFEMSW$, the  one-step {\rm AFEM}s. We then combine the quasi-optimal performances of $\GALERKIN$ and $\DATA$ to prove rate-optimality of the two-step {\rm AFEM} in Section \ref{S:optimality-AVEM-TS}. We conclude in Section \ref{S:Besov} upon bridging the gap between appproximation and regularity classes. In particular, we give sufficient conditions for functions in Besov, Sobolev and Lipschitz spaces to belong to the approximation classes.

%%%%%%%%%%%%%%%%%%%%%%%%%%%%%%%%%%%%%%%%%%%%%%%%%%%%%%%%%%%%%%%%%%%%%%%%%%%%%%%%
\subsection{Nonlinear approximation classes}\label{S:approx-classes}
%%%%%%%%%%%%%%%%%%%%%%%%%%%%%%%%%%%%%%%%%%%%%%%%%%%%%%%%%%%%%%%%%%%%%%%%%%%%%%%%
%

In Section \ref{S:approx-solution} we discuss approximation classes for functions in $\V$,
which are applicable to the solution $u$ of \eqref{weak-form}. In Section \ref{S:approx-data}
we turn our attention to approximation classes for functions in $\D$, which are in turn
applicable to data $\data$. We refer to \cite{DeVore:98}, as well as \cite{DeVoreLorentz:1993,BiDaDeVPe:02} for a discussion within nonlinear
approximation theory.

%-------------------------------------------------------------------------------
\subsubsection{Nonlinear approximation classes for functions in $\V$}\label{S:approx-solution}
%-------------------------------------------------------------------------------
%
For any $N \in {\mathbb N}$, $N\ge\#\grid_0$, we define the following collection of partitions within $\grids$
\[
\index{Meshes!$\grids_N$: set of all conforming refinement of $\grid_0$ with no more than $N$ elements} 
\grids_N = \big\{ \grid :  \grid \in\grids  \text{ satisfies }  \  \# \grid  \leq N \big\} \ .
\]
This is the set of {\it conforming} meshes generated from $\grid_0$ with at most $N-\#\grid_0$ bisections.
Given $v\in\V$ we let $\sigma_N(v)$ be the smallest approximation $H^1_0$-error incurred
on $v$ with {\it continuous} piecewise polynomial functions of degree $\le n$ over meshes $\grids_N$:
\begin{equation}\label{E:sigmaN}
\sigma_N(v) :=  \inf_{\grid \in \grids_N} \inf_{v_\grid \in \V_\grid} |v - v_\grid|_{H^1_0(\Omega)}.
\end{equation}
This is a theoretical measure of performance in that finding a mesh $\mesh\in\grids_N$ that realizes
$\sigma_N(v)$ has exponential complexity.
Proving a bound $|v-v_\grid|_{1,\Omega}\le C_1\sigma_{C_2 N}(v)$ for $\grid\in\grids_{N}$ with $C_2 \le 1 \le C_1$
independent of $N$, the so-called {\it instance optimality}, is rather difficult and beyond the scope of this survey. In fact, a function
$v \in \V_{\mesh}$ with $\grid\in\grids_N$ could be the solution of our model problem \eqref{weak-form}, because we allow forcing
$f \in H^{-1}(\Omega)$. Hence, we see that $\sigma_N(v) = 0$ and {\rm AFEM} should then capture $v$ exactly on a finer mesh
$\mesh \in \grids_{C_2^{-1}N}$. We refer to \cite{DieningKreuzerStevenson:16} for a proof of instance optimality
for a forcing $f \in L^2(\Omega)$ and the Laplace operator, namely for coefficients $\bA=\vec{I}$ and $c=0$.

We will instead be able to prove that the error $|v-v_\grid|_{H^1_0(\Omega)}$ for the
Galerkin solution $v_\grid$ for $\grid \in \grids_N$ decays in terms of $N$ with
the same rate $N^{-s}$ as $\sigma_N(v)$; we thus say that $\AFEM$ is \emph{rate-optimal}. We first note that
for $v \in H^{n+1}(\Omega)$ and $\grid\in\grids_N$ quasi-uniform, we expect to have
\begin{equation}\label{E:linear-scale}
\inf_{v_\grid \in \V_\grid} |v-v_\grid|_{H^1_0(\Omega)} \lesssim N^{-\frac{n}{d}} |v|_{H^{n+1}(\Omega)}
\end{equation}
because the global meshsize $h$ and $N$ satisfy $h\approx N^{-1/d}$. This error
estimate within the linear Sobolev scale provides the largest possible decay rate $-n/d$.

\begin{definition}[approximation class of $u$]\label{D:approx-class-u}
\index{Functional Spaces!$\mathbb A_s$: approximation class for $u$}
Given $0<s\le n/d$, the class $\As:=\As\big(H^1_0(\Omega);\grid_0\big)$, relative to the
partition $\grid_0$ and approximation in the $H^1_0$-norm by continuous piecewise polynomials
of degree $\le n$ on the forest $\grids$ emanating from $\grid_0$,
is the set of functions $v\in\V=H^1_0(\Omega)$ such that
\begin{subequations}\label{E:quasi-norm}
\begin{equation}\label{E:quasi-norm-def}
|v|_\As := \sup_{N\ge\#\grid_0} \big( N^s \sigma_N(v) \big)<\infty
\end{equation}
whence
\begin{equation}\label{E:quasi-norm-decay}
\sigma_N(v) \le |v|_{\As} N^{-s} \quad\forall N\ge\#\grid_0.
\end{equation}
\end{subequations}
\end{definition}
We also write $\As = \mathbb{A}_s^0$ to emphasize continuity of the discrete functions in
$\V_\grid=\mathbb{S}^{n,0}_\grid \cap \V$ with $\grid\in\grids$.
The quantity $|v|_\As$ is a quasi seminorm in $\As$, which is not a linear space but rather a nonlinear
class of functions. Notice that as $s$ increases, the cost of membership
to be in $\As$ increases, namely $\A_{s_1} \subset \A_{s_2}$ for $s_1 \ge s_2$.

We may as well consider approximating $v\in\V$ with {\it discontinuous} piecewise polynomials
$\mathbb{S}^{n,-1}_\grid$ of degree $\le n$, which is a richer space than $\mathbb{S}^{n,0}_\grid$. We
can likewise define the corresponding modulus of approximation
\begin{equation}\label{E:sigmaN-1}
\sigma_N^{(-1)}(v) := \inf_{\grid\in\grids_N} \inf_{v_\grid\in \mathbb{S}^{n,-1}_\grid}
|v-v_\grid|_{H^1_0(\Omega;\grid)}
\end{equation}
and approximation class $\mathbb{A}_s^{-1}:=\mathbb{A}_s^{-1}\big(H^1_0(\Omega);\grid_0\big)$
of functions $v\in H^1_0(\Omega)$ such that
\begin{equation}\label{E:approx-class-disc}
|v|_{\mathbb{A}_s^{-1}} := \sup_{N\ge\#\grid_0} \big( N^s \sigma_N^{(-1)}(v) \big)<\infty
\quad\Rightarrow\quad
\sigma_N^{(-1)}(v) \le |v|_{\mathbb{A}_s^{-1}} N^{-s}.
\end{equation}
It is obvious that $\sigma_N^{(-1)}(v)\le\sigma_N(v)$ for all $v \in H^1_0(\Omega)$ because
$\mathbb{S}^{n,0}_\grid \subset \mathbb{S}^{n,-1}_\grid$. However, we have the following equivalence
result taken from \cite{Veeser:2016}. The original proof, although more complicated and for a different notion of
error relevant to discontinuous Galerkin approximations, can be traced back to \cite[Proposition 5.2]{bonito2010quasi}; see Proposition~\ref{p:equiv_classes_dG}.

\begin{proposition}[equivalence of classes]\label{P:equivalence-classes}
Assume that all stars of meshes $\grid\in\grids$ are $(d-1)$-face-connected. Then, there exists a
constant $C_{\textrm{dG}}$ that depends on the shape regularity of $\grids$, the dimension $d$ and the
polynomial degree $n\ge1$, such that
\[
\sigma_N(v) \le C_{\textrm{dG}} \, \sigma_N^{(-1)}(v) \quad\forall v\in H^1_0(\Omega), \quad N\ge\#\grid_0.
\]
Moreover, the approximation classes coincide $\mathbb{A}_s^{0} = \mathbb{A}_s^{-1}$.
\end{proposition}
\begin{proof}
We simply resort to \eqref{eq:gradient-approx} of Proposition \ref{P:cont-vs-discont} (approximation of
gradients), namely for $v\in H^1_0(\Omega)$
\begin{equation*}
    1 \le \frac{\min_{w\in\mathbb {S}^{n,0}_\grid} |v-w|_{H^1_0(\Omega)}}
    {\min_{w\in\mathbb {S}^{n,-1}_\grid} |v-w|_{H^1_0(\Omega;\grid)}}
    \le C_\textrm{dG},
  \end{equation*}
 and use the definitions \eqref{E:sigmaN} and \eqref{E:sigmaN-1}. This completes the proof.
\end{proof}

In the rest of the paper we will make the following approximability assumption.

\begin{assumption}[approximability of $u$]\label{A:approx-u}
\index{Assumptions!Approximability of $u$}
The exact solution $u\in H^1_0(\Omega)$ of problem \eqref{strong-form} belongs to the
approximation class $\As\big(H^1_0(\Omega);\grid_0\big)$ with $s=s_u \in (0,\frac{n}{d}]$.
\end{assumption}

The following condition \eqref{E:equiv-to-As} is simpler to handle in practice than \eqref{E:quasi-norm}.

\begin{lemma}[membership in $\As$]\label{L:equiv-to-As}
Let $v \in \As$ and $\eps_0 := \inf_{v_{\grid_0}\in\V_{\grid_0}} |v-v_{\grid_0}|_{H^1_0(\Omega)}$.
Then for all $0<\eps\le\eps_0$ there exist $\grid_\eps \in \grids$ and
$v_\eps\in \V_{\grid_\eps}$ such that
\begin{equation}\label{E:equiv-to-As}
  |v - v_\eps|_{H^1_0(\Omega)} \le \eps,
  \quad
  \#\grid_\eps \le 1+ |v|_{\As}^{\frac{1}{s}} \eps^{-\frac{1}{s}}.
\end{equation}
\end{lemma}
\begin{proof}
Given $0<\eps\le\eps_0$, let $\grid_\eps\in\grids$ be a conforming refinement of $\grid_0$ with minimal cardinality
and $v_\eps \in \V_{\grid_\eps}$ such that
\[
|v-v_\eps|_{H^1_0(\Omega)} \le \eps.
\]
Therefore, if $\eps<\eps_0$ we deduce from the minimal property of $\grid_\eps$ that
\[
\inf_{v_\grid\in\V_\grid} |v-v_\grid|_{H^1_0(\Omega)} > \eps \quad\forall \grid\in\grids :
\quad \#\grid \le \# \grid_\eps-1.
\]
If $N := \#\grid_\eps - 1$, definition \eqref{E:quasi-norm} implies
\[
\eps < \inf_{\grid\in\grids_N} \inf_{v_\grid\in\V_\grid} |v-v_\grid|_{H^1_0(\Omega)} = \sigma_N(v)
\le |v|_{\A_s} N^{-s},
\]
whence
\[
\#\grid_\eps = 1 + N  \le 1+ |v|_{\A_s}^{\frac{1}{s}} \eps^{-\frac{1}{s}},
\]
as asserted in \eqref{E:equiv-to-As}. On the other hand, if $\eps=\eps_0$ we see that
\[
\eps_0 \le |v|_{\As} (\#\grid_0)^{-s}
\quad\Rightarrow\quad
\#\grid_0 \le |v|_{\As}^{\frac{1}{s}} \eps_0^{-\frac{1}{s}} < 1 + |v|_{\As}^{\frac{1}{s}} \eps_0^{-\frac{1}{s}}.
\]
This completes the proof.
\end{proof}

\begin{remark}\label{R:case-n=1}
  If $d=2$ and $n=1$, then Corollary \ref{C:optimal-mesh} (optimal convergence rate) shows
  that $W^2_p(\Omega) \subset \A^{1/2}$ for $p>1$. The space $W^2_p(\Omega)$ is much larger
  than $H^2(\Omega)$, fits within the nonlinear Sobolev scale, and delivers the same decay
  rate as \eqref{E:linear-scale}. We will investigate the connection between approximation
  classes $\As$ and regularity classes in any dimension $d$ and for any polynomial degree
  $n \ge1$ later in Section \ref{S:Besov}.
\end{remark}

%-------------------------------------------------------------------------------
\subsubsection{Nonlnear approximation classes for data in $\D$}\label{S:approx-data}
%-------------------------------------------------------------------------------
%

Given data $\data=(\bA,c,f)\in\D$ and a mesh $\grid\in\grids$, we consider the best approximation
of $\data$ by discrete (piecewise polynomial) data
$\widehat{\data}=(\widehat{\bA},\widehat{c},\widehat{f})\in\D_\grid$, where $\D$ and $\D_\grid$ are
defined in \eqref{E:space-data} and \eqref{E:space-discrete-data}. We measure the error in
the space $D(\Omega)$ defined in \eqref{E:space-DO-redefine}
with $q=2$ for $d<4$ or $q > \frac{d}{2}$ for $d\ge 4$. We now discuss the best approximation errors for the components of data in $D(\Omega)$, which are used to define the corresponding approximation classes. For the coefficients $(\bA,c)$, they are characterized by the quantities
\begin{equation}
\label{e:delta_A_c}
\delta_\grid (\bA)_r := \inf_{\wh{\bA}\in[\mathbb{S}_\grid^{n-1,-1}]^{d \times d}} \|\bA-\wh{\bA}\|_{L^r(\Omega)},
\quad
\delta_\grid (c)_q := \inf_{\wh{c}\in\mathbb{S}_\grid^{n-1,-1}} \|c-\wh{c}\|_{L^q(\Omega)}
\end{equation}
for $r,q\in[2,\infty]$ as above.
Note that $\hbA$ and $\hc$ in \eqref{e:delta_A_c} are unconstrained in the sense that they do not necessarily satisfy the structural assumption \eqref{E:structural-assumption-wh} and are thus not suited coefficient for the perturbed problem \eqref{E:perturbed-weak-form}. We define the best constrained approximation errors for $\bA\in M(\alpha_1,\alpha_2)$ and $c\in R(c_1,c_2)$ by
\begin{equation}
\label{e:wtdelta_A_c}
\begin{split}
\wt \delta_\grid (\bA)_r &:= \inf_{\wh{\bA}\in[\mathbb{S}_\grid^{n-1,-1}]^{d \times d} \cap M(\halpha_1,\halpha_2)} \|\bA-\wh{\bA}\|_{L^r(\Omega)},\\
\wt \delta_\grid (c)_q &:= \inf_{\wh{c}\in\mathbb{S}_\grid^{n-1,-1} \cap R(\hc_1,\hc_2)} \|c-\wh{c}\|_{L^q(\Omega)},
\end{split}
\end{equation}
where in view of \eqref{E:wh-alpha}
\begin{equation}\label{d:values_const}
\halpha_1 = \frac{\alpha_1}2, \ \halpha_2 = C_\textrm{ctr} \alpha_2, \ \hc_1 = - \frac{\alpha_1}{4 C_P^2}, \ \hc_2 = C_\textrm{ctr}(\alpha_1+c_1).
\end{equation}
We anticipate that in Section~\ref{S:data}, we prove the equivalences
\begin{equation}\label{e:wtdelta_to_delta}
\delta_\grid (\bA)_r \leq \wt \delta_\grid (\bA)_r \leq \Cdata \delta_\grid (\bA)_r, \qquad \delta_\grid (c)_q \leq \wt\delta_\grid (c)_q \leq \Cdata \delta_\grid (c)_q,
\end{equation}
for all $\bA \in M(\alpha_1,\alpha_2)$ and $c \in R(c_1,c_2)$; see Remarks~\ref{r:equivalent_class_A} and~\ref{r:equivalent_class_c}.
For the load function $f$, the definition \eqref{global-osc-with-elms} of $\osc_\grid(f)_{-1}$ suggests to consider
\begin{equation*}
\delta_\grid (f)_{-1} := \left(\sum_{T\in\grid} \inf_{\wh{f} \in \F_{\grid_{\omega_T}}}
\|f-\wh{f}\|_{H^{-1}(\omega_T)}^2\right)^{\frac12}, 
\end{equation*}

All these best approximation errors are hard
to evaluate and are thus replaced by the computable oscillations defined in \eqref{E:zeta-Ac-global}
and \eqref{E:zeta-f-global} in practice. We recall that they rely on the local $L^2$-projection
operator $\Pi_\grid$ for $(\bA,c)$ and the local $H^{-1}$-projection operator $P_\grid$ for $f$
to compute linear approximations $\wt\data$ of $\data$ to a desired accuracy. These projections are
later modified nonlinearly to give rise to $\hdata$ satisfying the side constraints \eqref{E:structural-assumption-wh} without compromising accuracy.
We recall that the \DATA module is assumed to construct approximations so that
\begin{equation*}
\index{Constants!$\Ldata$: \DATA quasi-optimality constant}
\osc_\grid(\bA)_r \le \Ldata \delta_\grid(\bA)_r,
\;\;
\osc_\grid(c)_q \le  \Ldata \delta_\grid(c)_q,
\;\;
\osc_\grid(f)_{-1} \le  \Ldata \delta_\grid(f)_{-1} \looseness=-1
\end{equation*}
with a mesh independent constant $\Ldata$; see Assumption~\ref{A:prop-data}.
We discuss in Section~\ref{S:DATA} practical realizations of \DATA. 

For the
purpose of assessing the cardinality of {\rm AFEM}, we do not need the specific form of $\wh\data$ but rather the decay of the best approximation errors in terms of degrees of freedom. Therefore, we postpone
until Section \ref{S:data-approx} the construction of $\wh\data$ for $n\ge1$ and to Section \ref{S:Besov}
the discussion of regularity
properties of $\data$ that guarantee membership in the following approximation classes.

\begin{definition}[approximation classes of $\bA$]\label{D:class-M}
\index{Functional Spaces!$\ACA_s$: approximation classes of $\bA$}
For $0<\alpha_1\leq \alpha_2$, $2\leq r \leq \infty$, let 
$\ACA_s :=  \ACA_s \big(L^r(\Omega;\R^{d\times d});\grid_0\big)$ be the set of matrix-valued functions
$\bA \in M(\alpha_1,\alpha_2)$ satisfying
\begin{equation}\label{E:zeta-A}
| \bA |_{{\ACA}_{s}}:=
\sup_{N \geq \# \grid_0}
\Big( N^s\inf _{\grid \in \mathbb{T}_N} \wt \delta_\grid (\bA)_r \Big) < \infty
\,\Rightarrow\,
 \inf _{\grid \in \mathbb{T}_N}\wt \delta_\grid (\bA)_r \le | \bA |_{{\ACA}_{s}} N^{-s}.
\end{equation}
\end{definition}

\begin{definition}[approximation classes of $c$] \label{D:class-C}
\index{Functional Spaces!$\ACc_s$: approximation classes of $c$}
The class $\ACc_{s} := \ACc_{s}\big(L^q(\Omega);\grid_0\big)$
is the set of functions $c\in R(c_1,c_2)$ such that
\begin{equation}\label{E:zeta-c}
| c |_{{\ACc}_{s}}:=
\sup_{N \geq \# \grid_0} 
\Big( N^s \inf _{\grid\in \mathbb{T}_N} \wt \delta_\grid (c)_q \Big) < \infty
\,\,\Rightarrow\,\, %\quad\Rightarrow\quad
\inf _{\grid \in \mathbb{T}_N} \wt \delta_\grid (c)_q \le | c |_{{\ACc}_{s}} N^{-s}. \looseness=-1
\end{equation}
\end{definition}
\begin{definition}[approximation classes of $f$] \label{D:class-F}
\index{Functional Spaces!$\ACf_s$: approximation classes of $f$}
The class $\ACf_{s} := \ACf_{s}\big(H^{-1}(\Omega);\grid_0\big)$
is the set of functions $f\in H^{-1}(\Omega)$ such that
\begin{equation}\label{E:zeta-f}
\hskip-0.35mm
| f |_{{\ACf}_{s}}:=
\sup_{N \geq \# \mesh_0}
\Big( N^s \inf_{\grid\in \mathbb{T}_N} \delta_\grid (f)_{-1} \Big) < \infty
\,\Rightarrow\, %\!\!\quad\!\Rightarrow\!\quad\!\!
\inf _{\grid \in \mathbb{T}_N} \delta_\grid (f)_{-1} \le | f |_{{\ACf}_{s}} N^{-s}. \looseness=-1
\end{equation}
\end{definition}

Since the polynomial degree of discrete coefficients $(\wh{\bA},\wh{c})$
in definition \eqref{E:space-discrete-data} is $n-1$,
we expect decay rates $s_A,s_c \le n/d$ according to nonlinear approximation theory.
The specific values of $(s_A,s_c)$ depend on the
regularity of $(\bA,c)$, a delicate topic that we further investigate in Sections \ref{S:Besov} and
\ref{S:data-approx}.
However, because $u$ and $\data = (\bA,c,f)$ satisfy the elliptic problem \eqref{strong-form}, the
above approximation classes are somewhat related. We now quantify this statement.

\begin{lemma}[relation between approximation classes]\label{L:approx-regularity}
Let $2 \le r,q \le \infty$ be so that $\frac{d}{2}<q$.
If $u\in\A_{s_u}\big(H^1_0(\Omega);\grid_0\big)$,
$\bA\in\ACA_{s_A} \big(L^r(\Omega;\R^{d\times d});\grid_0\big)$, and
$c\in\ACc_{s_c}\big(L^q(\Omega);\grid_0\big)$, with $0<s_u,s_A,s_c\le\frac{n}{d}$,
then $f\in\ACf_{s_f}\big(H^{-1}(\Omega);\grid_0\big)$ and
\begin{equation}\label{E:f-uAc}
|f|_{\ACf_{s_f}} \le C \Big( |u|_{\A_{s_u}}  +  |\bA|_{\ACA_{s_A}}   + |c|_{\ACc_{s_c}} \Big),
\quad
s_f = \min\{s_u,s_A,s_c\},
\end{equation}
where the constant $C>0$ depends $\|u\|_{W^1_p(\Omega)}$, $p=\frac{2r}{r-2}$,
and $\alpha_1, \alpha_2, c_1, c_2$. In particular, if $(\bA,c)$ are discrete in $\grid_0$, then
\begin{equation}\label{E:f-u}
|f|_{\ACf_{s_f}} \le C |u|_{\A_{s_u}}, \quad s_f = s_u.
\end{equation}
\end{lemma}
\begin{proof}
Let $L[u] := -\div (\bA\nabla u) + cu$ be the operator in \eqref{strong-form} and note that
$f=L[u] \in H^{-1}(\Omega)$ can be approximated by
$\wh{f} = -\div (\wh{\bA}\nabla \wh{v}) + \wh{c}\wh{v} \in \F_\grid$, where the
discrete space $\F_\grid$ is given in Definition \ref{D:discrete-functionals} and
$\wh{v}\in\mathbb{S}_\grid^{n,0}, \wh{\bA}\in(\mathbb{S}_\grid^{n-1,-1})^{d \times d},
\wh{c} \in \mathbb{S}_\grid^{n-1,-1}$. Let's now express $f-\wh{f}$ as follows
\begin{align*}
f-\wh{f} = - \div\big((\bA-\wh{\bA}) \nabla u \big) + (c-\wh{c}) u - \div\big( \wh{\bA} \nabla (u-\wh{v})\big)
+ \wh{c} (u-\wh{v}),
\end{align*}
and recall that we have to estimate $\|f-\wh{f}\|_{H^{-1}(\omega_T)}$ for every $T\in\grid$, rather than
a global norm in $H^{-1}(\Omega)$, to get an upper bound on $\delta_\grid(f)_{-1}$. 
%We thus resort
%to Lemma \ref{L:loc-of-H^-1-norm} (localization of $H^{-1}$-norm) and the finite overlapping property
%of the patches $\{\omega_T\}_{T\in\grid}$ with constant $C_\text{ovrl}$ depending only on shape regularity
%of $\grids$
%
%\[
%\sum_{T\in\grid} \|f-\wh{f}\|_{H^{-1}(\omega_T)}^2 \le C_\text{ovrl}^2 \|f-\wh{f}\|_{H^{-1}(\Omega)}^2.
%\]
%
Therefore, we proceed as in the proof of Lemma \ref{L:perturbation} (continuous
dependence on data) to obtain
\begin{align*}
\sum_{T\in\grid} \|f-\wh{f}\|_{H^{-1}(\omega_T)}^2 &\lesssim
\|\nabla u\|_{L^p(\Omega)}^2 \|\bA-\wh{\bA}\|_{L^r(\Omega)}^2
+ \|\nabla u\|_{L^2(\Omega)}^2 \|c-\wh{c}\|_{L^q(\Omega)}^2
\\& + \|\wh{\bA}\|_{L^\infty(\Omega)}^2 \|\nabla (u-\wh{v})\|_{L^2(\Omega)}^2
+ \|\wh{c}\|_{L^\infty(\Omega)}^2 \|u-\wh{v}\|_{L^2(\Omega)}^2,
\end{align*}
where $p=\frac{2r}{r-1}$ and $\|\nabla u \|_{L^p(\Omega)} < \infty$ according to \eqref{E:Wp-reg}
and $\frac{d}{2}<q\le\infty$. 
Note that thanks to \eqref{d:values_const}, $\| \wh A \|_{L^\infty(\Omega)} \leq \halpha_2 = \Cctr \alpha_2$ and $\| \wh c\|_{L^\infty(\Omega)} \leq \hc_2 = \Cctr (\alpha_1+c_2)$.
Moreover, since $\wh{v}, \wh{\bA}, \wh{c}$ can be chosen separately,
invoking \eqref{E:quasi-norm}, \eqref{E:zeta-A}, \eqref{E:zeta-c}, and \eqref{E:zeta-f},
we realize that
\begin{align*}
\inf_{\grid\in\grids_N} \delta_\grid(f)_{-1} 
& \lesssim \inf_{\grid\in\grids_N} \inf_{\wh{v}\in\mathbb{S}_\grid^{n,0}}\|\nabla (u-\wh{v})\|_{L^2(\Omega)}
\\& +  \inf_{\grid\in\grids_N} \inf_{\wh{\bA}\in[\mathbb{S}_\grid^{n-1,-1}]^{d\times d}} \|\bA - \wh{\bA}\|_{L^r(\Omega)}
\\& +  \inf_{\grid\in\grids_N} \inf_{\wh{c}\in\mathbb{S}_\grid^{n-1,-1}} \|c - \wh{c}\|_{L^q(\Omega)}
\\& \le  |u|_{\A_{s_u}} N^{-s_u} +  |\bA|_{\ACA_{s_A}} N^{-s_A} + |c|_{\ACc_{s_c}} N^{-s_c}
\end{align*}
gives \eqref{E:f-uAc} with $s_f = \min\{s_u,s_A,s_c\}$; \eqref{E:f-u} is a trivial consequence.
\end{proof}

Estimate \eqref{E:f-u} will be useful later in Theorem \ref{T:optimality-one-step} (rate-optimality
of one-step $\AFEM$s). 
It is important to realize that the multiplicative structure between solution $u$ and
coefficients $(\bA,c)$ is hidden in the constants $C$ in \eqref{E:f-uAc} and \eqref{E:f-u}.
Moreover, these estimates are possible due to the fact that the space $H^{-1}(\Omega)$ is the range of
the linear operator $L:H^1_0(\Omega) \to H^{-1}(\Omega)$ and that the discrete functions in
$\F_\grid$ are images by $L$ of functions in $\V_\grid$. This would not be true for $L^2$-weighted
surrogates of $\delta_\grid(f)_{-1}$ that typically overestimate the error in $H^{-1}(\Omega)$.

\begin{assumption}[approximability of data] \label{A:approx-data} There exist
\index{Assumptions!Approximability of data}
$s_A, s_c, s_f \in  (0, \frac{n}{d}]$ such that data $\data=(\bA,c,f)\in\D$ satisfies
$\bA \in \ACA_{\, s_{\! A}}$, $c \in \ACc_{s_c}$, $f \in \ACf_{s_f}$.
\end{assumption}

We recall that, if $\osc_\grid(\data)=\|\data-\wh\data\|_{D(\Omega)} > \Cdata \tau$ over a conforming refinement
$\grid\in\grids$ of $\grid_0$, then the call
\[
[\wh\grid,\wh\data] = \DATA \, (\grid,\data,\tau)
\]
produces a conforming refinement $\wh\grid$ of $\grid$ and  approximate data $\hdata = (\hbA,\hc,\hf)\in\D_{\wh\grid}$ over $\wh\grid$ that satisfies $\osc_{\wh\grid}(\data)\le\Ldata\delta_{\wh\grid}(\data)$ and for $r,q \in [2,\infty]$
\[
\osc_{\wh\grid}(\data) = \osc_{\wh\grid} (\bA)_r +  \osc_{\wh\grid} (c)_q + \osc_{\wh\grid} (f)_{-1} \le \Cdata \tau,
\]
as well as the
constraints $\hbA \in M(\halpha_1,\halpha_2)$ and $\hc \in R(\hc_1,\hc_2)$ defined in
\eqref{E:structural-assumption-wh}. We will show in Section \ref{S:data-approx} that the
routine responsible for reducing oscillations, namely $\GREEDY$, exhibits optimal
performance in the sense that the cardinalities $N_\grid(\bA), N_\grid(c), N_\grid(f)$
of the sets of elements necessary to reduce the individual oscillations of $(\bA,c,f)$
below the threshold $\Cdata\frac{\tau}{3}$ starting from any $\grid\ge\grid_0$ satisfy
\begin{equation}\label{eq:decay-data}
N_\grid(\bA) \lesssim  | \bA |_{\ACA_{\, s_{\! A}}}^{\frac{1}{s_A}} \tau^{-\frac{1}{s_A}} \,, \quad 
N_\grid(c) \lesssim  | c |_{\mathbb{C}_{s_c}}^{\frac{1}{s_c}} \tau^{-\frac{1}{s_c}} \,, \quad 
N_\grid(f) \lesssim  | f |_{\mathbb{F}_{s_f}}^{\frac{1}{s_f}} \tau^{-\frac{1}{s_f}} \,.
\end{equation}
Therefore, the cost of one call to $\DATA$ can be quantified by the total number
$N_\grid(\data)$ of elements marked, which obeys the relation
\begin{align*}
N_\grid(\data) &= N_\grid(\bA) + N_\grid(c) + N_\grid(f)
\\ &
\lesssim | \bA |_{\ACA_{\, s_{\! A}}}^{\frac{1}{s_A}} \tau^{-\frac{1}{s_A}}
+ | c |_{\mathbb{C}_{s_c}}^{\frac{1}{s_c}} \tau^{-\frac{1}{s_c}}
+ | f |_{\mathbb{F}_{s_f}}^{\frac{1}{s_f}} \tau^{-\frac{1}{s_f}}
\le |\data|_{\mathbb{A}_\data}^{\frac{1}{s_\data}} \tau^{-\frac{1}{s_\data}}
\end{align*}
with
\begin{equation}\label{E:class-data}
s_\data := \min\big\{s_A, s_c, s_f \big\}, \quad
| \data |_{\mathbb{A}_\data} := \left( | \bA |_{\ACA_{\, s_{\! A}}}^{\frac{1}{s_A}} +
| c |_{\mathbb{C}_{s_c}}^{\frac{1}{s_c}} +
| f |_{\mathbb{F}_{s_f}}^{\frac{1}{s_f}} \right)^{s_\data}.
\end{equation}
It is thus natural to make the following assumption on $\DATA$.

\begin{assumption}[quasi-optimality of $\DATA$]\label{A:optim-data}
\index{Assumptions!Quasi-optimality of \DATA}
The call $[\wh\grid,\wh\data] \!=\! \DATA \, (\grid,\data,\tau)$, from an arbitrary conforming
refinement $\grid$ of $\grid_0$ with tolerance $\tau$, marks the number of elements
$N_\grid(\data)$ to produce an approximation $\wh\data$ of $\data$ over $\wh\grid$
so that
\begin{equation}
\osc_{\wh\grid}(\data) = \|\data - \wh\data\|_{D(\Omega)} \le \Cdata\tau, \quad
N_\grid(\data) \lesssim | \data |_{\mathbb{A}_\data}^{\frac{1}{s_\data}} \, \tau^{-\frac{1}{s_\data}} 
\, .
\end{equation}
\end{assumption}

%%%%%%%%%%%%%%%%%%%%%%%%%%%%%%%%%%%%%%%%%%%%%%%%%%%%%%%%%%%%%%%%%%%%%%%%%%%%%%%%
\subsection{$\eps$-approximation of order $s$}\label{S:eps-approx}
%%%%%%%%%%%%%%%%%%%%%%%%%%%%%%%%%%%%%%%%%%%%%%%%%%%%%%%%%%%%%%%%%%%%%%%%%%%%%%%%
%
Inspection of the structure of algorithm $\AFEMTS$ (Algorithm~\ref{algo:AFEMTS}) reveals that
the approximate data $\data_k$ is fixed inside $\GALERKIN$. Therefore, the performance of $\GALERKIN$ is dictated by the regularity of the exact solution $\wh{u}_k=\wh{u}_k(\data_k)\in H^1_0(\Omega)$ with data $\data_k$, rather than the exact solution $u=u(\data)$ with data $\data$.  We know that $u\in \As$ and wonder what regularity is inherited by $\wh{u}_k$. This leads to the following concept introduced in \cite[Def.  3.1 and Lemma 3.2]{BonitoDeVoreNochetto:2013}.
\begin{definition}[$\varepsilon$-approximation of order $s$]\label{D:approx-order-s}
\index{Definitions!$\eps$-approximation of order $s$}
Given $u\in\As\big(H^1_0(\Omega);\grid_0\big)$ and $\varepsilon>0$,  a function $v\in H^1_0(\Omega)$ is said to be an $\varepsilon$-approximation of order $s$ to $u$ if $| u - v|_{H^1_0(\Omega)}\leq \varepsilon$ and there exists a constant $C>0$ independent of $\varepsilon$, $u$ and $v$ such that for all $\delta\geq \varepsilon$ there exists $N \geq \# \grid_0$ satisfying
\begin{equation}\label{E:eps-approx}
\sigma_N(v) \leq \delta
\qquad
N\leq 1 + C | u |_{\As}^{\frac{1}{s}}\delta^{-\frac{1}{s}}.
\end{equation}
\end{definition}

\begin{lemma}[$\varepsilon$-approximation of $u$ of order $s$]\label{L:eps-approx}
Let $u\in \As\big(H^1_0(\Omega);\grid_0\big)$ and $v\in H^1_0(\Omega)$ satisfy $| u-v |_{H^1_0(\Omega)}\leq \varepsilon$ for some $0<\varepsilon\le\eps_0$ with $\eps_0$ defined in Lemma \ref{L:equiv-to-As}. Then $v$ is a $2\varepsilon$-approximation of order $s$ to $u$.  
\end{lemma}
\begin{proof}
Let $\delta\geq 2\varepsilon$. By definition \eqref{E:sigmaN} of $\sigma_N(v)$, it suffices to invoke the triangle inequality to realize that
\[
\sigma_N(v) \le | u-v |_{H^1_0(\Omega)} + \sigma_N(u) \le \frac{\delta}{2} + \sigma_N(u).
\]
Since $u\in \As\big(H^1_0(\Omega)\big)$, in view of Lemma \ref{L:equiv-to-As}, there exist $N \geq \# \grid_0$
and $\mathcal{T}\in \mathbb{T}_N$
\[
\sigma_N(u) \le \frac{\delta}{2},
\qquad
N \le 1 + |u|_{\As}^{\frac{1}{s}} \Big(\frac{\delta}{2}\Big)^{-\frac{1}{s}}.
\]
The estimate \eqref{E:eps-approx} thus follows with constant $C=2^{\frac{1}{s}}$.
\end{proof}

This is a simple but crucial result to study $\AFEMTS$. It says that any function $v$ that is $\eps$-close
to a function $u\in\As(X;\grid_0)$ in the norm of the space $X$ defining the approximation class
$\As(X;\grid_0)$ can be approximated with a similar decay rate as $u$ in $X$ for as long as
the desired accuracy does not exceed $\eps$. In other words, the approximability of $u$ is inherited by
$v$ up to scale $\eps$. However, beyond the scale $\eps$ the approximability of $v$
may differ from that of $u$.
Note that neither the definition \eqref{E:sigmaN} of $\sigma_N(v)$
nor Lemma \ref{L:eps-approx} require $X = H^1_0(\Omega)$.

%%%%%%%%%%%%%%%%%%%%%%%%%%%%%%%%%%%%%%%%%%%%%%%%%%%%%%%%%%%%%%%%%%%%%%%%%%%%%%%%
\subsection{Properties of D\"orfler marking}\label{S:optim-mark}
%%%%%%%%%%%%%%%%%%%%%%%%%%%%%%%%%%%%%%%%%%%%%%%%%%%%%%%%%%%%%%%%%%%%%%%%%%%%%%%%
%

We follow the ideas of \cite{Stevenson:07}, \cite{CaKrNoSi:08} and \cite{Axioms:2014} to explore the insight \eqref{E:stevenson}
about D\"orfler marking.
Hereafter, we recall \eqref{E:grid-gridstar} and
consider two admissible partitions $\grid,\grid_*\in\grids$ such that  $\grid\le\grid_*$, i.e. the latter is a refinement of the former obtained by applying (newest-vertex) bisection to some of the elements of $\grid$.

In the sequel, we let $u = \wh{u} \in H^1_0(\Omega)$ be the exact solution with discrete
coefficients $(\wh{\bA},\wh{c})$ over a fixed mesh $\wh{\grid}\le\grid$ and forcing function $f$ that may
or may not be discrete. 
We rewrite the a posteriori error estimates \eqref{E:aposteriori-bounds-H1}
\begin{equation}\label{E:apost-discrete-data}
\index{Constants!$(C_L,C_U)$: a-posteriori lower and upper bounds constants}
C_L\,\estG\big(\uG,f\big) \le |\wh{u}-\uG|_{H^1_0(\Omega)} \le C_U\,\estG\big(\uG,f\big),
\end{equation}
where the total estimator $\estG\big(\uG,f\big)$ consists of the PDE estimator \rhn{$\eta_\grid(u_\grid,f)$}
and the oscillation $\osc_\grid(f)_{-1}$ and reads, according to \eqref{E:global-est-osc},
\[
\estG\big(\uG,f\big)^2 = \rhn{\eta_\grid(u_\grid,f)^2} + \osc_\grid(f)_{-1}^2.
\]
We also recall that when $f =P_\grid f  \in \mathbb F_{\grid}$ is discrete,  $\osc_\grid(f)_{-1}=0$ and $\estG\big(\uG,f\big)$ reduces to $\eta_\grid(u_\grid,f)$, \rhn{and that $P_\grid f$ is used within
$\eta_\grid(u_\grid,f)$ rather than $f$.}
The global nature of the elliptic boundary value problem \eqref{strong-form} prevents upper a posteriori
energy error estimates such as \eqref{E:apost-discrete-data} between the continuous and discrete solution
to be local. Remarkably, the situation for two Galerkin
solutions $u_\grid\in\V_\grid$ and $u_{\grid_*}\in\V_{\grid_*}$ is different, as stated in
Theorem \ref{T:ubd-corr} (upper bound for corrections) \looseness=-1
\begin{equation}\label{E:local-upper}
\index{Constants!$\wt{C}_U$: localized upper bound constant}
  |u_{\grid_*} - u_\grid|_{H^1_0(\Omega)} \le \wt{C}_U\,\estG(u_\grid,f,\refined[]),
\end{equation}
where $\refined[] = \grid \backslash \grid_*$ is the refined set defined in \eqref{E:refined-set}.
It thus turns out that $|u_{\grid_*} - u_\grid|_{H^1_0(\Omega)}$ is controlled by the estimator
$\estG(u_\grid,f,\refined[])$ on the set of elements $\refined[]$ where the meshes differ.
This crucial observation goes back to Stevenson \cite{Stevenson:07}; see also \cite{CaKrNoSi:08,NoSiVe:09}.

In the sequel we impose restrictions on the ranges of the D\"orfler marking parameter \eqref{E:theta} and
the threshold parameter $\omega$ for $\GALERKIN$ and $\AFEMSW$. We will impose a different restriction
later on $\omega$ for $\AFEMTS$.

\begin{assumption}[marking parameter]\label{A:theta-no-osc}
\index{Assumptions!Marking parameter $\theta$}
Let $\theta$ satisfy $\theta\in(0,\theta_0)$ with
\[
\theta_0 := \min\Big\{\big(2\Clip\wt{C}_U\big)^{-1}, 1\Big\},
\]
where $\Clip, \wt{C}_U$ are the constants in \eqref{E:lipschitz-estimator} and \eqref{E:local-upper}
respectively.
\end{assumption}

\begin{assumption}[restriction on $\omega$]\label{A:PDE-dominates}
\index{Assumptions!Restriction on $\omega$}
We assume $0\le\omega\le\omega_0<1$ with
\AB{
\[
\omega_0 :=  \sqrt{\frac{\theta_0^2-\theta^2}{2 +  \theta_0^2-\theta^2}}
\]
}
\end{assumption}

We now ready to make Stevenson's insight \eqref{E:stevenson} precise.

\begin{lemma}[D\"orfler marking]\label{L:dorfler-no-osc}
Let Assumptions \ref{A:theta-no-osc} and \ref{A:PDE-dominates} hold and $0<\mu\leq \frac 12$. 
Let $\grid \in \grids$, and let $\grid_*\in\grids$ be a refinement of $\grid$ with respective Galerkin solutions $u_\grid \in \V_\grid$ and
$u_{\grid_*} \in \V_{\grid_*}$; let ${\cal R}=\grid \setminus \grid_*$ be the refined set. 
Assume that the oscillation on $\grid$ is dominated by the total estimator
\begin{equation}\label{e:PDE_est_dominate}
\osc_\grid(f)_{-1} \le \omega \, \est_\grid(u_\grid,f)
\end{equation}
and that
\begin{equation}\label{E:mu-no-osc}
\rhn{\eta_{\grid_*}(u_{\grid_*},f) \leq \mu \, \eta_{\grid}(u_{\grid},f).}
\end{equation}
Then D\"orfler marking is valid for any $0 < \theta < \theta_0$
\begin{equation}\label{E:dorfler-no-osc}
\rhn{\theta \,\eta_\grid(u_\grid,f) \leq \eta_\grid(u_\grid,f,{\cal R}).}
\end{equation}
\end{lemma}
\begin{proof}
\rhn{We invoke Proposition \ref{P:est-reduction} (estimator reduction)
with $\delta=1$ along with the localized upper bound \eqref{E:local-upper}
to write
$$
\eta_{\grid}(u_{\grid},f)^2 \le 2\, \eta_{\grid_*}(u_{\grid_*},f)^2 + 2
\Clip^2 \left( \wt{C}_U^2 \est_\grid(u_\grid,f,\mathcal R)^2 + \osc_{\grid}(f)_{-1}^2\right).
$$
The last term accounts for the presence of $P_\grid f$ and $P_{\grid_*} f$ in the definitions 
of $\eta_{\grid}(u_{\grid},f)$ and $\eta_{\grid_*}(u_{\grid_*},f)$.
}
\AB{
In view of \eqref{e:PDE_est_dominate} and the definition \eqref{E:global-est-osc} of the total estimator, we have
\[
\osc_\grid(f)_{-1} \le \sigma \, \eta_\grid(u_\grid), \quad \sigma^2 := \frac{\omega^2}{1-\omega^2}
\]
so that
\rhn{
$$ 
\eta_{\grid}(u_{\grid},f)^2 \le 2\eta_{\grid_*}(u_{\grid_*},f)^2 + 2
\Clip^2 \wt{C}_U^2\left(  \eta_\grid(u_\grid,f,\mathcal R)^2 
+ 2 \sigma^2 \eta_\grid(u_\grid,f)^2\right).
$$

Using \eqref{E:mu-no-osc} and rearranging the above expression we obtain
\[
\big( \theta_0^2 - 2\sigma^2 \big) \eta_{\grid}(u_{\grid},f)^2  \le
\left( \frac{1-2\mu^2}{2\Clip^2 \wt{C}_U^2} - 2\sigma^2  \right)
\eta_{\grid}(u_{\grid},f)^2 \le \eta_\grid(u_\grid,f,{\cal R})^2
\]
}
provided $0<\mu\le\frac12$, because of the definition of $\theta_0$ in Assumption
\ref{A:theta-no-osc}. Finally, for any $\theta<\theta_0$ we realize that
$\omega_0$ from Assumption \ref{A:PDE-dominates} satisfies
\[
\sigma_0^2:=\frac{\AB{\omega_0^2}}{1-\omega_0^2}=\frac 1 2(\theta_0^2-\theta^2)
\quad\Rightarrow\quad
\theta^2 = \theta_0^2 - 2\sigma_0^2 \le \theta_0^2 -2\sigma^2,
\]
and D\"orfler marking is valid for ${\cal R}$ with parameter $\theta$.
}
\end{proof}

We remark that Lemma~\ref{L:dorfler-no-osc} requires that the oscillation on $\grid$ is dominated by the total estimator to guarantee a \AB{D\"orfler} marking property. This is always the case when $f$ is discrete as in Algorithm~\ref{A:GALERKIN} (\GALERKIN) because in that case $\osc_{\grid}(f)_{-1}=0$, or within Algorithm~\ref{A:one-step-switch} (one-step \AFEM with switch), which marks elements for refinement only if this property holds.

We also see that $\theta_0$ in Assumption \ref{A:theta-no-osc} corresponds to the choices
$\mu=\frac12$ and $\delta=1$. However, the proof reveals that for $\mu\to0$ we could obtain the largest
possible value $\theta_0 = (\Clip \wt{C}_U)^{-1}$, thereby the less restrictive.
Since this is just twice the value of $\theta_0$
in Assumption \ref{A:theta-no-osc}, the practical choice $\mu=\frac12$ is justified.

Lemma \ref{L:dorfler-no-osc} hinges on two ingredients: the Lipschitz property \eqref{E:lipschitz-estimator}
and the localized upper bound \eqref{E:local-upper} of the estimator. In particular, it does not
rely on the lower a posteriori error estimate in \eqref{E:apost-discrete-data}, as the original
proofs in \cite{Stevenson:07} and \cite{CaKrNoSi:08}, and easily extends
to discontinuous Galerkin methods in Section \ref{S:dg} and inf-sup stable methods in
Section \ref{S:conv-rates-infsup}. The original statement is, however, a bit more insighful:
if $\theta_0^2 = \frac{C_2^2}{2C_1^2}$, then for all $0<\theta<\theta_0$, $\omega^2 \le \theta_0^2-\theta^2$ 
\[
\enorm{u-u_{\grid_*}} \le\mu \enorm{u-u_\grid}
\quad\Rightarrow\quad
\eta_\grid(u_\grid,\refined[])\ge\theta\,\eta_\grid(u_\grid)
\]
provided $0<\mu\le 2^{-\frac12}$. We see that the threshold $\theta_0$ is related to the gap between
{\it reliability} constant $C_1$ and {\it efficiency} constant $C_2$
in the a posteriori bounds \eqref{E:aposteriori-bounds} in the energy norm; hence the
ratio $\frac{C_2}{C_1}\le 1$ is a quality measure of the estimator $\eta_\grid(u_\grid)$.
It is thus reasonable to be cautious about marking decisions if the constants $C_1$ and
$C_2$ are very disparate, and thus the ratio
$C_2/C_1$ is far from $1$. This justifies the constraint $\theta<\theta_0$.

%%%%%%%%%%%%%%%%%%%%%%%%%%%%%%%%%%%%%%%%%%%%%%%%%%%%%%%%%%%%%%%%%%%%%%%%%%%%%%%%%
\subsection{Rate-optimality of one-step {\rm \AFEM}s}\label{S:quasi-opt-one-step}
%%%%%%%%%%%%%%%%%%%%%%%%%%%%%%%%%%%%%%%%%%%%%%%%%%%%%%%%%%%%%%%%%%%%%%%%%%%%%%%%%
%

Recall that  $\marked_j$ is the output of the module $\MARK$ and that $\grid_j,u_j$ are the meshes and associated Galerkin solutions generated within Algorithms \ref{A:GALERKIN}
($\GALERKIN$) and \ref{A:one-step-switch} (one-step $\AFEM$ with switch). 
To express the cardinality $N_j(u)$ of $\marked_j$ in terms of $|u-u_j|_{H^1_0(\Omega)}$
we must relate the performance of these one-step $\AFEM$s with the approximation classes
$\As=\As\big(H^1_0(\Omega);\grid_0\big)$ for $u$ and $\F_s=\F_s\big(H^{-1}(\Omega);\grid_0\big)$
for $f$, which are never used in the design of these algorithms.
Even though this might appear like an undoable task, the key
to unravel this connection is given by Lemma
\ref{L:dorfler-no-osc} (D\"orfler marking) and the following assumption.
\begin{assumption}[cardinality of $\marked$]\label{A:cardinality}
\index{Assumptions!Cardinality of the marked set}
The module $\MARK$ selects a set $\marked$ in \eqref{E:theta} with \emph{minimal} cardinality.
\end{assumption}

According to the equidistribution principle \eqref{E:equidistribution}
and the local lower bound \eqref{e:equivalence_tmp} in the proof of Theorem \ref{T:modified-estimator} (modified residual estimator) for  discrete coefficients, i.e.,
\[
C_L \eta_{\grid_j}(u_j,T) \le C_L\est_{\grid_j}(u_j,f,T)\le |u-u_j|_{H^1(\omega_T)},
\]
it is
natural to mark elements with largest error indicators. This explains the choice of 
a minimal set $\marked$ in Assumption \ref{A:cardinality}.

We are now ready to bound the cardinality of $\marked_j$ in terms of $|u-u_j|_{H^1_0(\Omega)}$.

\begin{proposition}[cardinality of $\marked_j$]\label{P:cardinality-marked}
Let Assumptions \ref{A:theta-no-osc}, \ref{A:PDE-dominates} and \ref{A:cardinality} be valid.
If $u\in\As$ and
\begin{equation*}
\osc_{\grid_j}(f)_{-1} \le \omega \, \est_{\grid_j}(u_j,f),
\end{equation*}
then the cardinality $N_j(u)$ of $\marked_j$ satisfies
\[
N_j(u) \lesssim |u|_{\As}^{\frac{1}{s}} |u-u_j|_{H^1_0(\Omega)}^{-\frac{1}{s}}
\quad\forall \, j\ge0.
\]
\end{proposition}
\begin{proof}
Let $C_\textrm{cea}:=\sqrt{\frac{C_\B}{c_\B}}$ and let $\delta = \mu \frac{C_L}{\CCea} \eta_j(u_j)$ with $\mu\le\frac12$ the quasi-monotonicity constant in \eqref{E:quasi-monotonicity-coercive}.
We invoke \eqref{E:equiv-to-As} for $u\in\As$ to find
a mesh $\grid_\delta\in\grids$ and a Galerkin solution $u_\delta\in\V_{\grid_\delta}$  so that
\[
|u-u_{\grid_\delta}|_{H^1_0(\Omega)} \le \delta, \quad \#\grid_\delta \lesssim |u|_{\As}^{\frac{1}{s}}
\delta^{-\frac{1}{s}}.
\]
Since $\grid_\delta$ may be totally unrelated to $\grid_j$, we introduce the overlay
$\grid_*=\grid_j\oplus\grid_\delta$. We exploit that $\grid_*\ge\grid_\delta$, hence the
space nestedness $\V_{\grid_\delta}\subset \V_{\grid_*}$, along with the property that the 
Galerkin solution $u_{\grid_*}\in\V_{\grid_*}$ minimizes the energy error in $\V_{\grid_*}$
\[
\eta_{\grid_*}(u_{\grid_*}) \le \frac{1}{C_L} |u-u_{\grid_*}|_{H^1_0(\Omega)}
\le \frac{\CCea}{C_L} |u-u_{\grid_\delta}|_{H^1_0(\Omega)} \le \frac{\CCea}{C_L}\delta
= \mu \, \eta_j(u_j),
\]
because of the lower bound in \eqref{E:aposteriori-bounds-H1} and \eqref{E:quasi-monotonicity-coercive}.
Therefore, Lemma \ref{L:dorfler-no-osc} (D\"orfler marking) implies that the refined set
$\refined[] = \grid\backslash\grid_*$ satisfies D\"orfler marking with parameter $\theta<\theta_0$.
Since $\MARK$ delivers a minimal set $\marked_j$ with this property, according to
Assumption \ref{A:cardinality}, we deduce
\[
N_j(u) = \#\marked_j \le \#\refined[] \le \#\grid_* - \#\grid \le \#\grid_\eps - \#\grid_0 \lesssim
|u|_{\As}^{\frac{1}{s}} \delta^{-\frac{1}{s}},
\]
where we have used Lemma \ref{nv-L:overlay} (mesh overlay). The assertion
follows from $\osc_j(f)_{-1} \le \omega \est_j(u_j,f)$ in Assumption~\ref{A:PDE-dominates}
and the upper bound in \eqref{E:aposteriori-bounds-H1}
\[
|u-u_j|_{H^1_0(\Omega)} \le C_U \est_j(u_j,f) \le \frac{C_U}{\sqrt{1-\omega^2}} \eta_j(u_j),
\]
and completes the proof.
\end{proof}
 
We next prove rate optimality of the one-step $\AFEM$s of Algorithm \ref{A:GALERKIN} and
Algorithm \ref{A:one-step-switch}. To this end, we need an additional assumption.

\begin{assumption}[initial labeling]\label{A:initial-labeling}
\index{Assumptions!Initial labeling}
If the initial mesh $\grid_0$ is made of simplices, then let the initial labeling \eqref{nv-initial-label}
for $d=2$, or that of \cite[Section 4]{Stevenson:08} for $d>2$, be valid.
\end{assumption}

This assumption ensures %guarantees
the validity of Theorem \ref{nv-T:complexity-refine}
(complexity of $\REFINE$): if
$\marked_j\subset\mesh_j$ is a set of marked elements for a sequence $\{\grid_j\}_{j=0}^{k-1}$
of consecutive refinements of $\grid_0$, then the
cardinality of the $k$-th mesh satisfies
\begin{equation}\label{E:BDD}
\#\grid_k - \#\grid_0 \le \Ccompl \sum_{j=0}^{k-1} \#\marked_j,
\end{equation}
with a universal constant $\Ccompl$ depending only on $\grid_0$ and $d$. We always assume that
$\#\grid_k \ge \frac{3}{2}\#\grid_0$, whence $\#\grid_k-\#\grid_0\ge \frac{1}{3} \#\grid_k$
and, if $\wtCcompl=3 \Ccompl$, \eqref{E:BDD} reads instead
\begin{equation}\label{E:BDD-2}
\index{Constants!$\Ccompl2$@$\wtCcompl$: modified complexity of \REFINE constant}
\#\grid_k \le \wtCcompl \sum_{j=0}^{k-1} \#\marked_j,
\end{equation}

\begin{theorem}[rate optimality of one-step $\AFEM$s]\label{T:optimality-one-step}
For Algorithms \ref{A:GALERKIN} (with $\grid=\grid_0$) and \ref{A:one-step-switch},
let Assumptions \ref{A:theta-no-osc}, \ref{A:cardinality} and
\ref{A:initial-labeling} be valid, and in addition let the parameter $\omega>0$ satisfy Assumption
\ref{A:PDE-dominates} for Algorithm \ref{A:one-step-switch}. If $u\in\As$, then both
one-step {\rm AFEM}s give rise to sequences $\{\grid_k,\V_k,u_k\}_{k=0}^\infty$ such that
\begin{equation}\label{E:optimality-one-step}
|u-u_k|_{H^1_0(\Omega)} \lesssim |u|_{\As} (\#\grid_k)^{-s}.
\end{equation}
\end{theorem}
\begin{proof}
We consider first Algorithm \ref{A:GALERKIN}, for which the forcing $f\in\F_{\grid_0}$ is
discrete, whence $\osc_j(f)_{-1} = 0$ for all $j\ge 0$ and $\omega=\sigma=0$ in Assumption \ref{A:PDE-dominates}.
In view of \eqref{E:BDD-2}, we apply Proposition \ref{P:cardinality-marked} (cardinality of $\marked_j$)
to infer that
\[
\#\grid_k \le \wtCcompl \sum_{j=0}^{k-1} \#\marked_j \lesssim |u|_{\As}^{\frac{1}{s}}
\sum_{j=0}^{k-1} |u-u_j|_{H^1_0(\Omega)}^{-\frac{1}{s}}.
\]
We now recall the inequality $|u-u_k|_{H^1_0(\Omega)} \le C_* \alpha^{k-j} |u-u_j|_{H^1_0(\Omega)}$ from
Corollary \ref{C:linear-convergence} (linear convergence), and replace the sum above by
\[
\sum_{j=0}^{k-1} |u-u_j|_{H^1_0(\Omega)}^{-\frac{1}{s}} \le |u-u_k|_{H^1_0(\Omega)}^{-\frac{1}{s}}
\sum_{j=0}^k \alpha^{\frac{k-j}{s}} \le \frac{\alpha^{\frac{1}{s}}}{1-\alpha^{\frac{1}{s}}}
|u-u_k|_{H^1_0(\Omega)}^{-\frac{1}{s}}
\]
because $0<\alpha<1$ and the geometric series is summable.

We next deal with Algorithm \ref{A:one-step-switch}. If the algorithm calls $\MARK$, then
$\osc_j(f)_{-1} \le \omega \est_j(u_j,f)$ and 
the number of marked elements $N_j(u)$ obeys Proposition \ref{P:cardinality-marked}
\[
N_j(u) \lesssim |u|_{\As}^{\frac{1}{s}} |u-u_j|_{H^1_0(\Omega)}^{-\frac{1}{s}}.
\]
Instead, if the algorithm calls $\DATA$, then $\osc_j(f)_{-1} > \sigma_j=\omega\est_j(u_j,f)$ and $\DATA$
returns a mesh $\grid_{j+1}$ and reduces the oscillation $\osc_{j+1}(f)_{-1} \le \sigma_j$ with optimal
complexity. To quantify the cost, we recall that $u\in\As$ yields $f\in\F_s$ according to
Lemma \ref{L:approx-regularity} (relation between approximation classes) and $|f|_{\F_s}\lesssim |u|_{\As}$.
Therefore, the number of marked elements $N_j(f)$ to reduce $\osc_j(f)_{-1}$ to tolerance $\sigma_j$
satisfies
\[
N_j(f) \lesssim |f|_{\F_s}^{\frac{1}{s}} \sigma_j^{-\frac{1}{s}}
\lesssim |u|_{\As}^{\frac{1}{s}} \est_j(u_j,f)^{-\frac{1}{s}}
\lesssim |u|_{\As}^{\frac{1}{s}} |u-u_j|_{H^1_0(\Omega)}^{-\frac{1}{s}},
\]
because of \eqref{E:aposteriori-bounds-H1}. It thus remains to sum over $j$, apply again
\eqref{E:BDD-2}
\[
\#\grid_k \le \wtCcompl \sum_{j=0}^{k-1} \left( N_j(u) + N_j(f) \right)
\lesssim |u|_{\mathbb A_s}^{1/s} \sum_{j=0}^{k-1}|u-u_j|_{H^1_0(\Omega)}^{-\frac{1}{s}},
\]
and finally argue as before with the help of \AB{Corollary \ref{C:linear-convergence-error-sw} (linear convergence of error)}.
\end{proof}

%%%%%%%%%%%%%%%%%%%%%%%%%%%%%%%%%%%%%%%%%%%%%%%%%%%%%%%%%%%%%%%%%%%%%%%%%%%%%%%%
\subsection{Rate-optimality of two-Step {\rm \AFEM}}\label{S:optimality-AVEM-TS}
%%%%%%%%%%%%%%%%%%%%%%%%%%%%%%%%%%%%%%%%%%%%%%%%%%%%%%%%%%%%%%%%%%%%%%%%%%%%%%%%
%

The output of $[\wh{\grid}_k,\wh{\data}_k] = \DATA \, (\grid_k,\data,\omega\eps_k)$,
in the $k$-step of $\AFEMTS$ (Algorithm~\ref{algo:AFEMTS}), is fed to $[\grid_{k+1},u_{k+1}]=\GALERKIN \, (\wh{\grid}_k,\wh{\data}_k,\eps_k)$, which in turn iterates $J_k$ times.
We denote by $(\grid_{k,j}, \marked{}_{k,j}, u_{\grid_{k,j}})$ the triplets of grids, marked sets and discrete solutions computed within $\GALERKIN \, (\wh{\grid}_k,\wh{\data}_k,\eps_k)$ for $0\le j < J_k$. We further assume that
\[
\wh{\eps}_k := \eta_{\grid_{k,0}} (u_{\grid_{k,0}},\wh{\data}_k) > \eps_k
\]
for otherwise the module $\GALERKIN$ is skipped. In view of the lower a posteriori error estimate in
\eqref{E:apost-discrete-data} for discrete data $\wh{\data}_k$, we infer that
\[
|\wh{u}_k - u_{\grid_{k,0}}|_{H^1_0(\Omega)} \ge C_L \wh{\eps}_k > C_L\eps_k,
\]
where $\wh{u}_k \in H^1_0(\Omega)$ is the exact solution of \eqref{E:perturbed-weak-form}
with approximate data $\wh{\data}_k$. The module $\DATA$ guarantees \eqref{E:bound-data} and \eqref{E:error-DATA}, namely
\begin{equation}\label{E:omega-epsk}
\|\data - \wh{\data}_k \|_{D(\Omega)} \le \Cdata\omega\eps_k
\quad\Rightarrow\quad
| u - \wh{u}_k |_{H^1_0(\Omega)} \le C_D \omega\eps_k,
\end{equation}
where $u=u(\data)\in H^1_0(\Omega)$ is the exact solution of \eqref{weak-form}.
We see that the parameter $\omega$ controls the discrepancy between $u$ and $\wh{u}_k=\wh{u}_k(\data_k)$ relative to $\eps_k$. We now make an assumption on the appropriate size of $\omega$, which replaces Assumption \ref{A:PDE-dominates}
for $\AFEMSW$.

\begin{assumption}[size of $\omega$]\label{A:size-omega}
\index{Assumptions!Size of $\omega$}
The parameter $\omega$ in $\AFEMTS$ satisfies $\omega \in (0,\omega_0]$ where
$\omega_0 := \frac{\mu C_L}{2 C_D \CCea}$ with $C_\textrm{cea}$ as in \eqref{E:quasi-monotonicity-coercive} and the parameter $0<\mu\le\frac12$ appears in Lemma \ref{L:dorfler-no-osc}
(D\"orfler marking).
\end{assumption}

Consequently, if Assumption \ref{A:size-omega} is valid then \eqref{E:omega-epsk} yields for
$\omega\le\omega_0$
\begin{equation}\label{E:def-omega}
|u - \wh{u}_k|_{H^1_0(\Omega)} \le \frac{\mu C_L}{2\CCea} \eps_k.
\end{equation}

\begin{corollary}[cardinality of marked sets]\label{C:cardinality-no-osc}
Let Assumptions \ref{A:theta-no-osc}, \ref{A:cardinality}, and \ref{A:size-omega} hold. If
$u\in\As\big( H^1_0(\Omega);\grid_0 \big)$ and $\wh{\eps}_k>\eps_k$, then {\rm \GALERKIN} is called and there
exists a constant $C_0$ such that for all $0 \le j < J_k$
\begin{equation}\label{E:card-Mk-no-osc}
\#\marked_{k,j} \le C_0 |u|_{\As}^{1/s}\, |u-u_{\grid_{k,j}}|_{H^1_0(\Omega)}^{-1/s},
\end{equation}
and
\begin{equation}\label{E:card-Mk-no-osc2}
\#\marked_{k,j} \le C_0 |u|_{\As}^{1/s}\, \eps_k^{-1/s}.
\end{equation}
\end{corollary}
\begin{proof}
We argue as in Proposition \ref{P:cardinality-marked}.
Fix $0\leq j < J_k$ and set 
\[
\delta:= \mu \frac{C_L}{\CCea} \eta_{\grid_{k,j}}(u_{\grid_{k,j}})
\quad\Rightarrow\quad
\delta \geq \mu \frac{C_L}{\CCea} \eps_k. 
\]
Since $|u-\widehat{u}_k|_{H^1_0(\Omega)} \leq \frac{\delta}{2}$, by virtue of \eqref{E:def-omega},
we deduce that $\widehat{u}_k $ is an
$\delta$-approximation of order $s$ to $u$ according to Lemma \ref{L:eps-approx}
($\eps$-approximatiom of order $s$).  Therefore, there exists an admissible mesh
$\mathcal{T}_\delta\in\grids$ such that 
\[
|\wh{u}_k - u_{\grid_\delta}|_{H^1_0(\Omega)} \leq \delta, \qquad 
\#\mathcal{T}_\delta \lesssim | u |_{\As}^{\frac{1}{s}}\delta^{-\frac{1}{s}},
\]
and we proceed exactly as in Proposition \ref{P:cardinality-marked} to show that
\[
\# \marked_{k,j} \lesssim | u |_{\As}^{\frac{1}{s}} \delta^{-\frac{1}{s}}
\approx | u |_{\As}^{\frac{1}{s}} |u-u_{\grid_{k,j}}|_{H^1_0(\Omega)}^{-1/s}
\lesssim | u |_{\As}^{\frac{1}{s}} \eps_k^{-\frac{1}{s}}.
\]
This concludes the proof.
\end{proof}

\begin{corollary}[quasi-optimality of $\GALERKIN$]\label{C:quasi-optimality}
Let Assumptions \ref{A:approx-u}, \ref{A:theta-no-osc}, \ref{A:cardinality}, and \ref{A:size-omega} be
valid. Then, the number of marked elements $N_k(u)$ within the $k$-th call to {\rm \GALERKIN} satisfies
$$ 
N_k(u) \leq J C_0 | u |_{\As}^{\frac{1}{s}}\varepsilon_k^{-\frac{1}{s}},
$$
where $J\ge J_k$ is a uniform upper bound for the number of iterations of {\rm\GALERKIN}.
\end{corollary}
\begin{proof}
Use $N_k(u)=\sum_{j=0}^{J_k-1} \#\mathcal{M}_{k,j}$ and combine Corollary \ref{C:cardinality-no-osc}
(cardinality of marked sets) and Proposition \ref{P:cost-galerkin} (computational cost of $\GALERKIN$).
\end{proof}

We finally address the rate-optimality 
of the two-step algorithm $\AFEMTS$, by proving the announced bound \eqref{E:optimality-bound}.

\begin{theorem}[rate-optimality of $\AFEMTS$]\label{T:optimality-AFEM}
 Let Assumptions \ref{A:approx-u} (approximability of $u$),  \ref{A:approx-data} (approximability of data),
 \ref{A:optim-data} (quasi-optimality of the module {\rm \DATA}), 
 \ref{A:theta-no-osc} (marking parameter), \ref{A:size-omega} (size of $\omega$),
 and \ref{A:initial-labeling} (initial labeling)
 be valid. Then, {\rm \AFEMTS} gives rise to a sequence $\big(\grid_k,\V_{\grid_k},u_{\grid_k}\big)_{k=0}^{K}$
 such that
  \begin{equation*}
  |u -u_k|_{H^1_0(\Omega)} \leq C(u,\data) \big(\#\mesh_k\big)^{-s}  \quad 1\leq k \leq K,
  \end{equation*}
  where $0 < s = \min\{s_u, s_\data\} = \min\{s_u, s_A, s_c, s_f\} \leq\frac{n}{d}$ and
 \begin{equation*}
 C(u,\data) = C_* \Big( |u|_{\A_{s_u}}^{\frac{1}{s_u}} + |\bA|_{\ACA_{s_A}}^{\frac{1}{s_A}} + |c|_{\ACc_{s_c}}^{\frac{1}{s_c}} +
 |f|_{\ACf_{s_f}}^{\frac{1}{s_f}} \Big)^s
\end{equation*}
with constant $C_*>0$ independent of $u$ and $\data$.
  \end{theorem}
\begin{proof}
In view of Assumption \ref{A:approx-u}, 
Corollary \ref{C:quasi-optimality} implies that the number of marked elements $N_k(u)$ within the
$(k+1)$-th call to $\GALERKIN$ satisfies 
$$
N_k(u)\le C_3 | u |^{\frac{1}{s_u}}_{\A_{s_u}}\varepsilon_k^{-\frac{1}{s_u}}
$$ 
with $s_u\leq \frac{n}{d}$ and $C_3>0$ a suitable constant.  Moreover,  by Assumption  \ref{A:optim-data} the number of marked elements $N_k(\data)$ within the $(k+1)$-th call to $\DATA$ satisfies 
$$
N_k(\data)\le C_3 | \data |_{\mathbb{A}_{s_\data}}^{\frac{1}{s_\data}} \, \varepsilon_k^{-\frac{1}{s_\data}} 
$$ 
with $s_\data \leq \frac{n}{d}$. The total number of marked elements in the $(k+1)$-th loop of $\AFEMTS$ is thus
$$
N_k(\data)+N_k(u)\leq C_3 \Big( | u |_{{\A}_{s_u}}^{\frac{1}{s_u}} + | \data |_{\mathbb{A}_{s_\data}}^{\frac{1}{s_\data}}
\Big)\, \varepsilon_k^{-\frac{1}{s}} \,.
$$
Upon termination,  $\DATA$ and $\GALERKIN$ give 
\begin{equation*}
\begin{split}
|u-\widehat{u}_k|_{H^1_0(\Omega)}  & \le \frac{\mu C_L}{2\CCea} \eps_k \le \frac{C_L}{4\CCea} \eps_k, \\[5pt]
|\widehat{u}_k-u_{k+1}|_{H^1_0(\Omega)}
& \leq C_U \eta_{\grid_{k+1}}(u_{k+1},\wh{\data}_k)\leq C_U \eps_k \,,
\end{split}
\end{equation*}
because of \eqref{E:def-omega}, \eqref{E:apost-discrete-data} and the fact that $\mu<\frac12$.
This implies by triangle inequality
\begin{equation*}
| u - u_{k+1} |_{H^1_0(\Omega)} \leq \Big( \frac{C_L}{4\CCea} + C_U \Big) \eps_k = C_4 \eps_k.
\end{equation*}
Therefore, applying Theorem \ref{nv-T:complexity-refine} (complexity of $\REFINE$), the total amount of elements created by $k+1$ iterations within $\AFEMTS$, besides those in $\mathcal{T}_0$,  obeys the expression
\begin{eqnarray*}
\#\mathcal{T}_{k+1}
\leq \wtCcompl \sum_{j=0}^{k} \big( N_j(\data) + N_j(u) \big) \leq \wtCcompl C_3 
\Big(| u |_{{\A}_{s_u}}^{\frac{1}{s_u}} + | \data |_{\mathbb{A}_{s_\data}}^{\frac{1}{s_\data}}
\Big)\sum_{j=0}^{k}\varepsilon_j^{-\frac{1}{s}}.
\end{eqnarray*}
according to \eqref{E:BDD-2}.
Since $\varepsilon_j=2^{-j}\eps_0$ and
$
\sum_{j=0}^{k-1}(2^{-\frac{1}{s}})^j \leq \frac{1}{1-2^{-1/s}}
$
we obtain
\begin{equation*}\label{aux:2:opt-vem}
\#\mathcal{T}_{k+1} \leq C
\Big(| u |_{{\A}_{s_u}}^{\frac{1}{s_u}} + | \data |_{\mathbb{A}_{s_\data}}^{\frac{1}{s_\data}}
\Big)\, \varepsilon_k^{-\frac{1}{s}}
\end{equation*}
with $C=\frac{\wtCcompl C_3 \eps_0}{1-2^{-1/s}}$ provided $\#\grid_{k+1}\ge\frac32 \#\grid_0$.
This together with $| u - u_{k+1} |_{H^1_0(\Omega)} \le C_4 \eps_k$
gives the asserted estimate after $1 \le k+1 \le K$ loops.
\end{proof}

\begin{remark}\label{R:thresholds} 
The thresholds $\theta_0, \omega_0$ play no role in Proposition \ref{P:convergence-AFEM} (convergence of $\AFEMTS$) but are critical
in Theorem \ref{T:optimality-AFEM} (rate-optimality of $\AFEMTS$). The former takes care of the
discrepancy between error and estimator
\cite{Stevenson:07,CaKrNoSi:08,NoSiVe:09,bonito2010quasi,NochettoVeeser:2012}.
The latter guarantees that the perturbation
error \eqref{E:omega-epsk} is much smaller than $\varepsilon_k$ and enables $\GALERKIN$ to learn the
regularity of $u$ from $\widehat{u}_{\widehat{\mesh}_k}$ \cite{Stevenson:07,BonitoDeVoreNochetto:2013}.
\end{remark}

 \begin{remark}\label{R:s-sharp} {\rm
We claim that the convergence rate $s = \min\{s_u,s_\data\}$ cannot be improved to $s_u$ (the optimal
rate for approximations of $u\in\A_{s_u}\big(H^1_0(\Omega);\grid_0\big)$) when $s_\data < s_u$
by any algorithm that uses approximations $\wh\data = (\wh{\bA}, \wh{c}, \wh{f})$ of data $\data = (\bA, c, f)$.
In fact, given any $\delta>0$ consider the ball
\begin{equation}\label{E:ball-data}
B(\data,\delta) := \big\{ \wh{\data}\in\D: \|\data-\wh{\data}\|_{D(\Omega)} \le \delta \big\},
\end{equation}
where $D(\Omega)$ is defined in \eqref{E:space-DO-redefine}. If $u,\wh{u}\in H^1_0(\Omega)$ are the exact
solutions for data $\data,\wh{\data}$, then there are constants $0< c_* \le C_*$ such that
\[
c_* \delta \le \sup_{\wh{\data}\in B(\data, \delta)} |u-\wh{u} |_{H^1_0(\Omega)} \le C_* \delta.
\]
}
\end{remark}
The rightmost inequality is a consequence of Lemma \ref{L:perturbation} (continuous dependence on data). 
For the leftmost inequality, consider first a perturbation $\wh{f} = (1+\delta) f$ of the source term with coefficients $(\wh{\bA}, \wh{c})=(\bA,c)$, whence $\|\data-\wh{\data}\|_{D(\Omega)}=\delta$. Proceeding as in \eqref{nsv-enorm-Vnorm}, the coercivity and continuity of the bilinear form $\B$ imply
\[
c_\B \, |u-\wh{u} |_{H^1_0(\Omega)} \le \|f-\wh{f}\|_{H^{-1}(\Omega)} = \delta \le C_\B \, |u-\wh{u} |_{H^1_0(\Omega)}.
\]
On the other hand, if $\wh{f}=f$ and $(\wh{\bA},\wh{c}) = \frac{1}{\alpha} (\bA,c)$ with
$\alpha = 1+\frac{\delta}{\|\data\|_{D(\Omega)}}$, then
\[
\|\data-\wh{\data}\|_{D(\Omega)} < \delta,
\qquad
|u-\wh{u}|_{H^1_0(\Omega)} = \frac{|u|_{H^1_0(\Omega)}}{\|\data\|_{D(\Omega)}} \delta
\ge \frac{\|f\|_{H^{-1}(\Omega)}}{C_\B \, \|\data\|_{D(\Omega)}} \delta.
\]
This argument takes care of the multiplicative nature of $(\bA,c)$ in \eqref{strong-form}, which makes
$\wh{u}=\alpha u$, and proves our claim.

%%%%%%%%%%%%%%%%%%%%%%%%%%%%%%%%%%%%%%%%%%%%%%%%%%%%%%%%%%%%%%%%%%%%%%%%%%%%%%%%
\subsection{Rate-optimality of {\rm \AFEM} with other boundary conditions}\label{S:rate-other-bc}
%%%%%%%%%%%%%%%%%%%%%%%%%%%%%%%%%%%%%%%%%%%%%%%%%%%%%%%%%%%%%%%%%%%%%%%%%%%%%%%%
%
The key ingredient for rate-optimality of $\AFEM$, regardless of boundary conditions, is the validity
of Lemma \ref{L:dorfler-no-osc} (D\"orfler marking). This lemma provides a bridge between FEM meshes 
and optimal meshes and, in turn, hinges on \AB{three} properties of the PDE residual estimator
\rhn{$\eta_\grid(u_\grid,f)$: Theorem \ref{T:ubd-corr} (upper bound for corrections),
Lemma \ref{L:lipschitz} (Lipschitz property of the estimator), and Lemma~\ref{L:est-forcing} (estimator dependence on discrete forcing) to account for the possible change in the discrete
forcing $P_\grid f$.} Since  their 
proofs are insensitive to boundary conditions, because they do not alter the structure  of
\rhn{$\eta_\grid(u_\grid,f)$}, we conclude their validity as well as for Robin,
Neumann and non-homogeneous Dirichlet conditions.

Therefore, our three $\AFEM$s based on D\"orfler marking deliver the same asymptotic 
convergence rates associated with the approximations classes $\A_{s_u}$ for the solution
$u\in H^1(\Omega)$ and $\mathbb{A}_{s_\data}$ for data $\data = (\bA,c,p,\ell)$ for Robin and Neumann boundary conditions and $\data = (\bA,c,f,g)$ for Dirichlet boundary conditions. We need three
new approximation classes for $p\in \mathbb{P}_{s_p} (L^\infty(\partial\Omega);\grid_0)$,
$\ell\in \mathbb{L}_{s_\ell} (H^1(\Omega)^*;\grid_0)$ for Robin or Neumann conditions
and $g\in \mathbb{G}_{s_g}(H^{1/2}(\partial\Omega);\grid_0)$.

If coefficients $(\bA,c,p)$ are discrete for the Robin condition, then Lemma \ref{L:approx-regularity} (relation between approximation classes) extends and yields
\[
\langle \ell - \wh{\ell},v\rangle = \B [u-\wh{u},v]
\quad\Rightarrow\quad
\|\ell - \wh{\ell}\|_{\V^*} \le C \|u-\wh{u}\|_{\V},
\]
with $\V = H^1(\Omega)$ and $\B=\wh \B, \wh \ell$ given in \eqref{E:B-robin}. This in turn implies $|\ell|_{\mathbb{L}_{s_\ell}} \le C |u|_{\A_{s_u}}$, $s_\ell = s_u$ and the validity of Theorem \ref{T:optimality-one-step} (rate optimality of one-step $\AFEM$s). In \AFEMTS, \DATA approximates $\ell$ along with the other data and Theorem \ref{T:optimality-AFEM} (rate-optimality of $\AFEMTS$) is also valid for Robin and Neumann boundary conditions. We do not explore this matter any longer.

For non-homogeneous Dirichlet boundary conditions the analysis is simpler. If $g$ is discrete, then there is no difference with $g=0$. If not, we note that the solution map $g \mapsto u$ (all other data being fixed) is affine and that the error and augmented total estimator $\est_\grid(u_\grid,f,g):=\est_\grid(u_\grid,f)+\osc_{\grid}(g)_{1/2}$ are equivalent (Theorem~\ref{T:Dirichlet}). This indicates that the role of $g$ is similar to the role of $f$. Therefore, it suffices to replace $\est_\grid(u_\grid,f)$ by $\est_\grid(u_\grid,f,g)$ and $\osc_\grid(f)_{-1}$ by $\osc_\grid(f)_{-1}+\osc_\grid(g)_{1/2}$ in \AFEMSW. For \AFEMTS, the approximation of $g$ is handled by \DATA along with the other data. Hence, we again conclude that Theorems \ref{T:optimality-one-step} and \ref{T:optimality-AFEM} extend to non-vanishing Dirichlet conditions.

%%%%%%%%%%%%%%%%%%%%%%%%%%%%%%%%%%%%%%%%%%%%%%%%%%%%%%%%%%%%%%%%%%%%%%%%%%%%%%%%
\subsection{Rate-optimality of {\rm \AFEM} driven by alternative estimators}\label{S:rate-alternative}
%%%%%%%%%%%%%%%%%%%%%%%%%%%%%%%%%%%%%%%%%%%%%%%%%%%%%%%%%%%%%%%%%%%%%%%%%%%%%%%%
%
We recall the notation $\zeta_\grid(u_\grid)$ of Section \ref{S:conv-non-residual-est} for any of the
three alternative estimators in Section \ref{S:other_estimators} and the crucial local properties
\eqref{E:star-equiv-res} and \eqref{E:star-equiv-est}.

As already alluded to in Section \ref{S:rate-other-bc}, the key instrument for rate-optimality
is Lemma \ref{L:dorfler-no-osc} (D\"orfler marking). We now check the validity of
its \AB{three} main pillars: Theorem \ref{T:ubd-corr} (upper bound for corrections) and 
Lemma \ref{L:lipschitz} (Lipschitz property of the estimator) \rhn{and Lemma~\ref{L:est-forcing} (estimator dependence on discrete forcing) to account for possible change in the discrete forcing $P_\grid f$. It turns out that
if they were valid for $\zeta_\grid(u_\grid) = \zeta_\grid(u_\grid,f)$,}
then statements about rates of convergence similar to those for $\eta_\grid(u_\grid)$
would follow for $\zeta_\grid(u_\grid)$.

\begin{lemma}[localized discrete upper bound]\label{L:localized-upper-zeta}
Let $\grid,\grid_*\in\grids$ and $\grid_*$ be a refinement of $\grid$. Let the coefficients $(\wh\bA,\wh{c})$ be discrete over $\grid$ \AB{and $f\in H^{-1}(\Omega)$}. Then the error between the
corresponding Galerkin solutions $u_\grid\in\V_\grid$ and $u_{\grid_*} \in \V_{\grid_*}$ is bounded
by the indicator in the refined set $\overline{\mathcal{R}}$ plus data oscillation
\[
|u_\grid - u_{\grid_*}|_{H^1_0(\Omega)} \le
\overline{C}_U \Big(\big(\CLequiv\big)^{-1} \zeta_\grid(u_\grid,\overline{\mathcal{R}})^2
+\osc_\grid(f)_{-1}^2 \Big)^{\frac12},
\]
where $\overline{\mathcal{R}} := \big\{z\in\vertices: T \in \grid\backslash\grid_*, T\subset\omega_z   \big\}$
collects all vertices whose associated stars change from $\grid$ to $\grid_*$.
\end{lemma}
\begin{proof}
It suffices to realize that
$\grid\backslash\grid_* \subset \bigcup\{\omega_z: z\in\overline{\mathcal{R}}\}$, and appeal to
Theorem \ref{T:ubd-corr} and \eqref{E:star-equiv-est} to arrive at
\begin{align*}
|u_\grid - u_{\grid_*}|_{H^1_0(\Omega)}^2 &\le \overline{C}_U^2 \Big(\eta_\grid(u_\grid,\mathcal{R})^2
+\osc_\grid(f, \mathcal{R})_{-1}^2 \Big)
\\
& \le \overline{C}_U^2 \Big(\eta_\grid(u_\grid,\overline{\mathcal{R}})^2
+\osc_\grid(f)_{-1}^2 \Big)
\\
& \le \overline{C}_U^2 \Big(\big(\CLequiv\big)^{-2} \zeta_\grid(u_\grid,\overline{\mathcal{R}})^2
+\osc_\grid(f)_{-1}^2 \Big).
\end{align*}
This is the desired estimate.
\end{proof}  

\begin{lemma}[Lipschitz property of the estimator]\label{E:lipschitz-zeta}
    Let the coefficients $(\wh\bA,\wh{c})$ be discrete over $\grid$. There exists $\Clip$ such that
    \[
    \big| \zeta_\grid(v_1)- \zeta_\grid(v_2) \big| \le \Clip |v_1-v_2|_{H^1_0(\Omega)}
    \quad\forall \, v_1,v_2\in \V_\grid.
    \]
\end{lemma}
\begin{proof}
We resort to the star equivalence \eqref{E:star-equiv-res} between discrete residual and estimator.
It thus suffices to derive the Lipschitz property for $\|P_\grid R_\grid(v)\|_{H^{-1}(\omega_z)}$
with respect to $v\in\V_\grid$ for all $z\in\vertices$. Since 
$P_\grid R_\grid(v) = P_\grid f - \wh{\B}[v,\cdot]$, we get
\[
\big\langle P_\grid R_\grid(v_1) - P_\grid R_\grid(v_2),w \big\rangle = 
\int_{\omega_z} \nabla w \cdot\wh{\bA}\nabla (v_1-v_2) + \wh{c} (v_1-v_2)w 
\]
for all $w\in H^1_0(\omega_z)$. Therefore, Lemma \ref{L:Poincare} (first Poincar\'e inequality) yields
\[
\|P_\grid R_\grid(v)-P_\grid R_\grid(v_2)\|_{H^{-1}(\omega_z)}\le C(\wh{\bA},\wh{c}) 
\|v_1-v_2\|_{H^1(\omega_z)}
\]
where $C(\wh{\bA},\wh{c})$ depends on the $L^\infty$-norms of $(\wh{\bA},\wh{c})$. Finally, 
using the triangle inequality to accumulate over $z\in\vertices$ together with 
\eqref{E:star-equiv-res} gives the assertion. 
\end{proof}

Lemmas \ref{L:localized-upper-zeta} and \ref{E:lipschitz-zeta} lead to Lemma 
\ref{L:dorfler-no-osc} (D\"orfler marking) for $\zeta_\grid(u_\grid)$. If we further
choose a {\it minimal} set $\marked$ of vertices that satisfies D\"orfler property 
\eqref{E:doerfler-property}, then the previous rates of convergence for the 
three algorithms $\GALERKIN$, $\AFEMSW$, and $\AFEMTS$ but now driven by $\zeta_\grid(u_\grid)$
are valid provided $u\in\mathbb{A}_s$, the approximation class in Definition \ref{D:approx-class-u}. 
We do not restate these results.

%%%%%%%%%%%%%%%%%%%%%%%%%%%%%%%%%%%%%%%%%%%%%%%%%%%%%%%%%%%%%%%%%%%%%%%%%%%%%%%%
\subsection{Approximation vs regularity classes}\label{S:Besov}
%%%%%%%%%%%%%%%%%%%%%%%%%%%%%%%%%%%%%%%%%%%%%%%%%%%%%%%%%%%%%%%%%%%%%%%%%%%%%%%%
%
The purpose of this section is to reconcile the notion of approximation classes, discussed above, with that
of regularity classes. We recall the DeVore Diagram of Fig.~\ref{F:devore-diagram} that depicts the 
Sobolev line for the energy space $H^1_0(\Omega)$, namely
%corresponding to nonlinear Sobolev approximation of the energy space $H^1_0(\Omega)$, namely
%
\begin{equation*}
\sob (H^1_0) = \sob (W^s_p)
\quad\Rightarrow\quad
s - \frac{d}{p} = 1 - \frac{d}{2}.
\end{equation*}
The differentiability $s\ge1$ is only limited by the polynomial degree $n$, so $s\in [1,n+1]$. On the other
hand, the integrability $p$ is not restricted to be $p\ge1$ as is customary with Sobolev spaces. For example,
for $d=2$ and $s=n+1$, we get $p=\frac{2}{n+1} < 1$ provided $n\ge2$. Therefore, to take full advantage of
nonlinear approximation theory, we need to abandon the framework of Sobolev spaces $W^s_p(\Omega)$ and
deal with {\it Besov spaces} $B^s_{p,q}(\Omega)$ (frequently denoted $B^s_q(L^p(\Omega))$ or $B^s_q(L_p(\Omega))$ in the literature) with integrability index $p\in (0,\infty]$. The second index
$q \in (0,\infty]$ is useful to characterize special limiting cases; we will provide below a few interesting
examples but take $p=q$ most of the time. At this point, we only mention that when $s$ is non-integer and $1\leq p\leq \infty$, $B^s_{p,p}(\Omega)=W^s_p(\Omega)$ while when $r$ is integer $W^r_p(\Omega)$ for $p\not =2$ is not a Besov space but it is slightly smaller than $B^r_{p,\infty}(\Omega)$. The case $p=2$ is special since $B^s_{2,2}(\Omega)=H^s(\Omega)$ even when $s$ is an integer. 

This section
is devoted to the definition and properties of Besov and Lipschitz spaces, including their close relation to
approximation classes. Our presentation follows closely \cite{BiDaDeVPe:02} for $n=1$ and
\cite{GaspozMorin:14,GaspozMorin:17} for $n\ge1$, but it adds a few new ingredients. Since our results
involve three different type of spaces to account for the particular cases when the differentiability is integer, it is pertinent to introduce the following abstract
space $X^s_p(\Omega)$ with differentiability index $s \in (0, n+1]$ and integrability index $p \in (0,\infty]$
\begin{equation}\label{E:abstract-space}
\index{Functional Spaces!$X^s_p(\Omega)$: abstract functional spaces}
X^s_p(\Omega) :=
\begin{cases}
B^s_{p,p}(\Omega) & \quad  s \in (0, n+1), p \in (0,\infty],
\\
W^{n+1}_p(\Omega) & \quad  s = n+1, p\in [1,\infty],
\\
\textrm{Lip}^{n+1}_p(\Omega) & \quad s = n+1, p\in (0,1).
\end{cases}
\end{equation}
Here $\textrm{Lip}^{s}_p(\Omega)=\textrm{Lip}(s,L^p(\Omega))$, $s\in \mathbb N$, \index{Functional Spaces!$\textrm{Lip}^{n+1}_p(\Omega)$: Lipschitz spaces} are the Lipschitz spaces; see \eqref{E:lipschitz-Lp} below. For $s\in \mathbb N$ and $1< p < \infty$ the Sobolev spaces coincide with the Lipschitz spaces \cite[Theorem 10.55]{Leoni:2009}, i.e.,
\begin{equation}\label{e:Lip-Sob}
\textrm{Lip}^s_p(\Omega) = W^s_p(\Omega), \qquad s\in \mathbb N, \ 1 < p < \infty,
\end{equation}
while for $p=1$ we only have
\begin{equation}\label{e:Lip-Sob2}
 W^s_1(\Omega) \hookrightarrow \textrm{Lip}^s_1(\Omega), \qquad s\in \mathbb N.
\end{equation}

We use the following conventions: $X^s_p(\Omega):=L^p(\Omega)$ for $s=0$; $X^s_p(\Omega;\grid)$
is the space of functions with piecewise regularity $X^s_p$ over $\grid\in\grids$; $X^s_p(\Omega;\R^m)$ 
is the space $X^s_p(\Omega)$ of vector or matrix-valued functions.

We will prove in Section~\ref{S:global-approx} the following crucial approximation
results for functions in $L^q(\Omega)$ by discontinuous piecewise polynomials $\mathbb{S}^{n,-1}_\grid$
of degree $n \ge 1$ over conforming refinements $\grid$ of $\grid_0$. It turns out that this will also
allow us to deal with approximations in $H^1_0(\Omega)$ by continuous piecewise polynomials $\V_\grid$
of degree $n \ge 1$.

\begin{theorem}[regularity yields approximation]\label{T:direct}
Let $q\in[1,\infty], p\in (0,\infty]$, $s\in (0,n+1]$ and a function $g\in L^q(\Omega)$ satisfy $g\in X^s_p(\Omega)$ with
$s-\frac{d}{p} + \frac{d}{q}>0$. Then, there exists a constant $C=C(p,q,s,t,d,\Omega,\grid_0)$
such that
\begin{equation}\label{E:global-besov}
 E_n(g,\Omega)_q
 :=
 \adjustlimits \inf_{\grid\in\grids_N} \inf_{v \in\mathbb{S}^{n,-1}_\grid} |g-v|_{L^q(\Omega)}
 \le
 C |g|_{X^s_p(\Omega)} N^{-\frac{s}{d}}.
\end{equation}
Therefore, $g\in\A_{\frac{s}{d}}=\A_{\frac{s}{d}} \big(L^q(\Omega);\grid_0\big)$ and
\begin{equation}\label{E:approx-reg}
|g|_{\A_{\frac{s}{d}}} \le C |g|_{X^s_p(\Omega)}.
\end{equation}
\end{theorem}

We see that the decay rate $\frac{s}{d}$ in \eqref{E:global-besov} is proportional to
the difference of the differentiability
indices between the space $X^s_p(\Omega)$ and $L^q(\Omega)$ provided the Sobolev numbers
satisfy the relation
\[
\sob\big(X^s_p(\Omega)\big) > \sob\big(L^q(\Omega)\big),
\]
which implies
that the embedding of $X^s_p(\Omega)$ into $L^q(\Omega)$ is compact. The factor $d$ in the
denominator is a manifestation of the so-called {\it curse of dimensionality.}
The limiting case $s = n+1$ entails dealing with
Sobolev spaces $W^s_p(\Omega)$ and Lipschitz spaces $\textrm{Lip}^s_p(\Omega)$ depending on whether $p\geq 1$ or $p<1$.

%----------------------------------------------------------------------------------
\subsubsection{Modulus of smoothness}\label{S:modulus-smoothness}
%----------------------------------------------------------------------------------

%----------------------------------------------------------------------------------
\medskip\noindent
{\em Difference Operators.}
%----------------------------------------------------------------------------------
%
Since we intend to allow $p \in (0,1)$, the underlying functions in $B^s_{p,p}(\Omega)$ might not be locally
integrable, whence they might not be distributions in $\Omega$. Therefore, the notion of weak derivative
does not apply, which in turn has the drawback of being defined for integers and not for fractional
numbers. This leads to the most standard definition of Besov spaces $B^s_{p,q}(\Omega)$
using difference operators, which
only requires integrability in $L^p(\Omega)$ and is valid for any $s>0, p,q\in(0,\infty]$.
Other definitions, which provides equivalent results in the range $1\le p,q\le\infty$ can be found in
\cite{AdamsFournier:03,BerghLofstrom:76}.

Given a bounded Lipschitz domain $\Omega\subset\R^d$, and a vector $h \in\R^d$, we set
\[
\Omega_h := \big\{x\in\Omega: \quad [x,x+h] \subset \Omega  \big\}
\]
where $[x,x+h]$ denotes the closed segment connecting $x$ and $x+h$,
and define the first-order difference operator to be
\begin{equation}\label{E:first-difference}
\Delta_h^1 g(x) = \Delta_h^1 (g,x,\Omega) :=
\begin{cases}
g(x+h) - g(x) & x \in \Omega_h,
\\
0 & \textrm{otherwise}.
\end{cases}
\end{equation}
For $k\in\N$, $k\ge1$, we define the $k$-th difference operator by iteration
\begin{equation}\label{E:k-difference}
\Delta_h^k g(x) := \Delta_h^1 \big( \Delta_h^{k-1}\big)  g(x) \quad x \in \Omega_{kh}
\end{equation}
and observe that it has the explicit form
\[
\Delta_h^k g(x) =
%\left\lbrace \begin{array}{rl}
\begin{cases}
\sum\limits_{j=0}^k (-1)^{k+j}
\Big( \begin{matrix} k \\ j \end{matrix} \Big)  g(x+jh) & [x,x + kh] \subset\Omega,
\\
0 &  \textrm{otherwise}.
%\end{array} \right.
\end{cases}
\]
Note the property
\begin{equation}\label{E:poly-vanish}
p\in\P_k \quad\Rightarrow\quad \Delta_h^{k+1} p = 0 \quad \forall h.
\end{equation}
%

%----------------------------------------------------------------------------------
\medskip\noindent
{\em Smoothness.}
%----------------------------------------------------------------------------------
%
Given $p\in (0,\infty]$ and $t>0$, we define the {\it modulus of smoothness} of order $k$ in $L^p(\Omega)$ to be
\begin{equation}\label{E:modulus-smoothness}
\omega_k(g,t)_p = \omega_k(g,t,\Omega)_p := \sup_{|h|\le t} \| \Delta_h^k g \|_{L^p(\Omega)}.
\end{equation}
We note that if $\omega_k(g,t)_p = o(t^{n+1})$ as $t\to0$, then $g$ is a.e a polynomial in $\P_n$ and
\begin{equation}\label{E:non-degeneracy}
g \notin \P_n \quad\Rightarrow\quad
\omega_k(g,t)_p \ge C t^{n+1} \,\, 0 < t \le 1
\end{equation}
for some $C>0$ \cite[Proposition 7.4]{DeVoreLorentz:1993}. We also observe that the definition 
\eqref{E:modulus-smoothness} only requires $L^p$-integrability of $g$ and leads to the
following celebrated Whitney estimate
of the {\it best approximation error} 
\[
E_n(g,G)_p := \inf_{v\in\P_n} \|g-v\|_{L^p(G)}
\]
of $g$ by polynomials of degree $\le n$ in $G\subset\Omega$
\cite{BiDaDeVPe:02}, \cite[Theorem 1.4]{DekelLeviatan:04} \cite[Lemma 4.4]{GaspozMorin:14,GaspozMorin:17}.

\begin{lemma}[Whitney's lemma]\label{L:whitney}
Let $\grid\in\grids$ be an admissible grid, and let $T\in\grid$ be a generic element.
%Let $G$ be either $G=T$ or $G=\omega_\grid(T)$ the patch of $T$ over $\grid$.
If \ $0<p\le\infty$ and $n\ge 0$, then
\[
E_n(g,T)_p \le C \omega_{n+1}(g,h_T,T)_p \quad\forall g\in L^p(T),
\]
where $C = C(p,n,d,\grid_0)$ but is independent of $g$ and the size of $T$.
\end{lemma}
%

%----------------------------------------------------------------------------------
\subsubsection{Besov spaces}\label{S:besov}
%----------------------------------------------------------------------------------
%
Given $s > 0$ and $0< p,q \le \infty$, the Besov space $B^s_{p,q}(\Omega)$ is the set of all functions
$v\in L^p(\Omega)$ such that the following quantity is finite
\begin{equation}\label{E:besov-norm}
\index{Functional Spaces!$B^s_{p,q}(\Omega)$: Besov spaces}
|v|_{B^s_{p,q}(\Omega)} :=
\begin{cases}
\Big( \int_0^\infty \big[ t^{-s} \omega_{k}(v,t)_p \big]^q \frac{dt}{t}  \Big)^{\frac{1}{q}}
& 0< q < \infty,
\\
\sup_{t>0} \big[ t^{-s} \omega_{k} (v,t)_p \big] & q = \infty,
\end{cases}
\end{equation}
with $k = [s] + 1\in\N$ and $[s]$ stands for the integer part of $s$.
If we split the integral in \eqref{E:besov-norm} for $0<q<\infty$ in dyadic intervals, we obtain the
following equivalent expression for $|v|_{B^s_{p,q}(\Omega)}$
\begin{equation}\label{E:equiv-norm}
|v|_{B^s_{p,q}(\Omega)}^q = \sum_{m\in\Z} \int_{2^{-m-1}}^{2^{-m}} t^{-sq} \omega_{k}(v,t)^q_p \frac{dt}{t}
\approx \sum_{m\in\Z} 2^{msq} \omega_{k}(v, 2^{-m})_p^q;
\end{equation}
here we have used that both $\omega_{k}(v,t)_p$ and $t^{-s}$ are monotone functions of $t$. The hidden
constants depend on $s$ and $q$ but are otherwise independent of $v, k$ and $p$. Note that, with
obvious changes, \eqref{E:equiv-norm} is also valid for $q=\infty$
\begin{equation}\label{E:equiv-norm-q=infty}
|v|_{B^s_{p,\infty}(\Omega)}
\approx \sup_{m\in\Z} \, 2^{ms} \omega_{k}(v, 2^{-m})_p.
\end{equation}

We point out that $|v|_{B^s_{p,q}(\Omega)}$ is a semi-norm for $p,q\ge 1$ and is otherwise a semi-(quasi)norm
in that the triangle inequality is valid up to a constant larger than $1$; note that $|1|_{B^s_{p,q}(\Omega)}=0$.
The quasi-norm of $B^s_{p,q}(\Omega)$ is defined to be
\[
\| v \|_{B^s_{p,q}(\Omega)} := \|v\|_{L^p(\Omega)} + |v|_{B^s_{p,q}(\Omega)}.
\]
If an integer $k'>k$ is chosen in \eqref{E:besov-norm}, then the ensuing quasi-norms
$\| v \|_{B^s_{p,q}(\Omega)}$  are equivalent. This hinges on the {\it Marchaud inequality}
\cite[eq (2.6)]{DeVorePopov:1988}, \cite[Theorems 1 and 3]{Ditzian:1988}
\begin{equation}\label{E:marchaud}
\omega_{k}(v,t)_p \le C t^{k} \left( \|v\|_{L^p(\Omega)} + \left( \int_t^\infty \big(z^{-k} \omega_{k'}(v,z)_p\big)^p
\frac{dz}{z} \right)^{\frac{1}{p}} \right).
\end{equation}
The following lemma characterizes the precise blow-up of $|v|_{B^s_{p,p}(\Omega)}$ as $s\to n+1$.
\begin{lemma}[blow-up of $|v|_{B^s_{p,p}(\Omega)}$]\label{L:blow-up-besov}
{\it Let $s\in(0,n+1), p \in (0,\infty]$. Then
\[
|v|_{B^s_{p,p}(\Omega)} \le \big( p(n+1-s)\big)^{-\frac{1}{p}} \|v\|_{B^{n+1}_{p,p}(\Omega)}
\quad\forall v\in B^{n+1}_{p,p}(\Omega).
\]
}
\end{lemma}
\begin{proof}
We take $p<\infty$ and combine the definition \eqref{E:besov-norm} with \eqref{E:marchaud}, after replacing
the upper limit of integration by $\diam (\Omega) \approx 1$, to write
\begin{align*}
|v|_{B^s_{p,p}(\Omega)}^p & \approx \int_0^1 \Big(t^{-s} \omega_{n+1}(v,t)_p\Big)^p \frac{dt}{t} \lesssim I + II,
\end{align*}
with
\[
 I
 =
 \int_0^1 t^{(n+1-s)p} \int_t^1 \Big(z^{-(n+1)} \omega_{n+2}(v,z)_p\Big)^p \frac{dz}{z} \frac{dt}{t},
\,\,
II = \int_0^1 t^{(n+1-s)p} \|v\|_{L^p(\Omega)}^p \frac{dt}{t}.
\]
Exchanging the order of integration yields
\begin{align*}
 I
 &=
 \int_0^1 \Big(z^{-(n+1)} \omega_{n+2}(v,z)_p\Big)^p \Big(\int_0^z t^{(n+1-s)p-1} dt\Big) \frac{dz}{z}
\\
 &=
 \frac{1}{p(n+1-s)} \int_0^1 \Big(z^{-s} \omega_{n+2}(v,z)_p\Big)^p \frac{dz}{z}
 \le
 \frac{1}{p(n+1-s)} |v|_{B^{n+1}_{p,p}(\Omega)}^p.
\end{align*}
Since $II = \big(p(n+1-s)\big)^{-1} \|v\|_{L^p(\Omega)}^p$, the proof is thus complete.
\end{proof}

The following equivalence between
Sobolev and Besov spaces is valid for fractional differentiability $s$ 
\cite[Proposition 14.40]{Leoni:2009} (see also \cite[\S 6.4.4]{BerghLofstrom:76}, \cite[\S 7.33, \S 7.67]{AdamsFournier:03}):
for all $s\ge0, s\notin\N$ and $p\in[1,\infty]$
\begin{equation}\label{E:sobolev=besov}
B^s_{p,p} (\Omega) = W^s_p(\Omega).
\end{equation}
However, if $s\in\N$ is an integer, then $B^s_{p,q}(\Omega)$ is defined using $k=s+1$ differences whereas
$W^s_p(\Omega)$ involves $s$ weak derivatives in $L^p(\Omega)$ provided $p\in[1,\infty]$. It turns out that
for integer values $s \in\N$ the spaces differ 
\begin{equation}\label{E:sobolevbesov}
B^s_{p,p} (\Omega) \not = W^s_p(\Omega), \quad  p \not = 2,
\end{equation}
except for the exceptional case $p=2$ for which  $B^s_{2,2} (\Omega) = H^s(\Omega)$ \cite{DeVore:98}.

The Besov semi-norm is {\it sub-additive} in the following sense: if $\{T_i\}_{i=1}^N$ is a disjoint collection
of elements $T_i\in\grid$ and $\grid\in\grids$, $p \in (0,\infty]$ and $s>0$, then there exists
a constant $C$ depending on $p,s,d$ and $\grid_0$ but independent of $N$ such that
\begin{equation}\label{E:sub-additive}
\sum_{i=1}^N |v|_{B^s_{p,p}(T_i)}^p \le C |v|_{B^s_{p,p}(\Omega)}^p
\quad\forall v \in B^s_{p,p}(\Omega).
\end{equation}
The localization of Besov norms is more general than \eqref{E:sub-additive}. In fact, if $\omega_\grid(T)$ denotes
the patch of elements in $\grid$ around $T\in\grid$ (first ring), then the following is valid with equivalence
constants depending $p,s,d$ and $\grid_0$ but independent of $N$ \cite[Lemmas 4.3 and 4.4]{BiDaDeVPe:02}
\begin{equation}\label{E:localization-besov}
\sum_{T\in\grid} |v|_{B^s_{p,p}(\omega_\grid(T))}^p \approx C |v|_{B^s_{p,p}(\Omega)}^p
\quad\forall v \in B^s_{p,p}(\Omega).
\end{equation}

The following statements about embeddings between Besov spaces on bounded Lipschitz domains $\Omega$
will turn to be useful in the sequel \cite[\S 3.2.4, \S 3.3.1]{Triebel:10}:
if $0<p\le\infty$, $0<q_1,q_2\le\infty$, and $s_1,s_2,s>0$, then
\begin{equation}\label{E:besov-embeddings}
\begin{aligned}
s_1 > s_2 &\quad\Rightarrow\quad
B^{s_1}_{p,q_1}(\Omega) \hookrightarrow  B^{s_2}_{p,q_1}(\Omega),
\\
q_1 < q_2 &\quad\Rightarrow\quad
B^{s}_{p,q_1}(\Omega) \hookrightarrow  B^{s}_{p,q_2}(\Omega).
\end{aligned}
\end{equation}
Because of the second relation in \eqref{E:besov-embeddings}, statements valid for all second index $q$ are written for the largest space corresponding to $q=\infty$. 
In addition, for all $0< p,q,r \le \infty$ and $s>0$, the {\it discrepancy} between the spaces
$B^{s}_{p,r}(\Omega)$ and $L^q(\Omega)$ is the quantity
\begin{equation}\label{E:discrepancy}
\delta := s-\frac{d}{p} + \frac{d}{q}.
\end{equation}
The discrepancy $\delta$ governs the embedding between these two spaces 
\cite[Theorems 14.29 and 14.32]{Leoni:2009}, \cite{DeVore:98}, namely
\begin{equation}\label{E:embedding}
\delta > 0 \ \Rightarrow\
B^{s}_{p,\infty}(\Omega) \hookrightarrow  L^q(\Omega); \quad \delta =0 \ \Rightarrow\
B^{s}_{p,p}(\Omega) \hookrightarrow  L^q(\Omega), \ q \not = \infty,
\end{equation}
and the embedding is compact when $\delta>0$.
Notice that $\delta=0$
determines the Sobolev embedding line of the DeVore diagram in Fig.~\ref{F:devore-diagram}.

We stress that  when $\delta>0$ the third parameter $r$ in $B^{s}_{p,r}(\Omega)$ plays no role in \eqref{E:embedding} involving the largest space $B^s_{p,\infty}(\Omega)$, see \eqref{E:besov-embeddings}.
However, it turns
out to be useful to quantify regularity in extreme cases. For instance, the
characteristic function $\chi_G$ of a smooth set $G\varsubsetneqq\Omega$ satisfies
\[
\chi_G \in B^{\frac{1}{p}}_{p,\infty}(\Omega) \backslash B^{\frac{1}{p}}_{p,r}(\Omega)
\quad 0 < p, r < \infty.
\]
Moreover, the Lagrange basis functions $\{\phi_z\}_{z\in\nodes}$ of $\V_\grid$ satisfy for
any $0<p\le\infty$ \cite[Proposition 4.7]{GaspozMorin:14}
\begin{equation*}
\omega_{n+1}(\phi_z,t)_p =
\begin{cases}
|\supp \phi_z|^{\frac{d-1-p}{dp}} t^{1+\frac{1}{p}} & \quad 0 < t \le |\supp \phi_z|^{\frac{1}{d}},
\\
|\supp \phi_z| & \quad t > |\supp \phi_z|^{\frac{1}{d}}.
\end{cases}
\end{equation*}
This readily implies that for all $0< s < 1+\frac{1}{p}$ and $0<q<\infty$
\begin{equation}\label{E:Vtau-besov}
\V_\grid \subset B^s_{p,q}(\Omega),\quad
\V_\grid \subset B^{1+\frac{1}{p}}_{p,\infty}(\Omega).
\end{equation}

%----------------------------------------------------------------------------------
\subsubsection{Local approximation}\label{S:local-approx}
%----------------------------------------------------------------------------------
%

We are now in the position to prove a key approximation estimate. In finite element theory it comes
by the name of Bramble-Hilbert Lemma whereas in nonlinear approximation theory it is called
Jackson Theorem. We distinguish between the case $0<s<n+1$ and the limit integral
case $s=n+1$.

\begin{proposition}[Bramble-Hilbert for Besov spaces]\label{P:Bramble-Hilbert-Besov}
Let $\grid\in\grids$ and $T\in\grid$. Let $0<p,q\le\infty$, $0< s < n+1$, and either
$s - \frac{d}{p} + \frac{d}{q} \ge 0, q<\infty$ or $s>\frac{d}{p}, q=\infty$.
Set $r=\infty$ when $s-\frac d p + \frac d q>0$ and $r=p$ otherwise. Then we have
\begin{equation}\label{E:local-best-besov}
\inf_{P\in\P_n} \|v-P\|_{L^q(T)} \le C h_T^{s - \frac{d}{p} + \frac{d}{q}} |v|_{B^s_{p,r}(T)}
\quad\forall \, v\in B^s_{p,r}(T),
\end{equation}
where the constant $C = C(p,q,s,d,n,\grid_0)$ is independent of $v$ and $T$.
\end{proposition}
\begin{proof}
We first point out that we could use $k=n+1 \geq [s]+1$ in the definition of $|v|_{B^s_{p,r}(T)}$ according to \eqref{E:marchaud}. We next proceed in three steps.

\medskip\noindent
\step{1} Suppose first that $T$ is the master element, namely $|T| \approx 1$.
If $P\in\P_n$ is an arbitrary polynomial, using that the discrepancy $\delta=s - \frac{d}{p} + \frac{d}{q}\ge0$
yields
\[
E_n(v,T)_q \le \|v-P\|_{L^q(T)} \lesssim \|v-P\|_{B^s_{p,r}(T)}
= \|v-P\|_{L^p(T)} + |v-P|_{B^s_{p,r}(T)}
\]
due to the embedding \eqref{E:embedding}. Since the definition of $\omega_{n+1}(v,t)_p$ involves $n+1$
differences, we deduce $\Delta_h^{n+1} P = 0$ in view of \eqref{E:poly-vanish}, whence
$|v-P|_{B^s_{p,r}(T)} = |v|_{B^s_{p,r}(T)}$. We now take $P$ to be the best approximation of $v$ in
$L^p(T)$ to derive
\[
E_n(v,T)_q \lesssim E_n(v,T)_p + |v|_{B^s_{p,r}(T)}.
\]

\medskip\noindent
\step{2} We perform a scaling argument from the element $T\in\grid$ to the master element $\wh{T}$.
Let $\wh{x} = |T|^{-1/d} x$ be the change of variables and note that
\begin{align*}
\omega_{n+1}(v,t,T)_p & = \sup_{|h|\le t} \|\Delta_h^{n+1} v\|_{L^p(T)}
\\
& =  |T|^{\frac{1}{p}} \sup_{|h|\le t} \|\Delta_{h|T|^{-1/d}}^{n+1} \wh{v}\|_{L^p(\wh{T})}
= |T|^{\frac{1}{p}} \omega_{n+1}(\wh{v},\wh{t},\wh{T})_p
\end{align*}
with $\wh{t}=t|T|^{-1/d}$, whence
\begin{align*}
 |v|_{B^s_{p,r}(T)}^r
 &=
 \int_0^\infty \Big( t^{-s} \omega_{n+1} (v,t,T)_p \Big)^r \frac{dt}{t}
\\
 &=
 |T|^{\frac{r}{p}-\frac{sr}{d}}
  \int_0^\infty \Big( \wh{t}^{-s} \omega_{n+1} (\wh{v},\wh{t},\wh{T})_p \Big)^r \frac{d\wh{t}}{\wh{t}}
 =
 |T|^{\frac{r}{p}-\frac{sr}{d}} |\wh{v}|_{B^s_{p,r}(\wh{T})}^r.
\end{align*}
Therefore, since $E_n(v,T)_q = |T|^{\frac{1}{q}} E_n(\wh{v},\wh{T})_q$, we obtain
\[
E_n(v,T)_q \lesssim |T|^{\frac{1}{q}-\frac{1}{p}} E_n(v,T)_p +
|T|^{\frac{s}{d} -\frac{1}{p} + \frac{1}{q}} |v|_{B^s_{p,r}(T)}.
\]

\medskip\noindent
\step{3} It remains to estimate $E_n(v,T)_p$ which, in view of Lemma \ref{L:whitney} (Whitney's lemma),
satisfies $E_n(v,T)_ p \le C \omega_{n+1} (v,h_T,T)_p$ with $h_T \approx |T|^{1/d}\approx 2^{-m}$ for some
$m\in\Z$. Since $k=n+1 \geq [s]+1$, invoking the equivalent definition \eqref{E:marchaud} of $|v|_{B^s_{p,r}(T)}$ yields
\[
 E_n(v,T)_p^r
 \lesssim
 \omega_{n+1}(v,2^{-m},T)_p^r
 \lesssim
 h_T^{sr} \sum_{m\in\Z} 2^{msr} \omega_{n+1}(v,2^{-m},T)_p^r
 =
 h_T^{sr} |v|_{B^s_{p,r}(T)}^r.
\]
Inserting this estimate into that of Step 2 gives \eqref{E:local-best-besov}, as asserted.
\end{proof}

We now consider the integer case $s=n+1$. The first thing to notice is that \eqref{E:local-best-besov}
cannot possibly be valid: the definition \eqref{E:besov-norm} requires $k=[s]+1=n+2$, whence
any polynomial $g\in\P_{n+1}\setminus\P_n$ satisfies $E_n(g,T)_p>0$ as well as $\omega_{n+2}(g,T)_p=0$
according to \eqref{E:poly-vanish}. Lemma \ref{L:blow-up-besov} (blow-up of $|v|_{B^s_{p,p}(\Omega)}$)
reveals that replacing the semi-norm $|v|_{B^s_{p,p}(\Omega)}$ by the full norm $\|v\|_{B^{n+1}_{p,p}(\Omega)}$
is not a good idea either. To overcome this problem, we now introduce the space
$\textrm{Lip}^s_p(\Omega) := \textrm{Lip} \big(s,L^p(\Omega)\big)$ of $s$-Lipschitz functions with values in $L^p(\Omega)$, $0<p<\infty$ 
\cite[p.92]{DeVore:98}
\begin{equation}\label{E:lipschitz-Lp}
|g|_{\textrm{Lip}^s_p(\Omega)} := \sup_{t>0} \Big(t^{-s} \omega_{n+1}(g,t,\Omega)_p\Big).
\end{equation}
Comparing with \eqref{E:besov-norm} we realize that $\textrm{Lip}^s_p(\Omega) = B^s_{p,\infty}(\Omega)$ provided 
$s \not \in \mathbb N$ but $\textrm{Lip}^{s}_p(\Omega) \ne B^{s}_{p,\infty}(\Omega)$ when $s\in \mathbb N$.
Moreover,
\begin{equation}\label{E:compact-lipschitz}
 \delta
 =
 s -\frac{d}{p} + \frac{d}{q}
 >
 0, \ s \in \mathbb N 
\quad\Rightarrow\quad
 \textrm{Lip}^{s}_p(\Omega) \hookrightarrow  L^q(\Omega),
\end{equation}
with compact embedding.
If $\delta=0$ and $p\geq 1$, $q \not = \infty$, the above embedding is continuous in view of \eqref{e:Lip-Sob} and \eqref{E:embedding}.

\begin{proposition}[Bramble-Hilbert for Lipschitz spaces]\label{P:Bramble-Hilbert-Lipschitz}
Let $\grid\in\grids$ and $T\in\grid$. If $p \in (0,\infty)$, $q\in(0,\infty]$, $k\geq 0$ integer, and $k+1 - \frac{d}{p} + \frac{d}{q} \geq 0$, with strict inequality when $p<1$ or $q=\infty$, then we have
\begin{equation}\label{E:local-best-besov-n+1}
\inf_{P\in\P_k} \|v-P\|_{L^q(T)} \le C h_T^{k+1 - \frac{d}{p} + \frac{d}{q}} |v|_{\textrm{Lip}^{k+1}_p(T)}
\quad\forall \, v\in \textrm{Lip}^{k+1}_p(T),
\end{equation}
where the constant $C = C(p,q,d,k,\grid_0)$ is independent of $v$ and $T$.
\end{proposition}
\begin{proof}
In view of \eqref{E:compact-lipschitz}, we proceed as in the proof of
Proposition \ref{P:Bramble-Hilbert-Besov}, except for the following change in Step 3. For $h_T\approx 2^{-m}$,
we have instead
\begin{align*}
 E_k(v,T)_p
 &\lesssim
 \omega_{k+1}(v,2^{-m},T)_p
\\
 &\lesssim
 h_T^{k+1} \sup_{m\in\Z} \Big(2^{m(k+1)} \omega_{k+1}(v,2^{-m},T)_p \Big)
 =
 h_T^{k+1} |v|_{\textrm{Lip}^{k+1}_p(T)}.
\end{align*}
This concludes the proof.
\end{proof}

It is instructive to realize that Propositions \ref{P:Bramble-Hilbert-Besov} and \ref{P:Bramble-Hilbert-Lipschitz}
extend to Besov and Lipschitz spaces the usual Bramble-Hilbert lemma for Sobolev spaces
\cite[Lemma 4.3.8]{BrennerScott:08}.

\begin{proposition}[Bramble-Hilbert for Sobolev spaces]\label{P:Bramble-Hilbert-Sobolev}
{\it Let $\grid\in\grids$ and $T\in\grid$. For all $1 \le p, q, \le \infty$ and
$0< s\le n+1$ such that  $s - \frac{d}{p} + \frac{d}{q}\ge0$ with strict inequality when $q=\infty$, then
\begin{equation}\label{E:local-best-sobolev}
\inf_{P\in\P_n} \|v-P\|_{L^q(T)} \le C h_T^{s - \frac{d}{p} + \frac{d}{q}} |v|_{W^s_p(T)}
\quad\forall \, v\in W^s_p(T),
\end{equation}
where the constant $C = C(p,q,s,d,n,\grid_0)$ but is independent of $v$ and $T$.
}
\end{proposition}
\begin{proof}
When $s$ is fractional, $W^s_p(T)=B^s_{p,p}(T)$ in view of \eqref{E:sobolev=besov} and the result follows from Proposition~\ref{P:Bramble-Hilbert-Besov}.
Instead, when $s$ is integral, we invoke \eqref{e:Lip-Sob} and \eqref{e:Lip-Sob2} to deduce the result from Proposition~\ref{P:Bramble-Hilbert-Lipschitz}.
\end{proof}

%----------------------------------------------------------------------------------
\subsubsection{Global approximation: direct estimates}\label{S:global-approx}
%----------------------------------------------------------------------------------
%
We now collect local contributions from Propositions \ref{P:Bramble-Hilbert-Besov},
\ref{P:Bramble-Hilbert-Lipschitz} and \ref{P:Bramble-Hilbert-Sobolev}, depending on the range
of parameters $p,q,s$, to find global
error estimates for the solution $u$ as well as the coefficients $(\bA,c)$ of \eqref{weak-form}.
They are trivial consequences of Theorem \ref{T:direct}, which we prove first. The analysis of the
forcing function $f$ is somewhat different, due to the non-locality of the corresponding norm
$H^{-1}(\Omega)$, and is postponed to Section \ref{S:right-hand-side}.

\medskip\noindent
{\bf Proof of Theorem \ref{T:direct}.} Since the discrepancy $\delta=s - \frac{d}{p} + \frac{d}{q} > 0$,
the embedding $X^s_p(\Omega) \hookrightarrow  L^q(\Omega)$ is compact according to \eqref{E:embedding} and
\eqref{E:compact-lipschitz}. Given
$g \in X^s_p(\Omega)$, we consider the surrogate quantity $e_\grid(g,T) := C h_T^\delta |g|_{X^s_p(T)}$,
which satisfies
\begin{equation*}
E_n(g,T)_q = \inf_{v\in\P_n} \| g - v\|_{L^q(T)} \le e_\grid(g,T) \quad\forall T \in\grid
\end{equation*}
by virtue of Bramble-Hilbert Propositions \ref{P:Bramble-Hilbert-Besov}, \ref{P:Bramble-Hilbert-Lipschitz} and
\ref{P:Bramble-Hilbert-Sobolev}. We finally combine
Proposition \ref{P:abstract_greedy} (abstract greedy) with the subadditivity property
\eqref{E:sub-additive}
to deduce the desired estimate \eqref{E:sigmaN}. The remaining estimate \eqref{E:approx-reg} follows
from the definition of $|g|_{\A_{\frac{s}{d}}}$. This concludes the proof.

Inspection of this proof reveals that our estimate is
stronger than \eqref{E:global-besov}. In fact, we need the weaker regularity
\[
|g|_{X^s_p(\Omega;\grid)}^p = \sum_{T\in\grid} |g|_{X^s_p(T)}^p < \infty,
\]
which allows for piecewise Besov smoothness of $g$ very much in the spirit of \eqref{eq:localized-interp}.
This may accommodate singular behavior of $g$ aligned with the initial mesh $\grid_0$.

\begin{corollary}[approximation class of $u$]\label{C:approx-u}
Let the solution $u\in H^1_0(\Omega)$ of \eqref{weak-form} satisfy $u \in X^s_p(\Omega)$ with
$s \in (0, n+1], p\in (0,\infty]$ and $s - 1 - \frac{d}{p} + \frac{d}{2}>0$, where $X^s_p(\Omega)$ is
defined in \eqref{E:abstract-space}. Then
$u\in \A_{\frac{s-1}{d}}\big(H^1_0(\Omega);\grid_0\big)$ and
\begin{equation}\label{E:approx-u}
|u|_{\A_{\frac{s-1}{d}}} \lesssim |u|_{X^s_p(\Omega)}.
\end{equation}
Equivalently, $\sigma_N(u)$ defined in \eqref{E:sigmaN} satisfies
\begin{equation}\label{E:sigmaN-u}
\sigma_N(u) \lesssim |u|_{X^s_p(\Omega)} N^{-\frac{s-1}{d}} \quad N \ge \#\grid_0.
\end{equation}
\end{corollary}
\begin{proof}
In view of \eqref{eq:gradient-approx} of Proposition \ref{P:cont-vs-discont} (approximation of gradients), namely
\[
\inf_{v\in \mathbb{S}^{n,0}_\grid} \|\nabla (u-v)\|_{L^2(\Omega)} \lesssim
\inf_{v\in \mathbb{S}^{n,-1}_\grid} \|\nabla (u-v)\|_{L^2(\Omega)},
\]
we realize that it suffices to bound the element errors for $g=\nabla u \in L^2(\Omega;\R^d)$
by vector-valued discontinuous piecewise polynomials of degree $\le n-1$. Therefore, applying
Theorem \ref{T:direct} (regularity yields approximation) with $n$ replaced by $n-1$
gives the desired estimates \eqref{E:approx-u} and \eqref{E:sigmaN-u}.
\end{proof}

We now turn our attention to the coefficients $(\bA,c)$. 
Regarding $\bA$, Lemma \ref{L:perturbation} (continuous dependence on data) shows that the natural function space for $\bA$ is $L^\infty(\Omega,\R^{d\times d})$ provided $u\in H^1_0(\Omega)$. However, Lemma \ref{L:perturbation} also allows for $\bA \in L^r(\Omega,\R^{d\times d})$, 
$2 \le r < \infty$, provided $u \in W^1_p(\Omega)$ with $2\le p=\frac{2r}{r-2} < p_1$ which in turn is guaranteed by
Lemma \ref{L:Wp-regularity} ($W^1_p$-regularity). The latter permits discontinuities of $\bA$
within elements, which is of practical importance. Therefore, we consider the most general situation $2\le r \le \infty$ in the sequel. \looseness=-1

\begin{corollary}[approximation class of $A$]\label{C:approx-A}
For $0<\alpha_1\leq \alpha_2$ and $2 \le r \le \infty$ let the diffusion coefficient $\bA \in M(\alpha_1,\alpha_2)$ of
\eqref{weak-form} satisfy
$\bA \in X^s_p(\Omega;\R^{d\times d})$ with $s\in (0, n], p\in (0,\infty]$ and $s-\frac{d}{p} + \frac{d}{r} > 0$. Then
$\bA \in \ACA_{\frac{s}{d}} = \ACA_{\frac{s}{d}} \big((L^r(\Omega))^{d\times d};\grid_0\big)$ and
\begin{equation}\label{E:approx-A}
|\bA|_{\ACA_{\frac{s}{d}}} \lesssim |\bA|_{X^s_p(\Omega)}.
\end{equation}
Equivalently, $\wt \delta_\grid(\bA)_r$ defined in Section \ref{S:approx-data} satisfies
\begin{equation}\label{E:zeta-A-besov}
\inf _{\grid \in \mathbb{T}_N} \wt \delta_\grid(\bA)_r \lesssim |\bA|_{X^s_p(\Omega)} N^{-\frac{s}{d}}\quad
\forall \grid\in\grids_N, N \ge \#\grid_0,
\end{equation}
and this error decay is achieved by Algorithm \ref{algo:greedy} (greedy algorithm).
\end{corollary}
\begin{proof}
Simply recall the relation \eqref{e:wtdelta_to_delta} between the best constrained and unconstrained approximation errors and apply Theorem \ref{T:direct} (regularity yields approximation).
\end{proof}

Consider the special case $s=n$ and $r=\infty$ in Corollary \ref{C:approx-A}. We readily see that
$p > \frac{d}{n}$ which might be less than $1$ for $n>d$, hence the need for Besov spaces.

We finally deal with the reaction coefficient $c \in L^\infty(\Omega)$. According to Lemma
\ref{L:perturbation} (continuous dependence on data), and the discussion in Section \ref{S:DATA},
a natural space for $c$ is $L^q(\Omega)$ with $\frac{d}{2} < q \le \infty$, $s=0$; we could take $q=2$
for $d<4$. Section \ref{S:DATA} also reveals that the case $n=1$ is somewhat special 
in that we can exploit superconvergence in $W^{-1}_q(\Omega)$ with $q>d$. In fact, combining the argument following \eqref{E:c_super_conv} with \eqref{E:best-approx} yields
\[
\inf_{\wh{c} \in \mathbb{S}^{n-1,-1}_\grid} \|c-\wh{c}\|_{W^{-1}_q(\Omega)}^2 \lesssim 
\sum_{T\in\grid} h_T^{2t} \|c-\Pi_T{c}\|_{L^2(T)}^2 \lesssim \sum_{T\in\grid} h_T^{2t} \delta_\grid(c,T)_2^2 = \osc_\grid(c)_2^2
\]
with $0<t=1 - \frac{d}{2} + \frac{d}{q} < 2 - \frac{d}{2}$ provided $d<4$. This gives the following statement. We note that \eqref{E:c_super_conv-infty} could also be combined with \eqref{E:best-approx} for $n=1$ to obtain a similar result for $\osc_\grid(c)_\infty$ with $t=1$ and any $d\ge2$; however, we do not elaborate further.

\begin{corollary}[approximation of $c$]\label{C:approx-c}
Let $0\leq c_1 \leq c_2$ and the reaction coefficient $c \in R(c_1,c_2)$ satisfy $c\in X^s_p(\Omega)$ with
$s\in (0,n], p\in (0,\infty]$. If $n\ge 1$,  $q>\frac{d}{2}$, and $s -\frac{d}{p}+\frac{d}{q}>0$, then
$c \in \ACc_{\frac{s}{d}} = \ACc_{\frac{s}{d}} \big(L^q(\Omega);\grid_0\big)$ and
\begin{equation}\label{E:approx-c-n>1}
|c|_{\ACc_{\frac{s}{d}}} \lesssim |c|_{X^s_p(\Omega)}.
\end{equation}
If instead $n=1$, $q>d, s -\frac{d}{p}+\frac{d}{2}\ge 0, 0<t=1 - \frac{d}{2} + \frac{d}{q} <2-\frac{d}{2}$ and $d<4$, then
$c \in \ACc_{\frac{s+t}{d}} = \ACc_{\frac{s+t}{d}} \big(L^2(\Omega);\grid_0\big)$ and
\begin{equation}\label{E:approx-c}
|c|_{\ACc_{\frac{s+t}{d}}} \lesssim |c|_{X^s_p(\Omega)}.
\end{equation}
Equivalently for all $n\geq 1$,  $\wt \delta_\grid(c)_q$ defined in Section \ref{S:approx-data} satisfies
\[
\inf _{\grid \in \mathbb{T}_N} \wt \delta_\grid(c)_q \lesssim |c|_{X^s_p(\Omega)} N^{-\frac{s+t}{d}}\quad
\forall \grid\in\grids_N, N \ge \#\grid_0,
\]
with $t=0$ when $n>1$.
This error decay is achieved by Algorithm \ref{algo:greedy} (greedy algorithm).
\end{corollary}
\begin{proof}
In view of \eqref{e:wtdelta_to_delta}, inequality \eqref{E:approx-c-n>1} is a direct application of Theorem \ref{T:direct}
(regularity yields approximation). The superconvergence rate in \eqref{E:approx-c} is a
consequence of  \eqref{e:wtdelta_to_delta} and the proof of Proposition \ref{P:abstract_greedy} (abstract greedy) with $s$ replaced by $s+t$.
\end{proof}

We finally go back to the abstract space $X^s_p(\Omega)$, defined in \eqref{E:abstract-space},
and introduce the corresponding abstract approximation class 
$\ACX_{\frac{s}{d}}=\ACX_{\frac{s}{d}}\big(L^q(\Omega);\grid_0\big)$ of functions $v\in L^q(\Omega)$ such that
\begin{equation*}
|v|_{\ACX_{\frac{s}{d}}} = \sup_{N\ge \#\grid_0} \big(N^{\frac{s}{d}}
\inf_{\grid\in\grids_N} \osc_\grid(v)_q \big) < \infty
\quad\Rightarrow\quad
\inf_{\grid\in\grids_N} \osc_\grid(v)_q \le |v|_{\ACX_{\frac sd}} N^{-\frac{s}{d}}.
\end{equation*}
Consequently, Theorem \ref{T:direct} (regularity yields approximation) implies
\begin{equation}\label{E:aaprox-class-X}
X^s_p(\Omega) \subset \ACX_{\frac{s}{d}},
\quad
|v|_{\ACX_{\frac{s}{d}}} \lesssim |v|_{X^s_p(\Omega)}.
\end{equation}
We will utilize this abstract notation and estimates
in Section \ref{S:data-approx} while discussing the approximation 
of data $\data=(\bA,c,f)$ by a greedy algorithm.

%----------------------------------------------------------------------------------
\subsubsection{Global approximation: inverse estimates}\label{S:inverse-estimates}
%----------------------------------------------------------------------------------
%
Theorem \!\ref{T:direct} gives sufficient regularity
properties for a function $g\in L^q(\Omega)$ to belong
to an approximation class $\A_{\frac{s}{d}}\big(L^q(\Omega);\grid_0\big)$; this is called {\it direct estimate}.
Such regularity is written in terms of a Besov space $B^s_{p,p}(\Omega)$, except in the limiting case
$s=n+1$. The converse statement is also true and is called an {\it inverse estimate}: if $g$ belongs to an approximation class
$\A_{\frac{t}{d}}\big(L^q(\Omega);\grid_0\big)$, then it is a member of a Besov space
$\wh{B}^s_{p,p}(\Omega)$ provided $t>s$ and $0<s<n+1, s - \frac{d}{p} + \frac{d}{q} = 0$
\cite{BiDaDeVPe:02,GaspozMorin:14}.

Several comments are in order.
The Besov space $\wh{B}^s_{p,p}(\Omega)$ is defined via a multilevel decomposition of $L^p(\Omega)$
and coincides with $B^s_{p,p}(\Omega)$ only when $s<1+\frac{1}{p}$. This restriction of $s$ is natural
because $\V_\grid \subset \wh{B}^s_{p,p}(\Omega)$ for all $s$ but $\V_\grid \subset B^s_{p,p}(\Omega)$
requires $s<1+\frac{1}{p}$ according to \eqref{E:Vtau-besov}. The discrepancy between the
spaces $\wh{B}^s_{p,p}(\Omega)$ and $L^q(\Omega)$ is $\delta = s - \frac{d}{p} + \frac{d}{q} = 0$,
but the decay rate $t/d$ of $\A_{\frac{t}{d}}\big(L^q(\Omega);\grid_0\big)$ is larger than $s/d$.
This accounts for the embedding of $\A_{\frac{t}{d}}\big(L^q(\Omega);\grid_0\big)$
\[
 \sum_{n\in\N} \Big(\sigma_{2^n}(g) 2^{\frac{s}{d}n}\Big)^p
\le \sup_{n\in\N} \Big( \sigma_{2^n}(g) 2^{\frac{t}{d}n}\Big)^p \sum_{n\in\N} 2^{\frac{s-t}{d}pn}
\lesssim |g|_{\A_{\frac{t}{d}}}^p
\]
into a space with decay $s/d$ and summability $\ell^q$ that in turn embeds into $\wh{B}^s_{p,p}(\Omega)$
\cite{BiDaDeVPe:02,GaspozMorin:14}.
This reveals that there is no complete characterization of the approximation class $\A_{\frac{s}{d}}$ in
terms of Besov regularity.

%--------------------------------------------------------------------------------

\section{Data Approximation}\label{S:data-approx}
% \rhn{(AB $\longrightarrow$ AV)}
 
%\begin{itemize}
%\item
%  Tree approximation.
  
%\item
%  Greedy algoritms.
  
%\item
%  Discontinuous coefficients.

%\item
%  Approximation classes.
%\end{itemize}

This section focuses on the module $\DATA$ of Algorithms \ref{algo:AFEM-Kernel} ($\AFEMTS$) and
\ref{A:one-step-switch} (one-step $\AFEM$ with switch). According to Assumption \ref{A:optim-data}
(quasi-optimality of $\DATA$), the call
\begin{equation}\label{E:call-data}
[\wh\grid, \wh\data] = \DATA (\grid,\data,\tau)
\end{equation}
is meant to construct a quasi-optimal conforming refinement $\wh\grid$ of $\grid\in\grids$
and approximate piecewise polynomial
data $\wh\data = (\wh\bA, \wh{c}, \wh{f}) \in\D_{\wh\grid}$ over $\wh\grid$ that satisfies %
\begin{equation}\label{E:error-data}
\|\data - \wh\data \|_{D(\Omega)} \le \Cdata\tau
\end{equation}
as well as the constraints $\wh{\bA} \in M(\halpha_1,\halpha_2)$ and $\hc \in R(\hc_1,\hc_2)$ defined in
\eqref{E:structural-assumption}.

Sections \ref{ss:approx_diff} and \ref{ss:approx_reac} are devoted to the construction of $(\wh{\bA},\wh{c})$. To approximate the coefficients $(\bA,c)$ we proceed in two steps. First, we solve an {\it unconstrained} approximation problem upon computing the $L^2$-projection $(\wt{\bA},\wt{c})$ of $(\bA,c)$ onto the space of discontinuous piecewise polynomials of degree $\le n-1$; this step is linear, easily achieves the desired accuracy, but does not guarantee the monotonicity of oscillations with respect to refinement \eqref{E:monotoncicity-osc} and violates the constraints in \eqref{E:structural-assumption} unless $n=1$. Second, we resort to the nonlinear selection \eqref{E:wh-v} of the local $L^2$ approximation to force the resulting oscillations to be monotone. Third, we solve a {\it constrained} problem, which modifies $(\wt{\bA},\wt{c})$ locally into ($\wh{\bA},\wh{c})$ and restores \eqref{E:structural-assumption} without accuracy degradation; this is a delicate nonlinear procedure executed element by element, introduced and discussed in Section \ref{ss:positivity}.

The approximation of the right-hand side $f\in H^{-1}(\Omega)$ is a conceptually different linear process.
Without further structural assumptions on $f$ it is not possible to evaluate $\osc_{\grid}(f)_{-1}$ and reduce it. Hence we introduce surrogate estimators $\wt \osc_{\grid}(f)_{-1}$, which are larger than $\osc_{\grid}(f)_{-1}$,
but computable, for several classes of forcing functions $f$ relevant in practice. We discuss this in Section \ref{S:right-hand-side}.

We start in Section \ref{S:greedy} with a presentation and assessment of quasi-optimal \GREEDY algorithms to reduce the data error. 
An important consideration is that the
local error estimators $\{\osc_\grid(v,T)\}_{T\in\grid}$ may accumulate in $\ell^\infty$ as well as
in $\ell^q$ for $q<\infty$. Both are handled via a $\GREEDY$ algorithm similar to Algorithm~\ref{algo:greedy} but with different stopping criteria when the local errors accumulate in $\ell^q$ with $q<\infty$. 
The module $\DATA$ combines both:
its structure is displayed in Algorithm \ref{A:data} and its performance is elucidated in
Corollary \ref{C:data} below.

%-----------------------------------------------------------------------------------------
\subsection{Quasi-optimal {\rm \GREEDY} algorithms for data reduction}\label{S:greedy}
%-----------------------------------------------------------------------------------------
%
Algorithm \ref{algo:greedy} (greedy algorithm) is well suited for dealing with local error
estimators $\osc_\grid(v,T)_q$ that accumulate with respect to $T\in\grid$ in the space $\ell^\infty$. 
This is the framework for approximating coefficients $v = \bA,c$ in $L^\infty(\Omega)$, in which case
the local error estimators  $\osc_\grid(v,T)_\infty$ are defined in \eqref{E:zeta-Ac-loc}
for $r=q=\infty$. This requires that $v = \bA,c$ be piecewise uniformly continuous on $\grid$ for
$\osc_\grid(v;T)_\infty\to0$ as $h_T\to0$. However, for discontinuous $(\bA,c)$ and the forcing function $f$,
the accumulation of $\osc_\grid(v,T)_q$ for $v=\bA,c,f$ is in $\ell^q$ for $q < \infty$.
In this case, Algorithm \ref{algo:greedy} does not provide a direct
relation between a desired output tolerance $\tau$ for the total error
$$ 
\index{Error Estimators!$E_\grid(v)_q$: generic total error}
E_\grid (v)_q :=\| \{\osc_{\grid}(v,T)_q \}_{T \in \grid} \|_{\ell^q}= \left(\sum_{T \in \grid} \osc_{\grid}(v,T)_q^q\right)^{1/q}
$$
and the threshold $\delta$; recall that $\osc_\grid(v,T)_q := \|v-\wh{v}\|_{L^q(T)}$ for $T\in\grid$.

Another subtle difference from Algorithm~ \ref{algo:greedy} is that the algorithm $\GREEDY$ below does not start from $\grid_0$ but from any $\grid \in \grids$. Since \DATA and thus \GREEDY is called repeatedly within
\AFEM, it seems advantageous to exploit the mesh refinement already performed in the adaptive process rather than restarting from scratch; this thus improve the computational efficiency.  

\begin{algo}[$\GREEDY$]\label{A:greedy}
\index{Algorithms!\GREEDY: abstract greedy algorithm for \DATA}
Given a tolerance $\tau>0$, $0<q\leq \infty$, a number of bisection $b\geq 1$ performed per element to be refined, and an arbitrary conforming grid $\grid\in\grids$, not necessarily
$\grid_0$, \textsf{GREEDY} finds
a conforming refinement $\wh\grid\ge\grid$ of $\grid$ by bisection
and $\wh v\in \mathbb S^{n-1,-1}_\grid$ 
such that $E_\grid (v)_q \le\tau$:
%{\rm
\begin{algotab}
  \> $[\wh\grid,\wh v] = \textsf{GREEDY} \, (\grid,\tau,q,b,v)$\\
  \> \> $[\wh v] = \PROJECT(\grid,v)$\\
  \> \> while $E_\grid (v)_q > \tau$\\
  \> \> \> $[\marked] = \argmax\{\osc_\grid(v,T)_q \, : \, T\in\grid\}$\\
  \> \> \> $[\grid] = \REFINE(\grid,\marked,b)$\\
  \> \> \> $[\wh v] = \PROJECT(\grid,v)$\\
  \> \> return $\grid,\wh v$
\end{algotab}
%}
%
\end{algo}

In \GREEDY above, the element $T$ with largest error is refined as long as the total error $E_\grid (v)_q$ exceeds the target tolerance $\tau$. When the largest error is achieved by several elements, an ad-hoc criteria such as lexical order is used to break ties. 
We also recall that the routine \REFINE bisects all the elements in $\marked$ (in this case only one) $b$ times and performs additional refinements necessary to produce a conforming subdivision. \PROJECT computes the local approximations $\wh{v}$ of $v$ needed to evaluate $\osc_\grid(v,T)_q$; refer to  Section~\ref{S:DATA} and \eqref{E:wh-v} for the definition of $\wh{v}$.  The dependency on $\wh v$ in $\osc_\grid(v,T)_q$ and $E_\grid(v)$ is not indicated. 

To discuss the performances of the $\GREEDY$ algorithm, we recall that $X^s_p(\Omega;\grid_0)$ is the abstract space defined in \eqref{E:abstract-space} which satisfies
\begin{equation}\label{E:lp-sum-data}
\sum_{T\in \grid} | v |_{X^s_p(T)}^p \lesssim  | v |_{X_p^s(\Omega;\grid_0)}^p
\end{equation}
for all $\grid \in \grids$ and $v \in X^s_p(\Omega;\grid_0)$. 

The \GREEDY algorithm analyzed in Proposition~\ref{P:abstract_greedy} (abstract greedy) relies on the abstract assumptions \eqref{E:lq-summability}, \eqref{A:termination_abstract}, and \eqref{E:subadditive-norm}. With the aim of reducing the data oscillations, we  make these assumptions more concrete. 

\begin{assumption}[admissible set of parameters for \GREEDY]\label{ass:parameter-greedy}
\index{Assumptions!Admissible set of parameters for \GREEDY}
We say that the set of parameters $(v,s,t,p,q)$ is admissible for \GREEDY with local oscillations 
$\{\osc_\grid(v,T)_q\}_{T\in\grid}$ if $0<p,q \leq \infty$, $s,t \ge 0$ satisfy 
\begin{enumerate}
    \item $v\in X^s_p(\Omega;\grid_0)$;
    \item $t + s > 0, \;\; s - \frac{d}{p} + \frac{d}{q} \ge 0$ \;\; with strict inequality when $q=\infty$ or $s=n+1$, $p<1$;
    \item for $r := t + s - \frac{d}{p} + \frac{d}{q} > 0$
    \begin{equation}\label{E:decay-data}
\osc_\grid(v,T)_q \lesssim h_T^r |v|_{X^s_p(T)} \quad\forall \, T\in\grid, \grid\in \grids.
\end{equation}
\end{enumerate}
When the local oscillations considered are clear from the context, we say that $(v,s,t,p,q)$ is admissible for \GREEDY.
\end{assumption}

Relation ~\eqref{E:decay-data} replaces \eqref{A:termination_abstract}  and is a regularity assumption guaranteeing a convergence rate when approximating  $v$ by $\wh v = \PROJECT(\grid,v)$ (appearing in the definition of $\osc_\grid(v,T)_q$). We refer to Propositions~\ref{P:Bramble-Hilbert-Besov}, \ref{P:Bramble-Hilbert-Lipschitz}, and~\ref{P:Bramble-Hilbert-Sobolev} for examples where Assumption~\ref{ass:parameter-greedy} holds. Note that in view of \eqref{E:embedding} and \eqref{E:compact-lipschitz}, the conditions (ii) in Assumption~\ref{ass:parameter-greedy} guarantees $v\in X^s_p(\Omega) \subset L^q(\Omega)$. 
The parameter $t \ge 0$ reflects a possible additional power of $h$ in the oscillation term, see for e.g. \eqref{E:c_super_conv}, \eqref{E:c_super_conv-infty} and \eqref{E:weighted-L^2}.
Furthermore, in view of \eqref{E:lp-sum-data} assumption \eqref{E:sub-additive} is always satisfied by the $X^s_p(\Omega;\grid_0)$ semi-norms, and \eqref{E:lp-sum} is not needed any longer.

As alluded to above, the case $q<\infty$ is more complex to analyze and cannot solely rely on the decay property \eqref{E:decay-data} as in the proof of Proposition~\ref{P:abstract_greedy}. It requires the local oscillations to be monotone with respect to refinements.

\begin{assumption}[monotonicity of local oscillations]\label{a:monot_osc_data}
\index{Assumptions!Monotonicity of local oscillations}
We say that for $0<q\leq \infty$, the local errors satisfy the monotonocity property in $\ell^q$ if for any $v\in L^q(\Omega)$, any $\grid, \grid_* \in \grids$ with $\grid_* \geq \grid$ and any $T_* \in \grid_*$, $T\in \grid$ with $T_* \subset  T$, we have
\begin{equation}\label{e:osc_monotone_data}
\osc_{\grid_*}(v,T_*)_q \leq  \osc_{\grid}(v,T)_q.
\end{equation}
\end{assumption}

In view of Lemma \ref{L:monotonicity-osc} (monotonicity of oscillation),
Assumption~\ref{a:monot_osc_data} holds for the oscillations on $\bA$ and $c$ given in \eqref{E:zeta-Ac-loc} and \eqref{E:c_super_conv}, but not for the oscillations \eqref{E:zeta-star} of $f$. However, in Section~\ref{S:right-hand-side} below we derive computable surrogates for the local error $\osc_\grid(f,T)_{-1}$. These surrogates satisfy the monotonicity property and are used in turn to drive the \GREEDY algorithm. In passing, we note that we refrain from using the right-hand side of inequality \eqref{E:decay-data} as surrogate for the local oscillation. In fact, it is monotone with respect to refinements but at the expense of being difficult to evaluate because it involves the semi-norm
$|v|_{X^s_p(T)}$.

The following result is the counterpart of Proposition \ref{P:abstract_greedy} (abstract greedy) for \GREEDY with errors accumulating in $\ell^q$, $0<q \leq \infty$, and still starting from $\grid_0$.
We address the case where $\grid \not = \grid_0$ in Lemma \ref{l:greedy_start_grid} below.

\begin{proposition}[performance of $\GREEDY$]\label{P:greedy-data}
Let the initial subdivision $\gridk[0]$ of $\Omega\subset\R^d$
satisfy Assumption~\ref{A:initial-labeling} (initial labeling).
Let $\tau>0$ be the target tolerance and $b\geq 1$ be the number of bisections performed on each marked element.
Let $(v,s,t,p,q)$ satisfy Assumption~\ref{ass:parameter-greedy} (admissible set of parameters for {\rm\GREEDY}) with local errors $\{\osc_\grid(v,T)_q\}_{T\in\grid}$ which in turn verify 
Assumption~\ref{a:monot_osc_data} (monotonicity of local oscillations) in $\ell^{q}$. Then ${\rm \GREEDY}(\grid_0,\tau,q,b,v)$ terminates in a finite number of iterations
and
\begin{equation}\label{E:greedy_complexity}
E_\grid (v)_q \leq \tau \leq C |v|_{X^s_p(\Omega;\grid_0)} \big(\#\grid\big)^{-\frac{s+t}{d}}.
\end{equation}
with a constant $C = C(p,q,s,b,d,\Omega,\grid_0)$.
Furthermore, $v\in \mathbb X_{\frac{s+t}{d}}$
and $|v|_{\mathbb X_{\frac{s+t}{d}}} \lesssim |v|_{X^s_p(\Omega;\grid_0)}$.
Moreover, the estimate \eqref{E:greedy_complexity} is valid for tensor-valued functions $v$.
\end{proposition}
\begin{proof}
Since the proof is similar to that of Proposition \ref{P:abstract_greedy} (abstract greedy) with $e_\grid(v,T)_q=\osc_{\grid}(v,T)_q$, we only report the new ingredients. We recall that we use the convention $1/\infty=0$. Let $\grid_1, \dots, \grid_k$ be the sequence of refinements produced by $\GREEDY$, and  $T_1, \dots T_k$ be the sequence of marked elements. We need to estimate $\# \marked=k$ with $\marked = \{T_1,\dots, T_k\}$.
Set 
$$
\delta_i := \osc_{\grid_i}(v,T_i)_q,  \qquad (1 \leq i \leq k) \qquad \text{and }
\qquad \delta := \delta_{k-1} \,.
$$
Then, there holds
\begin{equation}\label{eq:pgreedy-1}
   E_{\grid_k}(v)_q  \leq \tau < E_{\grid_{k-1}}(v)_q  \leq \delta \, (\# \grid_{k-1})^\frac1q \\
\leq \delta \, (\# \grid_{k})^\frac1q .
\end{equation}
On the other hand, since $\REFINE$ does not increase the element estimators $\osc_{\grid_i}(v,T_i)$ thanks to \eqref{e:osc_monotone_data},
one has $\delta_i \geq \delta$ for any $i$, whence
$$
\osc_{\grid_i}(v,T_i)_q =  \delta_i \geq  \delta, \quad
\forall 1 \le i < k \,.
$$
Let us now partition $\marked$ into disjoint subsets ${\cal P}_j$ as in the proof of Proposition
\ref{P:abstract_greedy}. If $T_i \in {\cal P}_j$,
\eqref{E:decay-data} implies
$$
\delta \leq \osc_{\grid_i}(v,T_i)_q \lesssim h^r_T |v|_{X^s_p(T_i)} \leq 2^{-\frac{jr}2} |v|_{X^s_p(T_i)} \,,
$$
whence, exploiting the $\ell^p$ summability \eqref{E:lp-sum-data} gives
$$
\#\mathcal{P}_j\lesssim \delta^{-p}  2^{-\frac{jrp}{2}} |v|_{X^s_p(\Omega;\grid_0)}^p.
$$
which is similar to \eqref{nv-large}. Recalling \eqref{nv-small}, and proceeding 
as in the proof of Proposition \ref{P:abstract_greedy}, yields
$$
\delta \lesssim |v|_{X^s_p(\Omega;\grid_0)}\ ( \#\grid_k-\#\grid_0)^{-\frac{s+t}{d}-\frac{1}{q}} \,.
$$
We conclude the proof using \eqref{eq:pgreedy-1} and the bound $\#\grid_k \geq c_0 \#\grid_0$ for $c_0>1$.
\end{proof}

In contrast to Section \ref{S:tools}, and most of the existing literature, 
Algorithm~\ref{A:greedy} starts from a
refinement $\grid$ of $\grid_0$ rather than $\grid_0$ and thus exploits the mesh refinement already performed in the adaptive process. 
We now give a simple argument, updated from 
\cite{BonitoDeVoreNochetto:2013}, that shows that the number of elements $N(\grid,\tau,b,v)$ marked by \GREEDY starting from $\grid$ with target tolerance $\tau$ and refined $b \geq 1$ times is dominated by $N(\grid_0,\tau,1,v)$, namely
\begin{equation}\label{E:greedy-initialization}
    N(\grid,\tau,b,v) \le N(\grid_0,\tau,1,v).
\end{equation}
Estimate \eqref{E:greedy-initialization} is crucial because it avoids studying the cardinality of
$\GREEDY$ starting from $\grid \ne \grid_0$ directly, and simplifies the analysis. 
Even though \eqref{E:greedy-initialization} is plausible, the fact that the output of $\GREEDY(\grid_0,\tau, q, 1, v)$ is unrelated to $\grid$ makes it non-obvious.
In fact, note that we do not claim that $N(\grid,\tau,b,v) \le N(\grid_0,\tau,b,v)$, which is unclear. The proof presented below hinges on the fact that all the elements refined within $\GREEDY(\grid_0,\tau,q, 1,v)$ are either refined because they are marked by \GREEDY (and thus of largest oscillation) or because their refinement is necessary to guarantee conformity of the resulting subdivision. For our purposes, \eqref{E:greedy-initialization} suffices.  

\begin{lemma}[\GREEDY starting from $\grid$]\label{l:greedy_start_grid}
Let $\tau>0$ be a target tolerance and $b\geq 1$ be the number of bisections per marked element. Assume that the local errors employed by {\rm\GREEDY} satisfy Assumption~\ref{a:monot_osc_data} (monotonicity
of local oscillations) in $\ell^{q}$. Then, the  number of elements $N(\grid,\tau,b,v)$ marked by ${\rm \GREEDY}(\grid,\tau, q, b, v)$ satisfies \eqref{E:greedy-initialization} for any admissible refinement $\grid \in \grids$ of $\grid_0$.
\end{lemma}
\begin{proof}
    We simply write $\GREEDY (\grid_0,\tau,1)$ and $\GREEDY(\grid,\tau,b)$ because $v$ and $q$ are fixed. Let $N:=N(\grid_0,\tau,1,v)$, and 
    recall that the bisection rules define a unique forest $\grids$ emanating from $\grid_0$ and a unique sequence of elements $\{ T_i \}_{i=1}^{N}$ marked by $\GREEDY(\grid_0,\tau,1)$. We denote by 
    $\{\mathcal T^i\}_{i=1}^N$ the sequence of intermediate subdivisions built within $\GREEDY(\grid_0,\tau,1)$ starting with $\grid^0=\grid_0$: $T_i\in\grid^{i-1}$ is bisected once by
    $\REFINE$ which also produces the smallest conforming refinement $\grid^i$ of $\grid^{i-1}$ containing
    the two children of $T_i$. We thus say that $\GREEDY(\grid_0,\tau,1)$ satisfies the
    {\it minimality property} that all the elements refined are either marked elements because their error 
    is largest or necessary to guarantee conforming subdivisions. Notice that this is not true for $\GREEDY(\grid_0,\tau,b)$ when $b>1$.

    For any $\grid \in \grids$, we let $\Lambda_\grid$ be the set of indices $j\in \{ 1,...,N\}$ such that $T_j$ is never refined in the process to create $\grid$, i.e. $T_j$ is either an element of $\grid$ or a successor of an element of $\grid$. We show that 
    \begin{equation}\label{e:greedy_init_induc}
    N(\grid,\tau,b,v) \leq \#\Lambda_\grid
    \end{equation}
    by induction on $\# \Lambda_\grid$. If $\#\Lambda_\grid = 0$ then $\grid$ is a refinement of $\grid^{N}$, whence the monotonicity of the total error 
    $$
    E_{\grid}(v)_q \leq E_{\grid^{N}}(v)_q \leq \tau,
    $$ 
    guaranteed by \eqref{e:osc_monotone_data}, implies that $N(\grid,\tau,b,v)=0$; this satisfies \eqref{e:greedy_init_induc} as desired. 

    We now assume that \eqref{e:greedy_init_induc} is valid for any $\grid \in \grids$ such that $\# \Lambda_\grid \leq k$, a non-negative integer, and deduce it must also hold for any $\grid \in \grids$ such that $\# \Lambda_\grid \leq k+1$. 
    Let $\grid \in \grids$ be one such mesh, namely $\#\Lambda_\grid = k+1$. If $E_{\grid}(v)_q \leq \tau$ then $N(\grid,\tau,b,v)=0$ and 
    $N(\grid,\tau,b,v) \leq \#\Lambda_\grid$ holds trivially. 
    
    When instead $E_{\grid}(v)_q > \tau$, we let $j$ be the smallest index in $\Lambda_\grid$ and show that $T_j \in \grid$ using the minimality property of $\GREEDY(\grid_0,\tau,1)$.
    Assume by contradiction that $T_j \not \in \grid$ but $T_j$ belongs to a refinement $\widetilde \grid$ of
    $\grid$ and is thus a successor of an element $T \in \grid$. Note that $T$ is refined by $\GREEDY(\grid_0,\tau,1)$ to produce $T_j$ but was not marked because otherwise $T=T_i$ for some $i<j$ and $i\in \Lambda_\grid$, which would contradict the minimality of $j$. Hence, $T$ must have been refined by the \REFINE routine to guarantee conformity when bisecting a marked element $T_\ell$, $\ell<j$. Invoking the minimality of $j$ again yields that $\ell\notin\Lambda_\grid$
    and $T_\ell$ cannot be in $\grid$ because $T_\ell$ has been refined to get to $\grid$ by definition of
    $\Lambda_\grid$. Since \REFINE refines the minimal number of non-marked elements to guarantee conformity, and $\grid$ is conforming, $T$ must have been refined as well when refining 
    $T_\ell$ in the process of constructing $\grid$ and therefore cannot be in $\grid$. This is a contradiction and $T_j\in\grid$.
    
    Therefore, $\grid$ is a refinement of $\grid^{j-1}$ because all the elements marked or refined to ensure conformity by $\GREEDY(\grid_0,\tau,1)$ have been refined in the process of creating $\grid$. Moreover, $T_j\in\grid$ is the element with largest error $\osc_{\grid}(v,T_j)$ within $\grid$ (with ad-hoc criteria to break ties), because 
    $\osc_{\grid^{j-1}}(v,T_j)$ is largest in $\grid^{j-1}$ by definition of $T_j$
    and monotonicity of the local error \eqref{e:osc_monotone_data};
    hence $T_j$ must be the first element marked by $\GREEDY(\grid,\tau,b)$.
    Let $\grid^*$ be the subdivision obtained from $\grid$ upon bisecting $b$ times $T_j$. Notice that $\Lambda_{\grid^*}$ is a strict subset of $\Lambda_\grid$, because $j\notin\Lambda_{\grid^*}$, so that the induction assumption yields
    $$
    N(\grid,\tau,b,v) = 1 + N(\grid^*,\tau,b,v) \leq 1 + \# \Lambda_{\grid^*}\leq \#\Lambda_\grid.
    $$
    This proves \eqref{e:greedy_init_induc} and \eqref{E:greedy-initialization} follows immediately since 
    $
    \# \Lambda_\grid \leq N(\grid_0,\tau,1,v). 
    $
\end{proof}

Estimate \eqref{E:greedy-initialization} is critical to analyze the performances of \GREEDY starting from any admissible subdivision $\grid \in \grids$. We emphasize that the complexity estimate provided by Corollary~\ref{c:greed_start_grid} is expressed in terms of number of marked elements 
$N (\grid,\tau,q,b,v)$ and tolerance $\tau$ instead of error
and cardinality of $\grid$.
This is why $\textsf{GREEDY}$ can start from any mesh $\grid \in \grids$. 

\begin{corollary}[performance of greedy]\label{c:greed_start_grid}
Let the initial subdivision $\gridk[0]$ of $\Omega\subset\R^d$ satisfy Assumption~\ref{A:initial-labeling}
(initial labeling) and $\grid \in \grids$ be any admissible refinement of $\grid_0$.
Let $\tau>0$ be the target tolerance and $b\geq 1$ be the number of bisections performed on each marked element.
Let $(v,s,t,p,q)$ satisfy Assumption~\ref{ass:parameter-greedy} (admissible set of parameters for {\rm \GREEDY})
with local errors $\{\osc_\grid(v,T)_q\}_{T\in\grid}$ which in turn verify Assumption~\ref{a:monot_osc_data} (monotonicity of local oscillations) in $\ell^{q}$.
The number of marked elements $N(\grid,\tau,q,b,v)$ by ${\rm\GREEDY}(\grid,\tau,q,b,v)$ satisfies
\begin{equation}\label{E:greedy-estimate}
N(\grid,\tau,q,b,v) \le C |v|_{X^s_p(\Omega)}^{^{\frac{d}{s+t}}} \tau^{-\frac {d}{s+t}}.
\end{equation}
with a constant $C = C(p,q,s,b,d,\Omega,\grid_0)$.
Moreover, the estimate \eqref{E:greedy-estimate} is valid for tensor-valued functions $v$.
\end{corollary}

\begin{proof}
Invoking Proposition~\ref{P:greedy-data} (performance of $\GREEDY$), which gives rise to a mesh
$\wh\grid$, and Lemma~\ref{l:greedy_start_grid} ($\GREEDY$ starting from $\grid$), we readily deduce
\[
N(\grid,\tau,q,b,v) \le N(\grid_0,\tau,q,1,v) \le \#\wh{\grid} \le C 
|v|_{X^s_p(\Omega)}^{^{\frac{d}{s+t}}} \tau^{-\frac {d}{s+t}},
\]
which is the desired inequality \eqref{E:greedy-estimate}.
\end{proof}
%

%----------------------------------------------------------------------------------
\subsection{Constrained approximations}\label{ss:positivity}
%----------------------------------------------------------------------------------

 We discuss how the approximations produced by \GREEDY (see Corollary~\ref{c:greed_start_grid}) can be modified to satisfy the structural assumption \eqref{E:structural-assumption-wh} without sacrificing their accuracy.

\subsubsection{Constrained approximations of scalar functions}

The approximate data $\widetilde\data = (\widetilde{\vec{A}},\wtc,\wtf)$ constructed in the previous sections using the \GREEDY algorithm are not guaranteed to satisfy the necessary conditions for perturbed problem \eqref{E:perturbed-weak-form} with $\hdata = \widetilde \data$ to have a solution $\hu = \hu(\wh \data) \in H^1_0(\Omega)$. 
Recall that the data $\data=(\vec{A},c,f) \in D(\Omega)$ is assumed to satisfy the structural assumption \eqref{E:structural-assumption}, i.e., $\vec{A} \in M(\alpha_1,\alpha_2)$ and $c \in R(c_1,c_2)$ with $0<\alpha_1\leq \alpha_2$ and $0\leq c_1\leq c_2$.
It turns out that constructing approximate data $\hdata$ with the same constraints is a difficult task. We follow \cite{BonitoDeVoreNochetto:2013} and modify the data $\widetilde \data$ to obtain $\hdata=(\widehat{\vec{A}},\hc,\hf)$ in such a way that the approximation property of $\widetilde \data$ is preserved,
$$
\| \data - \hdata \|_{D(\Omega)} \le \Cdata \| \data - \widetilde \data\|_{D(\Omega)} \,,
$$
while ensuring that 
\begin{equation}\label{e:data_structural}
\widehat{\vec{A}} \in M\left(\frac{\alpha_1}2,C\halpha_2\right), \qquad \hc \in R\left(-\frac{\alpha_1}{4C_P^2},C \hc_2\right).
\end{equation}
Here $C$ is a constant independent of relevant quantities (we make this more precise below).  In particular, the data $\hdata$ satisfies the structural assumption \eqref{e:data_structural} which guarantees that the perturbed problem \eqref{E:perturbed-weak-form} has a unique solution.
Note that the general case is more subtle than when the data are approximated by piecewise constant approximations \eqref{E:meanvalues}, which are directly satisfying the structural assumption and used as motivation in Section~\ref{S:DATA}. 

We start by discussing a process modifying the approximation of a strictly positive scalar function $v \in L^\infty(\Omega)$, i.e. $v \in R(c_1,c_2)$ for some $0<c_1\leq c_2$; see \eqref{d:ReactionSpace}. Because the polynomial degree used to approximate the data might differ depending on the application, we use $m\in \mathbb N$ to denote a generic polynomial degree. We think of $\wt v \in \mathbb S^{m,-1}_\grid$ to be an approximation to $v$ not necessarily strictly positive. The following process modifies $\wt v$ locally to construct $\wh v \in \mathbb S^{m,-1}_\grid$. \index{Constants!$L$: threshold parameter for constrained approximation} It involves a parameter $L>2$ responsible for the truncation of $\wt v$ whenever it is too large, i.e., $\wt v \geq L c_2$. For $T\in \grid$, we set $\wh v|_T:=\wh v_T$ where
\begin{equation}\label{def:positivity}
\wh v_T:= \left\lbrace 
\begin{array}{ll}
c_2 & \ \textrm{when } \| \wt v \|_{L^\infty(T)}\geq L c_2, \\
\wt v|_T - \min_{x \in T} \wt v(x) + \frac {c_1} 2   & \ \textrm{when otherwise } \min_{x \in T} \wt v(x) < \frac{c_1} 2, \\ 
\wt v|_T & \  \textrm{otherwise.}
\end{array}
\right.
\end{equation}

Corollary~\ref{c:positivity} below is in essence Proposition~3 of \cite{BonitoDeVoreNochetto:2013} and states that the constructed $\wh v$ satisfies
$$
0<\frac{c_1}2 \leq \wh v \leq \left(\frac 1 2 + 2 L\right)c_2, \qquad \textrm{a.e. in }\Omega.
$$
This is at the expense of inflating the approximation error in $L^q$, $1\leq  q \leq \infty$, by a multiplicative constant $C$ only depending on $d$, $m$, $c_2/c_1$, $q$, $L$ and the shape regularity of $\grids$
$$
\| v - \wh v \|_{L^q(T)} \leq C \| v - \wt v \|_{L^q(T)}.
$$

In preparation for this result, we introduce the following notations.
We denote by $C_I$ the smallest constant such that for any $T\in \grid$ and any polynomial $P\in \P_m(T)$, the inverse inequality
\begin{equation}
\label{e:inverse}
\|\nabla P\|_{L^\infty(T)}\le C_I\|P\|_{L^\infty(T)}|T|^{-1/d}, 
\end{equation}
holds. The inverse inequality constant $C_I$ only depends on the shape regularity of $\grids$, $m$ and $d$. Note that for such polynomial $P \in \P_m(T)$, we have
$$
| P(x) - P(y)| \leq C_I \|P\|_{L^\infty(T)} |T|^{-1/d} |x-y|, \qquad \forall x,y\in T.
$$
Consequently, for any $\rho>0$ and $x \in T$, we define 
$$
T(x,\rho) := T \cap \overline{B(x,\rho |T|^{1/d}/C_I)},
$$
which is motivated by the fact that for $x \in T$ and $y \in T(x,\rho)$ we have
\begin{equation}\label{e:poly_variation}
| P(x) - P(y)| \leq C_I \|P\|_{L^\infty(T)} |T|^{-1/d} |x-y| \leq \rho \|P\|_{L^\infty(T)}.
\end{equation}
Critical for the analysis below is the existence of a constant $0<C_S(\rho) \leq 1$ depending on $\rho$ but also on $d$, $m$, and the shape regularity of $\grids$, such that 
\begin{equation}\label{e:Cs}
|T(x,\rho)| \geq C_S (\rho)|T| \qquad \forall x \in T, \quad T\in \grid.
\end{equation}
This constant $C_S(\rho)$ assesses the area of a subset of $T$ where the polynomial $P$ varies no more that $\rho \| P \|_{L^\infty(T)}$ away from $P(x)$. 

We are now in position to analyze the effect of the nonlinear correction \eqref{def:positivity}. We proceed locally over each $T\in \grid$ and start with the case where $\| \wt v \|_{L^\infty(T)}$ is large (Lemma~\ref{l:upper_bound_pres}). We then discuss the case where $\wt v(x)$ is small on $T$ (Lemma~\ref{l:lower_bound_pres}) while for the remaining case, the function $\wt v$ does not need to be modified on $T$. These three cases are collected in Corollary~\ref{c:positivity} for scalar valued functions and in Corollary~\ref{c:positivity-matrix} for matrix valued functions. In all the arguments below we used the convention $a^{1/\infty}=1$ for any $a>0$.

\begin{lemma}[locally enforcing constraints for large approximations]\label{l:upper_bound_pres}
    Let $\grid \in \grids$ be any conforming refinement of $\grid_0$ satisfying Assumption~\ref{A:initial-labeling} (initial labelling). Let $c_2>0$, $T\in \grid$ and $v_T\in L^\infty(T)$ satisfying $0<v_T \leq c_2$ a.e. in $T$.
    Furthermore, for $m\geq 0$ and $L>2$, assume that $\wt v_T \in \P_m(T)$ satisfies
    \begin{equation}\label{e:extra_upper_bound}
    \| \wt v_T \|_{L^\infty(T)}\geq  Lc_2.
    \end{equation}
    Then for the constant function $\wh v_T := c_2 \in \P_m(T)$ there holds
    $$
    \frac{c_1}2 <  \wh v_T < Lc_2 <\left(\frac 1 2 + 2 L\right) c_2.
    $$
    Moreover, for $1\leq q\leq \infty$, we have
    $$
    \| v_T - \wh v_T \|_{L^q(T)} \leq C^+_2\| v_T - \wt v_T \|_{L^q(T)},
    $$
    where $C^+_2 := \frac{4 C_S^{-1/q}}{L-2}$ and $C_S = C_S(1/2)$ is the constant appearing in \eqref{e:Cs} with $\rho=1/2$.
    \end{lemma}
\begin{proof}
Let $x_0 \in T$ and $\wt c_{2,T}$ defined by the relation
$$
 \wt c_{2,T}:=  |\wt v_T(x_0)|:= \| \wt v_T \|_{L^\infty(T)}.
$$
In view of the Lipschitz property \eqref{e:poly_variation} applied to $P=\wt v_T$ and with $\rho=\frac12$, we have
$$
|\wt v_T(x)-\wt v_T(x_0)|\le \frac {\wt c_{2,T}} 2,
$$
for $x \in T_0:=T(x_0,\frac 1 2) \subset T$. 
Recall \eqref{e:Cs}, which implies that $|T_0| \geq  \wt C_S |T|$ for some constant $\wt C_S := C_S(1/2)$ only depending $d$, $n$,  and the shape regularity of $\grids$.
On the one hand, this implies that $\wt v_T|_{T_0}$ is bounded below with $|\wt v_T(x)| \ge \frac {\wt c_{2,T}} 2$ for $x\in T_0$ and, on the other hand, $v_T$ is bounded from above by
$$
0 \leq v_T(x) \leq c_2 \leq  L^{-1} \wt c_{2,T}, \qquad x\in T.
$$
Consequently, for $x \in T_0$ and since $L > 2$, there holds 
$$
0 \leq v_T(x) \leq L^{-1} \wt c_{2,T}  \leq \frac {\wt c_{2,T}} 2  \leq |\wt v_T(x)|
$$ 
and thus
$$
|v_T(x)-\wt v_T(x)| \ge |\wt v_T(x)| - v_T(x) \geq \left(\frac 12 - \frac{1}{L}\right) \wt c_{2,T} = \frac{L-2}{2L}\wt c_{2,T},
$$
which indicates that $v_T$ and $\wt v_T$ are sufficiently far apart on a substantial portion $T_0$ of $T$.
Thus is responsible for the $L^q$-bound below. In fact, we have
\begin{equation}\label{e:diff_T0}
\|v_T-\wt v_T \|_{L^q(T)} \geq \| v_T-\wt v_T  \|_{L^q(T_0)} \geq  \frac{L-2}{2L} \wt c_{2,T} |T_0|^{1/q},
\end{equation}
whence, from the definition $\wh v_T := c_2$ and using \eqref{e:Cs}, we deduce
$$
\| v_T - \wh v_T \|_{L^q(T)} \leq 2 c_2 |T|^{1/q}    \leq 2 L^{-1} \wt{c}_{2,T} |T|^{1/q}
\leq \frac{4 C_S^{-1/q}}{L-2}  \|v_T-\wt v_T \|_{L^q(T)}
$$
as desired.
\end{proof}

\begin{lemma}[locally enforcing constraints for small approximations]\label{l:lower_bound_pres}
    Let $\grid \in \grids$ be any conforming refinement of $\grid_0$ satisfying Assumption~\ref{A:initial-labeling} (initial labelling). Let $0<c_1\leq c_2$,  $T\in \grid$ and $v_T\in L^\infty(T)$ satisfying $c_1<v_T \leq c_2$ a.e. in $T$.  Furthermore, for $m \geq 0$ and $L>2$ assume that $\wt v_T \in \P_m(T)$ satisfies
    \begin{equation}\label{e:assumption_other}
\| \wt v_T \|_{L^\infty(T)} \leq Lc_2
\end{equation}
and
    \begin{equation}\label{e:extra_lower_bound}
 \min_{x\in T} \wt v_T(x) < \frac{c_1}2.
 \end{equation}

    Then the function $\wh v_T := \frac {c_1} 2 + \wt v_T - \min_{x \in T} \wt v_T(x)   \in \P_n(T)$ is such that
    $$
    \frac {c_1} 2 \leq \wh v_T \leq 2 Lc_2 + \frac {c_1} 2 \leq \left( 2L+ \frac 1 2\right) c_2
    $$
    and
    $$
    \| v_T - \wh v_T \|_{L^q(T)} \leq C^+_1 \| v_T - \wt v_T \|_{L^q(T)},
    $$
    where $C^+_1:=(1+ C_S^{-1/q}(\rho)))$ and $C_S(\rho)$ is the constant appearing in \eqref{e:Cs} with $\rho =c_1/(2Lc_2)$.
    \end{lemma}
\begin{proof}
We define $x_0 \in T$, $\wt c_{1,T} \in \mathbb R$ by the relations
$$
\wt c_{1,T}:=\wt v_T(x_0):=\min_{x\in T}\wt v_T(x).
$$
From the Lipschitz property \eqref{e:poly_variation} and the assumption \eqref{e:assumption_other}, we find that
$$
|\wt v_T(x)-\wt v_T(x_0)|\le \frac{c_1}{2}
$$
for $x \in T_0:=T(x_0,\rho)$ with $\rho:= \frac{c_1}{2Lc_2}$.
Recall \eqref{e:Cs}, which implies that $|T_0| \geq  \widetilde{C}_S |T|$ for some constant $\widetilde{C}_S := C_S(\rho)$ only depending $d$, $m$, $c_2/c_1$, $L$,  and the shape regularity of $\grids$.

 For $x \in T_0$, we proceed by estimating the difference
 $$
 v_T(x) - \wt v_T(x) = v_T(x) - (\wt v_T(x) - \wt v_T(x_0)) - \wt v_T(x_0) \geq c_1 - \frac {c_1} 2 - \wt c_{1,T} = \frac {c_1} 2 - \wt c_{1,T}>0
 $$
 because $\wt c_{1,T} < c_1/2$ by assumption \eqref{e:extra_lower_bound}.
 This implies that
 $$
 |T_0|^{1/q} \left(\frac {c_1} 2 - \wt c_{1,T} \right) \leq \| v_T - \wt v_T\|_{L^q(T)},
 $$
 and $v_T$ and $\wt{v}_T$ are uniformly far apart in the substantial part $T_0$ of $T$.
 Therefore, $\wh v_T := \wt v_T + (\frac {c_1} 2 - \wt c_{1,T})$ satisfies
 $$
   \frac {c_1}2 \leq \wh v_T \leq 2Lc_2 + \frac {c_1} 2
  $$
  because $\wt{c}_{1,T} \ge -\|\wt{v}\|_{L^\infty(T)} \ge -Lc_2$ by assumption \eqref{e:assumption_other}, and
 \begin{equation*}
     \begin{split}
 \| v_T - \wh v_T \|_{L^q(T)} & \leq  \| v_T - \wt v_T \|_{L^q(T)} +  \left(\frac {c_1} 2 - \wt c_{1,T} \right) |T|^{1/q} 
 \\ 
 &\leq \| v_T - \wt v_T \|_{L^q(T)} +  \left(\frac {c_1} 2 - \wt c_{1,T} \right) \widetilde C_S^{-1/q} |T_0|^{1/q} \\
 & \leq 
 \big(1+ \widetilde C_S^{-1/q}\big)\| v_T - \wt v_T\|_{L^q(T)}.
\end{split}
\end{equation*}
This proves the assertions.
\end{proof}

\begin{corollary}[locally enforcing constraints]\label{c:positivity}
Let $\grid \in \grids$ be any conforming refinement of $\grid_0$ satisfying Assumption~\ref{A:initial-labeling} (initial labelling).
Let $0<c_1\leq c_2$, $T\in \grid$ and $v_T \in L^\infty(T)$ satisfying $c_1 \leq v_T \leq c_2$ a.e. in $T$. Then, for $m \geq 0$, $L>2$, and $\wt v_T \in \mathbb P_m(T)$, the function $\wh v_T  \in \mathbb P_m(T)$ defined in \eqref{def:positivity} satisfies
    $$
    \frac {c_1} 2 \leq \wh v_T \leq \left(\frac 12 +2L\right)c_2  \qquad \textrm{a.e. in }T.
    $$
    Moreover for $1\leq q \leq \infty$, we have
    $$
    \| v_T - \wh v_T \|_{L^q(T)} \leq \max(C^+_1,C^+_2) \| v_T - \wt v_T \|_{L^q(T)} \quad \forall T \in \grid,
    $$
    where $C^+_1$ and $C^+_2$ are the constants appearing in Lemmas~\ref{l:lower_bound_pres} and \ref{l:upper_bound_pres}, which only depend on $d$, $m$, $c_2/c_1$, $L$,  and the shape regularity of $\grids$.
\end{corollary}
\begin{proof}
The desired results follows from Lemma~\ref{l:upper_bound_pres} when 
$$
\| \wt v_T \|_{L^\infty(T)} \geq Lc_2
$$
and from Lemma~\ref{l:lower_bound_pres} when
$$
\| \wt v_T \|_{L^\infty(T)} < Lc_2 \qquad \textrm{and} \qquad \min_{x \in T} \wt v_T < \frac{c_1}{2}.
$$
In the remaining case 
$$
\| \wt v_T \|_{L^\infty(T)} < Lc_2 \qquad \textrm{and} \qquad \min_{x \in T} \wt v_T(x) \geq \frac{c_1}{2},
$$
since $\wh v_T=\wt v_T$ satisfies the desired constraints there is nothing to prove.
\end{proof}

%-----------------------------------------------------------------------------------
\subsubsection{Constrained approximation of the diffusion coefficients}\label{ss:approx_diff}
%-----------------------------------------------------------------------------------

For matrix-valued functions, the constraints are on the eigenvalues of the matrix rather than on the coefficients themselves. Although this requires a few adjustments, the process is similar to the scalar case. We recall that for $0<\alpha_1\leq \alpha_2$, $M(\alpha_1,\alpha_2) \subset L^\infty\big(\Omega;\mathbb R^{d\times d}_{\textrm{sym}}\big)$ denotes the class of symmetric matrix-valued functions whose eigenvalues lie between $\alpha_1$ and $\alpha_2$; see \eqref{d:DiffusionSpace}. 

Algorithm \APPLYCONSTRA is based on \eqref{def:positivity} and modifies approximations $\wt \bA \in (\mathbb S_{\grid}^{n_A,-1})^{d\times d}$ of $\bA\in M(\alpha_1,\alpha_2)$ to produce uniformly positive definite approximations $\wh \bA \in (\mathbb S_\grid^{n_A,-1})^{d\times d}$ of $\bA$. 

\begin{algo}[\APPLYCONSTRA]\label{A:const-project-A}
\index{Algorithms!\APPLYCONSTRA: modify an approximation of $\bA$ to satisfy the structural constrains}
Given a threshold parameter $L>2$, $0<\alpha_1 \leq \alpha_2$, a conforming refinement $\grid\in \grids$ of $\grid_0$, and $\wt \bA \in (\mathbb S_\grid^{n_A,-1})^{d\times d}$, this routine constructs a positive definite $\wh \bA \in (\mathbb S_\grid^{n_A,-1})^{d\times d}$.
\begin{algotab}
  \> $[\wh \bA] = \APPLYCONSTRA(\grid,\alpha_1,\alpha_2,L,\wt \bA)$\\
  \> \> For $T \in \grid$\\
  \> \> \> $\widetilde \alpha_{1,T}=\inf \{y^t \wt \bA(x) y, \ x\in T, \ |y|=1\} $\\
  \> \> \> $\wt \alpha_{2,T} = \sup \{ | y^t \wt \bA(x) y |, \ x \in T, \ | y|=1 \}$\\
  \> \> \> if $\wt \alpha_{2,T} \geq L \alpha_2$ \\
  \> \> \> \> $\wh \bA|_T = \alpha_2 \bI_d$ \\
  \> \> \> else if $\wt \alpha_{1,T} < \frac {\alpha_1} 2$ \\
  \> \> \> \> $\wh \bA|_T = \wt \bA|_T - \left(\frac {\alpha_1} 2 - \wtalpha_{1,T} \right)\vec{I}_d$ \\
  \> \> \> else \\
  \> \> \> \> $\wh \bA|_T = \wt \bA|_T$\\
  \> \> return $\wh \bA$
\end{algotab}
\end{algo}

Notice that \APPLYCONSTRA preserves symmetry, i.e, if $\wt \bA$ is symmetric so is the output $\wh \bA$.
In addition, when $n=1$ and $\wt \bA \in (\mathbb S_\grid^{0,-1})^{d\times d}$ is the piecewise constant local average of $\bA$, the output $\wh \bA$ of $\APPLYCONSTRA$ is $\wh \bA = \wt \bA$ since
in that case the parameters $\wtalpha_{1,T}$ and $\wtalpha_{2,T}$ satisfy
$$
\wtalpha_{1,T} \geq \alpha_1 > \frac {\alpha_1} 2 \qquad \textrm{and} \qquad \wtalpha_{2,T} \leq \alpha_2 < L\alpha_2, \qquad \forall T \in \grid.
$$
This is consistent with the observation made in Section~\ref{S:DATA}.

The next corollary hinges on Corollary~\ref{c:positivity} (locally enforcing constraints) to derive properties of \APPLYCONSTRA. In passing, we recall that for $\bA \in L^p(\Omega;\mathbb R^{d\times d})$ we write
$$
\| \bA \|_{L^p(\Omega)}:= \| | \bA | \|_{L^p(\Omega)},
$$
where for $x\in \Omega$, $|\bA(x)|$ is the spectral norm of $\bA(x)$.

\begin{corollary}[locally enforcing constraints for matrices]\label{c:positivity-matrix}
Let the threshold be $L>2$, $0<\alpha_1\leq \alpha_2$,  and $\bA \in M(\alpha_1,\alpha_2)$.
    Let $\grid \in \grids$ be any conforming refinement of $\grid_0$ and $\wt \bA \in (\mathbb S_\grid^{n_A,-1})^{d\times d}$ be a symmetric approximation of $\bA$. Then the output $[\wh \bA]=\APPLYCONSTRA(\grid,\alpha_1,\alpha_2,L,\wt \bA)$ is symmetric and satisfies 
    $$
    \frac {\alpha_1} 2 \leq \lambda_j(\wh \bA) \leq \left(\frac 12+ 2L\right) \alpha_2  \qquad \textrm{a.e. in }\Omega, \ 1\leq j \leq d.
    $$
    Moreover for $1\leq q \leq \infty$, we have 
    $$
    \| \bA - \wh \bA \|_{L^q(T)} \leq \Cdata \| \bA - \wt \bA \|_{L^q(T)} \qquad \forall T \in \grid,
    $$
    where $\Cdata:=\max(C^+_1,C^+_2)$ and $C^+_1$ and $C^+_2$ are the constants appearing in Lemmas~\ref{l:lower_bound_pres} and \ref{l:upper_bound_pres}, which only depends on $d$, $n_A$, $\alpha_2/\alpha_1$, $L$,  and the shape regularity of $\grids$. \looseness=-1
\end{corollary}
\begin{proof}
We observe that $\wt{\bA}$ is not assumed positive semi-definite.
We argue locally and fix $T\in \grid$. Let $\wt \alpha_{2,T} >0 $ and $y_0 \in \mathbb R^d$ be such that $|y_0|=1$ and
$$
\wt \alpha_{2,T}:= \sup_{x\in T} | y_0^t \wt \bA(x) y_0| := \sup_{x\in T} \sup_{y\in \mathbb R^d, |y|=1} |y^t \wtA(x) y|.
$$

We first consider the case $\wt \alpha_{2,T} \geq L \alpha_2$ for which $\wh \bA|_T := \alpha_2 \vec{I}_d$. For $x\in T$, we set
$$
a(x) := y_0^t \bA(x) y_0 \qquad \textrm{and} \qquad \wta_T(x) = y_0^t \wt \bA(x) y_0 \in \mathbb P_{n_A}(T).
$$
These notations allows us to reduce to the scalar case upon noting that
$$
\| a - \wta_T \|_{L^q(T)} \leq \| \bA - \wt \bA \|_{L^q(T)}
$$
and  $\alpha_1 \leq a \leq \alpha_2$ a.e. in $\Omega$.
Because $\wtalpha_{2,T} \geq L \alpha_2$, Lemma~\ref{l:upper_bound_pres} with $m=n_A$ guarantees that $\ha_T:= \alpha_2$ satisfies 
$$
\| a-\ha_T \|_{L^q(T)} \leq C^+_2 \|  a- \wta_T \|_{L^q(T)} \leq C^+_2 \| \bA - \wt \bA \|_{L^q(T)}.
$$
Consequently, the matrix-valued approximation $\wh \bA|_T := \alpha_2 \vec{I}_d$ satisfies
$$
\| \bA - \wh \bA \|_{L^q(T)} = \| a-\ha_T \|_{L^q(T)} \leq C^+_2 \| \bA - \wt \bA \|_{L^q(T)}.
$$
This proves the desired result when $\wtalpha_{2,T} \geq L \alpha_2$.

We now consider the case where $\wtalpha_{2,T} < L \alpha_2$ and define $\wtalpha_{1,T}  \in \mathbb R$, $y_1 \in \mathbb R^d$ with $|y_1|=1$ by the relations
$$
\wtalpha_{1,T} = \inf_{x\in T} y_1^t \wt \bA(x) y_1 = \inf_{x\in T} \inf_{|y|=1} y^t \wt \bA(x) y.
$$
We also redefine the associated scalar functions for $x\in T$ using $y_1$ instead of $y_0$
$$
a(x) := y_1^t \bA(x) y_1 \qquad \textrm{and} \qquad \wta_T(x) = y_1^t \wt \bA(x) y_1 \in \mathbb P_{n_A}(T).
$$
If $\wtalpha_{1,T} < \frac {\alpha_1} 2$ then $\wh \bA|_T = \wt \bA|_T + (\frac {\alpha_1} 2 - \wtalpha_{1,T})\vec{I}_d$. Lemma~\ref{l:lower_bound_pres} with $m=n_A$ ensures that $\ha_T = \wta_T + \frac {\alpha_1} 2 - \wtalpha_{1,T}$ satisfies
$$
\frac {\alpha_1} 2 \leq \ha_T  \leq \left(\frac 1 2 + 2M\right) \alpha_2
    $$
    and
    $$
    \| a - \ha_T \|_{L^q(T)} \leq C^+_1 \| a - \wta_T \|_{L^q(T)} \leq C^+_1 \| \bA - \wt \bA \|_{L^q(T)}.
    $$ 
Thus, $\wh \bA|_T$ satisfies the desired properties provided $\wt{\alpha}_{2,T} \ge L\alpha_2$ as well.

It the remaining case $\wtalpha_{2,T} < L \alpha_2$ and $\wtalpha_{1,T} \geq \frac {\alpha_1} 2$, the function $\wh \bA|_T=\wt \bA|_T$ satisfies the desired properties and there is nothing to prove.
\end{proof}

As a corollary, we report the complexity of an algorithm that concatenates the linear approximation of \GREEDY with the nonlinear correction into the constraint of \APPLYCONSTRA. We recall from Corollary \ref{C:approx-A} (approximation class of $\bA$) that the admissible set of parameters of $\bA$ for
$\GREEDY$ are $n_A\le n-1$
\[
s_A \in (0,n_A], \quad p_A \in (0,\infty], \quad 
q_A\in[2,\infty], \quad s_A - \frac{d}{p_A} + \frac{d}{q_A} > 0, \quad t_A=0.
\]

\begin{corollary}[complexity of constrained \GREEDY for $A$]\label{c:complex_consteained_A}
Let the initial mesh $\gridk[0]$ of 
 $\Omega\subset\R^d$
satisfy Assumption~\ref{A:initial-labeling} (initial labeling) and $\grid \in \grids$ be any admissible refinement of $\grid_0$.
Let $\tau>0$ be the target tolerance, $b\geq 1$ be the number of bisections performed on each marked element, and $L>2$ be a threshold parameter.
Furthermore, assume that $(\bA,s_A,t_A,p_A,q_A)$ satisfies Assumption~\ref{ass:parameter-greedy} (admissible set of parameters for {\rm \GREEDY})
with local oscillations $\{\|\bA-\wh\bA\|_{L^{q_A}(T)}\}_{T\in\grid}$ and, in addition, 
$\bA\in M(\alpha_1,\alpha_2)$ for some $0<\alpha_1\leq \alpha_2$.
The algorithm
\begin{algotab}
  \> $[\wh\grid,\wt \bA] = \GREEDY(\grid,\tau,q_A,b,\bA)$\\
  \> $[\wh \bA] = \APPLYCONSTRA(\wh\grid,\alpha_1,\alpha_2,L,\wt \bA)$
\end{algotab}
where ${\rm \GREEDY}$ is applied to the $d(d+1)/2$ distinct components of $\bA$,
marks $N$ elements of $\grid$ for refinement with
\begin{equation}\label{E:greedy-estimate-A}
N \leq C |\bA|_{X^{s_A}_{p_A}(\Omega;\grid_0)}^{^{\frac{d}{s_A}}} \tau^{-\frac {d}{s_A}}
\end{equation}
and $C = C(p_A,q_A,s_A,b,d,n_A,\alpha_2/\alpha_1,L,\Omega,\grid_0)$.
Moreover, $\wh \bA \in (\mathbb S_\grid^{n_A,-1})^{d\times d}$ satisfies
\begin{equation}\label{e:halphas}
\index{Constants!$(\halpha_1,\halpha_2)$: lower and upper bounds on the spectrum of $\hbA$}
\wh \bA \in M \left(\halpha_1,\halpha_2 \right):
\qquad\halpha_1 = \frac {\alpha_1} 2,\quad \halpha_2 = (1 +4L)\frac{\alpha_2}{2},
\end{equation}
and there is a constant $\Cdata>0$ such that
$$
\| \bA - \wh \bA \|_{L^q(\Omega)} \leq \Cdata \tau.
$$
\end{corollary}
\begin{proof}
    This result follows upon invoking Corollary~\ref{c:greed_start_grid} (performance of greedy) and Corollary~\ref{c:positivity-matrix} (locally enforcing constraints for matrices).
\end{proof}

\begin{remark}[constrained approximation class of matrices]\label{r:equivalent_class_A}
    As a consequence of Corollary~\ref{c:complex_consteained_A}, we realize that for $A\in M(\alpha_1,\alpha_2)$,
    $$
\delta_\grid (\bA)_r \leq \wt \delta_\grid (\bA)_r \leq \Cdata \delta_\grid (\bA)_r, 
    $$
    where the best approximation error $\delta_\grid (\bA)_r$ and best constrained approximation error $\wt \delta_\grid (\bA)_r$ are defined in \eqref{e:delta_A_c} and \eqref{e:wtdelta_A_c}.
\end{remark}

%----------------------------------------------------------------------------------
\subsubsection{Constrained approximation of the reaction coefficients}\label{ss:approx_reac}
%----------------------------------------------------------------------------------

If the reaction coefficient $c \in R(c_1,c_2)$ is strictly positive ($c_1>0$), then Corollary~\ref{c:positivity} (locally enforcing constraints) with $m=n_c$ directly applies to $v_T=c|_T$, $T\in \grid$, and guarantees that the approximate coefficient $\wh c \in \mathbb S_{\grid}^{n_c,-1}$ defined on $T\in \grid$ by $\wh c|_T:=\wh v_T$ satisfies 
$$
\wh c \in R\left( \wh{c}_1, \wh{c}_2 \right):\qquad
\wh{c}_1 = \frac{c_1}{2}, \quad \wh{c}_2 = (1+4L)\frac {c_2}2.
$$ 
However, reaction coefficients are not necessarily strictly positive on $\overline{\Omega}$ and Corollary~\ref{c:positivity} cannot be invoked directly.
Instead, we take advantage of the fact that the perturbed problem \eqref{E:perturbed-weak-form} 
is still well-posed provided $\wh{c} \ge -\halpha_1/(2 C_P^2)$ and the approximate diffusion coefficient $\wh \bA\in M(\halpha_1, \halpha_2)$ of $\bA\in M(\alpha_1,\alpha_2)$ satisfies
$\halpha_1 \ge \alpha_1/2$ according with \eqref{E:wh-alpha}; hence $\wh{c} \ge -\alpha_1/(4 C_P^2)$.
Therefore, we apply Corollary \ref{c:positivity} to the shifted reaction coefficient  $v=c+\frac{\halpha_1}{C_P^2}$, which satisfies
\begin{equation}\label{E:v1-v2}
v_1 := c_1 + \frac{\halpha_1}{C_P^2} \le v \le c_2 + \frac{\halpha_1}{C_P^2} =: v_2.
\end{equation}
Below is the proposed algorithm for the construction of $\wh c$ in the general case
$c_1\ge 0$.

\begin{algo}[\APPLYCONSTRC]\label{A:const-project-c}
\index{Algorithms!\APPLYCONSTRC: modify an approximation of $\bA$ to satisfy the structural constrains}
Given $L>2$, $\halpha_1 > 0$, a conforming refinement $\grid\in \grids$ of $\grid_0$, and $\wt c \in \mathbb S_\grid^{n_c,-1}$, this routine constructs $\wh c \in \mathbb S_\grid^{n_c,-1}$ as follows:
\begin{algotab}
  \> $[\wh c] = \APPLYCONSTRC(\grid,\halpha_1,L,\wt c)$\\
  \> \> $\wt v = \wt c + \frac{\halpha_1}{C_P^2}$\\
  \> \> For $T \in \grid$\\
  \> \> \> if $\| \wt v \|_{L^\infty(T)} \geq L v_2$\\
  \> \> \> \> $\wh v|_T = v_2$ \\
  \> \> \> else if $\min_{x\in T} \wt v(x) < \frac{v_1}{2}$ \\
  \> \> \> \> $\wh v|_T = \wt v|_T - \min_{x\in T} \wt v(x) + \frac {v_1}2$\\
  \> \> \> else \\
  \> \> \> \> $\wh v|_T = \wt v|_T$\\
  \> \> $\wh c = \wh v - \frac{\halpha_1}{C_P^2}$\\
  \> \> return $\wh c$
\end{algotab}
\end{algo}

We note that if $n_c=0$, then $\wt c$ is the piecewise average of $c$ and \APPLYCONSTRC does not modify $\wt c$ which already satisfies the structural assumption~\eqref{e:data_structural}.

The next result shows that the output $\hc$ of \APPLYCONSTRC is a modification of $\wt c$ which satisfies $\hc \in R(\wh c_1,\wh c_2)$, with 
\begin{equation}\label{e:def_hc}
\wh c_1:= \frac {c_1}{2}-\frac{\halpha_1}{2C_P^2} \qquad \textrm{and} \qquad \wh c_2:=(1 +4L)\frac{c_2}{2}+(4L-1)\frac{\halpha_1}{2C_P^2}
\end{equation}
without affecting the approximation of $c$ in $L^q$, $1\leq q \leq \infty$ (up to a multiplicative constant). In particular $\wh c_1 \geq -\frac{\halpha_1}{2C_P^2}$, which is necessary for the well-posedness of the perturbed problem~\eqref{E:perturbed-weak-form} when $\wh A \in M(\halpha_1,\halpha_2)$.

\begin{corollary}[locally enforcing constraints for nonnegative scalar functions]
\label{c:positivity-scalar}
Let $\bA\in M(\halpha_1,\halpha_2)$ with $0<\halpha_1\leq \halpha_2$, and $c \in R(c_1,c_2)$ with $0\leq c_1 \leq c_2$. Let $L>2$ and $v_1 \le v_2$ be defined in \eqref{E:v1-v2}.
Let $\grid \in \grids$ be any conforming refinement of $\grid_0$ and $\wt c \in \mathbb S_\grid^{n_c,-1}$. Then the output  $[\wh c]=\APPLYCONSTRC(\grid,\halpha_1,L,\wt c)$ satisfies 
    $$
    \wh c_1 \leq \wh c \leq \wh c_2 \qquad 
\textrm{a.e. in }\Omega,
    $$
    where $\wh c_1$ and $\wh c_2$ are given by \eqref{e:def_hc}. 
    Moreover for $0<q \leq \infty$, we have
    $$
    \| c - \wh c \|_{L^q(T)} \leq \Cdata \| c - \wtc \|_{L^q(T)} \qquad \forall T \in \grid,
    $$
    where $\Cdata$ is a constant only depending on $d$, $n$, $v_2/v_1$, $\Omega$, $L$,  and the shape regularity of $\grids$. 
    \end{corollary}
\begin{proof}
Set $\kappa:=\frac{\halpha_1}{C_P^2}$ and $v:=c+\kappa \in  R(c_1+\kappa,c_2+\kappa)$  so that $c_1+\kappa>0$.
On each $T\in \grid$, we invoke Corollary~\ref{c:positivity} (locally enforcing constraints) with $m=n_c$, $\wt v_T = \wt c|_T + \kappa$ and where $c_1$, $c_2$ are replaced by $c_1+\kappa$, $c_2+\kappa$ respectively. 
Hence, we deduce that the function $\wh v$  constructed within \APPLYCONSTRC satisfies 
$$
\frac{c_1+\kappa}{2} \leq \wh v  \leq (1 +4L)\frac{c_2+\kappa}{2}.
$$
and
\begin{equation}\label{e:apprx_c_data_tmp}
\| v - \wh v \|_{L^q(T)} \leq \Cdata \| v - \wt v\|_{L^q(T)} \qquad \forall T\in \grid,
\end{equation}
with a constant $\Cdata$ depending on $d,n, v_2/v_1, L$, and the shape regularity of $\grids$.
Shifting back, $c=v-\kappa$ and $\hc:= \wh v -\kappa$, we find that the approximation $\wh c$ constructed by \APPLYCONSTRC satisfies
$$
\frac{c_1+\kappa}{2}- \kappa \leq \wh c \leq (1+4L)\frac{c_2+\kappa}{2}-\kappa
$$
or equivalently
$$
\frac{c_1}{2} - \frac{k}{2} \leq \hc \leq  (1+4L)\frac{c_2}{2} +  (4L-1)\frac{\kappa}{2}.
$$
In view of \eqref{e:def_hc} and $\kappa=\frac{\halpha_1}{C_P^2}$, this is the first desired inequality in disguised. 

Furthermore, the second desired inequality follows from \eqref{e:apprx_c_data_tmp} because for $T\in \grid$ we have $c - \hc = v - \wh v$ and $c - \wtc = v - \wt v$.
\end{proof}

The next corollary combines the linear approximation of \GREEDY together with the nonlinear correction into the constraint of $\APPLYCONSTRC$. We recall from Corollary \ref{C:approx-c} (approximation class of $c$) that the admissible set of parameters of $c$ for $\GREEDY$ are $n_c\le n-1, s_c \in (0,n_c], p_c \in (0,\infty]$
\begin{align*}
n_c > 0 &\quad\Rightarrow\quad
q_c > \frac{d}{2}, \quad s_c - \frac{d}{p_c} + \frac{d}{q_c} >0, \quad t_c = 0;
\\
n_c = 0 &\quad\Rightarrow\quad
q_c =2, \quad s_c - \frac{d}{p_c} + \frac{d}{2} >0, \quad 0< t_c < 2 - \frac{d}{2}.
\end{align*}

\begin{corollary}[complexity of constrained \GREEDY for $c$]\label{c:complex_greedy_const_c}
Let the initial subdivision $\gridk[0]$ of $\Omega\subset\R^d$
satisfy Assumption~\ref{A:initial-labeling} (initial labeling) and $\grid \in \grids$ be any admissible refinement of $\grid_0$.
Let $\tau>0$ be the target tolerance, $b\geq 1$ be the number of bisections performed on each marked element, $L>2$ be the threshold parameter, and $\halpha_1>0$.
Furthermore, assume that $(c,s_c,t_c,p_c,q_c)$ satisfies Assumption~\ref{ass:parameter-greedy} (admissible set 
of parameters for {\rm \GREEDY}) with local oscillations
$\{\|c-\wh{c}\|_{L^{q_c}(T)} \}_{T\in\grid}$ and that $c\in R(c_1,c_2)$ for some $0\leq c_1 \leq c_2$.
The algorithm
\begin{algotab}
  \> $[\wh\grid,\wtc] = \GREEDY(\grid,\tau,q_c,b,c)$\\
  \> $[\wh c] = \APPLYCONSTRC(\wh\grid,\halpha_1,L,\wtc)$
\end{algotab}
marks $N$ elements of $\grid$ for refinement with
\begin{equation}\label{E:greedy-estimate-c}
N \leq C |c|_{X^{s_c}_{p_c}(\Omega;\grid_0)}^{^{\frac{d}{s_c+t_c}}} \tau^{-\frac {d}{s_c+t_c}}
\end{equation}
and a constant $C = C(p_c,q_c,s_c,b,d,n_c,v_2/v_1,L,\Omega,\grid_0)$ with $v_1\le v_2$ defined in \eqref{E:v1-v2} to construct $\wh\grid$.
The function $\wh c \in \mathbb S_{\wh\grid}^{n_c,-1}$ is a piecewise polynomial of degree $\le n_c$ over $\wh{\grid}$ and satisfies
$$
\wh c \in R(\wh c_1,\wh c_2),
$$
where $\wh c_1 \le \wh c_2$ are given by \eqref{e:def_hc}. Moreover, for $1 < q_c \le \infty$, there is a constant $\Cdata$ only depending on $d$, $n_c$, $v_2/v_1$, $\Omega$, $L$, and the shape regularity of $\grids$ such that
$$
\| c - \wh c \|_{L^{q_c}(\Omega)} \leq \Cdata \tau.
$$
\end{corollary}

\begin{proof}
 Simply apply Corollary~\ref{c:greed_start_grid} 
 (performance of greedy) and Corollary~\ref{c:positivity-scalar} (locally enforcing constraints for nonnegative scalar functions).
\end{proof}

\begin{remark}[constrained approximation class of scalars]\label{r:equivalent_class_c}
           Corollary~\ref{c:complex_greedy_const_c} implies that for $c\in R(c_1,c_2)$,
    $$
     \delta_\grid (c)_q \leq \wt\delta_\grid (c)_q \leq \Cdata \delta_\grid (c)_q,
    $$
        where the best approximation error $\delta_\grid (c)_q$ and best constrained approximation error $\wt \delta_\grid (c)_q$ are defined in \eqref{e:delta_A_c} and \eqref{e:wtdelta_A_c}.
\end{remark}

%-------------------------------------------------------------------------------
\subsection{Approximation of the load term $f$}\label{S:right-hand-side}
%-------------------------------------------------------------------------------

We now turn our attention to the question of designing a practical algorithm for reducing the global oscillation
\begin{equation}\label{e:datasec-osc}
E_\grid(f)_{-1}^2:= \sum_{T \in \grid}\| f - P_\grid f\|_{H^{-1}(\omega_T)}^2
\approx \sum_{z \in \vertices}\| f - P_\grid f\|_{H^{-1}(\omega_z)}^2,
\end{equation}
where the projection $P_\grid$ is defined in \eqref{E:defPmesh}.
The approximation of functionals in $H^{-1}(\Omega)$ is rather intricate and out of reach without assuming additional structure enabling practical evaluation of their actions on polynomial functions. 

We examine three cases of independent interest.
In Section~\ref{ss:data_f_Lp} we consider $f \in L^q(\Omega)$ for $q$ satisfying $2d/(d+2)< q \leq \infty$, which includes the most common setting $f \in L^2(\Omega)$. Sections~\ref{ss:data_f_dirac} and~\ref{ss:data_f_div} present examples of right-hand sides not in $L^1$. In Section~\ref{ss:data_f_dirac} we treat the case $f = g \delta_\Gamma$, where $\Gamma$ is an hyper-surface not necessarily captured by the faces of the  subdivisions and $g \in L^q(\Gamma)$, $q\geq 2$, while in Section~\ref{ss:data_f_div} we consider $f=\divo {\bg}$
for some ${\bg} \in L^2(\Omega;\R^d)$.
In all cases, the total error $E_{\grid}(f)_{-1}$ is estimated by a surrogate $\wt E_{\grid}(f)_{-1}$, namely $E_{\grid}(f)_{-1} \le \Cdata \wt{E}_{\grid}(f)_{-1}$ \looseness=-1
\index{Error Estimators!$\wt E_{\grid}(f)_{-1}^2$: generic surrogate estimator for the approximation of the load term}
$$
\wt E_{\grid}(f)_{-1}^2 := \sum_{T\in \grid} \widetilde \osc_\grid(f,T)_q^2
$$
with a definition of $\wt \osc_\grid(f,T)_q$ depending on the situation but local to $T\in \grid$ (and not on stars).
This allows Algorithm~\ref{A:greedy} ($\GREEDY$) to reduce $\wt{E}_\grid(f)_{-1}$.

Before starting, we recall relevant definitions and results from Section~\ref{S:aposteriori}
(a posteriori error analysis). For $z \in \vertices$, we denote by $\grid_z \subset \grid$ all the elements in $\omega_z$ and $\faces_z \subset \faces$ all the faces in $\omega_z$. For
$\ell \in H^{-1}(\Omega)$, the restriction $P_\grid \ell|_{\omega_z}$ belongs to the space $\mathbb F(\grid_z)=\mathbb F_{m_1,m_2}(\grid_z)$ made of functional whose action against $w \in H^1_0(\omega_z)$ reads \looseness=-1
\begin{equation}\label{e:data_local_F}
\langle \ell, w \rangle = \sum_{T \in \mesh_z} \int_T q_T  \, w + \sum_{F \in \faces_z} \int_F q_F \, w
\end{equation}
for some $q_F \in P_{m_1}(F)$, $F \in \faces_z$ and $q_T \in P_{m_2}(T)$, $T\in \grid_z$. The polynomial degrees are chosen to be $m_1 = n-1$ and $m_2=n-2$ but can be general in this discussion.

Corollary~\ref{C:local-near-best-approx-of-P} (local near-best approximation) guarantees $P_\grid \ell|_{\omega_z}$ is the quasi-best discrete functional in $\mathbb F(\grid_z)$, namely
\begin{equation}\label{e:PT_near_best_data}
    \| \ell - P_\grid \ell \|_{H^{-1}(\omega_z)}
    \leq
    C_P \inf_{\chi \in \F(\grid_z)} \| \ell - \chi \|_{H^{-1}(\omega_z)}.
\end{equation}
This will be used repeatedly to replace $P_\grid\ell$ by more tractable quantities and justify the use of \GREEDY algorithms to reduce \eqref{e:datasec-osc}.

%--------------------------------------------------------------------------------------
\subsubsection{The case $f \in L^q(\Omega)$}\label{ss:data_f_Lp}
%--------------------------------------------------------------------------------------

In this section, we show how to reduce the oscillation error \eqref{e:datasec-osc} when $f \in L^q(\Omega)$, with $q>\frac{2d}{d+2}$ to guarantee that $L^q(\Omega)$ compactly embeds in $H^{-1}(\Omega)$.  Note that this not only includes the most treated case in the literature $f \in L^2(\Omega)$ but also the more intricate cases $q<2$ originally analyzed in \cite{CohenDeVoreNochetto:2012}. 

If $\Pi_\grid f$ is the $L^2$-projection of $f$ into the space $\mathbb S_{\grid}^{n_f,-1}$ of
discontinuous piecewise polynomials of degree $n_f$, let $\wh f \in \mathbb S_{\grid}^{n_f,-1}$ be defined by \eqref{E:wh-v}; $n_f = n-1$ in some applications but not always. Since $\wh f|_{\omega_z} \in \mathbb F(\grid_z)$ by taking $q_F=0$ and $q_T = \wh f|_T$ in \eqref{e:data_local_F}, the local near-best approximation property \eqref{e:PT_near_best_data} of $P_\grid$ implies
$$
\| f - P_\grid f \|_{H^{-1}(\omega_z)} \le C_P  \| f-\wh f \|_{H^{-1}(\omega_z)}.
$$
Furthermore, for $v \in H^1_0(\omega_z)$ one has
$$
\langle  f - \wh f,v\rangle \leq \| f - \wh f\|_{L^q(\omega_z)}  \| v \|_{L^{\wt q}(\omega_z)}
$$
where $\frac 1 q + \frac 1 {\wt q} =1$.
Note that the restriction $q> \frac{2d}{d+2}$ guarantees that $1 \leq \wt q<\frac{2d}{d-2}$ and thus $\sob(H^1)>\sob(L^{\wt q})$. Therefore, Lemma \ref{L:Poincare} (first Poincar\'e inequality) yields
$$
\|   f - \wh f \|_{H^{-1}(\omega_z)} \lesssim  \textrm{diam}(\omega_z)^{1+d(\frac 1 2-\frac  1 q)} \| f - \wh f\|_{L^q(\omega_z)}.
$$
Returning to \eqref{e:datasec-osc}, after rearranging the terms element-wise and invoking the shape-regularity of $\grids$, we obtain $E_\grid(f)_{-1} \le \Cdata \wt E_\grid(f)_{-1}$, where
\begin{equation}\label{e:f-Lp-surrogate}
\wt E_\grid(f)_{-1}^2:=\sum_{T \in \grid} \wt \osc_\grid(f,T)_q^2,
\end{equation}
and $\widetilde \osc_{\grid}(f,T)_q: = h_T^t \| f - \wh  f \|_{L^{q}(T)}$ with $t:=1+d(\frac 1 2-\frac 1 q)>0$.

In view of the definition \eqref{E:wh-v}, the local oscillations $\wt \osc(f,T)_q$ satisfy 
Assumption~\ref{a:monot_osc_data} (monotonicity of local oscillations) in $\ell^2$ and we can now employ Algorithm~\ref{A:greedy} (\GREEDY) with local errors $\widetilde \osc_{\grid}(f,T)_q$ accumulating in $\ell^2$.
Recall that we use the convention $X^0_q(\Omega;\grid_0)=L^q(\Omega)$. 

\begin{corollary}[approximation class of $f\in L^q(\Omega)$]\label{C:p-greedy-Lp}
Let the initial subdivision $\gridk[0]$ of $\Omega\subset\R^d$
satisfy Assumption~\ref{A:initial-labeling} (initial labeling) and $\grid \in \grids$ be any admissible refinement of $\grid_0$.
Let $\tau>0$ be the target tolerance and $b\geq 1$ be the number of bisections performed on each marked element.
Let $2d/(d+2)< q \leq \infty$ and set $t=1+d(\frac 1 2-\frac 1 q)$.
Let $(f,s,t,p,2)$ satisfy Assumption~\ref{ass:parameter-greedy} (admissible set of
parameters for {\rm \GREEDY}) with local oscillations $\{\wt \osc(f,T)_q\}_{T\in\grid}$.
Then  $[\wh \grid, \wh f]=\GREEDY(\grid,\tau,2,b,f)$ terminates in a finite number of steps
with $\wt E_{\wh \grid}(f)_{-1} \leq \tau$, whence
$$
E_{\wh \grid}(f)_{-1} \leq \Cdata \tau.
$$
Moreover, the number $N$ of marked elements by ${\rm \GREEDY}$ satisfies
\begin{equation}\label{E:greedy-estimate-f}
N \lesssim |f|_{X^s_p(\Omega;\grid_0)}^{^{\frac{d}{s+t}}} \tau^{-\frac {d}{s+t}}.
\end{equation}
%with a constant $C = C(p,q,s,b,d,\Omega,\grid_0)$.
In particular, $f \in \mathbb F_{\frac{d}{s+t}}$ with 
$| f |_{\mathbb F_{\frac{d}{s+t}}} \lesssim \|f\|_{X^s_p(\Omega;\grid_0)}$.
\end{corollary}
\begin{proof}
Directly apply Corollary~\ref{c:greed_start_grid} (performance of \GREEDY).
\end{proof}

%--------------------------------------------------------------------------------------
\subsubsection{The case $f = g \delta_{\mathcal C}$}\label{ss:data_f_dirac}
%--------------------------------------------------------------------------------------

We now consider the case where the right-hand side data $f$ is a density supported on a Lipschitz hyper-surface $\mathcal C \subset \Omega$ in $\mathbb R^d$ with $(d-1)$-measure $|\mathcal C|<\infty$.

The intricate interactions between bulk and interface contributions on $P_\grid$ makes it difficult to analyze when $f = g \delta_{\mathcal C}$ with density $g \in L^q(\mathcal{C)}$. 
We take a simpler approach, likely suboptimal when $n>1$ and $d>2$, which discards $P_\grid$ in view of the near-best approximation property \eqref{e:PT_near_best_data} 
\begin{equation}\label{e:dirac_overestim}
\| f - P_\grid f\|_{H^{-1}(\omega_z)}\lesssim \| f \|_{H^{-1}(\omega_z)}.
\end{equation}
The right-hand side of the above estimate is the starting point of the analysis in \cite{CohenDeVoreNochetto:2012} assuming $n=1$ and $d=2$.

We start with the derivation of a first upper bound for the local error $\| f \|_{H^{-1}(\omega_z)}$. \begin{lemma}[local oscillation]\label{l:surrogate_dirac}
    Let $\grid \in \grids$, $z \in \nodes$, and $q>\frac{2(d-1)}{d}$. If $g\in L^q(\mathcal{C})$
    and $t:=\frac d 2- \frac 1 q (d-1)>0$, then there holds
    \begin{equation}\label{e:local_hm1_dirac}
    \| f \|_{H^{-1}(\omega_z)} \lesssim | \omega_z \cap \mathcal C|^{\frac{t}{d-1}} \| g \|_{L^q(\omega_z \cap \mathcal C)} \lesssim \sum_{T \subset \omega_z} h_T^{t} \| g \|_{L^q(T \cap \mathcal C)}.
    \end{equation}
\end{lemma}
\begin{proof}
    For $v \in H^1_0(\omega_z)$ and $\frac 1 q + \frac 1 {\wt q} =1$, we have
\begin{equation}\label{e:diract_cont}
    \langle f, v \rangle  = \int_{\omega_z \cap \mathcal C} g v  \leq \| g \|_{L^q(\omega_z \cap \mathcal C)} \| v \|_{L^{\wt q}(\omega_z \cap \mathcal C)}.
    \end{equation}
    We realize that $H^{1/2}(\omega_z \cap \mathcal C)$ compactly embeds in $L^{\wt q}(\omega_z \cap \mathcal C)$ because
    \begin{align*}
    t:= \sob(H^{1/2}(\omega_z \cap \mathcal C)) - \sob(L^{\wt q}(\omega_z \cap \mathcal C)) &= \frac 1 2 - (d-1) \left(\frac 1 2 - \frac 1 {\wt q}\right) \\
    & =\frac d 2 - \frac 1 q (d-1)> 0,
    \end{align*}
    provided $q > \frac{2(d-1)}{d}$. Consequently, we find that
    $$
      \| v \|_{L^{\wt q}(\omega_z \cap \mathcal C)}  \lesssim | \omega_z \cap \mathcal C|^{\frac{t}{d-1}} \| v \|_{H^{1/2}(\omega_z \cap \mathcal C)}.
    $$
    It remains to invoke the continuity \eqref{E:trace-1/2} of the trace operator to write
    $$
      \| v \|_{L^{\wt q}(\omega_z \cap \mathcal C)}  \lesssim | \omega_z \cap \mathcal C|^{\frac{t}{d-1}} \| v \|_{H^{1}(\omega_z)},
    $$
    which, together with \eqref{e:diract_cont}, yields
    the first estimate in \eqref{e:local_hm1_dirac}.
    To deduce the second estimate, it suffices to note that $| \omega_z \cap \mathcal C| \lesssim \textrm{diam}(\omega_z)^{d-1} \lesssim h_T^{d-1}$ for $T\subset \omega_z$ and that $\| g \|_{L^q(\omega_z \cap \mathcal C)} \leq \sum_{T\subset \omega_z} \| g \|_{L^q(T \cap \mathcal C)}$.
\end{proof}

Estimate \eqref{e:dirac_overestim} and Lemma~\ref{l:surrogate_dirac} provide a surrogate for data oscillation
\begin{equation}\label{e:dirac_surrogate}
\wt{E}_\grid(f)_{-1}^2 := \sum_{T\in\grid} \wt{\osc}_\grid (g,T)_q^2,
\quad
\wt \osc_\grid(g,T)_q := h_T^t \| g \|_{L^q(T\cap \mathcal C)},
\end{equation}
where $t=\frac{d}{2} - \frac{1}{q} (d-1)$.
The quantity $\wt \osc_\grid(g,T)_q$ verifies Assumption \ref{a:monot_osc_data} (monotonicity of local oscillations) with $\Omega$ replaced by $\mathcal C$ because of its element-wise structure. Therefore, Proposition
\ref{P:greedy-data} (performance of \GREEDY) states that Algorithm~\ref{A:greedy}
(\GREEDY) can reduce $\wt{E}_\grid(f)_{-1}$. This is in contrast with the star-wise \GREEDY algorithm analyzed in \cite{CohenDeVoreNochetto:2012}, which requires that all marked stars are refined $d$ times to ensure all the faces in the marked stars are refined.

We now discuss the performance of \GREEDY with local indicators $\wt \osc_\grid(g,T)_q$.

\begin{lemma}[approximation class of $f=g\delta_{\mathcal{C}}$]\label{l:p-greedy-dirac}
Let $\mathcal C \subset \Omega$ be Lipschitz hyper-surface.
Let the initial subdivision $\gridk[0]$ of $\Omega\subset\R^d$
satisfy Assumption~\ref{A:initial-labeling} (initial labeling) and $\grid \in \grids$ be any admissible refinement of $\grid_0$.
Let $\tau>0$ be the target tolerance and $b\geq 1$ be the number of bisections performed on each marked element, and $2(d-1)/d< q \leq \infty$.
Then  $[\wh \grid, \wh f]=\GREEDY(\grid,\tau,2,b,f)$ terminates in a finite number of steps
with surrogate estimator $\wt E_{\wh \grid}(f)_{-1} \leq \tau$ defined in \eqref{e:dirac_surrogate}, 
whence
$$
E_{\wh \grid}(f)_{-1} \leq \Cdata\tau.
$$
Moreover, the number $N$ of marked elements by {\rm \GREEDY} satisfies
\begin{equation}\label{E:greedy_complexity-dirac}
N \lesssim  \| g \|_{L^q(\mathcal C)}^{2(d-1)} \tau^{-2(d-1)}.
\end{equation}
%with a constant $C = C(q,s,b,d,\Omega,\grid_0)$.
In particular, $f=g \delta_{\mathcal C} \in \mathbb F_{\frac{1}{2(d-1)}}$ with 
$| f |_{\mathbb F_{\frac{1}{2(d-1)}}} \lesssim \|g\|_{L^q(\mathcal C)}$
\end{lemma}

\begin{proof}
This proof mainly follows the proof of Proposition~\ref{P:greedy-data} (performance of $\GREEDY$) but requires a few modifications to account for the geometry of the problem.
Since in turn the proof of Proposition~\ref{P:greedy-data} describes modifications to the proof of Proposition~\ref{P:abstract_greedy} (abstract greedy), we now provide a complete proof. We proceed in several steps. 
We first consider the call $\GREEDY(\grid_0,\tau,2,1,f)$ from $\grid_0$ with one bisection $b=1$
and accumulation in $\ell^2$, and discuss the general call from $\grid$ with $b\ge1$
in the last step of this proof.\\

\step{1} {\it Termination}.
Since $h_T$ decreases monotonically to $0$ with bisection, so does $\wt{\osc}_\grid(g,T)_q$. 
Consequently, \GREEDY terminates in finite number $k\ge 1$ of
iterations. Let $T_1, \dots T_k$ be the sequence of marked elements, with ${\cal M}= \{ T_1 , ..., T_k\}$ and $\grid_1,...,\grid_k$ be the sequence of refinements produced by \GREEDY starting from $\grid_0$. 
Upon termination, the surrogate error satisfies $\wt E_{\grid_k}(f)_{-1}\le \tau$,
whence
$
E_{\grid_k}(f)_{-1} \le \Cdata \tau.
$
\looseness=-1

\medskip
\step{2} {\it Counting}.
To estimate the cardinality of $\grid_k$, we need to count $\#\cal M$. Set   
$$
\delta_i := \wt{\osc}_{\grid_i}(g,T_i)_q, \quad  1 \leq i \leq k, \qquad \text{and }
\qquad \delta := \delta_{k-1} \,.
$$
Then, there holds
\begin{equation}\label{dirac:pgreedy-1}
\wt E_{\grid_k}(f)_{-1} \leq \tau < \wt E_{\grid_{k-1}}(f)_{-1} \leq \delta \, (\# \grid_{k-1})^\frac1 2
\leq \delta \, (\# \grid_{k})^\frac12 .
\end{equation}

We organize the elements in $\marked$ 
by size in such a way that allows for a counting argument. Let ${\mathcal P}_j$ be the set of elements $T$ of $\marked$ with size
\[
 2^{-(j+1)}
 \le
 |T|<2^{-j}
\quad\Rightarrow\quad
 2^{-\frac{j+1}{d}}\le h_T<2^{-\frac{j}{d}}.
\]
We first observe that all $T$'s in ${\mathcal P}_j$ are
\emph{disjoint}. This is because if $T_1,\,T_2\in {\mathcal P}_j$ and $\mathring{T}_1\cap\mathring{T}_2\neq\emptyset$, then one of them is contained in the other, say $T_1\subset T_2$, due to the bisection procedure which works in any dimension $d\ge1$; see Section \ref{S:bisection}. Hence,
\[
	|T_1|\le\frac{1}{2}\,|T_2|
\]
contradicting the definition of ${\mathcal P}_j$. 
On the one hand, this implies the first bound
\begin{equation}\label{dirac:nv-small}
 2^{-\frac{(j+1)(d-1)}{d}}\,\#{\mathcal P}_j
 \lesssim
 |\mathcal C|
\quad\Rightarrow\quad
 \#{\mathcal P}_j\lesssim |\mathcal C|\, 2^{\frac {(j+1)(d-1)}{d}},
\end{equation}
where we used that $h_T^{d-1} \approx |\omega_T \cap \mathcal C|$ since $T\cap \mathcal C \not = \emptyset$ for all marked elements. Recall that $\omega_T$ stands for the patch of elements around $T$.

On the other hand, the monotonicity of the local error indicators $\wt \osc_{\grid_i}(g,T)_q=h_T^t \| g \|_{L^q(T\cap \mathcal C)}$, implies that \REFINE does not increase $\wt \osc_{\grid_i}(g,T)_q$
and thus
$$
\delta \leq \delta_i = \wt \osc_{\grid_i}(g,T_i)_q,  \qquad 1 \leq i \leq k-1,
$$
where $t=\frac d2 - \frac1q (d-1)$.
In view of \eqref{e:dirac_surrogate}, if $T_i\in {\mathcal P}_j$, then we obtain
\[
	 \delta\le \wt \osc_{\grid_i}(g,T_i)_q \lesssim 2^{-\frac{jt}{d}}\|g\|_{L^q(T_i \cap \mathcal C)}.
\]
Therefore, accumulating these quantities in $\ell^q$ yields
\[
\delta^{q}\,\#{\mathcal P}_j\lesssim 2^{-\frac{jtq}{d}} \| g \|_{L^q(\mathcal C)}^q
\]
and gives rise to the second bound 
\begin{equation}\label{dirac:nv-large}
	\#{\mathcal P}_j\lesssim  \delta^{-q}\,2^{-\frac{jtq}{d}}\, \| g \|_{L^q(\mathcal C)}^q.
\end{equation}
\step{3} {\it Cardinality}. The two bounds for $\#{\mathcal P}$ in \eqref{dirac:nv-small} and
\eqref{dirac:nv-large} are complementary. The first one is good for $j$ small
whereas the second is suitable for $j$ large (think of $\delta\ll
1$). The crossover takes place for $j_0$ such that 
\[
2^{\frac{(j_0+1)(d-1)}{d}}|\mathcal C| \approx \delta^{-q}\,2^{-\frac{j_0 t q}{d}}\|g\|^q_{L^q(\mathcal C)}
\quad\Rightarrow\quad 
2^{j_0}\approx |\mathcal C|^{-\frac{2}{q}}  \delta^{-2} \|g\|^2_{L^q(\mathcal C)},
\]
upon using the expression for $t$.
We now compute 
\[
k=\#\marked=\sum_j\#{\mathcal P}_j\lesssim
\sum_{j\le j_0}2^{\frac{j(d-1)}{d}}|\mathcal C|+ \delta^{-q}\,\|g\|^q_{L^q(\mathcal C)}
\sum_{j>j_0}2^{-\frac{tq}{d}j}.
\]
Since 
\[
\sum_{j\le j_0}2^{\frac{j(d-1)}{d}}\approx 2^{\frac{j_0(d-1)}{d}},
\qquad \sum_{j>j_0}(2^{-\frac{tq}{d}})^j\lesssim 2^{-\frac{tqj_0}{d}},
\]
we can write 
\[
\#\marked\lesssim |\mathcal C|^{1-\frac{2(d-1)}{qd}}\big(\delta^{-1} \|g\|_{L^q(\mathcal C)}\big)^{\frac{2(d-1)}{d}}.
\]
We finally apply Theorem \ref{nv-T:complexity-refine} (complexity of \REFINE) to arrive at 
\[
\#\grid_k-\#\grid_0\lesssim \#\marked\lesssim |\mathcal C|^{1-\frac{2(d-1)}{qd}}\big(\delta^{-1} \|g\|_{L^q(\mathcal C)}\big)^{\frac{2(d-1)}{d}},
\]
or equivalently
\[
\delta \lesssim |\mathcal C|^{\frac{d}{d-1}-\frac 2 q}  \|g\|_{L^q(\mathcal C)} \big(\#\grid - \#\grid_0 \big)^{-\frac{d}{2(d-1)}}.
\]

We deduce from \eqref{dirac:pgreedy-1} that 
\[
\tau \lesssim \delta (\#\grid)^{\frac{1}{2}} \lesssim |\mathcal C|^{\frac{d}{d-1}-\frac 2 q}  \|g\|_{L^q(\mathcal C)} \big(\#\grid - \#\grid_0 \big)^{-\frac{d}{2(d-1)}+\frac 1 2}
\]
or equivalently
\begin{equation}\label{e:data_dirac_complexity_init}
\#\grid_k - \#\grid_0 \lesssim  \| g \|_{L^q(\mathcal C)}^{2(d-1)} \tau^{-2(d-1)},
\end{equation}
From this we conclude that $f = g \delta_{\mathcal C} \in \mathbb F_{\frac 1{2(d-1)}}$ with $|f|_{\mathbb F_{\frac 1{2(d-1)}}} \lesssim \| g \|_{L^q(\mathcal C)}$ as desired.

\medskip
\step{4} {\it Starting from $\grid$}.
To derive similar properties for \GREEDY starting from $\grid \in \grids$, we proceed as in the proof of Corollary~\ref{c:greed_start_grid} (performance of greedy). We distinguish the output  $[\wt{\grid},\wt f]=\GREEDY(\grid_0,\tau,2,1,f)$ starting from $\grid_0$ and performing $b=1$ bisection per marked element with $[\wh{\grid},\wh f]=\GREEDY(\grid,\tau,2,b,f)$ starting from $\grid \in \grids$ and performing $b\geq 1$ bisections per marked element. Lemma~\ref{l:greedy_start_grid} ($\GREEDY$ starting from $\grid$) guarantees that $\GREEDY(\grid,\tau,2,b,f)$ terminates with
$
\wt E_{\wh \grid}(f)_{-1} \leq \tau.
$
Moreover, Lemma~\ref{l:greedy_start_grid} also ensures that the number of marked elements satisfies
\[
N\le  \#\wt{\grid}-\# \grid_0 \lesssim  \| g \|_{L^q(\mathcal C)}^{2(d-1)} \tau^{-2(d-1)},
\]
where we used \eqref{e:data_dirac_complexity_init} to derive the last inequality.
This ends the proof.
\end{proof}

%-----------------------------------------------------------------------------------
\subsubsection{The case $f = \div \bg$ with $\bg \in L^2(\Omega;\R^d)$}\label{ss:data_f_div}
%-----------------------------------------------------------------------------------
%
A characterization of distributions in
$H^{-1}(\Omega)$ is given in \cite[Section 5.9.1]{Evans:10}: they are of the form
\[
f = f_0 + \div \bg
\]
with $f_0\in L^2(\Omega), \bg \in L^2(\Omega;\R^d)$. Since we have already treated separately
the ubiquitous case $\bg = \vec{0}$ in Section~\ref{ss:data_f_Lp}, we consider now the case $f_0=0$. Therefore,
\begin{equation}\label{E:divg-def}
\langle f ,v \rangle = -\int_\Omega \bg \cdot \nabla v
\quad\forall v \in H^1_0(\Omega)
\end{equation}
gives the action of $f$ on $v$ and its norm is \cite[Section 5.9.1]{Evans:10}

\begin{equation}\label{E:divg-norm}
\| f \|_{H^{-1}(\Omega)} = \inf \big\{\|\bg\|_{L^2(\Omega)}: \bg \in L^2(\Omega;\R^d)
\textrm{ satisfies } \eqref{E:divg-def} \big\} .
\end{equation}

Since adding the curl of a smooth vector field to $\bg$ does not change \eqref{E:divg-def},
we realize that the actual computation of \eqref{E:divg-norm} is problematic. We assume here
that $\bg$ is given and simply deal directly with $\bg$ thereby exploiting the relation
\begin{equation}\label{E:f->g}
\| f \|_{H^{-1}(\Omega)}  \le \|\bg\|_{L^2(\Omega)};
\end{equation}
this leads to a surrogate estimator.
We first approximate $\bg$ by discontinuous piecewise
polynomials of degree $n_f\le n-1$, namely we compute the $L^2$-projection $\bg_\grid = \Pi_\grid \bg$ onto
$[\mathbb{S}^{n_f,-1}_\grid]^d$, then we let $f_\grid := \div \bg_\grid \in \F_\grid \subset H^{-1}(\Omega)$
be the approximation of $f$: \looseness=-1
\begin{equation*}
\langle f_\grid , v \rangle = - \sum_{T\in\grid} \int_T \div \bg_\grid v - \sum_{F\in\faces} \int_F \jump{\bg_\grid}\cdot\vec{n}_F v
\quad\forall v \in H^1_0(\Omega).
\end{equation*}
We see that for $z \in \vertices$, $f_\grid|_{\omega_z}$ has the form of a functional in $\mathbb F(\grid_z)$ (see \eqref{e:data_local_F}) with $q_T = \div \bg_\grid|_T \in \P_{n_f-1},
q_F = \jump{\bg_\grid} \cdot\vec{n}_F \in \P_{n_f}$
for all $T\in\grid, F\in\faces$, but with smaller polynomial degree than functions in $\F(\grid_z)$.
We next exploit the local near-best approximation \eqref{e:PT_near_best_data} to replace $P_\grid f$ by $f_\grid$
\begin{equation}\label{E:surrogate-divg}
\|f-P_\grid f\|_{H^{-1}(\omega_z)} \le C_P \| f - f_\grid \|_{H^{-1}(\omega_z)}
\le C_P \| \bg - \bg_\grid \|_{L^2(\omega_z)}
\end{equation}
by virtue of \eqref{E:f->g} with $\Omega$ replaced by $\omega_z$.
This leads to the surrogate element-wise oscillation $\wt \osc_{\grid}(\bg,T)_2:=\| \bg - \bg_\grid \|_{L^2(T)}$, which satisfies Assumption~\ref{a:monot_osc_data} (monotonicity of local oscillation).
We thus have the global surrogate
$$
\wt{E}_{\grid}(f)_{-1}^2:= \sum_{T\in \grid} \wt \osc_{\grid}(\bg,T)_2^2.
$$
\begin{corollary}[approximation class of $\div\bg$]\label{C:divg}
Let the initial subdivision $\gridk[0]$ of $\Omega\subset\R^d$ satisfy 
Assumption~\ref{A:initial-labeling} (initial labeling) and $\grid \in \grids$ be
any admissible refinement of $\grid_0$.
Let $\tau>0$ be the target tolerance and $b\geq 1$ be the number of bisections performed on each marked element.
Let $(\bg,s,0,p,2)$ satisfy Assumption~\ref{ass:parameter-greedy} (admissible set of parameters for {\rm\GREEDY}) with local oscillations $\{\wt \osc(\bg,T)_2\}_{T\in\grid}$.
Then  $[\wh \grid, \wh f]=\GREEDY(\grid,\tau,2,b,f)$ terminates in a finite number of steps with 
$\wt E_{\wh \grid}(f)_{-1} \leq \tau$, whence $E_{\wh \grid}(f)_{-1} \leq \Cdata\tau$.
Moreover, the number $N$ of marked elements by {\rm \GREEDY} satisfies
$$
N \lesssim \| \bg \|_{X^s_p(\Omega)}^{\frac{d}{s}} \tau^{-\frac{d}{s}}.
$$
%
%with a constant $C = C(p,s,d,\Omega,\grid_0)$. 
In particular $f = \div \bg\in \F_{\frac{s}{d}}$ with
\[
|f|_{\F_\frac{s}{d}} \lesssim \| \bg \|_{X^s_p(\Omega)}.
\]
\end{corollary}
\begin{proof}
Apply Corollary~\ref{c:greed_start_grid} (performance of \GREEDY) with $q=2$ to $\bg$.
\end{proof}
 
%------------------------------------------------------------------------------------------
\subsection{{\rm \DATA} module}\label{S:data}
%------------------------------------------------------------------------------------------
%
We summarize now in one single algorithm, called $\DATA$, all the developments in Sections
\ref{ss:approx_diff}, \ref{ss:approx_reac}, and \ref{S:right-hand-side}. We first recall that 
Corollaries \ref{c:complex_consteained_A} (complexity of constrained \GREEDY for $\bA$) and \ref{c:complex_greedy_const_c} (complexity of constrained \GREEDY for $c$) deliver 
piecewise polynomial approximations $(\wh{\bA},\wh{c})$ of the coefficients $(\bA,c)$
over an admissible mesh $\wh\grid$ that satisfy both the global errors estimates
\[
E_{\wh \grid}(\bA)_{q_A} \leq \Cdata \tau,
\quad
E_{\wh \grid}(c)_{q_c} \leq \Cdata \tau,
\]
where $2 \le q_A, q_c \le \infty$ are the corresponding integrability indices, as well as
the structural constraint \eqref{E:structural-assumption-wh}.

The situation for the load $f$ is more intricate due to the evaluation of the nonlocal norm 
$H^{-1}(\Omega)$,
which requires further structure of $f$ besides regularity. Section~\ref{S:right-hand-side} provides
three examples of practical significance that allow for computable surrogate errors $\wt{E}_\grid(f)_{-1}$ larger than the desired oscillations $E_\grid(f)_{-1}$.
Since these examples have different requirements for the approximation procedure to work, we gather
the salient structural points in the following assumption.

\begin{assumption}[structure of $f$]\label{a:rhs-data}
  \index{Assumptions!Structure of $f$}  
  Let $(s_f,p_f)$ denote the additional
    regularity-integrability indices of $f$ beyond the basic $H^{-1}$-regularity, which are required
    by Assumption \ref{ass:parameter-greedy} (admissible set of parameters for \GREEDY). 
    Let $|f|_{\wt X^{s_f}_{p_f}(\Omega;\grid_0)}$ be a measure of piecewise regularity of 
    $f$ in $\grid_0$ expressed below in terms of surrogates. Assume that 
    either one of the following cases holds and note that all accumulate local oscillations in $\ell^2$.
    
    \begin{itemize}
        \item $f\in L^q(\Omega)$, with $\frac{2d}{d+2} < q \leq \infty$. Let $\widetilde \osc_{\grid}(f,T)_q = h_T^{t_f} \| f - \wh  f \|_{L^{q}(T)}$ be the local oscillation with $t_f =1+d(\frac 1 2- \frac 1 q) \ge 0$ and  $(f,s_f,t_f,p_f,2)$ satisfy Assumption \ref{ass:parameter-greedy}, and set
        $|f|_{\wt X^{s_f}_{p_f}(\Omega;\grid_0)}:=|f|_{X^{s_f}_{p_f}(\Omega;\grid_0)}$.
        
        \item $f = g \delta_{\mathcal C}$ where $\mathcal C \subset \Omega$ is a Lipschitz hyper-surface  and $g \in L^q(\mathcal C)$ with $\frac{2(d-1)}{d} < q \le \infty$. Let
        $\wt \osc_\grid(g,T)_q=h_T^{r} \| g \|_{L^q(T\cap \mathcal C)}$ be the local oscillation with $r = \frac d 2 - \frac 1 q (d-1)>0$. Set $s_f = 0$, $t_f=d/(2(d-1))$,  $p_f=q$, and $|f|_{\wt X^{s_f}_{p_f}(\Omega;\grid_0)}:=\|g\|_{L^q(\mathcal C)}$.
        
        \item  $f = \textrm{div } {\bg}$ with $\bg\in L^2(\Omega;\R^d)$. Let $\wt \osc_\grid(f,T)_2=\| {\bg} - \Pi_\grid {\bg}\|_{L^2(T)}$ be the local oscillation, $t_f=0$,  and $(\bg,s_f,t_f,p_f,2)$ satisfy
        Assumption~\ref{a:monot_osc_data}, and set $|f|_{\wt X^{s_f}_{p_f}(\Omega;\grid_0)}:=\|{\bg} \|_{X^{s_f}_{p_f}(\Omega;\grid_0)}$.
    \end{itemize}
    
\end{assumption}

In all these cases, \GREEDY algorithms with tolerance $\tau>0$ reduce the surrogate error $\wt E_\grid(f)_{-1}^2$ and eventually guarantee that
$$
E_\grid(f)_{-1} \leq  \Cdata \tau,
$$
where $\Cdata\ge1$ is the constant appearing in Corollary~\ref{C:p-greedy-Lp}, Lemma~\ref{l:p-greedy-dirac}, or Corollary~\ref{C:divg} depending on Assumption \ref{a:rhs-data} (structure of $f$).

\begin{algo}[$\DATA$] \label{A:data}
\index{Algorithms!$\DATA$: procedure to approximates the data $\data=(\bA,c,f)$}
Given a tolerance $\tau>0$ and an arbitrary conforming grid $\grid\in\grids$, not necessarily
$\grid_0$, $\DATA$ finds a conforming refinement $\wh\grid\ge\grid$ of $\grid$ and approximate
data $\wh\data = (\wh\bA, \wh{c}, \wh{f}) \in \D_{\wh\grid}$ over $\wh\grid$ such that 
$$
\| \data - \wh\data \|_{D(\Omega)}= E_{\wh\grid}(A)_{q_A} + E_{\wh\grid}(c)_{q_c} +  E_{\wh\grid}(f)_{-1}  \leq \Cdata \tau.
$$
\begin{algotab}
  \> $[\wh \grid,\wh\data] = \DATA \, (\grid,\tau,\data)$\\
  \>  \> $[\grid_{\bA},\wh{\bA}]=\GREEDY(\grid,\tau/3,q_A,b,\bA)$\\
  \>  \> $\wh \bA = \APPLYCONSTRA(\grid_A,\alpha_1,\alpha_2,L,\wt \bA)$\\
  \>  \> Set $\halpha_1 = \frac 1 2 \alpha_1$ and $\halpha_2 = (1+4L)\frac {\alpha_2} 2$\\
  \>  \> $[\grid_c,\wtc]=\GREEDY(\grid_A,\tau/3,q_c,b,c)$\\
  \>  \> $\wh c = \APPLYCONSTRC(\grid_c,\halpha_1,L,\wtc)$\\
  \>  \> $[\wh\grid,\wh f]=\GREEDY(\grid_c,\tau/3,2,b,f)$\\
  \>  \> return $\wh\grid,\wh\data$
\end{algotab}
\end{algo}
Note that \DATA depends on the threshold parameter $L>2$ used in \APPLYCONSTRA and \APPLYCONSTRC, although for simplicity it is not listed among the input parameters.
 
The next result summarizes the properties of $\DATA$. 

\begin{corollary}[performance of $\DATA$]\label{C:data}
Let the initial subdivision $\gridk[0]$ of $\Omega\subset\R^d$ satisfy Assumption~\ref{A:initial-labeling}
(initial labeling) and $\grid \in \grids$ be any admissible refinement of $\grid_0$.
Let $b\geq 1$ be the number of bisections performed on each marked element.
Let the assumptions of Corollaries \ref{c:complex_consteained_A} and \ref{c:complex_greedy_const_c}
for the coefficients $(\bA,c)$ be valid, and let $f$ satisfy Assumption \ref{a:rhs-data}.

For any target tolerance $\tau>0$ and any threshold parameter $L>2$, $[\wh \grid,\wh\data] = \DATA \, (\grid,\tau,\data)$ terminates in a finite number of iterations and outputs $\wh\data$, $\wh \grid \in \grids$ such that $\wh \bA$ is symmetric and 
$$
\wh \bA \in M(\halpha_1,\halpha_2), \qquad \wh c \in R(\wh c_1,\wh c_2),
$$
where $\halpha_1, \halpha_2$ are given by \eqref{e:halphas} while $\hc_1, \hc_2$ are given by \eqref{e:def_hc}.
Moreover, there is a constant $\Cdata \ge 1$ such that {\rm \DATA} terminates with
$$
\| \data - \wh\data\|_{D(\Omega)} \leq \Cdata \tau,
$$
and the number $N$ of elements marked to construct $\wh{\grid}$
satisfies
\begin{equation}\label{E:elements-doerfler}
N \lesssim |\data|_{\X_{\frac{s_\data}{d}}}^{\frac{d}{s_\data}} \tau^{-\frac{d}{s_\data}}
\end{equation}
with
$
s_\data := \min\{s_A,s_c+t_c,s_f+t_f\}$, 
and 
$$
|\data|_{\X_{\frac{s_\data}{d}}} = \Big(|\bA|_{X^{s_A}_{p_A}(\Omega,\grid_0)}^{\frac{d}{s_A}} + |c|_{X^{s_c}_{p_c}(\Omega,\grid_0)}^{\frac{d}{s_c}} + |f|_{\wt X^{s_f}_{p_f}(\Omega,\grid_0)}^{\frac{d}{s_f}} \Big)^{\frac{s_\data}{d}} .
$$
\end{corollary}
\begin{proof}
Since the local oscillations for $\bA$ and $c$ satisfy Assumption \ref{a:monot_osc_data} (monotonicity of local oscillations), we deduce that global oscillations do not increase upon refinement, namely
for $\wh\grid\ge \grid_c \ge \grid_{\bA}$
$$
E_{\wh\grid}(\bA)_{q_A} + E_{\wh\grid}(c)_{q_c} \leq E_{\grid_A}(\bA)_{q_A}+ E_{\grid_c}(c)_{q_c}. $$
In view of Corollaries~\ref{c:complex_consteained_A} and ~\ref{c:complex_greedy_const_c}, this in turn
implies
$$
E_{\wh\grid}(\bA)_{q_A} + E_{\wh\grid}(c)_{q_c} \leq  \Cdata \frac 2 3\tau.
$$

For the load term $f$, we invoke Corollary~\ref{C:p-greedy-Lp}, Lemma~\ref{l:p-greedy-dirac}, or Corollary~\ref{C:divg}, depending on Assumption \ref{a:rhs-data} (structure of $f$), to infer that
$$
E_{\wh\grid}(f)_{-1} \leq \Cdata\frac 1 3 \tau.
$$
Hence,
$$
\| \data - \wh \data \|_{D(\Omega)} = E_{\wh\grid}(\bA)_{q_A} + E_{\wh\grid}(c)_{q_c} +
E_{\wh\grid}(f)_{-1} \leq \Cdata \tau
$$
as desired.
The complexity estimate \eqref{E:elements-doerfler} directly follows from the complexity estimates given in Corollaries~\ref{c:complex_consteained_A} and \ref{c:complex_greedy_const_c} for $(\bA,c)$,
and Corollary~\ref{C:p-greedy-Lp}, Lemma~\ref{l:p-greedy-dirac}, or Corollary~\ref{C:divg} 
for $f$ depending on its structure.
\end{proof}

Similar ideas apply to approximate non-vanishing Dirichlet data or boundary flux condition
for Robin or Neumann problems, but we do not elaborate on this.

%--------------------------------------------------------------------------------

\section{Mesh Refinement: The Bisection Method} \label{S:mesh-refinement} %\cite{NochettoVeeser:2012}
%  \rhn{(CC $\longrightarrow$ RHN)}

%{\color{brown}  
%\begin{itemize}
%\item
%  Chains and labeling for $d=2$.

%\item
%  Recursive bisection.
  
%\item
%  Complexity of bisection for conforming meshes.

%\item
%  Complexity of bisection for non-conforming meshes.

%\end{itemize}
%}

This section is devoted to the complexity analysis of $\REFINE$ for $\Lambda$-admissible triangulations. Precisely, we prove the existence of a constant $\Ccompl>0$ such that
\[
\index{Constants!$\Ccompl$: complexity of \REFINE constant}
\#\gridk - \#\gridk[0] \le \Ccompl \sum_{j=0}^{k-1} \# \markedk[j]\,, \qquad k \geq 0\,.
\]
This kind of result holds for conforming meshes ($\Lambda=0$) and was stated in Theorem
\ref{nv-T:complexity-refine}, and for nonconforming meshes ($\Lambda >1$) as anticipated in Theorem
\ref{T:nonconforming-meshes}. 
The results
of Sections \ref{SS:conforming} and
\ref{SS:nonconforming} are valid for $d=2$ but the proofs of  the cited theorems extend easily to $d>2$. We
refer to the survey \cite{NoSiVe:09} for a full discussion for
$d\ge2$.

%-------------------------------------------------------------------------------------
\subsection{Conforming meshes}\label{SS:conforming}
%-------------------------------------------------------------------------------------

%-------------------------------------------------------------------------------------
\subsubsection{Chains and labeling for $d=2$}\label{nv-SS:chains-labeling}
\begin{figure}[ht]
	\centering
	\includegraphics[scale=0.65]{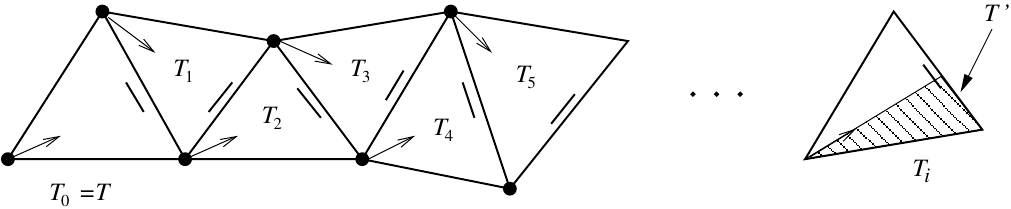}
        \caption{Typical chain $\mathcal{C}(T,\grid) = \{T_j\}_{j=0}^i$ emanating 
          from $T=T_0\in\grid$ with $T_j=F(T_{j-1}), j\ge1$.}\label{nv-F:dist-gen}
\end{figure}

%-------------------------------------------------------------------------------------
%
In order to study nonlocal effects of bisection for $d=2$ we introduce
now the concept of chain \cite{BiDaDeV:04}; this concept is not
adequate for $d>2$ \cite{NoSiVe:09,Stevenson:08}.  Recall that $E(T)$
denotes the edge of $T$ assigned for refinement.  To each $T\in\grid$
we associate the element $F(T)\in\grid$ sharing the edge $E(T)$ if
$E(T)$ is interior and $F(T)=\emptyset$ if $E(T)$ is on
$\partial\Omega$. A \emph{chain} $\mathcal{C}(T,\grid)$, with starting
element $T\in\grid$, is a sequence $\{T,F(T),\dots,F^m(T)\}$ with no
repetitions of elements and with
\[
F^{m+1}(T)=F^k(T)\text{ for some }k\in\{0,\dots,m-1\},\text{ or }F^{m+1}(T)=\emptyset;
\]
see Figure \ref{nv-F:dist-gen}. We observe that if an element $T$
belongs to two different grids, then the corresponding chains may be
different as well.
Two adjacent elements $T,T'=F(T)$ are \emph{compatibly divisible} (or
equivalently $T,T'$ form a \emph{compatible bisection patch}) if $F(T')=T$. Hence, 
$\mathcal{C}(T,\grid)=\{T,T'\}$ and a bisection of either $T$ or $T'$ does
not propagate outside the patch.

\medskip\noindent
\textbf{Example} (chains):
Let $\mathcal{F}=\{T_i\}_{i=1}^{12}$ be the forest of Figure \ref{nv-F:tree}. Then
$\mathcal{C}(T_6,\grid)=\{T_6,T_7\}, \mathcal{C}(T_9,\grid)=\{T_9\}$, and
$\mathcal{C}(T_{10},\grid)=\{T_{10},T_8,T_2\}$ are chains, but only
$\mathcal{C}(T_6,\grid)$ is a compatible bisection patch.

To study the structure of chains we rely on the initial labeling
\eqref{nv-initial-label} and the bisection rule of Section
\ref{S:bisection} (see Figure \ref{nv-F:labeling}):
\begin{equation}\label{nv-bisection-rule}
\begin{minipage}{0.8\linewidth}
{\it every triangle $T\in\grid$ with 
generation $g(T)=i$ receives the label $(i+1,i+1,i)$
with $i$ corresponding to the refinement edge $E(T)$, its side $i$ is
bisected and both new sides as well as the bisector
are labeled $i+2$ whereas the remaining labels do not change}.
\end{minipage}
\end{equation}
We first show that once the initial labeling and bisection rule are set,
the resulting master forest $\mathbb F$ is uniquely
determined: the label of an edge is independent of the elements
sharing this edge and no ambiguity arises in the recursion process. 

\begin{lemma}[labeling]\label{L:unique-labeling}
Let the initial labeling \eqref{nv-initial-label} for $\gridk[0]$ and
above bisection rule be enforced. 
If $\grid_0\le\grid_1\le\dots\le\grid_n$ are generated according to
\eqref{nv-bisection-rule},
then each side in $\grid_k$ has a unique label independent of the two 
triangles sharing this edge.
\end{lemma}
\begin{proof}
We argue by induction over $\grid_k$. For $k=0$ the assertion
is valid due to the initial labeling. Suppose the statement is true
for $\grid_k$. An edge $S$ in $\grid_{k+1}$ can be obtained in two
ways. The first is that $S$ is a bisector, and so a new edge, in which
case there is nothing to prove about its label being unique. The
second possibility is that $S$ was obtained by bisecting an edge
$S'\in\sides_k$. Let $T,\,T'\in\grid_k$ be the elements sharing $S'$,
and let us assume that $E(T')=S'$. Let $(i+1,i+1,i)$ be the label of
$T'$, which means that $S$ is assigned the label $i+2$. By induction
assumption over $\grid_k$, the label of $S'$ as an edge of $T$ is also
$i$. There are two possible cases for the label of $T$:
\begin{itemize}
	\item Label $(i+1,i+1,i)$: this situation is symmetric,
          $E(T)=S'$, and $S'$ is bisected with both halves getting
          label $i+2$. This is depicted in Fig. \ref{F:bisect-labeling-1}.
	
	\begin{figure}[ht]\label{F:bisect-labeling-1}
	\centering
	\includegraphics[scale=0.8]{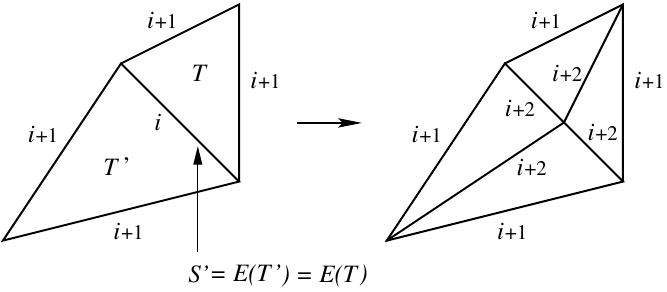}
	\caption{$T$ and $T'$ form a compatible patch, as they share the generation.}
	\end{figure}
	
	\item Label $(i,i,i-1)$: a bisection of side $E(T)$ with label
          $i-1$ creates a child $T''$ with label $(i+1,i+1,i)$ that
          is compatibly divisible with $T'$. Joining the new node of
          $T$ with the midpoint of $S'$ creates a conforming partition
          with level $i+2$ assigned to $S$. This is depicted in Fig. \ref{F:bisect-labeling-2}.
	
	\begin{figure}[ht]\label{F:bisect-labeling-2}
	\centering
	\includegraphics[scale=0.65]{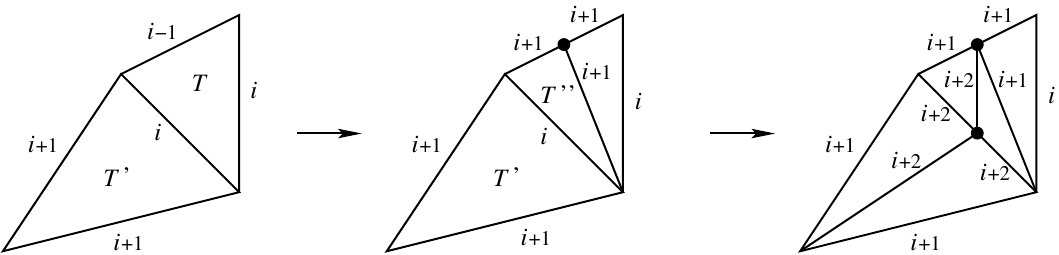}
	\caption{$T'$ form a compatible patch with the child $T''$ of $T$, indeed $T$ has a lower generation than $T'$.}
	\end{figure}

\end{itemize}
Therefore, in both cases the label $i+2$ assigned to $S$ is the same
from both sides, as asserted.
%\qed
\end{proof}

The two possible configurations displayed in the two figures above
lead readily to the following statement about generations.

\begin{corollary}[generation of consecutive elements]\label{nv-C:3_cor1}
For any $\grid\in\grids$ and $T,\,T'\in\grid$ with $T=F(T')$ we either have:
\begin{itemize}
	\item[(a)] $g(T)=g(T')$ and $T$, $T'$ are compatibly divisible, or
	\item[(b)] $g(T)=g(T')-1$ and $T'$ is compatibly divisible with a child of $T$.
\end{itemize}
\end{corollary}

\begin{corollary}[generations within a chain]\label{nv-C:3_cor2}
	For all $\grid\in\grids$ and $T\in\grid$, its chain
        $\mathcal{C}(T,\grid)=\{T_k\}_{k=0}^m$ with $T_k=F^k(T)$ have the
        property 
\[
  g(T_k)=g(T)-k\quad 0\le k\le m-1
\]
and $T_m=F^m(T)$ has generation $g(T_m)=g(T_{m-1})$ or it is a
boundary element with lowest labeled edge on $\partial\Omega$. 
In the first case, $T_{m-1}$ and $T_m$ are compatibly divisible.
\end{corollary}
\begin{proof} Apply Corollary \ref{nv-C:3_cor1} repeatly to
  consecutive elements of $\mathcal{C}(T,\grid)$. %\qed
\end{proof}

%-------------------------------------------------------------------------------------
\subsubsection{Recursive bisection}\label{nv-SS:recursive-bisection}
%-------------------------------------------------------------------------------------
%
Given an element $T\in\marked$ to be refined, the routine
$\textsf{REFINE\_RECURSIVE}~(\grid,T)$
recursively refines the chain $\mathcal{C}(T,\grid)$ of $T$, from the end 
back to $T$, and
creates a minimal conforming partition $\grid_*\ge\grid$ such that
$T$ is bisected once. This procedure reads as follows:
\begin{algotab}
  \> $[\grid_*]=\textsf{REFINE\_RECURSIVE}~(\grid,T)$\\
  \> \> if $g(F(T)) < g(T)$ \\
  \> \> \> $[\grid]= \textsf{REFINE\_RECURSIVE}~(\grid,F(T))$\\
  \> \> else\\
  \> \> \> bisect the compatible bisection patch $\mathcal{C}(T,\grid)$ \\
  \> \> \> update $\grid$ \\
  \> \> return~$\grid$
\end{algotab}

\noindent
We denote by
${\mathcal{C}_*(T,\grid)}\subset\grid_*$ the recursive refinement of
$\mathcal{C}(T,\grid)$ (or completion of $\mathcal{C}(T,\grid)$) caused by
bisection of $T$. Since $\textsf{REFINE\_RECURSIVE}$ refines solely
compatible bisection patches, intermediate meshes are always conforming.

\begin{figure}[ht]
	\centering
	\includegraphics[scale=0.7]{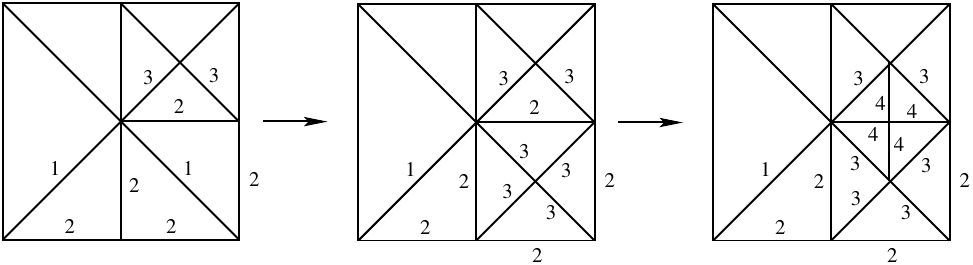}
\caption{Recursive refinement of $T_{10}\in\grid$ in Figure
\ref{nv-F:sequence} by $\textsf{REFINE\_RECURSIVE}$. This
entails refining the chain $\mathcal{C}(T_{10},\grid)=\{T_{10},T_8,T_2\}$,
starting from the last element $T_2\in\grid$, which form alone a compatible bisection
patch because its refinement edge is on the boundary, and continuing
with $T_8\in\grid$ and finally $T_{10}\in\grid$. Note that the successive meshes are
always conforming and that $\textsf{REFINE\_RECURSIVE}$ bisects
elements in $\mathcal{C}(T_{10},\grid)$ twice before getting back to $T_{10}$.}
\label{nv-F:recursive}
\end{figure}

We refer to Figure \ref{nv-F:recursive} for an example of recursive
bisection ${\mathcal{C}_*(T_{10},\grid)}$ of
$\mathcal{C}(T_{10},\grid)=\{T_{10},T_8,T_2\}$ in Figure
\ref{nv-F:sequence}: $\textsf{REFINE\_RECURSIVE}$ starts bisecting
from the end of $\mathcal{C}(T_{10},\grid)$, namely $T_2$, which is a
boundary element, and goes back the chain bisecting elements twice
until it gets to $T_{10}$.

We now establish a fundamental property of
$\textsf{REFINE\_RECURSIVE}~(\grid,T)$ relating the generation of elements
within  ${\mathcal{C}_*(T,\grid)}$ \cite{BiDaDeV:04}.

\begin{lemma}[recursive refinement]\label{nv-L:refpatch}
Let $\gridk[0]$ satisfy the labeling \eqref{nv-initial-label}, and let
$\grid\in\grids$ be a conforming refinement of $\gridk[0]$.
A call to {REFINE\_RECURSIVE}~$(\grid,T)$
terminates, for all $T$ in the set $\marked$ of marked elements, 
and outputs the smallest conforming refinement $\gridk[*]$
of $\grid$ such that $T$ is bisected. In addition, all newly created
$T'\in{\mathcal{C}_*(T,\grid)}$ satisfy
\begin{equation}\label{nv-generations}
	g(T')\le g(T)+1.
\end{equation}
\end{lemma}	
\begin{proof}
We first observe that $T$ has maximal generation within
$\mathcal{C}(T,\grid)$. So recursion is applied to elements with generation
$\le g(T)$, whence the recursion terminates.
We also note that this procedure creates children of $T$ and either
children or grandchildren of triangles
$T_k\in\mathcal{C}(T,\grid)=\{T_i\}_{i=0}^m$ 
with $k\ge1$. If $T'$ is a child of $T$ there is nothing to
prove. If not, we consider first $m=1$, in which case $T'$ is a
child of $T_1$ because $T_0$ and $T_1$ are compatibly
divisible and so have the same generation; thus
$g(T')=g(T_1)+1=g(T_0)+1$. 
Finally, if $m>1$, then $g(T_k)<g(T)$ and
we apply Corollary \ref{nv-C:3_cor2} to deduce
\[
g(T')\le g(T_k)+2\le g(T)+1,
\]
as asserted.	%\qed
\end{proof}

The following crucial lemma links generation and distance between $T$
and $T'\in{\mathcal{C}_*(T,\grid)}$, the latter being defined as \cite{BiDaDeV:04}
\[
\dist(T',T) := \inf_{x'\in T',x\in T}|x'-x|.
\]
\begin{lemma}[distance and generation]\label{nv-L:dist-gen}
Let $T\in\marked$. Any newly created $T'\in{\mathcal{C}_*(T,\grid)}$ by
$\textsf{REFINE\_RECURSIVE}~(\grid,T)$ satisfies 
\begin{equation}\label{nv-dist-gen}
\dist(T',T) \le D_2\frac{2}{\sqrt{2}-1}\,2^{-g(T')/2},
\end{equation}
where $D_2>0$ is the constant in \eqref{nv-diam-gen}.
\end{lemma}
\begin{proof} Suppose $T'\subset T_i\in\mathcal{C}(T,\grid)$ has been
created by subdividing $T_i$ (see Figure \ref{nv-F:dist-gen}). 
If $i\le1$ then $\dist(T',T)=0$ and there is nothing to prove. If
$i>1$, then we observe that $\dist(T',T_{i-1})=0$, whence
\begin{align*}
  \dist(T',T)&\le \dist(T_{i-1},T)+\diam(T_{i-1})\le\sum_{k=1}^{i-1}\diam(T_k)\\
  &\le D_2\,\sum_{k=1}^{i-1} 2^{-g(T_k)/2}<D_2\,\frac{1}{1-2^{-1/2}}\,2^{-g(T_{i-1})/2},
\end{align*}
because the generations decrease exactly by $1$ along the chain
$\mathcal{C}(T)$ according to Corollary \ref{nv-C:3_cor1}(b). 
Since $T'$ is a child or grandchild of $T_i$, we deduce
\[
g(T')\le g(T_i)+2 = g(T_{i-1})+1, 
\]
whence 
\[
\dist(T',T)< D_2\,\frac{2^{1/2}}{1-2^{-1/2}}\,2^{-g(T')/2}. 
\]
This is the desired estimate.	%\qed
\end{proof}

The recursive procedure $\textsf{REFINE\_RECURSIVE}$ is the core of 
the routine $\REFINE$ of Section \ref{S:bisection}:
given a conforming mesh $\grid\in\grids$ and a subset
$\marked\subset\grid$ of marked
elements, $\REFINE$ creates a conforming refinement 
$\gridk[*]\ge\grid$ of $\grid$ such that all elements of $\marked$ are
bisected at least once:
\begin{algotab}
  \> $[\grid_*]=\textsf{REFINE}~(\grid,\marked)$\\
  \> \> for all $\elm\in\marked\cap\grid$ do\\
  \> \> \> $[\grid]= \textsf{REFINE\_RECURSIVE}~(\grid,\elm)$\\
  \> \> return~$\grid$
\end{algotab}
It may happen that an element $T'\in\marked$ is scheduled prior to $T$ for 
refinement and $T\in\mathcal{C}(T',\grid)$. Since the call
$\textsf{REFINE\_RECURSIVE}~(\grid,T')$ bisects
$T$, its two children replace $T$ in $\grid$. This implies that
$T\notin\marked\cap\grid$, which prevents further refinement of $T$.

 In practice, one often likes to bisect selected elements several
  times, for instance each marked element is scheduled for $b\geq 1$
  bisections. This can be done by assigning the number $b(\elm)=b$ of
  bisections that have to be executed for each marked element $\elm$. If
  $\elm$ is bisected then we assign $b(\elm)-1$ as the number of
  pending bisections to its children and the set of marked elements is
  $\marked\definedas\{\elm\in\grid\mid b(\elm)>0\}$.

%-------------------------------------------------------------------------------------
\subsubsection{Complexity of bisection for conforming meshes}\label{nv-SS:proof-complexity-conforming}
%-------------------------------------------------------------------------------------
%
Figure \ref{nv-F:recursive} reveals that the issue of propagation of 
mesh refinement to keep conformity is
rather delicate. In particular, an estimate of
the form
\[
\#\gridk - \#\gridk[k-1] \le C \,  \# \markedk[k-1]
\]
is not valid with a constant $C$ independent of $k$; in fact the 
constant can be proportional to $k$ according to Figure
\ref{nv-F:recursive}.

Binev, Dahmen, and DeVore \cite{BiDaDeV:04} for $d=2$ 
and Stevenson \cite{Stevenson:08} for $d>2$ show that control of 
the propagation of refinement by bisection is possible when
considering the collective effect:
\begin{equation}\label{e:completion_2}
\#\gridk[k] - \#\gridk[0] \le D \sum_{j=0}^{k-1} \#\markedk[j].
\end{equation}
This can be heuristically motivated as follows.
Consider the set
\begin{math}
  \marked\definedas\bigcup_{j=0}^{k-1}\marked_j
\end{math}
used to generate the sequence 
\begin{math}
  \gridk[0]\leq\gridk[1]\leq\dots\leq\gridk[k]\asdefined\grid
\end{math}.
Suppose that each element $\elm_*\in\marked$ is assigned a fixed
amount $C_1$ of money to spend on refined elements in $\grid$, \ie
on $\elm\in\grid\setminus\gridk[0]$. Assume
further that $\lambda(\elm,\elm_*)$ is the portion of money spent by
$\elm_*$ on $\elm$. Then it must hold
\begin{subequations}\label{nsv-lambda}
  \begin{equation}\label{nsv-lambda.a}
    \sum_{\elm\in\grid\setminus\gridk[0]}\lambda(\elm,\elm_*)\leq C_1
    \qquad\text{for all } \elm_*\in\marked.
  \end{equation}
  In addition, we suppose that the investment of all elements in $\marked$ is fair in
  the sense that each $\elm\in\grid\setminus\gridk[0]$ gets at least a fixed
  amount $C_2$, whence
  \begin{equation}\label{nsv-lambda.b}
    \sum_{\elm_*\in\marked}\lambda(T,\elm_*)\ge C_2
    \qquad\text{for all } \elm\in\grid\setminus\gridk[0].
  \end{equation}
\end{subequations}
Therefore, summing up \eqref{nsv-lambda.b} and using the upper bound
\eqref{nsv-lambda.a} we readily obtain
\begin{displaymath}
  C_2(\#\grid-\#\gridk[0])\leq
  \sum_{\elm\in\grid\setminus\gridk[0]}\sum_{\elm_*\in\marked}\lambda(\elm,\elm_*)
  =
  \sum_{\elm_*\in\marked}\sum_{\elm\in\grid\setminus\gridk[0]}\lambda(\elm,\elm_*)
  \leq
  C_1\,\#\marked,
\end{displaymath}
which proves \eqref{e:completion_2} for $\grid$ and $\marked$.
In the remainder of this
section we design such an allocation function 
\begin{math}
  \lambda\colon\grid\times\marked\to\R^+
\end{math}
in several steps and prove that recurrent refinement by bisection yields
\eqref{nsv-lambda} provided $\gridk[0]$ satisfies \eqref{nv-initial-label}, thereby establishing Theorem \ref{nv-T:complexity-refine} (complexity of $\REFINE$).

\paragraph{Construction of the Allocation Function.}
The function $\lambda(T,T_*)$ is defined with the help of two sequences
$\big(a(\ell)\big)_{\ell=-1}^\infty$,
$\big(b(\ell)\big)_{\ell=0}^\infty\subset\R^+$ of positive numbers
satisfying
\begin{displaymath}
  \sum_{\ell\ge -1} a(\ell) = A < \infty,\qquad 
  \sum_{\ell\ge 0} 2^{-\ell/2}\,b(\ell) = B <\infty,
  \qquad
  \inf_{\ell\geq 1} b(\ell)\,a(\ell) = c_* > 0,
\end{displaymath}
and $b(0)\geq 1$.  Valid instances are $a(\ell)=(\ell+2)^{-2}$ and
$b(\ell)=2^{\ell/3}$.

With these settings we are prepared to define
\begin{math}
  \lambda\colon\grid\times\marked\to\R^+
\end{math}
by 
\begin{displaymath}
  \lambda(\elm,\elm_*)
  \definedas
  \begin{cases} 
    a(g(\elm_*)-g(\elm)), &\dist(\elm,\elm_*) < D_3\,B\,2^{-g(\elm)/d}
    \text{ and }g(\elm)\leq g(\elm_*)+1\\ 
    0,& \text{else},
  \end{cases}
\end{displaymath}
where
\begin{math}
%  D_3\definedas D_2\left(1+ \frac{2^{1/d}}{1-2^{-1/d}}\right).
  D_3\definedas D_2\big(1+ 2(\sqrt{2}-1)^{-1}\big).
\end{math}
Therefore, the investment of money by $\elm_*\in\marked$ is restricted to cells
$\elm$ that are sufficiently close and are of generation $g(\elm)\leq g(\elm_*)+1$. Only elements of these generations can be created during refinement of 
$\elm_*$ according to Lemma~\ref{nv-L:refpatch}.
We stress that except for the definition of $B$, this construction is
mutidimensional and we refer to \cite{NoSiVe:09,Stevenson:08} for details.

The following lemma shows that the total amount of money spend by
the allocation function $\lambda(T,T_*)$ per marked element $T_*$ is bounded.

\begin{lemma}[upper bound]\label{nv-lambda_upper}
  There exists a constant $C_1>0$ only depending on $\gridk[0]$ such that
  $\lambda$ satisfies \eqref{nsv-lambda.a}, \ie
  \begin{displaymath}
    \sum_{\elm\in\grid\setminus\gridk[0]}\lambda(\elm,\elm_*)\leq C_1
    \qquad\text{for all } \elm_*\in\marked.
  \end{displaymath}
\end{lemma}
\begin{proof} We proceed in two steps.

  \step{1} Given $\elm_*\in\marked$ we set $g_*=g(\elm_*)$ and we let
  $0\leq g\leq g_*+1$ be a generation of interest in the definition of
  $\lambda$. We claim that for such $g$ the cardinality of the set
  \begin{displaymath}
    \grid(\elm_*,g) = \{\elm\in\grid \mid
    \dist(\elm,\elm_*) < D_3\,B\,2^{-g/2} \text{ and }g(\elm) = g\}
  \end{displaymath}
  is uniformly bounded, \ie $\#\grid(\elm_*,g)\leq C$ with $C$ solely
  depending on $D_1,D_2,D_3,B$.

  From \eqref{nv-diam-gen} we learn that
  \begin{math}
    \diam(\elm_*)\leq 
    D_2 2^{-g_*/2}\leq 2D_2 2^{-(g_*+1)/2} \leq 2 D_2 2^{-g/2}
  \end{math}
  as well as 
  \begin{math}
    \diam(\elm) \leq D_2 2^{-g/2}
  \end{math}
  for any $\elm\in\grid(\elm_*,g)$. Hence, all elements of the set
  $\grid(\elm_*,g)$ lie inside a ball centered at the barycenter of $\elm_*$
  with radius $(D_3B + 3D_2)2^{-g/2}$. Again relying on \eqref{nv-diam-gen}
  we thus conclude
  \begin{displaymath}
    \#\grid(\elm_*,g) D_12^{-g} \leq \sum_{\elm\in\grid(\elm_*,g)} \vol{\elm}
    \leq c (D_3B + 3D_2)^2 2^{-g},
  \end{displaymath}
  whence $\#\grid(\elm_*,g)\leq c\,D_1^{-1}\,(D_3B + 3D_2)^2=:C$.

  \step{2} Accounting only for non-zero contributions 
  $\lambda(\elm,\elm_*)$ we deduce
  \begin{align*}
    \sum_{\elm\in\grid\setminus\gridk[0]}\lambda(\elm,\elm_*)
    &=
    \sum_{g=0}^{g_*+1} \sum_{\elm\in\grid(\elm_*,g)} a(g_* - g)
    \leq C \sum_{\ell=-1}^\infty a(\ell) = CA\asdefined C_1,
  \end{align*}
  which is the desired upper bound. %\qed
\end{proof}

The definition of $\lambda$ also implies that each refined element
receives a fixed amount of money. We show this next.

\begin{lemma}[lower bound]\label{nsv-lambda_lower}
  There exists a constant $C_2>0$ only depending on $\gridk[0]$ such that 
  $\lambda$ satisfies \eqref{nsv-lambda.b}, \ie
  \begin{displaymath}
    \sum_{\elm_*\in\marked}\lambda(\elm,\elm_*)\geq C_2\qquad\text{for all }
    \elm\in\grid\setminus\gridk[0].
  \end{displaymath}
\end{lemma}
\begin{proof} We proceed in several steps.

  \step{1} Fix an arbitrary $\elm_0\in\grid\setminus\gridk[0]$. Then
  there is an iteration count $1\leq k_0 \leq k$ such that
  $\elm_0\in\gridk[k_0]$ and $\elm_0\notin\gridk[k_0-1]$.  Therefore
  there exists an $\elm_1\in\markedk[k_0-1]\subset\marked$
  such that $\elm_0$ is generated during
  $\textsf{REFINE\_RECURSIVE}~(\gridk[k_0-1],\elm_1)$. Iterating this
  process we construct a sequence $\{\elm_j\}_{j=1}^J\subset\marked$
  with corresponding iteration counts $\{k_j\}_{j=1}^J$ such that
  $\elm_{j}$ is created by
  $\textsf{REFINE\_RECURSIVE}~(\gridk[k_j-1],\elm_{j+1})$. The sequence
  is finite since the iteration counts are strictly decreasing and
  thus $k_J=0$ for some $J>0$, or equivalently $\elm_J\in\gridk[0]$.

  Since $\elm_j$ is created during refinement of $\elm_{j+1}$ we infer
  from \eqref{nv-generations} that
  \begin{displaymath}
    g(\elm_{j+1})\geq g(\elm_j)-1.
  \end{displaymath}
  Accordingly, $g(\elm_{j+1})$ can decrease the previous value of
  $g(\elm_j)$ at most by $1$. Since $g(\elm_J)=0$ there exists a
  smallest value $s$ such that $g(\elm_s)=g(\elm_0)-1$. Note that for
  $j=1,\dots, s$ we have $\lambda(\elm_0,\elm_j)>0$ if
  $\dist(\elm_0,\elm_j)\leq D_3B g^{-g(\elm_0)/d}$.

  \step{2} We next estimate the distance $\dist(\elm_0,\elm_j)$. For
  $1\le j\le s$ and $\ell\geq 0$ we define the set
  \begin{displaymath}
    \grid(\elm_0,\ell,j)\definedas 
    \{\elm\in \{\elm_0,\dots,\elm_{j-1}\}\mid g(\elm)=g(\elm_0)+\ell\}
  \end{displaymath}
  and denote by $m(\ell,j)$ its cardinality. 
  The triangle inequality combined with an induction argument yields
  \begin{align*}
    \dist(\elm_0,\elm_j)
    &\leq
    \dist(\elm_0,\elm_1) + \diam(\elm_1) + \dist(\elm_1,\elm_j)\\
    &\leq
    \sum_{i=1}^j \dist(\elm_{i-1},\elm_i) + \sum_{i=1}^{j-1} \diam(\elm_i).
  \end{align*}
  We apply \eqref{nv-dist-gen} for the terms of the first sum
  and \eqref{nv-diam-gen} for the terms of the second sum to obtain
  \begin{align*}
    \dist(\elm_0,\elm_j)
    &<
    D_2\frac{2}{\sqrt{2}-1}\,\sum_{i=1}^j 2^{-g(\elm_{i-1})/2}
    + D_2 \sum_{i=1}^{j-1} 2^{-g(\elm_i)/2}\\
    &\le D_2\left(1+ \frac{2}{\sqrt{2}-1}\right)
    \sum_{i=0}^{j-1} 2^{-g(\elm_i)/2}\\
    &=
    D_3 \sum_{\ell=0}^\infty m(\ell,j)\,2^{-(g(\elm_{0})+\ell)/2}\\
    &=
    D_3 2^{-g(\elm_{0})/2} \sum_{\ell=0}^\infty m(\ell,j)\,2^{-\ell/2}.
  \end{align*}
  For establishing the lower bound we distinguish two cases
  depending on the size of $m(\ell,s)$. This is done next.
  
  \step{3} \emph{Case 1:} $m(\ell,s)\leq b(\ell)$ for all $\ell\geq
  0$. From this we conclude
  \begin{displaymath}
    \dist(\elm_0,\elm_s) <
    D_3 2^{-g(\elm_0)/2} \sum_{\ell=0}^\infty b(\ell)\,2^{-\ell/2}
    = D_3B\, 2^{-g(\elm_0)/2}
  \end{displaymath}
  and the definition of $\lambda$ then readily implies
  \begin{displaymath}
    \sum_{\elm_*\in\marked} \lambda(\elm_0,\elm_*) \geq
    \lambda(\elm_0,\elm_s) = a(g(\elm_s)-g(\elm_0)) = a(-1)>0.
  \end{displaymath}

  \step{4} \emph{Case 2:} There exists $\ell\ge 0$ such that $m(\ell,s)>
  b(\ell)$.  For each of these $\ell$'s there exists a smallest
  $j=j(\ell)$ such that $m(\ell,j(\ell))>b(\ell)$. We let $\ell^*$ be
  the index $\ell$ that gives rise to the smallest $j(\ell)$, and set
  $j^*=j(\ell^*)$. Consequently
  \begin{displaymath}
    m(\ell,j^*-1)\le b(\ell) \quad\text{for all } \ell\geq 0
    \qquad\text{and}\qquad m(\ell^*,j^*) > b(\ell^*).
  \end{displaymath}
  As in Case 1 we see
  \begin{math}
    \dist(\elm_0,\elm_i)< D_3B\, 2^{-g(\elm_0)/2}
  \end{math}
  for all $i\leq j^*-1$, or equivalently 
  \begin{displaymath}
    \dist(\elm_0,\elm_i)< D_3B\, 2^{-g(\elm_0)/2}
    \qquad \text{for all } \elm_i\in\grid(\elm_0,\ell,j^*).
  \end{displaymath}
  
  We next show that the elements in $\grid(\elm_0,\ell^*,j^*)$ spend
  enough money on $T_0$. We first consider $\ell^*=0$ and note that
  $\elm_0\in\grid(\elm_0,0,j^*)$. Since $m(0,j^*)>b(0)\geq 1$ we discover
  $j^*\geq 2$. Hence, there is an
  $\elm_i\in\grid(\elm_0,0,j^*)\cap\marked$, which yields the estimate
  \begin{displaymath}
    \sum_{\elm_*\in\marked} \lambda(\elm_0,\elm_*) \geq
    \lambda(\elm_0,\elm_i) = a(g(\elm_i)-g(\elm_0)) = a(0)>0.
  \end{displaymath}

  For $\ell^*>0$ we see that $\elm_0\not\in\grid(\elm_0,\ell^*,j^*)$, whence
  $\grid(\elm_0,\ell^*,j^*)\subset\marked$. In addition,
  $\lambda(\elm_0,\elm_i)=a(\ell^*)$ for all
  $\elm_i\in\grid(\elm_0,\ell^*,j^*)$. From this we conclude
  \begin{displaymath}
    \begin{aligned}
      \sum_{\elm_*\in\marked} \lambda(\elm_0,\elm_*) &\geq
      \sum_{\elm_*\in\grid(\elm_0,\ell^*,j^*)} \lambda(\elm_0,\elm_*)
      = m(\ell^*,j^*)\,a(\ell^*) \\
      & > b(\ell^*)\,a(\ell^*)
      \geq \inf_{\ell\geq 1} b(\ell)\,a(\ell)=c_*>0.
    \end{aligned}
  \end{displaymath}

 \step{5}  In summary we have proved the assertion 
 since for any $\elm_0\in\grid\setminus\gridk[0]$
  \begin{equation}%\tag*{\qed}
    \sum_{\elm_*\in\marked} \lambda(\elm_0,\elm_*) \geq \min\{a(-1), a(0), c_*\}
    \asdefined C_2 > 0.
  \end{equation}
  This completes the proof.
\end{proof}

\begin{remark}[complexity with $b>1$ bisections]\label{nv-R:complexity-b>1}
  To show the complexity estimate when $\REFINE$ performs $b>1$ bisections,
  the set $\markedk$ is to be understood as a sequence of \textit{single}
  bisections recorded in sets $\{\markedk(j)\}_{j=1}^b$, which belong to
  intermediate triangulations between $\gridk$ and $\gridk[k+1]$ with
  $\#\markedk(j)\leq 2^{j-1}\#\markedk$, $j=1,\dots,b$. Then we also 
  obtain Theorem~\ref{nv-T:complexity-refine}  because
  \begin{displaymath}
    \sum_{j=1}^b \# \markedk(j) \leq \sum_{j=1}^b 2^{j-1} \# \markedk = (2^b-1)\#\markedk.
  \end{displaymath}
  In practice, it is customary to take $b=d$ \cite{Siebert:12}.
\end{remark}

%-------------------------------------------------------------------------------------
\subsection{Nonconforming meshes}\label{SS:nonconforming}
%-------------------------------------------------------------------------------------
In this subsection, we consider two kinds of nonconforming meshes undergoing a refinement process: a) quadrilateral meshes with at most one hanging node per edge ($\Lambda =1$ in the definition of $\Lambda$-admissible meshes), and b) triangular meshes having global index bounded by a fixed, but arbitrary $\Lambda>1$.

%-------------------------------------------------------------------------------------
\subsubsection{Complexity of bisection for nonconforming quadrilateral meshes}
\label{nv-SS:proof-complexity-quad-nonconforming}
%-------------------------------------------------------------------------------------
%
We examine briefly the refinement process for quadrilaterals with
one hanging node per edge, which gives rise to the so-called {\it 1-meshes}.
The refinement of $T\in\grid$ might affect four elements
of $\grid$ for $d=2$ (or $2^d$ elements for any dimension $d\ge2$), all contained in
the {\it refinement patch} $R(T,\grid)$ of $T$ in $\grid$. The latter is defined as
\[
R(T,\grid) := \{T'\in\grid | ~ T' \textrm{ and } T \textrm{ share an
  edge and } g(T') \le g(T) \},
\]
and is called {\it compatible} provided $g(T')=g(T)$ for all $T'\in
R(T,\grid)$. The generation gap between elements sharing an edge, in particular 
those in $R(T,\grid)$, is always $\le 1$ for $1$-meshes, and
is $0$ if $R(T,\grid)$ is compatible. The element size satisfies
\begin{equation*}
h_T = 2^{-g(T)} h_{T_0}
\qquad\forall T\in\grid
\end{equation*}
where $T_0\in\grid_0$ is the ancestor of $T$ in the initial mesh
$\grid_0$. Lemma \ref{nv-L:diam-gen} is thus valid
\begin{equation}\label{nv-diam-gen-quad}
h_T < \bar{h}_T \le D_2 2^{-g(T)}
\qquad\forall T\in\grid.
\end{equation}
Given an element $T\in\marked$ to be refined, the routine
$\textsf{REFINE\_RECURSIVE}~(\grid,T)$ refines recursively
$R(T,\grid)$ in such a way that the intermediate meshes are always
$1$-meshes, and reads as follows:
\begin{algotab}
  \> $[\grid_*] = \textsf{REFINE\_RECURSIVE}~(\grid,T)$\\
  \> \> if $g=\min\{g(T''): T''\in R(T,\grid)\} < g(T)$ \\
  \> \> \> let $T'\in R(T,\grid)$ satisfy $g(T')=g$ \\
  \> \> \> $[\grid] = \textsf{REFINE\_RECURSIVE}~(\grid,T')$\\
  \> \> else\\
  \> \> \> subdivide $T$ \\
  \> \> \> update $\grid$ upon replacing $T$ by its children \\
  \> \> return~$\grid$
\end{algotab}
The conditional prevents the generation gap within $R(T,\grid)$
from getting larger than $1$.
If it fails, then the refinement patch $R(T,\grid)$ is
compatible and refinining $T$ increases the generation gap from $0$ 
to $1$ without violating the $1$-mesh structure.
This implies a variant of Lemma \ref{nv-L:refpatch}: $\textsf{REFINE\_RECURSIVE}~(\grid,T)$ creates a minimal $1$-mesh
$\grid_*\ge\grid$ refinement of $\grid$ so that for all newly created elements $T'\in\grid_*$
\begin{equation}\label{nv-generations-quad}
g(T') \le g(T) + 1
\end{equation}
and $T$ is subdivided only
{\it once}. This yields Lemma \ref{nv-L:dist-gen}: there exist a geometric
constant $D_g>0$ such that for all newly created elements $T'\in\grid_*$
\begin{equation}\label{nv-dist-gen-quad}
\dist(T,T') \le D_g 2^{g(T')}.
\end{equation}

The procedure $\textsf{REFINE\_RECURSIVE}$ is the core of
$\textsf{REFINE}$, which is conceptually identical to that in
Section \ref{nv-SS:recursive-bisection}.
Suppose that each marked element $T\in\marked$ is
to be subdivided $b\ge1$ times. We assign a flag $q(T)$ to each
element $T$ which is
initialized $q(T)=b$ if $T\in\marked$ and $q(T)=0$ otherwise. 
The marked set $\marked$ is then the set of
elements $T$ with $q(T)>0$, and every time $T$ is subdivided it is removed 
from $\grid$ and replaced by its children, which inherit the flag
$q(T)-1$.
This avoids the conflict of subdividing again an element that has been
previously refined by $\textsf{REFINE\_RECURSIVE}$. The procedure
$\textsf{REFINE}~(\grid,\marked)$ reads
\begin{algotab}
\index{Algorithms!\REFINE: refine all marked elements $b$ times and others necessary to produce a conforming mesh}
  \> $[\grid_*]=\textsf{REFINE}~(\grid,\marked)$\\
  \> \> for all $\elm\in\marked\cap\grid$ do\\
  \> \> \> $[\grid]= \textsf{REFINE\_RECURSIVE}~(\grid,\elm)$;\\
  \> \> end\\
  \> \> return~$\grid$
\end{algotab}
and its output is a minimal $1$-mesh 
$\grid_*\ge\grid$, refinement of $\grid$, so that all marked elements of
$\marked$ are refined at least $b$ times. Since $\grid_*$ has one
hanging node per side it is thus admissible in the sense of
\eqref{nv-level-nonconformity}.
However, the refinement may spread outside $\marked$ and the
issue of complexity of $\REFINE$ again becomes non-trivial.

With the  above ingredients in place, a statement similar to Theorem \ref{nv-T:complexity-refine} (complexity of $\REFINE$) for nonconforming quadrilateral meshes follows along the lines of Section 
\ref{nv-SS:proof-complexity-conforming}.

%-------------------------------------------------------------------------------------
\subsubsection{Complexity of bisection for $\Lambda$-admissible triangular meshes}
\label{nv-SS:proof-complexity-Ltria-nonconforming}
%-------------------------------------------------------------------------------------
%

Let $\grid \in \grids^\Lambda$ be a $\Lambda$-admissible simplicial mesh. Given any $T\in \mesh$, let us denote again by $E(T)$ the edge of $T$ assigned for refinement, i.e., the edge opposite to the newest vertex $v(T)$. Let us denote by $x(T)$ the midpoint of the edge $E(T)$.

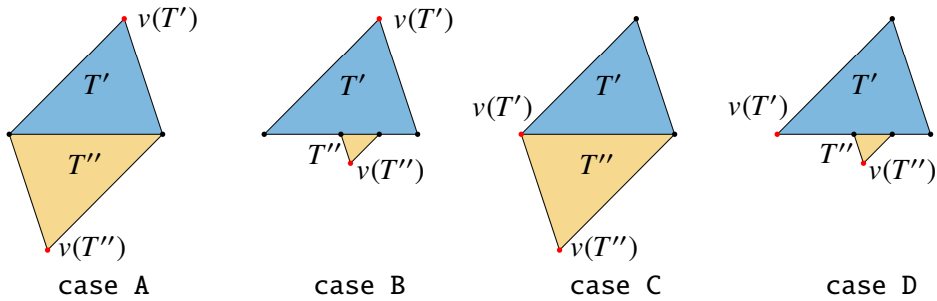
\begin{figure}[h!]
\begin{center}
\begin{overpic}[scale=0.19]{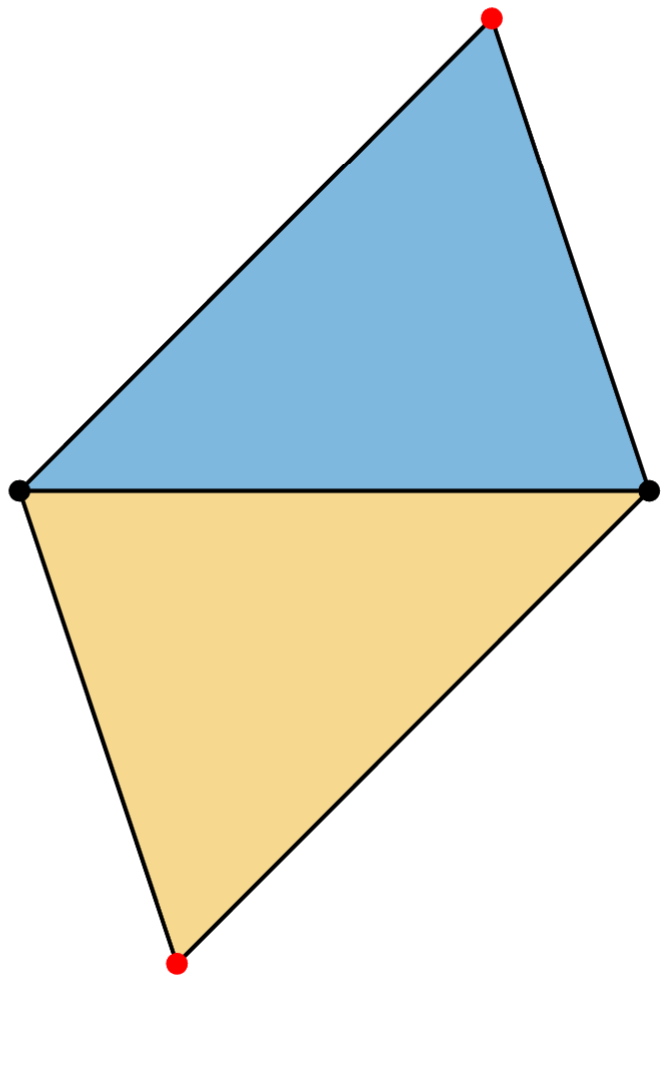}
\put(24,40){$T''$}
\put(30,70){$T'$}
\put(20,9){$v(T'')$}
\put(50,95){$v(T')$}
\put(20,-5){\texttt{case A}}
\end{overpic}
\qquad \quad
\begin{overpic}[scale=0.19]{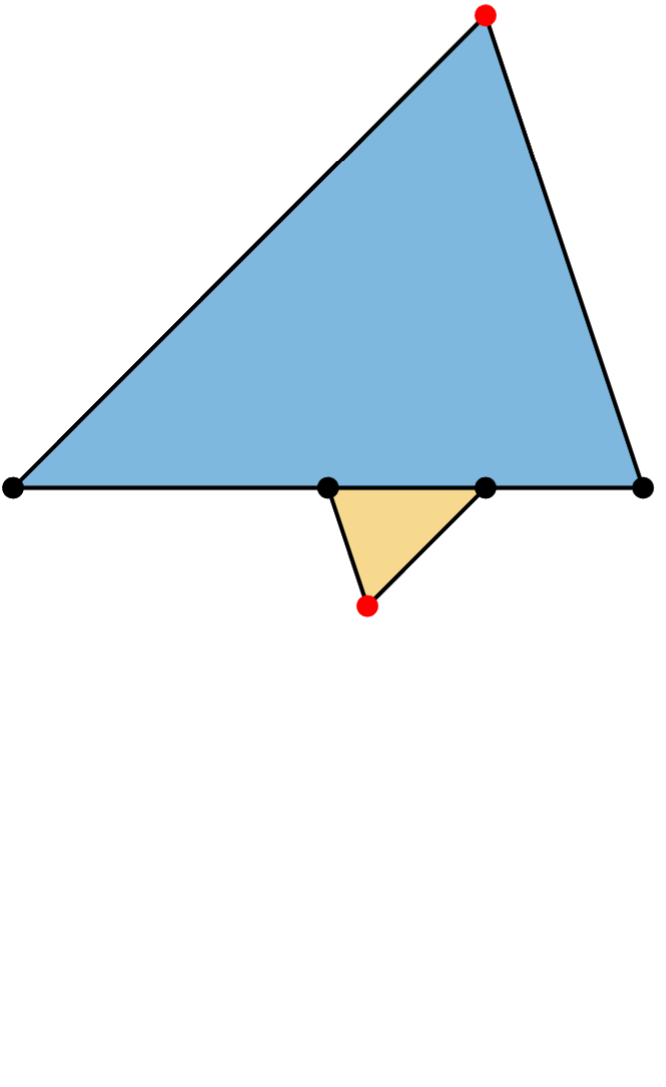}
\put(18,44){$T''$}
\put(30,70){$T'$}
\put(36,38){$v(T'')$}
\put(50,95){$v(T')$}
\put(20,-5){\texttt{case B}}
\end{overpic}
\qquad \quad
\begin{overpic}[scale=0.19]{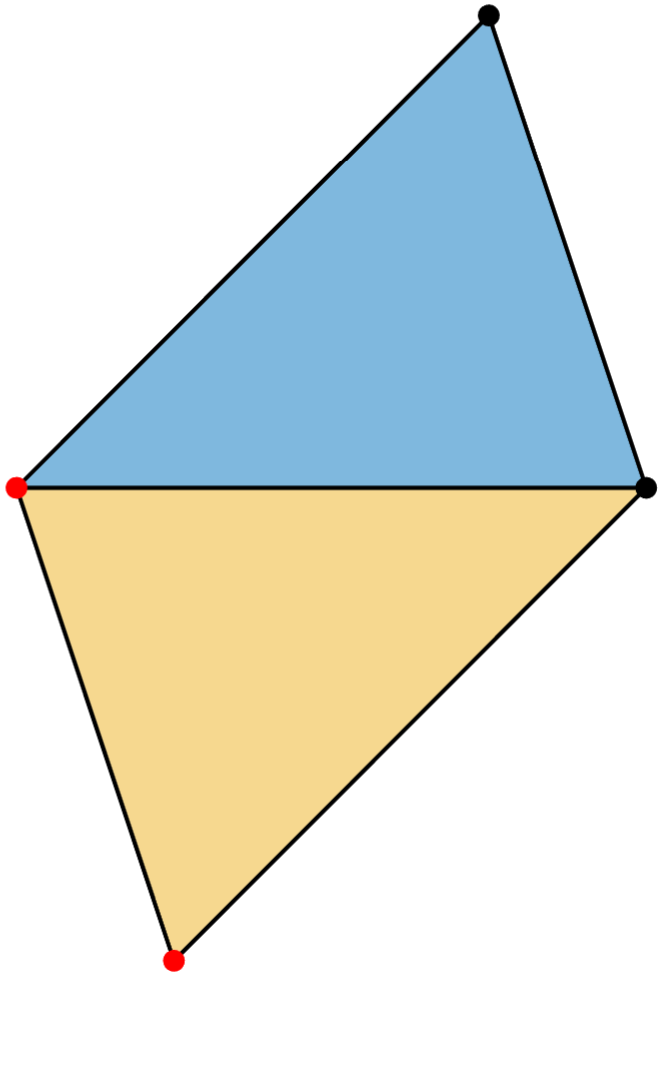}
\put(24,40){$T''$}
\put(30,70){$T'$}
\put(20,9){$v(T'')$}
\put(-17,62){$v(T')$}
\put(20,-5){\texttt{case C}}
\end{overpic}
\qquad \quad
\begin{overpic}[scale=0.19]{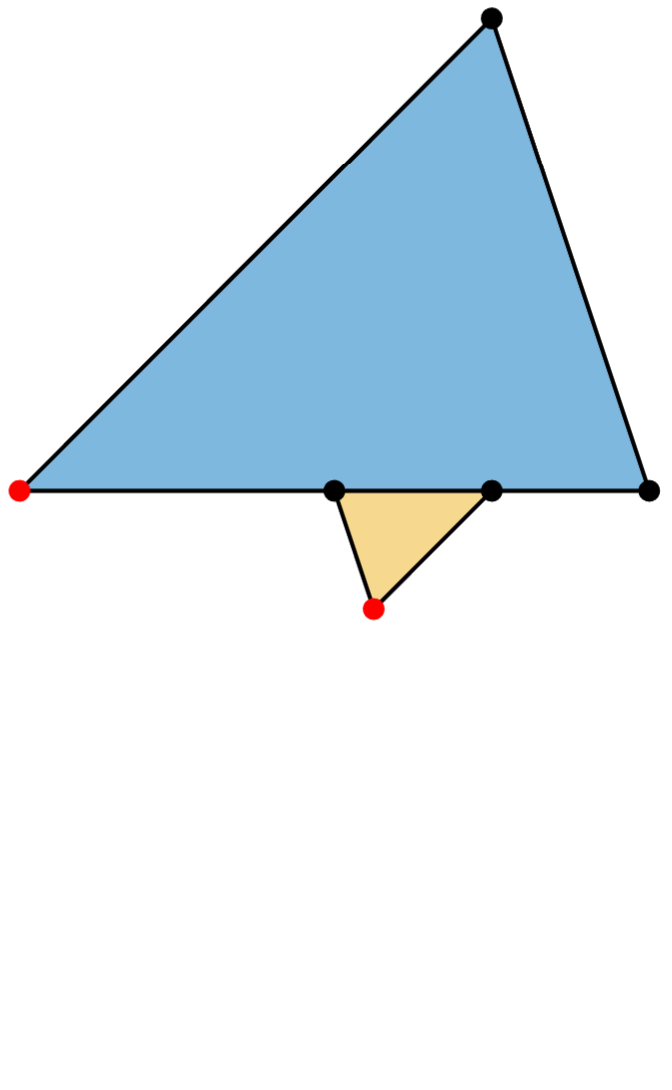}
\put(18,44){$T''$}
\put(30,70){$T'$}
\put(36,38){$v(T'')$}
\put(-17,62){$v(T')$}
\put(20,-5){\texttt{case D}}
\end{overpic}
\end{center}
\caption{ {The elements $T'$ and $T''$ are adjacent in cases A to D. They are compatible in cases A and B,
    and non-compatible in cases C and D.}}
\label{fig:compatible}
\end{figure}

Two elements $T', T'' \in \mesh$ are said {\em adjacent} if $E=T'\cap T''$ is an edge for at least one element, and are said {\em compatible} if they are adjacent and both $E(T')$ and $E(T'')$ belong to the same line (see Fig. \ref{fig:compatible}, cases A and B).

The following technical results will be helpful in the design of the refinement procedure.

\begin{lemma}[global index of a hanging node]\label{L:global-index}
Consider an edge $E=[{x}',{x}'']$ of the partition $\mesh$. If ${x} \in {\cal H} \cap {\rm int\, }E$ is generated by $m \geq 1$ bisections of $E$, then its global index $\lambda({x})$ satisfies
$$
\lambda({x})= \max (\lambda({x}'), \lambda({x}'')) +m \,.
$$
\end{lemma}
\begin{proof} If $m=1$, then $x=x_M$ is the midpoint of $E$, and the formula is just the Definition \ref{d-globalindex} of global index. If $m >1$, then ${x}$ is generated by bisecting some interval $[{z}',{z}''] \subset E$, and $\lambda({x})=\max (\lambda({z}'), \lambda({z}''))+1$. Exactly one between ${z}',{z}''$ has been generated by $m-1$ bisections, whereas the other one has been
generated by less than $m-1$ bisections. Hence, one concludes by induction.  \end{proof}

\begin{lemma}[reducing the global index of hanging nodes]\label{L:reduce-globindex}
Let ${\cal H} \cap {\rm int\, }E$ contain at least the midpoint ${x}_M$ of $E$. Assume that a bisection of some element in $\mesh$ transforms ${x}_M$ into a proper node, and let $\lambda_{\rm new}$ denote the new global-index mapping of the nodes in ${\cal H} \cap {\rm int\, }E$ after the bisection. Then there holds
$$
\lambda_{\rm new}({x}) \leq \lambda({x}) -1 \qquad \forall {x} \in {\cal H} \cap {\rm int\, }E \,.
$$
\end{lemma}
\begin{proof}
If ${x}={x}_M$, then trivially $\lambda_{\rm new}({x}) = 0 \leq \lambda({x}) -1$. If ${x} \in {\cal H} \cap {\rm int\, }E$ is contained, say, in $({x}',{x}_M)$ and has been generated by $m>1$ successive bisections of $E$, then it is generated by $m-1$ successive bisections of $[{x}',{x}_M]$. Thus, by applying Lemma \ref{L:global-index} we get
\begin{eqnarray*}
\lambda_{\rm new}({x}) &\leq& \max(\lambda_{\rm new}({x}'), \lambda_{\rm new}({x}_M) ) + m-1 \\
&=& \max(\lambda({x}'), 0 ) + m-1 \ = \ \lambda({x}') + m-1 \\
&\leq& \max((\lambda({x}'), \lambda({x}''))+m-1 \ = \ \lambda({x})-1\,.  \qquad \qquad 
\end{eqnarray*} 
This gives the desired estimate.
\end{proof}

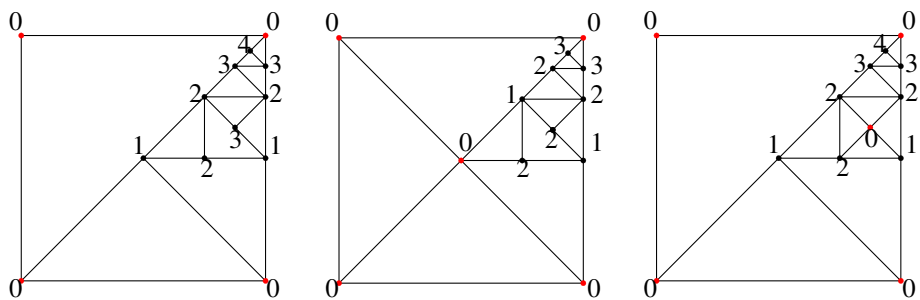
\begin{figure}[t!]
 {\small
\begin{center}
\begin{overpic}[scale=0.20]{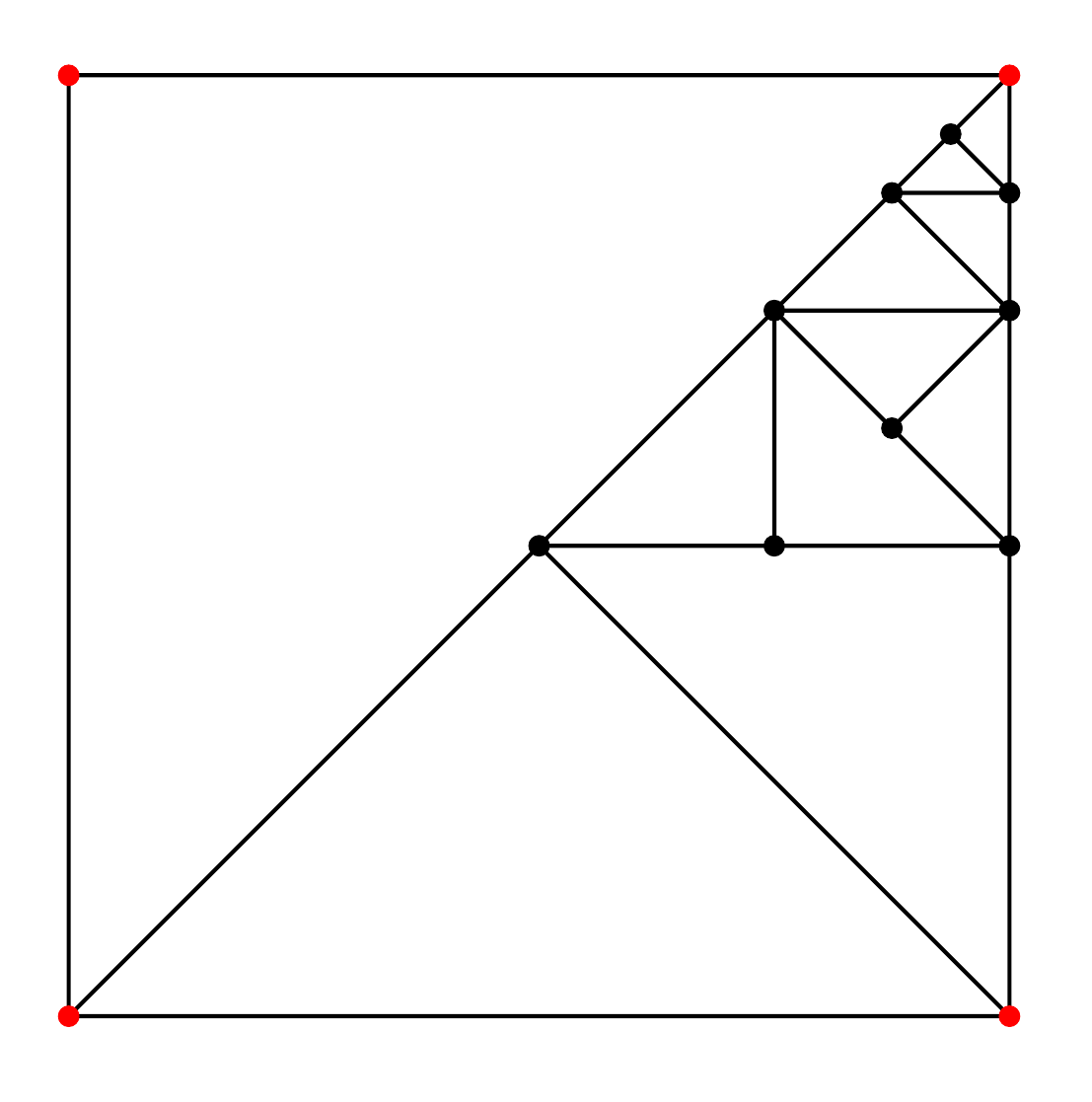}
\put( 2, 2){$0$}
\put( 2,95){$0$}
\put(92, 2){$0$}
\put(92,95){$0$}
\put(45,52){$1$}
\put(93,52){$1$}
\put(65,70){$2$}
\put(69,44){$2$}
\put(93,70){$2$}
\put(75,80){$3$}
\put(79,54){$3$}
\put(93,80){$3$}
\put(82,87){$4$}
\end{overpic}
\quad
\begin{overpic}[scale=0.20]{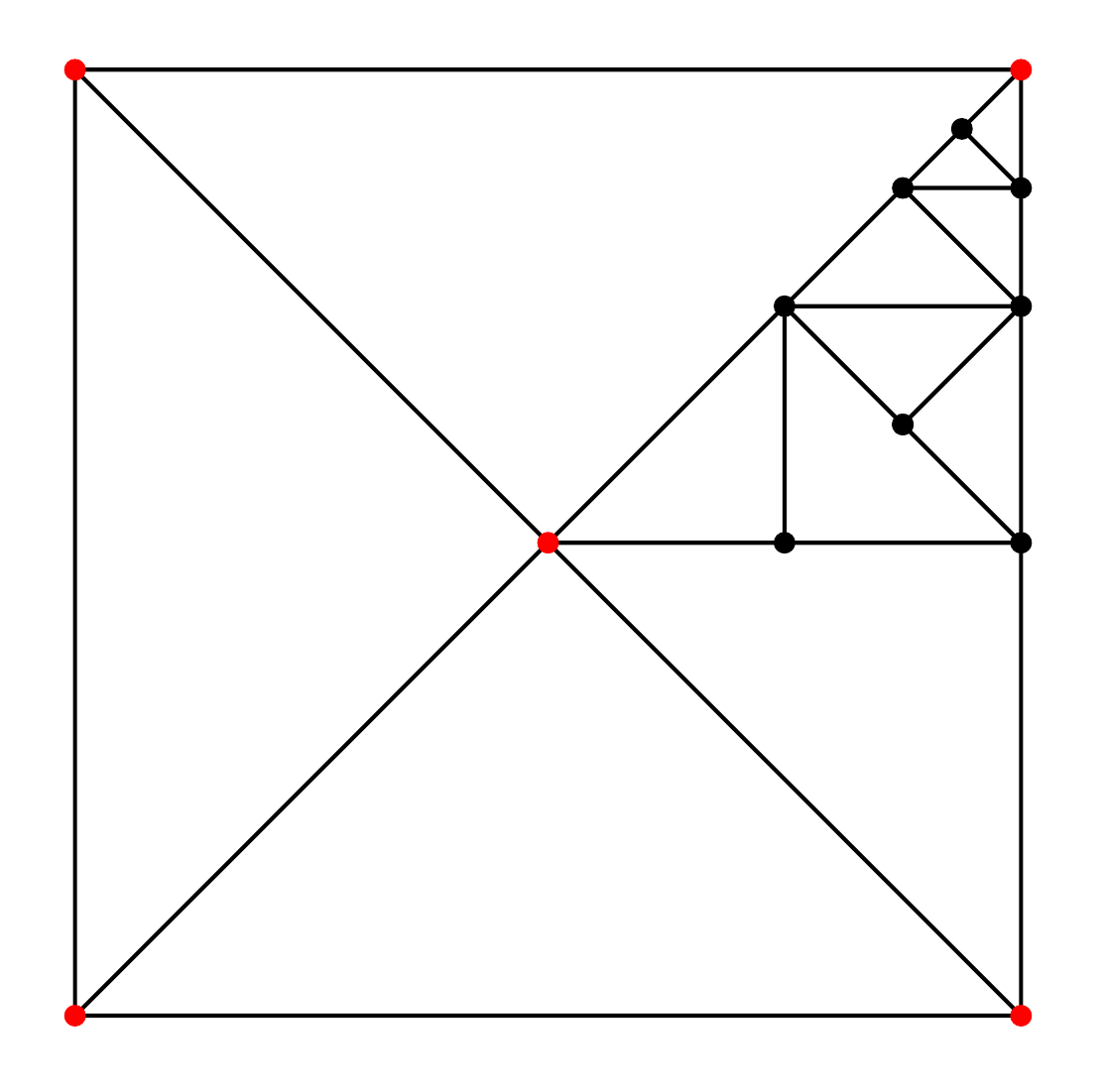}
\put( 2, 2){$0$}
\put( 2,95){$0$}
\put(94, 2){$0$}
\put(94,95){$0$}
\put(49,53){$0$}
\put(95,53){$1$}
\put(65,70){$1$}
\put(69,44){$2$}
\put(95,70){$2$}
\put(75,80){$2$}
\put(79,54){$2$}
\put(95,80){$3$}
\put(82,87){$3$}
\end{overpic}
\quad
\begin{overpic}[scale=0.20]{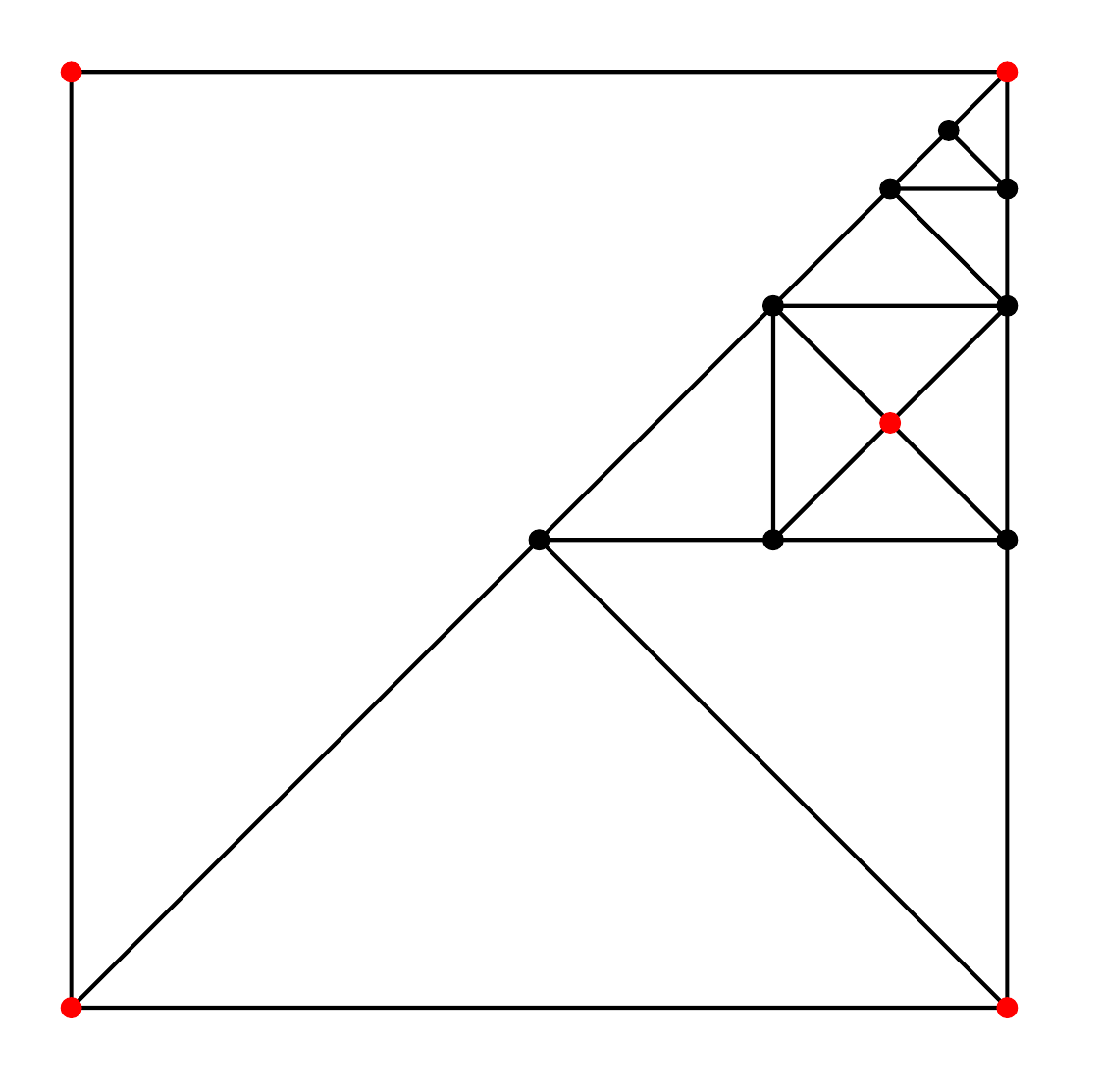}
\put( 2, 2){$0$}
\put( 2,95){$0$}
\put(92, 2){$0$}
\put(92,95){$0$}
\put(45,52){$1$}
\put(93,52){$1$}
\put(65,70){$2$}
\put(69,44){$2$}
\put(93,70){$2$}
\put(75,80){$3$}
\put(79,54){$0$}
\put(93,80){$3$}
\put(82,87){$4$}
\end{overpic}
\end{center}
\caption{Three examples of distributions of proper nodes (red) and hanging nodes (black), with associated global indices $\lambda$. The bisection added in the middle picture converts the centered node into proper, and induces nonlocal changes of global indices on chains associated with it; if $\Lambda=3$, the leftmost mesh is not admissible and this procedure is instrumental to restore admissibility. The right picture illustrates the creation of a proper node without nonlocal effects on global indices.
}
\label{fig:sample_element}
}
\end{figure}

The result just established is the motivation for the proposed refinement strategy, introduced in \cite{da2023adaptive}. Indeed, it assures that in order to reduce the global index of a hanging node sitting on an edge, it is enough to transform the midpoint of the edge into a proper node. The situation is well represented in Figure \ref{fig:sample_element}.

The following remark will be useful in the sequel.
\begin{remark}[facing element]\label{rem:facing}
{\rm
 Given a $\Lambda$-admissible mesh $\mesh$ and $T \in \mesh$, let $x(T)$ be the midpoint of $E(T)$, and suppose that $\lambda(x(T)) > \Lambda$. Then $x(T)$ is not a node of $\mesh$,  whence the edge $E(T)$ cannot contain any hanging node in its interior.  We conclude that there exists a unique adjacent element $\widetilde{T} \in \mesh$, $\widetilde{T} \not = T$, such that  $T \cap \widetilde{T} = E(T)$. This element will be called the element {\em facing} $T$, and denoted by $F(T)$.
}
\end{remark}

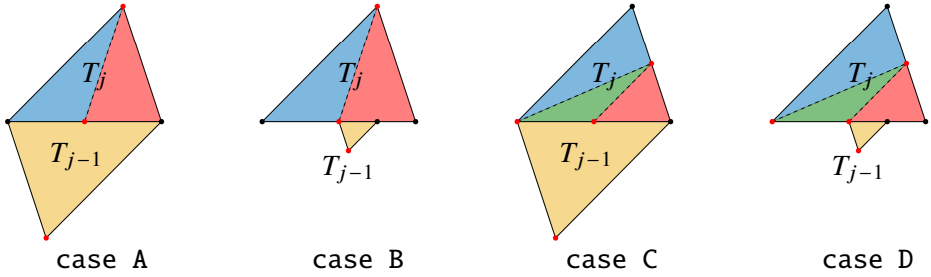
\begin{figure}[t!]
\begin{center}
\begin{overpic}[scale=0.19]{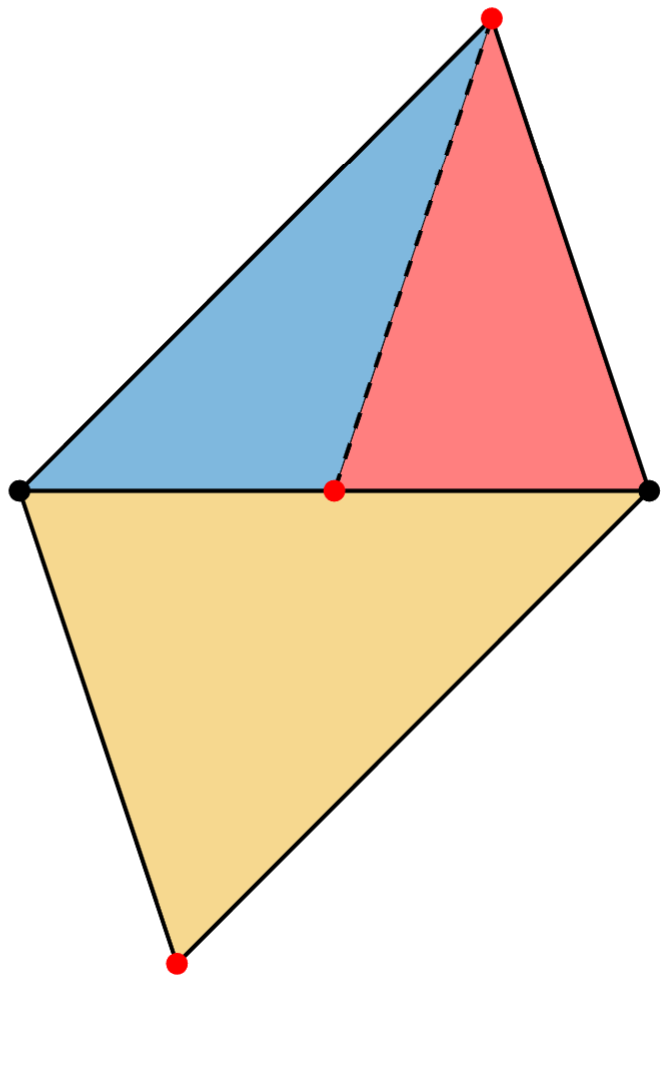}
\put(18,40){$T_{j-1}$}
\put(30,70){$T_j$}
%\put(20,10){$\bm{nv}(E_{k-1})$}
%\put(50,95){$\bm{nv}(E_k)$}
\put(20,0){\texttt{case A}}
\end{overpic}
\qquad \quad
\begin{overpic}[scale=0.19]{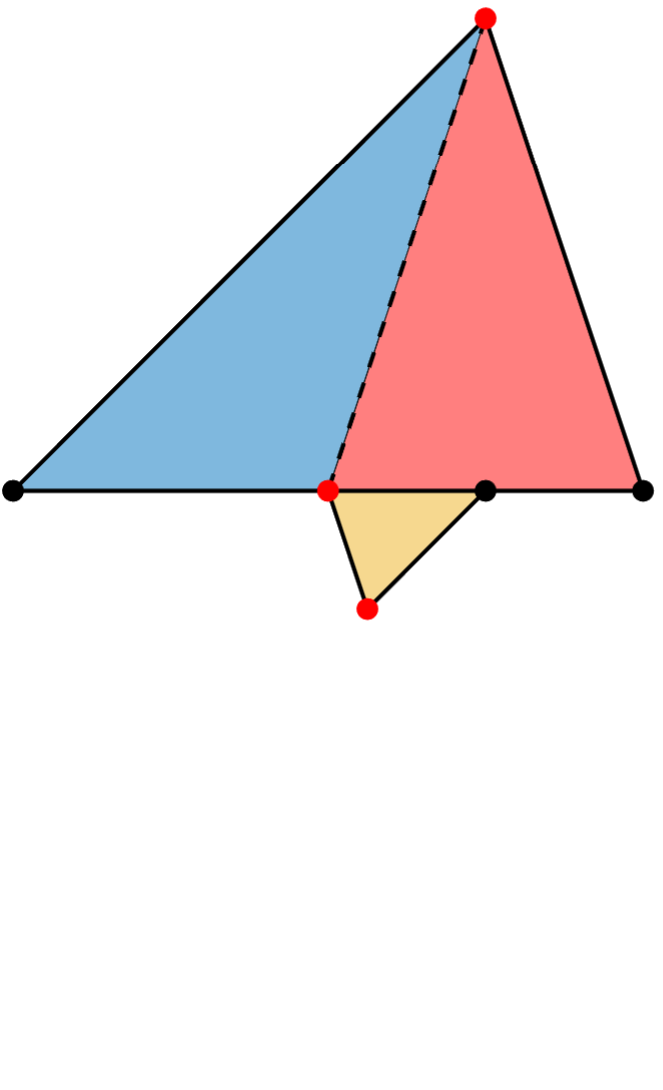}
\put(24,35){$T_{j-1}$}
\put(30,70){$T_j$}
%\put(38,42){$\bm{nv}(E_{k-1})$}
%\put(50,95){$\bm{nv}(E_k)$}
\put(20,0){\texttt{case B}}
\end{overpic}
\qquad \quad
\begin{overpic}[scale=0.19]{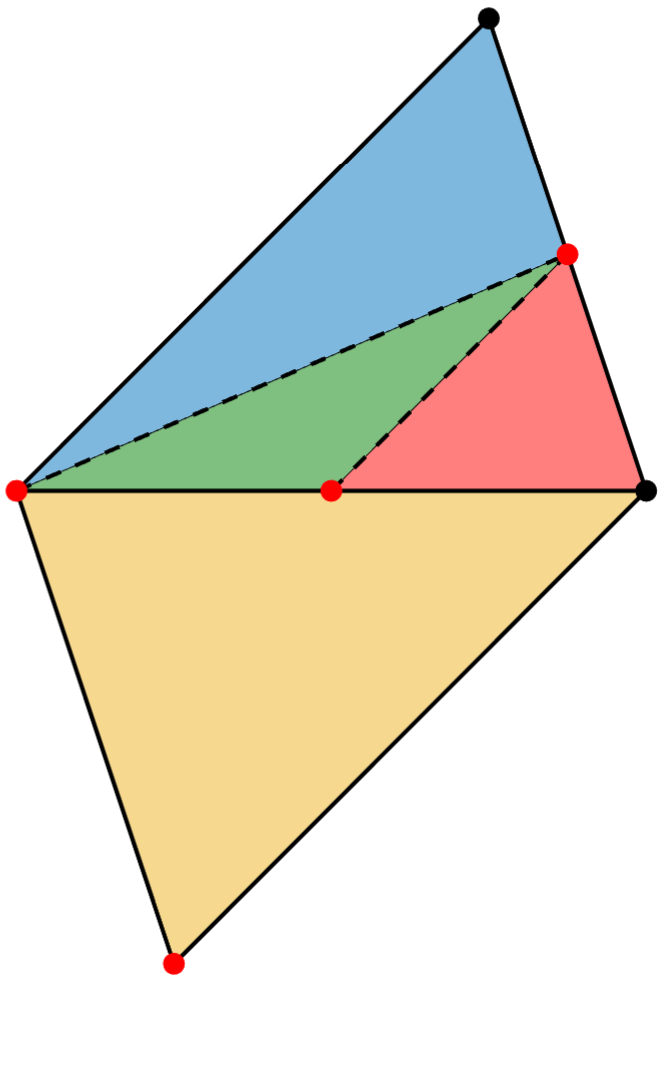}
\put(18,40){$T_{j-1}$}
\put(30,70){$T_j$}
%\put(20,10){$\bm{nv}(E_{k-1})$}
%\put(-18,62){$\bm{nv}(E_k)$}
\put(20,0){\texttt{case C}}
\end{overpic}
\qquad \quad
\begin{overpic}[scale=0.19]{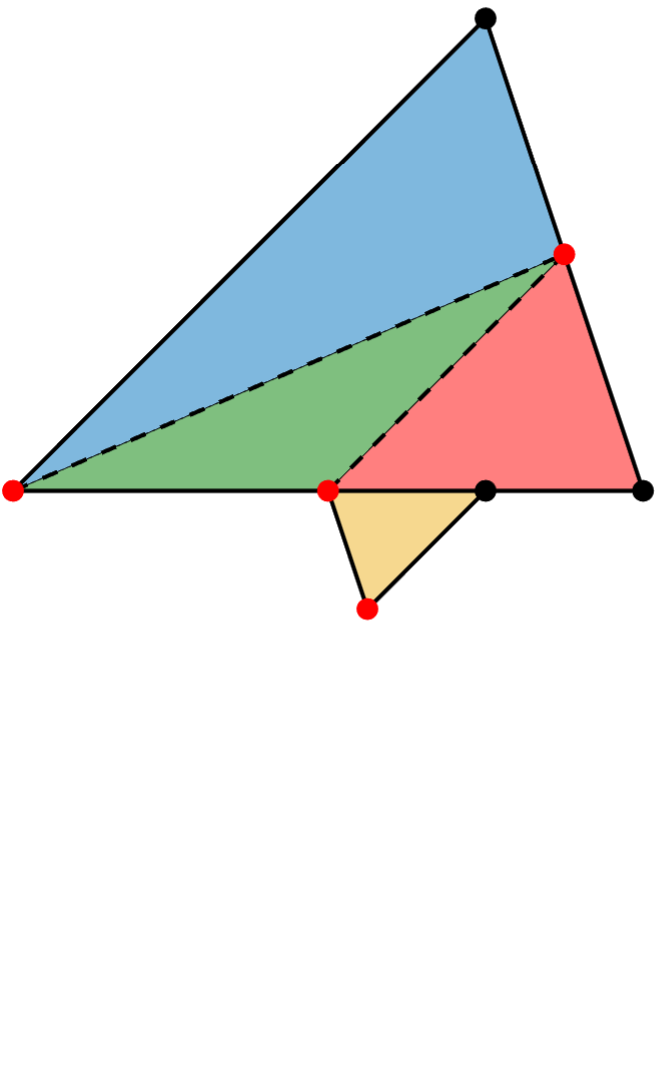}
\put(24,35){$T_{j-1}$}
\put(30,70){$T_j$}
%\put(38,42){$\bm{nv}(E_{k-1})$}
%\put(-18,62){$\bm{nv}(E_k)$}
\put(20,0){\texttt{case D}}
\end{overpic}
\end{center}
\caption{Two elements $T_{j-1}$ and $T_j$ in the chain $\RC(\mesh,T)$: $T_{j-1}$ can be bisected in a $\Lambda$-admissible way, only after $T_j$ is refined once (cases A and B), or twice (cases C and D)}
\label{fig:catena}
\end{figure}

Given an element $T \in \mesh$ which has been marked for refinement, we are ready to identify those elements  {in $\mesh$} that need be bisected with $T$ in order to create a $\Lambda$-admissible refinement of $\mesh$. Figure~\ref{fig:catena} illustrates the possible situations.

\begin{definition}[chain of elements to be refined]\label{D:ref-chain}
Define by recurrence the chain of elements starting at $T$
$$
\RC(T,\mesh)=\{T_0, T_1, \dots, T_k\}
$$
for some $k \geq 0$, as follows:  set first $T_0=T$ and, assuming to have defined $T_j$ for $j \geq 0$, then
\begin{enumerate}[label=(\roman*)]
\item if $\lambda(x(T_j)) \leq \Lambda$, set $k=j$ and stop;
\item if $\lambda(x(T_j)) = \Lambda+1$ and the facing element $F(T_j)$ is compatible with $T_j$, set $T_{j+1}=F(T_j)$, $k=j+1$ and stop;
\item if $\lambda(x(T_j)) = \Lambda+1$ and the facing element $F(T_j)$ is not compatible with $T_j$, set $T_{j+1}=F(T_j)$ and continue.
\end{enumerate}
\end{definition}

\begin{lemma}[properties of the chain of refinement]\label{L:chain}
 {The chain $\RC(T,\mesh)$ has finite length, precisely it holds $k \leq g(T)+1$, where $g(T)$ is the generation of $T$, defined in Sect. \ref{S:bisection}.
Furthermore, the sequence of element generations $\{ g(T_j) \}_{j = 0}^k$} is not increasing.
\end{lemma}
\begin{proof}
 {We claim that step (iii) in Definition \ref{D:ref-chain} reduces the generation by at least one. In fact, $T_j$ coincides with or is a refinement of a triangle $\hat{T} \in \mathbb{T}$ sharing a full edge with $T_{j+1}$; thus $g(T_j)\ge g(\hat{T})$. Such triangle $\hat{T}$ satisfies $g(\hat{T}) = g(T_{j+1})+1$, whence 
\begin{equation}\label{eq:bound-levels}
g(T_{j+1}) = g(\hat{T}) -1 \leq g(T_j)-1.
\end{equation}
Therefore, for as long as case (iii) is active, i.e. for all $j<k$, we have $g(T_j) \leq g(T_0) - j$ and
$$
0 \leq g(T_{k-1}) \leq g(T_0)-(k-1) \, ,
$$
which gives the first statement of the lemma. The monotonicity of $\{ g(T_j) \}_{j = 0}^k$
follows from \eqref{eq:bound-levels} and the fact that $g(T_{k-1})=g(T_k)$ in case (ii).
}
\end{proof}

Once the chain $\RC(\mesh,T)$ is defined, all its elements are refined, starting from the last one and proceeding backwards. This is accomplished in the following procedure.

\begin{algotab}
  \> $[\grid_*]=\textsf{REFINE\_RECURSIVE}~(\grid,T,\Lambda)$\\
  \> \> if $\lambda(x(T))\le \Lambda$ \\
    \> \> \> bisect $T$ \\
    \> \> \> update $\grid$ \\
    \> \> else if $F(T)$ is compatible with $T$ \\
    \> \> \>  bisect $F(T)$ and $T$ \\
    \> \> \>  update $\grid$ \\
    \> \>  else \\
    \> \> \>  $[\grid] =\textsf{REFINE\_RECURSIVE}~(\grid,F(T),\Lambda$)\\
   \> \> return $\grid$
\end{algotab}

\begin{proposition}[properties of $\textsf{REFINE\_RECURSIVE}$]\label{prop:level-bound}
   If $\mesh$ is $\Lambda$-admissible, %then
   the call $[\mesh_*] =\textsf{REFINE\_RECURSIVE}~(\grid,T,\Lambda)$ outputs the smallest $\Lambda$-admissible refinement $\gridk[*]$
of $\grid$ such that $T$ is bisected. In addition,  every element $T'\in\mesh_*$ generated by this call satisfies
\begin{equation}\label{eq:levels-chain}
g(T') \leq g(T)+1 \,.
\end{equation}
\end{proposition}
\begin{proof}
Let $\RC(T,\mesh) = \{T_j\}_{j=0}^k$ and
observe that, for $j \geq 1$, one or two bisections of $T_j $ convert the midpoint of the edge $E$ of $T_j$
shared with $T_{j-1}$ into a proper node. Therefore, Lemma \ref{L:reduce-globindex} (reducing the
global index of hanging nodes) implies that
the global indices of all interior nodes to $E$ decrease by at least $1$, and makes the
bisection of $T_{j-1}$ $\Lambda$-admissible as desired.
  
To prove \eqref{eq:levels-chain} we take $j \geq 1$ and consider the following two mutually
exclusive cases.
If $T_j$ and $T_{j-1}$ are compatible, then $T_j$ is replaced by two elements $T'\in\mesh_*$ of generation
\[
g(T') = g(T_j)+1 \leq  g(T)+1,
\]
according to Lemma \ref{L:chain} (properties of the chain of refinement). On the other hand, if $T_j$ and $T_{j-1}$ are not compatible, then $T_j$ is replaced by one element of generation $g(T_j)+1$ and two elements $T'\in\mesh_*$ of generation
\[
g(T') = g(T_j)+2 \leq g(T_{j-1})+1 \leq g(T)+1
\]
because of \eqref{eq:bound-levels}. Finally, the element $T_0=T$ is replaced by two elements of generation $g(T)+1$.
\end{proof}

If one considers the chains starting at any element $T \in \marked$, one obtains the     procedure
$\textsf{REFINE}~(\grid,\marked,\Lambda)$, which reads
\begin{algotab}
\index{Algorithms!\REFINE: refine marked elements and others necessary to produce a $\Lambda$-admissible mesh}
  \> $[\grid_*]=\textsf{REFINE}~(\grid,\marked, \Lambda)$\\
  \> \>  for all $\elm\in\marked\cap\grid$ do\\
  \> \> \> $[\grid]=\textsf{REFINE\_RECURSIVE}~(\grid,\elm,\Lambda)$\\
  \> \> return~$\grid$
\end{algotab}
and outputs a minimal $\Lambda$-admissible mesh 
$\grid_*\ge\grid$, refinement of $\grid$, so that all marked elements of $\marked$ are refined.

\medskip
{\it Proof of Theorem \ref{T:nonconforming-meshes} (complexity of $\REFINE$ for $\Lambda$-admissible meshes).} 
The arguments given in Sect. \ref{nv-SS:proof-complexity-conforming} for the conforming case can be adapted to the current situation. The two crucial properties needed are the relation \eqref{nv-dist-gen} between the distance of two elements in a chain and their generation, which is valid for bisection grids regardless of $\Lambda$-admissibility, and the relation \eqref{eq:levels-chain} between generations of elements. $\square$

%-----------------------------------------------------------------------------
\subsubsection{Mesh overlay and $\Lambda$-admissibility}\label{sec:admissible-overlay}
%-----------------------------------------------------------------------------
%
Given two partitions $\mesh_A$ and $\mesh_B$,   denote by $\mesh_A \oplus \mesh_B$  the {\em overlay} of $\mesh_A$ and $\mesh_B$, i.e., the partition whose associated tree is the union of the trees of $\mesh_A$ and $\mesh_B$.  The following property holds.
 \begin{proposition}[mesh overlay is $\Lambda$-admissible]\label{P:mesh-overlay-Lambda}
 If $\mesh_A$ and $\mesh_B$ are $\Lambda$-admissible, then $\mesh_A \oplus \mesh_B$ remains $\Lambda$-admissible.
 \end{proposition}
 \begin{proof}
 Denote here by ${\cal N}$ the set of all nodes obtained by newest-vertex bisection from the root partition $\mesh_0$. Let ${\cal N}_0$, ${\cal N}_A$, ${\cal N}_B$, ${\cal N}_{A + B}$, resp., be the set of nodes of the partitions  $\mesh_0$, $\mesh_A$, $\mesh_B$,  $\mesh_A \oplus \mesh_B$, resp.. It is easily seen that for each ${x} \in {\cal N}\setminus {\cal N}_0$ there exists a unique set ${\cal B}({x}) = \{{x}', {x}''\} \subset {\cal N}$ such that ${x}$ is generated by the bisection of the segment $[{x}', {x}'']$. Furthermore, if ${x} \in {\cal N}_{A+B}$ is a proper node of $\mesh_A$ (of $\mesh_B$, resp.), then it is also a proper node of $\mesh_A \oplus \mesh_B$.
 
 Let us denote by $\lambda_A$, $\lambda_B$, $\lambda_{A+B}$, resp., the global-index mappings defined on ${\cal N}_A$, ${\cal N}_B$, ${\cal N}_{A + B}$, resp.. It is convenient to extend the definition of $\lambda_A$ and $\lambda_B$ to the whole ${\cal N}_{A + B}$ by setting 
$$
\lambda_A({x}) = +\infty \quad \text{if } {x} \in {\cal N}_{A + B}\setminus {\cal N}_A \,, \qquad 
\lambda_B({x}) = +\infty \quad \text{if } {x} \in {\cal N}_{A + B}\setminus {\cal N}_B \,. 
$$
With these notations at hand, we are going to prove the inequality
\begin{equation}\label{eq:overlay-1}
\lambda_{A+B}({x}) \leq \min ( \lambda_A({x}),  \lambda_B({x}) ) \qquad \forall {x} \in {\cal N}_{A + B}\,,
\end{equation} 
from which the thesis immediately follows.

We proceed by induction on $k=\lambda_{A+B}(x)$, ${x} \in {\cal N}_{A + B}$. If $k=0$, the inequality is trivial since $\lambda_A({x}),  \lambda_B({x}) \geq 0$. So suppose \eqref{eq:overlay-1} hold up to some $k\geq 0$. If ${x} \in {\cal N}_{A + B}$ satisfies $\lambda_{A+B}({x})= k+1 >0$, then it is a hanging node of $\mesh_A \oplus \mesh_B$ by definition of global index, hence, it is a hanging node of $\mesh_A$ or $\mesh_B$; without loss of generality, suppose it is a hanging node of $\mesh_A$.
If ${x}$ is generated by the bisection of the segment $[{x}', {x}'']$, then again by definition of global index it holds
$$
k+1 = \lambda_{A+B}({x}) = \max( \lambda_{A+B}({x'}), \lambda_{A+B}({x}'') ) + 1 \,,
$$ 
which implies 
$$
\lambda_{A+B}({x'}) \leq k \,, \qquad  \lambda_{A+B}({x''}) \leq k \,.
$$
By induction,
$$
\lambda_{A+B}({x'}) \leq \min ( \lambda_A({x'}),  \lambda_B({x'}) ) \,, \qquad \lambda_{A+B}({x''}) \leq \min ( \lambda_A({x''}),  \lambda_B({x''}) ) \,, 
$$
from which we obtain
$$
\lambda_{A+B}({x}) \leq \max(\lambda_A({x'}), \lambda_A({x''}) ) +1 = \lambda_A({x})
$$ 
since ${x}$ is a hanging node of $\mesh_A$. On the other hand, either ${x} \in {\cal N}_B$ or ${x} \not \in {\cal N}_B$. In the latter case, $\lambda_B(x) = +\infty$, and \eqref{eq:overlay-1} is proven. In the former case, necessarily ${x}$ is a hanging node of $\mesh_B$, hence as above
$$
\lambda_{A+B}({x}) \leq \max(\lambda_B({x'}), \lambda_B({x''}) ) +1 = \lambda_B({x}) \,,
$$
and the thesis  is proven.    
\end{proof}

%--------------------------------------------------------------------------------

\section{Discontinuous Galerkin Methods }\label{S:dg}
% \rhn{(AB $\longrightarrow$ RHN)}

So far we have studied conforming finite element approximations. In this section we present and 
analyze a two-step $\AFEM$ for discontinuous Galerkin methods (dG). The core PDE routine $\GALERKIN$
is thereby replaced by $\GALERKINIP$, which hinges on the interior penalty discontinuous FEM.
We regard dG as a prototype non-conforming method of practical importance and thus the natural 
first step to investigate the effects of non-conformity within adaptivity.

Finite element functions, being discontinuous, allow for non-conforming meshes to support them.
We consider $\Lambda-$admissible subdivisions, according to Definition \ref{d-LambdaAdmissibility}, where $\Lambda \geq 0$ restricts the level of non-conformity, and denote by $\grids^\Lambda$ the collection of all $\Lambda-$admissible refinements of an initial subdivision $\grid_0$; we refer to Section~\ref{S:mesh-refinement} for details. However,
we further assume that $\grid_0$ is conforming to limit the level of technicalities.

There are several novel but characteristic aspects of dG. The most notable one is the appearance
of jumps in its formulation, to compensate for the lack of $H^1$-conformity, as well as in 
the a posteriori upper bounds and the comparison of Galerkin solutions on different meshes. The
lack of monotonicity of these jumps presents a formidable obstruction to the available proof teachniques in adaptivity.
However, we show in Lemma~\ref{l:jump} that they are controlled by the residual estimator,
thereby enabling us to loosely follow the roadmap of the conforming method, namely Sections~\ref{S:aposteriori},~\ref{S:convergence-coercive}, and \ref{S:conv-rates-coercive}.
Our approach is based on \cite{bonito2010quasi} for the one-step AFEM. 

The extra flexibility provided by non-conforming meshes, and corresponding discontinuous functions,
does not yield better asymptotic rate in $H^1$. An early manifestation of this fact, although written for conforming subdivisions, is Proposition~\ref{P:equivalence-classes} (equivalence of classes for $u$). We extend this result below for general $\Lambda-$admissible partitions.

One advantage of the two-step \AFEM is that its design and analysis allows for $f\in H^{-1}(\Omega)$ without added difficulties: the function $f$ is replaced by the discrete functional $\wh f = P_\grid f \in \mathbb F_\grid$, which applies to functions in $\mathbb S_\grid^{n,-1}$. This is in contrast with $f$, which cannot be applied to functions in $\mathbb S_\grid^{n,-1}$. We exploit this property and thereby extend the applicability of dG to load functions in $H^{-1}(\Omega)$. \looseness=-1

Our intend is to analyze the following algorithm for the approximation of the solution $u\in H^1_0(\Omega)$ to the coercive problem \eqref{weak-form}.
\begin{algo}[\AFEMTSIP]\label{algo:AFEM-TS-IP}
\index{Algorithms!\AFEMTSIP: interior penalty version of \AFEMTS}
Given an initial tolerance $\eps_0>0$, a target tolerance $\tol$ and initial mesh $\mesh_0$,
as well as a safety parameter $\omega \in (0,1]$, \AFEMTSIP is a two-step algorithm alternating between the resolution of data $\data$ and the Galerkin solution $u_\grid$:
\begin{algotab}
  \> $[\mesh, u_\mesh]=\AFEMTSIP (\mesh_0, \varepsilon_0, \omega, \tol)$\\
  \> \>  $\text{set }k=0$\\
  \> \> $\text{do}$  \\
  \> \> \> $[\widehat\mesh_{k},\widehat{\data}_{k}]=\DATA \, (\mesh_k, \data, \omega \, \varepsilon_k) $ \\
  \> \> \> $[\mesh_{k+1},u_{k+1}]=\GALERKINIP \, (\widehat{\mesh}_{k},\widehat{\data}_{k},\varepsilon_k)$ \\
  \> \> \> $\varepsilon_{k+1}=\tfrac12 {\varepsilon_k}$ \\
  \> \> \> $k \leftarrow k+1$ \\
  \> \> \text{while $\eps_{k-1} > \tol$} \\
  \> \> return $\grid_{k},u_{k}$
\end{algotab}
\end{algo}
In \AFEMTSIP, the module $\DATA \, (\grid,\data,\tau)$ is the same as described in Section~\ref{S:DATA} except that it produces approximate data $\wh \data \in \mathbb D_{\wh 
 \grid}$, defined in \eqref{E:space-discrete-data}, subordinate with a $\Lambda-$admissible refinement $\wh \mesh$ of $\grid_0$ for $\Lambda \geq 0$ rather than $\Lambda=0$ (conforming).  The discrete data $\wh \data$ also satisfies the structural assumption \eqref{E:structural-assumption-wh} as discussed in Section~\ref{S:data-approx}. It is worth pointing out that the projection $P_\grid$ used to approximate the right-hand side $f \in H^{-1}(\Omega)$ as well as all the results and algorithms presented in Section~\ref{S:data} are restricted to conforming subdivisions $\grids$. We briefly discuss in Section~\ref{ss:P_grid_nc} the extension of $P_\grid$ and $\DATA$ to $\Lambda$-admissible subdivisions.
 Algorithm~\ref{A:GALERKIN-IP} describes the module \GALERKINIP, the counterpart of \GALERKIN 
 for dG formulations.

In Section~\ref{ss:dG} we introduce notations and tools relevant for the characterization of discontinuous finite elements. Among them is the operator $\IdG$ that projects piecewise polynomial functions onto globally continuous piecewise polynomial functions. It is instrumental to derive a Poincar\'e inequality on the discontinuous spaces and guarantee that the approximation classes $\mathbb A_s^{-1}$ for the solution $u$ using discontinuous approximation on $\Lambda$-admissible subdivisions are equivalent to their conforming counter-parts $\mathbb A_s^0$ introduced in Section~\ref{S:conv-rates-coercive}. We present the discontinuous Galerkin method in Section~\ref{ss:IP}. We start with the standard symmetric interior penalty, discuss its drawbacks regarding the unnecessary regularity beyond $H^1_0(\Omega)$ imposed on the exact solution $u$, and describe a reformulation valid in $H^1_0(\Omega)$. The latter suffers from lack of consistency that needs to be accounted for. The a posteriori estimates for the perturbed problem \eqref{E:perturbed-weak-form} are derived in Section~\ref{ss:posteriori-dG}. Because the data is polynomial within \GALERKINIP, the a posteriori estimators are oscillation free.
The \GALERKINIP module is analyzed in Section~\ref{ss:GalerkinIP} while the discussion of rate-optimality of \AFEMTSIP is reserved for Section~\ref{ss:AFEMTSIP}.

%---------------------------------------------------------------------------------
\subsection{Discontinuous Galerkin setting}\label{ss:dG}
%---------------------------------------------------------------------------------

We start with an initial conforming subdivision $\grid_0$ made of simplices or hexahedra satisfying Assumption~\ref{A:initial-labeling} (initial labeling). Given $\Lambda>0$, the refinement procedure 
$\REFINE$ is designed to produce a $\Lambda$-admissible sequence of meshes $\grids^\Lambda$ obeying Theorem~\ref{T:nonconforming-meshes} ($\REFINE$ for $\Lambda$-admissible meshes). From now on, 
we do not specify the dependency on $\Lambda$ in the constants.

\subsubsection{Basic setting}
\index{Functional Spaces!$\mathbb V_\grid^{-1}$: non-conforming finite element space}
For $\grid \in \grids^\Lambda$, we denote by $\mathbb V_\grid^{-1}:=\mathbb S_{\grid}^{n,-1}:= \prod_{T \in \grid} \mathbb P_n(T)$ the space of piecewise polynomials of degree at most $n \geq 1$ subordinate to a partition $\grid$. In contrast with the conforming spaces
\[
\mathbb V_\grid^0 := \mathbb S_{\grid}^{n,0} \cap H^1_0(\Omega)
\]
considered earlier, the space $\mathbb V_\grid^{-1}$ consists of (possibly) discontinuous functions across the elements $T\in \grid$ and do not necessarily satisfy the vanishing boundary condition. Continuity across elements and vanishing boundary condition will be weakly imposed in the discontinuous Galerkin formulations.

We recall from Section \ref{S:nonconforming-meshes} that for a proper (interior) node $P \in {\cal P}$, the {\it domain of influence} $\omega_\grid(P)=\textrm{supp}(\psi_P)$ is the support of the Lagrange basis function $\psi_P\in\mathbb V_\grid^0$ associated with the node $P$; we refer to Fig.~\ref{f-basis-function}. Since the sequence of meshes is $\Lambda$-admissible, Proposition \ref{p-size-domain-influence}
(size of the domain of influence) shows that the number of elements $T\in \grid$ such that $T\subset \omega_\grid(P)$ is uniformly bounded for $\grid \in \grids^\Lambda$.

The set of faces associated with a subdivision $\grid \in \grids^\Lambda$ is denoted $\faceswbdy:=\faceswbdy(\grid)$, and it contains boundary faces as well as interior faces. The set of interior faces is denoted $\faces$.
For a face $F \in \faceswbdy$, we denote by $\meaninline{v}|_F$ and $\jump{.}|_F$ the average and jump operators across a face $F$.
To define them precisely, we associate for each face $F \in \faceswbdy$ one of the two unit normals $\vec{n}_F$. The choice of $\bn_F$ is fixed but irrelevant as long as the outward pointing normal to $\Omega$ is chosen for boundary faces. Let $T_{\pm}\in\grid$ be the elements that share the interior face $F$, namely
$F = T_-\cap T_+$, and $\mp \vec{n}_F$ be their outward pointing normals. Now, given $v \in \mathbb V_{\grid}^{-1}$, let $v_\pm := v|_{T_\pm}$ and define for an interior face $F$
\begin{equation}\label{e:average_jump}
\index{Meshes!$\meaninline{.}$: average on faces}
\index{Meshes!$\jump{.}$: jump across faces}
\meaninline{v}|_F := \frac 1 2 \big(v_- + v_+\big)|_F, \qquad \jump{v}|_F := \big(v_- - v_+\big)|_F, 
\end{equation}
By convention, we set $\meaninline{v}|_F := v_-$ and $\jump{v}|_F := v_-$ whenever $F$ is a boundary face. These definitions extend readily for vector valued functions. 

We use the subscript $\grid$ to denote the piecewise version of differential operators. For instance,
the {\it broken gradient} $\nabla_\grid$ is the piecewise gradient $\nabla_\grid v|_{\mathring{T}}= \nabla v|_{\mathring{T}}$ for $T\in \grid$ and $v \in \mathbb V_\grid^{-1}$.
For simplicity, we write 
$$
\| v \|_{L^2(\tau)}^2:= \sum_{T \in \tau} \| v \|_{L^2(T)}^2
$$
for any subset $\tau \subset \grid$ of elements and 
$$
\| v \|_{L^2(\sigma)}^2:= \sum_{F \in \sigma} \| v \|_{L^2(F)}^2
$$
for any subset $\sigma \subset \faceswbdy$ of faces.
We also define a meshsize function $h:=h_{\grid}:\overline{\Omega} \rightarrow (0,\infty)$ such that $h|_{\mathring{T}} \approx \textrm{diam}(T)$ for $T \in \grid$ and $h|_F \approx \textrm{diam}(F)$ for $F \in \faceswbdy$. 
With such notations at hand, the broken $H^1$ space
$$
\index{Functional Spaces!$\mathbb E_\grid$: broken $H^1$ space}
\mathbb E_\grid:= H^1(\Omega;\grid) = \prod_{T\in \grid} H^1(T),
$$
is endowed with the mesh dependent seminorm
\begin{equation}\label{E:dG-norm}
\index{Norms!$\dGnorm[a]{v}$: discontinuous Galerkin norm}
\dGnorm[a]{v}^2:= \| \nabla_\grid v \|_{L^2(\grid)}^2 +   a\| h^{-1/2} \jump{v}\|_{L^2(\faceswbdy)}^2,
\end{equation}
where $a$ is some positive parameter. We will prove below that this is indeed a norm.

With these notations, we can extend functionals $f \in \mathbb F_\grid$ in Definition~\ref{D:discrete-functionals} to $\mathbb V_\grid^{-1}$ for $\grid \geq \wh \grid$. 
Before doing so, recall that for $\wh f \in \mathbb F_{\wh \grid}$ and $v\in H^1_0(\Omega)$ we have
$$
\langle \wh f, v \rangle= \sum_{\wh T \in \wh \grid} \int_{\wh T} \wh f v + \sum_{\wh F \in \faces(\wh \grid)} \int_{\wh F} \wh f v,
$$
where, compared to Definition~\ref{D:discrete-functionals}, we slightly abused the notation 
$$
\wh f|_{\wh T} = \wh f_T \in \mathbb P_{2n-2}(\wh T) \qquad \wh f|_{\wh F} = \wh f_{\wh F} \in \mathbb P_{2n-1}(\wh F).
$$
In view of this, we can extend the duality pairing to $\mathbb V_\grid^{-1}$ by setting
\begin{equation}\label{e:ext_f-dG}
\langle \wh f, v \rangle_\grid := \sum_{T \in \grid} \int_T \wh f v + \sum_{F \in \faces} \int_F \wh f \mean{v}
\end{equation}
so that consistency in $H^1_0(\Omega)$ is preserved
\begin{equation}\label{e:consitency_f-dG}
\langle \wh f, v \rangle_\grid = \langle \wh f, v \rangle \qquad \forall v \in H^1_0(\Omega).
\end{equation}

%------------------------------------------------------------------------------------
\subsubsection{Interpolation operator $\IdG$}\label{S:IdG}
%------------------------------------------------------------------------------------
%
We shall need the interpolation operator $\IdG: \mathbb E_\grid \rightarrow \mathbb V_\grid^0$ from \cite{bonito2010quasi}. Its construction is based on an original idea of Cl\'ement \cite{clement1975approximation}, see also \cite{bernardi1998local} and other alternatives \cite{brenner2003poincare,bonito2021dg}.

Before embarking on the construction of  $\IdG$, we introduce a few notations. For an interior or boundary proper node
$P \in {\cal P}$ of the subdivision $\grid$, we denote by 
\begin{equation}\label{E:def-Vomega0P}
\mathbb{V}^0_{\omega_\grid(P)}:= H^1_0(\Omega) \cap \prod_{T\subset \mathcal
  \omega_\grid(P)} \mathbb P_n(T)
\end{equation}
  the space of continuous piecewise polynomial with support on the domain of influence $\omega_\grid(P)$ of $P$ and vanishing on $\partial \Omega$. When the underlying grid $\grid$ is clear from the context, we will simplify the notation and write $\mathbb{V}^0_P:=\mathbb{V}^0_{\omega_\grid(P)}$ and $\omega_P:= \omega_\grid(P)$; we refer to Figs.~\ref{nv-F:nonconforming-meshes} and \ref{fig:sample_element}.

 We now construct $\IdG$ in two steps.
 First, we define $V_P \in \mathbb{V}^0_P$ locally as satisfying
\begin{equation}\label{e:dg_interp}
\int_{\omega_P}(v-V_P)w = 0, \qquad 
\forall w \in \mathbb V^0_P.
\end{equation}
The value $V_P(P)$ is then used as the nodal value of $\IdG v$, namely
\begin{equation}\label{e:dg_int}
\index{Operators!$\IdG$: discontinuous Galerkin quasi-interpolant}
\IdG v  := \sum_{P\in {\cal P}}V_P(P)\psi_P,
\end{equation}
and we recall that $\{ \psi_P \}_{P \in {\cal P}}$ is a basis of $\mathbb V_\grid^0$ (see Section~\ref{S:tools}); note that $\IdG v=0$ on $\partial\Omega$ for all $v\in\mathbb E_\grid$. Moreover, including boundary proper nodes in the definition \eqref{e:dg_interp}-\eqref{e:dg_int} and replacing $H^1_0(\Omega)$ by $H^1(\Omega)$ in the definition \eqref{E:def-Vomega0P},
$\IdG$ easily extends to $\mathbb{S}^{n,0}_\grid$ without zero trace; we denote this operator
$\IdGp:\mathbb{E}_\grid\to\mathbb{S}^{n,0}_\grid$. An immediate property of $\IdG$
is {\it local invariance}:
\begin{equation}\label{E:IdG-invariance}
v\in \V^0_{\omega_T}
\quad\Rightarrow\quad
v=\IdG v \quad\textrm{in } T,
\end{equation}
where $\omega_T:= \bigcup \{\omega_P: P\in{\cal P},  T\subset\omega_P \}$; a similar property
is valid for $\IdGp$. We next gather a few more properties satisfied by $\IdG$. 

\begin{lemma}[interpolation operator]\label{L:IdG}
Let Assumption~\ref{A:initial-labeling} (initial labeling) hold and let $\grid \in \grids^\Lambda$.
For $v \in H^1_0(\Omega)$, there holds
\begin{equation}\label{e:dg_interpol_nonconfA}
\| v- \IdG v\|_{L^2(T)} \lesssim \|h \nabla v\|_{L^2(\omega_T)}, \quad \| \nabla \IdG v \|_{L^2(T)} \lesssim \|\nabla v\|_{L^2(\omega_T)},
\end{equation}
where $\omega_T$ is defined above. Instead, for $v \in \mathbb E_\grid$ there holds
\begin{equation}\label{e:dg_interpol_nonconfB}
\begin{aligned}
\| v - \IdG v \|_{L^2(T)} &+ \| h \nabla_\grid (v - \IdG v) \|_{L^2(T)} 
\\ & \lesssim \| h^{\frac{1}{2}} \jump{v}\|_{L^2(\faceswbdy \cap \omega_T)} + \|h\nabla_\grid (v-\Pi_\grid v) \|_{L^2(\omega_T)},
\end{aligned}
\end{equation}
where $\Pi_\grid$ is the $L^2$ projection operator onto $\V_\grid^{-1} = \mathbb{S}^{n,-1}_\grid$.
\end{lemma}
\begin{proof}
We start with \eqref{e:dg_interpol_nonconfA} and let $v\in H^1_0(\Omega)$. 
The definition \eqref{e:dg_interp} of the local projection $V_P\in\V_P^0$ yields for all $P \in {\cal P}$
$$
\|V_P\|_{L^2(\omega_P)}\le\|v\|_{L^2(\omega_P)}
\quad\Rightarrow\quad
\|V_P\|_{L^\infty(\omega_P)}\lesssim \diam \,(\omega_P)^{-\frac{d}{2}}\|v\|_{L^2(\omega_P)}.
$$
Proposition \ref{p-size-domain-influence} (size of the domain of influence) gives
$
\diam \omega_T \le C h_T,
$
whence the number of $\omega_P$ containing $T$ is uniformly bounded. Combining this with
the definition \eqref{e:dg_int} of $\IdG$ implies
\begin{equation}\label{e:dg_L2-stab}
\|\IdG v\|_{L^2(T)}  \lesssim \sum_{P\in {\cal P}:T\subset\omega_P} |V_P(P)| \, \|\psi_P\|_{L^2(T)}
\lesssim \|v\|_{L^2(\omega_T)} \quad \forall T \in \grid.
\end{equation}
Since $\IdG$ reproduces constants exactly locally, according to \eqref{E:IdG-invariance},
the first relation in \eqref{e:dg_interpol_nonconfA} follows from invoking the local $L^2$ stability property \eqref{e:dg_L2-stab} together with Proposition \ref{P:Bramble-Hilbert-Sobolev} (Bramble-Hilbert for Sobolev spaces). The second relation is proved using the same arguments and an inverse inequality
\begin{equation}\label{E:stab-grad}
\| \nabla \IdG v \|_{L^2(T)} \lesssim h_T^{-1} \inf_{v_0\in\R}\| \IdG (v-v_0) \|_{L^2(T)} \lesssim \| \nabla v \|_{L^2(\omega_T)}.
\end{equation}

We consider now $v\in\mathbb{E}_\grid$ and let $\wh{v}=\Pi_\grid v\in\V_\grid^{-1}$. We intend to prove \eqref{e:dg_interpol_nonconfB} by dealing with $v-\wh{v}$ and $\wh{v}$ separately and applying the 
triangle inequality. Since $v-\wh{v}$ has zero mean in $T$ according to \eqref{E:L2-projection-def},
we apply Lemma \ref{L:poincare-friedrichs} (second Poincar\'e inequality) to deduce
\[
\|v-\wh{v}\|_{L^2(T)} \lesssim h_T \| \nabla(v-\wh{v}) \|_{L^2(T)},
\]
whence, combining an inverse estimate with \eqref{e:dg_L2-stab}, we further infer that
\[
h_T\|\nabla \IdG (v-\wh{v})\|_{L^2(T)} \lesssim \|\IdG (v-\wh{v})\|_{L^2(T)} 
\lesssim \|v-\wh{v}\|_{L^2(\omega_T)} \lesssim h_T \|\nabla (v-\wh{v})\|_{L^2(\omega_T)}.
\]
This argument yields the inequality \eqref{e:dg_interpol_nonconfB} for $v-\wh{v}$.
It remains to deal with $\wh{v}$.

We scale $\omega_T$ to a reference domain with unit diameter. Estimate \eqref{e:domain_influence_control} on the size of the domains of influence guarantees that the number of such reference patches is uniformly finite over $\grids^\Lambda$. We relabel $\wh{v}$ as $v$ and examine the seminorm
$
\| \jump{v}\|_{L^2(\faceswbdy \cap \omega_P)}
$
on the space of discontinuous piecewise polynomials
$$
\left\lbrace v \in \Pi_{T\subset\omega_P}\mathbb P_n(T)| ~ V_P=0\right\rbrace,
$$
where $V_P$ is defined by \eqref{e:dg_interp}. If this seminorm vanishes then $v$ is continuous in
$\omega_P$ and thus $v\in\mathbb V_P^0$, whence the seminorm dominates any norm in this 
finite dimensional space. Consequently, scaling back gives
\begin{equation}\label{e:dg_scale}
\| v-V_P \|_{L^2(\omega_P)} + \| h \nabla_\grid(v-V_P) \|_{L^2(\omega_P)} 
\lesssim
\| h^{\frac12}\jump{v}\|_{L^2(\faceswbdy \cap \omega_P)}.
\end{equation}
We now deduce corresponding estimates for $\IdG$. 
For $T\in\mathcal T$ and $P,Q\in \omega_T \cap {\cal P}$, \eqref{e:dg_scale} implies that
\begin{equation}\label{e:dg_P-Q}
\| V_P - V_Q\|_{L^2(T)} + \| h \nabla_\grid(V_P - V_Q)\|_{L^2(T)} \lesssim  \|h^{\frac12} \jump{v}\|_{L^2(\faceswbdy \cap \omega_T)}.
\end{equation}
Consequently, the definition \eqref{e:dg_int} of $\IdG$ yields
\[
v - \IdG v = v - \sum_{P\in \omega_T\cap {\cal P}} V_P \psi_P
= (v-V_Q) - \sum_{P\in\omega_T \cap {\cal P}} (V_P-V_Q) \psi_P,
\]
which, combined with \eqref{e:dg_scale} and \eqref{e:dg_P-Q}, implies
$$
\| v - \IdG v\|_{L^2(T)} + \| h \nabla_\grid(v - \IdG v)\|_{L^2(T)} \lesssim \| h^{1/2} \jump{v} \|_{L^2(\faceswbdy \cap \omega_T)}.
$$
This is the desired estimate \eqref{e:dg_interpol_nonconfB} for $v=\wh{v}\in\V_\grid^{-1}$. To finish the proof we still need to express the right-hand side of the last inequality in terms of $v\in\mathbb{E}_\grid$. Applying the triangle inequality we are left with estimating $\|\jump{v-\wh{v}}\|_{L^2(F)}$ for any $F\in\faceswbdy \cap \omega_T$. If $T_F\in\grid$ is an element within $\omega_T$ that contains $F$ in its boundary, we employ the scaled trace inequality to arrive at
\[
h_T^{\frac{1}{2}} \|v-\wh{v}\|_{L^2(F)} \lesssim 
h_T \|\nabla(v-\wh{v})\|_{L^2(T_F)} + \|v-\wh{v}\|_{L^2(T_F)} 
\lesssim h_T \|\nabla(v-\wh{v})\|_{L^2(T_F)}.
\]
Finally, collecting all the estimates completes the proof.
\end{proof}

We discuss now consequences of Lemma \ref{L:IdG}. 
The first one is that jumps are solely responsible for controlling the discrepancy between
$v\in\V_\grid^{-1}$ and $\IdG v\in \V_\grid^0$:
\begin{equation}\label{e:dg_interpol_nonconfC}
\| v - \IdG v \|_{L^2(T)} + \| h \nabla_\grid (v - \IdG v) \|_{L^2(T)} 
\lesssim \| h^{\frac{1}{2}} \jump{v}\|_{L^2(\faceswbdy \cap \omega_T)},
\end{equation}
because $\Pi_\grid v = v$ in $\omega_T$. We next observe that \eqref{e:dg_interpol_nonconfB} is
also valid for $\IdGp$ with the same proof. We can thus apply \eqref{e:dg_interpol_nonconfB}
for $\IdGp$ to $w=v-\IdGp v$, use the invariance of $\IdGp$ in $\mathbb{S}^{n,0}_\grid$, and its continuity across internal faces in $\faces$, to deduce
\begin{equation}\label{e:dg_interpol_nonconfD}
\begin{aligned}
\| v - \IdGp v \|_{L^2(T)} &+ \| h \nabla_\grid (v - \IdGp v) \|_{L^2(T)} 
\lesssim \|h\nabla_\grid (v-\Pi_\grid v) \|_{L^2(\omega_T)}
\\ & + \| h^{\frac{1}{2}} \jump{v}\|_{L^2(\faces \cap \omega_T)} + 
\| h^{\frac{1}{2}} (v - \IdGp v) \|_{L^2(\partial\Omega \cap \omega_T)}.
\end{aligned}
\end{equation}
A third consequence of \eqref{e:dg_interpol_nonconfB} is the following Poincar\'e-type inequality 
on $\mathbb E_\grid$.

\begin{lemma}[Poincar\'e-type inequality on $\mathbb E_\grid$]\label{l:dG_Poincare}
Let $\grid \in \grids^\Lambda$ be a $\Lambda-$admissible refinement of $\grid_0$ satisfying Assumption~\ref{A:initial-labeling} (initial labeling). There exists $C_P = C_P(\Omega,\grid_0)$, such that for all $v \in \mathbb E_\grid$ there holds
\begin{equation}\label{E:dG_poincare}
\| v \|_{L^2(\Omega)} \leq  C_P \Big( \|\nabla_\grid v\|_{L^2\grid)} + \|h^{-\frac12} \jump{v} \|_{L^2(\faces)}
+ \|h^{-\frac12} v \|_{L^2(\partial\Omega)} \Big).
\end{equation}
In particular, if $v=0$ on $\partial\Omega$ then \eqref{E:dG_poincare} is a dG version of \eqref{E:Poincare}.
\end{lemma}
\begin{proof}
We argue locally with \eqref{e:dg_interpol_nonconfB}. First we realize that an argument similar to
\eqref{E:stab-grad} yields 
$\|\nabla_\grid(v-\Pi_\grid v)\|_{L^2(\omega_T)} \lesssim \|\nabla_\grid v\|_{L^2(\omega_T)}$,
whence adding over $T\in\grid$
\[
\| v - \IdG v \|_{L^2(\Omega)} + \| \nabla \IdG v \|_{L^2(\Omega)}
\lesssim \|\nabla_\grid v \|_{L^2(\Omega)}
+ \| h^{-\frac{1}{2}} \jump{v}\|_{L^2(\faceswbdy)}.
\]
It thus suffices to write
$$
\| v \|_{L^2(\Omega)} \leq \| \IdG  v \|_{L^2(\Omega)} + \| v - \IdG  v \|_{L^2(\Omega)},
$$
and invoke \eqref{E:Poincare} for $\IdG v\in H^1_0(\Omega)$ together with the preceding inequality.
\end{proof}

Another important property obtained using the interpolation operator $\IdG$ is that the approximation classes $\mathbb A_s^0 := \mathbb A_s(H^1_0(\Omega);\grid_0)$ defined using globally continuous piecewise polynomial approximations of degree $\le n$ on conforming subdivisions are equivalent to those without global continuity on $\Lambda$-admissible subdivisions $\grid\in\grids^\Lambda$, provided $\dGnorm[1][\grid]{\,.\,}$ (defined in \eqref{E:dG-norm}) is used as norm on $\mathbb E_\grid$. We define
\begin{equation}\label{E:sigmaN-1-IP}
\sigma_N^{n,-1}(v) := \inf_{\grid\in\grids^\Lambda_N} \inf_{v_\grid\in \mathbb{S}^{n,-1}_\grid}
\dGnorm[1][\grid]{v-v_\grid}
\end{equation}
and \index{Functional Spaces!$\mathbb A_s^{-1}$: approximation classes for $v$ for the discontinuous Galerkin norm}
$\mathbb A_s^{-1}:=\mathbb A_s^{-1}(H^1_0(\Omega);\grid_0)$ to be the class of functions $v \in H^1_0(\Omega)$ such that
$$
| v |_{\mathbb A_s^{-1}}:= \sup_{N \geq \# \grid_0}\left( N^s \sigma_N^{n,-1}(v) \right) < \infty \quad \Rightarrow \quad \sigma_N^{n,-1}(v) \leq | v |_{\mathbb A_s^{-1}}N^{-s}.
$$
Note the scaling parameter $a$ for jumps in the definition of $\sigma_N^{n,-1}$ is just $a=1$.

The following result can be traced back to \cite{bonito2010quasi}. 

\begin{proposition}[equivalence of classes for $u$]\label{p:equiv_classes_dG}
Let $\grid_0$ be an initial conforming subdivision satisfying Assumption~\ref{A:initial-labeling} (initial labeling). There are two constants $m\in \mathbb N$ and $C \geq 1$ such that for all $N \geq \#\grid_0$ and all $v \in H^1_0(\Omega)$
$$
\sigma_{N}^{n,-1}(v) \leq \sigma_{N}^{n,0}(v) \qquad \textrm{and} \qquad \sigma_{mN}^{n,0}(v) \leq C \sigma_{N}^{n,-1}(v).
$$
In particular, the approximation classes coincide $\mathbb A_s^0 \equiv \mathbb A^{-1}_s$, $s \geq 0$. 
\end{proposition}
\begin{proof}
The start with the first inequality.
For $v \in H^1_0(\Omega)$ and $N \geq \# \grid_0$, we let $\grid \in \grids_N$ be a conforming subdivision of $\grid_0$ and $v^0_{\grid} \in \mathbb V_{\grid}^0 \subset \mathbb V_{\grid}^{-1}$ be such that
$$
\sigma_N^{n,0}(v) = | v - v^0_{\grid} |_{H^1_0(\Omega)}.
$$
Because $v-v^0_{\grid} \in H^1_0(\Omega)$, we have $\dGnorm[1][\grid]{v-v^0_{\grid}} = | v - v^0_{\grid} |_{H^1_0(\Omega)}$ and thus
$$
\sigma_N^{n,-1}(v) \leq | v - v^0_{\grid} |_{H^1_0(\Omega)} = \sigma_N^{n,0}(v).
$$

We now prove the second inequality. For
$v \in H^1_0(\Omega)$ and $N \geq \#\grid_0$, let  $\mathcal T \in \mathbb T^\Lambda_N$ be a $\Lambda$-admissible mesh with $N$ elements and $v_{\grid} \in \mathbb V_{\grid}^{-1}$ be so that
$$
\dGnorm[1][\grid]{v-v_{\grid}} = \sigma_{N}^{n,-1}(v).
$$
We first show that $\IdG[\grid]v_{\grid} \in \mathbb V_{\grid}^0$ satisfies
$$
|v-\IdG[\grid]v_{\grid}|_{H^1_0(\Omega)} \lesssim  \sigma_{N}^{n,-1}(v).
$$
Indeed, using the triangle inequality we obtain
$$
\dGnorm[1][\grid]{v-\IdG[\grid]v_{\grid}}
\leq \dGnorm[1][\grid]{v-v_{\grid}} + \dGnorm[1][\grid]{v_{\grid}-\IdG[\grid]v_{\grid}}.
$$
Interpolation estimate \eqref{e:dg_interpol_nonconfC} yields
\begin{equation}\label{e:atojump}
\dGnorm[1][\grid]{v_{\grid}-\IdG[\grid]v_{\grid}} \lesssim
\| h^{-1/2} \jumpinline{v_{\grid}} \|_{L^2(\faceswbdy)},
\end{equation}
because $v_{\grid}-\IdG[\grid]v_{\grid} \in \V^{-1}_{\grid}$, whence
\begin{equation*}
| v-\IdG[\grid] v_{\grid}|_{H^1_0(\Omega)} = \dGnorm[1][\grid]{v-\IdG[\grid]v_{\grid}} \leq C \sigma_N^{n,-1}(v)
\end{equation*}
as claimed for a constant $C\ge1$ independent of $v$ and $N$. 
To assert an estimate on $\sigma_N^{n,0}(v)$, we now exhibit a conforming refinement $\overline{\grid}$ of $\grid$ with a comparable number of elements. To do this, we note that because $\grid\in \grids^\Lambda$ is $\Lambda$-admissible, it is the product of successive calls 
$[\grid_j]=\REFINE(\grid_{j-1}, T_{j-1})$, $j=1,...,J$, where  $\grid_j$ is the smallest {\it $\Lambda-$admissible} refinement of $\grid_{j-1}$ such that the element $T_{j-1}\in\grid_{j-1}$ 
is bisected once. We now let $\overline{\grid}\in\grids$ be the conforming subdivision obtained from the successive calls $[\overline{\grid}_j]=\REFINE(\overline{\grid}_{j-1}, \{T_{j-1}\}\cap \overline{\grid}_{j-1})$ with $\overline{\grid}_0=\grid_0$ but where this time $\REFINE$ produces the smallest {\it conforming} refinement of $\overline{\grid}_{j-1}$ where the element of $T_{j-1}$ is bisected once if $T_{j-1} \in \overline{\grid}_{j-1}$ or otherwise $\overline{\grid}_{j}=\overline{\grid}_{j-1}$. A simple induction argument, exploiting the
minimality of the meshes generated by $\REFINE$, reveals that
$\overline{\grid}_{j}\ge\grid_{j}$ for $0\le j \le J$. Consequently,
Theorem~\ref{nv-T:complexity-refine} (complexity of \REFINE) guarantees that
$$
\# \overline{\grid} - \#\grid_0 \leq \Ccompl \sum_{j=0}^{J-1} \# \big(\{T_{j-1}\} \cap \overline{\grid}_{j-1} \big) \leq \Ccompl J \leq \Ccompl \big(\# \grid - \#\grid_0 \big)
$$
whence $\# \overline{\grid} \leq \Ccompl \# \grid \leq m N$ with $m:=\lceil \Ccompl \rceil$ because $\Ccompl\ge1$. 

Therefore, $\mathbb V_{\grid}^0 \subset \mathbb V_{\overline{\grid}}^0$ because $\overline{\grid}$ is a conforming refinement of $\grid$. Since $\# \overline{\grid} \leq m N$ and $\IdG[\grid]v_{\grid} \in \mathbb V_{\grid}^0$, we deduce
$$
\sigma_{mN}^{n,0}(v) \le | v-\IdG[\grid] v_{\grid}|_{H^1_0(\Omega)} \leq C \sigma_N^{n,-1}(v),
$$
which is the desired inequality. Finally, the equivalence of moduli of approximation yields $\mathbb A_s^0 \equiv \mathbb A^{-1}_s$ and completes the proof.
\end{proof}

\begin{remark}[equivalent classes for $\data$]\label{r:equiv_class_data}
The approximation classes for data $\data=(\bA,c,f)$, namely $\mathbb M_{s}((L^r(\Omega))^{d\times d});\grid_0$, $\mathbb C_s(L^q(\Omega);\grid_0)$ and $\mathbb F_s(H^{-1}(\Omega);\grid_0)$, are defined for conforming subdivisions in Section~\ref{S:conv-rates-coercive}. However, repeating the construction of the smallest conforming refinement $\overline{\grid}$ of any $\Lambda$-admissible subdivision $\grid$, and using the fact that $\# \overline{\grid} \approx \# \grid$ proved above, we deduce that these classes are equivalent to their counter-parts on non-conforming meshes. 
Therefore, we do not repeat the proof here and use from now on the same notation to denote the approximation classes on $\Lambda$-admissible subdivisions. 
\end{remark}

%---------------------------------------------------------------------------------
\subsection{Discontinuous Galerkin formulation}\label{ss:IP}
%---------------------------------------------------------------------------------
%
This section discusses the \SOLVE routine at the core of the module \GALERKINIP.
Recall that within the two-step method \AFEMTSIP, data $\data=(\bA,c,f)$ is approximated by $\wh \data=(\wh \bA,\wh c,\wh f) \in \mathbb D_{\wh \grid}$ subordinate to a partition $\wh \grid \in \grids^\Lambda$. For a subdivision $\grid \in \grids^\Lambda$, $\grid \geq \wh \grid$, the Galerkin solution $[u_\grid]=\SOLVE(\grid)$ is constructed to approximate $\wh u = u(\wh\data)\in H^1_0(\Omega)$, the exact weak solution of the perturbed problem \eqref{E:perturbed-weak-form} with approximate data $\wh \data=(\wh{\bA},\wh c,\wh f)$ constructed using Algorithm~\ref{A:data} ($\DATA$). Corollary~\ref{C:data} (performance of \DATA) guarantees that the output $[\wh \data,\wh \grid]$ of  $\DATA$ satisfies the structural assumption
 \begin{equation}\label{e:structural-IP}
 \wh{\bA} \in M(\halpha_1,\halpha_2), \qquad \wh c \in R(\wh c_1,\wh c_2)
 \end{equation} with $0<\halpha_1 \leq \halpha_2$ and $-\frac {\halpha_1}{2 C_P^2} \leq \wh c_1 \leq \wh c_2$
upon replacing the Poincar\'e constant $C_P$ by the larger constant $C_P$ appearing in Lemma~\ref{l:dG_Poincare} in Algorithm~\ref{A:const-project-c} (\APPLYCONSTRC).
We do not specify the dependency on $\halpha_1$, $\halpha_2$, $\wh c_1$, and $\wh c_2$ of the constants appearing in the analysis below. We also emphasize that the constants involved in \eqref{e:structural-IP} do not depend on $\wh \grid$ and are thus uniform among all the discrete data constructed within $\AFEMTSIP$. 

Relation \eqref{e:structural-IP} not only ensures the existence and uniqueness of a solution $\wh u \in H^1_0(\Omega)$ satisfying the perturbed problem \eqref{E:perturbed-weak-form} but also, as we shall see in Corollary~\ref{c:dG_cont_coer}, the existence and uniqueness of its discontinuous Galerkin approximation.   
We first present the standard symmetric interior penalty method and point out its consistency requires the exact solution $u\in H^{s}(\Omega)$, $s>3/2$. To circumvent this rather restrictive assumption, we introduce lifting operators allowing a reformulation valid in $H^1(\Omega)$.
However, this reformulation is only consistent on the conforming subspace $\mathbb V_\grid^0 = \mathbb V_\grid^{-1} \cap H^1_0(\Omega)$ and requires our analysis to decompose the discrete space $\mathbb V_\grid^{-1}$ into $\mathbb V_\grid^0$ and its complement $\V_\grid^\perp$ with respect to an appropriate scalar product. 

%---------------------------------------------------------------------------------
\subsubsection{The symmetric interior penalty method}
%---------------------------------------------------------------------------------
%
The symmetric interior penalty (SIP) formulation is the most standard discontinuous Galerkin method. For $\grid \in \grids^\Lambda$, it consists in finding $u_\grid \in \mathbb V_\grid^{-1}$ satisfying
\begin{equation}\label{e:IP}
\bilin[\grid]{u_\grid}{v} = \langle \wh f, v \rangle_\grid, \qquad \forall v \in \mathbb V_\grid^{-1},
\end{equation}
where $\mathcal B_\grid:\mathbb V_\grid^{-1} \times \mathbb V_\grid^{-1} \rightarrow \mathbb R$ is the bilinear form defined by
\begin{equation}\label{e:bilin_IP}
\begin{split}
\bilin[\grid]{w}{v}:= &\int_\Omega (\nabla_\grid v \cdot \wh{\bA}\nabla_\grid w + \wh c wv)- \sum_{F\in \faceswbdy} \int_F \jump{v}\vec{n}_F \cdot \meaninline{\wh{\bA} \nabla_\grid w} \\ 
&-\sum_{F \in \faceswbdy} \int_F \jump{w}\vec{n}_F \cdot \meaninline{\wh{\bA} \nabla_\grid v} + \kappa \sum_{F \in \faceswbdy}\int_F h_F^{-1}\jump{w}\jump{v}.
\end{split}
\end{equation}
The parameter $\kappa>0$ is responsible for keeping the discontinuity of  the Galerkin solution under control and its value is discussed below. Unless specified otherwise, all the constant appearing in the discussion below are independent of $\kappa$ and the notation $A \lesssim B$ signifies $A \leq C B$ with a constant $C$ independent of the discretization parameters and $\kappa$.

A few comments regarding the weak formulation \eqref{e:IP} are in order. An integration by parts reveals that the method is consistent whenever the exact solution satisfies the additional regularity $u \in H^{s}(\Omega)$, $s>3/2$. However, we do not make this assumption in the analysis below but rather extend the formulation to the energy space $\mathbb E_\grid \supset \mathbb V_\grid^{-1}$ using lifting operators. The same integration by parts also indicates that the term $\sum_{F \in \faceswbdy}\int_F \jump{w}\vec{n}_F \cdot \meaninline{\wh{\bA} \nabla v}$ is not necessary but included to achieve a symmetric formulation. Recall that $\wh{\bA}$ constructed by $\DATA$ is symmetric.
In addition, the presence of $\langle \wh f, v \rangle_\grid$ is not standard but allows for  right-hand sides $\wh f \in \mathbb F_{\wh \grid}$ and in turn for $f \in H^{-1}(\Omega)$ within the $\AFEMTSIP$ algorithm.

%--------------------------------------------------------------------------------
\subsubsection{Lifting operators}
%--------------------------------------------------------------------------------

The interior penalty bilinear form \eqref{e:bilin_IP} includes inter-element terms 
\begin{equation}\label{e:inter-element-terms}
\sum_{F\in \faceswbdy} \int_F \jump{v}\vec{n}_F \cdot \meaninline{\wh{\bA} \nabla_\grid w} +\sum_{F \in \faceswbdy} \int_F \jump{w}\vec{n}_F \cdot \meaninline{\wh{\bA} \nabla_\grid v},
\end{equation}
which are not defined on $H^1(\Omega)$ but on $H^{s}(\Omega)$, $s>3/2$.  In turn, the method is consistent when $u \in H^{s}(\Omega)$, $s>3/2$.
The key ingredient to extend $\bilin[\grid]{w}{v}$ to $\mathbb E_\grid\times \mathbb E_\grid$ without additional regularity is a lifting operator \cite{brezzi2000discontinuous,arnold2002unified,perugia2003,houston2004mixed,houston2007energy,bonito2010quasi} introduced in this section.

For $n'>0$, we define $\mathcal L^{n'}_\grid: \mathbb E_\grid \rightarrow [\mathbb S_\grid^{n',-1}]^d$ by the relations
\begin{equation}\label{d:lift}
\int_\Omega \mathcal L^{n'}_\grid[v] \cdot \wh{\bA}\vec{w}= \sum_{F \in \faceswbdy} \int_F \jump{v} \vec{n}_F \cdot \meaninline{\wh{\bA}\vec{w}}, \qquad \forall \vec{w} \in [\mathbb S_\grid^{n',-1}]^d.
\end{equation}
From this definition, we easily deduce a $L^2$ stability estimate.

\begin{lemma}[stability of lift]\label{l:estimate_lift}
Let $\grid \in \grids^\Lambda$ be a $\Lambda-$admissible subdivision of $\grid_0$ satisfying Assumption~\ref{A:initial-labeling} (initial labeling). Assume $\wh{\bA} \in M(\halpha_1,\halpha_2)$ with $0<\halpha_1 \leq \halpha_2$. For $n'\geq 0$ and all  $v \in \mathbb S_\grid^{n',-1}$, there holds
\begin{equation}\label{e:estimate_lift}
\|\mathcal L_\grid^{n'}[v] \|_{L^2(\Omega)} \leq C \| h^{-1/2} \jump{v}\|_{L^2(\faceswbdy)},
\end{equation}
where $C=C(\halpha_2/\halpha_1,n',\grid_0)$.
\end{lemma}
\begin{proof}
Let $v\in \mathbb S_\grid^{n',-1}$ and set $\vec{w}=\mathcal L_\grid^{n'}[v]$ in \eqref{d:lift} to write
\begin{equation*}
\begin{split}    
\|\wh{\bA}^{1/2}\mathcal L^{n'}_{\grid}[v]\|_{L^2(\Omega)}^2 &=
\int_\Omega \mathcal L_{\grid}^{n'}[v] \cdot \wh{\bA} \mathcal L_{\grid}^{n'}[v] \\
& = \sum_{F \in \faceswbdy} \int_F h^{-1/2} \jump{v} \vec{n}_F \cdot h^{1/2} \meaninline{\wh{\bA} \mathcal L_\grid^{n'}[v]} \\
&\leq 
\| h^{-1/2}\jump{v}\|_{L^2(\faceswbdy)} \|h^{1/2} \meaninline{\wh{\bA} \mathcal L_\grid^{n'}[v]} \|_{L^2(\faceswbdy)}.
\end{split}
\end{equation*}
A local inverse estimate along with the eigenvalue bounds for $\wh{\bA} \in M(\halpha_1,\halpha_2)$ yields
$$
\|h^{1/2} \meaninline{\wh{\bA} \mathcal L_\grid^{n'}[v]} \|_{L^2(\faceswbdy)}\leq C \halpha_2 \| \mathcal L_\grid^{n'}[v]\|_{L^2(\Omega)},
$$
where $C$ only depends on $n'$ and on the shape regularity constant of $\grid_0$.
Combining the above two inequalities and taking advantage again of the assumption $\wh{\bA} \in M(\halpha_1,\halpha_2)$ implies \eqref{e:estimate_lift}. 
\end{proof}

We record two estimates based on \eqref{e:estimate_lift} and used multiple times in the analysis below. Combining the estimate \eqref{e:estimate_lift} on the lifting operator with assumption \eqref{e:structural-IP}  and a Cauchy-Schwarz inequality we find that
\begin{equation}\label{e:estimate_lift2}
\int_\Omega \mathcal L_\grid^{n'}[v] \cdot \wh{\bA} \nabla_\grid w \leq C_{\textrm{lift}} \| h^{-1/2} \jump{v} \|_{L^2(\faceswbdy)} \| \nabla_\grid w\|_{L^2(\grid)}, \ \  \forall v,w \in \mathbb E_\grid,
\end{equation}
for a constant $C_{\textrm{lift}}=C_{\textrm{lift}}(\halpha_1,\halpha_2,n',\grid_0)$ and in particular independent of the discretization parameters and $\kappa$. 
This, together with a Young inequality, yields for any $\epsilon>0$ the second estimate for all 
$v,w \in \mathbb E_\grid$
\begin{equation}\label{e:estimate_lift3}
\int_\Omega \mathcal L_\grid^{n'}[v] \cdot \wh{\bA} \nabla_\grid w \leq  \frac{C^2_{\textrm{lift}}}{2\epsilon}\| h^{-1/2} \jump{v} \|_{L^2(\faceswbdy)}^2 + \frac \epsilon 2  \| \nabla_\grid w\|_{L^2(\grid)}^2.
\end{equation}

We now return to the SIP weak formulation \eqref{e:IP} and take advantage of the lifting operators to deduce an equivalent expression of the bilinear form $\mathcal B_\grid$ on $\mathbb V_\grid^{-1}$, which is well defined on $\mathbb  E_\grid$. The problematic inter-elements terms \eqref{e:inter-element-terms} are equivalently rewritten as
\begin{equation}\label{e:dG_lift_term}
\int_\Omega \mathcal L^{n'}_\grid[v] \cdot \wh{\bA}\nabla_\grid w + \int_\Omega \mathcal L^{n'}_\grid[w] \cdot \wh{\bA}\nabla_\grid v  
\end{equation}
provided 
$$
\nabla_\grid \mathbb V_\grid^{-1} \subset [\mathbb S_\grid^{n',-1}]^d.
$$
The above condition is satisfied when $n' \geq n-1$ for subdivisions $\grid$ made of simplices and $n' \geq n$ for hexahedra. To continue with an analysis incorporating both case, we set $n'=n$ and write $\mathcal L_\grid := \mathcal L_\grid^{n}$. With this choice, the bilinear form $\mathcal B_\grid$ in the symmetric interior penalty method \eqref{e:IP} reads
\begin{equation}\label{e:bilinIP2}
\begin{split}
\B_\grid[w,v] = &\ a_\grid[w,v] - \int_\Omega \mathcal L_\grid[v] \cdot \wh{\bA}\nabla_\grid w \\
&- \int_\Omega \mathcal L_\grid[w] \cdot \wh{\bA}\nabla_\grid v + \kappa \sum_{F \in \faceswbdy}\int_F h_F^{-1}\jump{w}\jump{v},
\end{split}
\end{equation}
for all $w,v \in \mathbb V_\grid^{-1}$ and where we used
\begin{equation}\label{d:dg_bilin_a}
    a_\grid[w,v]:= \int_\Omega \nabla_\grid v \cdot \wh{\bA} \nabla_\grid w + \wh c wv
 \end{equation}
 to denote the bilinear form related to the conforming method.

Expression \eqref{e:bilinIP2} is well defined for $w,v \in \mathbb E_\grid$ and the weak formulation \eqref{e:IP} is well-posed. These two claims follow from Corollary~\ref{c:dG_cont_coer} below, which in turn is a consequence of the next result focusing on the bilinear form $a_\grid$; we recall 
\eqref{E:wh-alpha}.

\begin{lemma}[properties of $a_\grid$]\label{l:prop_a_grid}
Let $\grid \in \grids^\Lambda$ be a $\Lambda-$admissible refinement of $\grid_0$ satisfying Assumption~\ref{A:initial-labeling} (initial labeling). Furthermore, assume that $\wh{\bA}$ and $\wh c$ satisfy the structural assumption \eqref{e:structural-IP}. Then, we have
\begin{equation}\label{e:dG_cont_a}
a_\grid[w,v] \leq (\halpha_2+|\wh c_2| C_P^2) \dGnorm[1]{v} \dGnorm[1]{w}, \qquad \forall v,w \in \mathbb E_\grid
\end{equation}
and
\begin{equation}\label{e:dG_partial_coer_a}
a_\grid[v,v] \geq \frac{\halpha_1}{2} \| \nabla_\grid v\|_{L^2(\Omega)}^2 +\min(0, \wh c_1) C_P^2\| h^{-1/2} \jump{v} \|_{L^2(\faceswbdy)}^2, \qquad \forall v \in \mathbb E_\grid,
\end{equation}
where $C_P$ is the constant in Lemma~\ref{l:dG_Poincare} 
(Poincar\'e-type inequality in $\mathbb{E}_\grid$).
\end{lemma}
\begin{proof}
We start with the continuity estimate \eqref{e:dG_cont_a}. The assumption on the discretized coefficients implies that for $v,w \in \mathbb E_\grid$ there holds
$$
a_\grid[w,v] \leq \alpha_2 \| \nabla_\grid w \|_{L^2(\grid)} \| \nabla_\grid v \|_{L^2(\grid)} + |c_2| \| w \|_{L^2(\Omega)}\| v \|_{L^2(\Omega)}.
$$
It remains to invoke Lemma~\ref{l:dG_Poincare} (Poincar\'e-type inequality on $\mathbb{E}_\grid$)
to deduce \eqref{e:dG_cont_a}.

Similarly for the partial coercivity estimate \eqref{e:dG_partial_coer_a}, we have
$$
a_\grid[v,v] \geq \halpha_1 \| \nabla_\grid v\|_{L^2(\Omega)}^2 + \wh c_1 \| v \|_{L^2(\Omega)}^2 \geq \halpha_1 \| \nabla_\grid v\|_{L^2(\Omega)}^2 +\min(0, \wh c_1) C_P^2 \dGnorm[1][\grid]{v}^2
$$
and the desired estimate follows from the assumption $-\frac{\halpha_1}{2 C_P^2} \leq \wh c_1$. 
\end{proof}

For the next result, we recall that the discrete norm $\|.\|_{\kappa,\grid}$ is defined in \eqref{E:dG-norm}.
\begin{corollary}[properties of $\cal B_\grid$]\label{c:dG_cont_coer}
    Let $\grid \in \grids^\Lambda$ be a $\Lambda-$admissible refinement of $\grid_0$ satisfying Assumption~\ref{A:initial-labeling} (initial labeling). Furthermore assume that  $\wh{\bA}$ and $\wh c$ satisfy the structural assumption \eqref{e:structural-IP}. There exist a constant $C_\textrm{cont}$ such that    
    \begin{equation}\label{e:cont_Bdg}
    \bilin[\grid]{w}{v} \leq C_\textrm{cont} \dGnorm{v} \dGnorm{w} \qquad \forall v,w \in \mathbb E_\grid.
    \end{equation}
    Moreover, there are constants $\overline{\kappa}_{\textrm{stab}}, C_\textrm{coer}>0$ such that for all $\kappa > \overline{\kappa}_{\textrm{stab}}$, there holds
    \begin{equation}\label{e:coerc_Bdg}
    C_\textrm{coer}  \dGnorm{v}^2 \leq \bilin[\grid]{v}{v} \qquad \forall v \in \mathbb E_\grid.
    \end{equation}
    In particular, the Galerkin formulation \eqref{e:IP} has a unique solution $u_\grid \in \mathbb V_\grid^{-1}$.
\end{corollary}
\begin{proof}
   The continuity estimate is a direct consequence of the continuity estimate \eqref{e:dG_cont_a}, estimate \eqref{e:estimate_lift2} for the lifting terms, and Cauchy-Schwarz inequality
   $$
   \kappa \sum_{F \in \faceswbdy} \int_F h_F^{-1} \jump{w}\jump{v} \leq \kappa \| h^{-1/2} \jump{w}\|_{L^2(\faceswbdy)} \| h^{-1/2} \jump{v}\|_{L^2(\faceswbdy)},
   $$
   which holds for all $v,w \in \mathbb E_\grid$. 
   
   We now focus on the coercivity estimate \eqref{e:coerc_Bdg} and start from \eqref{e:dG_partial_coer_a}, which we write for $v\in \mathbb E_\grid$ as 
\begin{equation}
    \label{e:dg_coercivity_rel}
    \begin{split}
\frac{\halpha_1}{2} \| \nabla_\grid v \|_{L^2(\grid)}^2 - \max\big(0,-c_1  \wh C_P^2\big) \| h^{-1/2} \jump{v} \|_{L^2(\faceswbdy)}^2 \leq  a_\grid[v,v] .
    \end{split}
\end{equation}
Furthermore, the terms involving the lifting operators in the definition \eqref{e:bilinIP2} of the bilinear form $\mathcal B_\grid$ reduce to
$
-2 \int_{\Omega} \mathcal L_\grid[v] \cdot \wh{\bA} \nabla_\grid v
$
when $w=v$. 
Hence, the estimate \eqref{e:estimate_lift3} with $\epsilon = \halpha_1/4$ implies that
$$
2\left| \int_{\Omega} \mathcal L_\grid[v] \cdot \wh{\bA} \nabla_\grid v \right| \leq \frac{\halpha_1}{4} \| \nabla_\grid v\|_{L^2(\grid)}^2 + \frac{4 C_{\mathcal L}^2}{\halpha_1}\| h^{-1/2} \jump{v}\|_{L^2(\faceswbdy)}^2.
$$
Gathering the above inequalities and recalling definition \eqref{e:bilinIP2} of $\mathcal B_\grid$, we find that
\begin{equation*}
%\label{e:dG_coercivity_plus}
\frac{\alpha_1}{4}\| \nabla_\grid v\|_{L^2(\grid)}^2 + (\kappa-\overline{\kappa}_{\textrm{stab}}) \| h^{-1/2} \jump{v}\|_{L^2(\faceswbdy)}^2 \leq \bilin[\grid]{v}{v}
\end{equation*}
with 
\begin{equation*}
\overline{\kappa}_{\textrm{stab}}:= \frac{4 C_{\mathcal L}^2}{\alpha_1}+\max\big(0,-\wh c_1 C_P^2\big).
\end{equation*}
The desired coercivity estimate directly follows provided $\kappa > \overline{\kappa}_{\textrm{stab}}$.
\end{proof}

%-----------------------------------------------------------------------------------
\subsubsection{Partial consistency and role of the conforming Galerkin solution}\label{ss:dg_conforming}
%-----------------------------------------------------------------------------------
%
From now on, we shall use the expression of $\mathcal B_\grid$ in \eqref{e:bilinIP2} extending $\mathcal B_\grid$ to $\mathbb E_\grid \times \mathbb E_\grid$.
This reformulation comes at the price of a partial consistency. Since $\mathcal L_\grid[v]=0$ whenever $v \in H^1_0(\Omega)$ and the duality product $\langle \cdot,\cdot\rangle_\grid$ satisfies the consistency \eqref{e:consitency_f-dG}, we have
\begin{equation}\label{e:dg_partialconsistency}
\bilin[\grid]{\wh u}{v} = a_\grid[\wh{u},v] =  \langle \wh f,v \rangle_\grid \qquad \forall v \in H^1_0(\Omega),
\end{equation}
which indicates that the reformulation \eqref{e:bilinIP2} using lifts is consistent on $H^1_0(\Omega)$. However, \eqref{e:dg_partialconsistency} does not hold for all $v \in \mathbb V_\grid^{-1}$. 

This suggests splitting $\mathbb V_\grid^{-1}$ into a conforming space where the consistency holds and its orthogonal complement. 
We decompose the discontinuous space as
\begin{equation}\label{e:conf-nonconf-space-decomp}
\mathbb V_\grid^{-1} = \mathbb V_\grid^{0} \oplus \mathbb V_\grid^{\perp},
\end{equation}
where $\mathbb V_\grid^{0} = \mathbb V_\grid^{-1}\cap H^1_0(\Omega)$ is the finest conforming subspace of $\mathbb V_\grid^{-1}$ and $\mathbb V_\grid^{\perp}$ is the orthogonal complement with respect to the $\bilin[\grid]{\cdot}{\cdot}$ scalar product. 
Note that the later is well defined provided the assumption on the penalty parameter $\kappa > \overline{\kappa}_{\textrm{stab}}$, required by Corollary~\ref{c:dG_cont_coer}, is satisfied. From now on, we assume this is the case and point out that although the constants appearing in the analysis below do not depend on $\kappa$, they may depend on $\overline{\kappa}_{\textrm{stab}}$. 

We also emphasize that there might not be a conforming subdivision associated with $\mathbb V_\grid^0$. The latter is the span of the basis functions associated with proper nodes; see Figure~\ref{f-basis-function} for an illustration and refer to Section~\ref{S:mesh-refinement} for more details. Consequently, the analysis provided  below relies on the decomposition \eqref{e:conf-nonconf-space-decomp} of the space $\mathbb V_\grid$ rather than on a subdivision $\grid$. It is also worth pointing out that the conforming part $u_\grid^0 \in \mathbb V_\grid^0$ of the Galerkin solution $u_\grid \in \mathbb V_\grid^{-1}$ satisfies
\begin{equation}\label{e:decomp_conf_eq}
\bilin[\grid]{u_\grid^{0}}{v} = \langle \wh f,v \rangle, \qquad \forall v \in \mathbb V_\grid^{0}.
\end{equation}
Hence, $u_\grid^0$ is the conforming Galerkin approximation on $\mathbb V_\grid^{0}$.
As we shall see this finest coarser conforming Galerkin solution plays a critical role in the convergence of \AFEMTSIP. This justifies the orthogonal decomposition \eqref{e:conf-nonconf-space-decomp} associated with the $\cal B_\grid$ scalar product. 

Another advantage of using the $\mathcal B_\grid$-orthogonal decomposition \eqref{e:conf-nonconf-space-decomp} is that it offers a control on the non-conforming component of $v \in \mathbb V_\grid^{-1}$ by its scaled jumps. 
To achieve this, the operator $\IdG $ defined by \eqref{e:dg_int} is instrumental.

\begin{lemma}[control of non-conformity]\label{l:control_perp}
Let $\grid \in \grids^\Lambda$ be a $\Lambda-$admissible refinement of $\grid_0$ satisfying Assumption~\ref{A:initial-labeling} (initial labeling). Assume that $\wh{\bA}$ and $\wh c$ satisfy the structural assumption \eqref{e:structural-IP}. For  $\kappa > \overline{\kappa}_{\textrm{stab}}$, if $v = v^0 + v^\perp \in \mathbb V_\grid^{-1}$ according to \eqref{e:conf-nonconf-space-decomp}, then
$$
\dGnorm{v^\perp} \lesssim \kappa^{1/2}\| h^{-1/2} \jump{v^\perp} \|_{L^2(\faceswbdy)} = \kappa^{1/2} \| h^{-1/2} \jump{v} \|_{L^2(\faceswbdy)}.
$$
\end{lemma}
\begin{proof}
    Because $\IdG  v \in \mathbb V_\grid^{0}$, the orthogonal decomposition \eqref{e:conf-nonconf-space-decomp} implies that
    $$
    \bilin[\grid]{v^\perp}{v^\perp} \leq \bilin[\grid]{v-\IdG  v}{v-\IdG  v}.
    $$
    The desired result follows from the coercivity \eqref{e:coerc_Bdg} and continuity \eqref{e:cont_Bdg} of $\mathcal B_\grid$ along with the interpolation estimate \eqref{e:dg_interpol_nonconfB}.
\end{proof}

%----------------------------------------------------------------------------------
\subsection{A posteriori error estimates}\label{ss:posteriori-dG}
%----------------------------------------------------------------------------------
%
We derive a residual error estimate for the discontinuous Galerkin method. Because the data $\wh \data \in \mathbb D_{\grid}$ is discrete, the analysis is free from data oscillation. In the notation introduced in Section~\ref{S:aposteriori}, this means $\est_{\grid} = \eta_\grid$ where for $v \in \mathbb E_\grid$
$$
\eta_\grid^2(v):= \sum_{T \in \grid} \eta_\grid(v,T)^2, 
$$
and
\begin{equation*}
 \eta_\grid(v,T)^2
 := h_T  \sum_{F \subset \partial T \setminus \partial \Omega} \| j_\grid(v)-\wh f\|_{L^2(F)}^2+
 h_T^2 \| r_\grid(v) \|_{L^2(T)}^2
 \end{equation*}
 with $j_\grid(v)|_F:=\vec{n}_F \cdot \jumpinline{\wh{\bA} \nabla_\grid v}$ and $r_\grid(v)|_T:=\wh f - \wh c v + \div_\grid ( \wh{\bA} \nabla_\grid v)$.
 
We start with a result mimicking the conforming argument, discuss its drawback and improvement next.

\begin{lemma}[a posteriori error estimates]\label{L:aposteriori-dg}
Let $\grid \in \grids^\Lambda$ be a $\Lambda-$admissible refinement of $\grid_0$ satisfying Assumption~\ref{A:initial-labeling} (initial labeling). Assume that $\wh{\bA}$ and $\wh c$ satisfy the structural assumption \eqref{e:structural-IP}.
If $\kappa > \overline{\kappa}_{\textrm{stab}}$, then there holds
\begin{equation}\label{e:dg_upper}
\dGnorm{\wh u - u_\grid}^2 \lesssim  \eta_\grid(u_\grid)^2 + \kappa \| h^{-1/2} \jump{u_\grid} \|_{L^2(\faceswbdy)}^2
\end{equation}
and
\begin{equation}\label{e:dg_lower}
C_\textrm{L} \, \eta_\grid(u_\grid)   \leq  \dGnorm{\wh u-u_\grid},
\end{equation}
for some constant $C_{\textrm{L}}$.
\end{lemma}
\begin{proof}
We start with the upper bound \eqref{e:dg_upper}. To exploit the consistency \eqref{e:dg_partialconsistency} in $\mathbb V_\grid^0$, we decompose the error $e:= \wh u-u_\grid \in \mathbb E_\grid$ into a conforming part $e^0 := \wh u-u_\grid^0 \in H^1_0(\Omega)$ and a non-conforming part $e^\perp:=-u_\grid^\perp \in \mathbb V_\grid^{-1}$ according to \eqref{e:conf-nonconf-space-decomp}.
The proof thus relies on techniques used in the conforming theory coupled with Lemma~\ref{l:control_perp} (control of non-conformity). We denote by $C$ a generic constant independent of the discretization and $\kappa$ but possibly depending on $\overline{\kappa}_{\textrm{stab}}$. 

Using the coercivity~\eqref{e:coerc_Bdg} and partial consistency \eqref{e:dg_partialconsistency}, we get
\begin{equation}\label{e:upper_bound_split}
C_\textrm{coer}\dGnorm{e}^2 \leq \bilin[\grid]{e}{e} = \bilin[\grid]{e}{e^0 - \mathcal I_\grid  e^0}-\bilin[\grid]{e}{u_\grid^\perp},
\end{equation}
where $I_\grid$ is the Scott-Zhang interpolant provided in Proposition~\ref{P:quasi-interp-bc}. 
For the first term, we note that since $e^0 - I_\grid e^0 \in H^1_0(\Omega)$ we have
$$
\bilin[\grid]{e}{e^0 - I_\grid e^0} = a_\grid[e,e^0-I_\grid e^0]- \int_\Omega \mathcal L_\grid[e] \cdot \bA \nabla (e^0 - I_\grid e^0).
$$
For the term involving the bilinear form $a_\grid$ we proceed as in the conforming case to arrive at
\begin{equation*}
    a_\grid[e,e^0-I_\grid e^0] \lesssim   \eta_\grid(u_\grid) \| \nabla e^0 \|_{L^2(\Omega)}.
\end{equation*}
This, combined with estimate \eqref{e:estimate_lift2} on the lifting operators and the $H^1$ stability of the Scott-Zhang interpolant, yield
$$
\bilin[\grid]{e}{e^0 - I_\grid e^0} \lesssim \left( \eta_\grid(u_\grid)^2+  \| h^{-1/2} \jump{u_\grid}\|^2_{L^2(\faceswbdy)} \right)^{1/2} \| \nabla e^0 \|_{L^2(\Omega)}.
$$
We rewrite $e^0 = e + u_\grid^\perp$ and use the estimate on the non-conforming component provided by Lemma~\ref{l:control_perp} along with a Young inequality to write
$$
\bilin[\grid]{e}{e^0 - I_\grid e^0} \leq C \left(\eta_\grid(u_\grid)^2+\kappa \| h^{-1/2} \jump{u_\grid}\|^2_{L^2(\faceswbdy)} \right) + \frac{C_\textrm{coer}}{4}   \| \nabla e \|_{L^2(\Omega)}^2.
$$
For the second term in \eqref{e:upper_bound_split}, the continuity \eqref{e:cont_Bdg} of the bilinear form $\mathcal B_\grid$, Lemma~\ref{l:control_perp} again, and a Young inequality yield
$$
\bilin[\grid]{e}{u_\grid^\perp} \leq \frac{C_\textrm{coer}}{4}\dGnorm{e}^2 + C \kappa  \| h^{-1/2} \jump{u_\grid}\|^2_{L^2(\faceswbdy)}.
$$
Returning to \eqref{e:upper_bound_split}, we find that
$$
\| e \|_\grid^2 \lesssim \eta_\grid(u_\grid)^2  + \kappa \| h^{-1/2} \jump{u_\grid}\|^2_{L^2(\faceswbdy)},
$$
which is the desired upperbound.

We finally deal with the lower bound \eqref{e:dg_lower}. For $T\in \grid$ and $v \in H^1_0(T)$, we get
$$
\int_T (-\textrm{div}_\grid(\wh{\bA}\nabla_\grid u_\grid) + \wh c u_\grid - \wh f) v = \int_T \nabla v \cdot \wh{\bA} \nabla_\grid (\wh u-u_\grid) v + \wh c (\wh u-u_\grid) v.
$$
For an interior face $F \in \faces$, $v \in H^1_0(\omega_F)$ with $\omega_F:= \{ T \in \grid \ : \ T\cap F \not = \emptyset \}$, we have \looseness=-1
\begin{equation*}
    \begin{split}
\int_F (\jumpinline{\wh{\bA}\nabla_\grid u_\grid}-\wh f) v  = & \int_{\omega_F}  (-\textrm{div}_\grid(\wh{\bA}\nabla_\grid u_\grid) + \wh c u_\grid -\wh f) v \\
&- \int_{\omega_F} \nabla v \cdot \wh{\bA}\nabla_\grid (\wh u-u_\grid)  - \wh c (\wh u-u_\grid) v.
\end{split}
\end{equation*}
The desired lower bound follows from the same arguments as in the conforming case;
we refer to Proposition~\ref{P:lw-bd-with-std-res-est} (partial lower bound).
\end{proof}

Upper bound \eqref{e:dg_upper} may suggest adding
the jump term $\kappa^{1/2}  \|  h^{-1/2}\jump{u_\grid}\|_{L^2(\faceswbdy)}$
to the residual estimator $\eta_\grid(u_\grid)$.
This would result in a clean upper bound but, because of the presence of negative powers of the meshsize, would be at the expense of destroying the monotonicity property of the estimator; see e.g. Proposition~\ref{P:est-reduction} (estimator reduction). The later is instrumental in the analysis provided below.

The next result mitigates the effect of the additional jump term by showing that $\|h^{-1/2}\jump{u_\grid}\|_{L^2(\faceswbdy)}$ can be bounded by $\eta_\grid(u_\grid)/\kappa$ and thus can be absorbed by the estimator in the upper bound provided $\kappa$ is sufficiently large.
We follow the proof provided in \cite{bonito2010quasi} and refer to \cite[eq. (3.20)]{karakashian2007convergence} for an alternative (original) proof.

\begin{lemma}[discontinuity control]\label{l:jump}
Let $\grid \in \grids^\Lambda$ be a $\Lambda-$admissible refinement of $\grid_0$ satisfying Assumption~\ref{A:initial-labeling} (initial labeling). Assume that $\wh{\bA}$ and $\wh c$ satisfy the structural assumption \eqref{e:structural-IP}.
There exists a constant $\overline{\kappa}_{\textrm{jump}} \geq \overline{\kappa}_{\textrm{stab}} > 0$ such that
if  $\kappa\geq \overline{\kappa}_{\textrm{jump}}$ then there holds
$$
\| h^{-1/2}\jump{u_\grid}\|_{L^2(\faceswbdy)} \lesssim \kappa^{-1} \eta_\grid(u_\grid).
$$
\end{lemma}
\begin{proof}
For $v_\grid^0 \in \mathbb V^0_\grid$ we realize that because $\jumpinline{v_\grid^0}=0$, the coercivity estimate \eqref{e:coerc_Bdg} implies that
\begin{equation}\label{e:dg_jumpcoerc}
C_\textrm{coerc}\kappa \| h^{-1/2}\jump{u_\grid}\|_{L^2(\faceswbdy)}^2 \leq \bilin[\grid]{u_\grid - v_\grid^0}{u_\grid - v_\grid^0}.
\end{equation}
We now rewrite the right-hand side of \eqref{e:dg_jumpcoerc} to produce residual terms.
Since $u_\grid$ solves \eqref{e:IP}, we have
\begin{equation}\label{e:dg_jump_A}
\bilin[\grid]{u_\grid-v_\grid^0}{u_\grid-v_\grid^0}= \langle \wh f, u_\grid - v_\grid^0\rangle_\grid -\bilin[\grid]{v_\grid^0}{u_\grid-v_\grid^0}.
\end{equation}
We concentrate for the moment on the second term.
Since $\jumpinline{v_\grid^0}=0$, the stabilization term vanishes as well:
$$
\kappa \sum_{F\in \faceswbdy} \int_F \jump{v_\grid^0} \jump{u_\grid - v_\grid^0} = 0.
$$
Hence, writing $v_\grid^0 = u_\grid + (v_\grid^0-u_\grid)$, we deduce that
\begin{equation*}
\begin{split}
    \bilin[\grid]{v_\grid^0}{u_\grid-v_\grid^0} &= a_\grid[u_\grid,u_\grid-v_\grid^0] 
    \\ & \quad - a_\grid[u_\grid-v_\grid^0,u_\grid-v_\grid^0]     
    - \int_\Omega \nabla v_\grid^0 \cdot \wh{\bA} \mathcal L_\grid[u_\grid] ,
\end{split}
\end{equation*}
where we invoked again the property $\jumpinline{v_\grid^0}=0$ to infer that $\mathcal L_\grid[v_\grid^0]=0$. Integrating by parts the first term on the right-hand side, adding it to
the first term on the right-hand side of \eqref{e:dg_jump_A}, and using the extended definition
\eqref{e:ext_f-dG} of the duality pairing leads to the following expression involving the residuals
$r_\grid(u_\grid), j_\grid(u_\grid)$:
\begin{equation*}
\begin{split}
    \bilin[\grid]{u_\grid-v_\grid^0}{u_\grid-v_\grid^0} =  &  \int_\Omega r_\grid(u_\grid) (u_\grid - v_\grid^0) + \int_{\faceswbdy} \big(j_\grid(u_\grid)-\wh f\big) \meaninline{u_\grid-v_\grid^0}\\
    &+ a_\grid[u_\grid-v_\grid^0,u_\grid-v_\grid^0]- \int_\Omega \nabla (u_\grid-v_\grid^0) \cdot \wh{\bA} \mathcal L_\grid[u_\grid].
\end{split}
\end{equation*}
We point out that we have also employed the definition \eqref{d:lift} of lift to rewrite the resulting 
face terms.
Inserting this estimate back in \eqref{e:dg_jumpcoerc}, together with the bound \eqref{e:estimate_lift2} for lifts and the continuity estimate \eqref{e:dG_cont_a} of $a_\grid$,  gives
\begin{equation*}
\begin{split}
 \kappa\| h^{-1/2} & \jump{u_\grid} \|_{L^2(\faceswbdy)}^2 \lesssim 
\dGnorm[1]{u_\grid - v_\grid^0}^2 \\
& + \eta_\grid(u_\grid) \left(\| h^{-1} (u_\grid-v_\grid^0)\|_{L^2(\Omega)} + \|h^{-1/2}\meaninline{u_\grid-v_\grid^0}\|_{L^2(\faceswbdy)}\right).
\end{split}
\end{equation*}
Note that the presence of $\dGnorm[1]{\,.\,}$ rather than $\dGnorm[\kappa]{\,.\,}$ on the right-hand side of the above estimate is critical for the argument below.
The former is independent of $\kappa$ and can thus be absorbed on the left-hand side for sufficiently large $\kappa$ provided $v_\grid^0=\IdG  u_\grid$.
In fact, the interpolation estimates \eqref{e:dg_interpol_nonconfC} in turn imply
\begin{equation*}
    \begin{split}
\| h^{-1} (u_\grid-\IdG  u_\grid)\|_{L^2(\Omega)} &+ \|h^{-1/2}\meaninline{u_\grid-\IdG  u_\grid}\|_{L^2(\faceswbdy)} 
\\ & + \dGnorm[1]{u_\grid - \IdG  u_\grid} 
\lesssim \| h^{-1/2} \jump{u_\grid} \|_{L^2(\faceswbdy)}.        
    \end{split}
\end{equation*}
Hence, a Young's inequality yields
$$
(\kappa-\overline{\kappa}_{\textrm{stab}}) \| h^{-1/2} \jump{u_\grid} \|_{L^2(\faceswbdy)}^2 \lesssim \kappa^{-1} \eta_\grid(u_\grid)^2  + \| h^{-1/2} \jump{u_\grid} \|_{L^2(\faceswbdy)}^2,
$$
and the desired estimate follows provided $\kappa$ is sufficiently large.
\end{proof}

As a direct consequence of the previous Lemma, 
we obtain a simpler practical upper bound.
\begin{corollary}[stabilization-free a posteriori upper bound]
\label{c:dg_upper2}
Under the assumptions of Lemma~\ref{l:jump}, there exists a constant $C_\textrm{U}$ such that for all $\kappa \geq \overline{\kappa}_{\textrm{jump}}$ we have
\begin{equation}\label{e:dg_upper2}
\dGnorm{\wh u - u_\grid} \leq C_\textrm{U} \, \eta_\grid(u_\grid).
\end{equation}
\end{corollary}
\begin{proof}
Combine the upper bound \eqref{e:dg_upper} and Lemma~\ref{l:jump}.
\end{proof}

The partial consistency \eqref{e:dg_partialconsistency} leads to {\it partial Galerkin orthogonality}
\begin{equation}\label{E:partial-galerkin-ortho}
\B_\grid[\wh{u}-u_\grid,v]=0 \quad\forall v\in\V_\grid^0.
\end{equation}
This would suggest that a quasi-best approximation (Cea's lemma) result in the full space
$\V_\grid^{-1}$ is questionable. 
However, the lack of consistency is built into the jump terms which are in turn controlled by the estimator weighted by a negative power of the penalty parameter $\kappa$. It thus remains to resort to the lower bound to return to the error and derive a quasi-best approximation estimate for 
sufficienlty large $\kappa$. We prove this result next, which expresses the important fact that
dG is quasi-optimal with respect to the norm $\dGnorm[\kappa][\grid]{\,\cdot \,}$ defined in 
\eqref{E:dG-norm}. This has two significant consequences: first it leads to quasi-monotonicity of
the error upon refinement (see Corollary \ref{c:quasi-monotonicity-dG} below), and second it dictates
the approximation class for dG already alluded to in Proposition \ref{p:equiv_classes_dG}
(equivalence of classes for $u$).

\begin{corollary}[Cea's lemma]\label{c:cea-dG}
Under the assumptions of Lemma~\ref{l:jump}, there is $\overline{\kappa}_\textrm{Cea} \geq \overline{\kappa}_\textrm{jump}$ such that for $\kappa \geq \overline{\kappa}_{\textrm{Cea}}$, there holds
\begin{equation}\label{e:cea-dG}
\dGnorm[\kappa][\grid]{\wh u-u_{\grid}} \leq C_\textrm{Cea} \inf_{v_\grid \in \mathbb V_\grid^{-1}} \dGnorm[\kappa][\grid]{\wh u-v_{\grid}}.
\end{equation}
\end{corollary}
\begin{proof}
    We combine the orthogonal decomposition \eqref{e:conf-nonconf-space-decomp} and the partial Galerkin orthogonality \eqref{E:partial-galerkin-ortho} to write
\begin{equation*}
    \bilin[\grid]{\wh u - u_\grid}{\wh u - u_\grid} =  \bilin[\grid]{\wh u - u_\grid}{\wh u - v_\grid} -     \bilin[\grid]{\wh u - u_\grid}{u_\grid^\perp} +   \bilin[\grid]{\wh u - u_\grid}{v_\grid^\perp}
\end{equation*}
    for all $v_\grid \in \mathbb V_\grid$. 
    Invoking the coercivity and continuity of $\B_\grid$ in Lemma~\ref{c:dG_cont_coer}) (properties of
    $\B_\grid$) in conjunction with Lemma~\ref{l:control_perp} (control of discontinuity) yields
    \begin{align*}
      \dGnorm[\kappa][\grid]{\wh u - u_\grid} \lesssim \dGnorm[\kappa][\grid]{\wh u - v_\grid} +\kappa^{1/2} \| h^{-1/2} \jump{u_\grid} \|_{L^2(\faceswbdy)} + \kappa^{1/2} \| h^{-1/2} \jump{v_\grid} \|_{L^2(\faceswbdy)}.
    \end{align*}
    Applying now Lemma~\ref{l:jump} (discontinuity control) results in
    \[
    \| h^{-1/2} \jump{u_\grid} \|_{L^2(\faceswbdy)} \lesssim \kappa^{-1} \eta_\grid(u_\grid)
    \lesssim \kappa^{-1} \dGnorm[\kappa][\grid]{\wh u - u_\grid}
    \]
    because of the lower bound \eqref{e:dg_lower}. We thus end up with
    \[
    \dGnorm[\kappa][\grid]{\wh u - u_\grid} \lesssim \dGnorm[\kappa][\grid]{\wh u - v_\grid}
    + \kappa^{-1/2} \dGnorm[\kappa][\grid]{\wh u - u_\grid},
    \]
    which for $\kappa$ sufficiently large gives the desired bound.
\end{proof}

With this best approximation result, we deduce the following crucial property.

\begin{corollary}[quasi-monotonicity]\label{c:quasi-monotonicity-dG}
Under the assumptions of Lemma~\ref{l:jump}, there is a constant $C_\textrm{Mo}$ independent of the discretization parameters and $\kappa$ such that for all $\kappa \geq \overline{\kappa}_{\textrm{jump}}$ and all $\grid_* \geq \grid$ we have
\begin{equation}\label{e:quasi-monotone-dG}
\dGnorm[\kappa][\grid_*]{\wh u-u_{\grid_*}} \leq C_\textrm{Mo} \dGnorm[\kappa][\grid]{\wh u-u_{\grid}}.
\end{equation}
\end{corollary}
\begin{proof}
We rely on the orthogonal decomposition \eqref{e:conf-nonconf-space-decomp} to write $u_\grid = u_\grid^0 + u_\grid^\perp$ and on Corollary~\ref{c:cea-dG} (Cea's lemma). Since 
$u_\grid^0 \in \mathbb V_\grid^0 \subset \mathbb V_{\grid_*}^{-1}$, we see that
$$
C_\textrm{Cea}^{-1} \dGnorm[\kappa][\grid_*]{\wh u-u_{\grid_*}} \le \dGnorm[\kappa][\grid_*]{\wh u - u_\grid^0}
= | \wh u - u_\grid^0 |_{H^1_0(\Omega)} = \dGnorm[\kappa][\grid]{\wh u - u_\grid^0}.
$$
Therefore, adding and subtracting $u_\grid^\perp$ and making use of Lemma~\ref{l:control_perp} (control of non-conformity) together with Lemma~\ref{l:jump} (discontinuity control) implies
$$
\dGnorm[\kappa][\grid_*]{\wh u-u_{\grid_*}} \lesssim \dGnorm[\kappa][\grid]{\wh u - u_\grid} + \kappa^{-1/2} \eta_{\grid}(u_\grid).
$$
It remains to invoke the lower bound \eqref{e:dg_lower} to deduce the desired result.
\end{proof}

Corollary \ref{c:quasi-monotonicity-dG} assumes the same data.
In estimating the cost of \GALERKINIP we need a variant of this result that allows for different data. We establish this next.

\begin{corollary}[quasi-monotonicity with different data]\label{c:quasi-monotonicity-dG-data}
Let $\grid_*\ge\grid$ and $\wh{\data}_*, \wh{\data}$ be discrete data on these meshes. Let
$\wh{u}_* = u(\wh{\data}_*), \wh{u}=u(\wh{\data}) \in H^1_0(\Omega)$ and 
$u_{\grid_*}\in\V_{\grid_*}^{-1}, u_\grid \in \V_\grid^{-1}$ be the corresponding exact and 
Galerkin solutions.
Under the assumptions of Lemma~\ref{l:jump}, for $\kappa \geq \overline{\kappa}_{\textrm{jump}}$ there holds
\begin{equation}\label{e:quasi-monotone-dG-data}
\dGnorm[\kappa][\grid_*]{\wh{u}_* - u_{\grid_*}} \leq C_\textrm{Mo} \big(\dGnorm[\kappa][\grid]{\wh u-u_{\grid}} + | \wh{u}_* - \wh u |_{H^1_0(\Omega)} \big).
\end{equation}
\end{corollary}
\begin{proof}
   We proceed as in the proof of Corollary \ref{c:quasi-monotonicity-dG} with $\wh{u}_*$,
   but in the last step use the fact $u_\grid$ and $\wh{u}$ are the functions associated with the
   same data $\wh{\data}$ and thereby satisfy $C_L \eta_\grid(u_\grid) \le \|\wh{u} - u_\grid\|_{\kappa,\grid}$ according to \eqref{e:dg_lower}. Applying the triangle inequality and the 
   property $\dGnorm[\kappa][\grid]{\wh{u}_* - \wh{u}} = |\wh{u}_* - \wh{u}|_{H^1_0(\Omega)}$
   concludes the proof.
\end{proof}

We end this section with the dG counterpart of Theorem~\ref{T:ubd-corr} (upper bound for corrections).
One striking difference is that the lack of consistency prevents the discrete lower bound in the dG context from localizing to the refined set $\grid \setminus \grid_*$ for $\grid_* \geq \grid$. Rather, it contains a global jump term that expresses the lack of conformity and vanishes as $\kappa \to \infty$ in view of Lemma \ref{l:jump} (discontinuity control). This is consistent with the upper bound \eqref{e:dg_upper}.
We use the notation $\omega_\grid(\tau)$ for a set of elements $\tau \in \grid$ to denote $\tau$ augmented by one layer of elements
$$
\omega(\tau):=\omega_\grid(\tau)= \bigcup_{T \in \tau} \omega_\grid(T).
$$

\begin{lemma}[quasi-localized discrete upper bound]\label{l:upper-discrete-dG}
Let $\grid, \grid_* \in \grids^\Lambda$, with $\grid_* \geq \grid$, be two $\Lambda-$admissible refinements of $\grid_0$ satisfying Assumption~\ref{A:initial-labeling} (initial labeling). Assume that $\wh{\bA}$ and $\wh c$ satisfy the structural assumption \eqref{e:structural-IP}, $\wh{f}\in\F_\grid$, and denote by $u_\grid \in \mathbb V_{\grid}^{-1}$, $u_{\grid_*} \in \mathbb V_{\grid_*}^{-1}$ the two Galerkin solutions associated with $\grid$, $\grid_*$ respectively. 
There is a constant $C_\textrm{LU}$ such that for all $\kappa > \overline{\kappa}_{\textrm{stab}}$ we have
$$
\dGnorm[\kappa][\grid]{u_{\grid_*}^0-u_\grid}^2 \leq C_{\textrm{LU}}^2 \left(\eta_\grid^2\big(u_\grid,\omega(\grid \setminus \grid_*)\big) + \kappa \| h^{-1/2} \jump{u_\grid}\|_{L^2(\faceswbdy)}^2\right),
$$
where $u_{\grid_*}=u_{\grid_*}^0+u_{\grid_*}^\perp$ is the orthogonal decomposition according to \eqref{e:conf-nonconf-space-decomp}.
Moreover, if $\kappa \geq \overline{\kappa}_{\textrm{jump}}$, there holds
\begin{equation}\label{e:upper-discrete-dG}
\dGnorm[\kappa][\grid]{u_{\grid_*}^0-u_\grid}^2 \leq C_{\textrm{LU}}^2 \left(\eta_\grid^2\big(u_\grid,\omega(\grid \setminus \grid_*)\big) + \kappa^{-1} \eta_\grid^2(u_\grid) \right).
\end{equation}
\end{lemma}
\begin{proof}
    We decompose $u_{\grid_*} = u_{\grid_*}^0+u_{\grid_*}^\perp$ according to \eqref{e:conf-nonconf-space-decomp}, exploit the partial consistency \eqref{e:dg_partialconsistency} for $u_{\grid_*}^0$
    with $v^0 \in \V_\grid^0$, and $\wh{f} \in \F_\grid\subset\F_{\grid_*}$ to obtain
\begin{equation*}
    \B_{\grid_*}[u_{\grid_*},v^0] = a_{\grid_*} [u_{\grid_*},v^0] =
    a_\grid [u_{\grid_*},v^0] = \B_\grid [u_{\grid_*},v^0] = \langle \wh{f},v^0 \rangle_\grid.
\end{equation*}
Since $\B_\grid[u_\grid,v^0]=\langle \wh{f},v^0 \rangle_\grid$, we readily see that
\[
\bilin[\grid]{u_{\grid_*}^0-u_\grid}{v^0}=0 \qquad \forall v^0 \in \V_\grid^0.
\]
We rely on this reduced form of Galerkin orthogonality to prove the assertions. To this end,
we write $u_{\grid_*}^0-u_\grid= e_*^0-u_\grid^\perp$, with $e_*^0:= u_{\grid_*}^0 - u_\grid^0$.
Using the coercivity estimate \eqref{e:coerc_Bdg} for $v=u_{\grid_*}^0-u_\grid\in\mathbb{E}_\grid$ yields for $\kappa \geq \overline{\kappa}_{\textrm{stab}}$ 
    \begin{equation}\label{e:local_upper-dG-tmp}
    \begin{split}
    \dGnorm[\kappa][\grid]{u_{\grid_*}^0-u_\grid}^2 & \lesssim \bilin[\grid]{u_{\grid_*}^0-u_\grid}{u_{\grid_*}^0-u_\grid}\\
    & =  \bilin[\grid]{u_{\grid_*}^0-u_\grid}{e_*^0} - \bilin[\grid]{u_{\grid_*}^0-u_\grid}{u_\grid^\perp}.
    \end{split}
    \end{equation}
    Note that the last term cannot be localized and accounts for the lack of consistency of the dG method. However, it can be made arbitrarily small by increasing the penalty parameter $\kappa$. In fact, combining the continuity \eqref{e:cont_Bdg} with Lemma~\ref{l:control_perp} (control of non-conformity) gives
$$
\bilin[\grid]{u_{\grid_*}^0-u_\grid}{u_\grid^\perp} \lesssim \kappa^{1/2} \| h^{-1/2} \jump{u_\grid}\|_{L^2(\faceswbdy)} \dGnorm[\kappa][\grid]{u_{\grid_*}^0-u_\grid}.
$$

    To localize $\bilin[\grid]{u_{\grid_*}^0-u_\grid}{e_*^0}$, we choose $v^0=\IdG e_*^0$, where the interpolation operator $\IdG$ is given by \eqref{e:dg_int}, and exploit the reduced Galerkin orthogonality. Since $e_*^0 - \IdG e_*^0 \in H^1_0(\Omega)$, the decomposition \eqref{e:bilinIP2} of the bilinear form $\mathcal B_\grid$ reads
    \begin{align*}
    \bilin[\grid]{u_{\grid_*}^0 &-u_\grid}{e_*^0} = \bilin[\grid]{u_{\grid_*}^0-u_\grid}{e_*^0-\IdG e_*^0}
    \\
    & = a_\grid[u_{\grid_*}^0-u_\grid,e_*^0-\IdG e_*^0] - \int_\Omega \mathcal L_\grid[u_{\grid_*}^0-u_\grid] \cdot \wh{\bA}\nabla_\grid (e_*^0-\IdG e_*^0).
    \end{align*}
    We handle the first term as in the conforming case (Theorem~\ref{T:ubd-corr}), namely
    $$
    a_{\grid}(u_{\grid_*}^0-u_\grid,e_*^0) \lesssim \eta_\grid \big(u_\grid,\omega(\grid\setminus \grid_*)\big) \, | u_{\grid_*}^0 - u_\grid^0 |_{H^1_0(\Omega)}.
    $$
    Note that the interpolation estimate \eqref{e:dg_interpol_nonconfA} for $\IdG$ is responsible
    for the appearance of $\omega(\grid\setminus \grid_*)$ rather than the smaller set $\grid\setminus \grid_*$.
    
    For the second term, we use the lift estimate \eqref{e:estimate_lift2} along with the $H^1$-stability of $\IdG$ and $\jumpinline{u_{\grid_*}^0}=0$ to write
    $$
    \int_\Omega \mathcal L_\grid[u_{\grid_*}^0-u_\grid] \cdot \wh{\bA}\nabla_\grid e_*^0 \lesssim  \| h^{-1/2} \jump{u_\grid}\|_{L^2(\faceswbdy)} \, | u_{\grid_*}^0 -u_\grid^0|_{H^1_0(\Omega)}.
    $$
    Inserting the estimates into \eqref{e:local_upper-dG-tmp}, and recalling that $1\lesssim \overline{\kappa}_{\textrm{stab}} \leq \kappa$, we find that
    \begin{align*}
    \dGnorm[\kappa][\grid]{u_{\grid_*}^0-u_\grid}^2 & \lesssim 
    \kappa^{1/2} \| h^{-1/2} \jump{u_\grid}\|_{L^2(\faceswbdy)} \dGnorm[\kappa][\grid]{u_{\grid_*}^0-u_\grid}
    \\
    & + \left(\eta_\grid (u_\grid,\grid\setminus \grid_*) + \kappa^{1/2} \| h^{1/2} \jump{u_\grid} \|_{L^2(\faceswbdy)} \right) | u_{\grid_*}^0 -u_\grid^0|_{H^1_0(\Omega)}.
    \end{align*}
    Notice that $u_{\grid_*}^0 -u_\grid^0 = u_{\grid_*}^0 -u_\grid + u_\grid^\perp$ so that in view of Lemma~\ref{l:control_perp} (control of non-conformity), we have
    $$
    | u_{\grid_*}^0 -u_\grid^0|_{H^1_0(\Omega)} = \dGnorm[\kappa][\grid]{u_{\grid_*}^0 -u_\grid^0}\lesssim  \dGnorm[\kappa][\grid]{u_{\grid_*}^0 -u_\grid}  +\kappa^{1/2} \| h^{-1/2} \jump{u_\grid}\|_{L^2(\faceswbdy)}.
    $$
    The first desired inequality follows from the last two estimates. For the second inequality, it suffices to further invoke Lemma~\ref{l:jump} (discontinuity control).
\end{proof}

%--------------------------------------------------------------------------------
\subsection{Module \GALERKINIP}\label{ss:GalerkinIP}
%--------------------------------------------------------------------------------
%
The main ingredients for the a posteriori estimation have been derived in the previous section and we can now turn our attention to the adaptive method. In essence, it is the same as in the conforming case (Algorithm~\ref{A:GALERKIN}) but accounting for the perturbation arising from the non-conforming setting. Compared to Algorithm~\ref{A:GALERKIN}, $\SOLVE \, (\grid)$ determines the discontinuous Galerkin solution to \eqref{e:IP} and $\REFINE(\grid,\cal M)$ produce the smallest $\Lambda$-admissible refinement of $\grid$ where all the marked elements $\cal M$ are refined at least $b\geq 1$ times.

\begin{algo}[$\GALERKINIP$]\label{A:GALERKIN-IP}
\index{Algorithms!\GALERKINIP: discontinuous Galerkin version of \GALERKIN}
Let $\wh\grid \ge \grid_0$ be a $\Lambda$-admissible refinement, $\Lambda \geq 0$, of a suitable initial mesh $\grid_0$.
Let data $\wh \data=(\wh{\bA},\wh c,\wh f) \in \mathbb D_{\wh{\grid}}$ be discrete on $\wh \grid$ and  $\varepsilon>0$ be
a stopping tolerance. The following routine creates a $\Lambda$-admissible refinement
$\grid\ge \wh {\grid}$ and discontinuous Galerkin solution $u_\grid\in\V_\grid^{-1}$ for data $\wh{\data}$
such that $\eta_{\grid}(u_\grid)\le \varepsilon$.
\begin{algotab}
  \> $[\grid,u_\grid] = \GALERKINIP \, (\wh \grid, \wh \data, \eps)$ \\
  \> \> $\text{set } j=0, \grid_0 = \wh{\grid} \text{ and do}$ \\
  \> \> \>$[u_j] = \SOLVE \, (\grid_j)$;\\
  \> \> \>$[\{\eta_j(u_j,T)\}_{T\in\grid_j}] = \ESTIMATE(u_j,\grid_j,\wh \data)$;\\
  \> \> \>$\text{if } \eta_j(u_j) \le \eps;$\\
  \> \> \> \>$\text{ return } (\grid_j,u_j)$ \\
  \> \> \>$[\marked_j] = \MARK \, \big(\{\eta_j(u_j,T)\}_{T\in\grid_j},\grid_j,\theta\big)$;\\
  \> \> \>$[\grid_{j+1}] = \REFINE(\grid_j,\marked_j)$; \\
  \> \> \>$j\leftarrow j+1$; \\
  \> \> $\text{while true}$
\end{algotab}
\end{algo}

We start the analysis of \GALERKINIP by investigating how the energy norm $\bilin[\grid]{v}{v}^{1/2}$ changes upon refining $\grid$. Note that in the conforming case, Lemma~\ref{L:Pythagoras} (Pythagoras)
directly provides the relation $\enorm{\wh u-u_{\grid_*}}\le\enorm{\wh u-u_\grid}$.
In the non-conforming setting, the constant on the right-hand side is no longer 1 and jumps terms are present in the estimate.
Regardless, it is possible to assess the effect of refinement in the energy norm and, in turn, compare two consecutive Galerkin solutions $u_\grid$ and $u_{\grid_*}$ where $\grid \in \grids^\Lambda$ and $\grid_* = \REFINE(\grid,\marked)$ for some $\marked \subset \grid$. 
This is the subject of the next three results but before embarking on this path, we mention a key ingredient for this comparison to hold: {\it the routine {\rm \REFINE} does not refine elements in $\grid$ more than $d$ times for $b=1$}; see Corollary~\ref{C:Lambda-refined}. This implies for any $F \in \faceswbdy$ and $F_* \in \faceswbdy_*$ with $F_* \subset F$, one has
\begin{equation}\label{e:comparison_faces_dG}
h_{F} \lesssim h_{F_*}.
\end{equation}

\begin{lemma}[mesh perturbation]\label{l:coarse_to_fine}
Let $\grid \in \grids^\Lambda$ be a $\Lambda-$admissible refinement of $\grid_0$ satisfying Assumption~\ref{A:initial-labeling} (initial labeling), $\marked \subset \grid$, and $\grid_*=\REFINE(\grid,\marked)$. Assume that $\wh{\bA}$ and $\wh c$ satisfy the structural assumption \eqref{e:structural-IP}.  There is a constant $C$ such that for 
$0<\varepsilon<1$ and all $v \in \mathbb E_\grid$ there holds
\begin{equation}\label{e:coarse_to_fine}
\bilin[\grid_*]{v}{v}\leq  (1+\varepsilon)\bilin[\grid]{v}{v} 
+ C \varepsilon^{-1} \kappa \| h^{-1/2}\jump{v}\|_{L^2(\faceswbdy)}^2.
\end{equation}
\end{lemma}
\begin{proof}
Because $\nabla_{\grid_*} v = \nabla_\grid v$ when $v \in \mathbb E_\grid$, we directly deduce that
\begin{equation}\label{e:dg_comp_tmp}
    \begin{split}
\bilin[\grid_*]{v}{v} = & \bilin[\grid]{v}{v} + 2\int_\Omega (\mathcal L_{\grid}[v]-\mathcal L_{\grid_*}[v]) \cdot \wh{\bA} \nabla_\grid v \\
&+\kappa \| h_*^{-1/2} \jump{v}\|_{L^2(\faceswbdy_*)}^2-\kappa \| h^{-1/2} \jump{v}\|_{L^2(\faceswbdy)}^2,
    \end{split}
\end{equation}
where $h_*:= h_{\grid_*}$ denotes the mesh size of $\grid_*$.

Unlike the broken gradients, the lifting operators are affected by refinements. However, this effect is controlled by the scaled jumps as we show now. 
Using estimate \eqref{e:estimate_lift3} twice with $\epsilon = \frac{\varepsilon C_\textrm{coer}}{2}$ 
yields
\begin{equation*}
    \begin{split}
&2\int_\Omega (\mathcal L_{\grid}[v]-\mathcal L_{\grid_*}[v]) \cdot \wh{\bA} \nabla_\grid v
\\
& \qquad \leq  \varepsilon C_{\textrm{coer}} \|  \nabla_\grid v\|_{L^2(\grid)}^2+ \frac{2 C^2_{\textrm{lift}}}{\varepsilon C_\textrm{coer}} \Big(\| h_*^{-1/2} \jump{v}\|_{L^2(\faceswbdy_*)}^2+\| h^{-1/2} \jump{v}\|_{L^2(\faceswbdy)}^2\Big).
    \end{split}
\end{equation*}
Hence, the coercivity estimate \eqref{e:coerc_Bdg} gives 
\begin{equation*}
\begin{split}
&2\int_\Omega (\mathcal L_{\grid}[v]-\mathcal L_{\grid_*}[v]) \cdot \wh{\bA} \nabla_\grid v
\\
& \qquad \leq \epsilon \bilin[\grid]{v}{v} +  \frac{2 C^2_{\textrm{lift}}}{\varepsilon C_\textrm{coer}}  \left(\| h^{-1/2} \jump{v}\|_{L^2(\faceswbdy)}^2+\| h_*^{-1/2} \jump{v}\|_{L^2(\faceswbdy_*)}^2\right).
\end{split}
\end{equation*}
Inserting this back into \eqref{e:dg_comp_tmp}, and using the fact that the jumps of $v$
occur only on $\faceswbdy\subset\faceswbdy_*$, the meshsize relation \eqref{e:comparison_faces_dG} proves the desired estimate.
\end{proof}

\begin{lemma}[comparison of solutions]\label{L:sol_comp}
Let $\grid \in \grids^\Lambda$ be a $\Lambda-$admissible refinement of $\grid_0$ satisfying Assumption~\ref{A:initial-labeling} (initial labeling), $\marked \subset \grid$, and $\grid_*=\REFINE(\grid,\marked)$.
Assume that $\wh{\bA}$ and $\wh c$ satisfy the structural assumption \eqref{e:structural-IP}. Let $\wh u = u (\wh{\data}) \in H^1_0(\Omega)$ be the solution of the perturbed problem \eqref{E:perturbed-weak-form} with discrete data $\wh \data$ and denote by $u_\grid \in \mathbb V_\grid^{-1}$, $u_{\grid_*} \in \mathbb V_{\grid_*}^{-1}$ the Galerkin solutions associated to $\grid$, $\grid_*$ respectively with data $\wh{\data}$. Let $\overline{\kappa}_{\textrm{jump}}$ be as in Lemma~\ref{l:jump}. There exists a constant $C_\mathrm{comp}$ such that for all $\kappa \geq \overline{\kappa}_{\textrm{jump}}$ and all $0<\epsilon <1$ we have
\begin{equation*}
\begin{aligned}
\bilin[\grid_*]{\wh u-u_{\grid_*}}{\wh u-u_{\grid_*}}  
&\leq (1+\epsilon)\bilin[\grid]{\wh u-u_{\grid}}{\wh u-u_{\grid}}\\
& \hskip-1.5cm
-\frac{C_\mathrm{coer}}{2}\| \nabla_{\grid_*} (u_{\grid_*}-u_\grid)\|_{L^2(\grid_*)}^2
+ \frac{C_\mathrm{comp}}{\epsilon \kappa}\left(
\eta_{\grid}(u_\grid)^2+\eta_{\grid_*}(u_{\grid_*})^2\right).
\end{aligned}
\end{equation*}
\end{lemma}
\begin{proof}
We invoke the partial Galerkin orthogonality \eqref{E:partial-galerkin-ortho} of $\wh u-u_{\grid_*}$
upon testing with $v^0:=u_{\grid_*}^0 - u_{\grid}^0 \in \V_{\grid_*}^0$:
$$
\bilin[\grid_*]{\wh u-u_{\grid_*}}{\wh u-u_{\grid_*}} = \bilin[\grid_*]{\wh u-u_{\grid_*}+v^0}{\wh u-u_{\grid_*}+v^0} - \bilin[\grid_*]{v^0}{v^0}.
$$
Note that 
$$
\wh u-u_{\grid_*}+v^0 = \wh u-u_\grid + u_\grid^\perp - u_{\grid_*}^\perp
$$
and $\dGnormstar{v^0} = \| \nabla_{\grid_*} v^0 \|_{L^2(\grid_*)}$, which is critical for the argument below. Hence, from the coercivity and continuity of $\mathcal B_{\grid_*}$ (Corollary~\ref{c:dG_cont_coer}), we deduce that
\begin{equation*}
    \begin{split}
& \bilin[\grid_*]{\wh u-u_{\grid_*}}{\wh u-u_{\grid_*}}\\
 & \qquad \leq  \bilin[\grid_*]{\wh u-u_{\grid}}{\wh u-u_{\grid}} +  2C_\textrm{cont}^{1/2} \bilin[\grid_*]{\wh u-u_{\grid}}{\wh u-u_{\grid}}^{1/2} \dGnormstar{u_{\grid}^\perp - u_{\grid_*}^\perp}\\
 & \qquad \qquad + C_\textrm{cont} \dGnormstar{u_\grid^\perp - u_{\grid_*}^\perp}^2 - C_\textrm{coer} \| \nabla_{\grid_*} (u_{\grid_*}^0 - u_{\grid}^0) \|_{L^2(\grid_*)}^2.
    \end{split}
\end{equation*}
We now apply the reverse triangle inequality and Young's inequality
$$
\| \nabla_{\grid_*} (u_{\grid_*}^0 - u_{\grid}^0)\|_{L^2(\grid_*)}^2 \geq \frac  1 2 \| \nabla_{\grid_*} (u_{\grid_*} - u_{\grid})\|_{L^2(\grid_*)}^2 - \| \nabla_{\grid_*} (u_{\grid_*}^\perp - u_{\grid}^\perp)\|_{L^2(\grid_*)}^2
$$
to deduce that for any $0<\varepsilon<1$
\begin{equation*}
    \begin{split}
\bilin[\grid_*]{\wh u-u_{\grid_*}}{\wh u-u_{\grid_*}} \leq & (1+\varepsilon)\bilin[\grid_*]{\wh u-u_{\grid}}{\wh u-u_{\grid}} 
\\
& - \frac{C_\textrm{coer}}{2} \| \nabla_{\grid_*} (u_{\grid_*} - u_{\grid})\|_{L^2(\grid_*)}^2
+ \frac{C}{\varepsilon}  \dGnormstar{u_{\grid_*}^\perp - u_{\grid}^\perp}^2,
    \end{split}
\end{equation*}
where $C$ is for the remainder of this proof a constant independent of the discretization parameters and $\kappa$.

To bound the last term we recall Lemma~\ref{l:control_perp} (control of non-conformity)
\begin{equation*}
\dGnormstar{u_{\grid_*}^\perp - u_{\grid}^\perp} \lesssim \kappa^{1/2} \| h_*^{-1/2} \jumpinline{u_{\grid_*}}\|_{L^2(\faceswbdy_*)} + \kappa^{1/2} \| h_*^{-1/2} \jump{u_{\grid}}\|_{L^2(\faceswbdy_*)},
\end{equation*}
and notice that the last integral over $\faceswbdy_*$ has weights relative to the local meshsize of
$\grid_*\ge\grid$. Since for consecutive meshes the local meshsizes are comparable, according to
\eqref{e:comparison_faces_dG}, we can write $\| h^{-1/2} \jump{u_{\grid}}\|_{L^2(\faceswbdy)}$ instead.
Inserting these expressions in the preceding estimate, and using Lemma~\ref{l:coarse_to_fine} (mesh perturbation) to replace $\B_{\grid_*}$ with $\B_\grid$ on the right-hand side, yields 
\begin{equation*}
    \begin{split}
\bilin[\grid_*]{\wh u-u_{\grid_*}}{\wh u-u_{\grid_*}} \leq & (1+\varepsilon)\bilin[\grid]{\wh u-u_{\grid}}{\wh u-u_{\grid}} - \frac{C_\textrm{coer}}{2} \| \nabla_{\grid_*} (u_{\grid_*} - u_{\grid})\|_{L^2(\grid_*)}^2 \\
&+ C \frac{\kappa}{\varepsilon} \left( \| h_*^{-1/2} \jump{u_{\grid_*}}\|_{L^2(\faceswbdy_*)}^2 +\| h^{-1/2} \jump{u_{\grid}}\|_{L^2(\faceswbdy)}^2\right),
    \end{split}
\end{equation*}
where $2\eps$ has been relabeled $\eps$. Finally, to derive the desired estimate, it remains to invoke Lemma~\ref{l:jump} (discontinuity control).
\end{proof}

Combining Lemma \ref{L:sol_comp} (comparison of solutions) with Lemma \ref{L:aposteriori-dg} (a posteriori error estimate), we derive the following dG version of 
Lemma \ref{L:Pythagoras} (Pythagoras).

\begin{corollary}[quasi-orthogonality of dG errors]\label{C:quasi-ortho-dg}
If the assumptions of Lemma \ref{L:sol_comp} hold, then
for all $\kappa\ge \kappa_{\textrm{QO}} := \frac{Ccomp}{\eps^2 C_L}$ and $0<\eps\le\frac14$ there holds
\[
\dGnorm[\kappa][\grid_*]{\wh u-u_{\grid_*}}^2 \le (1+4\eps) \,\dGnorm[\kappa][\grid]{\wh u-u_{\grid}}^2
-\frac{C_{\textrm{coer}}}{2} \,\|\nabla_{\grid_*} (u_{\grid_*} - u_\grid) \|_{L^2(\grid_*)}^2.
\]
\end{corollary}
\begin{proof}
We make use of the lower bound \eqref{e:dg_lower}, and set $D:=\frac{C_\textrm{comp}}{\eps C_L}$,
to rewrite the estimate of Lemma \ref{L:sol_comp} as follows:
\[
\Big(1 - \frac{D}{\kappa} \Big)\dGnorm[\kappa][\grid_*]{\wh u-u_{\grid_*}}^2 \le
\Big(1+\eps+\frac{D}{\kappa} \Big) \dGnorm[\kappa][\grid]{\wh u-u_{\grid}}^2
- \frac{C_{\textrm{coer}}}{2} \|\nabla_{\grid_*} (u_{\grid_*} - u_\grid) \|_{L^2(\grid_*)}^2.
\]
For $\kappa\ge \kappa_{\textrm{QO}}:= \frac{D}{\eps}$ this
inequality implies
\[
\dGnorm[\kappa][\grid_*]{\wh u-u_{\grid_*}}^2 \le \frac{1+2\eps}{1-\eps} \dGnorm[\kappa][\grid]{\wh u-u_{\grid}}^2 - \frac{C_{\textrm{coer}}}{2(1-\eps)}
\|\nabla_{\grid_*} (u_{\grid_*} - u_\grid) \|_{L^2(\grid_*)}^2.
\]
It remains to realize that $\frac{1+2\eps}{1-\eps} \le 1 + 4\eps$ provided $\eps\le\frac14$.
\end{proof}

The last ingredient to prove convergence of \GALERKINIP is a dG version of 
Proposition \ref{P:est-reduction} (estimator reduction) \AB{with $f=f_* \in \mathbb F_{\grid}$}. It turns out that the same estimate and
proof are valid for dG except that the $H^1_0$-seminorm is to be replaced by the broken $H^1_0$-seminorm.
We thus state the result without proof.

\begin{proposition}[estimator reduction]\label{P:est-reduction-dG}
Given $\grid\in\grids^\Lambda$ and a subset $\marked\subset\grid$
of elements marked for refinement, let $\grid_*=\REFINE\big(\grid,\marked\big)$.
\rhn{If $f=P_\grid f\in\F_\grid$, then}
there is a constant $\Clip$ such that for all $v\in\V_\grid$, $v_*\in\V_{\grid_*}$ and any $\delta>0$
\begin{equation*}
      \eta_{\grid_*}(v_*,\grid_*)^2 \le
      (1+\delta)\big(\eta_\grid(v,\grid)^2 - \lambda\,\eta_\grid(v,\marked)^2\big)
      +
      (1+\delta^{-1})\,\Clip^2 \, \| \nabla_{\grid _*}(v_*-v)\|_{L^2(\grid_*)}^2.
\end{equation*}
\end{proposition}

We are now in a position to prove a contraction property between two consecutive iterations of the adaptive loop \GALERKINIP. 
\begin{theorem}[contraction property]\label{t:dG_conv}
 Let $\wh \grid$ be a $\Lambda-$admissible refinement of $\grid_0$ satisfying Assumption~\ref{A:initial-labeling} (initial labeling). Let $\wh \data \in \mathbb D_{\wh \grid}$ such that $\wh{\bA}$ and $\wh c$ satisfy the structural assumption \eqref{e:structural-IP}. Let $\theta\in (0,1]$ be the D\"orfler marking parameter used in the {\rm \MARK} module and $\{\grid_j,\V_j,u_j\}_{j=0}^J$
  be a sequence of conforming meshes, finite element spaces and
  discrete solutions $u_j = u_j(\wh\data)\in\V_j$ created within {\rm \GALERKINIP}.
  If $\wh u\in H^1_0(\Omega)$ is the exact solution of \eqref{E:perturbed-weak-form} with discrete data $\wh \data$,
  then there exist constants $\overline{\kappa}_{\textrm{conv}}\geq 0$, $\gamma>0$, and $0<\alpha<1$ independent of the discretization parameters and $\kappa$, such that
  for all $\kappa \geq \overline{\kappa}_{\textrm{conv}}$ and $0 \leq j < J$
  \begin{equation}
  \label{e:dG-contraction}
B_{j+1}^2
+\gamma \eta_{\grid_{j+1}}^2(u_{j+1})
\leq \alpha^2 \Bigl(B_j^2 +\gamma \eta_{\grid_j}^2(u_j)\Bigr),
  \end{equation}
  where $B_{j}:=\big(\bilin[\grid_{j}]{\wh u-u_{j}}{\wh u-u_{j}}\big)^{1/2}$ is the dG norm of
  $\wh u-u_{j}$.
\end{theorem}
\begin{proof}
    In essence, we proceed as in Theorem~\ref{T:contraction} (general contraction property) 
    for the conforming case but with minor changes that account for non-conformity. We only 
    explain the differences below.
    For $j\geq 0$, we shorten the notations and write $\eta_{j} := \eta_{\grid_j}(u_j)$
    and $E_j := \| \nabla_{\grid_{j+1}}(u_{j+1}-u_j) \|_{L^2(\grid_{j+1})}$.

Corollary \ref{C:quasi-ortho-dg} (quasi-orthogonality of dG errors) gives for any $0<\eps\le\frac14$
\[
B_{j+1}^2 \le \big(1+4\eps\big) B_j^2 - \frac{\Ccoer}{2} E_j^2.
\]
Combining Proposition~\ref{P:est-reduction-dG} (estimator reduction), written in terms of $\grid = \grid_j$, $\grid_*=\grid_{j+1}$, $v=u_j$ and $v_*=u_{j+1}$, with D\"orfler marking
$\eta_j(u_j,\mathcal M_j) \ge \theta \eta_j$ yields
$$
\eta_{j+1}^2 \leq \big(1+\delta\big)\big(1-\lambda\theta^2\big) \eta_j^2
+ (1+\delta^{-1})  \Clip^2 E_j^2
$$
for any $\delta>0$. We now multiply this inequality by $\gamma>0$ and add it to the previous
one with the following choice of parameters:
\[
\delta = \frac{1-\frac{\lambda\theta^2}{2}}{1-\lambda\theta^2} - 1,
\qquad
\gamma = \frac{\Ccoer}{2\Clip^2 (1+\delta^{-1})}.
\]
Consequently, the terms involving $E_j^2$ cancel out and we end up with
\begin{align*}
B_{j+1}^2 + \gamma \eta_{j+1}^2 &\le \big(1+4\eps\big) B_j^2 + \gamma \big(1+\delta\big) \big(1-\lambda\theta^2\big) \eta_j^2
\\
& = \Big(1+4\eps - \frac{\gamma\lambda\theta^2}{4} \Big) B_j^2
+ \gamma \Big( 1 - \frac{\lambda\theta^2}{4} \Big) \eta_j^2.
\end{align*}
We finally choose $\eps := \frac{\gamma\lambda\theta^2}{32}$ to obtain \eqref{e:dG-contraction}
with $\alpha^2 = \max\big\{ 1- \frac{\gamma\lambda\theta^2}{8}, 1 - \frac{\lambda\theta^2}{4} \big\}$,
and conclude the proof.
\end{proof}

\begin{corollary}[linear convergence]\label{c:convergence-IP}
    Under the assumptions of Theorem~\ref{t:dG_conv}, and if $0<\alpha<1$, $\gamma >0$, $\overline{\kappa}_{\textrm{conv}}>0$ are the constants in \eqref{e:dG-contraction}, then for all $\kappa \geq \overline{\kappa}_{\textrm{conv}}$ there holds
    $$
    \dGnorm[\kappa][\grid_k]{\wh u-u_k} \leq C_* \alpha^{k-j} \dGnorm[\kappa][\grid_j]{\wh u-u_j},
    $$
    for some constant $C_*$ independent of the discretization parameters and $\kappa$.
\end{corollary}
\begin{proof}
Let $e_k^2= \dGnorm[\kappa][\grid_k]{\wh u-u_k}^2$ and $B_k^2 = \bilin[\grid_k]{\wh u-u_k}{\wh u-u_k}^2$, and use the coercivity estimate \eqref{e:coerc_Bdg}, the contraction property \eqref{e:dG-contraction}, the continuity estimate \eqref{e:cont_Bdg} and the lower bound \eqref{e:dg_lower} to arrive at
\[
\Ccoer e_k^2 \le B_k^2 \le \alpha^{2(k-j)} \big( B_j^2 + \gamma \eta_j^2 \big)
\le \alpha^{2(k-j)} \Big( \Ccont + \frac{\gamma}{C_L^2} \Big) e_j^2.
\]
This is the desired estimate in disguised with $C_*:= \frac{1}{\Ccoer}\big(\Ccont+\frac{\gamma}{C_\textrm{L}^2}\big)^{1/2}$.
\end{proof}

We end the discussion on \GALERKINIP by deriving the optimality property of the D\"orlfer marking strategy. We mimic the proof of Lemma~\ref{L:dorfler-no-osc} (D\"orfler marking) but directly use the optimal parameter $\mu=\frac 1 2$ to simplify the argument. We refer to the discussion after Lemma~\ref{L:dorfler-no-osc} for the role of $\mu$ and its influence on $\theta_0$. 
Notice that $\theta_0$ depends on $\kappa^{-1}$ because of its appearance in the perturbed localized upper bound \eqref{e:upper-discrete-dG}. It plays a similar role to $\sigma$ in Assumption~\ref{A:PDE-dominates} (restriction on $\omega$) in the presence of oscillations (one step method with switch).

\begin{lemma}[D\"orfler marking]\label{L:dorfler-IP}
Let $\grid_* \geq \grid$ be two $\Lambda-$admissible refinements of $\grid_0$ satisfying Assumption~\ref{A:initial-labeling} (initial labeling). Let $\wh \data \in \mathbb D_{\wh \grid}$ such that $\wh{\bA}$ and $\wh c$ satisfy the structural assumption \eqref{e:structural-IP}. 
Let $u_{\grid}\in\V_{\grid}^{-1}$, $u_{\grid_*}\in\V_{\grid_*}^{-1}$ denote the Galerkin solutions associated with $\grid$, $\grid_*$ respectively and $\wh u \in H^1_0(\Omega)$ denotes the solution to \eqref{E:perturbed-weak-form} with discrete data $\wh \data$. 
Assume $\kappa > \overline{\kappa}_{D}:=\max(\overline{\kappa}_{\textrm{stab}},4\Clip^{2}C_\textrm{LU}^{2})$.
If 
\begin{equation}\label{E:mu-no-osc-IP}
\eta_{\grid_*}(u_{\grid_*}) \leq \frac{1}{2} \eta_\grid(u_\grid)
\end{equation}
then the refined set $\grid \setminus \grid_*$ satisfies the D\"orfler property 
\begin{equation}\label{E:dorfler-IP}
\eta_\grid\big(u_\grid,\omega(\grid \setminus \grid_*)\big)\ge \theta_0\,\eta_\grid(u_\grid),
\end{equation}
with $0<\theta_0^2:=\theta_0^2(\kappa):= \frac{1-4\Clip^2 C_\textrm{LU}^2\kappa^{-1}}{4\Clip^2 C_\textrm{LU}^2}<\frac{1}{4\Clip^2 C_\textrm{LU}^2}$. 
\end{lemma}
\begin{proof}
To relate $\eta_{\grid}$ with $\eta_{\grid_*}$, we invoke Proposition~\ref{P:est-reduction-dG}
(estimator reduction) with $\delta=1$, along with the localized upper bound \eqref{e:upper-discrete-dG}, to write
$$
\eta_{\grid}(u_\grid)^2 \leq 2 \eta_{\grid_*}(u_{\grid_*})^2 +  2  C_\textrm{Lip}^2 C_{\textrm{LU}}^2\left( \eta_{\grid}\big(u_\grid,\omega(\grid \setminus \grid_*)\big)^2 + \kappa^{-1} \eta_\grid(u_\grid)^2\right). 
$$
This, combined with \eqref{E:mu-no-osc-IP} yields
$$
\left(\frac 1 2 - 2C_\textrm{Lip}^2 C_{\textrm{LU}}^2\kappa^{-1} \right) \eta_{\grid}(u_\grid)^2 \leq 2 C_\textrm{Lip}^2 C_{\textrm{LU}}^2 \eta_{\grid}\big(u_\grid,\omega(\grid \setminus \grid_*)\big)^2
$$
for $\kappa \geq \overline{\kappa}_\textrm{D}$. 
This is the desired results in disguised. 
\end{proof}

%--------------------------------------------------------------------------------
\subsection{Convergence of {\rm \AFEMTSIP}}\label{ss:AFEMTSIP}
%--------------------------------------------------------------------------------
%
Algorithm~\ref{algo:AFEM-TS-IP} (\AFEMTSIP) relies on two modules: \GALERKINIP and \DATA.
We have analyzed the performance of \GALERKINIP in the previous section and showed in Section \ref{S:data-approx} that the output $[\wh \grid_k, \wh \data_k] = \DATA(\grid_k,\data,\omega \eps_k)$ satisfies

\begin{equation}\label{E:omega-epsk-IP}
\|\data - \wh{\data}_k \|_{D(\Omega)} \le \omega\eps_k
\quad\Rightarrow\quad
|u - \wh{u}_k|_{H^1_0(\Omega)} \le C_D \omega\eps_k.
\end{equation}
Recall that $u=u(\data), \wh{u}_k = u(\wh{\data}_k)\in H^1_0(\Omega)$ are the exact solutions to \eqref{weak-form} with exact data $\data$ and discrete data $\wh{\data}_k$, respectively.
We also recall that $\wh \data_k$ satisfies the structural assumption \eqref{e:structural-IP} uniformly in $k$ and thus $C_D$ does not depend on $k$. 

We start with a result guaranteeing that the cost of \GALERKINIP does not depend on the iteration 
counter $k$ within \AFEMTSIP.

\begin{lemma}[computational cost of \GALERKINIP]\label{c:cost_dg}
For any $\kappa \geq \overline{\kappa}_{\textrm{conv}}$ and any $k\in \mathbb N$, the number of sub-iterations $J_k$ inside a call of {\rm \GALERKINIP} at iteration $k$ of Algorithm~\ref{algo:AFEM-TS-IP} ({\rm \AFEMTSIP}) is bounded independently of $k$.
\end{lemma}
\begin{proof}
We proceed as in the proof of Proposition \ref{P:cost-galerkin} (computational cost of $\GALERKIN$)
for the conforming case, and focus on the essential differences. We fix the iteration counter $k\ge1$, recall that the output of the $(k-1)$-th loop of \AFEMTSIP is $[\grid_{k},u_k]=\GALERKINIP \, (\wh \grid_{k-1},\wh \data_{k-1},\eps_{k-1})$, and denote by $\grid_{k,j}$ and $\wh u_{k,j} \in \mathbb V^{-1}_{\grid_{k,j}}$ the $j$-th mesh and Galerkin solution to \eqref{e:IP} with data $\wh \data_{k}$ in the $k$-th loop of \AFEMTSIP. The exact solution to the perturbed problem \eqref{E:perturbed-weak-form} with discrete coefficient $\wh \data_{k}$ is $\wh u_k = u(\wh{\data}_k)$.

We recall that $\grid_{k,0} = \wh{\grid}_k$ is the mesh produced by $\DATA$, and assume that
$u_{k,0}\in\V_{k,0}$ satisfies $\eta_{\grid_{k,0}}(u_{k,0})>\eps_k$ because otherwise $J_k=0$ and there is nothing to prove. In view of Corollary~\ref{c:convergence-IP} (linear convergence), all we need to prove is that the error $\dGnorm[\kappa][\grid_{k,0}]{u_{k,0}-\wh u_k}$ entering \GALERKINIP is bounded
by $\eps_k$. We resort to Corollary~\ref{c:quasi-monotonicity-dG-data} (quasi-monotonicity with different data) to write
\[
\dGnorm[\kappa][\grid_{k,0}]{u_{k,0}-\wh u_k} \leq C_\textrm{Mo} \left(\dGnorm[\kappa][\grid_{k}]{u_{k}-\wh u_{k-1}} + | \wh u_k - \wh u_{k-1}|_{H^1_0(\Omega)}\right).
\]
The appearance of the last term is the only difference with respect to Proposition \ref{P:cost-galerkin}.
However, in view of property \eqref{E:omega-epsk-IP} of $\DATA$, we infer that
\[
|\wh u_k - \wh u_{k-1}|_{H^1_0(\Omega)} \le
|\wh u_k - u|_{H^1_0(\Omega)} + |\wh u_{k-1} - u|_{H^1_0(\Omega)} \le C_D \omega(\eps_k+\eps_{k-1})
= 3 C_D \omega\eps_k.
\]
Moreover, the stabilization-free upper bound \eqref{e:dg_upper2} implies
\[
\dGnorm[\kappa][\grid_{k}]{u_k-\wh u_{k-1}} \leq C_{\textrm{U}} \, \eta_{\grid_k}(u_k) \leq C_{\textrm{U}} \eps_{k-1} = 2 C_{\textrm{U}} \eps_k,
\]
which, combined with the lower bound \eqref{e:dg_lower}, further yields
\[
\dGnorm[\kappa][\grid_{k,0}]{u_{k,0}-\wh u_k}
\le \big( 3C_D\omega + 2 C_U \big) \eps_k = \Lambda \eps_k
\quad\Rightarrow\quad
\eta_{\grid_{k,0}} \le C_L^{-1} \Lambda \eps_k.
\]
This is the requisite estimate. In fact, recalling Corollary~\ref{c:convergence-IP},
we see that
\[
\eta_{\grid_{k,j}}(u_{k,j}) \leq C_{\textrm{L}}^{-1}  \dGnorm[\kappa][\grid_{k,j}]{ u_{k,j} - \wh u_k} \leq C_{\textrm{L}}^{-1} C_*  \alpha^j \dGnorm[\kappa][\grid_{k,0}]{ u_{k,0} - \wh u_k}
\le C_{\textrm{L}}^{-1} C_* \Lambda \eps_k \alpha^j.
\]
Since \GALERKINIP stops when $\eta_{\grid_{J_k,k}}(u_{J_k,k})\le\eps_k$, we finally
conclude as in the proof of Proposition \ref{P:cost-galerkin} that $J_k$ is independent of $k$.
\end{proof}

The proof of convergence of $\AFEMTSIP$ is identical to the proof of  Proposition~\ref{P:convergence-AFEM} upon replacing the semi-norm $|.|_{H^1_0(\Omega)}$ by the appropriate dG norm. It is therefore not repeated here.

\begin{proposition}[convergence of \AFEMTSIP]
For any $\kappa \geq \overline{\kappa}_{\textrm{conv}}$ and $k\ge0$, the $(k+1)$-th iteration
of {\rm \AFEMTSIP} terminates and requires a finite number of inner iterations of $\GALERKINIP$
independent of $k$.
Moreover, if $u\in H^1_0(\Omega)$ denotes the solution to \eqref{weak-form}, there exists a constant $C_*$ such that the output of $[\mesh_{k+1},u_{k+1}]=\GALERKINIP \, (\widehat{\mesh}_{k},\widehat{\data}_{k},\varepsilon_k)$ satisfies
$$
\dGnorm[\kappa][\grid_{k+1}]{{u} -u_{k+1}} \leq C_* \varepsilon_k, \qquad \forall k \geq 0.
$$
Therefore, {\rm \AFEMTSIP} stops after $K<2+\frac{\log\frac{\eps_0}{\tol}}{\log2}$ iterations and delivers
\begin{equation*}
  \dGnorm[\kappa][\grid_{K}]{{u} -u_{K+1}} \leq C_* \tol.
\end{equation*}   
\end{proposition}
%

%--------------------------------------------------------------------------------
\subsection{Rate-optimality of $\AFEMTSIP$}\label{S:rate-optimal-dg}
%--------------------------------------------------------------------------------
%
To derive rates of convergence for the discontinuous Galerkin method, we proceed similarly to
Section~\ref{S:conv-rates-coercive} for the conforming case . 
Recall that in the $k$th-step of Algorithm~\ref{algo:AFEM-TS-IP} ($\AFEMTSIP$), the output of $[\wh{\grid}_k,\wh{\data}_k] = \DATA \, (\grid_k,\data,\omega\eps_k)$ is fed to $[\grid_{k+1},u_{k+1}]=\GALERKINIP \, (\wh{\grid}_k,\wh{\data}_k,\eps_k)$, which in turn iterates $J_k$ times. Lemma~\ref{c:cost_dg} shows that $J_k$ is uniformly bounded in $k$,
and we assume that $J_k \geq 1$ for otherwise the module $\GALERKINIP$ is skipped altogether. 
We denote by $(\grid_{k,j}, \marked_{k,j}, u_{\grid_{k,j}})$ the triplets of grids, marked sets and discrete solutions computed within $\GALERKINIP \, (\wh{\grid}_k,\wh{\data}_k,\eps_k)$ for $0\le j < J_k$. 
Note that 
\[
\wh{\eps}_{k,j} := \eta_{\grid_{k,j}} (u_{\grid_{k,j}},\wh{\data}_k) > \eps_k, \qquad 0\leq j <J_k
\]
so that together with the lower a posteriori error estimate \eqref{e:dg_lower}, we infer that
\[
\dGnorm[\kappa][\grid_{k,j}]{\wh{u}_k - u_{\grid_{k,j}}} \ge C_\textrm{L} \wh{\eps}_{k,j} > C_\textrm{L}\eps_k,
\]
where $\wh{u}_k = u (\wh{\data}_k)\in H^1_0(\Omega)$ is the exact solution with approximate data $\wh{\data}_k$. 

The module $\DATA$ guarantees \eqref{E:omega-epsk-IP}, and the parameter $\omega$ modulates the discrepancy between $u$ and $\wh{u}_k$ relative to $\eps_k$. The error due to data approximation can be made small relative to the finite element approximation by choosing $\omega$ much smaller than $1$. 
In addition, we have established Lemma~\ref{L:dorfler-IP} (D\"orfler marking) for $\theta_0<1$ which implies a D\"orfler property for any $0< \theta \leq \theta_0$. 
The restrictions on the parameters $\kappa$, $\omega$, and $\theta$ are gathered in the following assumption.

\begin{assumption}[restrictions on $\kappa$, $\omega$ and $\theta$]\label{a:omega_theta-IP}
\index{Assumptions! Restrictions on $\kappa$, $\omega$, and $\theta$}
Assume that $\kappa > \max(\overline{\kappa}_D,\overline{\kappa}_\textrm{conv})$, that  $0<\omega \leq  \frac 1 4 C_\textrm{Mo}^{-1}C_\textrm{L}C_D^{-1}$ and that $0 < \theta \leq \theta_0(\kappa)$, where $\kappa_D$ and $\theta_0$ are defined in Lemma~\ref{L:dorfler-IP} (D\"orfler marking).
\end{assumption}

Note that if Assumption \ref{a:omega_theta-IP} is valid then
\begin{equation}\label{E:def-omega-IP}
| u - \wh{u}_k|_{H^1_0(\Omega)} \le \frac{1}{4} C_\textrm{Mo}^{-1}C_\textrm{L}\eps_k.
\end{equation}

The next results rely on Assumption~\ref{A:approx-u} (approximability of $u$) and Assumption~\ref{A:approx-data} (approximability of data). They are stated and proved for conforming meshes and continuous approximations of $u$. However, Proposition~\ref{p:equiv_classes_dG} (equivalence of classes for $u$) and  Remark~\ref{r:equiv_class_data} (equivalence of classes for $\data$)
show that these classes coincide with the conforming case.

\begin{proposition}[cardinality of marked sets]\label{P:cardinality-IP}
Let Assumptions \ref{A:approx-u} (approximability of $u$),  \ref{A:cardinality} (cardinality of $\mathcal M$), and~\ref{a:omega_theta-IP} (restrictions on $\kappa$, $\omega$ and $\theta$) hold.
If $\wh{\eps}_{k,0}>\eps_k$, then {\rm \GALERKINIP} at iteration $k$ of {\rm \AFEMTSIP} is called and
the cardinality $N_{k,j}(u)$ of the marked set $\marked_{k,j}$ satisfies
\begin{equation}\label{nv-card-Mk-no-osc}
N_{k,j}(u) \lesssim |u|_{s}^{1/s}\,\dGnorm[\kappa][\grid_{k,j}]{u-u_{\grid_{k,j}}}^{-1/s}
\quad\fa 0 \le j < J_k.
\end{equation}
\end{proposition}
\begin{proof}
Fix $0\leq j < J_k$ and set
$$
\delta:= \frac 1 2 C_\textrm{Mo}^{-1}C
_\textrm{L} \, \eta_{\grid_{k,j}}(u_{\grid_{k,j}}) \geq \frac 1 2 C_\textrm{Mo}^{-1} C_\textrm{L}\eps_k,
$$
because $\eta_{\grid_{k,j}}(u_{\grid_{k,j}}) > \eps_k$ for $j<J_k$. 
Thanks to \eqref{E:def-omega-IP}, $\widehat{u}_k $ is an
$(\frac 1 2 C_\textrm{Mo}^{-1} C_\textrm{L}\eps_k)$-approximation of order $s$ to $u$ according to Lemma \ref{L:eps-approx} ($\eps$-approximation of order $s$).  
Therefore, there exists a conforming mesh $\mathcal{T}_\delta\in\grids^\Lambda$ and $ u_{\grid_\delta}^0 \in \mathbb V_{\grid_\delta}$ such that 
$$ 
\dGnorm[\kappa][\grid_\delta]{\wh{u}_k - u_{\grid_\delta}^0} = |\wh{u}_k - u_{\grid_\delta}^0|_{H^1_0(\Omega)} \leq \delta, \qquad 
\#\mathcal{T}_\delta \lesssim | u |_{\As}^{\frac{1}{s}}\delta^{-\frac{1}{s}}.
$$
To compare $\grid_\delta$ with $\mathcal{T}_{k,j}$ we consider the overlay 
$\mathcal{T}_*=\mathcal{T}_{k,j}\oplus\mathcal{T}_\delta$, which satisfies
$$
\# \mathcal{T}_*\leq \#\mathcal{T}_{k,j} + \# \mathcal{T}_\delta - \#\mathcal{T}_0;
$$
see Proposition \ref{P:mesh-overlay-Lambda} (mesh overlay is $\Lambda$-admissible).
Let $u_{\grid_*}\in \mathbb{V}_{\grid_*}^{-1}$ be the Galerkin solution on the subspace $\mathbb{V}_{\grid_*}^{-1}$ and invoke Corollary~\ref{c:cea-dG} (Cea's lemma) to write
$$
\eta_{\grid_*}(u_{\grid_*}) \leq C_{\textrm{L}}^{-1} \dGnorm[\kappa][\grid_*]{\hu_k - u_{\grid_*}} \leq C_{\textrm{L}}^{-1} C_{\textrm{Mo}} \dGnorm[\kappa][\grid_*]{\hu_k - u_{\grid_\delta}^0} = C_{\textrm{L}}^{-1} C_{\textrm{Mo}} |\hu_k - u_{\grid_\delta}^0|_{H^1_0(\Omega)},
$$
whence $\eta_{\grid_*}(u_{\grid_*}) \leq  C_{\textrm{L}}^{-1} C_{\textrm{Mo}}\delta$ and
$$
\eta_{\grid_*}(u_{\grid_*}) \leq \frac 1 2 \eta_{\grid_{k,j}}(u_{\grid_{k,j}}).
$$

Applying Lemma \ref{L:dorfler-IP} (D\"orfler marking) to $\mathcal{T}_*$ and $\mathcal{T}_{k,j}$ we infer that the enlarged refined set $\omega(\mathcal{T}_{k,j}\setminus \mathcal{T}_*)$ satisfies the D\"orfler marking property
$$
 \eta_{\grid_{k,j}}\big(u_{\grid_{k,j}},\omega(\mathcal{T}_{k,j}\setminus \mathcal{T}_*)\big) \geq \theta \, \eta_{\grid_*}(u_{\grid_*}) 
$$
since $0<\theta\leq \theta_0$ by Assumption \ref{a:omega_theta-IP}.
The D\"orfler marking involves a minimal set $\mathcal{M}_{k,j}$ according to Assumption~\ref{A:cardinality}, which thus implies
$$
N_{k,j}(u)\leq \# \omega(\mathcal{T}_{k,j}\setminus \mathcal{T}_*) \lesssim 
\# (\mathcal{T}_{k,j}\setminus \mathcal{T}_*) \le
\# \mathcal{T}_\delta - \#\mathcal{T}_0
\lesssim | u |_{\As}^{\frac{1}{s}}\delta^{-\frac{1}{s}}\lesssim
| u |_{\As}^{\frac{1}{s}}\varepsilon_k^{-\frac{1}{s}}.
$$
because $\#(\mathcal{T}_{k,j}\setminus \mathcal{T}_*) \le \#\mathcal{T}_* - \#\mathcal{T}_{k,j}$.
This concludes the proof.
\end{proof}

\begin{corollary}[quasi-optimality of \GALERKINIP]\label{C:quasi-optimality-IP}
Let Assumptions \ref{A:approx-u} (approximability of $u$),  \ref{A:cardinality} (cardinality of $\mathcal M$), and~\ref{a:omega_theta-IP} (restrictions on $\kappa$, $\omega$ and $\theta$) hold. Assume $\kappa \geq \max(\overline{\kappa}_\textrm{conv},\overline{\kappa}_\textrm{D})$. Then, the total number of marked elements $N_k(u)$ within a call to {\rm \GALERKINIP} satisfies
$$ 
N_k(u) \leq J C_0 | u |_{\As}^{\frac{1}{s}}\varepsilon_k^{-\frac{1}{s}},
$$
where $J\ge J_k$ is a uniform upper bound for the number of iterations of {\rm \GALERKINIP}
according to Lemma \ref{c:cost_dg} (computational cost of \GALERKINIP).
\end{corollary}
\begin{proof}
Use that $N_k(u)=\sum_{j=0}^{J_k-1} N_{k,j}(u)$ and combine Propositions \ref{P:cardinality-IP}
(cardinality of marked sets) and \ref{c:cost_dg} (computational cost of $\GALERKINIP$).
\end{proof}

We finally address the rate-optimality of the two-step algorithm \AFEMTSIP.

\begin{theorem}[rate-optimality of $\AFEMTSIP$]\label{T:optimality-AFEM-IP}
 Let Assumptions \ref{A:approx-u} (approximability of $u$), \ref{A:approx-data} (approximability of data), \ref{A:optim-data} (quasi-optimality of $\DATA$), \ref{A:cardinality} (cardinality of $\mathcal M$), \ref{A:initial-labeling} (initial labeling), and~\ref{a:omega_theta-IP} (restrictions on $\kappa$, $\omega$ and $\theta$) 
 hold. Then, {\rm \AFEMTSIP} gives rise to a sequence $\big(\grid_k,\V_{\grid_k}^{-1},u_{\grid_k}\big)_{k=0}^{K+1}$
 such that
  \begin{equation*}
   \dGnorm[\kappa][\grid_k]{u-u_{\grid_k}}  \leq C(u,\data) \big(\#\mesh_k\big)^{-s}  \quad 1\leq k \leq K+1,
  \end{equation*}
  where $0 < s = \min\{s_u, s_\data\} = \min\{s_u, s_A, s_c, s_f\} \leq\frac{n}{d}$ and
 \begin{equation*}
 C(u,\data) = C_* \Big( |u|_{\A_{s_u}}^{\frac{1}{s_u}} + |\bA|_{\ACA_{s_A}}^{\frac{1}{s_A}} + |c|_{\ACc_{s_c}}^{\frac{1}{s_c}} +
 |f|_{\ACf_{s_f}}^{\frac{1}{s_f}} \Big)^s
\end{equation*}
with constant $C_*>0$ independent of $u$ and $\data$.
  \end{theorem}
\begin{proof}
Assumptions \ref{A:approx-u}, \ref{A:cardinality}, and \ref{a:omega_theta-IP} combined with
Corollary \ref{C:quasi-optimality-IP} for $u$, and Assumptions \ref{A:approx-data} and \ref{A:optim-data}
for $\data$, imply the existence of a constant $C_\#$ such that the total number of marked elements within one loop of \AFEMTSIP is
\[
N_k(u) + N_k(\data) \leq C_\# \Big( | u |_{{\A}_{s_u}}^{\frac{1}{s_u}} + | \data |_{\mathbb{A}_{s_\data}}^{\frac{1}{s_\data}}
\Big)\, \varepsilon_k^{-\frac{1}{s}} \,,
\]
with $s_u, s_\data \le \frac{n}{d}$. Moreover, upon termination $\DATA$ and $\GALERKINIP$ give 
\begin{equation*}
\begin{split}
| u-\widehat{u}_k|_{H^1_0(\Omega)}  & \le \frac 1 4 C_\textrm{Mo}^{-1} C_\textrm{L} \eps_k, \\[5pt]
\dGnorm[\kappa][\grid_{k+1}]{\widehat{u}_k-u_{\grid_{k+1}} }
& \leq C_\textrm{U} \eta_{\mathcal{T}_{k+1}}(u_{\grid_{k+1}})\leq C_\textrm{U} \eps_k \,,
\end{split}
\end{equation*}
because of \eqref{E:def-omega-IP} and \eqref{e:dg_upper2}.
This implies by triangle inequality
\begin{equation*}
\dGnorm[\kappa][\grid_{k+1}]{ u - u_{\grid_{k+1}} } \leq \left( \frac 1 4 C_\textrm{Mo}^{-1} C_\textrm{L}+C_\textrm{U}\right) \varepsilon_k.
\end{equation*}
We finally conclude as in Theorem~\ref{T:optimality-AFEM} (rate-optimality of $\AFEMTS$).
\end{proof}

Remark \ref{R:thresholds} about the role of $\omega$, $\theta^*$
and Remark \ref{R:s-sharp} about the optimality of the result, written after Theorem~\ref{T:optimality-AFEM} for the conforming case, remain valid for the non-conforming case and are not repeated here.

%--------------------------------------------------------------------------------
\subsection{Operator $P_\grid$ and routine \DATA on $\Lambda$-admissible partitions}\label{ss:P_grid_nc}
%--------------------------------------------------------------------------------
%
In this section we have used extensively the notion of $\Lambda$-admissible meshes for the design
and study of dG methods, including forcing $f\in H^{-1}(\Omega)$. To this end, as well as for the 
design of the two-step AFEM for dG, namely $\AFEMTSIP$, the construction of the
local projection $P_\grid f \in \F_\grid$ is critical. We discuss this now.

Recall that for a conforming partition $\grid \in \grids$, $P_\grid f$ is defined as a projection to $\mathbb F_\grid$; see Definition~\ref{D:Pgrid} (projection onto discrete functionals).
The definition and subsequent properties of $P_\grid$ hinge on extensions $E_F$ for $F \in \faces$,
studied in Lemma \ref{L:extending-from-faces} (extending from faces),
as well as on bubble functions $\phi_T$, $T\in \grid$, 
 and $\phi_F$, $F\in \faces$ satisfying Assumption \ref{A:abstract-cut-off} (abstract cut-off).

The definition of the element bubble functions $\phi_T$ in \eqref{elm-bubble} is local to $T$ and is thus unchanged on non-conforming subdivisions. The situation is different for faces.
If $F$ is a conforming face, we have the conforming definitions of $E_F$ and $\phi_F$. Instead, if $F$ is a non-conforming face, i.e. $F=T \cap T_*$ with $g(T_*) > g(T)$, and use a virtual conforming refinement of $\omega_F$ to define $E_F$ and $\phi_F$ as in \eqref{face-bubble}. Recall that $g(T)$ is the generation of $T\in\grid$, and $\grid \in \grids$ is a uniform refinement of $\grid_0$ if and only if $g$ is constant on $\grid$.
Let $\overline{\grid}$ be the uniform refinement of $\grid_0$ containing $T_*$, whence
$g(T)=g(T_*)$ for all $T\in\overline{\grid}$; $\overline{\grid}$ is conforming thanks to Assumption~\ref{A:initial-labeling} (initial labeling) on $\grid_0$. Let $\overline{T} \in \overline{\grid}$ be the element sharing $F$ with $T_*$ (and thus contained in $T$) and let $\overline \omega_F:= T_* \cup  \overline{T}$ be the virtual conforming patch around $F$. We now proceed by defining $E_F$ via \eqref{e:E_F} with $\omega_F$ replaced by $\overline{\omega}_F$ and $\phi_F$ as in \eqref{face-bubble} using the basis functions $\phi_z$, $z \in \vertices \cap F$, associated with $\prod_{T \subset\overline{\omega}_F} \mathbb P_n(T) \cap H^1_0(\overline{\omega}_F)$. Note that because $\grid$ is $\Lambda$-admissible, Proposition~\ref{p-size-domain-influence} guarantees that the diameters of $\overline{T}$, $T$, $T_*$, $\omega_F$ and $\overline{\omega}_F$ are all comparable with constants depending on the initial mesh $\grid_0$ and $\Lambda$.

Assumption \ref{A:abstract-cut-off} is an important ingredient in the analysis of $P_\grid$ and it holds true with $\omega_F$ replaced by $\overline{\omega}_F$ when $F$ is a non-conforming face. Therefore, Remark~\ref{R:Computation-of-P} (local computation) and Corollary~\ref{C:local-near-best-approx-of-P} (local near-best approximation) are valid for $\Lambda-$admissible partitions as well.
Consequently, all the algorithms and results presented in Section~\ref{S:data-approx} (data approximation) readily extend to $\Lambda$-admissible subdivisions as well. We do not dwell on this
matter any longer.

\section{AFEMs for Inf-Sup Stable Problems }\label{S:conv-rates-infsup} 

%\cite{Feischl:2019,Feischl:2022}
% \rhn{(RHN $\longrightarrow$ CC)}

% \begin{itemize}
% \item
%   Relaxed quasi-orthogonality.

% \item
%   Convergence for discrete data, estimator reduction property.

% \item
%   Two-step algorithm and convergence for general data.

% \item
%   Convergence rates 

% \item
%   Growth of the stability constant: the LU approach.
  
% \end{itemize}

We go back to the functional framework introduced in Section \ref{S:inf-sup}. Precisely,
let the bilinear form $\B: \V \times \W \rightarrow  \R$ be continuous and {\em inf-sup stable} (i.e., it satisfies one of the equivalent conditions stated in Theorem \ref{nsv-T:ex1-var-prob} (Ne\v cas)). Given $f\in \Wd$, let $u \in \V$ be the unique solution of the variational problem
\begin{equation}\label{E:bilinear-VW}
     u\in\V:\qquad  \bilin{u}{w} = \dual{f}{w} \qquad  \forall w \in \W.
\end{equation}
 Let $\Vj \subset \V$, $\Wj \subset \W$ be finite dimensional subspaces depending on an integer parameter $j\geq 0$, such that 
$$
\dim{\Vj}= \dim{\Wj}= n_j, \qquad \Vj \subset \Vk[j+1], \qquad \Wj \subset \Wj[j+1].
$$
(Note that the notation has changed with respect to Section \ref{S:Galerkin}, where $\V_N$ was a subspace of dimension $N$. Here $\Vj$ may stand for $\V_{\grid_j}$, where $\grid_j$ is the $j$-th mesh generated by an adaptive algorithm.)

We assume $\B$ to satisfy a {\em uniform} discrete inf-sup condition on any product of subspaces $\Vj\times\Wj$, i.e., there exists a constant $\beta>0$ such that for all $j$
\begin{equation}\label{E:inf-sup}
    \inf_{v\in \Vj}\sup_{w\in \Wj}\frac{\bilin{v}{w}}{\|v\|_\V \|w\|_\W}\ge \beta.
\end{equation}
Let $u_j\in \Vj$ be the solution of the (Petrov-)Galerkin problem
\begin{equation}\label{E:bilinear-VkWk}
    u_j\in \Vj: \qquad \bilin{u_j}{w}= \dual{f}{w} \qquad \forall w \in \Wj.
\end{equation}

The first part of this section, which is mostly based on the recent work by M. Feischl \cite{Feischl:2022}, is devoted to studying the convergence of this approximation. Convergence and rate-optimality of different {\AFEM}s will be  discussed next in Section \ref{S:infsup-optimality}.
Applications will be given to the Stokes problem (see Section~\ref{S:rates-stokes}) and the mixed formulation of a scalar diffusion problem (see Section~\ref{S:rates-mixed-fems}).

%-------------------------------------------------------------------------------------------
\subsection{Linear convergence of inf-sup stable methods} 
%-------------------------------------------------------------------------------------------

We make the following key assumptions that guarantee the convergence of the sequence $u_j$ to $u$ in the $\V$-norm, and comment about them afterwards. The first assumption is a relaxed form of the {\em general quasi-orthogonality} property introduced in \cite{Axioms:2014} as part of an abstract set of axioms of adaptivity.
\begin{assumption}[relaxed quasi-orthogonality]\label{A:relaxed-quasi-optimality}  
\index{Assumptions!Relaxed quasi-orthogonality}
For each $N \in \N$ there exists a nondecreasing constant $C=C(N)$ such that
    \begin{equation}\label{E:relaxed-quasi-optimality}
        \sum^{j+N}_{k=j} \|u_{k+1}- u_k\|^2_\V\le C(N)\|u - u_j\|^2_\V \qquad j \ge 0,
    \end{equation}
    and 
    $$
    C(N)=o(N),\qquad \text{as } N\rightarrow \infty.
    $$
    \end{assumption}
    \begin{assumption}[equivalence of error and estimator]\label{A:equivalence-estimator} 
\index{Assumptions!Equivalence of error and estimator}
    There exist constants $ C_U \geq C_L >0$ and, for each $j \geq 0$, an error estimator $\et[j]=\et[j](u_j)$, such that 
    \begin{equation}\label{E:equivalence-estimator}
     C_L \et[j] \leq \|u - u_j\|_\V \leq C_U \et[j]   \qquad j \ge 0.
    \end{equation}
    %\todo[inline]{RHN: Th estimator is called $\est_k$ earlier. We used the notation $\eta_k$ for the scaled sum of error and estimator in Section \S \ref{S:convergence-coercive}.}
    %  
\end{assumption}
    \begin{assumption}[estimator reduction]\label{A:estimator-reduction} 
    \index{Assumptions!Estimator reduction}
There exist constants $0<\ro[1]<1$ and $C_1>0$ such that
    \begin{equation}\label{E:estimator-reduction}
        \et[j+1]^2 \le \ro[1] \et[j]^2 + C_1 \|u_{j+1} - u_{j}\|^2_\V \qquad j \ge 0.
    \end{equation}
\end{assumption}

\begin{remark}
Assumptions \ref{A:equivalence-estimator} and \ref{A:estimator-reduction} are abstract and allow
for a general convergence theory. In the context of our model problems of Section \ref{S:ex-bvp}, they are valid for discrete data, i.e., if the coefficients of the linear operator corresponding to the bilinear form $\B$ are piecewise polynomials on the adopted meshes, and if
$f \in \F_\mesh$ (see Section \ref{S:discretize-functionals}). We make this concrete in Sections
\ref{S:rates-stokes} and \ref{S:rates-mixed-fems} below.
\end{remark}
%\todo[inline]{RHN: We should include mixed methods and say that $(A,c)$ are pw polynomials.}

\begin{remark}
We comment on the significance of Assumption \ref{A:relaxed-quasi-optimality} upon considering two extreme cases.
\begin{itemize}
    \item[1.] Assumption \ref{A:relaxed-quasi-optimality} with $C(N)=O(1)$ is precisely the general quasi-orthogonality property of \cite{Axioms:2014}. It is valid with $C(N)=1$ for $\V= \W$ and $\B$ symmetric and coercive. Indeed,  
    $$
    \bilin{u_{k+1}-u_k}{u -u_{k+1}}=0\qquad \text{(Galerkin orthogonality)},
    $$
    whence
    $$
    \enorm{u_{k+1}- u_k}^2+  \enorm{u -u_{k+1}}^2= \enorm{u - u_k}^2,
    $$
    where $\enorm{\cdot}$ is the energy norm induced by $\B$. Adding upon $k$ and using telescopic cancellation yields
    \begin{equation*}
    \begin{split}
    \sum_{k=j}^{j+ N}\enorm{u_{k+1}- u_k}^2 &= \sum_{k=j}^{j+ N} \enorm{u - u_k}^2 -  \enorm{u -u_{k+1}}^2 \\&= \enorm{u -u_j}^2 -\enorm{u- u_{j+N+1}}^2 \le  \enorm{u- u_j}^2.
    \end{split}
    \end{equation*}
    Finally, the equivalence \eqref{nsv-enorm-Vnorm} of the norms $\|.\|_\V$ and $\enorm{\cdot}$ yields the result.
    
    \item[2.] Assumption \ref{A:relaxed-quasi-optimality} trivially holds with $C(N)=O(N)$ for $\B$ continuous and inf-sup stable. 
    %
%     This stems from the following abstract result, which will also be useful later on. 
% \begin{lemma}[quasi-monotonicity property]\label{L:quasi-monotonicity}
% Let the bilinear form $\B$ be continuous on $\V \times \W$, and inf-sup stable on some $\V_a \times \W_a \subset \V \times \W$, with inf-sup constant $\beta>0$. If $u_a \in \V_a$ is the Galerkin $\B$-projection of some $u\in \V$ upon $V_a$, defined by
% \[
% \B[u-u_a,w_a]=0 \qquad \forall w_a \in \W_a,
% \]
% and if $\V_b$ is any subspace of $\V_a$, then
% \begin{equation}\label{E:quasi-optimal}
% \| u-u_a \|_\V \leq \left(1+\frac{\|\B\|}{\beta}\right) \| u-v_b \|_\V \qquad \forall v_b \in \V_b.
% \end{equation}
% \end{lemma}
% %
% \todo[inline]{RHN (12/27/23): Please check the new Corollary 4.3. Should we remove this lemma?}
% %
% \begin{proof}
% Let us write $\| u-u_a \|_\V \leq \| u-v_b \|_\V +\| u_a - v_b \|_\V$ for any $v_b\in \V_b$. By inf-sup stability in $V_a$ and Galerkin orthogonality, we get
% \begin{equation*}
% \begin{split}
% \| u_a - v_b \|_\V & \leq \frac1\beta \sup_{w_a \in \W_a}\frac{\B[u_a - v_b,w_a]}{\|w_a\|_\W} \\ &= 
% \frac1\beta \sup_{w_a \in \W_a}\frac{\B[u - v_b,w_a]}{\|w_a\|_\W} \leq \frac{\|\B\|}{\beta} \| u-v_b \|_\V,
% \end{split}
% \end{equation*}
% whence the result.
% \end{proof}
%
Indeed, choosing in Corollary \ref{C:quasi-monotonicty} (quasi-monotonicity) $\V_N=\V_k$ 
or $\V_{k+1}$ and $\V_M=\V_j$
for $j \leq k$, and using the triangle inequality gives
    $$
    \|u_{k+1}- u_k\|^2_\V\lesssim \|u_{k+1}- u\|^2_\V + \|u_k - u\|^2_\V\lesssim \|u - u_j\|^2_\V.
    $$
    Adding, we get 
    $$
    \sum_{k=j}^{j+N}\|u_{k+1}- u_k\|^2_\V\le C \sum_{k=j}^{j+N}\|u- u_j\|^2_\V = C N \|u- u_j\|^2_\V.
    $$
    However, the relation $C(N)= O(N)$ is not enough for the subsequent analysis. In fact, we need $C(N)=o(N)$.
\end{itemize}
\end{remark}

We now prove that the stated assumptions guarantee the linear convergence of the sequence of Petrov-Galerkin solutions \eqref{E:bilinear-VkWk}. This result is similar to the convergence result for the estimators given in \cite{Feischl:2022}, and exploits the equivalence \eqref{E:equivalence-estimator} between errors and estimators.
\begin{theorem}[linear convergence]\label{T:Feischl-theorem} Under Assumptions \ref{A:relaxed-quasi-optimality}, \ref{A:equivalence-estimator} and \ref{A:estimator-reduction}, the discretization \eqref{E:bilinear-VkWk} is convergent; precisely, there exist constants $0<\rho<1$ and $c>0$ such that 
\begin{equation}\label{E:Feischl-theorem}
e_{j+i}\le c \rho^i e_j \qquad \forall i, j \in \N,
\end{equation}   
where  $e_j:= \|u- u_j\|_\V$.
\end{theorem}
\begin{proof}
The proof is divided into several steps. Firstly, we set
$$
E_k := \|u_k - u_{k-1}\|_\V.
$$

   \step{1} We start by iterating \eqref{E:estimator-reduction} $1 \leq n \leq k$ times to obtain
    \begin{equation*}
    \begin{split} 
      \et[k]^2&\le \rho_1 \et[k-1]^2 + C_1 E^2_k\\ 
              &\le \rho_1 \left( \rho_1 \et[k-2]^2 + C_1 E^2_{k-1}\right) + C_1 E^2_k \le \rho^2_1 \et[k-2]^2 + C_1 \left(E^2_k + E^2_{k-1}  \right)\\ 
              &\le \rho_1^n \et[k-n]^2 + C_1 \sum_{\ell = k-n+1}^k E^2_\ell.
    \end{split}
    \end{equation*}
    We now invoke Assumption \ref{A:equivalence-estimator} to state the upper bound 
    $$
    e^2_k \le c_1 \et[k]^2
    $$
    and the lower bound
    $$
    \et[k]^2 \le c_2 e_k^2
    $$
    (with $c_1=C_U^2$ and $c_2=C_L^{-2}$). This yields
    \begin{equation}\label{E:linear_convergence_step1.1.01}
    e^2_k \le c_1 \et[k]^2\le c_1 c_2 \rho_1^n e^2_{k-n} + c_1 C_1 \sum_{\ell = k-n+1}^k E_\ell^2.
    \end{equation}
    Let $n\in \N$ be sufficiently large such that
    $$
    \rho_2= c_1 c_2 \rho_1^n<1, 
    $$
    and let us relabel $c_1 C_1$ as $C_1$ to get
 \begin{equation}\label{E:linear_convergence_step1.1}
        e^2_k \le \rho_2 e^2_{k-n} + C_1 \sum_{\ell = k-n+1}^k E_\ell^2.
    \end{equation}
    This shows that the reduction property \eqref{E:estimator-reduction} of the estimator is valid for the error after $n$ iterations. We cannot expect \eqref{E:linear_convergence_step1.1} to hold for $e_k$ with  $n=1$, not even in the coercive case: see Example \ref{EX:interior-node} (lack of strict monotonicity) for $\bA=\vec{I}$ and $f =1$. It is thus convenient to rewrite \eqref{E:linear_convergence_step1.1} as follows:
    \begin{equation}\label{E:linear_convergence_step1.2}
        e^2_{k n} \le \rho_2 e^2_{(k-1)n} + C_1 \sum_{\ell = (k-1)n+1}^{k n} E_\ell^2.
    \end{equation}

    \step{2} Sum up \eqref{E:linear_convergence_step1.2} from $k = j+1$ to $k = j +N$ to get
    $$
        \sum_{k= j+1}^{j + N}e^2_{k n} \le \rho_2 \sum_{k= j+1}^{j + N}e^2_{(k-1)n} + C_1 \sum_{\ell = j n+1}^{(j + N) n} E_\ell^2.
    $$
    Using \eqref{E:relaxed-quasi-optimality} we see that
    $$
        \sum_{\ell = j n+1}^{(j + N) n} E_\ell^2 = \sum_{\ell = j n}^{(j + N) n-1} E_{\ell+1}^2\le C(N n -1) \, e^2_{j n} \le C (N n) \, e^2_{j n},
    $$
    whence
    \begin{equation*}
    \begin{split}
        \sum_{k=j+1}^{j + N}e^2_{k n}&\le \rho_2 \sum_{k= j +1}^{j+N}e^2_{(k-1)n} + C_1 C (N n) \, e^2_{j n} \\&\le \rho_2 \left(\sum_{k = j+1}^{j + N} e^2_{k n} +e^2_{j n}\right) + C_1 C (N n) \, e^2_{j n}. 
    \end{split}
    \end{equation*}
    This implies 
    $$
    (1 - \rho_2)\sum_{k = j+1}^{j + N} e^2_{k n} \le \big(\rho_2 + C_1 C(N n)\big) \,e^2_{j n} \, ,
    $$
    or equivalently
    \begin{equation}\label{E:linear_convergence_step2.1}
        \frac{1 -\rho_2}{\rho_2 + C_1 C(N n)}\sum_{k = j+1}^{j +N}e^2_{k n} \le e^2_{j n}.
    \end{equation}
    Let us add the quantity $\sum_{k = j+1}^{j+ N} e^2_{k n}$ to both sides to arrive at
    $$
     \left( 1 + \frac{1 -\rho_2}{\rho_2 + C_1 C(N n)}\right)\sum_{k = j+1}^{j +N}e^2_{k n} \le \sum_{k = j}^{j+ N} e^2_{k n}.
    $$
    We can rewrite this inequality as follows:
    \begin{equation}\label{E:linear_convergence_step2.2}
        \sum_{k = j+1}^{j +N}e^2_{k n} \le \rho(N) \sum_{k = j}^{j+ N} e^2_{k n},
    \end{equation}
    where
    \begin{equation*}
    \begin{split}
        \rho(N) &= \frac{1}{1 + \frac{1 -\rho_2}{\rho_2 + C_1 C(N n)}}= \frac{\rho_2 + C_1 C(N n)}{1+ C_1 C(N n)}\\
                & =  1 -  \frac{1 -\rho_2}{1 + C_1 C(N n)} = 1 - \frac{1}{D(N)}
    \end{split}
   \end{equation*}
   with
   $$
   D(N) = \frac{1 + C_1 C(N n)} {1 -\rho_2}\rightarrow \infty, \qquad N \rightarrow \infty,
   $$
   whenever $C(N)$ diverges.
   Therefore, \eqref{E:linear_convergence_step2.2} is a contraction for the quantity $\sum_{k = j+1}^{j+ N} e^2_{k n}$ with a constant $\rho(N)$ uniform in $j$ that may degenerate to 1 as $N\rightarrow\infty$.

   \medskip
   
    \step{3} We iterate \eqref{E:linear_convergence_step2.2} and exploit that the left-hand side has one fewer term than the right-hand side. Take 
   $$
   j \rightarrow j +1 \qquad N \rightarrow N-1
   $$
   to get
   $$
   \sum_{k = j+2}^{j +N}e^2_{k n} \le \rho(N-1) \sum_{k = j+1}^{j+ N} e^2_{k n},
   $$
   whence
   $$
   \sum_{k = j+2}^{j +N}e^2_{k n} \le \rho(N-1) \rho(N)\sum_{k = j}^{j+ N} e^2_{k n}.
   $$
   Iterating, we get
   $$
   e^2_{(j+N)n}= \sum_{k= j+N}^{j+ N}e^2_{k n} \le \rho(1)\dots \rho(N-1) \rho(N) \sum_{k = j}^{j+N}e^2_{k n}.
   $$
   We now need to bound the sum on the right-hand side by a single term. To this end, we resort to \eqref{E:linear_convergence_step2.1}
   $$
  \sum_{k = j+1}^{j+ N} e^2_{k n}\le \frac{\rho_2 + C_1 C(N n)}{1 -\rho_2}e^2_{j n}
   $$
   and add $e^2_{j n}$ on both sides
   $$
   \sum_{k = j}^{j+ N} e^2_{k n}\le \left (1 + \frac{\rho_2 + C_1 C(N n)}{1 -\rho_2}\right)e^2_{j n}= D(N) \, e^2_{j n}.
   $$
   Altogether, we arrive at
   $$
   e^2_{(j+N)n}\le \rho(1)\dots \rho(N) D(N) \, e^2_{j n} = D(N)\prod^N_{k=1}\left(1 - \frac{1}{D(k)}\right)e^2_{j n}.
   $$
   
   \step{4} We estimate the factor on the right-hand side. For $N_0 >0$ to be chosen later, set
   $$
   \rho_0 :=  D(N_0)\prod^{N_0}_{k=1}\left(1 - \frac{1}{D(k)}\right)
   $$
   and compute the logarithm of $\rho_0$
   $$
   \log(\rho_0)= \log(D(N_0))+ \sum_{k=1}^{N_0}\log \left(1 - \frac{1}{D(k)}\right)\le \log(D(N_0)) - \sum_{k=1}^{N_0}\frac{1}{D(k)},
   $$
   because the $\log$ is concave and $\log(1+x)\le x$.
   Since we assume in $\eqref{E:relaxed-quasi-optimality}$ that $D(k)\simeq C(k N_0)=o(k)$, the series diverges and we see that
   $$
   \log(\rho_0)<0
   $$
   for $N_0$ sufficiently large. Summarizing, there exist $N_0>0$ and $0 < \rho_0 <1$ such that
   \begin{equation}\label{E:linear_convergence_step4}
       e^2_{(j + N_0) n} \le \rho_0 \, e^2_{j n} \qquad \forall j \in \N.
   \end{equation} 
   
   \step{5} For any $j,i \in \N$, we now find $c>0$ and $0<\rho<1$ such that the inequality
   \[
   e_{j+i} \le c \rho^i e_j
   \]
   holds. We decompose $j$ and $j+i$ in terms of integers $k,m$
   \begin{align*}
       j &= (k-1)n+\wh{j}, \quad k\ge1, \quad 0\le \wh{j} < n,
       \\
       j+i &= (k+m) n + \wh{i}, \quad m \ge -1, \quad 0\le \wh{i} < n,
   \end{align*}
   and examine first the case $m\ge0$. We further decompose
   \[
   m = a N_0 + b, \quad a,b\in\N, \quad 0\le b < N_0 
   \quad\Rightarrow\quad
   a = \frac{m}{N_0} - \frac{b}{N_0}.
   \]
   Note that
   \[
      kn = j -\wh{j} + n > j, \quad
      m = \frac{i}{n} -  \Big( \frac{\wh{i} - \wh{j}}{n} + 1 \Big), \quad
      a > \frac{i}{n N_0} - \frac{2+N_0}{N_0}.
   \]
   Therefore, invoking Corollary~\ref{C:quasi-monotonicty} (quasi-monotonicity)
   %
   % \todo[inline]{RHN (01/21/24): Claudio note that we refer to a Lemma \ref{L:quasi-monotonicity}, that comes later, rather than a quasi-monotonicity of norms. This statement has been commented out. Please check!}
   $$
   e_{j_2}\le \frac{\|\B\|}{\beta}\, e_{j_1} = C_\ast e_{j_1} \qquad j_2 \ge j_1\ge 0,
   $$
   in conjunction with \eqref{E:linear_convergence_step4}, yields
   \[
   e_{j+i} \le C_* e_{(k+aN_0)n} \le C_* \rho_0^{\frac{a}{2}} e_{kn} \le C_*^2 \rho_0^{\frac{a}{2}} e_j
   < C_*^2 \rho_0^{-\frac{2+N_0}{2 N_0}} \big( \rho_0^{\frac{1}{2 n N_0}} \big)^i e_j.
   \]
   This is the desired estimate with $c= C_*^2 \rho_0^{-\frac{2+N_0}{2 N_0}}$ and $\rho = \rho_0^{\frac{1}{2 n N_0}}$ for $m\ge0$. We finally consider $m=-1$ and use again the error quasi-monotonicity to write
   \[
   e_{j+i} \le C_* e_j = \frac{C_*}{\rho^i} \rho^i e_j < \frac{C_*}{\rho^n}\rho^i e_j.
   \]

This concludes the proof with $c = \max \big\{ C_*^2 \rho_0^{-\frac{2+N_0}{2 N_0}}, \frac{C_*}{\rho^n}\big\}$.
\end{proof}

\begin{remark}[improving on Assumption \ref{A:relaxed-quasi-optimality}]\label{R:improving quasi-optimality}
It is worth underlining that, while linear convergence \eqref{E:Feischl-theorem} is established in Theorem \ref{T:Feischl-theorem} using Assumption \ref{A:relaxed-quasi-optimality}, the same property combined with uniform inf-sup stability tell us a posteriori that the constant $C(N)$ in \eqref{E:relaxed-quasi-optimality} can be made independent of $N$, i.e. $C(N)=O(1)$. To see this,
we apply the linear convergence bound
$$
\Vert u - u_k \Vert_\V \leq  c \rho^{k-j}  \Vert u - u_j \Vert_\V  
$$
in conjunction with the triangle inequality
$$
\Vert u_{k+1}-u_k \Vert_\V  \leq \Vert u - u_{k+1} \Vert_\V + \Vert u-u_k \Vert_\V  \leq   2c \rho^{k-j}  \Vert u - u_j \Vert_\V; 
$$
summation of a geometric series gives
$$
\sum^{j+N}_{k=j} \|u_{k+1}- u_k\|^2_\V  \le C \|u - u_j\|^2_\V 
$$
with $
C = 4c^2 \sum_{\ell=0}^\infty \rho^{2\ell} < +\infty
$.
This suggests that Assumption \ref{A:relaxed-quasi-optimality} might be too pessimistic.
\end{remark}
%\todo[inline]{RHN: We have to use the bookstrapping argument to improve on the Assumption \ref{A:relaxed-quasi-optimality}.}

%-------------------------------------------------------------------------------------------------------
%-------------------------------------------------------------------------------------------------------
\subsection{Inf-sup stability implies quasi-orthogonality}\label{S:Quasi-Orthogonality}
%-------------------------------------------------------------------------------------------------------
%-------------------------------------------------------------------------------------------------------
We aim at proving the following key result in this section.
\begin{theorem}[sufficient condition for Assumption \ref{A:relaxed-quasi-optimality}]\label{T:Feischl-ass} The Assumption \ref{A:relaxed-quasi-optimality} (relaxed quasi-orthogonality)
is valid if the bilinear form  $\B: \V \times \W \rightarrow  \R$ is continuous and uniformly inf-sup stable on the sequence of subspaces $\Vj \times \Wj$, $j \geq 0$.   \looseness=-1
\end{theorem}
To accomplish this task, we proceed in two steps. Using variational techniques, we first
establish an intermediate result formally similar to \eqref{E:relaxed-quasi-optimality} (see Corollary \ref{C:quasi-orthogonality}), but involving the norm of a matrix $\Ub$ related to the form $\B$. Next, we rely on 
algebraic techniques to estimate such norm (see Theorem \ref{T:bound-matrices}) and close the proof of the desired result.

% We devote this and the subsequent subsection to the proof of this theorem. Precisely, here we establish an intermediate result formally similar to \eqref{E:relaxed-quasi-optimality} (see Corollary \ref{C:quasi-orthogonality}), but involving the norm of a matrix $\Ub$ related to the form $\B$. In Section \ref{S:growth}, we use algebraic techniques to estimate such norm and close the proof of the theorem.

In order to perform the first step, we introduce orthonormal bases of the finite element spaces $\Vj$, $\Wj$, $0\le j\le N$,
and next we biorthogonalize them. This procedure turns out to be crucial. 

We start with some notation. Let
$$
n_j = \dim \Vj = \dim \Wj.
$$
Let $\Vkort[j-1]$ and $\Wkort[j-1]$ denote the orthogonal complements of $\Vj[j-1]$ and $\Wj[j-1]$ within $\Vj$ and $\Wj$, respectively. Let 
$$
d_j = n_j - n_{j-1} =\dim \Vkort[j-1] = \dim \Wkort[j-1]
$$
be the dimension of the Galerkin update to augment the space $\Vk[j-1]$ into the next space $\Vj$, and likewise with the space $\Wk[j-1]$ and $\Wj$.

We consider orthonormal bases
\begin{equation}\label{E:orthonormal_basis}
    \vb = \left\{\vb(j)\right\}_{j=0}^N\;\subset \Vk[N], \qquad \wb = \left\{\wb(i)\right\}_{i=0}^N\;\subset \Wk[N]
\end{equation}
partitioned into blocks for $1\le j\le N$
\begin{equation}\label{E:orthogonal_basis_blocks}
    \vb(j) =\left(v_k\right)^{n_j}_{k = n_{j-1}+1}\subset \Vkort[j-1], \qquad \wb(i) =\left(w_k\right)^{n_i}_{k = n_{i-1}+1}\subset \Wkort[i-1],
\end{equation}
and $\vb(0) \subset \Vk[0]$, $\wb(0) \subset \Wk[0]$.
In other words, $\left(\vb(j), \wb(j)\right)$ represent the $d_j$ new directions added by Galerkin to the current spaces $\left(\Vk[j-1],\Wk[j-1]\right)$ for $1\le j\le N$.

We recall that the bilinear form $\B:\Vk[N]\times \Wk[N] \rightarrow \R$ satisfies the following uniform properties for all $0\le j\le N$
\begin{enumerate}[label=(P\arabic*)]
    \item {\it continuity}: 
    \begin{equation}\label{E:prop-B-infsup-1}
    \big|\bilin {v}{w} \big| \le \|{\B}\|\|v\|_\V \|w\|_\W \qquad\forall v \in \Vk[j], \; w \in \Wk[j];
    \end{equation}
    \item {\it inf-sup condition}: 
    \begin{equation}\label{E:prop-B-infsup-2}
    \beta \|v\|_\V \le \sup_{w \in \Wk[j]}\frac{\bilin{v}{w}}{\|w\|_\W}\qquad\forall v \in \Vk[j].
    \end{equation}
\end{enumerate}
The block bases $\vb$ and $\wb$ given in \eqref{E:orthonormal_basis} induce a block matrix
$$
\Bb := \left(\Bb (i,j)\right)_{i,j=0}^N \; \in \R^{n_N\times n_N}
$$
defined by
\begin{equation}\label{E:matrixB}
 \Bb (i,j) = \bilin{\vb(j)}{\wb(i)}.   
\end{equation}
Note that the actual size of $\Bb$ is $n_N = \dim \Vk[N]>\!>N$, and that the following analysis entails expressing important quantities in terms of the number of blocks $N$ rather than the dimension $n_N$.

\begin{figure}[ht]
\includegraphics[width=0.90\textwidth]{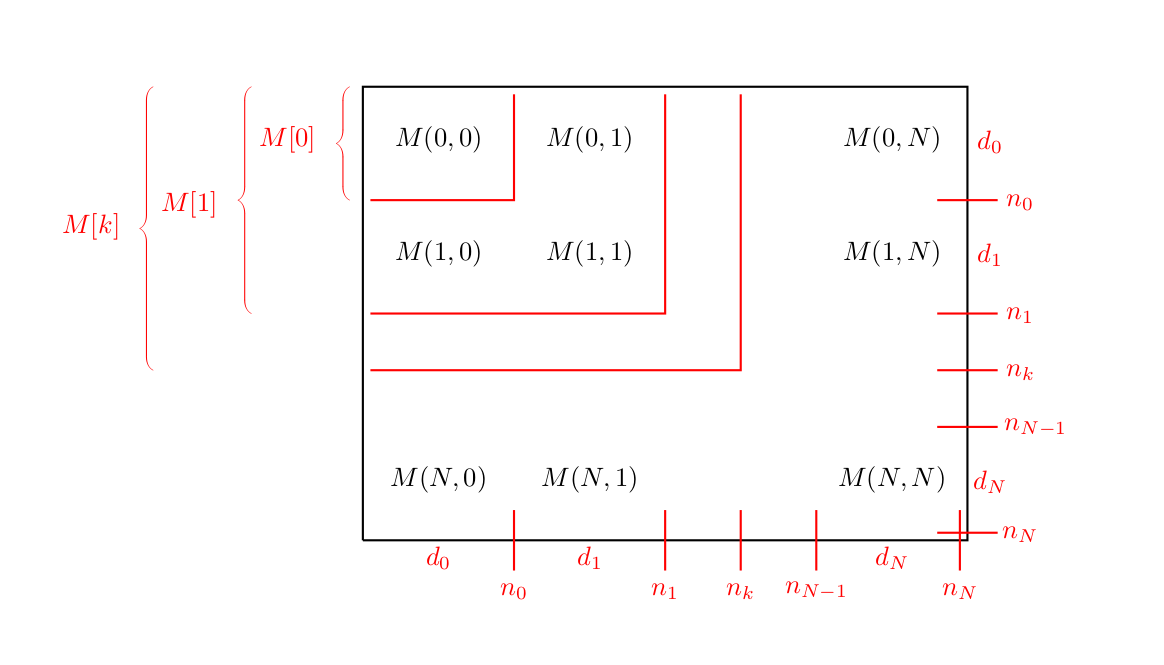}
\caption{Block partition of a matrix $\Mb \in \R^{n_N\times n_N}$ with $(N+1)\times(N+1)$ blocks $\Mb(i,j)\in \R^{d_i\times d_j}$ and principal $k$-th block $\Mb[k] \in \R^{n_k\times n_k}$ with $0\le i,j,k \le N.$}
\label{F:M-block-partition}
\end{figure}

We will use this block decomposition for a generic matrix 
$$
\Mb = \left(\Mb (i,j)\right)^N_{i, j = 0}  \in \R^{n_N \times n_N}
$$
and we denote by $\Mb[k]=\left(\Mb (i,j)\right)^k_{i, j = 0}$ the {\em principal $k$-th block} of $\Mb$. Fig.~\ref{F:M-block-partition} shows schematically what this means. 

%%%%%%%%%%%%
%%%%%%%%%%%%

We stress that (P2) implies that $\Bb[k]$ is uniformly invertible with
\begin{equation}\label{E:invertibleM}
    \|\Bb[k]^{-1}\|_2\le \frac{1}{\beta} \qquad\forall \, 0\le k \le N.
\end{equation}
In fact, (P2) with $j$ replaced by $k$ can be rephrased as follows in terms of the coordinates $\vbb\in\R^{n_k}$ relative to the orthonormal basis $\{\vb(j)\}_{j=0}^k$ of $\Vk[k]$ of a generic vector in $\Vk[k]$
$$
\beta \| \vbb \|_2 \leq \| \Bb[k] \vbb \|_2 \qquad \forall \vbb \in \mathbb{R}^{n_k} \,,
$$
i.e., setting $\zbb =\Bb[k] \vbb $,
$$
\beta \, \| \Bb[k]^{-1} \zbb \|_2 \leq \| \zbb \|_2 \qquad \forall \zbb\in \mathbb{R}^{n_k} \,,
$$
which is precisely \eqref{E:invertibleM}.

A fundamental linear algebra theorem of Gaussian elimination guarantees the existence of a unique normalized block $\Lb \Ub$ decomposition of $\Bb$ without pivoting due to \eqref{E:invertibleM}:
\begin{equation}\label{E:M=LU}
    \Bb =  \Lb \Ub,
\end{equation}
with block partitioning 
\begin{equation}\label{E:Ldefinition}
    \Lb(i,j)\in\R^{d_i\times d_j},\qquad \Lb(i,j)=\boldsymbol{0} \text{ \ for \ } j>i,\qquad \Lb(i,i)=\Ib(i,i);
\end{equation}
\begin{equation}\label{E:Udefinition}
    \Ub(i,j)\in\R^{d_i\times d_j},\qquad \Ub(i,j)=\boldsymbol{0} \text{ \ for \ } i>j.
\end{equation}
\subsubsection{Matrix representation} 
The $k$-th Galerkin solution $u_k$ satisfies:
$$
u_k \in \Vk: \qquad \bilin{u_k}{w} = \scp{f}{w}, \qquad  \forall w \in \Wk.
$$
Equivalently, if $\{\gammab(j)\}_{j=0}^k \in \R^{n_k}$ are the coordinates of $u_k$ with respect to the orthonormal basis $\{\vb(j)\}_{j=0}^k$
$$
u_k = \sum_{j=0}^k \gammab(j)\cdot \vb(j),
$$
then
$$
\sum_{j=0}^k \gammab(j)\cdot \bilin{\vb(j)}{\wb(i)} = \scp{f}{\wb(i)}, \qquad \forall \;0\le i\le k,
$$
or using matrix notation
\begin{equation}\label{E:Bmatrixformulation}
    \sum_{j=0}^k \Bb(i,j) \, \gammab(j) = \fb(i)=\scp{f}{\wb(i)}, \qquad \forall  \;0\le i\le k.
\end{equation}
If we further write
$$
\gammab_k = \left(\gammab(j)\right)_{j=0}^k \in \R^{n_k}, \qquad \fb_k =(\fb(i))^k_{i=0} \in \R^{n_k},
$$
then \eqref{E:Bmatrixformulation} reduces to
\begin{equation}\label{E:Bmatrixformulation-reduced}
    \Bb[k]\gammab_k = \fb_k.
\end{equation}
In view of the definition of $\fb_k$ we realize that the $k$-th section $\fb_N[k]\in\R^{n_k}$ of $\fb_N$ coincides with $\fb_k$:
$$
\fb_N[k]= \big(\fb(i)\big)_{i=0}^k = \fb_k.
$$
However, this statement is {\it not} true for the solution $\gammab_k$ of \eqref{E:Bmatrixformulation-reduced}, namely
$$
\gammab_N[k] \not= \gammab_k.
$$

%------------------------------------------------------------------------------------------
\subsubsection{Block biorthogonal bases} 
%------------------------------------------------------------------------------------------
We define biorthogonal bases $\vbt \subset \Vk[N]$ and $\wbt \subset\Wk[N]$ as follows:
\begin{equation}\label{E:definition-vtilde}
   \vbt :=\Ub^{-T}\vb \quad\Longrightarrow \quad  \vbt(j) = \sum_{m=0}^j \Ub^{-T}(j,m)\vb(m) \quad 0 \leq j \leq N,
\end{equation}
\begin{equation}\label{E:definition-wtilde}
  \wbt := \Lb^{-1} \wb \quad \Longrightarrow \quad \wbt(i) =\sum_{m=0}^i \Lb^{-1}(i,m)\wb(m) \quad 0 \leq i \leq N.
\end{equation}
We will see below that these bases are convenient to represent the Galerkin solution $u_k \in \V_k$. We start with a list of properties.

\begin{lemma}[span of new bases] The vectors $\vbt$ and $\wbt$ are bases of $\V_N$ and $\W_N$, respectively, and satisfy
$$
\textrm{span}\{\vbt(j)\}^k_{j=0}=\textrm{span}\{\vb(j)\}^k_{j=0},
$$
$$
\textrm{span}\{\wbt(i)\}^k_{i=0}=\textrm{span}\{\wb(i)\}^k_{i=0}.
$$ 
\end{lemma}   
\begin{proof}
    This relies on the fact that $\Ub^{-T}$ and $\Lb^{-1}$ are lower triangular and the diagonal blocks are non-singular (i.e., both $\Lb$ and $\Ub$ are invertible).
\end{proof}

Consider now the matrix $\Bbt$ induced by $(\vbt,\wbt)$, namely
\begin{equation}\label{E:definitionBtilde}
    \Bbt := \bilin{\vbt}{\wbt}\in \R^{n_N\times n_N}\,.
\end{equation}
\begin{lemma}[biorthogonality] The block matrix $\Bbt$ is equal to the identity, namely
$$
\Bbt(i,j) = \Ib(i,j) \qquad\forall\;0\le i, j \le N.
$$
\end{lemma}
\begin{proof}
    We simply combine the definition \eqref{E:definitionBtilde} with \eqref{E:definition-vtilde} and \eqref{E:definition-wtilde} to deduce, for all $0\le i,j\le N$, that
    \begin{equation*}
    \begin{split}
        \Bbt(i,j) &= \bilin{\vbt(j)}{\wbt(i)}\\
             &= \B \Big[\sum_{m=0}^j \Ub^{-T}(j,m)\vb(m),\sum_{k=0}^i \Lb^{-1}(i,k)\wb(k)\Big]\\
             %&= \bilin{ \sum_{m=0}^j \Ub^{-T}(j,m)\vb(m)}{\sum_{k=0}^i \Lb^{-1}(i,k)\wb(k)}\\
             &= \sum_{m=0}^j \sum_{k=0}^i \Lb^{-1}(i,k)\bilin{ \vb(m)}{\wb(k)} \Ub^{-T}(j,m) \\
             &= \sum_{m=0}^j \sum_{k=0}^i \Lb^{-1}(i,k)\Bb(k,m) \Ub^{-1}(m,j)\\
             &= \left(\Lb^{-1}\Bb\Ub^{-T}\right)(i,j)\\
             &= \left(\Lb^{-1}\left(\Lb \Ub\right)\Ub^{-1}\right)(i,j)\\
             &=  \Ib(i,j),
    \end{split}
    \end{equation*}
    as asserted.
\end{proof}

Generic functions $v \in \Vk[N]$ and $w \in \Wk[N]$ can be represented as follows in terms of the old and new bases:
\begin{equation}\label{E:representation-of-v-via-vt}
    v= \sum_{j=0}^N \gammab(j)\cdot \vb(j)= \sum_{j=0}^N \gammabt(j) \cdot \vbt(j),
\end{equation}
\begin{equation}\label{E:representation-of-w-via-wt}
    w =\sum^N_{i=0}\alphab(i)\cdot\wb(i)= \sum_{i = 0}^N \alphabt (i) \cdot\wbt(i).
\end{equation}
The following lemma relates the coordinates in the two systems.

\begin{lemma}[change of basis]
    The coordinates $\gammab =\left( \gammab(j)   \right)^N_{j=0}$ and $\alphab=\left(\alphab(i)\right)^N_{i=0}$ satisfy
    \begin{equation}\label{E:change-of-basis}
     \gammab = \Ub^{-1}\gammabt, \qquad \alphab = \Lb^{-T}\alphabt.   
    \end{equation}
\end{lemma}
\begin{proof}
    Write \eqref{E:representation-of-v-via-vt} in vector form and use \eqref{E:definition-vtilde} to obtain
    $$
    v =\gammabt^T \vbt= \gammabt^T\left(\Ub^{-T}\vb\right) = \left(\Ub^{-1}\gammabt\right)^T \vb = \gammab^T \vb,
    $$
    whence
    $$
    \gammab = \Ub^{-1}\gammabt.
    $$
    Similarly, combining \eqref{E:representation-of-w-via-wt} and \eqref{E:definition-wtilde} yields
    $$
    w = \alphabt^T \wbt = \alphabt^T \left(\Lb^{-1}\wb\right)= \left(\Lb^{-T}\alphabt\right)^T \wb = \alphab^T \wb
    $$
    and
    $$
    \alphab = \Lb^{-T} \alphabt.
    $$
    This concludes the proof.
\end{proof}
\subsubsection{An intermediate inequality}
We intend to prove the following crucial estimate. This result is in \cite{Feischl:2022}, but we give a different proof based on variational arguments. 
\begin{proposition}[quasi-orthogonality - I]\label{P:quasi-orthogonality}
    If $u_k\in\Vk$ denotes the $k$-th Galerkin solution of \eqref{E:bilinear-VkWk}, then there holds
    \begin{equation}\label{E:quasi-orthogonality}
        \sum_{k=0}^{N-1}\|u_{k+1}- u_k\|^2_\V\le \frac{\|\Ub\|_2^2}{\beta^2}\|u_N - u_0\|^2_\V.
    \end{equation}
\end{proposition}
\begin{proof}
    We proceed in several steps.
    
    \step{1}  {\it Estimate of $\|u_{k-1}- u_k\|_\V$.} Galerkin orthogonality yields
        \begin{equation}\label{E:quasi-orthogonality-step1}
            \bilin{u_{k+1}- u_k}{w}=0\quad \forall \, w\in \Wk.
        \end{equation}
        The uniform discrete inf-sup property (P2) implies the existence of $\wover \in \Wk[k+1]$ with $\|\wover\|_\W=1$ such that
        \begin{equation}\label{E:quasi-orthogonality-step1_1}
            \beta \|u_{k+1} - u_k\|_\V\le\bilin{u_{k+1}- u_k}{\wover}.
        \end{equation}
        We decompose $\wover$ orthogonally as follows:
        $$
        \wover = \wover_k +\woverort_k, \qquad \wover_k\in\Wk, \quad \woverort_k\in \Wkort,
        $$
        with $\|\woverort_k\|_\W \le 1$. In view of \eqref{E:quasi-orthogonality-step1}, \eqref{E:quasi-orthogonality-step1_1} also reads 
        $$
        \beta \|u_{k+1}-u_k\|_\V \le \bilin{u_{k+1}- u_k}{\woverort_k} \le \bilin{u_{k+1}- u_k}{\frac{\woverort_k}{\|\woverort_k\|_\W}}.
        $$
        We now let $\what_{k+1}:= \frac{\woverort_k}{\|\woverort_k\|_\W}\in\Wkort \subset\Wk[k+1]$ and decompose it along the oblique subspaces $\Wk =  \textrm{span}\{\wbt(j)\}^k_{j=0}$ and $\textrm{span}\{\wbt(k+1)\}$, as illustrated in Fig.~\ref{F:oblique-decomposition}.
        \begin{figure}
        \begin{center}
        \begin{tikzpicture}[scale = 0.9]
        %\draw [help lines] (-4,0) grid (4,4);
        \draw [->] (-4,-1/2) -- (4,1/2);
        \draw [->] (1/8,-1) -- (-1/2,4);
        \draw [->, blue] (-0.7,-0.7) -- (3,3);
        \draw [->,red, thick] (0,0) -- (-1/4,2);
        \draw [dashed,red] (-1/4,2) -- (16/7,16/7);
        \draw [dashed,red] (-18/7,-18/56) -- (-1/4,2);
        \draw [->,red, dashed,thick] (0,0) -- (16/7,16/7);
        \draw [->,red, dashed,thick] (0,0) -- (-18/7,-18/56);
        \coordinate [label= {[blue,scale =0.85]$\textrm{span}\{\wbt(k+1)\}$}] (z1) at (4.5,3);
        \coordinate [label= {[scale =0.85]$\Wkort = \textrm{span}\{\wb(k+1)\}$}] (z1) at (1,4);
        \coordinate [label= {[scale =0.85]$\Wk = \textrm{span}\{\wbt(j)\}^k_{j=0}$}] (z1) at (5,-0.55);
        \coordinate [label= {[red,scale =0.85]$\what_{k+1}$}] (z1) at (-0.55,1);
        \end{tikzpicture}
        \end{center}
        \caption{Oblique decomposition of the space $\Wk[k+1]$ into the subspaces $\Wk =  \textrm{span}\{\wbt(j)\}^k_{j=0}$ and $\textrm{span}\{\wbt(k+1)\}$.}
        \label{F:oblique-decomposition}
        \end{figure}
        Since $\Wkort= \textrm{span}\{\wb(k+1)\}$ and $\wb(k +1) =\left(w_j\right)_{j=n_k +1}^{n_{k+1}}$ is an orthonormal basis, the function $\what_{k+1}\in \Wkort$ can be written uniquely as
        $$
        \what_{k+1}= \alphab(k+1) \cdot \wb(k+1)
        $$
        with $\alphab(k+1)\in\R^{d_{k+1} }$ satisfying
        $$
        \|\alphab(k+1)\|_2 = 1 =  \|\what_{k+1}\|_\W.
        $$
        Invoking \eqref{E:definition-wtilde}, we can express $\wb(k+1)$ in terms of $\{\wbt(j)\}_{j=0}^k$ as follows:
        \begin{equation*}
        \begin{split}
          \wb(k+1)&=\sum_{j=0}^{k+1}\Lb(k+1,j) \,\wbt(j)\\
                  &= \wbt(k+1)+ \sum_{j=0}^k \Lb(k+1,j) \,\wbt(j)
          \end{split} 
        \end{equation*}
        because $\Lb(k+1,k+1)=\Ib(k+1,k+1) \in \R^{d_{k+1}\times d_{k+1}}$. Consequently
        $$
        \bilin{u_{k+1}-u_k}{\what_{k+1}}= \bilin{u_{k+1}-u_k}{\alphab(k+1) \cdot \wbt(k+1)}
        $$
        because \eqref{E:quasi-orthogonality-step1} implies
        $$
        \bilin{u_{k+1}-u_k}{\wbt(j)}=0 \qquad \forall\; 0\le j\le k.
        $$
        In addition, the biorthogonality of $\wbt(k+1)$ with respect to $\vbt(j)$ for $0\le j\le k$ translates into
        $$
        \bilin{u_k}{\wbt(k+1)}=0= \bilin{u_0}{\wbt(k+1)}.
        $$
        Moreover, Galerkin orthogonality yields
        $$
        \bilin{u_{k+1}}{\wbt_{k+1}}= \scp{f}{\wbt_{k+1}}= \bilin{u_N}{\wbt_{k+1}},
        $$
        and collecting the preceding expressions we obtain
        \begin{equation}\label{E:quasi-orthogonality-step1_end}  
        \|u_{k+1}- u_k\|\le \frac{1}{\beta} \alphab(k+1)\cdot \bilin{u_N-u_0}{\wbt(k+1)}.
        \end{equation}
        \step{2} {\it Estimate of $\bilin{u_N-u_0}{\wbt(k+1)}$.} We exploit the biorthogonality between $\{\vbt(j)\}^N_{j=0}$ and $\{\wbt(j)\}^N_{j=0}$. In fact, we write
        $$
        u_N - u_0 = \sum_{j=0}^N\gammabt(j)\cdot \vbt(j)
        $$
        and substitute into the right-hand side of \eqref{E:quasi-orthogonality-step1_end} to arrive at
        $$
        \bilin{u_N- u_0}{\wbt(k+1)}= \sum_{j=0}^N\gammabt(j)\cdot \bilin{\vbt(j)}{\wbt(k+1)}=\gammabt(k+1).
        $$
        Therefore, \eqref{E:quasi-orthogonality-step1_end} gives
        $$
        \|u_{k+1}- u_k\|_\V\le \frac{1}{\beta}\alphab(k+1)\cdot \gammabt(k+1),
        $$
        whence
        $$
        \|u_{k+1}-u_k\|_\V\le \frac{1}{\beta}\|\gammabt(k+1)\|_2
        $$
        because $\|\alphab(k+1) \|_2=1$.

        \medskip
        \step{3} {\it Final estimate.} Compute
        $$
                \sum_{k=0}^{N-1}\|u_{k+1}- u_k\|^2_\V \le \frac{1}{\beta^2}\sum_{k=0}^{N-1}\|\gammabt(k+1)\|_2^2
                \le \frac{1}{\beta^2}\|\gammabt\|_2^2\le \frac{\|\Ub\|_2^2}{\beta^2}\|\gammab\|_2^2,
       $$
        according to \eqref{E:change-of-basis}. Since $\{\vb(j)\}^N_{j=0}$ are orthonormal, we get
        $$
        u_N - u_0 = \sum^N_{j=0}\gammab(j)\cdot\vb(j) \quad \Rightarrow\quad \|u_N-u_0\|_\V =\|\gammab\|_2
        $$
        and
        $$
        \sum_{k=0}^{N-1}\|u_{k+1}-u_k\|_\V^2 \le \frac{\|\Ub\|^2_2}{\beta^2}\|u_N-u_0\|_\V^2,
        $$
        as asserted. This concludes the proof. 
\end{proof}

In order to get the quasi-orthogonality estimate, we still need to compare the errors $\|u_N - u_0\|_\V$ and $\|u- u_0\|_\V$. The following is a variant of \eqref{E:discrete-stability}.

\begin{lemma}[stability]\label{L:quasi-monotonicity}
There holds
    $$
    \|u_N - u_0\|_\V\le\frac{\|\B\|}{\beta}\|u-u_0\|_\V.
    $$
\end{lemma}
\begin{proof}
    We use \eqref{E:prop-B-infsup-2} and \eqref{E:prop-B-infsup-1}, in this order, together with Galerkin orthogonality to deduce
    \begin{equation*}
    \begin{split}
     \beta \|u_N- u_0\|_\V &\le \sup_{w\in \Wk[N]}\frac{\bilin{u_N-u_0}{w}}{\| w\|_\W}\\
                           &= \sup_{w\in\Wk[N]}\frac{\bilin{u-u_0}{w}}{\|w\|_\W}\le \|\B\| \|u- u_0\|_\V.
    \end{split}
    \end{equation*}
    This completes the proof.
\end{proof}
\begin{corollary}[quasi-orthogonality - II]\label{C:quasi-orthogonality}
Let $u_k\in\Vk[k]$ be the $k$-th Galerkin solution of \eqref{E:bilinear-VkWk}. Then, for all
$0\le j\le N$, we have
$$
\sum_{k=j}^{j+N-1}\|u_{k+1}- u_k\|^2_\V\le \frac{\|\B\|^2}{\beta^4}\|\Ub\|^2_2\|u - u_j\|^2_\V.
$$
\end{corollary}
\begin{proof}
Combining Proposition \ref{P:quasi-orthogonality} with Lemma \ref{L:quasi-monotonicity} yields
$$
 \sum_{k=0}^{N-1}\|u_{k+1}- u_k\|^2_\V\le \frac{\|\B\|^2}{\beta^4}\|\Ub\|^2_2\|u - u_0\|^2_\V.
$$
Finally, replacing $u_0\in \Vk[0]$ by the $j$-th Galerkin solution $u_j\in\V_j$, we obtain the desired estimate.
\end{proof}

 This corollary says that in order to prove Theorem \ref{T:Feischl-ass}, i.e., to check the validity of Assumption \ref{A:relaxed-quasi-optimality}, it is enough to investigate the growth of the block triangular factor $\Ub$ introduced in \eqref{E:M=LU}, and precisely to prove that
$$
\|\Ub\|_2=o(N^{1/2}).
$$
This is the second step of our analysis. Actually, we prove something more, which is expressed by the following result. 
\begin{theorem}[bound of block matrices $\Lb$ and $\Ub$]\label{T:bound-matrices} 
     There exist constants $C_{L U}>0$ and $p>2$ such that the block $\Lb \Ub$ factors of $\Bb$ satisfy
     \begin{equation}\label{T:inf-sup_bound_LU}
         \|\Ub\|_2 +\|\Lb\|_2 + \|\Ub^{-1}\|_2 + \|\Lb^{-1}\|_2 \le C_{L U} \; N^{1/p}.
     \end{equation}
     \end{theorem}

The proof of this theorem is lenghty and very technical; it involves subtle linear algebra arguments, which may not be familiar to many readers.
For such reasons, we prefer to postpone it to the end of this section (see Section \ref{S:growth}).

%-----------------------------------------------------------------------------------------
%-----------------------------------------------------------------------------------------
\subsection{Convergence rates of AFEMs for inf-sup stable methods}\label{S:infsup-optimality}
%-----------------------------------------------------------------------------------------
In this section, we discuss {\AFEM}s to solve a boundary-value problem admitting a variational formulation of the form 
 \begin{equation}\label{var-prob-gener}
    u\in\V:\quad \bilin{u}{w}=\dual{f}{w} 
    \quad\forall \, w\in\W,
 \end{equation} 
 in which the bilinear form $\B$ on $\V \times \W$ is continuous and inf-sup stable, with inf-sup constant $\beta>0$, and $f \in \W^*$. We will consider the one-step \AFEM given by Algorithm \ref{A:GALERKIN} (\GALERKIN) when all data are discrete,  the one-step \AFEM with switch given by Algorithm \ref{A:one-step-switch} (\AFEMSW) when the operator coefficients are discrete but the forcing term is not (as in the Stokes problem), and the general two-step \AFEM given by Algorithm \ref{algo:AFEM-Kernel} (\AFEMTS).

\medskip
\noindent
\subsubsection{\AB{ Algorithm \ref{A:GALERKIN} (\GALERKIN\!)}} For $j \geq 0$, let us denote by $(\grid_j,\V_j, u_j)$ with $u_j\in \V_j = \V_{\grid_j}$, the sequence of meshes, subspaces and Galerkin approximations to \eqref{var-prob-gener} generated by \GALERKIN. Let 
\begin{equation}\label{E:abstract-estimator}
\eta_{\grid_j}(v) \AB{ = \eta_{\grid_j}(v,f)} = \Big(\sum_{T \in \grid_j} \eta_{\grid}(v,T)^2 \Big)^{1/2}
\end{equation}
be the PDE error estimators used in the loop. If such estimators fulfil Assumptions \ref{A:equivalence-estimator} (equivalence of error and estimator) and \ref{A:estimator-reduction} (estimator reduction), then Theorem \ref{T:Feischl-theorem} (linear convergence) applies and the following result holds.
\begin{proposition}[convergence and termination of \GALERKIN]
\label{P:term-GAL} The module \\ {\rm \GALERKIN} produces a sequence $\{u_j\}$ converging linearly 
to $u\in\V$
\[
\Vert u-u_{j+i} \Vert_\V \leq C\, \rho^i \Vert u-u_{j} \Vert_\V \qquad \forall j,i \geq 0, \qquad 0<\rho<1,
\]
thereby reaching any prescribed accuracy $\| u - u_j \|_\V \leq\varepsilon$ in a finite number of iterations.
\end{proposition}
%\begin{proof}
%To terminate the iteration, it is enough to impose the condition
%\[
%\Vert u-u_{j} \Vert_\V \leq C\, \rho^j \Vert u-u_0 \Vert_\V \leq \varepsilon.
%\]
%\end{proof}

\subsubsection{\AB{Algorithm \ref{A:one-step-switch} (\AFEMSW).}}\label{S:AFEMSW-infsup}

\AB{This algorithm applies to the situation in which the operator coefficients are discrete, whereas the forcing $f \in \W^*$ is not. Then, the PDE estimator $\eta_{\grid}(v,f)$ depends on $f$ via a projection $P_\grid f$ upon a finite dimensional subspace of $\W^*$. Inspired by Lemma \ref{L:loc-of-H^-1-norm} (localization of $H^{-1}$-norm) we denote by $\W^*_\grid$ a suitable decomposition of $\W^*$ subordinate to $\grid$ with norm $\|f\|_{\W^*_\grid}$. In this part of the discussion, we prefer to make the dependence of $\eta_\grid$ upon $P_\grid f$ explicit to avoid confusion, so we will write $\eta_{\grid}(v,P_\grid f)$ rather than $\eta_{\grid}(v,f)$ as usual.}

\AB{Let us begin by stating} two assumptions on the estimator \eqref{E:abstract-estimator} to be used in the sequel.
\begin{assumption}[Lipschitz continuity of 
 estimator]\label{A:Lipschitz-estimator}
\index{Assumptions!Lipschitz continuity of estimator}
There exists a constant \\ $C_{\text{Lip}}>0 $ such that for any $\grid \in \grids$, \AB{ any $v,w \in \V_{\grid}$ and any $f,g \in \W^*$, we have 
\[
|\eta_{\grid}(v,P_\grid f) - \eta_{\grid}(w,P_\grid g)| \leq  C_{\text{Lip}} \big(\| v - w \|_\V + \|P_\grid f-P_\grid g \|_{\W^*_\grid} \big).
\]
%where $\| F \|_{\V^*(\grid)}^2 :=\sum_{T \in \grid}\| F \|_{\V^*(\omega_T)}^2$ and $\V^*(\omega_T)$ %denotes a local version of the dual space $\V^*$ on the star $\omega_T$.
}
\end{assumption}
\begin{assumption}[monotonicity of estimator]\label{A:monotonicity-estimator}
\index{Assumptions!Monotonicity of estimator}
If $\grid \in \grids$ and $\grid_*$ is a refinement of $\grid$, then \AB{the projection operator satisfies $P_{\grid_*} P_\grid = P_\grid$ and
\[
\eta_{\grid_*}(v,P_\grid f,T) \leq \eta_{\grid}(v,P_\grid f,T) \qquad \forall T \in \grid \quad \text{and} \quad \forall v \in \V_{\grid}, \ \forall f \in \W^*.
\]
}
\end{assumption}

It is useful for the subsequent applications to have explicit criteria which guarantee the fulfilment of Assumption \ref{A:estimator-reduction}. This is the purpose of the following result.\looseness=-1
\AB{
\begin{proposition}[estimator reduction under D\"orfler marking]\label{P:local-estimator-reduction}
Let the estimator $\eta_\grid(v,P_\grid f)$ in \eqref{E:abstract-estimator} be used in {\rm \GALERKIN}. Let Assumptions \ref{A:Lipschitz-estimator} and \ref{A:monotonicity-estimator} be valid.
Let
$\grid_*$ be a refinement of 
$\grid$, with estimator $\eta_{\grid_*}(v,P_{\grid_*}f)$, obtained by bisecting the elements $T \in \marked$ marked in {\rm \MARK}, using a D\"orfler condition on the estimator $\eta_\grid(u_\grid, P_\grid f)$ for the Galerkin solution $u_\grid \in \V_\grid$. 
Suppose that there exists $\lambda \in (0,1)$ such that
\begin{equation}\label{E:estim-red-abstract}
\eta_{\grid_*}(u_\grid,P_\grid f,T)^2 \leq \lambda \,\eta_\grid(u_\grid,P_\grid f,T)^2 \qquad \forall T \in \marked .
\end{equation}
Then, there exists $0 < \rho <1$ and $C>0$ such that
for all $v_{\grid_*} \in \V_{\grid_*}$
\[
\eta_{\grid_*}(v_{\grid_*},P_{\grid_*} f)^2 \leq \rho \, \eta_{\grid}(u_{\grid},P_\grid f)^2 +C \left(\|v_{\grid_*} - u_\grid\|_\V^2 + \|P_{\grid_*} f - P_\grid f \|_{\W^*_{\grid_*}}^2 \right).
\]
\end{proposition}
\begin{proof}
By Assumption \ref{A:Lipschitz-estimator} applied to $\grid_*$, we have 
$$
\eta_{\grid_*}(v_{\grid_*}, P_{\grid_*} f) \leq \eta_{\grid_*}(u_\grid, P_\grid f) + C_\text{Lip} \left( \| v_{\grid_*} - u_\grid\|_\V + \|P_{\grid_*} f - P_\grid f \|_{\W^*_{\grid_*}} \right).
$$
Using Assumption \ref{A:monotonicity-estimator} while extending Proposition \ref{P:est-reduction} to the current abstract setting, we have for any $\delta>0$
\begin{equation*}
\begin{split}
      \eta_{\grid_*}(v_{\grid_*},P_{\grid_*} f)^2 &\le
      (1+\delta)\big(\eta_\grid(u_\grid, P_\grid f)^2 
      - (1-\lambda)\,\eta_\grid(u_\grid, P_\grid f, \marked)^2\big) \\
      & \qquad +
      2(1+\delta^{-1})\, C_\text{Lip}^2 \left( \| v_{\grid_*} - u_\grid\|_\V^2 + \|P_{\grid_*} f - P_\grid f \|_{\W^*_{\grid_*}}^2 \right).
\end{split}      
\end{equation*}  
We conclude using D\"orfler condition $\eta_\grid(u_\grid,\marked) \geq \theta \eta_\grid(u_\grid)$ and choosing $\delta$ small enough.
\end{proof}

Before proceeding further, let us introduce the quantity
$$
\osc_\grid(f):= 
%\Big( \sum_{T \in \grid} \osc_\grid(f,T)_{-1}^2 \Big)^{1/2} \simeq 
\| f - P_\grid f \|_{\W^*_\grid},
$$
which is a measure of the oscillation of the data $f$. If $u_\grid\in\V_\grid$ is the solution of $\AFEMSW$, we 
let ${\cal E}_\grid(u_\grid,f)$ 
indicate the full estimator defined by
\begin{equation}\label{E:def-full-estimator}
{\cal E}_\grid(u_\grid,f)^2 :=\eta_\grid(u_\grid, P_\grid f)^2 + \osc_\grid(f)^2.
\end{equation}
We formulate the following assumption on the data oscillation.
}
%
%\todo[inline]{RHN (12/8/23): Data oscillation is not defined globally but as a sum of local contributions. This is because we do not know that $P_\grid$ is globally stable.}
%
%
\begin{assumption}[quasi-monotonicity of oscillation]\label{A:monotonicity-oscillation} 
\index{Assumptions!Quasi-monotonicity of oscillation}
There exists a constant $C_{\textrm{osc}}> 0$ such that for any $\grid \in \grids$ and any admissible refinement $\grid_* \geq \grid$, we have
 
\[
\osc_{\grid_*}(f) \leq C_{\osc} \osc_{\grid}(f).
\]   
\end{assumption}

\AB{
A consequence of this assumption is the bound
\[
\|P_{\grid_*} f - P_\grid f \|_{\W^*_{\grid_*}} \leq 
\|f- P_{\grid_*} f \|_{\W^*_{\grid_*}} +
\|f - P_\grid f \|_{\W^*_{\grid_*}} \leq (1+C_{\osc}) \osc_{\grid}(f),
\]
which, inserted into the reduction estimate of Proposition \ref{P:local-estimator-reduction}, gives the existence of $0 < \rho_0 <1$ and $C_0>0$ independent of $\grid$ such that
\begin{equation}\label{E:newbound-est-reduct}
\eta_{\grid_*}(u_{\grid_*},P_{\grid_*} f)^2 \leq \rho_0 \, \eta_{\grid}(u_{\grid},P_\grid f)^2 +C_0 \left(\|u_{\grid_*} - u_\grid\|_\V^2 + \osc_{\grid}(f)^2 \right).
\end{equation}
}

We aim at establishing a linear convergence result similar to Theorem \ref{T:Feischl-theorem} for the sequence $\{u_j\}_{j=0}^\infty$ generated by \AFEMSW. To this end, we introduce as usual the short-hand notation $e_j = \|u-u_j\|_\V$, $E_{j+1}=\|u_{j+1}-u_j\|_\V$, \AB{$\eta_j=\eta_{\grid_j}(u_j,P_{\grid_j} f)$, $\osc_j=\osc_{\grid_j}(f)$}, ${\cal E}_j = {\cal E}_{\grid_j}(u_j,f)$, and we also introduce the scaled estimator
\begin{equation}\label{E:def-full-estimator-scaled}
\zeta_j^2 := \eta_j^2 + \gamma \osc_j^2,
\end{equation}
where the parameter $\gamma>0$ is to be found. Note that at this point we have three parametes, $\omega \in (0,1)$, $\xi \in (0,1)$, and $\gamma > 0$ to play with, and the idea is to find conditions on them such that an inequality similar to \eqref{E:linear_convergence_step1.1} in the proof of Theorem \ref{T:Feischl-theorem} holds true. The following result is an intermediate step.
\begin{lemma}[linear estimator reduction]\label{L:lin-estim-red}  Let Assumptions \ref{A:estimator-reduction} (estimator reduction), \ref{A:Lipschitz-estimator} (Lipschitz continuity of estimator), \ref{A:monotonicity-estimator} (monotonicity of estimator) and \ref{A:monotonicity-oscillation} (quasi-monotonicity of oscillation) be valid. There exists \AB{ $\omega_0 >0 $ such that, for any choice of parameters $0 <\omega \leq   \omega_0$ and $0 < \xi \leq 1/\sqrt{2}$} in {\rm \AFEMSW}, 
there exist constants $0 < \rho <1$ , $\Lambda>0$, $\gamma \geq 1$ for which 
\begin{equation}\label{E:lin-estim-red-1}
\zeta_k^2 \leq \rho^{k-j} \zeta_j^2 + \Lambda \sum_{\ell=j+1}^k E_\ell^2, \qquad k \geq j \geq 0.
\end{equation}
\end{lemma}
\begin{proof}
We discuss separately the two alternatives in Algorithm \ref{A:one-step-switch} ($\AFEMSW$).

\step{1} {\em Case $\osc_j \leq \omega {\cal E}_j$}. We use \AB{\eqref{E:newbound-est-reduct} to get
\[
\eta_{j+1}^2 \leq \rho_0 \eta_j^2 + C_0 E_{j+1}^2 + C_0 \osc_j^2
\]
}
and Assumption \ref{A:monotonicity-oscillation} to write
\[
\osc_{j+1}^2 \leq C_{\osc}^2 \osc_j^2.
\]
From
\[
\osc_j^2 \leq \omega^2 {\cal E}_j^2 = \omega^2(\eta_j^2+\osc_j^2) \leq \omega^2(\eta_j^2+ \gamma\osc_j^2) = \omega^2 \zeta_j^2
\]
provided $\gamma \geq 1$, we deduce \AB{
\begin{equation}
\begin{split}    
\zeta_{j+1}^2 &= \eta_{j+1}^2+\gamma\osc_{j+1}^2 \leq \rho_0 \eta_j^2 +C_0 E_{j+1}^2 + (C_0+ \gamma C_{\osc}^2) \osc_j^2 \\ 
&\leq \rho_0 \eta_j^2  +C_0 E_{j+1}^2 +(C_0+\gamma C_{\osc}^2) \omega^2 (\eta_j^2+\gamma \osc_j^2) \\
& = [\rho_0 + (C_0+\gamma C_{\osc}^2)\omega^2]\eta_j^2 +
[(C_0+\gamma C_{\osc}^2)\omega^2]\gamma\osc_j^2 +C_0 E_{j+1}^2 \\
& \leq [\rho_0 + (C_0+\gamma C_{\osc}^2)\omega^2]\zeta_j^2 + C_0 E_{j+1}^2.
\end{split}
\end{equation}
Below we will impose 
\begin{equation}\label{E:lin-estim-red-2}
\rho_1 := \rho_0 + (C_0+\gamma C_{\osc}^2) \omega^2 < 1,
\end{equation}
which will yield the desired bound
\begin{equation}\label{E:lin-estim-red-3}
\zeta_{j+1}^2 \leq \rho_1 \zeta_j^2 + C_0 E_{j+1}^2.
\end{equation}
}

\step{2} {\em Case $\osc_j > \omega {\cal E}_j$}. We use ${\cal E}_j^2=\eta_j^2 + \osc_j^2 > \eta_j^2$ to get
\[
\eta_j^2 < \frac1{\omega^2}\osc_j^2.
\]
\AB{ Proceeding as in the proofs of Proposition \ref{P:local-estimator-reduction} (now with ${\cal M}=\emptyset$) and of \eqref{E:newbound-est-reduct}, we obtain for any $\delta>0$
\[
\eta_{j+1}^2 \leq (1+\delta)\eta_j^2 + C_\delta (E_{j+1}^2+\osc_j^2) \leq (1-\delta)\eta_j^2 + C_\delta E_{j+1}^2 +\left( \frac{2\delta}{\omega^2} +C_\delta \right) \osc_j^2
\]
with $C_\delta = C_0^2(1+\delta^{-1})$.} On the other hand, since $\osc_{j+1}$ is computed after a call to \DATA, it satisfies
\[
\osc_{j+1}^2 \leq \xi^2 \sigma_j^2 = \xi^2 \omega^2 {\cal E}_j^2 < \xi^2 \osc_j^2 =
\frac{1+\xi^2}2\osc_j^2 \ - \  \frac{1-\xi^2}2 \osc_j^2.
\]
Combining the two last equations, we obtain \AB{
\begin{equation*}
\begin{split}
\zeta_{j+1}^2 = \eta_{j+1}^2 + \gamma \osc_{j+1}^2 &\leq (1-\delta)\eta_j^2 +\gamma\frac{1+\xi^2}2\osc_j^2 
 \ +\  C_\delta E_{j+1}^2 \\
 & \qquad \qquad \qquad +\Big(\frac{2\delta}{\omega^2} +C_\delta - \gamma \frac{1-\xi^2}2 \Big) \osc_j^2.
\end{split}
\end{equation*}
Below we will enforce
\begin{equation}\label{E:lin-estim-red-5m1}
\Gamma := \frac{2\delta}{\omega^2} +C_\delta - \gamma \frac{1-\xi^2}2 \leq 0, 
\end{equation}
which will guarantee 
\begin{equation}\label{E:lin-estim-red-5}
\zeta_{j+1}^2 \leq \rho_2 \zeta_j^2 + C_\delta E_{j+1}^2,
\end{equation}
with $\rho_2 := \max(1-\delta, \frac{1+\xi^2}2) < 1$.

\step{3} {\it Choice of parameters:}
Summarizing, in both cases \step{1} and \step{2} we have obtained
\begin{equation}\label{E:lin-estim-red-6}
\zeta_{j+1}^2 \leq \rho \zeta_j^2 + \Lambda E_{j+1}^2
\end{equation}
with $\rho :=\max(\rho_1,\rho_2)<1$ and $\Lambda:=\max(C_0,C_\delta)$, which holds under the conditions \eqref{E:lin-estim-red-2} and \eqref{E:lin-estim-red-5m1}.
Iterating \eqref{E:lin-estim-red-6}, we obtain the desired bound \eqref{E:lin-estim-red-1}.

To fulfil \eqref{E:lin-estim-red-2}, we write $\omega^2$ in the form $\omega^2 = \frac{\sigma_0}\gamma$, which gives 
\[
\rho_1 = \rho_0 + (C_0+\gamma C_{\osc}^2)\frac{\sigma_0}\gamma = \rho_0 + \big(\frac{C_0}\gamma + C_{\osc}^2 \big)\sigma_0 \leq \rho_0 + ({C_0} + C_{\osc}^2 )\sigma_0
\]
since $\gamma \geq 1$, and we pick a $\sigma_0 >0 $ small enough to make $\rho_0 + ({C_0} + C_{\osc}^2 )\sigma_0<1$.

To fulfil \eqref{E:lin-estim-red-5m1}, we use $\xi \leq 1/\sqrt{2}$ and again $\omega^2 = \frac{\sigma_0}\gamma$ to write
\[
\Gamma =  C_\delta + \left(\frac{2\delta}{\sigma_0} -\frac{1-\xi^2}2 \right)\gamma \leq C_\delta + \left(\frac{2\delta}{\sigma_0} -\frac14 \right)\gamma. 
\]
Choosing $\delta=\delta_0 =\sigma_0/16$, yields
\[
\Gamma \leq C_{\delta_0}- \frac18 \gamma \leq 0 \qquad \text{provided} \quad \gamma \geq 8 C_{\delta_0}.
\]
In conclusion, setting $\gamma_0=\max(1, 8 C_{\delta_0})$ and $\omega_0 = \sqrt{\sigma_0/\gamma_0}$, we fulfil both conditions \eqref{E:lin-estim-red-2} and \eqref{E:lin-estim-red-5m1} for any $0 < \omega \leq \omega_0$, by choosing the scaling parameter $\gamma=\sigma_0/\omega^2 \geq \gamma_0$.
This completes the proof.
}
\end{proof}

Before establishing the linear convergence result for Algorithm \ref{A:one-step-switch} (one-step \AFEM with switch), we need to extend Assumption \ref{A:equivalence-estimator} (equivalence of error and estimator) to the present situation, in which the estimator $\eta_\grid$ is replaced by the full estimator ${\cal E}_\grid$ defined in \eqref{E:def-full-estimator}; see Theorem \ref{T:modified-estimator}
(modified residual estimator).

    \begin{assumption}[equivalence of error and full estimator]\label{A:equivalence-estimator-infsup} 
    \index{Assumptions!Equivalence of error and full estimator}
    There are constants $ C_U \geq C_L >0$  such that 
    \begin{equation}\label{E:equivalence-estimator-infsup}
     C_L {\cal E}_j \ \leq \|u - u_j\|_\V \leq C_U {\cal E}_j   \qquad j \ge 0,
    \end{equation}
    where ${\cal E}_j = {\cal E}_{\grid_j}(u_{\grid_j},f)$.
    \end{assumption}

\begin{theorem}[linear convergence for \AFEMSW]\label{T:lin-convergence-switch} Let Assumptions \ref{A:equivalence-estimator-infsup} \ \\ (equivalence of error and full estimator), \ref{A:estimator-reduction} (estimator reduction), \ref{A:Lipschitz-estimator} (Lipschitz continuity of estimator) and \ref{A:monotonicity-oscillation} (quasi-monotonicity of oscillation) be valid.
There exists $\omega_0 \in (0,1]$ such that, for any choice of parameters $0<\omega < \omega_0$ and $0 < \xi, \theta < 1$ in {\rm \AFEMSW}, constants $0 < \rho <1$ and $c>0$ exist for which 
\begin{equation}\label{E:lin-convergence-switch}
e_{j+1} \leq c \rho^i e_j \qquad \forall i,j \in \N,
\end{equation}
where $e_j:=\|u-u_j\|_\V$.
\end{theorem}
\begin{proof}
    By Assumption \ref{A:equivalence-estimator-infsup} and $\gamma \geq 1$ in \eqref{E:def-full-estimator-scaled}, we get the equivalence of error and scaled estimator
\[
\frac{C_L^2}\gamma \zeta_j^2 = \frac{C_L^2}\gamma (\eta_j^2+\gamma \osc_j^2) \leq C_L^2 {\cal E}_j^2 \leq e_j^2 \leq C_U^2 {\cal E}_j^2 \leq C_U^2 (\eta_j^2+\gamma \osc_j^2) = C_U^2 \zeta_j^2.
\]
Invoking \eqref{E:lin-estim-red-1} yields
\begin{equation*}
\begin{split}
 e_k^2 \ \leq \ C_U^2 \zeta_k^2 &\ \leq \ C_U^2 \rho^{k-j} \zeta_j^2 \ + \ C_U^2 \Lambda \sum_{\ell=j+1}^k E_\ell^2, \\
 &\leq \gamma \frac{C_U^2}{C_L^2} \rho^{k-j} e_j^2 \ + \ C_U^2 \Lambda \sum_{\ell=j+1}^k E_\ell^2.
\end{split}
\end{equation*}
This inequality is similar to the expression \eqref{E:linear_convergence_step1.1.01} obtained in Step \step{1} in the proof of Theorem \ref{T:Feischl-theorem}. Therefore, we can finally proceed as in that proof and obtain the desired result.
\end{proof}

\subsubsection{\AB{ Algorithm~\ref{algo:AFEMTS} (\AFEMTS)}}
%whose definition is repeated here for convenience.
% \begin{algo}[two-step AFEM]\label{A:two-step-infsup}
% Given an initial tolerance $\eps_0>0$, a target tolerance $\tol$ and initial mesh $\mesh_0$,
% as well as a safety parameter $\omega \in (0,1]$, {\rm AFEM} consists of the two-step algorithm:
% \medskip
% \begin{algotab}
%   \>  $[\mesh, u_\mesh]=\AFEMTS \, (\mesh_0, \varepsilon_0, \omega, \tol)$
%   \\
%   \>  \> $\text{set } k= 0 \text{ and do }$ 
%   \\
%   \>  \> \> $[\widehat\mesh_{k},\widehat{\data}_{k}]=\DATA \, (\mesh_k, \data, \omega \, \varepsilon_k) $ 
%   \\
%   \>  \> \> $[\mesh_{k+1},u_{k+1}]=\GALERKIN \, (\widehat{\mesh}_{k},\widehat{\data}_{k},\varepsilon_k)$ 
%   \\
%   \>  \> \> $\varepsilon_{k+1}=\tfrac12 {\varepsilon_k}$
%   \\
%   \>  \> \> $k \leftarrow k+1$
%   \\
%   \>  \> $\text{while } \eps_{k-1}>\tol$
%   \\
%   \>  \> $\text{return} \, \mesh_k, u_k$
% \end{algotab}
% \end{algo}
As usual, \DATA produces discrete data $\widehat\data_{k}$ on a mesh $\widehat\mesh_{k} \geq \grid_k$, whereas \GALERKIN produces an approximation $u_{k+1}$ on a mesh $\grid_{k+1} \geq \widehat\mesh_{k}$ to the exact solution $\widehat{u}_k$ of the boundary-value problem of interest with data $\widehat\data_{k}$. Its kernel is given in Algorithm \ref{A:GALERKIN} (\GALERKIN).

\smallskip
In order to proceed, we need some notation and some assumptions. Let us denote by $D(\Omega)$ the space of admissible data $\data$ for the boundary-value problem at hand; let $\|\data\|_{D(\Omega)}$ be a (quasi-)norm on $D(\Omega)$. If $\data$ collects all data of problem \eqref{var-prob-gener}, we write $\B =\B(\data)$ and ${\cal F}={\cal F}(\data) = \dual{f(\data)}{\cdot}$ to highlight the dependence of the bilinear and linear forms on the chosen data; similarly, we write $u=u(\data)$ for the corresponding solution. A perturbation $\widehat{\data}$ of $\data$ generates perturbed bilinear and linear forms $\widehat{\B}=\B(\widehat{\data})$ and $\widehat{\cal F}={\cal F}(\widehat{\data})=\dual{f(\widehat{\data})}{\cdot}$, and a perturbation $\widehat{u}=u(\widehat{\data})$ of $u$, which satisfies 
\begin{equation}\label{prob-var-perturb}
    \widehat{u}\in\V:\qquad \widehat{\B}[\widehat{u},{w}]=\widehat{\cal F}[w]
    \qquad\text{for all } w\in\W.
\end{equation} 

We assume, as in Section \ref{S:DATA}, that a call  $[\widehat{\grid},\widehat{\data}]=\DATA(\grid,\data, \tau)$ generates an admissible refinement $\widehat{\grid}$ of $\grid$ and discrete data $\widehat{\data}$ over $\widehat{\grid}$, such that
\begin{equation}\label{E:bound-data-infsup}
\|\data - \widehat{\data}\|_{D(\Omega)}  \leq \Cdata\tau,  
\end{equation}
where $\Cdata>0$ depends on data (see Section \ref{ss:positivity}).
Finally, we associate to any admissible refinement $\grid$ of $\grid_0$, two finite dimensional spaces $\V_\grid \subset \V$ and $\W_\grid \subset \W$ of equal dimension, made of piecewise polynomial functions on $\grid$ (typically, this is accomplished by choosing a type of finite element compatible with the pair $(\V, \W)$ and adopting it in any $\grid \in \grids$).

We are ready to state the assumptions which will rule our forthcoming analysis of \AFEMTS.

\begin{assumption}[perturbation estimate]\label{A:perturbation inf-sup}
For any $\widetilde{\grid} \in \grids$ and $\varepsilon \leq \varepsilon_0$, let $[\widehat{\grid},\widehat{\data}]=\DATA(\widetilde{\grid},\data, \varepsilon)$ and let $\widehat{u}=u(\widehat{\data})$ be the solution of \eqref{prob-var-perturb}.
There exists a constant $C_D>0$, independent of $\widetilde{\grid}$ and $\varepsilon$,  such that
\begin{equation}\label{E:perturbation inf-sup}
\| u - \widehat{u} \|_\V \leq C_D \|\data - \widehat{\data}\|_{D(\Omega)}.
\end{equation}
\end{assumption}
Note that, concatenating this inequality with
\eqref{E:bound-data-infsup} for $\tau=\varepsilon$, we can quantify the effect of a call to \DATA on the perturbation of the exact solution; we have indeed
\begin{equation}\label{E:perturbation-bound}
\| u - \widehat{u} \|_\V \leq C_D \Cdata\varepsilon.
\end{equation}

\begin{assumption}[uniform continuity constant]\label{A:continuity inf-sup}
For any $\widetilde{\grid} \in \grids$ and $\varepsilon \leq \varepsilon_0$, let $[\widehat{\grid},\widehat{\data}]=\DATA(\widetilde{\grid},\data, \varepsilon)$ and let $\widehat{\B}=\B(\widehat{\data})$ be the associated bilinear form.
There exists a constant $C_B \geq \| \B \|$, independent of $\widetilde{\grid}$ and $\varepsilon$,  such that
\begin{equation}\label{E:continuity inf-sup}
\| \widehat{\B} \| \leq C_B.
\end{equation}
\end{assumption}

\begin{assumption}[uniform inf-sup constant]\label{A:uniform inf-sup}
For any $\widetilde{\grid} \in \grids$ and $\varepsilon \leq \varepsilon_0$, let $[\widehat{\grid},\widehat{\data}]=\DATA(\widetilde{\grid},\data, \varepsilon)$, let $\widehat{\B}=\B(\widehat{\data})$ be the associated bilinear form, let $\grid$ be either $\widehat{\grid}$ or an admissible refinement of $\widehat{\grid}$, and finally let $\V_{\grid} \subset \V$, $\W_{\grid} \subset \W$ be the subspaces built on $\grid$ as above.
There exists a constant $0<\bar{\beta} \leq \beta$, independent of $\widetilde{\grid}$, $\varepsilon$, and $\grid$,  such that
\begin{equation}\label{E:uniform inf-sup}
\infimum_{w\in\W_{\grid}}\sup_{v\in \V_\grid} \frac{\widehat{\B}[v,w]}{\|v\|_\V \|w\|_\W}
      \geq \bar{\beta}. 
\end{equation}
\end{assumption}

The last assumption guarantees the well-posedness of all the discrete variational problems 
\begin{equation}\label{prob-var-perturb-discr}
    u_\grid\in\V_\grid:\quad \widehat{\B}[u_\grid,{w}_\grid]=\wh{\cal F}[w_\grid]
    \quad\forall \, w_\grid\in\W_\grid,
\end{equation}
associated with the successive refinements of the initial mesh $\grid_0$ performed by \AFEMTS.

We want to prove, 
as in the coercive case (see Proposition \ref{P:cost-galerkin}), that the number of iterations performed in any call to \GALERKIN inside \AFEMTS (which is finite by Proposition \ref{P:term-GAL}) is indeed uniformly bounded.

\begin{proposition}[computational cost of $\GALERKIN$]\label{P:cost-galerkin-infsup}
~~ Let Assumptions \ref{A:equivalence-estimator}, \ref{A:estimator-reduction}, \ref{A:perturbation inf-sup}, \ref{A:continuity inf-sup}, and \ref{A:uniform inf-sup}
be valid.
For any $k\in\N$, the number of subiterations $J_k$ inside a call to {\rm \GALERKIN} at
iteration $k$ of {\rm \AFEMTS} is bounded by a constant $J$ independent of $k$.
\end{proposition}
\begin{proof} 
% The proof is similar to that of Proposition \ref{P:cost-galerkin}, but now we use the linear convergence of the Galerkin error 
%  (Theorem \ref{T:Feischl-theorem}) rather than the contraction of the quasi-error.
Denote by $\grid_{k,j}$ the successive refinements of $\widehat{\grid}_k$ defined in \GALERKIN at iteration $k$, and let $u_{k,j} \in \V_{k,j}=\V_{\grid_{k,j}}$ be the corresponding Galerkin solutions, which are approximations of the solution $\widehat{u}_k \in \V$ of the perturbed problem \eqref{prob-var-perturb} with forms $\widehat{\B}=\widehat{\B}_k=\B(\widehat{\data}_k)$ and $\widehat{f}=\widehat{f}_k=f(\widehat{\data}_k)$.
Note as well that we use a posteriori estimators $\eta_{k,j}=\eta_{k,j}(v)$ defined on $\V_{k,j}$, which depend on $\wh{\data}_k$ via the coefficients of the equation. However, in reference to Assumptions \ref{A:equivalence-estimator} and \ref{A:estimator-reduction}, we always suppose the constants in the bounds \eqref{E:equivalence-estimator} and \eqref{E:estimator-reduction} to be  independent of $k$ and $j$.

Let us pick $j:=J_k-1$. By definition of stopping criterion in \GALERKIN, and by \eqref{E:equivalence-estimator} and \eqref{E:Feischl-theorem}, we 
get
\begin{equation}\label{E:Gal-cost-infsup-1}
\varepsilon_k < \eta_{k,j}(u_{k,j}) \leq \frac1{C_L} \| \widehat{u}_k - u_{k,j} \|_\V \leq \frac{c}{C_L} \rho^j \| \widehat{u}_k - u_{k,0} \|_\V. 
\end{equation}
The norm on the right-hand side can be bounded via Corollary \ref{C:quasi-monotonicty} (quasi-monotoni-city), applied to $\B:=\widehat{\B}_k$, $u:=\widehat{u}_k \in \V$, $u_N:=u_{k,0} \in \V_N:=\V_{\wh{\grid}_k}$, and $v: =u_k \in \V_M :=\V_{\grid_k} \subseteq \V_{\wh{\grid}_k}$ (the output of \GALERKIN at iteration $k-1$).  Using Assumptions \ref{A:continuity inf-sup} (uniform continuity constant) and \ref{A:uniform inf-sup} (uniform inf-sup constant),
we thus have
\begin{equation}\label{E:Gal-cost-infsup-2}
\| \widehat{u}_k - u_{k,0} \|_\V \leq 
\lambda \| \widehat{u}_k - {u}_k  \|_\V, 
\end{equation}
with $\lambda=\frac{C_B}{\bar{\beta}}$.
%
% \todo[inline]{RHN (01/21/24): Claudio, the quote above to Lemma \ref{L:quasi-monotonicity} seems incorrect. Moreover, this lemma does not seem the statement we need, which is
% %
% \[
% \| \wh{u}_k - u_{k,0} \|_\V \le \Big(1+\frac{\|\B\|}{\beta}\big) \|\wh{u}_k - u_k\|_\V
% \]
% %
% This is a consequence of the inf-sup over $\V_{k,0}$ and the fact that $u_k\in\V_k\subset\V_{k,0}$. Do you agree? If so, then maybe this can be written explicitly within the proof or before this proposition. What do you think? 

% Note that $u_k$ corresponds to data $\wh{\data}_{k-1}$ whereas $\wh{u}_k$ is the exact solution for data $\wh{\data}_k$, but this does not matter. This DOES matter for dG though.
% }
At last, we use again the triangle inequality to get
\begin{eqnarray*}
 \begin{split}
\| \widehat{u}_k - {u}_k\|_\V &\leq  \| \widehat{u}_k - \widehat{u}_{k-1}  \|_\V + \| \widehat{u}_{k-1} - {u}_k  \|_\V \\
& \leq \| u- \widehat{u}_k \|_\V + \| u- \widehat{u}_{k-1} \|_\V + \|\widehat{u}_{k-1} - u_k \|_\V;
\end{split}
\end{eqnarray*}
then, Assumption \ref{A:perturbation inf-sup} (perturbation estimate) yields 
$\| u- \widehat{u}_k \|_\V \leq C_D \Cdata\, \omega \varepsilon_k$ and $\| u- \widehat{u}_{k-1} \|_\V \leq C_D\Cdata \, \omega \varepsilon_{k-1}$, whereas the termination test for \GALERKIN at iteration $k-1$ yields $\|\widehat{u}_{k-1} - u_k \|_\V \leq C_U \, \varepsilon_{k-1}$. Hence, recalling $\varepsilon_{k-1}=2\varepsilon_k$ and $\omega \leq 1$, we get
\begin{equation}\label{E:Gal-cost-infsup-3}
\| \widehat{u}_k - {u}_k\|_\V \leq \sigma \, \varepsilon_k
\end{equation}
with $\sigma = 3C_D\Cdata + 2C_U$. Finally, concatenating \eqref{E:Gal-cost-infsup-1}, \eqref{E:Gal-cost-infsup-2} and \eqref{E:Gal-cost-infsup-3}, we obtain
\[
\rho^j \geq \frac{C_L}{c \lambda \sigma},
\]
which implies $J_k \leq 1 + \big(\log \frac{C_L}{c \lambda \sigma}\big) (\log \rho)^{-1} =: J$.
\end{proof}

The remaining of this section is devoted to investigate the rate-optimality of {\AFEM}s for inf-sup stable problems. Precisely, we aim at establishing the analogue of bound \eqref{E:optimality-bound} for such problems, i.e.,
\begin{equation}\label{E:optimality-bound-infsup}
\|u-u_\grid\|_{\V} \leq C(u,\data) \, \big( \#\mesh \big)^{-s} \,.
\end{equation}
To this end, we have to introduce approximation classes for the solution and the data, and to study the quasi-optimality properties of mesh refinement and \GALERKIN.

%----------------------------------------------------------------------------------
\subsubsection{Nonlinear approximation classes}
%----------------------------------------------------------------------------------

The definition of the approximation class $\A_s=\A_s(\V;\grid_0)$ for functions in $\V$ is identical to that given in Section \ref{S:approx-solution} for functions in $H^1_0(\Omega)$ (see Definition \ref{D:approx-class-u}), provided the norm $|v|_{H^1_0(\Omega)}$ is replaced by the norm $\| v\|_\V$ at all occurrences.

In the rest of the section we will make the following regularity assumption.

\begin{assumption}[approximability of $u$]\label{A:approx-u-infsup}
The exact solution $u\in \V$ of problem \eqref{var-prob-gener} belongs to the
approximation class $\As(\V;\grid_0)$ for some $s=s_u \in (0,\frac{n}{d}]$. 
\end{assumption}

The approximation classes of a data $\data \in D(\Omega)$ are defined via discrete approximations $\data_\grid \in \D_\grid$ subordinate to a partition $\grid \in \grids$, which produce the oscillation
\[
\osc_\grid(\data)= \infimum_{\data_\grid \in \D_\grid} \| \data - \data_\grid \|_{D(\Omega)}.
\]

\begin{definition}[approximation classes of $\data$] Let
$\D_s := \D_s \big(D(\Omega);\grid_0\big)$ be the set of data $\data \in D(\Omega)$
satisfying
\begin{equation}\label{E:osc-D}
| \data |_{\D_{s}}:=
\sup_{N \geq \# \grid_0}
\left( N^s\inf _{\grid \in \mathbb{T}_N} \osc_\mesh (\data) \right) < \infty
\,\,\Rightarrow\,\,
 \inf _{\grid \in \mathbb{T}_N}\osc_\grid (\data) \le | \data |_{\D_{s}} N^{-s}.
\end{equation}
\end{definition}

The following assumptions on the data of our boundary-value problem will be valid in the rest of the section.

\begin{assumption}[approximability of $\data$]\label{A:approx-data-infsup}
The data $\data\in D(\Omega)$ of problem \eqref{var-prob-gener} belongs to the
approximation class $\D_s\big(D(\Omega);\grid_0\big)$ for some $s=s_\data \in (0,\frac{n}{d}]$. 
\end{assumption}

\begin{assumption}[quasi-optimality of $\DATA$]\label{A:optim-data-infsup}
A call $[\wh\grid,\wh\data] = {\tt DATA} \, (\grid,\data,\varepsilon)$ marks a set of elements
${\cal M}_{\data}$ whose cardinality $N(\data)=\#{\cal M}_{\data}$ obeys
\begin{equation}\label{e:assumption_marked_data}
N(\data) \lesssim | \data |_{\D_s}^{\frac{1}{s_\data}} \, \varepsilon^{-\frac{1}{s_\data}} 
\, .
\end{equation}
\end{assumption}

The concept of {\em $\varepsilon$-approximation of order $s$} of $u\in \A_s(\V;\grid_0)$ is identical to the one given in Definition \ref{D:approx-order-s}, and so is the proof of the following result. %which will prove useful in the sequel.
\begin{lemma}[$\varepsilon$-approximation of $u$ of order $s$]\label{L:eps-approx-infsup}
Let $u\in \As(\V;\grid_0)$ and $v\in \V$ satisfy $\| u-v \|_{\V}\leq \varepsilon$ for some $0<\varepsilon\le\eps_0$. Then $v$ is a $2\varepsilon$-approximation of order $s$ to $u$.  
\end{lemma}

\subsubsection{Rate-optimality of {\rm \GALERKIN}}
To estimate the growth of the cardinality of the meshes produced inside a call to \GALERKIN, 
which always deals with discrete data, and to relate it to the approximation class of the exact solution $u$, we need an additional assumption of the estimators $\eta_\grid$. In the sequel, for any subset ${\cal R}\subset \grid$, we define $\eta_\grid(v,{\cal S})$ by
\[
\eta_\grid(v,{\cal S})^2 = \sum_{T \in {\cal S}} \eta_\grid(v,T)^2.
\]

\begin{assumption}[discrete reliability of the estimator]\label{A:discrete-rel-estimator}
There exists a constant $c_2>0 $ such that for any $\grid \in \grids$ and any refinement $\grid_* \geq \grid$, if ${\cal R}={\cal R}_{\grid \to \grid_*}=\grid \setminus\grid_*$ is the set of refined elements of $\grid$, it holds
\[
\| u_{\grid_*}-u_\grid\|_\V \leq c_2 \eta_\grid(u_\grid,{\cal R}),
\]
where $u_\grid$ and $u_{\grid_*}$ (resp.) are the Galerkin solutions in $\V_\grid$ and $\V_{\grid_*}$ (resp.). 
\end{assumption}

We recall that the module \MARK in \GALERKIN implements D\"orfler's strategy, i.e., for a fixed $\theta \in (0,1]$, it identifies a subset ${\cal M} \subseteq \grid$ of elements undergoing bisection by the condition
\begin{equation}\label{E: doerfler-infsup}
\eta_\grid(u_\grid,{\cal M}) \geq \theta \,\eta_\grid(u_\grid).
\end{equation}

The following property is the analogue of the one stated in Lemma \ref{L:dorfler-no-osc} for coercive problems. Since the proof is similar, we omit it.

%\todo[inline]{RHN (12/28/23): Lemma \ref{L:dorfler-no-osc-infsup} is now included in Lemma \ref{L:dorfler-no-osc}, and may be erased or briefly recalled.}

\begin{lemma}[D\"orfler marking]\label{L:dorfler-no-osc-infsup}
Let Assumptions \ref{A:Lipschitz-estimator} and \ref{A:discrete-rel-estimator} be valid.
For all $0<\mu<\frac12$ there exists $0<\theta_0<1$ such that, if $\grid \in \grids$ and $\grid_*$ is a refinement of $\grid$ with refined set ${\cal R}=\grid \setminus \grid_*$, and if the Galerkin solutions $u_\grid \in \V_\grid$ and $u_{\grid_*} \in \V_{\grid_*}$ satisfy
\[
\eta_{\grid_*}(u_{\grid_*}) \leq \mu \, \eta_{\grid}(u_{\grid}),
\]
then
\[
\theta_0 \,\eta_\grid(u_\grid) \leq \eta_\grid(u_\grid,{\cal R}). 
\]
%
%for all $0 < \theta \leq \theta_0$.
\end{lemma}
%
% \begin{proof}
%  Write
%  \[
%  \eta_\grid^2(u_\grid) = \eta_\grid^2(u_\grid,{\cal R}) + \eta_\grid^2(u_\grid, \grid \cap \grid_*)
%  \]
% and note that, using the Assumptions,
% \begin{equation*}
% \begin{split}
% \eta_\grid^2(u_\grid, \grid \cap \grid_*) &\leq (1+\delta)\, \eta_\grid^2(u_{\grid_*}, \grid \cap \grid_*) + c_1^2\, (1+\delta^{-1}) \,\| u_{\grid_*}-u_\grid\|_\V \\
% & \leq (1+\delta)\, \mu^2 \, \eta_\grid^2(u_\grid) +
% c_1^2(1+\delta^{-1}) \, c_2^2 \,\eta_\grid^2(u_\grid, {\cal R}).
% \end{split}
% \end{equation*}
% %
% \todo[inline]{RHN (12/26/23): It seems to me that we do not have to split into $\grid\cap\grid_*$ and $\refined$. In view of Proposition 5.36, how about saying:

% $$\eta_\grid(u_\grid)^2 \le (1+\delta) \eta_{\grid_*} (u_{\grid_*})^2 + C_1^2 (1+\delta^{-1}) \|u_\grid - u_{\grid_*}\|_V^2?$$

% Moreover, the upper bound of $\| u_{\grid_*}-u_\grid\|_\V$ involves also $\osc_\grid(f,\cal R)$ that may
% not vanish unless $f$ is discrete.}
% %
% This implies 
% \[
% \big(1-(1+\delta)\,\mu^2 \big) \,\eta_\grid^2(u_\grid) \leq c_1^2 c_2^2 \, (1+\delta^{-1})\, \eta_\grid^2(u_\grid, {\cal R}),
% \]
% or equivalently
% \[
% \frac{1-(1+\delta)\,\mu^2}{c_1^2 c_2^2 \, (1+\delta^{-1})} \, \eta_\grid^2(u_\grid) \leq  \eta_\grid^2(u_\grid, {\cal R}).
% \]
% Thus, we may set
% \begin{equation}\label{E:def-theta0}
% \theta_0^2 := \frac{1-(1+\delta)\,\mu^2}{c_1^2 c_2^2 \, (1+\delta^{-1})}
% \end{equation}
% and choose $\delta$ small enough so that $0 <\theta_0 <1$. 
% \end{proof}

We are ready to investigate the rate-optimality of the $k$-th call to \GALERKIN in the two-step \AFEM (see Definition \ref{algo:AFEM-Kernel}). We denote by $\marked_{k,j} \subseteq \grid_{k,j}$ the marked set at the $j$-th iteration inside \GALERKIN (hereafter, we refer to the notation in the proof of Proposition \ref{P:cost-galerkin-infsup}). 
To achieve quasi-optimality, the following assumption is fundamental.
\begin{assumption}[minimality of marked sets]\label{A:cardinality-infsup}
The module $\MARK$ selects a set $\marked_{k,j}$ with \emph{minimal} cardinality among those satisfying D\"orfler's condition
\[
\eta_{k,j}(u_{k,j},{\cal M}) \geq \theta \,\eta_{k,j}(u_{k,j}) \quad\forall \, k,j.
\]
\end{assumption}

% \todo[inline]{RHN (12/29/23): Please check Proposition \ref{P:cardinality-marked} and Corollary
% \ref{C:cardinality-no-osc}. They might greatly simplified these results because the proofs are very similar now.}

\begin{proposition}[cardinality of marked sets]\label{P:optim-GAL-infsup}
Let Assumptions \ref{A:equivalence-estimator}, \ref{A:perturbation inf-sup},
\ref{A:continuity inf-sup}, \ref{A:uniform inf-sup}, \ref{A:approx-u-infsup}, \ref{A:Lipschitz-estimator}, \ref{A:discrete-rel-estimator} and \ref{A:cardinality-infsup}  hold true. There exists a constant $C_0>0$ independent of $k$ and $j$ such that the cardinality $N_{k,j}(u)$ of $\marked_{k,j}$ satisfies
\begin{equation}
N_{k,j}(u) \leq C_0 \, |u|_{\A_s}^{1/s} \varepsilon_k^{-1/s} 
\end{equation}
and
\begin{equation}
N_{k,j}(u) \leq C_0 \, |u|_{\A_s}^{1/s} \| u - u_{k,j}\|_\V^{-1/s}.
\end{equation}    
\end{proposition}
\begin{proof}
The proof can be easily obtained by slightly adapting to the current abstract setting the proof of Corollary \ref{C:cardinality-no-osc}, keeping also into account Proposition \ref{P:cardinality-marked}.
\end{proof}

Let $\marked_k$ denote the set of marked elements in \GALERKIN at iteration $k$ of \AFEM. Since the cardinality $N_k(u)= \#\marked_k$ of $\marked_k$ satisfies $N_k(u) = \sum_{j=0}^{J_k-1} N_{k,j}(u)$, we can estimate its cardinality by combining Propositions \ref{P:cost-galerkin-infsup} and \ref{P:optim-GAL-infsup}.  

\begin{corollary}[rate-optimality of \GALERKIN]\label{C:opt-GAL-infsup} Under the assumptions of Propositions \ref{P:cost-galerkin-infsup} and \ref{P:optim-GAL-infsup}, the total number of marked elements $\marked_k$ in {\rm \GALERKIN} at iteration $k$ of {\rm \AFEM} satisfies
\[
N_{k}(u) \leq J \, C_0 \, |u|_{\A_s}^{1/s} \varepsilon_k^{-1/s}.  
\]   
\end{corollary}

%------------------------------------------------------------------------------------
\subsubsection{Rate-optimality of {\rm \AFEMTS}}
%------------------------------------------------------------------------------------
At last, we focus on the two-step {\rm \AFEM} in Definition \ref{algo:AFEM-Kernel} ($\AFEMTS$), and prove its rate-optimality, in relation to the nonlinear approximation classes of the solution $u$ and the problem data $\data$.

%\todo[inline]{RHN (12/29/23): This theorem might be a simple consequence of Theorem \ref{T:optimality-AFEM}. Please check.}

\begin{theorem}[rate-optimality of \AFEMTS]\label{T:opt-AFEM-infsup}
Under the same assumptions of Proposition \ref{P:optim-GAL-infsup}, plus Assumptions \ref{A:approx-data-infsup} and \ref{A:optim-data-infsup}, there exists a constant $C_*$ independent of $u$ and $\data$ such that the sequence $(\grid_k, \V_{\grid_k},u_{\grid_k})$, $k \geq 0$, produced by {\rm \AFEMTS} satisfies
\[
\| u - u_{\grid_k} \|_\V \leq C_* \left( |u|_{\A_{s_u}}^{1/s_u}+ |\data|_{\D_{s_\data}}^{1/s_\data} \right)^s \big(\# \grid_k\big)^{-s},
\]
with $0<s=\min(s_u,s_\data) \leq \frac{n}d$.
\end{theorem}
\begin{proof}
Let us denote by $\marked_\ell^u$ ($\marked_\ell^\data$, resp.) the set of elements marked by \GALERKIN (by \DATA, resp.) at iteration $\ell$ of \AFEM. By Corollary \ref{C:opt-GAL-infsup} and Assumption \ref{A:optim-data-infsup}, there exist constants $D_1, D_2$ independent of $u$, $\data$ and $k$ such that
\[
N_\ell(u) \leq D_1 \, |u|_{\A_{s_u}}^{1/s_u} \varepsilon_\ell^{-1/s_u}, \qquad 
N_\ell(\data) \leq D_2 \, | \data |_{\D_{s_\data}}^{1/s_\data} \, \varepsilon_\ell^{-1/s_\data} .
\]
Then, we conclude as in the proof of Theorem \ref{T:optimality-AFEM}.
% Hence, there exists $D_3$ such that
% \[
% \#\grid_{k} - \#\grid_0 \leq D_3 \left( |u|_{\A_{s_u}}^{1/s_u} + | \data |_{\D_{s_\data}}^{1/s_\data} \right)\sum_{\ell=0}^{k-1} \varepsilon_\ell^{-1/s}.
% \]
% Now we use $\sum_{\ell=0}^{k-1} \varepsilon_\ell^{-1/s} \simeq \varepsilon_k^{-1/s}$ and $\| u - u_{\grid_k} \|_\V \lesssim \varepsilon_{k-1} \lesssim \varepsilon_k$ to get the existence of $D_4$ such that
% \[
% \#\grid_{k} - \#\grid_0 \leq D_4 \left( |u|_{\A_{s_u}}^{1/s_u} + | \data |_{\D_{s_\data}}^{1/s_\data} \right) \| u - u_{\grid_k} \|_\V^{-1/s}.
% \]
% We conclude by observing that it is not restrictive to assume  $\#\grid_k \geq c \, \#\grid_0$ for some $c>1$ independent of $k$.
\end{proof}

\subsection{The Stokes problem}\label{S:rates-stokes}
%-----------------------------------------------------------------------------------------

Here we consider the Stokes problem
\begin{equation}\label{strong-stokes-1}
\begin{alignedat}{2}
-\Delta \bu + \nabla p &= \bF &\qquad&\text{in }\Omega,
\\
\div\bu & = 0 &&\text{in }\Omega,
\\
\bu &= 0 &&\textrm{on }\partial\Omega,
\end{alignedat}
\end{equation}
already introduced in Section~\ref{S:ex-bvp}. Assuming $\vec{f} \in H^{-1}(\Omega;\R^d)=\V^*$, its weak formulation is given in \eqref{weak-stokes} or, equivalently, in \eqref{weak-stokes-B}, where the bilinear form $\B$ is continuous and inf-sup stable, as a consequence of Brezzi's Theorem \ref{T:Brezzi} (Brezzi); see Section~\ref{S:inf-sup}.

%\todo[inline]{Should we discuss specific mixed FEMs (Taylor-Hood and mini element) here or we present theanalysis abstractly? Perhaps we could do the latter here and discuss examples in Section \ref{S:aposteriori}.} 

%\todo[inline]{Should we discuss stabilized FEMs? Is there a connection with dG? What is the state-of-the-art regarding a posteriori error analysis for this type of methods? This question belongs to Section \ref{S:aposteriori}.}

A Galerkin discretization of this problem, based on finite dimensional subspaces $\V_\grid \subset \V=H^1_0(\Omega;\R^d)$ and $\Q_\grid \subset \Q=L^2_0(\Omega)$, reads as follows: 
find $(\bu_\grid,p_\grid)\in\V_\grid\times\Q_\grid$ such that 
\begin{equation}\label{E:weak-stokes-Gal}
\begin{alignedat}{3}
  &a[\vec{u}_\grid,\vec{v}] + b[p_\grid,\vec{v}] &&= \scp{\vec{f}}{\vec{v}} 
&\qquad &\forall \, \vec{v}\in\V_\grid,
  \\
  &b[q,\vec{u}_\grid] &&= 0 &&\forall \, q\in\Q_\grid,
\end{alignedat}
\end{equation}
or equivalently,
\begin{equation*}
  (\vec{u}_\grid,p_\grid)\in \V_\grid\times\Q_\grid: \quad
  \bilin{(\vec{u}_\grid,p_\grid)}{(\vec{v},q)} =  \scp{\vec{f}}{\vec{v}} \quad\forall \,
  (\vec{v},q)\in \V_\grid\times\Q_\grid.
\end{equation*}
We assume that the pair $(\V_\grid,\Q_\grid)$ is {\em uniformly inf-sup stable} for the form $b$, i.e., there exists a constant $\beta>0$, independent of the refinement $\grid$, such that
\begin{equation}\label{E:unif-infsup}
\infimum_{q\in\Q_\grid}\sup_{\vec{v}\in \V_\grid}\frac{b[q,\vec{v}]}{\|\vec{v}\|_\V \|w\|_\W}
  \geq \beta \,.
\end{equation}  
This condition is equivalent to the uniform inf-sup stability of the bilinear form $\B$ on the product space $\mathbb X_\grid :=\V_\grid \times \Q_\grid$.
Then, applying a discrete form of Brezzi's Theorem, we obtain the existence and uniqueness of the solution of \eqref{E:weak-stokes-Gal}, which satisfies the stability bound
\begin{equation}\label{E:stability-Stokes}
\| \vec{u}_\grid \|_\V + \| p_\grid \|_\Q \leq C \, \| \vec{f} \|_{\V^*}\,,  
\end{equation}
where $C$ only depends on the continuity constant $\|a\|$ and the coercivity constant $\alpha$ of the form $a$, and the inf-sup constant $\beta$. Furthermore, we have the quasi-best approximation bounds (\cite[Proposition 8.2.1]{BoffiBrezziFortin})
\begin{eqnarray}
\| \vec{u}-\vec{u}_\grid \|_\V &\leq & C_{11} \, \min_{\vec{v} \in \V_\grid} \| \vec{u}-\vec{v} \|_\V \ + \ C_{12}\, \min_{q \in \Q_\grid} \| p-q \|_\Q \\
\| p-p_\grid \|_\V &\leq& C_{21} \, \min_{\vec{v} \in \V_\grid} \| \vec{u}-\vec{v} \|_\V \ + \ C_{22}\, \min_{q \in \Q_\grid} \| p-q \|_\Q
\end{eqnarray}
where the constants $C_{ij}$, $1 \leq i,j \leq 2$, only depend on the quantities $\|a\|$, $\|b\|$, $\alpha$ and $\beta$.

There are many families of finite element spaces that are uniformly inf-sup stable for the Stokes problem; see \cite[Chapter 8]{BoffiBrezziFortin}. Among them, we consider  here the Taylor-Hood element \cite{TaylorHood:73} and its higher-order versions. They all use continuous discrete pressures, so they fit in the general form
\begin{eqnarray*}
\V_\grid &=&\{\vec{v} \in H^1_0(\Omega; \R^d) : \  \vec{v}|_{T} \in \vec{V}_T, \ T \in \grid\} \\
\Q_\grid &=&\{ q \in L^2_0(\Omega)\cap C^0(\bar{\Omega}) : \ q|_{T} \in Q_T, \ T \in \grid\},
\end{eqnarray*}
where $\vec{V}_T$ and $Q_T$ are spaces of polynomials on the element $T$. Considering simplicial elements, we have for $n\geq 2$
\begin{equation}
\vec{V}_T = ({\mathbb P}_n(T))^d, \qquad Q_T = {\mathbb P}_{n-1}(T).
\end{equation}

\smallskip
The convergence and optimality of an adaptive algorithm for the Stokes problem based on Taylor-Hood elements was first established by Feischl \cite{Feischl:2019} (see also \cite[Section 6]{Feischl:2022}). We aim at deriving a similar result using the abstract framework presented in this section.

We start by developing the a posteriori error analysis, and for this we introduce the weak residual 
\[
\langle \boldsymbol{{\cal R}}_\grid, (\vec{v},q) \rangle := \langle \vec{f}, \vec{v} \rangle - \bilin{(\vec{u}_\grid,p_\grid)}{(\vec{v},q)} \qquad \forall (\vec{v},q)\in \V \times \Q,
\]
which we represent as $\boldsymbol{{\cal R}}_\grid = (\boldsymbol{{\cal R}}_\grid^m, {{\cal R}}_\grid^c)$ according to the momentum and continuity equations; note that ${{\cal R}}_\grid^c = \div\bu_\grid$. The continuity and inf-sup stability properties of the exact Stokes form ${\cal B}$ yield the equivalence
\begin{equation}\label{E:equiv-stokes}
\| \vec{u} - \vec{u}_\grid \|_\V + \| p - p_\grid \|_\Q \approx \| \boldsymbol{{\cal R}}_\grid \|_{\V^* \times \Q^*} \approx 
\| \boldsymbol{{\cal R}}_\grid^m \|_{\V^*} + \| \div\bu_\grid \|_{L^2(\Omega)}.
\end{equation}
We now apply Corollary \ref{C:residual-loc} (star localization of residual norm) to each component of the momentum residual $\boldsymbol{{\cal R}}_\grid^m$ to get
\[
\| \boldsymbol{{\cal R}}_\grid^m \|_{\V^*}^2 \approx \sum_{z \in {\cal V}} \| \boldsymbol{{\cal R}}_\grid^m \|_{(H^{-1}(\omega_z))^d}^2,
\]
whereas Lemma \ref{L:abstract_bds} 
(splitting of local residual norm) yields the equivalence
\[
\| \boldsymbol{{\cal R}}_\grid^m \|_{(H^{-1}(\omega_z))^d}^2 \approx \| P_\grid \vec{f} + \Delta \vec{u}_\grid - \nabla p_\grid \|_{(H^{-1}(\omega_z))^d}^2 + \| \vec{f}-P_\grid \vec{f} \|_{(H^{-1}(\omega_z))^d}^2.
\]
In view of mesh refinement, we recall Lemma \ref{L:local-re-indexing} (localization re-indexing) and we express the error indicator in terms of elements $T \in \grid$ rather than stars $\omega_z$, in analogy with the scalar case \eqref{mod-res-est}. To this end, define
\begin{equation}
\label{E:mod-res-est-stokes}
\begin{aligned}
\eta_\grid^m(\vec{u}_\grid,p_\grid,T)^2
 &\definedas
 h_T^2 \, \| P_T \vec{f} + \Delta \vec{u}_\grid - \nabla p_\grid\|_{(L^2(T))^d}^2  
\\
 &\qquad
+ h_T \!\!\! \sum_{F \subset \partial T \setminus \partial\Omega}
  \| \jump{(\nabla \vec{u}_\grid) \vec{n}_F} - P_F \vec{f} \|_{(L^2(F))^d}^2 \,,  
\\
 \osc_\grid(\vec{f},T)_{-1}^2
 &\definedas
  \| \vec{f} - P_\grid \vec{f} \|_{(H^{-1}(\omega_T))^d}^2. 
\end{aligned}
\end{equation}
Note that the jump term $\jump{(\nabla \vec{u}_\grid) \vec{n}_F}$ does not contain the pressure contribution, since discrete pressures in $\Q_\grid$ are globally continuous. We thus have
\begin{equation}\label{E:bound-Resm-stokes}
\| \boldsymbol{{\cal R}}_\grid^m \|_{\V^*}^2 \approx \sum_{T \in \grid} 
\Big(\eta_\grid^m(\vec{u}_\grid,p_\grid,T)^2 + \osc_\grid(\vec{f},T)_{-1}^2   \Big).
\end{equation}
Recalling \eqref{E:equiv-stokes}, the full local PDE residual indicator could be defined as
\[
\eta_\grid(\vec{u}_\grid,p_\grid,T)^2 :=
\eta_\grid^m(\vec{u}_\grid,p_\grid,T)^2 + \| \div\bu_\grid \|_{L^2(T)}^2.
\]
However, such a quantity is not guaranteed to strictly reduce under mesh refinement, due to the presence of the divergence term, which is not scaled by a positive power of the meshsize. The following result provides an equivalent expression of the last term, which does reduce.
We recall the definition \eqref{e:average_jump} of jumps across faces.
\begin{lemma}[norm equivalence for divergence]\label{L:div-jumps} It holds
\begin{equation*}
\| \divo \vec{u}_\grid \|_{L^2(\Omega)}^2 \ \approx \
  \sum_{T \in \grid} \sum_{F \subset \partial T\setminus \partial\Omega} h_F \, \| \jump{\divo \vec{u}_\grid}\|_{L^2(F)}^2 . 
\end{equation*}    
\end{lemma}
\begin{proof}
The result follows from applying to $\varphi=\divo \vec{u}_\grid$ the equivalence
\[
\| \varphi - \Pi_\grid \varphi \|_{L^2(\Omega)}^2 \approx  \sum_{T \in \grid} \sum_{F \subset \partial T\setminus \partial\Omega} h_F \, \| \jump{\varphi}\|_{L^2(F)}^2 \qquad \forall \varphi \in 
{\mathbb S}^{n-1,-1}_\grid
\]
(where $\Pi_\grid$ is the $L^2$-orthogonal projection upon ${\mathbb S}^{n-1,0}_\grid$), after observing that $\Pi_\grid \varphi = 0$ since
$\vec{u}_\grid$ is discretely divergence-free, i.e., it satisfies the second set of equations in \eqref{E:weak-stokes-Gal}). To prove the equivalence for arbitrary $\varphi \in {\mathbb S}^{n-1,-1}_\grid$, we use the quasi-interpolation operator $\IdG$ introduced in Section \ref{S:IdG}, which leaves ${\mathbb S}^{n-1,0}_\grid$ invariant. Then, it is easily seen that
\[
\| \varphi - \Pi_\grid \varphi \|_{L^2(\Omega)}^2 \approx \| \varphi - \IdG \varphi \|_{L^2(\Omega)}^2 ,
\]
so, it is enough to prove the equivalence with $\Pi_\grid$ replaced by $\IdG$. But this calculation can be done on patches $\omega_T$ since $\IdG$ is quasi-local:
\[
\| \varphi - \IdG \varphi \|_{L^2(T)}^2 \lesssim  \sum_{F \subset \omega_T} h_F \, \| \jump{\varphi}\|_{L^2(F)}^2 \lesssim \| \varphi - \IdG \varphi \|_{L^2(\omega_T)}^2
\quad \forall \varphi \in 
{\mathbb S}^{n-1,-1}_\grid.
\]
The first inequality follows from 
\eqref{e:dg_interpol_nonconfB}
(see also \cite[Proposition 5.4]{BanschMorinNochetto:02}). The second inequality results from the fact that if the rightmost term vanishes on $\omega_T$,
then $\varphi = \IdG \varphi$, whence $\varphi$ is continuous in $\omega_T$. This yields $\jump{\varphi}|_F = 0$ for all internal faces F of $\omega_T$ and the middle term vanishes.
\end{proof}

Applying Lemma \ref{L:div-jumps}, we are led to define the elemental residual indicator 
\begin{equation}\label{E:def-eta-stokes}
\eta_\grid(\vec{u}_\grid,p_\grid,T)^2 :=
\eta_\grid^m(\vec{u}_\grid,p_\grid,T)^2 + 
h_T \!\! \sum_{F \subset \partial T\setminus \partial\Omega} \| \jump{\divo \vec{u}_\grid}\|_{L^2(F)}^2.
\end{equation}

Concatenating \eqref{E:equiv-stokes}, \eqref{E:bound-Resm-stokes} and Lemma \ref{L:div-jumps}, we fulfill Assumption \ref{A:equivalence-estimator-infsup} (equivalence of error and full estimator). The precise result is as follows.
\begin{proposition}[a posteriori error analysis for Stokes] \label{P:equiv-error-estimator-stokes} There exists constants $C_U \geq C_L >0$ such that
\[
C_L {\cal E}_\grid(\vec{u}_\grid, p_\grid,\vec{f}) \leq 
\| \vec{u} - \vec{u}_\grid \|_\V + \| p - p_\grid \|_\Q \leq 
C_U {\cal E}_\grid(\vec{u}_\grid, p_\grid,\vec{f}),
\]
where the full estimator is defined by
\[
{\cal E}_\grid(\vec{u}_\grid, p_\grid,\vec{f})^2 := \sum_{T\in \grid}
{\cal E}_\grid(\vec{u}_\grid, p_\grid,\vec{f},T)^2
\]
with \ 
$
{\cal E}_\grid(\vec{u}_\grid, p_\grid,\vec{f},T)^2 := \eta_\grid(\vec{u}_\grid,p_\grid,T)^2 + \osc_\grid(\vec{f},T)_{-1}^2
$ introduced in \eqref{E:def-eta-stokes} and \eqref{E:mod-res-est-stokes} and 
$\osc_\grid(f,T)_{-1}=\|f-P_\grid f\|_{H^{-1}(\omega_T)}$.
\endproof
\end{proposition}

Since the Stokes problem has constant coefficients but variable forcing, it is natural to resort to Algorithm \ref{A:one-step-switch} (\AFEMSW), the one-step \AFEM with switch, for its adaptive discretization. 
\AB{With respect to the functional setting
of Section~\ref{S:AFEMSW-infsup}, the ambient space ${\mathbb W}$ is $\V \times \Q$ and the data projection operator $P_\grid$ is  
$$
\boldsymbol{P}_\grid:= ((P_\grid)^d, \Pi_\grid^{n-1}) : {\mathbb W}^* \to  ({\mathbb F}_\grid)^d \times {\mathbb S}^{n-1,-1},
$$
where ${\mathbb F}_\grid$ is the scalar discrete space introduced in Definition \ref{D:discrete-functionals}, $P_\grid$ is here the scalar projection operator introduced in Definition \ref{D:Pgrid}, and $\Pi_\grid^{n-1}$ is the $L^2$-orthogonal projection upon ${\mathbb S}^{n-1,-1}$. Furthermore, the norm used to measure data perturbations is $\| (\boldsymbol{f},g) \|_{{\mathbb W}^*_\grid}^2 = \sum_{T \in \grid} \big( \| \boldsymbol{f}\|_{(H^{-1}(\omega_T))^d}^2 +\| g \|_{L^2(T)}^2 \big)$. 

It is easily seen that $\eta_\grid(\vec{u}_\grid,p_\grid,T)$ satisfies Assumptions \ref{A:Lipschitz-estimator} (Lipschitz continuity of estimator) and \ref{A:monotonicity-estimator} (monotonicity of estimator) as well as the hypotheses of  Proposition \ref{P:local-estimator-reduction} (estimator reduction under D\"orfler marking): the estimator is clearly Lipschitz continuous and monotone, and it satisfies condition \eqref{E:estim-red-abstract} since all its addends appear multiplied by a positive power of the meshsize. In addition, the oscillation $\osc_\grid((\boldsymbol{f},g))=\| (\boldsymbol{f},g) - \boldsymbol{P}_\grid (\boldsymbol{f},g) \|_{{\mathbb W}^*_\grid}$ fulfills Assumption \ref{A:monotonicity-oscillation} (quasi-monotonicity of oscillation). 
}

 Theorem \ref{T:lin-convergence-switch} (linear convergence for \AFEMSW) provides sufficient conditions for the linear convergence of the algorithm, and these conditions have been verified along the previous discussion. Therefore, we may summarize our findings in the following theorem.

\begin{theorem}[linear convergence for Stokes]\label{T:lin-conv-Sokes}
Consider the Galerkin discretization \eqref{E:weak-stokes-Gal} of the Stokes problem which uses Taylor-Hood elements of order $n\ge2$, and let the a posteriori estimator be given in Proposition \ref{P:equiv-error-estimator-stokes}. Then, Theorem \ref{T:Feischl-theorem} guarantees the linear convergence of Algorithm \ref{A:one-step-switch} (\rm \AFEMSW) applied to this problem, i.e., it holds for some $c> 0$ and $0 \leq \rho <1 $
$$
e_{j+i}\le c \rho^i e_j \qquad \forall i, j \in \N,
$$    
with  $e_j:= \| \nabla(\vec{u}-\vec{u}_j) \|_{\Omega} +
\| p - p_j \|_{\Omega}$. \endproof
\end{theorem}

In order to assess the optimality of the discretization, we specify the definition of approximation classes for the solution of the Stokes problem. Precisely, given $(\vec{v},q)\in\V\times\Q$ we let $\sigma_N(\vec{v},q)$ be the smallest approximation error incurred
on $(\vec{v},q)$ with elements in $\V_\grid \times \Q_\grid$ over meshes belonging to $\grids_N$:
\begin{equation}\label{E:sigmaN-stokes}
\sigma_N(\vec{v},q) :=  \inf_{\grid \in \grids_N} \inf_{(\vec{v}_\grid,q_\grid) \in \V_\grid \times \Q_\grid} (\|\nabla(\vec{v} - \vec{v}_\grid)\|_{\Omega}^2 +\|q-q_\grid  \|_{\Omega}^2 )^{1/2}.
\end{equation}
For $0<s\le n/d$, the class $\As=\As(\V\times\Q;\grid_0)$, relative to the
partition $\grid_0$ is the set of functions $(\vec{v},q)\in\V\times \Q$ such that
\begin{equation}\label{E:quasi-norm-stokes}
|(\vec{v},q)|_\As := \sup_{N\ge\#\grid_0} \big( N^s \sigma_N(\vec{v},q) \big)<\infty
\end{equation}
By adapting the arguments used in the proof of Theorem \ref{T:optimality-one-step} (rate optimality of one-step \AFEM{s}), we can prove the following result.
\begin{theorem}[rate optimality of \AFEMSW for Stokes]\label{T:optimality-one-step-Stokes}
Let the assumptions of Theorem \ref{T:lin-conv-Sokes} be valid. If $(\vec{u},p)\in\As$, then the sequence $\{\grid_k,\V_k,(\vec{u}_k,p_k)\}_{k\geq 0}$ generated by {\rm \AFEMSW} are such that
\begin{equation}\label{E:optimality-one-step-Stokes}
\|\nabla(\vec{u} - \vec{u}_k)\|_{L^2(\Omega)} +\|p-p_k  \|_{L^2(\Omega)} \lesssim |(\vec{u},p)|_{\As} (\#\grid_k)^{-s}, \quad  k\geq 0.
\end{equation}
\end{theorem}
\smallskip

%\todo[inline]{RHN: We could mention here the one-step AFEM for $\vec{f}$ with piecewise polynomial of degree $\le n$ and the one-step AFEM with switch for more general $\vec{f}$.}

\begin{remark}[limits of the analysis] Other inf-sup stable elements, such as the Mini element  or the Crouzeix-Raviart element (see e.g. \cite{BoffiBrezziFortin}), do not fit in the present setting of analysis, since their velocities contain elementwise bubble components (which are indeed responsible for the stability of the elements). Unfortunately, a bubble on an element does not restricts to two bubbles when the element is bisected, preventing  the nestedness condition $\V_{\grid} \subset \V_{\grid_*}$ to be satisfied when $\grid_*$ is a refinement of $\grid$.    
\end{remark}

%-----------------------------------------------------------------------------------------
\subsection{Mixed FEMs for scalar elliptic PDEs}\label{S:rates-mixed-fems}
%-----------------------------------------------------------------------------------------

% \todo[inline]{Should we discuss specific mixed FEMs (Raviart-Thomas and BDM) here or we present the 
% analysis abstractly? Perhaps we could do the latter here and discuss examples in Section \ref{S:aposteriori}.} 

The diffusion-reaction problem \eqref{strong-form} can be formulated in mixed form as follows:
\begin{equation}
\label{E:strong-form-2}
\begin{aligned}
\bA^{-1} \bsigma &= \nabla u  &&\text{ in }\Omega,
\\
 -\div \bsigma + c u &= f &&\text{ in }\Omega,
\\
 u &= 0 &&\text{on }\partial\Omega.
\end{aligned}
\end{equation}
Introducing the porosity matrix $\bK:=\bA^{-1}$, we assume hereafter that data 
\begin{equation}\label{E:defAc-mixed}
{\cal D} = (\bK,c,f) \in D(\Omega) := M(\alpha_1,\alpha_2) \times R(c_1,c_2) \times L^2(\Omega),
\end{equation}
where $M(\alpha_1,\alpha_2)$ and $R(c_1,c_2)$ are defined in \eqref{d:DiffusionSpace} and \eqref{d:ReactionSpace}, respectively. Note that the current parameters $\alpha_1,\alpha_2$ are
the reciprocals of $\alpha_2,\alpha_1$ in \eqref{d:DiffusionSpace}, but to avoid complicating the 
notation further, we relabel them hereafter.

\medskip\noindent
{\em Weak formulation.}
 To write the weak formulation of these equations, we introduce the Hilbert space
\begin{equation}\label{E:def-Hdiv}
H(\divo;\Omega):= \big\{\btau \in L^2(\Omega;\R^d) : \divo \btau \in L^2(\Omega)\big\}    
\end{equation}
equipped with the norm $\|\btau \|_{H(\div;\Omega)}^2 := \|\btau \|_{\Omega}^2 + \|\divo \btau \|_{\Omega}^2$. Then, we multiply the first equation in \eqref{E:strong-form-2} by $\btau \in H(\div;\Omega)$ and the second equation by $v \in L^2(\Omega)$, we integrate over $\Omega$ and apply partial integration to the term containing $\nabla u$, keeping into account the Dirichlet boundary condition. In this way, we obtain the following variational problem: find $(\bsigma,u) \in \V:= H(\div;\Omega) \times  L^2(\Omega)$ such that
\begin{equation}\label{E:weak-mixedDR}
\begin{alignedat}{3}
  &\int_\Omega \bK \bsigma \cdot \btau + \int_\Omega u \divo \btau  &&= 0 
&\qquad &\forall \, \btau \in H(\div;\Omega),
  \\
  & \int_\Omega v \divo \bsigma  - \int_\Omega c\, u \, v  &&= - \int_\Omega f \, v &&\forall \, v\in L^2(\Omega).
\end{alignedat}
\end{equation}
This can be written as: 
find $(\bsigma, u)\in V\times Q$ such that 
\begin{equation}\label{E:weak-mixedDR2}
\begin{alignedat}{3}
  &a[\bsigma,\btau] + b[u,\btau] &&= 0
&\qquad &\forall \, \bsigma\in V,
  \\
  &b[v,\bsigma] + m[u,v]&&= -\langle f,v \rangle  &&\forall \, v\in Q,
\end{alignedat}
\end{equation}
if we set $V:=H(\div;\Omega)$, $Q:=L^2(\Omega)$, and we introduce the continuous bilinear  forms 
$a\colon V\times V\to\R$, 
$b\colon Q\times V\to\R$ and 
$m\colon Q\times Q\to\R$ by
\begin{equation*}
     a[\bsigma,\btau] = \int_\Omega \bK \bsigma \cdot \btau\,, \qquad
     b[v,\btau] = \int_\Omega v \divo \btau\,, \qquad 
     m[u,v] = - \int_\Omega c\, u \, v \,,
\end{equation*}
and the linear form $\scp{f}{v} = \int_\Omega f v$.
An equivalent formulation, similar to \eqref{weak-stokes-B}, is as follows:
\begin{equation}\label{weak-Dirichlet}
  (\bsigma,u)\in  V\times Q: \quad
  \bilin{(\bsigma, u)}{(\btau,v)} =  -\scp{f}{v} \quad\forall \,
  (\btau,v)\in  V\times Q,
\end{equation}
with
\[
\bilin{(\bsigma, u)}{(\btau,v)} 
\definedas
a[\bsigma, \btau] + b[u,\btau]
+ b[v,\bsigma] + m[u,v].
\]

Formulation \eqref{E:weak-mixedDR2} is a generalization of the classical saddle point problem considered in Section \ref{S:inf-sup}, given by the presence of the third bilinear form $m$. According to \cite[Theorem 4.3.1]{BoffiBrezziFortin}, the well-posedness of such a problem can be derived from the three following conditions:
\begin{enumerate}%[label=(\roman*)]
    \item the form $a$ is coercive on $V_0=\{\btau \in V : b[v,\btau]=0 \text{ for all }v \in Q\}$, 
    \item the form $b$ satisfies an inf-sup condition on $V \times Q$,
    \item the form $m$ is non-positive on $Q$, i.e., $m[v,v]\le 0$ for all $v\in Q$.
\end{enumerate}
These conditions are easily checked for our mixed formulation of the Dirichlet problem.

\medskip\noindent
{\em Discretization.}
To define a finite element discretization of this problem, we consider partitions $\grid\in\grids$ obtained by conforming bisection refinements of an initial partition $\grid_0$, and let $V_\grid \subset V$ and $Q_\grid \subset Q$ be finite dimensional subspaces made of piecewise polynomial functions on $\grid$. 
Among the families of {\em uniformly inf-sup stable} finite element spaces for this problem, we consider the Raviart-Thomas-N\'ed\'elec family \cite{RaviartThomas, nedelec:80}, and the Brezzi-Douglas-Marini family \cite{BrezziDouglasMarini:85} on simplicial elements. They fit into the general definition
\begin{eqnarray*}
V_\grid &=&\{\btau \in H(\div;\Omega) : \  \btau|_{T} \in \vec{V}_T, \ T \in \grid\}, \\
Q_\grid &=&\{ q \in L^2(\Omega) : \ q|_{T} \in Q_T, \ T \in \grid\}.
\end{eqnarray*}
For the Raviart-Thomas-N\'ed\'elec (RTN) family we have
\begin{eqnarray*}
\vec{V}_T = (\mathbb{P}_{n-1}(T))^d +{\vec{x}} \, \mathbb{P}_{n-1}(T), \qquad Q_T= \mathbb{P}_{n-1}(T), \qquad \quad n \geq 1,
\end{eqnarray*}
where ${\vec{x}}=(x_1, \dots, x_d)$ is the coordinate vector, whereas for the Brezzi-Douglas-Marini (BDM) family we have
\begin{eqnarray*}
\vec{V}_T = (\mathbb{P}_n(T))^d , \qquad Q_T= \mathbb{P}_{n-1}(T), \qquad \quad n \geq 1.
\end{eqnarray*}
Note that for any face $F$ of the triangulation it holds $\btau|_F \cdot \vec{n}_F \in \mathbb{P}_{n-1}(F)$ for the RTN family, and $\btau|_F \cdot \vec{n}_F \in \mathbb{P}_{n}(F)$ for the BDM family; furthermore, $\div \vec{V}_T=Q_T$. We refer to \cite[Sects. 2.3.1 and 7.1.2]{BoffiBrezziFortin} for more details.

Due to the presence of variable data, it is natural to perform the adaptive discretization of the problem by adopting Algorithm \ref{algo:AFEMTS} (\AFEMTS), the two-step \AFEM.
The procedure $[\widehat{\grid},\widehat{\data}]=\DATA(\grid,\data, \tau)$ generates an admissible refinement $\widehat{\grid}$ of $\grid$ and discrete data 
\[
\widehat{\data} = (\widehat{\bK},\widehat{c}, \widehat{f}) \in \D_{\widehat{\grid}} :=
 \big[ {\mathbb S}^{n-1,-1}_{\wh{\grid}} \big]^{d\times d}
\times {\mathbb S}^{n-1,-1}_{\wh{\grid}} \times {\mathbb S}^{n-1,-1}_{\wh{\grid}} 
\]
 over $\widehat{\grid}$, such that
$\widehat{\bK} \in M(\widehat{\alpha}_1,\widehat{\alpha}_2)$, $\widehat{c} \in R(\widehat{c}_1,\widehat{c}_2)$ (see Sects. \ref{ss:approx_diff} and \ref{ss:approx_reac}), and  
% \begin{equation*}  %\label{E:bound-data-infsup-bis}
% \|\data - \widehat{\data}\|_{D(\Omega)} :=
% \| \bK-\widehat{\bK} \|_{L^\infty(\Omega)} +\| c - \widehat{c} \|_{L^\infty(\Omega)} + \|f - \widehat{f} \|_{L^2(\Omega)} \leq \Cdata\tau.  
% \end{equation*}
\begin{equation*}  %\label{E:bound-data-infsup-bis}
\|\data - \widehat{\data}\|_{\wh{D}(\Omega)}  \leq \Cdata\tau, 
\end{equation*}
where the space $\wh{D}(\Omega)$ is defined in \eqref{E:space-DO}.

The Galerkin discretization with these discrete data reads:
find $({\bsigma}_\grid,{u}_\grid)\in V_\grid\times Q_\grid$ such that 
\begin{equation}\label{E:weak-mixed-Gal}
\begin{alignedat}{3}
  &\widehat{a}[{\bsigma}_\grid,\btau] + b[{u}_\grid,\btau] &&= 0 
&\qquad &\forall \, \btau\in V_\grid,
  \\
  &b[v,{\bsigma}_\grid] \,+  \widehat{m}[{u}_\grid,v]&&= -\langle \widehat{f}, v \rangle &&\forall \, v\in Q_\grid,
\end{alignedat}
\end{equation}
with
\begin{equation*}
     \widehat{a}[\bsigma,\btau] = \int_\Omega \widehat{\bK} \bsigma \cdot \btau\,,  \qquad 
     \widehat{m}[u,v] = - \int_\Omega \widehat{c}\, u \, v \,, \qquad
     \langle \widehat{f}, v \rangle = \int_\Omega \widehat{f} \, v
     \,,
\end{equation*}
or equivalently,
\begin{equation*} ({\bsigma}_\grid,{u}_\grid)\in V_\grid\times Q_\grid: \quad
  \widehat{\cal B}[({\bsigma}_\grid,{u}_\grid),(\btau,v)] =  - \scp{\widehat{f}}{v} \quad\forall \,
  (\btau,v)\in  V_\grid\times Q_\grid.
\end{equation*}

\medskip\noindent
{\em A posteriori error estimator.}
Let us denote by $(\wh{\bsigma},\widehat{u})\in  V\times Q$ the exact solution of the perturbed problem
\begin{equation*}  
  \widehat{\cal B}[(\wh{\bsigma},\widehat{u}),(\btau,v)] =  - \scp{\widehat{f}}{v} \quad\forall \,
  (\btau,v)\in  V\times Q;
\end{equation*}
note that the forcing $\wh{f}$ appears with a negative sign.
Then, by continuity and uniform inf-sup stability of the form $\widehat{\cal B}$, we know that the error
\[
\| \wh{\bsigma} - {\bsigma}_\grid \|_{H(\divo;\Omega)} + \| \widehat{u}-{u}_\grid \|_{L^2(\Omega)}
\]
is equivalent to the quantity
\begin{equation*}
\sup_{(\btau,v)\in V\times Q}   \frac{\widehat{\cal B}[(\wh{\bsigma} -{\bsigma}_\grid, \widehat{u}-{u}_\grid),(\btau,v)]}{\| \btau \|_{H(\divo;\Omega)} + \| v \|_{L^2(\Omega)}}
= 
\sup_{(\btau,v)\in V\times Q}
\frac{
\langle \widehat{f}, v \rangle + \widehat{\cal B}[({\bsigma}_\grid ,{u}_\grid),(\btau,v)]}{\| \btau \|_{H(\divo;\Omega)} + \| v \|_{L^2(\Omega)}} .
\end{equation*}
By Galerkin orthogonality, the numerator is equal to  
\[
\langle \widehat{f}, v-v_\grid \rangle + \widehat{\cal B}[({\bsigma}_\grid,{u}_\grid),(\btau-\btau_\grid,v-v_\grid)] \qquad
\forall (\btau_\grid,v_\grid) \in 
V_\grid\times Q_\grid,
\]
that we now proceed to estimate. The term
\[
\widehat{\cal B}[({\bsigma}_\grid,{u}_\grid),(\btau-\btau_\grid,0)]
\]
can be analyzed as in  Carstensen \cite{Carstensen:97} (see also \cite[Section 4.8]{Verfuerth:13}), by resorting to a stable decomposition of $H(\divo;\Omega)$: precisely, given $\btau \in H(\divo;\Omega)$, there exist $\boldsymbol{\Phi} \in (H^1(\Omega))^d$ and ${\boldsymbol u}\in (H^1(\Omega))^d$ such that
\begin{equation}\label{E:decomp-btau}
\btau= \boldsymbol{\Phi} + \curl\, {\boldsymbol u}
\end{equation}
with $\| \boldsymbol{\Phi} \|_{(H^1(\Omega))^d} + \| {\boldsymbol u} \|_{(H^1(\Omega))^d} \lesssim \|\btau \|_{H(\divo;\Omega)}$ (see \cite[Section 5.1.3]{xu-chen-nochetto2009}). Note that if $\Omega$ is convex, then \eqref{E:decomp-btau} is the Helmholtz decomposition of $\btau$, with $\boldsymbol{\Phi} = \nabla G$ for some $G \in (H^2(\Omega))^d$. Using \eqref{E:decomp-btau} and a suitable choice of $\btau_\grid$, one can show that
\begin{equation}\label{E:split-mixed-1}
|\widehat{\cal B}[({\bsigma}_\grid,{u}_\grid),(\btau-\btau_\grid,0)]| \lesssim \eta_{\grid,1}((\bsigma_\grid,u_\grid)) \, \|\btau \|_{H(\divo;\Omega)}
\end{equation}
with 
\[
\eta_{\grid,1}((\bsigma_\grid,u_\grid))^2=\sum_{T \in \grid} \eta_{\grid,1}((\bsigma_\grid,u_\grid),T)^2
\]
and
\begin{equation}
\label{E:mod-res-mixed}
\begin{split}
\eta_{\grid,1}((\bsigma_\grid,u_\grid),T)^2
 :=& 
   \ h_T^2 \| \widehat{\bK} \bsigma_\grid - \nabla u_\grid \|_{L^2(T)}^2 
   +h_T^2 \| \curl (\widehat{\bK} \bsigma_\grid) \|_{L^2(T)}^2 \\
  &+  h_T \!\!\!\!\!\! \sum_{F \subset \partial T \setminus \partial\Omega}
  \| \jump{(\widehat{\bK} \bsigma_\grid)_t} \|_{L^2(F)}^2 %\\
  + h_T \!\!\!\!\!\! \sum_{F \subset \partial T \cap \partial\Omega}
  \| (\widehat{\bK} \bsigma_\grid)_t  \|_{L^2(F)}^2 ,
\end{split}
\end{equation}
where $\boldsymbol{\phi}_t= \boldsymbol{\phi} - (\boldsymbol{\phi}\cdot \vec{n}_F) \vec{n}_F$ denotes the tangential component of the vector field $\boldsymbol{\phi}$ on $F$.
On the other hand, the term
\[
\langle \widehat{f}, v-v_\grid \rangle + \widehat{\cal B}[(\bsigma_\grid,u_\grid),(0,v-v_\grid)] = \int_\Omega (\widehat{f} + \div \bsigma_\grid - \widehat{c} \, u_\grid)(v - v_\grid)
\]
can be bounded as follows. For any $T \in \grid$, let $\Pi_T = \Pi_T^{n-1}$ be the $L^2$-orthogonal projection upon $Q_T=\mathbb{P}_{n-1}(T)$,  and let us choose $(v_\grid)|_{T}= \Pi_T v$. Then, noticing that $\widehat{f} + \div \bsigma_\grid \in Q_T$, we have
\begin{equation}\label{E:mixed-identity}
\begin{split}
\int_\Omega (\widehat{f} + \div \bsigma_\grid - \widehat{c} \, u_\grid)(v - v_\grid) &= 
\sum_{T \in \grid} 
\int_T (\widehat{f} + \div \bsigma_\grid - \Pi_T(\widehat{c} \, u_\grid))(v - \Pi_T v) \\
& \qquad 
-\sum_{T \in \grid} 
\int_T (\widehat{c} \, u_\grid - \Pi_T(\widehat{c} \, u_\grid))(v - \Pi_T v) \\
& = 
-\sum_{T \in \grid} 
\int_T (\widehat{c} \, u_\grid - \Pi_T(\widehat{c} \, u_\grid)) v \,,
\end{split}
\end{equation}
whence
\[
| \langle \widehat{f}, v-v_\grid \rangle + \widehat{\cal B}[(\bsigma_\grid,u_\grid),(0,v-v_\grid)] | \leq \sum_{T \in \grid}
\|   
\widehat{c} \, u_\grid - \Pi_T(\widehat{c} \, u_\grid)
\|_{L^2(T)} \| v \|_{L^2(T)} \,.
\]
Conversely, it is easily checked that \eqref{E:mixed-identity} implies the bound
\[
\|   
\widehat{c} \, u_\grid - \Pi_T(\widehat{c} \, u_\grid)
\|_{L^2(T)} \lesssim
\| \div \bsigma - \div \bsigma_\grid \|_{L^2(T)} + \| \widehat{c} \|_{L^\infty(T)} \| \widehat{u}-u_\grid \|_{L^2(T)}.
\]

The choice $n=1$ yields $\widehat{c} \, u_\grid \in \mathbb{P}_0(T)$, hence, $\widehat{c} \, u_\grid - \Pi_T(\widehat{c} \, u_\grid)=0$ . For $n \geq 2$, we could define as a (squared) local error indicator the quantity
\[
\eta_{\grid,1}((\bsigma_\grid,u_\grid),T)^2 + 
\|   
\widehat{c} \, u_\grid - \Pi_T(\widehat{c} \, u_\grid)
\|_{L^2(T)}^2 \,,
\]
but the second addend is not guaranteed to reduce under refinement, since it is not scaled by a positive power of the meshsize.
However, there is an equivalent quantity which does reduce, as stated in the following result.
\begin{lemma}[equivalence of local error indicators]\label{L:equivalence-mixed} 
Assume $\widehat{c}, u_\grid \in \mathbb{P}_{n-1}(T)$, for $n \geq 2$. Let $\Pi_T^j$ be the $L^2$-orthogonal projection upon $\mathbb{P}_j(T)$. Then,
\begin{equation}\label{E:equivalence-mixed}
\|   
\widehat{c} \, u_\grid - \Pi_T(\widehat{c} \, u_\grid)
\|_{L^2(T)} \approx 
h_T \sum_{j=1}^n \|u_\grid - \Pi_T^{n-1-j} u_\grid \|_{L^2(T)} 
\| \nabla \widehat{c} - \Pi_T^{j-2}\nabla \widehat{c} \|_{L^\infty(T)},
\end{equation}
where the constants hidden in the symbol $\approx$ are independent of $\widehat{c}$, $u_\grid$, and $T$.
\end{lemma}
\begin{proof}
By the Bramble-Hilbert theorem,
\[
\| \widehat{c} \, u_\grid - \Pi_T(\widehat{c} \, u_\grid)
\|_{L^2(T)} \lesssim h_T^n | \widehat{c} \, u_\grid |_{H^n(T)}
\]
and
\[
| \widehat{c} \, u_\grid |_{H^n(T)} \lesssim \sum_{j=1}^{n-1} 
|\widehat{c} |_{W^{j}_{\infty}(T)} 
|u_\grid |_{H^{n-j}(T)}.
\]
Moreover,
\[
|\widehat{c} |_{W^{j}_{\infty}(T)} = 
|\nabla \widehat{c} |_{W^{j-1}_{\infty}(T)} =
|\nabla \widehat{c} - \Pi_T^{j-2}\nabla \widehat{c} |_{W^{j-1}_{\infty}(T)}
\]
and
\[
|u_\grid |_{H^{n-j}(T)} =
|u_\grid - \Pi_T^{n-1-j} u_\grid |_{H^{n-j}(T)}.
\]
Applying inverse estimates for semi-norms, we obtain the $\lesssim$ inequality in \eqref{E:equivalence-mixed}. 

To get the opposite inequality, it is enough to check that the vanishing of the left-hand side implies the vanishing of the right-hand side, since both quantities are defined on finite-dimensional spaces and their scaling with respect to the element size is the same. Now,
$\widehat{c} \, u_\grid = \Pi_T(\widehat{c} \, u_\grid)$ implies $\widehat{c} \, u_\grid \in \mathbb{P}_{n-1}(T)$. Let us assume that $u_\grid \in \mathbb{P}_{n-1-k}(T)$ for some $0 \leq k \leq n-1$, and consequently $\widehat{c} \in \mathbb{P}_k(T)$, i.e., $\nabla \widehat{c} \in \mathbb{P}_{k-1}(T)$. Then, $\Pi_T^{n-1-j}u_\grid = u_\grid$ for any $j \leq k$, hence, the corresponding differences in the summation on the right-hand side of \eqref{E:equivalence-mixed} vanish. Conversely, for $j>k$ we have $j-2 \geq k-1$, which implies 
$\Pi_T^{j-2}\nabla \widehat{c} = \nabla \widehat{c}$, that is the corresponding differences in the summation on the right-hand side vanish. In conclusion, all terms in the summation in \eqref{E:equivalence-mixed} vanish, and the thesis is proven.
\end{proof}

Summarizing, we have obtained the following result.
\begin{proposition}[a posteriori error estimator for mixed methods]\label{P: apost-mixed} For 
every $T\in\grid$, the following local quantity
\begin{equation}
\label{E:mod-res-mixed-2}
\begin{split}
\eta_{\grid}((\bsigma_\grid ,u_\grid)&,T)^2
 :=
   \ h_T^2 \| \widehat{\bK} \bsigma_\grid - \nabla u_\grid \|_{L^2(T)}^2 
   +h_T^2 \| \curl (\widehat{\bK} \bsigma_\grid) \|_{L^2(T)}^2 \\
  &+  h_T \!\!\!\!\!\! \sum_{F \subset \partial T \setminus \partial\Omega}
  \| [(\widehat{\bK}\bsigma_\grid)_t] \|_{L^2(F)}^2 %\\
  + h_T \!\!\!\!\!\! \sum_{F \subset \partial T \cap \partial\Omega}
  \| (\widehat{\bK} \bsigma_\grid)_t  \|_{L^2(F)}^2 \\
  &+ h_T^2 \sum_{j=1}^n \|u_\grid - \Pi_T^{n-1-j} u_\grid \|_{L^2(T)}^2 
\| \nabla \widehat{c} - \Pi_T^{j-2}\nabla \widehat{c} \|_{L^\infty(T)}^2
\end{split}
\end{equation}
is a (squared) a posteriori error indicator for the mixed problem \eqref{E:weak-mixed-Gal}, which gives rise to a global a posteriori error estimator $\eta_\grid({\bsigma}_\grid,{u}_\grid)$ that satisfies Assumption \ref{A:equivalence-estimator}.
\endproof
\end{proposition}

Finally, Assumption \ref{A:estimator-reduction} follows from Proposition \ref{P:local-estimator-reduction} 
 (estimator reduction under D\"orfler marking), since the estimator is clearly Lipschitz continuous and monotone, and it satisfies condition \eqref{E:estim-red-abstract} since all its addends are scaled by positive powers of the meshsize.

As a consequence, the \GALERKIN step in \AFEMTS converges linearly by Theorem \ref{T:Feischl-theorem}, and the number of sub-iterations in the $k$-th call to \GALERKIN is bounded by a constant $J$ independent of $k$ (Proposition \ref{P:cost-galerkin-infsup}). Furthermore, Theorem \ref{T:opt-AFEM-infsup} guarantees the quasi-optimality of the two-step \AFEM.

\begin{theorem}[quasi-optimality of \AFEMTS for mixed methods]\label{T:quasi-optim-mixed}
Let the exact solution $(\bsigma,u)$ of the mixed problem \eqref{E:weak-mixedDR} belong to the approximation class $\A_{s_u}(V;\grid_0)$ and let the data $(\bK,c,f)$ belong to the approximation class $\D_{s_\data}(D(\Omega);\grid_0)$.
Let Assumptions
 \ref{A:theta-no-osc} (marking parameter), \ref{A:size-omega} (size of $\omega$),
 and \ref{A:initial-labeling} (initial labeling)
 be valid. Consider the Galerkin discretization \eqref{E:weak-mixed-Gal}  based on one of the Raviart-Thomas-N\'ed\'elec or Brezzi-Douglas-Marini finite element pairs. There exists a constant $C_*$ independent of the exact solution $(\bsigma,u)$ and the data $\data=(\bK,c,f)$ such that the sequence $\{\big(\grid_k, \V_{\grid_k}\times \Q_{\grid_k},(\bsigma_{\grid_k},u_{\grid_k})\big)\}_{k\geq 0}$ produced by {\rm\AFEMTS} satisfies for $k \geq 0$
\[
\| \bsigma - \bsigma_{\grid_k} \|_{H(\divo;\Omega)}
+ \| u - u_{\grid_k} \|_{L^2(\Omega)} \leq C_* \left( |(\bsigma,u)|_{\A_{s_u}}^{1/s_u}+ |\data|_{\D_{s_\data}}^{1/s_\data} \right)^s \big(\# \grid_k\big)^{-s},
\]
with $0<s=\min(s_u,s_\data) \leq \frac{n}d$. \endproof
\end{theorem}

\begin{remark}
Another family of uniformly inf-sup stable spaces is the Brezzi-Douglas-Fortin-Marini's (see \cite{BreDougFortMar:87}), where $\V_T= \{ \btau \in (\mathbb{P}_{n}(T))^d : \btau \cdot \vec{n}_F \in \mathbb{P}_{n-1}(F)  \   \forall F \in {\cal F}\cap \partial T \}$ and $Q_T = \mathbb{P}_n(T)$, $n \geq 1$.  However, the imposed condition on the normal component of vector fields on each face of $\grid$ prevents the inclusion of $\V_{\grid}$ into $\V_{\grid_*}$ to hold if $\grid_*$ is a bisection refinement of $\grid$.   
\end{remark}

%-----------------------------------------------------------------------------------------

%-------------------------------------------------------------------------------------------------------
%-------------------------------------------------------------------------------------------------------
\subsection{Proof of Theorem \ref{T:bound-matrices}}\label{S:growth}
%-------------------------------------------------------------------------------------------------------
%-------------------------------------------------------------------------------------------------------
This section is devoted to establishing Theorem \ref{T:bound-matrices}, which in turn contributes with Corollary \ref{C:quasi-orthogonality} to the proof of Theorem \ref{T:Feischl-theorem}.

It is important to notice that the growth of $\|\Ub\|_2$ is dictated by the number of blocks $N$ rather than the actual dimension $n_N\gg N$ of $\Ub$.
Therefore, we again use the block notation from Section \ref{S:Quasi-Orthogonality}
$$
\Bb =\left(\Bb(i,j)\right)^N_{i,j=0}\in \R^{n_N\times n_N}, 
$$
with lower and upper triangular factors
$$
\Lb =\left(\Lb(i,j)\right)^N_{i,j=0}\in \R^{n_N\times n_N}, \qquad \Ub =\left(\Ub(i,j)\right)^N_{i,j=0}\in \R^{n_N\times n_N}.
$$
We also set
$$
 \qquad \Ab = \left(\Ab(i,j)\right)_{i,j=0}^N:=\sqrt{\alpha}\Bb
$$
for a suitable parameter $\alpha>0$ defined below.
%

%----------------------------------------------------------------------------------------
\subsubsection{Representation of block inverse matrices}
%----------------------------------------------------------------------------------------

We first show that it suffices to derive the estimates
\begin{equation}\label{E:estimate_of_Uinv}
 \|\Ub^{-1}\|_2\lesssim N^{-1/p}\qquad p>2 \,, 
\end{equation}
\begin{equation}\label{E:estimate_of_Utinv}
 \|\Ubt^{-1}\|_2\lesssim N^{-1/p}\qquad p>2  \,,
\end{equation}
\begin{equation}\label{E:estimateD}
    \|\Dib\|_2\lesssim 1,
\end{equation}
where $\Bb^T = \Lbt \Ubt$ is the normalized block triangular decomposition of $\Bb^T$, and $\Dib\in \R^{n_N\times n_N}$ stands for the block diagonal part of $\Ub$,

In fact, in view of property (P1) (continuity of $\B$), we see that
$$
\|\Bb\|_2 = \|\Bb^T\|_2\le \|\B\|,
$$
whence
$$
\Lb =\Bb \Ub^{-1} \Rightarrow \|\Lb\|_2 \le \|\Bb\|_2 \|\Ub^{-1}\|_2\le \|\B\| \|\Ub^{-1}\|_2 
$$
and, similarly, 
$$
\|\Lbt\|_2 \le \|\B\| \|\Ubt^{-1}\|_2 \,.
$$
On the other hand, from
$$
\Bb = \Lb \Dib \left(\Dib^{-1}\Ub\right)\quad \Rightarrow \quad \Bb^T =\left(\Ub^T \Dib^{-1}\right)\Dib \Lb^T
$$
we infer that
\begin{equation}\label{E:definitions_Lbt_Ubt}
    \begin{split}
    \Lbt = \Ub^T \Dib^{-1}\quad &\Rightarrow\quad \Ub = \Dib \Lbt^T\\
    \Ubt = \Dib \Lb^{T} \quad &\Rightarrow \quad \Lb^{-1}= \Dib \Ubt^{-T}\,,
    \end{split}
\end{equation}
which implies
\begin{equation*}
    \begin{split}
     \|\Ub\|_2 &\le \|\Dib\|_2\|\Lbt^T\|_2 ,\\
     \|\Lb^{-1}\|_2 &\le \|\Dib\|_2\|\Ubt^{-T}\|_2.
    \end{split}
\end{equation*}
Therefore, we can focus on proving \eqref{E:estimate_of_Uinv} and \eqref{E:estimateD}, since the proof of \eqref{E:estimate_of_Utinv} is identical to that of \eqref{E:estimate_of_Uinv}. We proceed in several steps. The most delicate estimate is \eqref{E:estimate_of_Uinv}.

\medskip
    \step{1} {\it $j$-th column of $\Ub^{-1}$.} To prove \eqref{E:estimate_of_Uinv}, it turns out to be convenient to get first an explicit expression for the $j$-th column of $\Ub^{-1}$. We achieve this next.
    \begin{lemma}[representation of the $j$-th column of $\Ub^{-1}$] We have
    \begin{equation}\label{E:representation_column_Uinv}
    \Ub^{-1}(i,j)= \Bb[j]^{-1}(i,j)\qquad\forall \; 0\le i\le j \le N.    
    \end{equation}
    \end{lemma}
    \begin{proof}
        We compute the $(i,j)$ block of $\Bb[j]^{-1} = \Ub[j]^{-1}\Lb[j]^{-1}$
        $$
        \Bb[j]^{-1}(i,j)= \sum^{j}_{k=0}\Ub^{-1}[j](i,k)\Lb^{-1}[j](k,j) = \Ub[j]^{-1}(i,j)
        $$
        because $\Lb^{-1}(k,j)=\boldsymbol{0}$ for $k<j$ and $\Lb^{-1}(j,j)=\Ib(j,j)$. Moreover, we claim that
        $$
        \Ub^{-1}[j](i,j) = \Ub^{-1}(i,j)\qquad i\le j
        $$
        because $\Ub^{-1}$ is block upper triangular. To see this, let
        $$
        \xb(:,j) = \Ub^{-1}(:,j)\in \R^{n_N\times d_j}
        $$
        be the $j$-th block column of $\Ub^{-1}$, which satisfies
        $$
        \Ub \xb(:,j) = \Ib(:,j)\in \R^{n_N\times d_j}.
        $$
        Since $\Ib(i,j)=\boldsymbol{0}$ for $i>j$ and $\Ub$ is block upper triangular, we have $\xb(i,j)=\boldsymbol{0}$ for $i>j$. Therefore, the matrix
        $$
        \xbt(:,j) = \left(\xb(i,j)\right)^j_{i=0}\in \R^{n_j \times d_j}
        $$
        with the first $j$ blocks of $\xb(:,j)$ satisfies the reduced system
        \begin{equation*}
            \begin{split}
            \Ub[j](j,j) \, \xbt(j,j)&= \Ib(j,j)\\
            \sum^j_{k=i}\Ub[j](i,k) \, \xbt(k,j)&=\boldsymbol{0} \qquad 0 \le i \le j.
            \end{split}
        \end{equation*}
        We thus deduce that $\xbt(:,j)= \Ub[j]^{-1}(:,j)$, as asserted.
    \end{proof}
    This lemma justifies dealing with $\Bb[j]^{-1}$.

    \medskip
    \step{2} {\it Representation of $\Bb[j]^{-1}$.} We resort to the Neumann series expansion. We first consider the uniform SPD matrix
    $$
    \Bb[j]\Bb[j]^T\in \R^{n_j\times n_j}
    $$
    for which there exists $\alpha>0$ such that
    $$
    \|\Ib[j] -\alpha \Bb[j]\Bb[j]^T\|_2 < 1
    $$
    uniformly in $j$. In fact, note that for $\xb\in\R^{n_j\times n_j}$
    $$
    \|\xb-\alpha \Bb[j]\Bb[j]^T\xb\|_2^2=\|\xb\|^2_2 - 2 \alpha \scp{\xb}{\Bb[j]\Bb[j]^T \xb} +\alpha^2\|\Bb[j] \Bb[j]^T \xb \|^2_2
    $$
    as well as
    $$
     \scp{\xb}{\Bb[j]\Bb[j]^T \xb}= \|\Bb[j]^T \xb\|_2^2\ge \beta^2\|\xb\|_2^2
    $$
    in view of property (P2) (discrete inf-sup) and \eqref{E:discrete-inf-sup}, and
    $$
    \| \Bb[j]\Bb[j]^T \xb\|_2\le \|\Bb[j]\|^2_2\|\xb\|_2\le \|\B\|^2 \|\xb\|_2. 
    $$
    Consequently
    $$
    \| \xb -\alpha \Bb[j]\Bb[j]^T \xb\|^2_2 \le \left(1 -2 \alpha \beta^2 + \alpha^2 \|\B\|^4\right)\|\xb\|^2_2.
    $$
    The quadratic polynomial in $\alpha$ on the right-hand side is minimized by $\alpha = \frac{\beta^2}{\|\B\|^4}$, and gives
    \begin{equation}\label{E:growth_definition_rho}
        \|\Ib[j]- \alpha \Bb[j]\Bb[j]^T\|^2_2 \le 1- \frac{\beta^4}{\|\B\|^4}=: \rho^2.
    \end{equation}
    From now on, we fix this value of $\alpha$ and assume the uniform bound \eqref{E:growth_definition_rho}. Let
    $$
    \Ab[j]= \sqrt{\alpha} \Bb[j]\;\in \R^{n_j\times n_j}
   $$
   be the $j$-th principal section of the matrix $\Ab$ introduced previously, and let
   \begin{equation}\label{E:growth_definition_D}
       \Db[j]:= \Ib[j] - \Ab[j] \Ab[j]^T \; \in \R^{n_j\times n_j}.
   \end{equation}
   \begin{lemma}[representation of ${\Bb[j]^{-1}}$] The following expression is valid
   \begin{equation}\label{E:representation_Binv}
       \Bb[j]^{-1}= \alpha \Bb[j]^T \sum_{m=0}^\infty \Db[j]^m \qquad \forall \, 0\le j \le N.
   \end{equation}
   \end{lemma}
   \begin{proof}
       Since $\|\Db[j]\|_2\le \rho <1$ according to \eqref{E:growth_definition_rho}, the Neumann series theorem guarantees that
       $$
       \Ab[j]\Ab[j]^T = \Ib[j] - \Db[j]
       $$
       is invertible and the inverse reads
       $$
       \Ab[j]^{-T}\Ab[j]^{-1} = \sum^\infty_{m=0} \Db[j]^m,
       $$
       where $\Db[j]^0 =\Ib[j]$. Multiplying on the left by $\Ab[j]^T$, we obtain
       $$
       \frac{1}{\sqrt{\alpha}}\Bb[j]^{-1} = \sqrt{\alpha}\Bb[j]\sum^\infty_{m=0}\Db[j]^m,
       $$
       which yields the assertion.
   \end{proof}

   \medskip
   \step{3} {\it Representation of $\Ub^{-1}$.} In order to obtain a representation of $\Ub^{-1}$, we now build on \eqref{E:representation_column_Uinv}, which gives a formula for the $j$-th column of $\Ub^{-1}$ in terms of $\Bb[j]^{-1}$, and \eqref{E:representation_Binv}, which provides a series representation of $\Bb[j]^{-1}$. To this end, we introduce the block upper triangular matrix $\Db_m \in \R^{n_N\times n_N}$ given by
   $$
   \Db_m (i,j):= 
   \begin{cases}
       \Db[j]^m(i,j) & i \le j\\
       0             & i > j
   \end{cases}
   $$
   for $m\ge 1$ and $\Db_0 = \Ib$. Hence,
    $$
   \Ub^{-1}(i,j) = 
   \begin{cases}
       \displaystyle{\alpha \left( \Bb[j]^T \sum^\infty_{m=0}\Db_m[j]\right)(i,j)} & i \le j\\
       0             & i > j.
   \end{cases}
   $$  
   To write this expression in compact form, it is convenient to introduce the {\it block upper triangular truncation operator} $\Ucal: \R ^{n_N\times n_N} \rightarrow \R ^{n_N\times n_N}$ defined by 
   $$
   \Ucal(\Mb)(i,j) := 
   \begin{cases}
       \Mb(i,j) &i\le j\\
       0 & i>j
   \end{cases} \qquad \forall \Mb \in \R ^{n_N\times n_N}.
   $$
   \begin{lemma}[representation of $\Ub^{-1}$]\label{L:infsup-representation-Uinv} There holds
   \begin{equation}\label{E:growth_representation_Uinv}
     \Ub^{-1} = \alpha \; \Ucal \left(\Bb^T \sum^\infty_{m=0}\Db_m\right).
   \end{equation}
   \end{lemma}
   \begin{proof} Since $\Db_m$ is block upper triangular for all $m\ge 0$, so is the series $\sum^\infty_{m=0}\Db_m$. It thus suffices to check that
   $$
   \left(\Bb^T \sum^\infty_{m=0}\Db_m\right)(i,j) =  \left(\Bb[j]^T \sum^\infty_{m=0}\Db_m[j]\right)(i,j) \qquad i\le j.
   $$
   This shows the desired relation \eqref{E:growth_representation_Uinv}.
   \end{proof}
   
   \step{4} {\it Recursion.} In order to estimate $\Db_m$, it is useful to relate $\Db_m$ with $\Db_{m-1}$. We start with a simple property of the operator $\Ucal$: for $\Ab, \Bb \in \R^{n\times n}$ and $1\le i\le j \le n$, there holds
   $$
   \big(\Ab\; \Ucal(\Bb)\big)_{i j} = \sum^n_{k=1}\Ab_{i k} \;\Ucal (\Bb)_{k j} = \sum^j_{k=1}\Ab_{i k}\;\Bb_{k j}= \big(\Ab[j] \Bb[j]\big)(i,j).
   $$
   \begin{lemma}[recursion]\label{L:infsup-recursion} The following is valid for all $m\ge 1$
       \begin{equation}\label{E:growth_recursion1}
           \Db_m = \Db_{m-1} - \Ucal \left(\Ab\; \Ucal(\Ab^T\Db_{m-1})\right),
       \end{equation}
       with $\Db_0=\Ib$. Therefore, the $j$-th column of $\Db_m$ reads
       \begin{equation}\label{E:growth_recursion2}
           \Db_m(0:j,j) = \Db[j]\Db_{m-1}(0:j,j) \qquad 0\le j \le N.
       \end{equation}
   \end{lemma}
   \begin{proof}
       Take first $m=1$ and apply the proceeding relation for $0\le i \le j \le N$ to obtain
       \begin{equation*}
           \begin{split}
            \left( \Db_{0} - \Ucal \left(\Ab\; \Ucal(\Ab^T\Db_{0})\right)\right)(i,j) &= \Ib(i,j) - \left(\Ab \; \Ucal (\Ab^T)\right)(i,j)\\
            &= \Ib(i,j) - \Ab[j] \Ab[j]^T(i,j)\\
            & = \Db[j](i,j) = \Db_1(i,j),
           \end{split}
       \end{equation*}
       in light of \eqref{E:growth_definition_D}.
       Take next $m>1$ and $0\le i \le j$ to arrive at
       \begin{equation*}
           \begin{aligned}
            &\left( \Db_{m-1} - \Ucal \left(\Ab\; \Ucal(\Ab^T\Db_{m-1})\right)\right)(i,j) \\
            & \qquad\quad = \Db_{m-1}(i,j) - \sum_{k, \ell =1}^j \Ab (i,k) \, \Ab^T(k,\ell) \, \Db_{m-1}(\ell,j)\\
            & \qquad\quad = \sum_{\ell=1}^j\left(\Ib[j]- \Ab[j]\Ab[j]^T\right) (i,\ell) \, \Db_{m-1}(\ell,j)\\
            & \qquad\quad = \sum_{\ell=1}^j \Db[j](i,\ell) \, \Db[j]^{m-1}(\ell,j)\\
            & \qquad\quad = \Db[j]^m(i,j) = \Db_m(i,j).
           \end{aligned}
       \end{equation*}
    This is the asserted equality \eqref{E:growth_recursion1}. The remaining relation \eqref{E:growth_recursion2} follows from the last equality upon realizing that
    $$
    \Db_m(0:j,j) = \left(\Db_m(i,j)\right)^j_{i=0} = \Db[j]\Db[j]^{m-1}(0:j,j)=\Db[j]\Db_{m-1}(0:j,j).
    $$
    This completes the proof.
   \end{proof}
   
   \step{5} {\it Schatten norms.} In the view Lemmas \ref{L:infsup-recursion} (recursion) and  \ref{L:infsup-representation-Uinv} (representation of $\Ub^{-1}$), we intend to estimate $\|\Ub^{-1}\|_2$ in terms of suitable norms of $\Db_m$ that depend on the number $N$ of blocks rather than the dimension $n_N$, because $n_N\gg N$. These special norms are called {\it block Schatten norms.}

   However, for the sake of clarity, we start with the definition and properties of the usual Schatten norms. They include the operator 2-norm, the Frobenius norm, and satisfy a H\"older inequality.
   \begin{definition}[Schatten norms] 
       Given $\Mb\in \R^{n \times n}$ let
       $$
       \sigma_1(\Mb)\ge \sigma_2(\Mb)\ge \dots \ge \sigma_n(\Mb)\ge 0
       $$
       be the singular values of $\Mb$. Given $1\le p\le \infty$, let the $p$-Schatten norm be
       $$
       |\Mb|_p:= \left(\sum_{m=1}^n \sigma_m(\Mb)^p\right)^{1/p}.
       $$ 
    \end{definition}
    \begin{remark} Note that if $p=\infty$ the Schatten norm reduces to the 2-norm, i.e. 
       $$
       |\Mb|_\infty = \sigma_1(\Mb)= \|\Mb\|_2,
       $$
       and if $p= 2$ it is equivalent to the Frobenius norm,
       $$
       |\Mb|_2 = \left(\sum^n_{m=1}\sigma_m(\Mb)^2\right)^{1/2} = \left(\sum^n_{m=1}\Mb_{i j}^2\right)^{1/2} =\|\Mb\|_F.
       $$
    \end{remark}
    
    We now list a number of useful properties of these norms.
    
    \begin{lemma}[properties of $|\cdot|_p$]
        The following properties hold for $1\le p\le\infty$:
        \begin{enumerate}
            \item \label{i:prop_1} $\sigma_i (\Mb^T \Mb) = \sigma_i(\Mb)^2\qquad\Rightarrow\qquad |\Mb^T\Mb|_p = |\Mb|^2_{2 p}$
            \item \label{i:prop_2} $\sigma_i (\Mb) = \sigma_i(\Mb^T)\qquad\Rightarrow\qquad |\Mb|_p = |\Mb^T|_{p}$\\
            \item \label{i:prop_3} H\"older inequality: for $\frac{1}{r} = \frac{1}{p} + \frac{1}{q}$, with $r,p,q\in [1, \infty]$,
            $$
            |\Mb_1 \;\Mb_2|_r\le |\Mb_1|_p \; |\Mb_2|_q
            $$
            \item \label{i:prop_4} $|\Ucal(\Mb)|_\infty \lesssim \log(n) \;|\Mb|_\infty$
            \item \label{i:prop_5} $|\Ucal(\Mb)|_{2^j}\le 2^{j-1}|\Mb|_{2^j}$.
        \end{enumerate}
    \end{lemma}
    \begin{remark} Properties \ref{i:prop_1} and \ref{i:prop_2} are trivial. We refer to \cite[Lemma XI.9.20]{DunfordSchwartz:88} for property \ref{i:prop_3}, to \cite[Eq (15)]{Bhatia:00} for property \ref{i:prop_4}, and to \cite{Davies:88} and \cite[Lemma 17]{Feischl:2022} for property \ref{i:prop_5}.
    \end{remark}
    
    To define the block Schatten norms, we consider the subspace $\Dcalb$ of $\R^{n_N \times(N+1)}$ of matrices of the form
    \begin{equation*}
    \Xb = \begin{bmatrix} \Xb_0 &      &      &     \\
                                &\Xb_1 &      &     \\
                                &      &\ddots&     \\
                                &      &      &\Xb_N
    \end{bmatrix}, \quad \Xb_j\in\R^{d_j}\quad 0\le j\le N,
    \end{equation*}
     or equivalently
     $$
     \Xb \in \Dcalb \quad \Longleftrightarrow \quad \Xb_{i j}=\boldsymbol{0}\qquad \forall \;i \neq n_{j-1}+1, \dots, n_j.
     $$
     We can represent $\Xb$ using block notation as follows:
     $$
     \Xb =\left(\Xb(i,j)\right)^N_{i,j = 0}, \qquad \Xb(i,j)\in\R^{d_j\times 1},
     $$
     where
     $$
     \Xb(i,j)=\begin{cases}
         \Xb_j\quad &i=j\\
         \boldsymbol{0}\quad &i\not = j.
     \end{cases}
     $$
     Given a block matrix $\Mb =\left(\Mb(i,j)\right)^N_{i,j=0}\in\R^{n_N\times n}$, we consider
     $$
     \Mb \Xb = \left(\Mb(i,j)\Xb_j\right)^N_{i,j=0}\in \R^{n_N\times(N+1)},
     $$
     namely the $j$-th block column of $\Mb \Xb$ is
     $$
     \left(\Mb(i,j)\Xb_j\right)^N_{i=0}\in \R^{n_N}.
     $$
     % \todo{RHN: Should we relate these matrices to the Galerkin problem and motivate the block Shatten norms.}
     \begin{definition}[block Schatten norms]\label{D:block-shatten} For $1\le p\le \infty$, let
     $$
     |\Mb|_{b,p}:= \sup_{\Xb \in\Dcalb,\;|\Xb|_\infty \le 1} |\Mb \Xb|_p \qquad \forall\; \Mb \in \R^{n_N\times n_N}.
     $$       
     \end{definition}
     Note the unusual norm $|\Xb|_\infty$ instead of $|\Xb|_p$ in this definition of operator norm $|\Mb|_{b,p}$. This choice is deliberate and will be useful later; see Remark \ref{usual-block}. We now list important properties of the block Schatten norms; see \cite[Lemmas 15, 16, 17]{Feischl:2022} for proofs.
     \begin{lemma}[properties of $|\cdot|_{b, p}$]\label{L:growth-properties_norm}
         The following properties hold for all $\Mb, \Mb_1, \Mb_2 \in \R^{n_N\times n_N}$ and $1\le p\le \infty$:
         \begin{enumerate}
             \item $|\Mb|_{b, p}\le (N+1)^{1/p}|\Mb|_\infty = (N+1)^{1/p}\|\Mb\|_2$
             \item $|\Mb_1 \Mb_2|_{b, p}\le |\Mb_1|_\infty |\Mb_2|_{b,p}$
             \item $|\Mb|_\infty\le |\Mb|_{b,p}$, \quad $|\Mb|_{b, \infty}=|\Mb|_\infty$
             \item  If $\Mb_1\in\R^{n_N\times n_N}$ is block triangular with $j$-th block column
             $$
             \Mb_1(0:j,j)= \Pb_j\Mb_2(0:j,j), \qquad \Pb_j \in \R^{n_j\times n_j}
             $$
             for $0\le j\le N$, then
             $$
             |\Mb_1|_{b,2}\le \max_{0\le j\le N}|\Pb_j|_\infty \, |\Mb_2|_{b,2}
             $$
             \item $|\Ucal(\Mb)|_{b, 2^k}\le 2^{k-1}|\Mb|_{b,2^k}$ \qquad $k=1,2$
             \item $|\Ucal(\Mb)|_\infty\le \left(\lceil\log_2(N)\rceil+1\right)|\Mb|_\infty$.
         \end{enumerate}
     \end{lemma}
     \begin{remark}\label{usual-block}
         To understand the significance of Definition \ref{D:block-shatten}, we examine the growth of the usual and block $p$-Schatten norm relative to the $\infty$-Schatten norm for $1\le p<\infty$. Given $\Mb\in\R^{n_N\times n_N}$, we have for the usual $p$-norm
         $$
         |\Mb|_p = \left(\sum^{n_N}_{i=1}\sigma_i(\Mb)^p\right)^{1/p}\le n_N^{1/p} \sigma_1(\Mb) = n_N^{1/p} |\Mb|_\infty= n_N^{1/p} \|\Mb\|_2,
         $$
         whereas for the block $p$-norm we get
         $$
         |\Mb|_{b,p}\le (N+1)^{1/p}|\Mb|_\infty = (N+1)^{1/p}\|\Mb\|_2,
         $$
         according to Lemma \ref{L:growth-properties_norm} (i). In fact, given $\Xb\in\Dcalb$ with $|\Xb|_\infty=\|\Xb\|_2=1$, we first note that
         $$
         |\Mb\Xb|_p = \left( \sum_{j=0}^N \sigma_j(\Mb\Xb)^p\right)^{\frac{1}{p}}
         \le \sigma_0(\Mb\Xb) (N+1)^{\frac{1}{p}} = \|\Mb\Xb\|_2 (N+1)^{\frac{1}{p}},
         $$
         and also that
         $$
         \|\Mb\Xb\|_2 = \sup_{\xb\in\R^{N+1}} \frac{\|\Mb\Xb\xb\|_2}{\|\xb\|_2}
         \le \|\Mb\|_2 \sup_{\xb\in\R^{N+1}}  \frac{\|\Xb\xb\|_2}{\|\xb\|_2} \le \|\Mb\|_2
         $$
         because $\|\Mb\|_2=|\Mb|_\infty=1$. On the one hand, this explains why it is convenient
         to have the norm $|\Xb|_\infty$ rather than $|\Xb|_p$ in Definition \ref{D:block-shatten}.
         On the other hand, this calculation reveals the key point that
         $$
         |\Mb|_{b,p}\ll |\Mb|_p
         $$
         because the growth of $|\Mb|_{b,p}$ is dictated by the number of blocks $N+1$ whereas that of $|\Mb|_p$ is proportional to the dimension $n_N$ of $\Mb$ and $n_N\gg N$. This property is essential in the estimate of $\|\Ub^{-1}\|_2$ below.
     \end{remark}

     \step{6}  {\it Estimate of $\|\Ub^{-1}\|_2$.} We are now in a position to prove the desired bound \eqref{E:estimate_of_Uinv}.
     \begin{proposition}[estimate of $\|\Ub^{-1}\|_2$]\label{P:estimate_norm_Uinv} Let $\Bb\in\R^{n_N\times n_N}$ be a block matrix such that
     $$
     \|\Bb\|_2\le \|\B\|, \qquad \max_{0\le j \le N}\|\Bb[j]^{-1}\|_2\le \frac{1}{\beta}.
     $$
    Then there exist constants $C_{LU}$ and $p>2$ such that the block upper triangular factor $\Ub$ of $\Bb$ satisfies
    \begin{equation}\label{E:growth-estimate_norm_Uinv}
        \|\Ub^{-1}\|_2\le C_{LU} N^{1/p}.
    \end{equation}
     \end{proposition}
     \begin{proof}
         We recall \eqref{E:growth_representation_Uinv} of Lemma \ref{L:infsup-representation-Uinv} (representation of $\Ub^{-1}$)
         $$
         \Ub^{-1} = \alpha\;\Ucal\left(\Bb^T \sum^\infty_{m=0}\Db_m\right),
         $$
         along with \eqref{E:growth_recursion1} of Lemma \ref{L:infsup-recursion} (recursion)
         $$
         \Db_m = \Db_{m-1} -\Ucal \left(\Ab\; \Ucal(\Ab^T\Db_{m-1})\right) \qquad m\ge 1
         $$
         and \eqref{E:growth_recursion2} of Lemma \ref{L:infsup-recursion}
         $$
          \Db(0:j,j) = \Db[j]\Db_{m-1}(0:j,j) \qquad 0\le j \le N,
         $$
         with $\Db_0=\Ib$. We use these expressions in conjunction with Lemma \ref{L:growth-properties_norm} (properties of $|\cdot|_{b,p}$) to prove \eqref{E:growth-estimate_norm_Uinv}. We proceed in several steps.
         \begin{enumerate}
             \item {\it Bound for $|\Db_m|_{b,2}$.} In light of \eqref{E:growth_definition_rho} and \eqref{E:growth_definition_D}
             $$
             |\Db[j]|_\infty =\|\Db[j]\|_2 \le \rho = \sqrt{1-\frac{\beta^4}{\|\B\|^4}}<1 \qquad 0\le j\le N.
             $$
             Applying Lemma \ref{L:growth-properties_norm} (iv) to $\Db_m$ yields
             \begin{equation*}
              \begin{split}
                  |\Db_m|_{b,2}&\le \max_{0\le j\le N}|\Db[j]|_\infty \;|\Db_{m-1}|_{b,2}\\
                  &\le \rho \;|\Db_{m-1}|_{b,2}\le \rho^m\; |\Ib|_{b,2}.
             \end{split} 
             \end{equation*}
             Recalling Lemma \ref{L:growth-properties_norm} (i)
             $$
             |\Ib|_{b,2}\le (N+1)^{1/2}\|\Ib\|_2 =(N+1)^{1/2},
             $$
             whence
             $$
             |\Db_m|_{b,2} \le \rho^m\;(N+1)^{1/2}.
             $$
             We observe that this bound is not good enough for our purposes because it scales like $N^{1/2}$ instead of $N^{1/p}$ for $p>2$. We next improve upon this bound.

             \medskip
             \item {\it Bound for $|\Db_{m}|_{b,4}$.} We take $k=2$ in Lemma \ref{L:growth-properties_norm} (v) and use the triangle inequality to arrive at
             \begin{equation*}
                 \begin{split}
                  |\Db_m|_{b,4}&\le |\Db_{m-1}|_{b,4} + |\Ucal \left(\Ab\; \Ucal(\Ab^T\Db_{m-1})\right)|_{b,4}\\
                  &\le|\Db_{m-1}|_{b,4} + 2 |\Ab\; \Ucal(\Ab^T\Db_{m-1})|_{b,4}.
                 \end{split}
             \end{equation*}
            We further apply Lemma \ref{L:growth-properties_norm} (ii) and (v) to obtain
            $$
            |\Ab\; \Ucal(\Ab^T\Db_{m-1})|_{b,4} \le 2 |\Ab|_\infty \;|\Ab^T\Db_{m-1}|_{b,4} \le 2 |\Ab|^2_\infty \;|\Db_{m-1}|_{b,4}.
            $$
            Therefore,
            $$
            |\Db_m|_{b,4}\le \left(1+ 4 |\Ab|^2_\infty\right)|\Db_{m-1}|_{b,4}
            $$
            but the prefactor on the right-hand side is greater than 1 and so not suitable for iteration. We still have
            $$
            |\Db_m|_{b,4}\le \left(1+ 4 |\Ab|^2_\infty\right)^m|\Ib|_{b,4}.
            $$

            \medskip
            \item {\it Bound for $|\Db_m|_\infty$.} We combine the estimates from steps 1 and 2 to exploit their relative merits. Recall from Lemma \ref{L:growth-properties_norm} (iii) that
            $$
            |\Db_m|_\infty\le|\Db_m|_{b,p}\qquad \forall 1\le p\le \infty.
            $$
            Take $p = 2, 4$ and $0< t<1$ to be chosen later, and write
            \begin{equation*}
                \begin{split}
                |\Db_m|_\infty &\le |\Db_m|^{1-t}_{b,2}|\Db_m|^t_{b,4}\\
                &\le \left[\rho^{1-t}\left(1+ 4|\Ab|^2_\infty\right)^t\right]^m |\Ib|^{1-t}_{b,2} |\Ib|^t_{b,4}.
                \end{split}
            \end{equation*}
            Consequently, there exists $0<t_0<1$ such that 
            $$
            q:=\rho^{1-t}\left(1+4|\Ab|^2_\infty\right)^t<1\qquad0<t<t_0
            $$
            and
            $$
            |\Db_m|_\infty \le q^m |\Ib|^{1-t}_{b,2}|\Ib|^t_{b,4}.
            $$
            We now estimate the two terms on the right-hand side relying on Lemma \ref{L:growth-properties_norm} (i), namely
            \begin{equation*}
                \begin{split}
                    |\Ib|_{b,2}&\le (N+1)^{1/2}\|\Ib\|_2 = (N+1)^{1/2},\\
                    |\Ib|_{b,4}&\le (N+1)^{1/4}\|\Ib\|_2 = (N+1)^{1/4}.
                \end{split}
            \end{equation*}
            We thus obtain
            $$
            |\Db_m|_\infty \le q^m (N+1)^{1/\wt{p}}
            $$
            with $\frac{1}{\wt{p}}= \frac{1-t}{2} + \frac{t}{4}<\frac{1}{2}$ for $0<t<t_0$.

            \medskip
            \item {\it Estimate of $\|\Ub^{-1}\|_2$.} Recalling the expression
            $$
            \Ub^{-1} = \alpha \; \Ucal \left(\Bb^T \sum^\infty_{m=0}\Db_m\right),
            $$
            and applying Lemma \ref{L:growth-properties_norm} (vi), (ii) and (iii), we see that
            \begin{equation*}
                \begin{split}
                \|\Ub^{-1}\|_2 = |\Ub^{-1}|_\infty &\lesssim \log(N)\;|\Bb|_\infty \sum^\infty_{m=0}|\Db_m|_\infty\\
                &\lesssim |\Bb|_\infty (N+1)^{1/\wt{p}}\log(N) \sum^\infty_{m=0}q^m\\
                &\lesssim \|\Bb\|_2 (N+1)^{1/\wt{p}} \log(N).
                \end{split}
            \end{equation*}
            Finally, for any $2< p < \wt{p}$, we can absorb the logarithm thereby getting
            $$
            \|\Ub^{-1}\|_2\lesssim\|\Bb\|_2 (N+1)^{1/p},
            $$
            which is the desired estimate \eqref{E:growth-estimate_norm_Uinv}.
         \end{enumerate}
         This concludes the proof.
     \end{proof}

    % \todo[inline]{Note that we obtain a factor $\|\Bb\|_2$ on the right-hand side instead of the factor
    % $\|\Bb\|_2^3$ derived in \cite{Feischl:2022}; see Theorem 18 in \cite{Feischl:2022} and bottom of 
    % p.37 of my handwritten notes. This is due to the fact thatwe have an expression for $\Ub^{-1}$ that
    % involves the sum $\sum_{m=0}^\infty \Db_m$ starting from $m=0$. Is it correct that our estimate be
    % independent of $|\Bb|_\infty = \|\Bb\|_2$? This is to be checked!}

     \medskip
     \step{7} {\it Estimate of block diagonal $\Dib$.} We recall that $\Dib = \text{diag }\Ub\in \R^{n_N\times n_N}$ is the block diagonal of $\Ub$.
     We consider the block partitioning of $\Bb[j]$
     $$
     \Bb[j] = \begin{bmatrix}
         \Bb[j-1]   &\Rb_1\\
         \Rb_2^T    &\Rb_3
     \end{bmatrix} \; \in \R^{n_j\times n_j}
     $$
     where
     \begin{equation*}
         \begin{split}
          \Rb_1 &= \Bb[j](1:j-1,j)\;\in \Rb^{n_{j-1}\times d_j},\\
          \Rb_2^T &= \Bb[j](j,1:j-1)\;\in \Rb^{d_j\times n_{j-1}},\\
          \Rb_3 &=\Bb[j](j,j) \;\in\R^{d_j\times d_j}.
         \end{split}
     \end{equation*}
     \begin{lemma}[bound of $\|\Dib\|_2$]\label{L:inf-sup_bound_normD} There holds
     \begin{equation}\label{E:inf-sup_definition_CD}
         \|\Dib\|_2\le \|\B\| + \frac{\|\B\|^2}{\beta} = C_D.
     \end{equation}
     \end{lemma}
     \begin{proof}
         Compute the $\Lb \Ub$ factorization of $\Bb[j]$
         $$
         \Bb[j]=\begin{bmatrix}
             \Ib[j-1]                &0\\
             \Rb^T_2 \Bb[j-1]^{-1}   &1
         \end{bmatrix}
         \begin{bmatrix}
             \Bb[j-1]     &\Rb_1 \\
             0            &\Rb_3 -\Rb_2^T\Bb[j-1]^{-1}\Rb_1
         \end{bmatrix}
         $$
         and realize that
         $$
         \Ub(j,j) = \Dib(j,j) = \Rb_3 -\Rb_2^T\Bb[j-1]^{-1}\Rb_1\;\in\R^{d_j\times d_j}. 
         $$
         Since
         $$
         |\Rb_i|_\infty =\|\Rb_i\|_2 \le \|\Bb[j]\|_2 \le \|\Bb\|_2 =\|\B\| \qquad i= 1, 2
         $$
         $$
         |\Rb_3|_\infty  = \|\Rb_3\|_2 \le \|\Bb[j]\|_2 \le \|\B\|,
         $$
         and
         $$
         |\Bb[j-1]^{-1}|_\infty = \|\Bb[j-1]\|_2\le \frac{1}{\beta},
         $$
         according to properties (P1) and (P2) of the bilinear form $\B$, we deduce
         $$
         |\Dib(j,j)|_\infty =\|\Dib(j,j)\|_2 \le \|\B\| + \frac{\|\B\|^2}{\beta}
         $$
         as asserted.
     \end{proof}

     \medskip
     \step{8} {\it Bound of $\Lb \Ub$ factors.} We are finally in the position to prove Theorem \ref{T:bound-matrices}.
         We combine Proposition \ref{P:estimate_norm_Uinv} (estimate of $\|\Ub^{-1}\|_2$) and $\Lb = \Bb \Ub^{-1}$ to obtain
         $$
         \|\Lb\|_2\le \|\Bb\|_2\|\Ub^{-1}\|_2\le \|\B\| C_{L U} N^{1/p}.
         $$
         Then, invoking \eqref{E:definitions_Lbt_Ubt} in conjunction with Proposition \ref{P:estimate_norm_Uinv} and Lemma \ref{L:inf-sup_bound_normD} (bound of $\|\Db\|_2$) as well as the bounds of $\|\Ub^{-1}\|_2$ and $\|\Lb\|_2$, yields
         \begin{equation}
             \begin{split}
             \|\Ub\|_2 &\le \|\Dib\|_2\|\Lbt^T\|_2 \le C_D\|\B\| C_{L U} N^{1/p}\\
             \|\Lb^{-1}\|_2 &\le \|\Dib\|_2\|\Ubt^{-T}\|_2\le C_D C_{L U}N^{1/p}
             \end{split}
         \end{equation}
         with $C_D$ being the constant in \eqref{E:inf-sup_definition_CD}. This completes the proof.    \endproof

%--------------------------------------------------------------------------------
\begin{comment}
\section{Goal Oriented Adaptivity }\label{S:goal-oriented}
 \rhn{(AB $\longrightarrow$ CC)}  \cite{Stevenson:2009,BeckerPraetorius:2023}

\begin{itemize}
\item
Drawbacks of the usual approach

\item
Mommer-Stevenson approach

\item

\end{itemize}

\input{goal.tex}

%--------------------------------------------------------------------------------
\section{AFEMs for Eigenvalue Problems }\label{S:eigenvalues}
\rhn{CC $\longrightarrow$ AB}  \cite{BoffiGallistl:2017}, Gallistl and Bonito-Demlow for cluster robust.

\begin{itemize}
\item

\item

\item

\end{itemize}

\input{eigen.tex}

\newpage
\ 
\newpage

%\input{data-tmp.tex}
\end{comment}
%--------------------------------------------------------------------------------
\appendix
%\section{Notations}
%\input{notations.tex}

%\indexprologue{Text explaining the index}
\printindex

%--------------------------------------------------------------------------------
%%% BIBLIOGRAPHY
\clearpage 
\addcontentsline{toc}{section}{References}
\bibliography{afem_bib}
\label{lastpage}
\end{document}